\documentclass[11pt]{amsbook}

\usepackage{amsthm, amsmath, amssymb, amscd}
\usepackage{graphicx}
\usepackage{color}
\usepackage{tikz-cd}
\usepackage{import}
\usepackage{microtype}

\usepackage[hidelinks,pagebackref]{hyperref}

\renewcommand*{\backref}[1]{}
\renewcommand*{\backrefalt}[4]{
  \ifcase #1
  [No citations.]
  \or [#2]
  \else [#2]
  \fi }

\newcommand\blfootnote[1]{%
  \begingroup
  \renewcommand\thefootnote{}\footnote{#1}%
  \addtocounter{footnote}{-1}%
  \endgroup
}

\newtheorem{theorem}{Theorem}[chapter]
\newtheorem{corollary}[theorem]{Corollary}
\newtheorem{lemma}[theorem]{Lemma}
\newtheorem{conjecture}[theorem]{Conjecture}
\newtheorem{proposition}[theorem]{Proposition}

\theoremstyle{definition}
\newtheorem{definition}[theorem]{Definition}
\newtheorem{question}[theorem]{Question}
\newtheorem{example}[theorem]{Example}
\newtheorem{notation}[theorem]{Notation}
\newtheorem{exercise}{Exercise}[chapter]

\theoremstyle{remark}
\newtheorem{remark}[theorem]{Remark}

\numberwithin{section}{chapter}
\numberwithin{equation}{chapter}
\numberwithin{figure}{chapter}    
\numberwithin{table}{chapter}
\numberwithin{theorem}{chapter}

\newcommand{\refthm}[1]{theorem~\ref{Thm:#1}}
\newcommand{\reflem}[1]{lemma~\ref{Lem:#1}}
\newcommand{\refprop}[1]{proposition~\ref{Prop:#1}}
\newcommand{\refcor}[1]{corollary~\ref{Cor:#1}}
\newcommand{\refrem}[1]{remark~\ref{Rem:#1}}

\newcommand{\refconj}[1]{conjecture~\ref{Conj:#1}}
\newcommand{\refeqn}[1]{equation~\eqref{Eqn:#1}}
\newcommand{\refitm}[1]{\eqref{Itm:#1}}
\newcommand{\refdef}[1]{definition~\ref{Def:#1}}
\newcommand{\refex}[1]{exercise~\ref{Ex:#1}}
\newcommand{\refsec}[1]{section~\ref{Sec:#1}}
\newcommand{\refsubsec}[1]{subsection~\ref{Subsec:#1}}
\newcommand{\reffig}[1]{figure~\ref{Fig:#1}}
\newcommand{\refchap}[1]{chapter~\ref{Chap:#1}}

\newcommand{\refexamp}[1]{example~\ref{Example:#1}}
\newcommand{\refexample}[1]{example~\ref{Example:#1}}

\newcommand{\reftable}[1]{table~\ref{Table:#1}}
\newcommand{\refnot}[1]{notation~\ref{Not:#1}}

\newcommand{\Refthm}[1]{Theorem~\ref{Thm:#1}}
\newcommand{\Reflem}[1]{Lemma~\ref{Lem:#1}}
\newcommand{\Refprop}[1]{Proposition~\ref{Prop:#1}}
\newcommand{\Refcor}[1]{Corollary~\ref{Cor:#1}}

\newcommand{\Refex}[1]{Exercise~\ref{Ex:#1}}

\newcommand{\Reffig}[1]{Figure~\ref{Fig:#1}}
\newcommand{\Refchap}[1]{Chapter~\ref{Chap:#1}}

\newcommand{\HH}{{\mathbb{H}}}
\newcommand{\RR}{{\mathbb{R}}}
\newcommand{\ZZ}{{\mathbb{Z}}}
\newcommand{\CC}{{\mathbb{C}}}
\newcommand{\QQ}{{\mathbb{Q}}}
\newcommand{\EE}{{\mathbb{E}}}

\newcommand{\calA}{\mathcal{A}}
\newcommand{\calT}{\mathcal{T}}
\newcommand{\calF}{\mathcal{F}}
\newcommand{\calL}{\mathcal{L}}
\newcommand{\calD}{\mathcal{D}}
\newcommand{\calH}{\mathcal{H}}

\newcommand{\Ahyp}{{A_{\rm Hyp}}}
\newcommand{\AhypT}{{A^{\calT}_{\rm Hyp}}}
\newcommand{\APSL}{{A_{\rm PSL}}}

\newcommand{\mat}[2][cccc]{\left(\begin{array}{#1} #2\\ \end{array}\right)}

\newcommand{\from}{\colon\thinspace} % As in ``f maps _from_ X _to_ Y''.

\newcommand{\bdy}{\partial}
\newcommand{\vol}{\operatorname{vol}}
\newcommand{\guts}{\operatorname{guts}}
\newcommand{\PSL}{\operatorname{PSL}}
\newcommand{\SL}{\operatorname{SL}}
\newcommand{\SO}{\operatorname{SO}}
\newcommand{\tr}{\operatorname{tr}}

\newcommand{\area}{\operatorname{area}}
\newcommand{\link}{\operatorname{link}}
\newcommand{\Id}{{\mathrm{Id}}}
\newcommand{\Arg}{{\mathrm{Arg}}}
\newcommand{\tw}{{\operatorname{tw}}}
\newcommand{\injrad}{{\operatorname{injrad}}}
\newcommand{\Isom}{{\operatorname{Isom}}}
\newcommand{\half}{{\frac{1}{2}}}
\newcommand{\voct}{{v_{\rm{oct}}}}
\newcommand{\vtet}{{v_{\rm tet}}}
\newcommand{\cut}{{\backslash \backslash}}

\DeclareMathOperator{\arctanh}{arctanh}

\makeindex

\begin{document}

\frontmatter

\title{Hyperbolic Knot Theory}
\author{Jessica S. Purcell}
\address{School of Mathematics \\
  9 Rainforest Walk, Room 401 \\
  Monash University, VIC 3800 \\
  Australia }
\email{jessica.purcell@monash.edu}
%% \thanks{Supported in part by the Australian Research Council.}
\date{\today}

%    The 2010 edition of the Mathematics Subject Classification is
%    the current definitive version.


\maketitle

\tableofcontents

\chapter*{Introduction}

Knots appear in scientific literature as early as 1771, in work of Vandermonde. In approximately 1833, Gauss developed the linking number of two knots, and his student Listing published work on alternating knots in 1847. Tait was one of the first to try to classify knots up to equivalence, creating the first knot tables in the 1870s and 1880s.
%% By 1900, Tait and Little had successfully enumerated all prime knots through ten crossings, although without proof, and with a duplication that was not discovered until work of Perko in 1974~\cite{Perko}. 
For more on the history of knots, see for example the detailed article by Epple~\cite{Epple}, or the survey articles by Przytycki~\cite{przytycki:knots} and Silver~\cite{silver:knothistory}.

Since the early work of Tait, knot theory has been influenced by and influential in the mathematical fields of topology, algebra, quantum field theory, and in geometry. There are several books that investigate knots from topological, algebraic, and quantum perspectives; some of my favorites are those of Rolfsen~\cite{rolfsen}, Burde and Zieschang~\cite{BurdeZieschang}, Murasugi~\cite{murasugi}, and Lickorish~\cite{lickorish:KnotTheory}. This book focuses on knots from a geometric perspective, particularly hyperbolic geometry, and overlaps more with books on hyperbolic geometry than knot theory, particularly in the early chapters that develop prerequisites in hyperbolic geometry. See, for example, \cite{benedetti-petronio, ratcliffe, thurston:book}.

The study of the geometry of knots, particularly hyperbolic geometry, began with work of Robert Riley in the 1970s~\cite{Riley:Fig8},
and developed further in the late 1970s and early 1980s, with work of William Thurston~\cite{thurston}.

By Thurston's work, a knot complement has one of three forms: Either it is a \emph{torus knot}, which can be drawn as an embedded curve on a Heegaard torus in the 3-sphere, or a \emph{satellite knot}, which lies in a tubular neighborhood of another knot, or it is \emph{hyperbolic}~\cite{thurston:bulletin}. Torus knots are relatively well-understood, and satellite knots are often studied by considering other knots. 
Hyperbolic knots, however, are not well-understood in general, and yet they are extremely common.  
For example, of all prime knots up to 16 crossings,
classified by Hoste, Thistlethwaite, and Weeks,
13 are torus knots, 20 are satellite knots, and the remaining 1,701,903 are hyperbolic~\cite{htw}. Of all prime knots up to 19 crossings, 15 are torus knots, 380 are satellite knots, and the remaining 352,151,858 are hyperbolic~\cite{Burton:KnotEnum}.

Moreover, 
if a knot complement admits a hyperbolic structure, then that structure is unique, by work of Mostow and Prasad in the 1970s~\cite{mostow, prasad}.
More carefully, Mostow showed that if there is an isomorphism between the fundamental groups of two closed hyperbolic 3-manifolds, then there is an isometry taking one to the other. Prasad extended this work to 3-manifolds with torus boundary, including knot complements. Thus if two hyperbolic knot complements have isomorphic fundamental group, then they have exactly the same
hyperbolic structure. Finally, Gordon and Luecke showed that two knot complements with the same fundamental group are equivalent~\cite{gordon-luecke}
(up to mirror reflection).

Thus a hyperbolic structure on a knot complement is a complete invariant of the knot. If we could completely understand hyperbolic structures on knot complements, we could completely classify hyperbolic knots. This book is an introduction to the mathematics involved.

\chapter*{Preface}

\section*{Why I wrote this book}

This book is an introduction to hyperbolic geometry in three dimensions, with motivations and examples coming from the field of knots. It is also an introduction to knot theory, with tools, techniques, and topics coming from geometry. As I write, I believe it is the only book that attempts to be both. 

To be clear, there are dozens of excellent books on knot theory, available from undergraduate to graduate levels, many of them classics that I learned from and continue to learn from. There are also several excellent books on hyperbolic geometry, particularly from the three-dimensional viewpoint. The aim of this book is to fill in a gap between them: to feature the contributions of hyperbolic geometry to knot theory, and the contributions of knot theory to hyperbolic geometry. It also aims to put techniques and tools from both fields into one place. 

In recent years, the field of hyperbolic 3-manifolds has matured, with many open conjectures resolved in the early 2000s. The result is that we now have better insight than ever into the structure of hyperbolic manifolds. This insight can be applied to broad classes of 3-manifolds, including many knot and link complements. On the other hand, the area of knot theory has also ballooned in recent years, with new tools arising from algebra, homology theory, quantum topology, representation theory, as well as geometry. As new knot and link invariants arise, and new applications of knot theory to other fields develop, it is natural to ask how such invariants interact. In particular, how do these invariants interact with hyperbolic geometry, which contains some of the strongest information on knots and links? There are many open questions and conjectures about the interaction of hyperbolic geometry with other knot invariants, and many mathematicians are interested in learning hyperbolic geometry specifically as it applies to knot theory. This book is a more direct introduction to the hyperbolic geometry of knots.

Hyperbolic geometry was first applied to the study of knots and their complements in the 1970s. 
Since then, hyperbolic geometry has played an important role in the classification of knots, with invariants such as volume and canonical decomposition developing directly from geometry.

However, the contribution of knot theory to hyperbolic geometry should not be understated. Complements of knots and links have been the playground of the 3-dimensional hyperbolic geometer for decades, aided by diagrams and topology, and by computational software such as SnapPea by Weeks, to find hyperbolic structures on knots. Many conjectures in hyperbolic geometry are based upon geometric properties that were first observed in knots. Many results in hyperbolic geometry have been proved first by restricting to families of knots, especially twist knots, two-bridge knots, and alternating knots, all of which feature prominently in this text.

This book is a hands-on introduction to this mixing of fields, geometry and knots.

\section*{How I structured the book}

The book starts with an introductory chapter giving basic definitions required from knot theory, and motivating some of the problems discussed in this book.

The first example of a hyperbolic knot, identified by Riley, is the unique prime knot with crossing number four, known as the figure-8 knot. In \refchap{Fig8Decomp}, we give an introduction to the complement 
of the figure-8 knot, and describe how to decompose it into two polyhedra. The exercises outline a generalization of this decomposition to all knots, and lead the reader through complications that arise when generalizing. This decomposition, particularly for the figure-8 knot, will then serve as a running example for later chapters.

In chapters two through six, we develop the basics of geometric structures on manifolds, particularly in dimensions two and three. Much of this material overlaps with other texts on hyperbolic geometry. Here, we try to keep our presentation heavily illustrated by examples, especially examples from knot theory. More specifically, \refchap{IntroHyp} gives an introduction to the hyperbolic plane and hyperbolic 3-space, and gives properties and examples of calculations that we will need in the text. It is purposely brief, as it is not meant to be a comprehensive introduction to these spaces, but only to equip the reader with the tools required to calculate and compute in hyperbolic geometry.
\Refchap{Geometric} introduces geometric structures on manifolds, and works through examples in two dimensions, including careful examples of the torus and the 3-punctured sphere. \Refchap{GluingCompleteness} returns to 3-manifolds and knots, building the first examples of hyperbolic structures on knot complements by way of triangulations. The chapter covers Thurston's gluing and completeness equations, again using the figure-8 knot as a running example. \Refchap{Margulis} delves a little more deeply into properties of hyperbolic isometries, with a main goal of proving the thick-thin decomposition of hyperbolic 3-manifolds. This decomposition implies that thin parts of hyperbolic 3-manifolds can always be identified with knots or links in some 3-dimensional space. Finally, in \refchap{CompletionDehnFilling}, incomplete structures on hyperbolic 3-manifolds are investigated carefully. The main result is that such structures can often be viewed as Dehn fillings of hyperbolic manifolds.

Chapters seven through twelve focus on families of knots and links that have been particularly amenable to study through hyperbolic geometry, and to tools used to study these knots and links, including tools coming from more general 3-manifold topology. \Refchap{TwistKnots}, just after the chapter on hyperbolic Dehn filling, discusses knots described by Dehn filling links in the 3-sphere; many of these links have very explicit hyperbolic geometry. This chapter explores consequences of Dehn fillings for these families. \Refchap{Essential} then provides an interlude, with results from 3-manifold topology, defining essential surfaces, normal surfaces, and returning to hyperbolic geometry via angle structures. \Refchap{AngleStruct} develops the powerful tool of angle structures and volumes of 3-manifolds. The main result in the chapter is a proof of the theorem of Casson and Rivin relating volumes of angle structures to hyperbolic geometry.
Angle structures have had great success as applied to the family of two-bridge knots, and this is the subject of \refchap{TwoBridge}. The chapter develops topological descriptions of the knots as gluings of tetrahedra, and works through a proof that these tetrahedra are geometric using the theorems of \refchap{AngleStruct}.
In \refchap{Alternating}, we study alternating links. This chapter gives a proof, using properties of these knots, of the theorem of Menasco that a prime alternating knot with more than one twist region is hyperbolic. \Refchap{Quasifuchsian} discusses the geometry of surfaces embedded in knot and link complements, including three and four punctured spheres, and checkerboard surfaces.

The final chapters, chapters thirteen through fifteen, explore some of the more important knot and link invariants arising from hyperbolic geometry. One of the most important geometric invariants of a hyperbolic knot is its volume, and \refchap{Volume} is devoted to volumes of knots and links. It contains methods to bound the volume of a knot.
\Refchap{Canonical} discusses the Ford domain and canonical polyhedral decomposition, also called the Epstein--Penner decomposition of a manifold. This decomposition provides a tool that can be used to identify when two 3-manifolds are isometric; for example it is used by the software SnapPea (and SnapPy). \Refchap{Character} gives a brief introduction to the overlap of hyperbolic geometry and algebraic geometry, introducing gluing and character varieties of knots, and the $A$-polynomial, which is a polynomial invariant directly related to the hyperbolic geometry of a knot.

\section*{Prerequisites and notes to students}
I have tried to keep prerequisites to a minimum. A basic course in topology is required, as well as some knowledge of basic algebraic topology, particularly the fundamental group and covering spaces. Occasionally, experience with Riemannian geometry will be helpful, but it is not required, with one exception: we assume standard results from a first course in Riemannian geometry in parts of \refchap{Volume}. We also occasionally assume basic results in differential topology, such as the fact that smooth manifolds admit tubular neighborhoods, and that submanifolds can be isotoped to meet transversely. 

Also, this book is written to be interactive, with examples and exercises. I hope you work through the examples as they are presented, and generalize them in exercises. Many important results are saved for exercises.

\section*{Acknowledgments}
The first form of this book appeared as lecture notes for a unit at Brigham Young University (BYU). The subject was inspired by my participation in a workshop on interactions between hyperbolic geometry, quantum topology, and number theory held at Columbia University in 2009. I have also given related graduate student workshops at Iowa in 2014, at Melbourne in 2016, and at Luminy in 2018. I thank the organizers of these workshops for inviting me, and various agencies for supporting the workshops, and for supporting fundamental research in mathematics. 

I learned much of the material in the first part of this book as a graduate student under the direction of Steve Kerckhoff, reading notes of William Thurston from the 1970s that were ghost-written by Kerckhoff and Bill Floyd~\cite{thurston}. Learning along with me were fellow graduate students David Futer and Henry Segerman. Many of their insights and elucidations helped me develop my own understanding; those insights are contained in this book, and I thank Steve, David, and Henry for them. I also owe thanks to Henry Segerman and Saul Schleimer for figures, particularly figures~\ref{Fig:SchleimerSegerman}, \ref{Fig:52SchleimerSegerman}, \ref{Fig:63SchleimerSegerman}, and~\ref{Fig:Fig8Fiber3D}.
Thanks to Saul Schleimer for figures~\ref{Fig:Fig8Fiber2D}, and to David Bachman, Saul Schleimer, and Henry Segerman for figure~\ref{Fig:Fig8Fiber3D2}. 
Discussions with David Futer about various parts of this book, especially two-bridge knots, have been invaluable. I also thank Jim Cannon, who attended my course on this material, and provided ideas for helping students get involved with exercises, and I thank Kenneth Perko, Abhijit Champanerkar, Ilya Kofman, Yi Wang, and Yair Minsky for feedback on drafts of the book. 

I owe the most thanks to the many people who have worked through various drafts and incarnations of this book, especially Emma Turner and Mark Meilstrup, who gave great feedback during the original 2010 BYU course, and Sophie Ham, Max Jolley, Josh Howie, Emily Thompson, John Stewart, and Ensil Kang, who read through many chapters carefully and helped me fix exposition and errors.
I also thank students who have worked through drafts of this book with others, including students at the University of Warwick, at Michigan State University, Temple University, and at Oklahoma State University. 

Remaining errors are, of course, my own. Please tell me about them.  

\aufm{Jessica S. Purcell}

\mainmatter

\setcounter{chapter}{-1} %% Introduction to knot theory is chapter 0
\chapter{A Brief Introduction to Hyperbolic Knots}\label{Chap:KnotIntro}
\blfootnote{Jessica S. Purcell, Hyperbolic Knot Theory}

This book gives an introduction to knots, links, and hyperbolic geometry. Before we begin, we need to define carefully what we mean by knots and links, and that is done in this chapter. We also introduce classical problems in knot theory, and problems motivated by geometry, especially hyperbolic geometry. This chapter is meant to motivate future chapters, and it has many references to content covered in more detail later in the book, where we address some of these problems. Many of the questions described in this chapter have partial answers, and many are still wide open.

\section{An introduction to knot theory}\label{Sec:KnotIntro}

The earliest study of knots seems to be by Gauss, Listing, and especially Tait, who published several papers on knot theory in the years 1876 through 1885. In a preface to his work on knot theory, republished in his 1898 Scientific papers~\cite{Tait:SciPapers}, Tait writes:
\begin{quote}
  ``The subject [knot theory] is a very much more difficult and intricate one than at first sight one is inclined to think, and I feel that I have not succeeded in catching the key-note.''
\end{quote}

Since Tait's work, advances in knot theory have come through applications of topology, algebra, and invariants arising in quantum field theory, but no single mathematical field has led to simple tools that apply to all knots. In other words, perhaps mathematicians still have not succeeded in catching the ``key-note.'' Perhaps there is no ``key-note'' in knot theory. However, there are definitely mathematical techniques that work well when applied to particular problems or particular families. This book introduces techniques arising from geometry.

%%%%%%%%%%%%%%%%%%%%%%%%%%%%%%%%%%%%%%%%%%%%%%%%%%%%%%%%%%%%%%%%%
\subsection{Basic terminology}

To begin, we need careful definitions of the objects involved.

\begin{definition}\label{Def:knot}
A \emph{knot}\index{knot, definition} $K \subset S^3$ is a subset of points homeomorphic to a circle $S^1$ under a piecewise linear (PL) homeomorphism. We may also
think of a knot as a PL embedding $K\from S^1 \to S^3$. We will use the same symbol $K$ to refer to the map and its image $K(S^1)$.

More generally, a \emph{link}\index{link, definition} is a subset of $S^3$ PL homeomorphic to a disjoint union of copies of $S^1$. Alternatively, we may think of a link as a PL embedding of a disjoint union of copies of $S^1$ into $S^3$.
\end{definition}

A PL homeomorphism of $S^1$ is one that takes $S^1$ to a finite number of linear segments. Restricting to such homeomorphisms allows us to assume that a knot $K\subset S^3$ has a regular tubular neighborhood,\index{tubular neighborhood} that is there is an embedding of a solid torus $S^1\times D^2$ into $S^3$ such that $S^1\times\{0\}$ maps to $K$. An embedding of $S^1$ into $S^3$ that cannot be made piecewise linear defines an object called a \emph{wild knot}.\index{wild knot} Wild knots may have very interesting geometry, but we will only be concerned with the classical knots of \refdef{knot} here.

In fact, rather than working with PL embeddings and homeomorphisms, we obtain the same results working with smooth ones. That is, we could require instead in \refdef{knot} that a knot be a smooth embedding of $S^1$ into $S^3$, and we obtain an equivalent theory. We will assume this fact throughout the book, working with both PL and smooth maps, with very little mention of this fact. 

\begin{definition}\label{Def:knotequiv}
  We will say that two knots (or links) $K_1$ and $K_2$ are equivalent if they are \emph{ambient isotopic}\index{ambient isotopic}, that is, if there is a (PL or smooth) homotopy $h\from S^3 \times [0,1] \to S^3$ such that $h(*, t) = h_t\from S^3 \to S^3$ is a homeomorphism for each $t$, and
  \[ h(K_1, 0) = h_0(K_1) = K_1 \quad \mbox{ and } \quad
  h(K_1, 1) = h_1(K_1) = K_2. \]
  Such a map $h$ is called an \emph{ambient isotopy}.\index{ambient isotopy}
\end{definition}

A PL (or smooth) embedding of $S^1$ into $S^3$ defines two 3-manifolds, one open and one compact, as in the following definition. 

\begin{definition}\label{Def:knotcomplement}
  For a knot $K$, let $N(K)$ denote an open regular neighborhood of $K$ in $S^3$. The \emph{knot exterior}\index{knot exterior, definition} is the manifold $S^3-N(K)$. Notice that it is a compact 3-manifold with boundary homeomorphic to a torus.

  The \emph{knot complement}\index{knot complement, definition} is the open manifold $S^3-K$, homeomorphic to the interior of $S^3-N(K)$.

  Similarly, if $L$ is a link the \emph{link exterior}\index{link exterior, definition} is $S^3-N(L)$, and the \emph{link complement}\index{link complement, definition} is $S^3-L$. 
\end{definition}

It was an open question for many years as to whether two knots with homeomorphic complements must be equivalent (up to reflection). This was proved in the affirmative by Gordon and Luecke in 1989~\cite{gordon-luecke}.

\begin{theorem}[Gordon--Luecke Theorem]\label{Thm:GordonLuecke}\index{Gordon--Luecke knot complement theorem}
  If two knots have complements that are homeomorphic by an orientation-preserving homeomorphism, then the knots are equivalent. 
\end{theorem}

The complement of a knot and the complement of its reflection are homeomorphic, by the orientation-reversing reflection homeomorphism. However, the knot itself may not be equivalent to its reflection. In fact, hyperbolic geometry tools do not distinguish knots and their reflections, and so we often only consider knots up to reflection in this book. If we disregard reflections, the Gordon--Luecke theorem states that knots are determined by their complements.

The same is not true for links. There are infinitely many inequivalent links whose complements are homeomorphic. However, the ways in which such links can be constructed are relatively well-understood; see, for example~\cite{gordon:links}.

\begin{definition}\label{Def:knot-diagram}
A \emph{knot diagram}\index{knot diagram, definition} (or \emph{link diagram})\index{link diagram, definition} is a 4-valent graph with over/under crossing information at each vertex. The diagram is embedded in a plane $S^2\subset S^3$ called the \emph{projection plane},\index{projection plane, definition} or \emph{plane of projection}.\index{plane of projection, definition}
\end{definition}

\begin{figure}[h]
\includegraphics{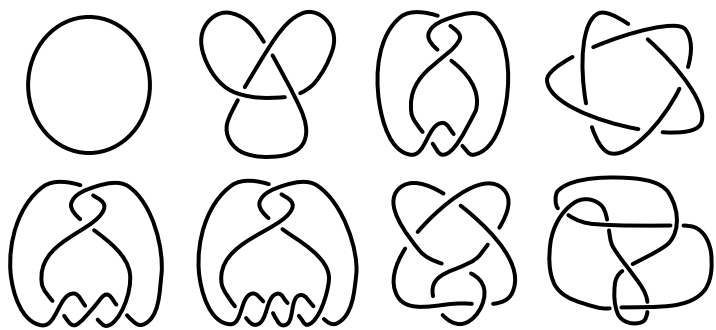}
\caption{Knots with at most six crossings.}
\label{Fig:KnotTable}
\end{figure}

\Reffig{KnotTable} shows diagrams of the eight knots with at most six crossings.
Classically, a knot has been described by a diagram.
Tait's works give many diagrams. In modern work, knots also appear without diagrams, for example when they arise as periodic orbits of a dynamical system~\cite{BirmanWilliams:Lorenz}, or from a gluing of polyhedra~\cite{CallahanDeanWeeks, ckp, ckm}. 

However, many open problems in knot theory still concern knot diagrams. 
One goal of \refchap{Fig8Decomp}, and then the next few chapters, is to give a method to pass from a knot or link \emph{diagram} to a topological and then geometric description of the knot or link \emph{complement}. That is, we start with a 4-valent graph describing a knot or link $K$, and obtain a mathematically rigorous decomposition of the 3-manifold $S^3-K$ into simple 3-dimensional pieces, which will be useful for applying tools from geometry and 3-manifold topology.

%%%%%%%%%%%%%%%%%%%%%%%%%%%%%%%%%%%%%%%%%%%%%%%%%%%%%%%%%%%%%%%%%
\section{Problems in knot theory}

There are many open problems in knot theory, and as new mathematical fields are brought to bear upon these problems, new questions and problems arise. This section gives a few highlights of the most classical problems, and also problems that seem most amenable to geometric techniques. Probably the most long-standing problem, and also one of the most broad, is the following.

\subsection{The classification problem}

When do two different descriptions of knots yield equivalent knots? When do they have homeomorphic complements?

\begin{figure}
  \includegraphics{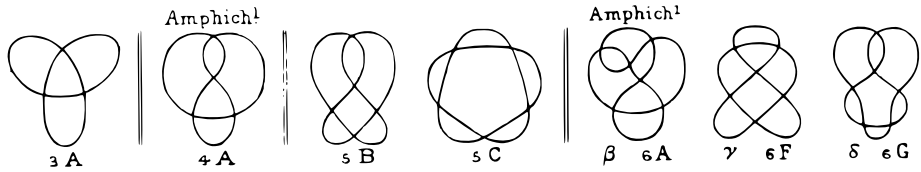}
  \caption{A very small portion of P.~Tait's 1884 tables of knot diagrams, from~\cite{Tait:OnKnotsII}. The original contains a full page with such diagrams, with additional pages of diagrams in~\cite{Tait:OnKnotsIII}.}
  \label{Fig:Tait}
\end{figure}

When the description of a knot is given by a diagram, this is the problem that Tait encountered while trying to list all knots with a fixed number of crossings. See \reffig{Tait}, which is modified from the 1884 paper~\cite{Tait:OnKnotsII}. 

\begin{figure}[h]
%% Accompanies image file 'F0-03-Nuga.eps' (pdf, eps, ps)
\begingroup%
  \makeatletter%
  \providecommand\color[2][]{%
    \errmessage{(Inkscape) Color is used for the text in Inkscape, but the package 'color.sty' is not loaded}%
    \renewcommand\color[2][]{}%
  }%
  \providecommand\transparent[1]{%
    \errmessage{(Inkscape) Transparency is used (non-zero) for the text in Inkscape, but the package 'transparent.sty' is not loaded}%
    \renewcommand\transparent[1]{}%
  }%
  \providecommand\rotatebox[2]{#2}%
  \newcommand*\fsize{\dimexpr\f@size pt\relax}%
  \newcommand*\lineheight[1]{\fontsize{\fsize}{#1\fsize}\selectfont}%
  \ifx\svgwidth\undefined%
    \setlength{\unitlength}{280.02919579bp}%
    \ifx\svgscale\undefined%
      \relax%
    \else%
      \setlength{\unitlength}{\unitlength * \real{\svgscale}}%
    \fi%
  \else%
    \setlength{\unitlength}{\svgwidth}%
  \fi%
  \global\let\svgwidth\undefined%
  \global\let\svgscale\undefined%
  \makeatother%
  \begin{picture}(1,0.24352916)%
    \lineheight{1}%
    \setlength\tabcolsep{0pt}%
    \put(0,0){\includegraphics[width=\unitlength]{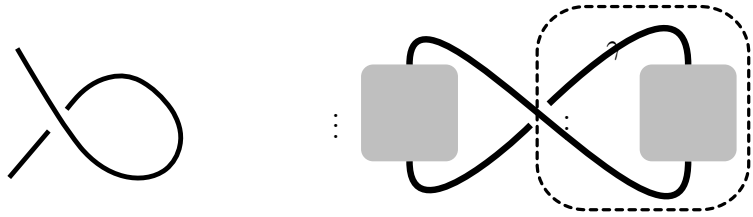}}%
    \put(0.80610283,0.21604877){\color[rgb]{0,0,0}\makebox(0,0)[lt]{\lineheight{1.25}\smash{\begin{tabular}[t]{l}$\gamma$\end{tabular}}}}%
    \put(0.43640286,0.10231692){\color[rgb]{0,0,0}\makebox(0,0)[lt]{\lineheight{1.25}\smash{\begin{tabular}[t]{l}$\vdots$\end{tabular}}}}%
    \put(0.74733461,0.09950607){\color[rgb]{0,0,0}\makebox(0,0)[lt]{\lineheight{1.25}\smash{\begin{tabular}[t]{l}$\vdots$\end{tabular}}}}%
  \end{picture}%
\endgroup%

  \caption{On the left, a nugatory crossing. On the right, a more general reducible crossing.}
  \label{Fig:Nugatory}
\end{figure}

There are a few moves that can be performed on a diagram that do not change the equivalence class of the underlying knot. For example, if the diagram contains a single crossing that forms a loop, as shown on the left of \reffig{Nugatory}, that loop can be untwisted to simplify the diagram.

\begin{definition}\label{Def:Nugatory}
  A single crossing forming a loop, as on the left of \reffig{Nugatory}, is called a \emph{nugatory crossing}\index{nugatory crossing}.

  More generally, a \emph{reducible crossing}\index{reducible crossing} is a crossing through which we may draw a circle $\gamma$ on the plane of projection such that $\gamma$ meets the diagram only at one point, at the crossing. See \reffig{Nugatory}, right.

  A diagram is \emph{reduced}\index{reduced diagram} if it contains no reducible crossings.
\end{definition}

Note that reducible crossings can be removed by isotoping the diagram. We typically will assume that our knot diagrams are reduced.

\begin{figure}[h]
  \includegraphics{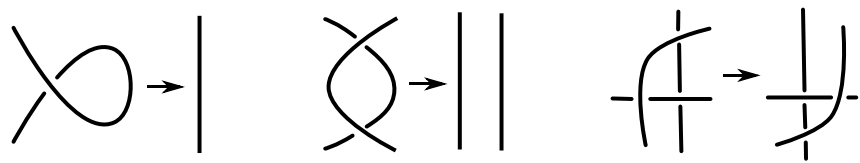}
  \caption{Three Reidemeister moves do not change knot equivalence.}
  \label{Fig:Reidemeister}
\end{figure}

There are other well-known moves to change a diagram into an equivalent diagram. These include the three moves shown in \reffig{Reidemeister}, called Reidemeister moves.\index{Reidemeister moves}

The Reidemeister moves appear in work of Maxwell in the 1800s (see, for example,~\cite{Epple}). In the 1920s, Reidemeister~\cite{Reidemeister} and Alexander and Briggs~\cite{AlexanderBriggs} independently gave rigorous proofs that two equivalent diagrams can always be related by a sequence of such moves.

\begin{figure}[h]
  \includegraphics{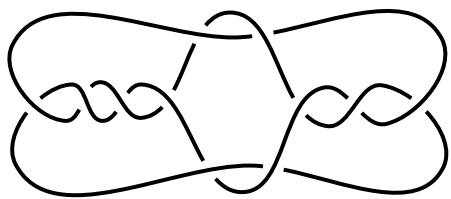}
  \caption{This diagram of the unknot was discovered in 1934 by Goeritz.}
  \label{Fig:Goeritz}
\end{figure}

The \emph{crossing number}\index{crossing number} of a knot is the minimal number of crossings in all diagrams of the knot. A minimal crossing diagram will necessarily be reduced. However, a reduced diagram is not necessarily a minimal crossing diagram. For example, \reffig{Goeritz} shows the reduced diagram of a knot that can, with a little work, be simplified to the \emph{unknot}\index{unknot}, i.e.\ the simple circle with no crossings. This diagram was discovered by Goeritz in 1934~\cite{Goeritz}. 

In fact, the diagram of \reffig{Goeritz} is an example of a knot diagram that cannot be simplified by Reidemeister moves without first increasing the number of crossings of the diagram. 

\begin{figure}[h]
  \includegraphics{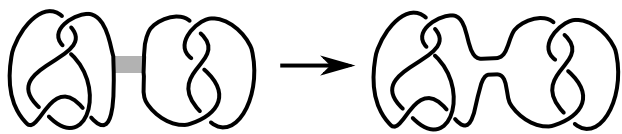}
  \caption{The knot sum of two knots.\index{knot sum}}
  \label{Fig:KnotSum}
\end{figure}

In addition to attempting to remove crossings, other moves can be performed on diagrams to simplify the classification problem. For example, there is a way of joining two simple diagrams into one more complicated diagram, shown in \reffig{KnotSum}.

Starting with two diagrams side-by-side, take a rectangle embedded in the plane of projection that has one side on one diagram, avoiding crossings, an opposite side on the other diagram, again avoiding crossings, and the final two sides disjoint from the two diagrams. Form the new diagram by removing the two edges of the rectangle that lie on the knots, and joining the knots along the two opposite sides of the rectangle. The resulting knot is called the \emph{knot sum}\index{knot sum}. It is also sometimes called the \emph{connected sum of the knots}\index{connected sum}. 

Given a knot sum of two knot diagrams, consider the embedded curve $\gamma$ in the plane of projection of the diagram that encircles exactly one of the original diagrams, cutting through the rectangle in the definition of the knot sum. This curve $\gamma$ meets the diagram of the knot sum in exactly two points, and it bounds disks on both sides (thinking of the projection plane as $S^2\subset S^3$), and both discs contain crossings.
We say that a diagram is \emph{prime}\index{prime}\index{prime!diagram} if no such curve $\gamma$ exists. That is, a knot or link diagram is \emph{prime} if, for every simple closed curve $\gamma$ in the plane of projection, if $\gamma$ meets the knot exactly twice transversely away from crossings, then $\gamma$  bounds a region of the diagram with no crossings. 

Curves such as $\gamma$ above detect knot sums. When knots are classified by diagram, listed according to crossing number, typically only prime diagrams are included.

The problem of listing all knots by crossing number, without duplicates, is a difficult one. There are 1,701,936 prime knots with at most 16 crossings, classified by Hoste, Thistlethwaite, and Weeks in 1998~\cite{htw}. More recently, Burton classified 352,152,252 prime knots up to 19 crossings~\cite{Burton:KnotEnum}. These knots can be downloaded with the 3-manifold software Regina~\cite{regina}. In both instances, the knots are only classified up to reflection in the plane of projection. 

\begin{definition}\label{Def:knotinvt}
A \emph{knot invariant}\index{knot invariant, definition} is a function from the set of knots to some other set whose value depends only on the equivalence class of the knot. A \emph{link invariant}\index{link invariant} is defined similarly.
\end{definition}

The crossing number of a knot is an example of a knot invariant.

Knot and link invariants are used to prove that two knots or links are distinct, or to measure the complexity of the link in various ways. We will revisit examples of knot invariants below, particularly geometric ones.

Notice that the number of knots with a given crossing number grows very rapidly. There does not seem to be a natural way of enumerating knots within a fixed class of crossing number. And while the crossing number was one of the first knot invariants to be studied by knot theorists, it does not seem to relate well to other knot invariants, particularly those that arise in geometry. For these reasons and others, other ways of classifying knots have arisen over the years, which we will discuss further below.

In this book we will apply geometry to the problem of the classification of knots. It has been known since the early 1980s, due to work of Thurston~\cite{thurston:bulletin}, that the complement of a knot decomposes into pieces, each admitting a 3-dimensional geometry. By using geometric properties of knot complements, we can often distinguish knots. This brings us to the second problem in knot theory that we discuss here.

\subsection{The problem of determining geometry of the complement}

Briefly, the complement of a knot is \emph{hyperbolic}\index{hyperbolic knot or link} if and only if it admits a complete metric with all sectional curvatures equal to $-1$. We will give other equivalent definitions of hyperbolic knots in later chapters, which will often be more useful for calculations, computations, and examples.

For now, it is known that when a knot complement is hyperbolic, its hyperbolic metric is unique. That is, hyperbolic knot complements that are homeomorphic must also be isometric under any hyperbolic metrics placed upon their complements.
Moreover, a large number of knots are hyperbolic, and many that are not hyperbolic decompose into hyperbolic pieces.

More precisely, consider the following families of knots.

\begin{figure}[h]
  \includegraphics{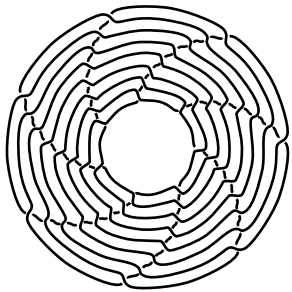}
  \caption{A torus knot}
  \label{Fig:TorusKnot0}
\end{figure}

\begin{definition}\label{Def:TorusKnot0}
  A \emph{torus knot}\index{torus knot} is a knot that can be embedded (without crossings) on the surface of an unknotted torus in $S^3$. See \reffig{TorusKnot0}. 
\end{definition}

By an unknotted torus, we mean the neighborhood of an unknot in $S^3$, with no crossings. 

\begin{figure}[h]
  \includegraphics{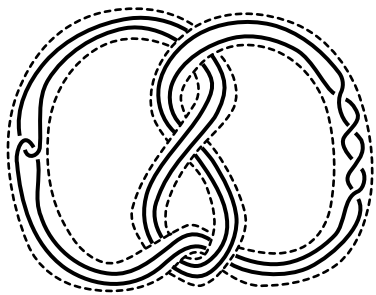}
  \caption{An example of a satellite knot. The dotted line forms the boundary of a neighborhood of a different knot, and the satellite lives inside that neighborhood.}
  \label{Fig:Satellite0}
\end{figure}

\begin{definition}\label{Def:Satellite0}
  A \emph{satellite knot}\index{satellite knot} is a knot that can be embedded in a regular neighborhood of another knot in $S^3$. See \reffig{Satellite0}.
\end{definition}

The complement of a torus knot admits a 3-dimensional geometry that is not hyperbolic, due to work of Thurston~\cite{thurston:bulletin}. He also showed that the complement of a satellite knot cannot be hyperbolic, but can be cut along a torus to decompose into pieces that admit 3-dimensional geometry, which could possibly be hyperbolic. For example, the knot complement in \reffig{Satellite0} can be cut along the dashed solid torus into two hyperbolic pieces, as we will see later in this book. 

Thurston showed that every knot in $S^3$ that is neither a torus knot nor a satellite knot must have hyperbolic complement~\cite{thurston:bulletin}. 

Thus hyperbolic geometry can be a useful tool in the classification problem of knots --- in theory.

In practice, we need tools and techniques to determine when a knot complement is hyperbolic. For example, if a knot is given by a messy diagram, how does one determine whether or not it is equivalent to a torus or satellite knot? How can we determine whether its complement is hyperbolic? And if it is hyperbolic, how can we find a hyperbolic metric?

Thurston outlined a procedure for finding a hyperbolic metric using the diagram of the figure-8 knot in his 1979 lecture notes~\cite{thurston}. This process was generalized by others, for example~\cite{menasco:links}, and even made algorithmic, in Weeks' 1985 PhD thesis~\cite{Weeks:Thesis}. There is now software that determines, given a knot diagram, whether or not the knot complement is hyperbolic. This is the computer program SnapPy, which is freely available~\cite{SnapPy}.\index{SnapPy}

Indeed, using computational tools, Burton has determined that of all prime knots with up to 19 crossings, $352,151,858$ are hyperbolic, and only $395$ are not hyperbolic~\cite{Burton:KnotEnum}. These are split into $14$ torus knots and $380$ satellite knots.

The next four chapters of this book concern the problem of determining a hyperbolic metric on a knot complement. We will step carefully through the necessary definitions and procedures, using Thurston's decomposition of the figure-8 knot complement as an example. This will give our first potential method to find a hyperbolic metric.

Chapters~\ref{Chap:Margulis} and~\ref{Chap:CompletionDehnFilling} give additional methods and tools from hyperbolic geometry to find or deform a hyperbolic metric. These first six chapters form the foundation required to discuss hyperbolic geometry and knots in more detail. 

Of course, these chapters require some work. The fact that software exists that can compute hyperbolic geometry of knots begs the question, why work through such computations by hand at all? Why not just work with the computer? There are many reasons, related to additional open problems. One reason is the next problem.

\subsection{The problem of determining geometry for families of knots.}
A computer program computes hyperbolic geometry for one knot at a time, or for a finite number of knots. But what can be said about infinite families of knots? For example, how does one determine the hyperbolic geometry of knots with descriptions given by infinite classes of diagrams? If two knots in a family are ``similar'' is their geometry also similar?

Potential answers to such questions seem to depend very heavily on the family of knots given. For example, for fixed $c$, it does not seem to be the case that the (finite) family of knots with crossing number $c$ have very similar hyperbolic geometry. 

On the other hand, certain infinite families of knots do exist with very similar hyperbolic geometry, and others at least seem to have geometry that reflects properties of the diagrams. We will discuss such knots and their properties, for example in chapters~\ref{Chap:TwistKnots}, \ref{Chap:TwoBridge}, and~\ref{Chap:Alternating}, with careful proofs. For now, we will present a definition of one such family. 

\begin{definition}\label{Def:Bigon}
  A \emph{bigon}\index{bigon} is a region of a graph bounded by exactly two edges and exactly two vertices. 
\end{definition}

\begin{figure}[h]
  \includegraphics{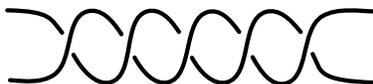}
  \caption{A twist region of a diagram}
  \label{Fig:TwistRegion}
\end{figure}

For example, \reffig{TwistRegion} shows several bigons connected end-to-end in a portion of a diagram graph of a knot. 

\begin{definition}\label{Def:TwistRegion}
  A \emph{twist region}\index{twist region} of a diagram of a knot is a maximal portion of the knot diagram where two strands twist around each other, as in \reffig{TwistRegion}.

More precisely, recall that a diagram of a knot is a 4-valent graph with over/under crossing information at each vertex. A twist region is a string of bigon\index{bigon} regions in the diagram graph, arranged end-to-end at their vertices, which is maximal in the sense that there are no additional bigon regions meeting the vertices on either end. A single crossing adjacent to no bigons is also a twist region. We will further restrict so that all twist regions are alternating, meaning crossings alternate over and under while following a strand of the twist region. If not, the second Reidemeister move\index{Reidemeister move} applied to the diagram removes two crossings from the twist region.
\end{definition}

The condition that twist regions be maximal ensures that there is only one way to put together exactly two twist regions in a diagram. 

\begin{definition}\label{Def:TwistKnot}
  The \emph{twist knot} $J(2,n)$\index{twist knot} is the knot with a diagram consisting of exactly two twist regions, one of which contains two crossings. The other twist region contains $n\in\ZZ$ crossings. The direction of crossing depends on the sign of $n$.
\end{definition}

\begin{figure}[h]
  \includegraphics{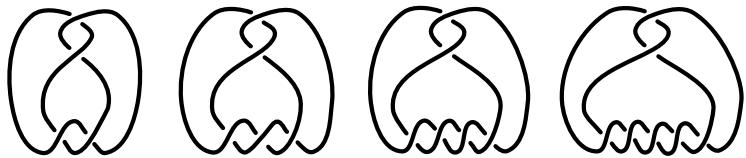}
  \caption{Twist knots $J(2,2)$ (the figure-8 knot), $J(2,3)$ (the $5_2$ knot), $J(2,4)$ (the $6_1$ or Stevedore knot), and $J(2,5)$}
  \label{Fig:TwistKnots}
\end{figure}

Twist knots $J(2,2)$, $J(2,3)$, $J(2,4)$, and $J(2,5)$ are shown in \reffig{TwistKnots}. 

\begin{figure}[h]
  \includegraphics{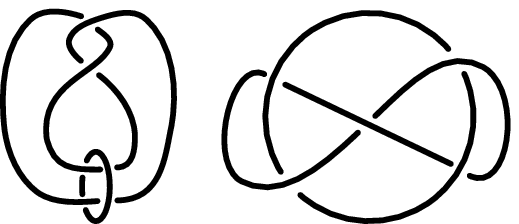}
  \caption{Two diagrams of the Whitehead link.\index{Whitehead link}}
  \label{Fig:Whitehead0}
\end{figure}

The family of twist knots $J(2,n)$ has very nice hyperbolic geometry, which we discuss in \refchap{TwistKnots}. In particular, as $n$ approaches infinity, we will see that the hyperbolic geometry of twist knot complements limits, in a precise sense, to the hyperbolic geometry of the Whitehead link complement; the Whitehead link is shown in \reffig{Whitehead0}.\index{Whitehead link}

More generally, any family of knots containing higher and higher numbers of crossings in a twist region will have complements converging to a link with a simple circle encircling that twist region. Knots with high numbers of crossings in twist regions are called \emph{highly twisted}.\index{highly twisted} Again these are discussed in \refchap{TwistKnots}. 

Given a diagram of a link, we can combine twist regions by performing a sequence of moves on the diagram called \emph{flypes}.

\begin{figure}[h]
  \includegraphics{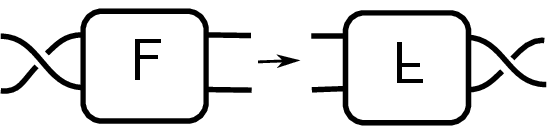}
  \caption{A flype.}
  \label{Fig:Flype0}
\end{figure}

\begin{definition}\label{Def:Flype}
  Let $\gamma$ be a simple closed curve meeting the diagram of $K$ transversely exactly four times away from crossings, with two intersections adjacent to a crossing on the outside of $\gamma$. A \emph{flype}\index{flype} is a move on the diagram that rotates the region inside $\gamma$ by $180^\circ$, moving the crossing adjacent to $\gamma$ to become a crossing adjacent to $\gamma$ but between the opposite two strands. 
  See \reffig{Flype0}.
\end{definition}

Now, suppose a simple closed curve $\gamma$ in the plane of projection meets a diagram transversely exactly four times away from crossings, and suppose also that the curve is adjacent to crossings on both sides. Then we can perform a flype to move one of the crossings to the opposite side of the curve, to form a bigon.\index{bigon} If the bigon is not alternating, remove both crossings, producing a diagram with fewer crossings. Otherwise, there are two cases. Either the curve $\gamma$ encloses only bigons on one side to begin with, and the flype produces a diagram that is unchanged, or the flype has moved a crossing out of one twist region, on one side of $\gamma$, into a distinct twist region on the other side of $\gamma$. Performing the same flype a finite number of times will move all crossings in the twist region on one side of $\gamma$ into the twist region on the other side, thus reducing the number of twist regions of the diagram. Thus by performing a finite number of flypes, we obtain a diagram with a minimal number of twist regions. Such a diagram is called \emph{twist-reduced}.\index{twist-reduced}

Every knot has a twist-reduced diagram with some number of twist regions. On the other hand, for a fixed positive integer $T$, there are only finitely many ways of combining twist regions to form a twist-reduced diagram with $T$ twist regions. The collection of twist-reduced diagrams with $T$ twist regions forms an infinite family of diagrams. Two highly twisted\index{highly twisted} diagrams with the same pattern of twist regions will have similar hyperbolic geometry, in ways that can be quantified. Thus rather than classifying knots by crossing number, from a geometric perspective it may make more sense to classify knots by number of twist regions in a twist-reduced diagram, or \emph{twist-number}.\index{twist-number} This brings us to another (broad and vaguely-worded) problem.

\subsection{The problem of enumerating knots by geometry}
Enumerating knots by twist region may make more geometric sense than enumerating by crossing number, because highly twisted\index{highly twisted} knots have diagrams that relate well to their geometry, in a sense that will be made precise in \refchap{TwistKnots}. Given any knot, is there always a diagram that encodes hyperbolic geometry?

\begin{figure}[h]
  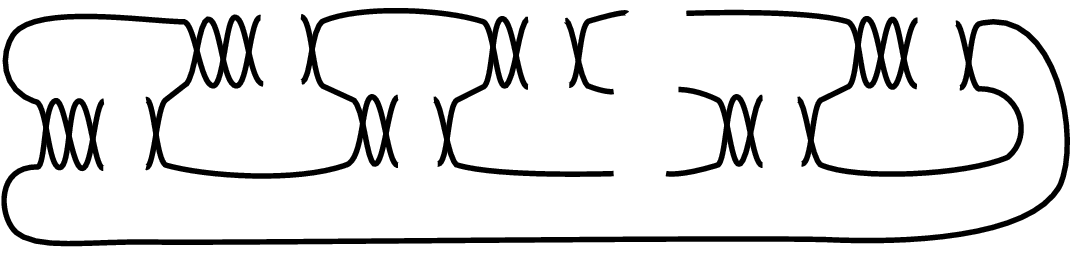
  \caption{A general form of a 2-bridge knot.\index{2-bridge knot or link}}
  \label{Fig:2BridgeDiagram0}
\end{figure}

Schubert considered a family of knots in 1956~\cite{Schubert}. He called the knots \emph{2-bridge knots}\index{2-bridge knot or link}. They can be described diagrammatically by taking four parallel strands, and twisting pairs of the strands into sequences of twist regions, then capping off either end with two ``bridges.'' A general form of such a diagram is shown in \reffig{2BridgeDiagram0}; see also \refchap{TwoBridge}. 

Although Schubert's work pre-dates the first work on the hyperbolic geometry of knots by nearly two decades, his 2-bridge knots\index{2-bridge knot or link} turn out to be very amenable to hyperbolic geometry techniques. We will see early on in this book that any knot exterior $S^3-N(K)$ can be decomposed into a collection of truncated tetrahedra.\index{truncated tetrahedron} Equivalently, $S^3-K$ is formed by gluing tetrahedra whose vertices have been removed. This is called an \emph{ideal triangulation}\index{ideal triangulation, definition} of the knot exterior, or sometimes simply a \emph{triangulation}.\index{triangulation, definition}

In the case of 2-bridge knots,\index{2-bridge knot or link} we will see that a triangulation of the knot complement can be read easily off the diagram. Not only that, we will see in \refchap{TwoBridge} that the edges and faces of the triangulation can be made totally geodesic under the hyperbolic metric, and the tetrahedra can be straightened simultaneously to be convex, with piecewise geodesic boundaries. Thus the combinatorics of the diagram of a 2-bridge knot gives a combinatorial method of describing the geometry of the 2-bridge knot.\index{2-bridge knot or link} This is very powerful.

It would be great to be able to extend these techniques to all knots, and some progress has been made with applications to other families, such as $n$-bridge knots for higher $n$. However, few families seem to be quite as nice as 2-bridge knots.\index{2-bridge knot or link}

There is still much ongoing work on triangulating knot exteriors and determining geometric properties of triangulations. We will discuss some of the techniques and applications in \refchap{AngleStruct}.

We have mentioned above that any knot exterior can be triangulated. In fact, any 3-manifold with torus boundary components can be decomposed into truncated tetrahedra. When the tetrahedra are convex hyperbolic tetrahedra, we say the triangulation is \emph{geometric}.\index{geometric triangulation} The software SnapPy has a census of orientable manifolds built up of at most nine geometric tetrahedra~\cite{SnapPy}. Some of these are knot complements.

This leads to a new way of classifying hyperbolic knots: by the number of geometric tetrahedra required to triangulate their exterior. This method of enumerating knots has been employed in~\cite{CallahanDeanWeeks, ckp, ckm}. 

\begin{figure}[h]
  \includegraphics{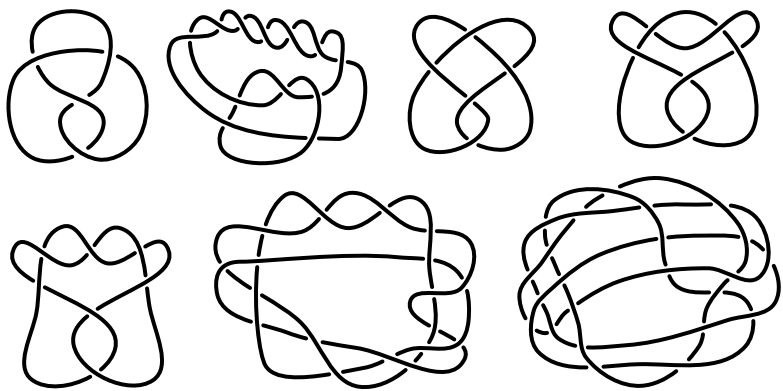}
  \caption{The seven simplest hyperbolic knots, built of at most four geometric tetrahedra.}
  \label{Fig:SimplestHyp}
\end{figure}

To date, 502 hyperbolic knots, built of at most eight geometric tetrahedra, have been classified. The diagrams of these knots often have large numbers of crossings. The knots built of at most four tetrahedra are shown in \reffig{SimplestHyp}.

Classifying knots by triangulations of their exteriors seems to be more difficult than classifying them by diagrams. This is because, given a triangulation of a 3-manifold with torus boundary, it is not obvious that the underlying space is a knot complement for a knot $K$ in $S^3$. We will discuss some techniques to detect whether such a manifold is a knot complement in \refchap{Essential}.

\subsection{The problem of finding geometric diagrams}
Twist knots and 2-bridge knots\index{2-bridge knot or link} have standard diagrams that encode a great deal of information about the geometry of the knot. Does every knot have such a diagram? (Probably not.) Does every knot have a diagram from which we may read some geometric information?

Alternating knots are another family of knots that seem to be amenable to hyperbolic geometric techniques.

\begin{definition}\label{Def:Alternating}
An \emph{alternating diagram}\index{alternating diagram!definition}\index{alternating diagram} is a diagram of a knot or link that has an orientation such that, when following the knot in the direction of the orientation, the crossings alternate between over and under. An \emph{alternating knot or link}\index{alternating knot or link}\index{alternating knot or link!definition} is a knot or link that has an alternating diagram.
\end{definition}

We will see in the exercises in \refchap{Fig8Decomp} that alternating knot complements decompose into pieces with the same combinatorics of the diagram. In chapters~\ref{Chap:Alternating}, \ref{Chap:Quasifuchsian}, and~\ref{Chap:Volume} we will use this decomposition to determine some geometric information on the knot complement.

How useful is this work broadly? All knots with at most seven crossings have alternating diagrams.\index{alternating diagram} Tait began his work~\cite{Tait:SciPapers} by assuming diagrams were alternating (although he did publish diagrams of eight- and ten-crossing non-alternating examples in 1877). 
However, the proportion of alternating knots in diagrams enumerated by crossing number rapidly drops to zero~\cite{SundbergThistlethwaite, Thistlethwaite:Tangles}. As for knots enumerated by geometric triangulations,\index{geometric triangulation} non-alternating examples seem to be even more common; a non-alternating example appears as the second knot on the list in \reffig{SimplestHyp}. Thus unfortunately, alternating knots\index{alternating knot or link} and links are not very common.

An open research question is, how many of the techniques presented in these chapters for determining geometry of alternating links generalize to other knots and links? There has been much work in recent years in extending this work to other families of knots, and some success. We are far from using such techniques to find hyperbolic geometry of all knots, though.

\subsection{The problem of determining geometric invariants}
One way of distinguishing knots is to compute invariants for each of them. If the invariants disagree, then the knots cannot be equivalent.

Several knot invariants arise classically, such as the crossing number that we encountered above. Many additional knot invariants arise through geometry. One aim of this book is to discuss such invariants, and give tools to calculate them. 

One of the most straightforward knot invariants that arises in geometry is the volume of a knot. We will show in \refchap{Margulis} that any knot complement that admits a hyperbolic structure has finite volume. Thus volumes of knots give knot invariants.

For those knots whose diagrams are particularly amenable to geometric techniques, such as twist knots, 2-bridge knots,\index{2-bridge knot or link} and alternating knots,\index{alternating knot or link} there are known methods to estimate volume using the combinatorics of the diagram. This is discussed along the way, but especially in \refchap{Volume}, where we bring to bear several tools in geometry to give two-sided bounds on volumes.

How powerful is volume as a knot invariant? It can be easy to calculate numerically, using the software SnapPy~\cite{SnapPy}, for example. Such computations can be rigorously verified to lie in a fixed error range using interval arithmetic, as in~\cite{HIKMOT}. Thus computing volume is a useful tool for distinguishing knots with distinct volume. However, there are many distinct knots that cannot be distinguished by volume; they have the same volume. We give some methods of constructing such knots and links in \refchap{Quasifuchsian}.

Then, is there a better geometric knot invariant than volume to distinguish knots? In \refchap{Canonical}, we describe the \emph{canonical decomposition}\index{canonical decomposition} of a hyperbolic knot complement. This is a decomposition consisting of convex polyhedra. We will show that when two knots have the same canonical decomposition, they must necessarily have homeomorphic complements, and thus by the Gordon--Luecke theorem, they must be equivalent (up to reflection). Thus the canonical decomposition is a complete invariant for hyperbolic knots. Unfortunately, it is not easy to compute in general, and provable forms of canonical decompositions are only known for a few infinite families of knots, including 2-bridge knots~\cite{Gueritaud:thesis}.\index{2-bridge knot or link} Canonical decompositions of alternating knots\index{alternating knot or link} are still unknown in general, for example.

Finally, we discuss very briefly one polynomial invariant. 
In most standard books on knot theory, there will be chapters on polynomial invariants, particularly the Alexander polynomial and the Jones polynomial. We will not treat such polynomials here; they arise from techniques that do not use hyperbolic geometry. There is one polynomial invariant of knots that depends heavily on hyperbolic geometry, however. This is the $A$-polynomial.\index{$A$-polynomial} We devote \refchap{Character} to a discussion of the $A$-polynomial, its definition and computation for a few examples. We will see that it relates to hyperbolic structures on a knot complement and the deformations of such structures.

\subsection{The problem of relating geometric invariants to other invariants}

What of the invariants that are being omitted from this book? We mentioned above Alexander and Jones polynomials. There are also more modern algebraic knot invariants, such as Khovanov homology and Floer homologies, and quantum invariants such as colored Jones polynomials. 

Many open problems in knot theory, driving much of the ongoing research in the field, concern relating invariants of knots arising from other fields of mathematics to hyperbolic geometry and hyperbolic knot invariants. We will not discuss in detail these open problems, because defining non-hyperbolic invariants will take us too far afield. However, one motivating factor for writing this book was to help mathematicians, particularly students, get up to speed with their hyperbolic geometry, in order to investigate the relations of geometry to other invariants in knot theory.

%%%%%%%%%%%%%%%%%%%%%%%%%%%%%%%%%%%%%%%%%%%%%%%%%%%%%%%%%%%%%%%%%
\section{Exercises}

\begin{exercise}
  Find a sequence of isotopy moves of the diagram of the Goeritz knot, \reffig{Goeritz}, that reduces it to the standard diagram of the unknot with no crossings.
\end{exercise}

\begin{exercise}
  Download and install the software SnapPy~\cite{SnapPy}. Use it to sketch diagrams of a few knots, and determine whether the knot is hyperbolic. Do this for at least one hyperbolic knot and at least one non-hyperbolic knot.
\end{exercise}

\begin{exercise}
  Convince yourself by drawing several examples that every 4-valent planar graph can be assigned over/under crossing information at each vertex to obtain an alternating knot.\index{alternating knot or link} Now try to prove this fact. (This may require some graph theory.)
\end{exercise}

\begin{exercise}
  Show that a connected sum\index{connected sum} of two knots is always a satellite knot.\index{satellite knot}
\end{exercise}

% Part I: Foundations of hyperbolic structures
%% Ch01_Fig8Decomp.tex
%% Chapter 1 of Purcell Hyperbolic Knot Theory

\part{Foundations of Hyperbolic Structures}\label{Part:Foundations}

\chapter{Decomposition of the Figure-8 Knot}\label{Chap:Fig8Decomp}
\blfootnote{Jessica S. Purcell, Hyperbolic Knot Theory}

In this chapter, we begin developing tools to work with knots and links and the 3-manifolds they define. We give a geometric method, explained carefully by example, to decompose a knot or link complement into simple pieces. The methods here are an introduction to topological techniques in 3-manifold geometry and topology, and an introduction to some of the tools used in the field. 

One goal of this chapter is to present a method that will allow us to pass from a knot or link \emph{diagram} to a description of the knot or link \emph{complement}. That is, we start with a 4-valent graph describing a knot or link $K$, and obtain a mathematically rigorous decomposition of the 3-manifold $S^3-K$ into simple 3-dimensional pieces, which pieces will be useful for applying tools from geometry and 3-manifold topology. 

\section{Polyhedra}

\begin{figure}[h]
\includegraphics{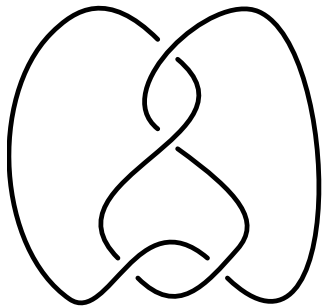}
\caption{A diagram of the figure-8 knot.}
\label{Fig:Fig8Diagram}
\end{figure}

Sometimes it is easier to study manifolds, including knot complements, if we split them into smaller, simpler pieces, for example 3-balls. We are going to decompose the figure-8 knot complement into two carefully marked 3-balls, namely ideal polyhedra. The diagram of the figure-8 knot that we use is shown in \reffig{Fig8Diagram}. The decomposition we describe appears in Thurston's notes~\cite{thurston}, and with a little more explanation in~\cite{thurston:book}. The procedure has been generalized to all link complements, for example in~\cite{menasco:links}. This work is essentially what we present below in the text and in exercises. 

\begin{definition}\label{Def:polyhedron}
A \emph{polyhedron}\index{polyhedron, definition} is a closed 3-ball whose boundary is labeled with a finite graph, containing a finite number vertices and edges, so that complementary regions, which are called \emph{faces},\index{face} are simply connected.

An \emph{ideal polyhedron}\index{ideal polyhedron, definition} is a polyhedron with all vertices removed. That is, to form an ideal polyhedron, start with a regular polyhedron and remove the points corresponding to vertices.
\end{definition}

We will cut $S^3 - K$ into two ideal polyhedra. We will then have a description of $S^3 - K$ as a gluing of two ideal polyhedra. That is, given a description of the polyhedra, and gluing information on the faces of the polyhedra, we may reconstruct the knot complement $S^3- K$. Although we use the example of the figure-8 knot, in the exercises, you will walk through the techniques below to determine decompositions of other knot complements into ideal polyhedra, and to generalize to all knots. 

\subsection{Overview}

\begin{figure}[h]
%% Accompanies image file 'F1-02-F8Face.eps' (pdf, eps, ps)
\begingroup%
  \makeatletter%
  \providecommand\color[2][]{%
    \errmessage{(Inkscape) Color is used for the text in Inkscape, but the package 'color.sty' is not loaded}%
    \renewcommand\color[2][]{}%
  }%
  \providecommand\transparent[1]{%
    \errmessage{(Inkscape) Transparency is used (non-zero) for the text in Inkscape, but the package 'transparent.sty' is not loaded}%
    \renewcommand\transparent[1]{}%
  }%
  \providecommand\rotatebox[2]{#2}%
  \newcommand*\fsize{\dimexpr\f@size pt\relax}%
  \newcommand*\lineheight[1]{\fontsize{\fsize}{#1\fsize}\selectfont}%
  \ifx\svgwidth\undefined%
    \setlength{\unitlength}{119.25789642bp}%
    \ifx\svgscale\undefined%
      \relax%
    \else%
      \setlength{\unitlength}{\unitlength * \real{\svgscale}}%
    \fi%
  \else%
    \setlength{\unitlength}{\svgwidth}%
  \fi%
  \global\let\svgwidth\undefined%
  \global\let\svgscale\undefined%
  \makeatother%
  \begin{picture}(1,0.94545795)%
    \lineheight{1}%
    \setlength\tabcolsep{0pt}%
    \put(0,0){\includegraphics[width=\unitlength]{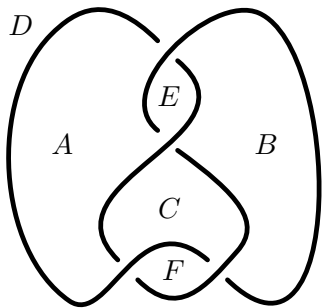}}%
    \put(0.49051237,0.08497803){\color[rgb]{0,0,0}\makebox(0,0)[lt]{\lineheight{1.25}\smash{\begin{tabular}[t]{l}$F$\end{tabular}}}}%
    \put(0.14165845,0.48550332){\color[rgb]{0,0,0}\makebox(0,0)[lt]{\lineheight{1.25}\smash{\begin{tabular}[t]{l}$A$\end{tabular}}}}%
    \put(0.47861463,0.2830131){\color[rgb]{0,0,0}\makebox(0,0)[lt]{\lineheight{1.25}\smash{\begin{tabular}[t]{l}$C$\end{tabular}}}}%
    \put(0.00635975,0.86214996){\color[rgb]{0,0,0}\makebox(0,0)[lt]{\lineheight{1.25}\smash{\begin{tabular}[t]{l}$D$\end{tabular}}}}%
    \put(0.47597072,0.63805362){\color[rgb]{0,0,0}\makebox(0,0)[lt]{\lineheight{1.25}\smash{\begin{tabular}[t]{l}$E$\end{tabular}}}}%
    \put(0.78312542,0.48550332){\color[rgb]{0,0,0}\makebox(0,0)[lt]{\lineheight{1.25}\smash{\begin{tabular}[t]{l}$B$\end{tabular}}}}%
  \end{picture}%
\endgroup%

\caption{Faces for the figure-8 knot complement.}
\label{Fig:Faces}
\end{figure}

Start with a diagram of the knot. There will be two polyhedra in our decomposition. These can be visualized as two balloons: One balloon expands above the diagram, and one balloon expands below the diagram. As the balloons continue expanding, they will bump into each other in the regions cut out by the graph of the diagram. Label these regions. In \reffig{Faces}, the regions are labeled $A$, $B$, $C$, $D$, $E$, and $F$.  These will correspond to faces of the polyhedra.

\begin{figure}[h]
%% Accompanies image file 'F1-03-CrFold.eps' (pdf, eps, ps)
\begingroup%
  \makeatletter%
  \providecommand\color[2][]{%
    \errmessage{(Inkscape) Color is used for the text in Inkscape, but the package 'color.sty' is not loaded}%
    \renewcommand\color[2][]{}%
  }%
  \providecommand\transparent[1]{%
    \errmessage{(Inkscape) Transparency is used (non-zero) for the text in Inkscape, but the package 'transparent.sty' is not loaded}%
    \renewcommand\transparent[1]{}%
  }%
  \providecommand\rotatebox[2]{#2}%
  \newcommand*\fsize{\dimexpr\f@size pt\relax}%
  \newcommand*\lineheight[1]{\fontsize{\fsize}{#1\fsize}\selectfont}%
  \ifx\svgwidth\undefined%
    \setlength{\unitlength}{183bp}%
    \ifx\svgscale\undefined%
      \relax%
    \else%
      \setlength{\unitlength}{\unitlength * \real{\svgscale}}%
    \fi%
  \else%
    \setlength{\unitlength}{\svgwidth}%
  \fi%
  \global\let\svgwidth\undefined%
  \global\let\svgscale\undefined%
  \makeatother%
  \begin{picture}(1,0.60655738)%
    \lineheight{1}%
    \setlength\tabcolsep{0pt}%
    \put(0,0){\includegraphics[width=\unitlength]{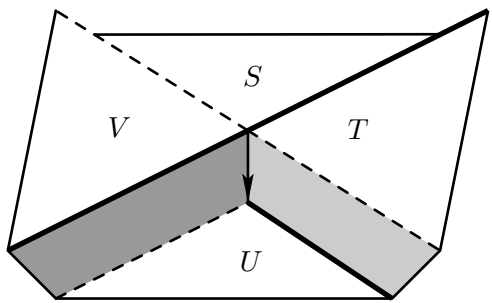}}%
    \put(0.48524595,0.06316217){\color[rgb]{0,0,0}\makebox(0,0)[lt]{\lineheight{1.25}\smash{\begin{tabular}[t]{l}$U$\end{tabular}}}}%
    \put(0.21354829,0.33987263){\color[rgb]{0,0,0}\makebox(0,0)[lt]{\lineheight{1.25}\smash{\begin{tabular}[t]{l}$V$\end{tabular}}}}%
    \put(0.70781735,0.33485977){\color[rgb]{0,0,0}\makebox(0,0)[lt]{\lineheight{1.25}\smash{\begin{tabular}[t]{l}$T$\end{tabular}}}}%
    \put(0.49226391,0.44614547){\color[rgb]{0,0,0}\makebox(0,0)[lt]{\lineheight{1.25}\smash{\begin{tabular}[t]{l}$S$\end{tabular}}}}%
  \end{picture}%
\endgroup%

  \caption{The knot runs along diagonals. Faces labeled $U$ and $T$ meet
    at the edge shown, marked by an arrow.}
\label{Fig:3dFolded}
\end{figure}

The faces meet up in edges. There is one edge for each crossing. It runs vertically from the knot at the top of the crossing to the knot at the bottom (or the other way around). The balloon expands until faces meet at edges. \Reffig{3dFolded} shows how the top balloon would expand at a crossing. The edge is drawn as an arrow from the top of the crossing to the bottom. Faces labeled $T$ and $U$ meet across the edge. Rotating the picture $180^\circ$ about the edge, we would see an identical picture with $S$ meeting $V$.

\begin{figure}
%% Accompanies image file 'F1-04-CModel.eps' (pdf, eps, ps)
\begingroup%
  \makeatletter%
  \providecommand\color[2][]{%
    \errmessage{(Inkscape) Color is used for the text in Inkscape, but the package 'color.sty' is not loaded}%
    \renewcommand\color[2][]{}%
  }%
  \providecommand\transparent[1]{%
    \errmessage{(Inkscape) Transparency is used (non-zero) for the text in Inkscape, but the package 'transparent.sty' is not loaded}%
    \renewcommand\transparent[1]{}%
  }%
  \providecommand\rotatebox[2]{#2}%
  \newcommand*\fsize{\dimexpr\f@size pt\relax}%
  \newcommand*\lineheight[1]{\fontsize{\fsize}{#1\fsize}\selectfont}%
  \ifx\svgwidth\undefined%
    \setlength{\unitlength}{147bp}%
    \ifx\svgscale\undefined%
      \relax%
    \else%
      \setlength{\unitlength}{\unitlength * \real{\svgscale}}%
    \fi%
  \else%
    \setlength{\unitlength}{\svgwidth}%
  \fi%
  \global\let\svgwidth\undefined%
  \global\let\svgscale\undefined%
  \makeatother%
  \begin{picture}(1,0.75510204)%
    \lineheight{1}%
    \setlength\tabcolsep{0pt}%
    \put(0,0){\includegraphics[width=\unitlength]{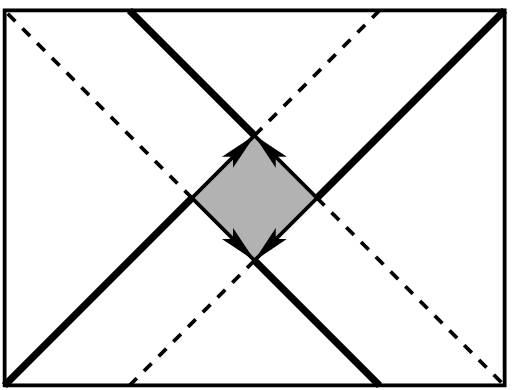}}%
    \put(0.45114145,0.05759257){\color[rgb]{0,0,0}\makebox(0,0)[lt]{\lineheight{1.25}\smash{\begin{tabular}[t]{l}$U$\end{tabular}}}}%
    \put(0.07439035,0.34395536){\color[rgb]{0,0,0}\makebox(0,0)[lt]{\lineheight{1.25}\smash{\begin{tabular}[t]{l}$V$\end{tabular}}}}%
    \put(0.79429697,0.37035192){\color[rgb]{0,0,0}\makebox(0,0)[lt]{\lineheight{1.25}\smash{\begin{tabular}[t]{l}$T$\end{tabular}}}}%
    \put(0.4735386,0.62711862){\color[rgb]{0,0,0}\makebox(0,0)[lt]{\lineheight{1.25}\smash{\begin{tabular}[t]{l}$S$\end{tabular}}}}%
  \end{picture}%
\endgroup%

  \caption{Cut out the shaded square. Start with a pair of parallel lines. Fold the thick part of the line in a direction opposite that	of the dashed part of the line. Fold parallel thick and dashed lines in opposite directions. Correct folding results in a model that looks like \reffig{3dFolded}.}
\label{Fig:3dEdge}
\end{figure}

It may be helpful to examine the meeting of faces at an edge by 3-dimensional model. Henry Segerman has come up with a paper model to illustrate the phenomenon of \reffig{3dFolded}. Start with a sheet of paper labeled as in \reffig{3dEdge}. Cut out the shaded square in the middle. Now fold the paper until it looks like that in \reffig{3dFolded}. By rotating the paper model,
we can see how faces meet up.

Stringing crossings such as this one together, we obtain the complete polyhedral decomposition of the knot. This is the geometric intuition behind the polyhedral expansion. We now explain a combinatorial method to describe the polyhedra.

\subsection{Step 1.}  Sketch faces and edges into the diagram.

Recall a diagram is a 4-valent graph lying on a plane, the plane of projection. The regions on the plane of projection that are cut out by the graph will be the faces, including the outermost unbounded region of the plane of projection. We start by labeling these, as in \reffig{Faces}.  

\begin{figure}
  \includegraphics{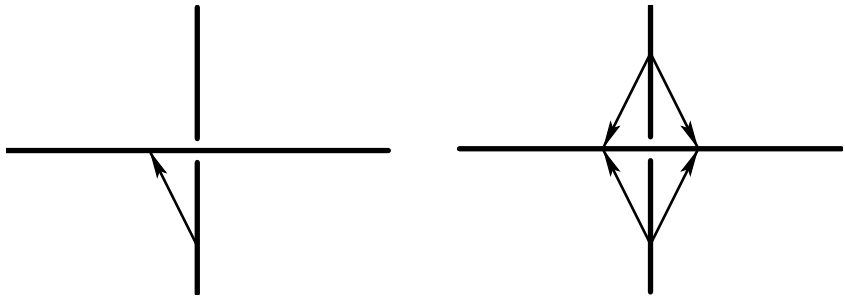}
  \caption{A single edge.}
  \label{Fig:Edges}
\end{figure}

Edges come from arcs that connect the two strands of the diagram at a crossing. These are called \emph{crossing arcs}.\index{crossing arc} For ease of explanation, we are going to draw each edge four times, as follows. Shown on the left of \reffig{Edges} is a single edge corresponding to a crossing arc. Note that the edge is ambient isotopic\index{ambient isotopic} in $S^3$ to the three additional edges shown on the right in \reffig{Edges}.

\begin{figure}
%% Accompanies image file 'F1-06-8Edge.eps' (pdf, eps, ps)
\begingroup%
  \makeatletter%
  \providecommand\color[2][]{%
    \errmessage{(Inkscape) Color is used for the text in Inkscape, but the package 'color.sty' is not loaded}%
    \renewcommand\color[2][]{}%
  }%
  \providecommand\transparent[1]{%
    \errmessage{(Inkscape) Transparency is used (non-zero) for the text in Inkscape, but the package 'transparent.sty' is not loaded}%
    \renewcommand\transparent[1]{}%
  }%
  \providecommand\rotatebox[2]{#2}%
  \newcommand*\fsize{\dimexpr\f@size pt\relax}%
  \newcommand*\lineheight[1]{\fontsize{\fsize}{#1\fsize}\selectfont}%
  \ifx\svgwidth\undefined%
    \setlength{\unitlength}{119.25789642bp}%
    \ifx\svgscale\undefined%
      \relax%
    \else%
      \setlength{\unitlength}{\unitlength * \real{\svgscale}}%
    \fi%
  \else%
    \setlength{\unitlength}{\svgwidth}%
  \fi%
  \global\let\svgwidth\undefined%
  \global\let\svgscale\undefined%
  \makeatother%
  \begin{picture}(1,0.94545795)%
    \lineheight{1}%
    \setlength\tabcolsep{0pt}%
    \put(0,0){\includegraphics[width=\unitlength]{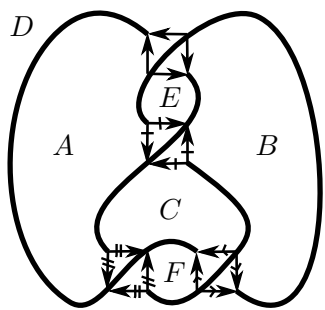}}%
    \put(0.49051237,0.08497803){\color[rgb]{0,0,0}\makebox(0,0)[lt]{\lineheight{1.25}\smash{\begin{tabular}[t]{l}$F$\end{tabular}}}}%
    \put(0.14165845,0.48550332){\color[rgb]{0,0,0}\makebox(0,0)[lt]{\lineheight{1.25}\smash{\begin{tabular}[t]{l}$A$\end{tabular}}}}%
    \put(0.47861463,0.2830131){\color[rgb]{0,0,0}\makebox(0,0)[lt]{\lineheight{1.25}\smash{\begin{tabular}[t]{l}$C$\end{tabular}}}}%
    \put(0.00635975,0.86214996){\color[rgb]{0,0,0}\makebox(0,0)[lt]{\lineheight{1.25}\smash{\begin{tabular}[t]{l}$D$\end{tabular}}}}%
    \put(0.47597072,0.63805362){\color[rgb]{0,0,0}\makebox(0,0)[lt]{\lineheight{1.25}\smash{\begin{tabular}[t]{l}$E$\end{tabular}}}}%
    \put(0.78312542,0.48550332){\color[rgb]{0,0,0}\makebox(0,0)[lt]{\lineheight{1.25}\smash{\begin{tabular}[t]{l}$B$\end{tabular}}}}%
  \end{picture}%
\endgroup%

\caption{Edges of the figure-8 knot complement.}
\label{Fig:Fig8Edges}
\end{figure}

The reason for sketching each edge four times is that it allows us to visualize easily which edges bound the faces we have already labeled. In \reffig{Fig8Edges}, we have drawn four copies of each of the four edges we get from crossing arcs of the diagram of the figure-8 knot.  Note that the face labeled $A$, for example, will be bordered by three edges, one with two tick marks, one with a single tick mark, and one with no tick marks.

\begin{remark}
Orientations on the edges can be chosen to run in either direction; that is, arrows on the edges can run from overcrossing to undercrossing or vice versa, as long as we are consistent with orientations corresponding to the same edge. We have chosen the orientations in \reffig{Fig8Edges} to simplify a later step, and to match a figure in \refchap{GluingCompleteness}. The opposite choice for any edge is also fine.
\end{remark}

\subsection{Step 2} Shrink the knot to ideal vertices on the top polyhedron. 

Now we come to the reason for using \emph{ideal} polyhedra, rather than regular polyhedra. Notice that the edges stretch from a part of the knot to a part of the knot. However, the manifold we are trying to model is the knot complement, $S^3 - K$. Therefore, the knot $K$ does not exist in the manifold. An edge with its two vertices on $K$ must necessarily be an ideal edge; that is, its vertices are not contained in the manifold $S^3- K$.

Since the knot is \emph{not} part of the manifold, we will shrink strands of the knot to ideal vertices. That is, retract each knot strand to a single point.
This may cause some confusion at first, because the strand of the knot is not homeomorphic to a single point. However, we are considering the \emph{complement} of the strand. The complement of the strand on the boundary of the ball is homeomorphic to the complement of a single point on the boundary of the ball, so we replace strands by ideal vertices (single removed points). 

Focus first on the polyhedron on top. Each component of the knot we ``see'' from inside the top polyhedron will be shrunk to a single ideal vertex. These visible knot components correspond to sequences of overcrossings of the diagram. Compare to \reffig{3dFolded} --- note that at an undercrossing, the component of the knot ends in an edge, but at an overcrossing the knot continues on. Moreover, note that at an overcrossing, the knot passes the same edge twice, once on each side.

\begin{figure}[h]
%% Accompanies image file 'F1-07-8Top.eps' (pdf, eps, ps)
\begingroup%
  \makeatletter%
  \providecommand\color[2][]{%
    \errmessage{(Inkscape) Color is used for the text in Inkscape, but the package 'color.sty' is not loaded}%
    \renewcommand\color[2][]{}%
  }%
  \providecommand\transparent[1]{%
    \errmessage{(Inkscape) Transparency is used (non-zero) for the text in Inkscape, but the package 'transparent.sty' is not loaded}%
    \renewcommand\transparent[1]{}%
  }%
  \providecommand\rotatebox[2]{#2}%
  \newcommand*\fsize{\dimexpr\f@size pt\relax}%
  \newcommand*\lineheight[1]{\fontsize{\fsize}{#1\fsize}\selectfont}%
  \ifx\svgwidth\undefined%
    \setlength{\unitlength}{119.25789642bp}%
    \ifx\svgscale\undefined%
      \relax%
    \else%
      \setlength{\unitlength}{\unitlength * \real{\svgscale}}%
    \fi%
  \else%
    \setlength{\unitlength}{\svgwidth}%
  \fi%
  \global\let\svgwidth\undefined%
  \global\let\svgscale\undefined%
  \makeatother%
  \begin{picture}(1,0.94545795)%
    \lineheight{1}%
    \setlength\tabcolsep{0pt}%
    \put(0,0){\includegraphics[width=\unitlength]{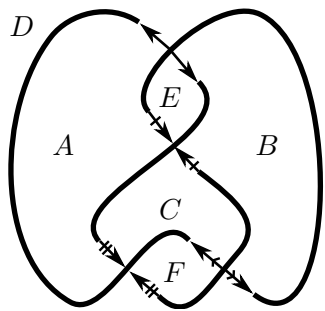}}%
    \put(0.49051237,0.08497803){\color[rgb]{0,0,0}\makebox(0,0)[lt]{\lineheight{1.25}\smash{\begin{tabular}[t]{l}$F$\end{tabular}}}}%
    \put(0.14165845,0.48550332){\color[rgb]{0,0,0}\makebox(0,0)[lt]{\lineheight{1.25}\smash{\begin{tabular}[t]{l}$A$\end{tabular}}}}%
    \put(0.47861463,0.2830131){\color[rgb]{0,0,0}\makebox(0,0)[lt]{\lineheight{1.25}\smash{\begin{tabular}[t]{l}$C$\end{tabular}}}}%
    \put(0.00635975,0.86214996){\color[rgb]{0,0,0}\makebox(0,0)[lt]{\lineheight{1.25}\smash{\begin{tabular}[t]{l}$D$\end{tabular}}}}%
    \put(0.47597072,0.63805362){\color[rgb]{0,0,0}\makebox(0,0)[lt]{\lineheight{1.25}\smash{\begin{tabular}[t]{l}$E$\end{tabular}}}}%
    \put(0.78312542,0.48550332){\color[rgb]{0,0,0}\makebox(0,0)[lt]{\lineheight{1.25}\smash{\begin{tabular}[t]{l}$B$\end{tabular}}}}%
  \end{picture}%
\endgroup%

  \caption{Isotopic edges in top polyhedron identified.}
\label{Fig:8Top}
\end{figure}

In terms of the four copies of the edge in \reffig{Edges}, when we consider the polyhedron on top, we may identify the two edges which are isotopic along an overstrand, but not those isotopic along understrands. See \reffig{8Top}.

\begin{figure}[h]
%% Accompanies image file 'F1-08-8TopP.eps' (pdf, eps, ps)
\begingroup%
  \makeatletter%
  \providecommand\color[2][]{%
    \errmessage{(Inkscape) Color is used for the text in Inkscape, but the package 'color.sty' is not loaded}%
    \renewcommand\color[2][]{}%
  }%
  \providecommand\transparent[1]{%
    \errmessage{(Inkscape) Transparency is used (non-zero) for the text in Inkscape, but the package 'transparent.sty' is not loaded}%
    \renewcommand\transparent[1]{}%
  }%
  \providecommand\rotatebox[2]{#2}%
  \newcommand*\fsize{\dimexpr\f@size pt\relax}%
  \newcommand*\lineheight[1]{\fontsize{\fsize}{#1\fsize}\selectfont}%
  \ifx\svgwidth\undefined%
    \setlength{\unitlength}{135.19239807bp}%
    \ifx\svgscale\undefined%
      \relax%
    \else%
      \setlength{\unitlength}{\unitlength * \real{\svgscale}}%
    \fi%
  \else%
    \setlength{\unitlength}{\svgwidth}%
  \fi%
  \global\let\svgwidth\undefined%
  \global\let\svgscale\undefined%
  \makeatother%
  \begin{picture}(1,0.85106607)%
    \lineheight{1}%
    \setlength\tabcolsep{0pt}%
    \put(0,0){\includegraphics[width=\unitlength]{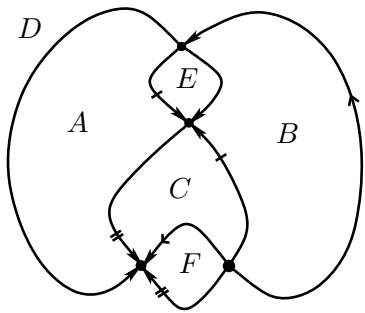}}%
    \put(0.16789609,0.50357645){\color[rgb]{0,0,0}\makebox(0,0)[lt]{\lineheight{1.25}\smash{\begin{tabular}[t]{l}$A$\end{tabular}}}}%
    \put(0.7505247,0.47267256){\color[rgb]{0,0,0}\makebox(0,0)[lt]{\lineheight{1.25}\smash{\begin{tabular}[t]{l}$B$\end{tabular}}}}%
    \put(0.45239336,0.31906224){\color[rgb]{0,0,0}\makebox(0,0)[lt]{\lineheight{1.25}\smash{\begin{tabular}[t]{l}$C$\end{tabular}}}}%
    \put(0.03064662,0.76807703){\color[rgb]{0,0,0}\makebox(0,0)[lt]{\lineheight{1.25}\smash{\begin{tabular}[t]{l}$D$\end{tabular}}}}%
    \put(0.47057214,0.62810076){\color[rgb]{0,0,0}\makebox(0,0)[lt]{\lineheight{1.25}\smash{\begin{tabular}[t]{l}$E$\end{tabular}}}}%
    \put(0.47966151,0.11000675){\color[rgb]{0,0,0}\makebox(0,0)[lt]{\lineheight{1.25}\smash{\begin{tabular}[t]{l}$F$\end{tabular}}}}%
  \end{picture}%
\endgroup%

  \caption{Top polyhedron, viewed from the inside.}
  \label{Fig:8TopPoly}
\end{figure}

Shrink each overstrand to a single ideal vertex. The result is pattern of faces, edges, and ideal vertices for the top polyhedron, shown in \reffig{8TopPoly}. Notice that the face $D$ is a disk, containing the point at infinity.

\subsection{Step 3} Shrink the knot to ideal vertices for the bottom polyhedron.

\vspace{.1in}

Notice that underneath the knot, the picture of faces, edges, and vertices will be slightly different. In particular, when finding the top polyhedron, we collapsed overstrands to a single ideal vertex. When you put your head underneath the knot, what appear as overstrands from below will appear as understrands on the usual knot diagram.

\begin{figure}[h]
%% Accompanies image file 'F1-09-CrFldB.eps' (pdf, eps, ps)
\begingroup%
  \makeatletter%
  \providecommand\color[2][]{%
    \errmessage{(Inkscape) Color is used for the text in Inkscape, but the package 'color.sty' is not loaded}%
    \renewcommand\color[2][]{}%
  }%
  \providecommand\transparent[1]{%
    \errmessage{(Inkscape) Transparency is used (non-zero) for the text in Inkscape, but the package 'transparent.sty' is not loaded}%
    \renewcommand\transparent[1]{}%
  }%
  \providecommand\rotatebox[2]{#2}%
  \newcommand*\fsize{\dimexpr\f@size pt\relax}%
  \newcommand*\lineheight[1]{\fontsize{\fsize}{#1\fsize}\selectfont}%
  \ifx\svgwidth\undefined%
    \setlength{\unitlength}{183bp}%
    \ifx\svgscale\undefined%
      \relax%
    \else%
      \setlength{\unitlength}{\unitlength * \real{\svgscale}}%
    \fi%
  \else%
    \setlength{\unitlength}{\svgwidth}%
  \fi%
  \global\let\svgwidth\undefined%
  \global\let\svgscale\undefined%
  \makeatother%
  \begin{picture}(1,0.60655738)%
    \lineheight{1}%
    \setlength\tabcolsep{0pt}%
    \put(0,0){\includegraphics[width=\unitlength]{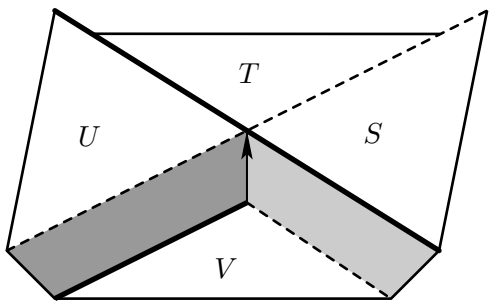}}%
    \put(0.74171388,0.32736424){\color[rgb]{0,0,0}\makebox(0,0)[lt]{\lineheight{1.25}\smash{\begin{tabular}[t]{l}$S$\end{tabular}}}}%
    \put(0.15432886,0.32268761){\color[rgb]{0,0,0}\makebox(0,0)[lt]{\lineheight{1.25}\smash{\begin{tabular}[t]{l}$U$\end{tabular}}}}%
    \put(0.48543443,0.45176266){\color[rgb]{0,0,0}\makebox(0,0)[lt]{\lineheight{1.25}\smash{\begin{tabular}[t]{l}$T$\end{tabular}}}}%
    \put(0.4349268,0.04957229){\color[rgb]{0,0,0}\makebox(0,0)[lt]{\lineheight{1.25}\smash{\begin{tabular}[t]{l}$V$\end{tabular}}}}%
  \end{picture}%
\endgroup%

  \caption{3-dimensional model, opposite side as in \reffig{3dFolded}. Now faces $V$ and $U$ meet along an edge.}
  \label{Fig:3dBottom}
\end{figure}

One way to see this difference is to take the 3-dimensional
model constructed in \reffig{3dEdge}. \Reffig{3dFolded} shows the view of the faces meeting at an edge from the top. If you turn the model over to the opposite side, you will see how the faces meet underneath. \Reffig{3dBottom}
illustrates this. Note $U$ now meets $V$, and $S$ meets $T$.

In terms of the combinatorics, edges of \reffig{Edges} that are isotopic by sliding an endpoint along an understrand are identified to each other on the bottom polyhedron, but edges only isotopic by sliding an endpoint along an overstrand are not identified.

\begin{figure}[h]
%% Accompanies image file 'F1-10-8BotP.eps' (pdf, eps, ps)
\begingroup%
  \makeatletter%
  \providecommand\color[2][]{%
    \errmessage{(Inkscape) Color is used for the text in Inkscape, but the package 'color.sty' is not loaded}%
    \renewcommand\color[2][]{}%
  }%
  \providecommand\transparent[1]{%
    \errmessage{(Inkscape) Transparency is used (non-zero) for the text in Inkscape, but the package 'transparent.sty' is not loaded}%
    \renewcommand\transparent[1]{}%
  }%
  \providecommand\rotatebox[2]{#2}%
  \newcommand*\fsize{\dimexpr\f@size pt\relax}%
  \newcommand*\lineheight[1]{\fontsize{\fsize}{#1\fsize}\selectfont}%
  \ifx\svgwidth\undefined%
    \setlength{\unitlength}{134.73947525bp}%
    \ifx\svgscale\undefined%
      \relax%
    \else%
      \setlength{\unitlength}{\unitlength * \real{\svgscale}}%
    \fi%
  \else%
    \setlength{\unitlength}{\svgwidth}%
  \fi%
  \global\let\svgwidth\undefined%
  \global\let\svgscale\undefined%
  \makeatother%
  \begin{picture}(1,0.82720976)%
    \lineheight{1}%
    \setlength\tabcolsep{0pt}%
    \put(0,0){\includegraphics[width=\unitlength]{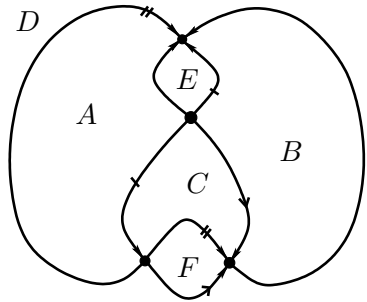}}%
    \put(0.18829424,0.49468102){\color[rgb]{0,0,0}\makebox(0,0)[lt]{\lineheight{1.25}\smash{\begin{tabular}[t]{l}$A$\end{tabular}}}}%
    \put(0.0201103,0.75532211){\color[rgb]{0,0,0}\makebox(0,0)[lt]{\lineheight{1.25}\smash{\begin{tabular}[t]{l}$D$\end{tabular}}}}%
    \put(0.46918785,0.59330194){\color[rgb]{0,0,0}\makebox(0,0)[lt]{\lineheight{1.25}\smash{\begin{tabular}[t]{l}$E$\end{tabular}}}}%
    \put(0.75800628,0.39253793){\color[rgb]{0,0,0}\makebox(0,0)[lt]{\lineheight{1.25}\smash{\begin{tabular}[t]{l}$B$\end{tabular}}}}%
    \put(0.49824579,0.30272238){\color[rgb]{0,0,0}\makebox(0,0)[lt]{\lineheight{1.25}\smash{\begin{tabular}[t]{l}$C$\end{tabular}}}}%
    \put(0.47182949,0.0676171){\color[rgb]{0,0,0}\makebox(0,0)[lt]{\lineheight{1.25}\smash{\begin{tabular}[t]{l}$F$\end{tabular}}}}%
  \end{picture}%
\endgroup%

  \caption{Bottom polyhedron, from the outside.}
  \label{Fig:8BotPoly}
\end{figure}

As above, collapse each knot strand corresponding to an understrand to a single ideal vertex. The result is \reffig{8BotPoly}.

\vspace{.1in}

One thing to notice: we sketched the top polyhedron with our heads inside the ball on top, looking out. If we move the face $D$ away from the point at infinity, then it wraps \emph{above} the other faces shown in \reffig{8TopPoly}.

On the other hand, we sketched the bottom polyhedron with our heads outside the ball on the bottom. If we move the face $D$ away from the point at infinity, it wraps \emph{below} the other faces shown in \reffig{8BotPoly}.

\subsection{Rebuilding the knot complement from the polyhedra} Figures~\ref{Fig:8TopPoly} and~\ref{Fig:8BotPoly} show two ideal polyhedra that we obtained by studying the figure-8 knot complement. We claim that they glue to give the figure-8 knot complement. That is, attach face $A$ on the bottom polyhedron to the face labeled $A$ on the top polyhedron, ensuring that the edges bordering face $A$ match up. Similarly for the other faces.

This process of gluing faces and edges gives exactly the complement of the knot. By construction, faces glue to give the faces illustrated in \reffig{Fig8Edges}, and edges glue to give the edges there, except when we have finished, all four edges in an isotopy shown in that figure have been glued together.

\section{Generalizing: Exercises}

This polyhedral decomposition works for any knot or link diagram, to give a polyhedral decomposition of its complement.

\begin{figure}[h]
  \begin{tabular}{ccccc}
    \includegraphics{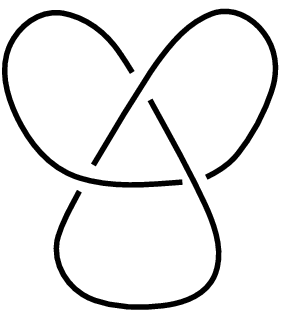} & \hspace{.2in} & 
    \includegraphics{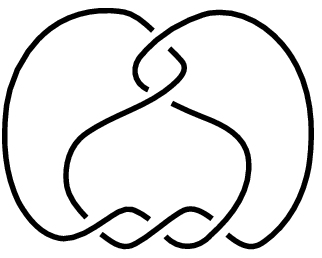} & \hspace{.2in} &
    \includegraphics{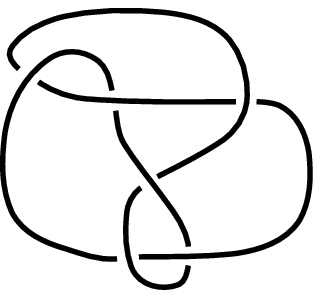} \\
    (a) Trefoil. & \hspace{.2in} & (b)  The $5_2$ knot.  &
    \hspace{.2in} & (c) The $6_3$ knot. \\
  \end{tabular}
  \caption{Three examples of knots.}
  \label{Fig:SimpleKnots}
\end{figure}

\begin{exercise}\label{Ex:WarmUp}
  As a warm-up exercise, determine the polyhedral decomposition for one (or more) of the knots shown in \reffig{SimpleKnots}. Sketch both top and bottom polyhedra.

  Your solution should consist of two ideal polyhedra, i.e.\ marked graphs on the surface of a ball, with faces and edges marked according to the gluing pattern. For example, the complete diagrams in Figures~\ref{Fig:8TopPoly} and~\ref{Fig:8BotPoly} form the solution for the figure-8 knot. 
\end{exercise}

 \begin{exercise}
The examples of knots we have encountered so far are all alternating,\index{alternating knot or link} as in \refdef{Alternating}. The diagram of the knot $8_{19}$ in \reffig{8-19} is not alternating. In fact, the knot $8_{19}$ has no alternating diagram.\index{alternating diagram!example with no alternating diagram}

\begin{figure}[h]
  \includegraphics{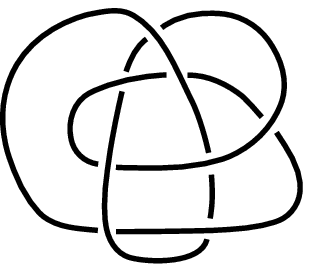}
  \caption{The knot $8_{19}$, which has no alternating diagram.\index{alternating diagram}\index{alternating diagram!example with no alternating diagram}}
  \label{Fig:8-19}
\end{figure}

Determine the polyhedral decomposition for the given diagram of the knot $8_{19}$. Note: as above, many ideal vertices are obtained by shrinking overstrands to a point. However, you will have to use, for example, \reffig{3dFolded} to determine what happens between two understrands.
\end{exercise}

\begin{exercise}
Recall that the \emph{valence}\index{valence} of a vertex in a graph is the number of edges that meet that vertex. The valence of an ideal vertex is defined similarly. 

  \begin{enumerate}
  \item[(a)] If a knot diagram is alternating,\index{alternating knot or link}\index{alternating knot or link!polyhedral decomposition} we obtain a very special ideal polyhedron. In particular, all ideal vertices will have the same valence. What is it?  Show that the ideal vertices for an alternating knot all have this valence.

  \item[(b)] What are the possible valences of ideal vertices in general, i.e.\ for non-alternating knots?  For which $n\geq 0 \in \ZZ$ is there a knot diagram whose polyhedral decomposition yields an ideal vertex of valence $n$?  Explain your answer, with (portions of) knot diagrams.
  \end{enumerate}
\end{exercise}

\begin{exercise}
  In the polyhedral decomposition for alternating knots,\index{alternating knot or link!polyhedral decomposition} the polyhedra are given by simply labeling each ball with the projection graph of the knot and declaring each vertex to be ideal.
  \begin{enumerate}
  \item Prove this statement for any alternating knot. That is, prove that the decomposition gives polyhedra whose edges match the projection graph of the diagram. 
  \item Show that for non-alternating knots, this is false. That is, the decomposition does not give polyhedra whose edges match the projection graph of the diagram. 
  \end{enumerate}
\end{exercise}

\begin{figure}[h]
\includegraphics{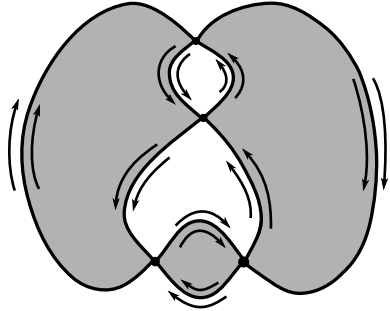}
\caption{Checkerboard coloring and ``gear rotation'' for the figure-8 knot.}
\label{Fig:Gears}
\end{figure}

\begin{exercise}
  A graph admits a \emph{checkerboard coloring}\index{checkerboard coloring} if all the complementary regions can be colored either white or shaded, with white faces meeting shaded faces across the edges. Any 4-valent graph can be checkerboard colored, particularly projection graphs of knot diagrams.

  In the case of an alternating knot,\index{alternating knot or link!polyhedral decomposition} faces are identified from the top polyhedron to the identical face on the bottom polyhedron, and the identification is by a \emph{gear rotation}:\index{gear rotation} white faces on the top are rotated once counter-clockwise and then glued to the corresponding face on the bottom; shaded faces on the top are rotated once clockwise and then glued. This is shown for the figure-8 knot in \reffig{Gears}. Prove that for the decomposition of any alternating knot, faces are identified by a gear rotation.
\end{exercise}

\begin{exercise}
The diagrams we have encountered so far are all reduced, as in \refdef{Nugatory}, but we can follow the above procedure for non-reduced diagrams. For example, we can obtain a polyhedral decomposition for diagrams which contain a nugatory crossing\index{nugatory crossing}.

Show that the polyhedral decomposition of a knot diagram will contain a monogon, i.e.\ a face whose boundary is a single edge and a single vertex, if and only if the diagram has a simple nugatory crossing.
\end{exercise}

\begin{exercise}\label{Ex:CollapseBigons}
Recall that a bigon\index{bigon} is a region of a graph bounded by exactly two edges and exactly two vertices. Note that when a bigon appears in our polyhedral decomposition, the two edges of the bigon must be isotopic to each other. Hence, we sometimes will remove bigon faces from the polyhedral decomposition, identifying their two edges.

\begin{quote}
  Let bigons\index{bigon} be bygone. --- William Menasco
\end{quote}

For the figure-8 knot, sketch the two polyhedra we get when bigon faces are removed. How many edges are there in this new, bigon-free decomposition? The resulting polyhedra are well known solids in this case. What are they?

For each of the polyhedra obtained in \refex{WarmUp}, sketch the resulting polyhedra with bigons removed.
\end{exercise}

\begin{exercise}
Suppose we start with an alternating knot\index{alternating knot or link} diagram with at least two crossings, and do the polyhedral decomposition above, collapsing bigons at the last step. What are possible valences of vertices? Sketch the diagram of a single alternating knot that has all possible valences of ideal vertices in its polyhedral decomposition.

What valences of vertices can you get if you don't require the diagram to be alternating but collapse bigons? Can you find 1-valent vertices? For any $n>4 \in \ZZ$, can you find $n$-valent vertices?  
\end{exercise}

%% Ch02_IntroHyp.tex
%% 

\chapter{Calculating in Hyperbolic Space}\label{Chap:IntroHyp}
\blfootnote{Jessica S. Purcell, Hyperbolic Knot Theory}

We will need to manipulate objects in 2 and 3-dimensional hyperbolic space. This chapter provides a very brief introduction to the tools that will be needed in the future, the objects that will be studied (lines, triangles, tetrahedra, metric properties), and examples of calculations that will appear.

We will use terminology and calculations from standard elementary Riemannian geometry. The reader who is not as comfortable with Riemannian geometry might find it helpful to follow along in the first few chapters of an introductory Riemannian geometry text, such as do~Carmo \cite[Chapter~1]{docarmo}. We will not provide all the details to all the statements given. The idea is that we want to begin calculating on knot complements and other 3-manifolds immediately, without getting lost early in details. Thus our aim is to provide just enough information here to start calculating in future chapters. Many more details and results can be found in other books, including full books on hyperbolic geometry. Anderson gives a very nice introduction to 2-dimensional hyperbolic geometry~\cite{anderson}. More details in all dimensions appear in Ratcliffe~\cite{ratcliffe}. The book by Marden includes more on groups of isometries of hyperbolic space, including results on infinite volume hyperbolic 3-manifolds~\cite{marden}. 
An introduction to hyperbolic geometry that includes a discussion of its visualization is also given by Thurston~\cite{thurston:book}. 

\section{Hyperbolic geometry in dimension two}

We start with hyperbolic 2-space, $\HH^2$.

There are several models of hyperbolic space. Here, we will work with the upper half plane model. In this model, hyperbolic 2-space $\HH^2$ is defined to be the set of points in the upper half plane:
\[
\HH^2 = \{ x+ i\,y \in \CC \mid y>0 \},
\]
equipped with the metric whose first fundamental form is given by
\[
ds^2 = \frac{dx^2 + dy^2}{y^2}.
\]
That is, start with the usual Euclidean metric on $\RR^2$, whose first fundamental form is $dx^2+dy^2$. To obtain the metric on the hyperbolic plane, rescale the usual Euclidean metric by $1/y$, where $y$ is height in the plane.

Note that a point in $\HH^2$ can either be thought of as a complex number $x+i\,y \in \CC$ or as a point $(x,y)\in\RR^2$. Both perspectives are useful: $\RR^2$ leads more easily to coordinates and calculations, and $\CC$ works seamlessly with our definition of isometries below. Changing perspectives does not affect our results, so we will regularly switch between the two without comment. 

Our first task is to explore the meaning of the hyperbolic metric, and how it affects measurements.

\subsection{Hyperbolic 2-space and Riemannian geometry}

In this subsection, we briefly review how the metric and the space $\HH^2$ described above fit into a more general picture of Riemannian geometry. We also describe tools from Riemannian geometry we will use to do calculations. If you are not yet familiar with Riemannian geometry, feel free to skim this section, noting equations \eqref{Eqn:ArcLength}, \eqref{Eqn:VolForm}, and \eqref{Eqn:Area}. This section was primarily written for a student who has seen some Riemannian geometry, but may have difficulty applying abstract concepts of that field to the specific metric of hyperbolic geometry. In the author's experience, a few key equations will be enough to get started. 

In Riemannian geometry, a \emph{Riemannian metric}\index{Riemannian metric} on a manifold $M$ is defined to be a correspondence associating to each point $p\in M$ an inner product $\langle\cdot, \cdot\rangle_p$ on the tangent space $T_pM$. This inner product gives us a way of measuring the lengths of vectors tangent to $M$ at $p$, as well as computing areas, angles between curves, etc. The \emph{first fundamental form}\index{first fundamental form} is defined by $\langle v, v\rangle_p$ for $v\in T_pM$. 

In our case, the Riemannian manifold we consider is $\HH^2$, and we have natural local coordinates on the manifold given by $x+i\,y \in \CC$, or $(x,y)\in\RR^2$, for $y>0$. We may use these coordinates to describe the Riemannian metric. In particular, at the point $(x,y)\in\HH^2$, a tangent vector $v\in T_{(x,y)}\HH^2$ can also be described by coordinates $v = v_x \frac{\partial}{\partial x} + v_y\frac{\partial}{\partial y}$, and we write it as a vector
\[ v = \left( \begin{array}{c} v_x \\ v_y \end{array} \right). \]
Then the metric on $\HH^2$ is given by a matrix
\[ \langle v, w \rangle = (v_x, v_y) \mat{1/y^2 & 0 \\ 0 & 1/y^2} \left(\begin{array}{c} w_x \\ w_y \end{array} \right). \]

One of the simplest geometric measurements we can compute using the definition of the metric is the arc length of a curve. If $\gamma(t)$ is a (differentiable) curve in $\HH^2$, for $t \in [a,b]$, then we obtain a tangent vector $\gamma'(t)$ at each point of $\gamma(t)$ in $\HH^2$, called the velocity vector. The \emph{arc length}\index{arc length} of $\gamma$ for $t\in[a,b]$ is defined to be
\[ |\gamma| = \int_a^b \sqrt{\langle \gamma'(s), \gamma'(s) \rangle}\, ds. \]
When considering $\HH^2$, we will have coordinates $\gamma(t)=(\gamma_x(t), \gamma_y(t))$, and $\gamma'(t) = (\gamma'_x(t), \gamma'_y(t))^T$. Thus the arc length will be
\begin{equation}\label{Eqn:ArcLength}
|\gamma| = \int_a^b \sqrt{ (\gamma'_x(s))^2 + (\gamma_y'(s))^2 }\, \frac{1}{\gamma_y(s)}\,ds.
\end{equation}
We will use this formula to compute examples in the next subsection.

Another piece of geometric information we can compute with a metric is the volume of a region, which we typically call ``area'' in two dimensions. In the most general setting, if $R\subset M$ is contained in a coordinate neighborhood of the Riemannian manifold $M$, with coordinates $(x_1, \dots, x_n)$ and metric given by the matrix $g_{ij}$ in these coordinates, then we can compute the volume of $R$ to be
\begin{equation}\label{Eqn:VolForm}
 \vol(R) = \int_R d\vol = \int_R \sqrt{\det(g_{ij})}\,dx_1 \dots dx_n.
\end{equation}
The form $d\vol$ is the \emph{volume form}\index{volume form}. 
Thus in our setting, with $M=\HH^2$ and metric as above,
\begin{equation}\label{Eqn:Area}
  \area(R) = \int_R \,\frac{1}{y^2}\, dx\,dy.
\end{equation}

\subsection{Computing arc lengths and areas}

Now we will use the formulas obtained above to do calculations, in order to better understand the hyperbolic space $\HH^2$. 

\begin{example}\label{Example:HorizontalLine}
Fix a height $h>0$, and consider first a horizontal line segment between points $(0,h) = i\,h$ and $(1,h) = 1+i\,h$ in $\HH^2$. We may parameterize the line segment by $\gamma(t) = (t, h)$, for $t\in[0,1]$. Using \refeqn{ArcLength}, we find the arc length\index{arc length} of $\gamma$ is $|\gamma| = 1/h$. Note that because $h$ is fixed, the arc length of $\gamma$ is just its usual Euclidean length rescaled by $1/h$. Thus when $h=1$, the length of $\gamma$ is $1$. When $h$ becomes large, the arc length becomes very small. In other words, points with the same height become very close together as their heights increase. On the other hand, as $h$ approaches $0$, the length of $\gamma$ approaches infinity. In fact, points near the real line $\RR=\{(x,0)\in\RR^2\}$ can be very far apart.
\end{example}

\begin{example}\label{Example:VerticalLine}
Consider now a vertical line between points $(x,a)$ and $(x,b)$, for $x,a,b$ fixed in $\RR$, $0<a<b$. Such a line can be parameterized by $\zeta(t) = (x,t)$ for $t\in [a,b]$. So $\zeta'(t) = (0,1)$. Thus its arc length\index{arc length} is given by
\[
|\zeta| = \int_a^b \sqrt{0+1}\,\frac{1}{s}\, ds = \log\left(\frac{b}{a}\right).
\]
If we set $b=1$ and let $a$ approach $0$, note that the arc length\index{arc length} of $\zeta$ gets arbitrarily large, approaching infinity. Similarly setting $a=1$ and letting $b$ approach infinity gives arbitrarily long lengths.
\end{example}

The real line $\RR = \{(x,0)\in \RR^2\}$ along with the point at infinity $\infty$ play an important role in the geometry of $\HH^2$, although these points are not contained in $\HH^2$.

\begin{definition}\label{Def:BdryatInfty}
We call $\RR\cup\{\infty\}$ the \emph{boundary at infinity}\index{boundary at infinity} for $\HH^2$. Note it is homeomorphic to a circle $S^1$, and hence is sometimes called the \emph{circle at infinity}\index{circle at infinity}. It is denoted by $S^1_\infty$, $\bdy\HH^2$, and sometimes $\bdy_\infty\HH^2$.
\end{definition}

Areas behave quite differently in hyperbolic space than in Euclidean space.

\begin{example}\label{Example:AreaHorocycle}
In this example, we will compute the area of the region $R$ of $\HH^2$ that is the intersection of the half-plane lying to the left of the line $x=1$, the half-plane to the right of the line $x=0$, and the plane lying above $y=1$. The region $R$ is shown in \reffig{example-region}. 

\begin{figure}
  \includegraphics{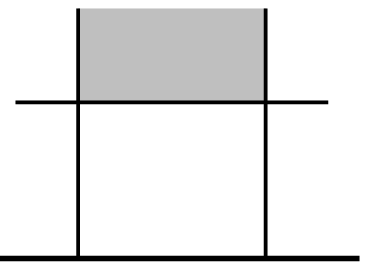}
  \caption{The region of \refexamp{AreaHorocycle}}
  \label{Fig:example-region}
\end{figure}

Using \refeqn{Area}, we see that the area of the region is given by
\begin{align*}
  \area(R) & = \int_R \frac{1}{y^2}\,dx\,dy \\
  & = \int_0^1 \int_1^\infty \frac{1}{y^2}\, dy\, dx \\
  & = \int_0^1 1\, dx = 1
\end{align*}

This example shows that regions with infinite Euclidean area may have finite hyperbolic area. 
\end{example}

%%%%%%%%%%%%%%%%%%%%%%%%%%%%%%%%%%%%%%%%%%%%%%%%%%%%%%%%%%%%%%%%%
\subsection{Geodesics and isometries}

Recall that a \emph{geodesic}\index{geodesic} between points $p$ and $q$ is a length minimizing curve between those points. An \emph{infinite geodesic}\index{geodesic!infinite} is a curve $\gamma$ from $\RR$ to a Riemannian manifold such that for any $s< t\in \RR$, the curve $\gamma([s,t])$ minimizes the distance between $\gamma(s)$ and $\gamma(t)$. 

\begin{theorem}\label{Thm:GeodesicsH2}
The infinite geodesics in $\HH^2$ consist of vertical straight lines and semi-circles with center on the real line.\qed
\end{theorem}

\begin{figure}
%% Accompanies image file 'F2-02-HypGeo.eps' (pdf, eps, ps)
\begingroup%
  \makeatletter%
  \providecommand\color[2][]{%
    \errmessage{(Inkscape) Color is used for the text in Inkscape, but the package 'color.sty' is not loaded}%
    \renewcommand\color[2][]{}%
  }%
  \providecommand\transparent[1]{%
    \errmessage{(Inkscape) Transparency is used (non-zero) for the text in Inkscape, but the package 'transparent.sty' is not loaded}%
    \renewcommand\transparent[1]{}%
  }%
  \providecommand\rotatebox[2]{#2}%
  \newcommand*\fsize{\dimexpr\f@size pt\relax}%
  \newcommand*\lineheight[1]{\fontsize{\fsize}{#1\fsize}\selectfont}%
  \ifx\svgwidth\undefined%
    \setlength{\unitlength}{182.00000954bp}%
    \ifx\svgscale\undefined%
      \relax%
    \else%
      \setlength{\unitlength}{\unitlength * \real{\svgscale}}%
    \fi%
  \else%
    \setlength{\unitlength}{\svgwidth}%
  \fi%
  \global\let\svgwidth\undefined%
  \global\let\svgscale\undefined%
  \makeatother%
  \begin{picture}(1,0.40659339)%
    \lineheight{1}%
    \setlength\tabcolsep{0pt}%
    \put(0,0){\includegraphics[width=\unitlength]{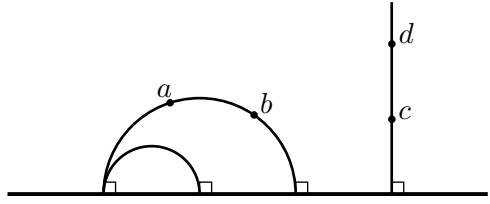}}%
    \put(0.31314728,0.20483647){\color[rgb]{0,0,0}\makebox(0,0)[lt]{\lineheight{1.25}\smash{\begin{tabular}[t]{l}$a$\end{tabular}}}}%
    \put(0.52696291,0.17284829){\color[rgb]{0,0,0}\makebox(0,0)[lt]{\lineheight{1.25}\smash{\begin{tabular}[t]{l}$b$\end{tabular}}}}%
    \put(0.81373397,0.15994079){\color[rgb]{0,0,0}\makebox(0,0)[lt]{\lineheight{1.25}\smash{\begin{tabular}[t]{l}$c$\end{tabular}}}}%
    \put(0.81429511,0.31707562){\color[rgb]{0,0,0}\makebox(0,0)[lt]{\lineheight{1.25}\smash{\begin{tabular}[t]{l}$d$\end{tabular}}}}%
  \end{picture}%
\endgroup%

  \caption{Some geodesics and points in $\HH^2$.}
  \label{Fig:HypGeodesics}
\end{figure}

Note these are exactly the circles and lines in the upper half plane that meet $S^1_\infty$ at right angles. See \reffig{HypGeodesics}. Observe that between any two points in the upper half plane, there is a unique vertical line or semi-circle between them. Thus a geodesic between points $p$ and $q$ in $\HH^2$ is a segment of a semi-circle or a vertical straight line. An infinite geodesic can also be viewed as the unique semi-circle or vertical straight line between two points on the boundary at infinity\index{boundary at infinity} of $\HH^2$. We will typically drop the word ``infinite'' to describe geodesics between points on the boundary at infinity. Thus we use the same word ``geodesic'' to describe both infinite or bounded arcs, depending on context. 

The proof of \refthm{GeodesicsH2} is left as an exercise in Riemannian geometry. The simplest way to prove the theorem uses coordinates and a bit more Riemannian geometry than we have reviewed so far. The interested reader can work through the details. The fact that these are the geodesics of $\HH^2$ is all we will need going forward.

An \emph{isometry}\index{isometry, definition} between Riemannian manifolds $M$ and $N$ is a diffeomorphism $f\from M\to N$ such that
\[
\langle v, w \rangle_p = \langle df_p(v), df_p(w) \rangle_{f(p)} \quad \mbox{ for all } p\in M, v, w \in T_pM.
\]

Isometries preserve lengths, angles, and other geometric information. We are most interested in orientation preserving isometries from hyperbolic space to itself, i.e.\ diffeomorphisms $\phi\from\HH^2\to\HH^2$ that preserve the metric and orientation on $\HH^2$. All such isometries form a group acting on $\HH^2$. We will assume the following theorem; see also \refex{Reflection}. 

\begin{theorem}\label{Thm:IsomH2}
The full group of isometries of $\HH^2$ is generated by reflections in geodesics in $\HH^2$.

The group of orientation preserving isometries of $\HH^2$ is the group of \emph{linear fractional transformations}\index{linear fractional transformation}
\[ z\mapsto \frac{az + b}{cz+d}, \]
where $a,b,c,d\in\RR$, and $ad-bc > 0$.\qed
\end{theorem}

By taking the quotient of $a$, $b$, $c$, and $d$ by $\sqrt{ad-bc}$, the linear fractional transformation is equivalent to an element of $\PSL(2,\RR)$, the group of projective 2 by 2 matrices with real coefficients and determinant $1$. That is, we may view $A\in\PSL(2,\RR)$ as given by a matrix
\[ A = \pm\mat{a&b\\c&d}, \]
where $a$, $b$, $c$, $d \in \RR$ and $ad-bc=1$. The sign in front reflects the fact that it is \emph{projective}; it is well-defined only up to multiplication by $\pm \Id$. On the other hand, $A$ acts on $\HH^2$ via
\[ Az = \frac{az + b}{cz + d}. \]
Note the action is unaffected when we multiply $a$, $b$, $c$, and $d$ by the same real constant, thus it is necessary to take projective matrices. 

Linear fractional transformations take circles and lines to circles
and lines, so they map geodesics to geodesics. For more information on these transformations, see for example \cite[pp~76--89]{Ahlfors}.

The following lemma is very useful.

\begin{lemma}\label{Lem:3points}
Given any three distinct points $z_1$, $z_2$, and $z_3$ in $\bdy\HH^2$, there exists an orientation preserving isometry of $\HH^2$ taking $z_3$ to $\infty$, and taking $\{z_1,z_2\}$ to $\{0,1\}$.  It follows that there exists an isometry of $\HH^2$ taking any three distinct points on $\bdy\HH^2$ to any other three distinct points, with appropriate orientation.
\end{lemma}

\begin{proof}
This is a standard fact of linear fractional transformations. We need to take some care to preserve orientation. If necessary, switch $z_1$ and $z_2$ so that the sequence $z_1,z_2,z_3$ runs in counterclockwise order around $\bdy\HH^2 \cong S^1$.

If none of $z_1$, $z_2$, and $z_3$ are infinity, then a linear fractional transformation sending $z_1$ to $1$, $z_2$ to $0$, and $z_3$ to $\infty$ is given by
\[ z\mapsto \frac{z-z_2}{z-z_3}\frac{z_1-z_3}{z_1-z_2}.\]
Note that the determinant of this transformation is 
\[(z_1-z_3)(z_1-z_2)(z_2-z_3).\]
Because the sequence $z_1, z_2, z_3$ is in counterclockwise order, this is positive. Thus it gives the desired orientation preserving isometry.

If $z_1=\infty$, $z_2=\infty$, or $z_3=\infty$, then the isometry is given by
\[
z\mapsto \frac{z-z_2}{z-z_3}, \quad z\mapsto\frac{z_1-z_3}{z-z_3}, \quad z\mapsto\frac{z-z_2}{z_1-z_2}
\]
respectively. One can check that again, because we ensured the sequence $z_1,z_2,z_3$ is in counterclockwise order, the determinant of each transformation is positive. 
\end{proof}

Many metric calculations in $\HH^2$ can be simplified greatly by applying an appropriate isometry, including the use of \reflem{3points}. For example, the following lemma is easily proved using an isometry. 

\begin{lemma}\label{Lem:IntGeodesics}
Two distinct geodesics $\ell_1$ and $\ell_2$ in $\HH^2$ either
\begin{enumerate}
\item intersect in a single point in the interior of $\HH^2$,
\item intersect in a single point on the boundary $\bdy \HH^2$, or
\item are completely disjoint in $\HH^2 \cup \bdy \HH^2$.
\end{enumerate}
In the third case, there is a unique geodesic $\ell_3$ that is perpendicular to both $\ell_1$ and $\ell_2$.
\end{lemma}

\begin{proof}
We may apply an isometry $g$ of $\HH^2$, taking endpoints of $\ell_1$ to $0$ and $\infty$, and taking one of the endpoints of $\ell_2$ to $1$. The image of the second endpoint of $\ell_2$ under $g$ is then some point $w$ in $\bdy\HH^2 = \RR\cup\{\infty\}$. Note that $g(\ell_1)$ is the vertical line from $0$ to $\infty$ in $\HH^2$. The point $w$ determines the image of $g(\ell_2)$. 

If $w=0$ or if $w=\infty$, then we are in the second case, and $g(\ell_2)$ is a semi-circle with endpoints $0$ and $1$, or a vertical line from $1$ to $\infty$.

If $w \in \RR$ is less than zero, then we are in the first case. The two endpoints of $g(\ell_2)$ are separated by the line $g(\ell_1)$, so the geodesics must meet.

Finally, if $w\in\RR$ is greater than zero, then we are in the third case, and the geodesics are disjoint. One way to see that there is a unique geodesic perpendicular to both is to apply another isometry $h$, taking $\sqrt{w}$ to $0$ and $-\sqrt{w}$ to $\infty$. That is, let $h\from \HH^2\to \HH^2$ be given by
\[ h(z) = \frac{z-\sqrt{w}}{z+\sqrt{w}}. \]
Note that $h(0) = -1$, $h(\infty) = 1$, so $h(g(\ell_1))$ is the geodesic that is a semi-circle with endpoints at $-1$ and $1$. Also,
\[ h(1) = \frac{1-\sqrt{w}}{1+\sqrt{w}}\quad \mbox{and} \quad h(w)=\frac{w-\sqrt{w}}{w+\sqrt{w}}=-\frac{1-\sqrt{w}}{1+\sqrt{w}}. \]
So $h(g(\ell_2))$ is the geodesic that is a semi-circle with endpoints $h(1)$ and $-h(1)$. Thus images of both geodesics are semi-circles with center at $0$. The geodesic from $0$ to $\infty$ is therefore perpendicular to both, and it is the unique such geodesic. Set $\ell_3$ to be the image of the line from $0$ to $\infty$ under the composition $g^{-1}\circ h^{-1}$. 
\end{proof}

In the previous proof, knowing which isometry $h$ to apply in the last step required a calculation. However, once that isometry was applied, the existence and uniqueness of the geodesic $\ell_3$ was clear.

Computing lengths of geodesics is also simplified by applying isometries. 

\begin{example}\label{Example:length} Length computation.

  Suppose you wish to compute the length of a segment, or the distance between two points in $\HH^2$.  One strategy for doing this is to apply an isometry taking the two points to a simpler picture.  For example, in \reffig{HypGeodesics}, we may find an isometry taking the geodesic containing points $a$ and $b$ to the vertical geodesic from $0$ to $\infty$. Then under this isometry, the points $a$ and $b$ map to points of the form $(0,t_1)$ and $(0, t_2)$. 

In \refexamp{VerticalLine}, we already computed the length of the vertical segment between $(0, t_1)$ and $(0, t_2)$; its length is $\log(t_1/t_2)$ (assuming here that $t_2<t_1$, otherwise take the negative of the log). This gives the distance between $a$ and $b$.
\end{example}

\subsection{Triangles and horocycles}

\begin{definition}\label{Def:ideal-triangle}
An \emph{ideal triangle}\index{ideal triangle} in $\HH^2$ is a triangle with three geodesic edges, with all three vertices on $\bdy\HH^2$.
\end{definition}

There is an isometry of $\HH^2$ taking any ideal triangle\index{ideal triangle} to the ideal triangle with vertices $0$, $1$, and $\infty$, by \reflem{3points}.  Hence all ideal triangles in $\HH^2$ are isometric. In fact, we will see that they all have finite area. Thus all ideal triangles have the same area!

\begin{figure}
  \includegraphics{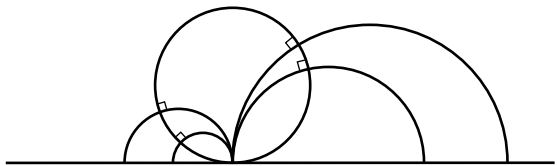}
\caption{A horocycle}
\label{Fig:horocycle}
\end{figure}

\begin{definition}\label{Def:horocycle}
A \emph{horocycle}\index{horocycle} centered at an ideal point $p \in \bdy\HH^2$ is defined as a curve perpendicular to all geodesics through $p$. When $p$ is a point on $\RR \subset \bdy\HH^2 = \RR \cup \{\infty\}$, a horocycle is a Euclidean circle tangent to $p$, as in \reffig{horocycle}.  When $p$ is the point $\infty$, a horocycle at $p$ is a line parallel to $\RR$.  That is, in this case the horocycle consists of points of the form $\{ (x,y) \mid y=c \}$ where $c>0$ is constant.
\end{definition}

\begin{definition}\label{Def:horoball}
  A \emph{horoball}\index{horoball} is the region of $\HH^2$ interior to a horocycle.
\end{definition}

Note a horoball will either be a Euclidean disk tangent to $\RR\subset\bdy\HH^2$ or a region consisting of points of the form $\{(x,y) \mid y>c \}$.

In \refexamp{AreaHorocycle}, we computed the area of a portion of a horoball, and we observed it was finite. Using this, we can show that the area of an ideal triangle is finite.

\begin{lemma}\label{Lem:AreaTriangleFinite}
  The area of an ideal triangle is finite.\index{ideal triangle}
\end{lemma}

\begin{proof}
Given any ideal triangle in $\HH^2$, we may apply an isometry taking its vertices to $0$, $1$, and $\infty$. Let $T$ denote this ideal triangle. Consider the intersection of $T$ with the horoball about infinity of height $1$. This is the region $R$ of \refexamp{AreaHorocycle}.

Note that the isometries
\[ z \mapsto \frac{z-1}{z} \quad \mbox{and} \quad z\mapsto\frac{-1}{z-1} \]
take the horoball about infinity to horoballs of Euclidean diameter $1$ centered at $1$ and at $0$, respectively, and take $T$ to $T$. Thus the intersections of $T$ with these horoballs also have areas $1$.

Finally, note that the complement of these horoballs in $T$ is a closed and bounded region $B$, lying below the line $y=1$, above the horocycle of Euclidean diameter $1$ centered at $0$, and above the horocycle of Euclidean diameter $1$ centered at $1$. The region $B$ lies in the rectangle $[0,1]\times[\half,1]$. It follows that the area of $B$ is at most the area of the rectangle, which is finite.

Thus the area of $T$ is $3$ plus the area of $B$, which is finite. 
\end{proof}

From the lemma, we see that the area of an ideal triangle is larger than~$3$. 
The exercises lead you through a calculation showing that the area of an ideal triangle is in fact $\pi$.

\section{Hyperbolic geometry in dimension three}

Hyperbolic 3-space is defined as follows:
\[
\HH^3 = \{ (x+iy, t) \in \CC \times \RR \mid t>0\}, \]
under the metric with first fundamental form
\begin{equation}\label{Eqn:3dHypMetric}
ds^2 = \frac{dx^2 + dy^2 + dt^2}{t^2}.
\end{equation}

We have the following theorems, which we will assume. Their proofs can be found in texts on hyperbolic geometry. 

\begin{theorem}\label{Thm:GeodesicsH3}
The geodesics in $\HH^3$ consist of vertical lines and semicircles orthogonal to the boundary $\bdy\HH^3 = \CC \cup \{\infty\}$.\index{boundary at infinity} Totally geodesic planes are vertical planes and hemispheres centered on $\CC$.\qed
\end{theorem}

\begin{theorem}\label{Thm:IsomH3}
The full group of isometries of $\HH^3$ is generated by reflections in geodesic planes.

The group of orientation preserving isometries of $\HH^3$ is $\PSL(2,\CC)$.  Its action on the boundary\index{boundary at infinity} $\bdy \HH^3 = \CC \cup \{\infty\}$ is the usual action of $\PSL(2,\CC)$ on $\CC \cup \{\infty\}$, via M\"obius transformation.\index{M\"obius transformation}\qed
\end{theorem}

An element $A\in\PSL(2,\CC)$ can be represented by a matrix, up to multiplication by $\pm \Id$. \Refthm{IsomH3} states that if
\[ A=\pm\mat{a&b\\c&d} \in \PSL(2,\CC),\]
then the action of $A$ on $\bdy\HH^3$ is given by
\[A(z) = \frac{az+b}{cz+d}, \mbox{ for } z\in\bdy\HH^3.\]

The action of an element of $\PSL(2,\CC)$ extends to the interior of hyperbolic 3-space, and there is a unique way to extend. Marden works through it carefully in \cite[Chapter~1]{marden}. However, we will not need the formula, and it is complicated, so we omit it here. 

\begin{theorem}\label{Thm:ClassifyPSL(2,C)}
Apart from the identity, any element of $\PSL(2,\CC)$ is exactly one of the following:
\begin{enumerate}
\item \emph{elliptic},\index{elliptic} which has two fixed points on $\bdy\HH^3$ and rotates about the geodesic axis between them in $\HH^3$, fixing the axis pointwise,
\item \emph{parabolic},\index{parabolic} which has a single fixed point on $\bdy\HH^3$, 
\item \emph{loxodromic},\index{loxodromic} which has two fixed points on $\bdy\HH^3$, and dilates and rotates about the axis between them.
\end{enumerate}
\end{theorem}

For example, the element $\mat{\exp(i\theta)& 0 \\ 0& \exp(-i\theta)}\in \PSL(2,\CC)$ is elliptic:\index{elliptic} it fixes the points $0$ and $\infty$, and the axis between them, and rotates about that axis by angle $2\theta$.

The element $\mat{1&1\\0&1} \in \PSL(2,\CC)$ is parabolic.\index{parabolic} It fixes the point $\infty$ only. Its action on $\bdy\HH^3$ is $z\mapsto z+1$, which extends to Euclidean translation by~$1$ in the interior of hyperbolic 3-space. 

Finally, the element $\mat{\rho & 0\\0& \rho^{-1}}$ is loxodromic\index{loxodromic} whenever $\rho$ is a complex number with $|\rho|>1$. It fixes the points $0$ and $\infty$ in $\bdy \HH^3$, but translates along the axis between them, and rotates and translates points in the interior of $\HH^3$ that do not lie on the axis.

In fact, after conjugating by an appropriate element of $\PSL(2,\CC)$, any element of $\PSL(2,\CC)$ actually becomes one of these three examples. This is stated as \reflem{MoreClassifyPSL} in \refchap{Margulis}. As a warm up for that theorem and \refthm{ClassifyPSL(2,C)}, \refex{IsomH2} works through a similar classification for isometries of $\HH^2$.

\begin{definition}\label{Def:ideal-tetr}
An \emph{ideal tetrahedron}\index{ideal tetrahedron, definition} is a tetrahedron in $\HH^3$ with all four vertices on $\bdy\HH^3$, and with geodesic edges and faces. 
\end{definition}

\begin{figure}
%% Accompanies image file 'F2-04-IdealT.eps' (pdf, eps, ps)
\begingroup%
  \makeatletter%
  \providecommand\color[2][]{%
    \errmessage{(Inkscape) Color is used for the text in Inkscape, but the package 'color.sty' is not loaded}%
    \renewcommand\color[2][]{}%
  }%
  \providecommand\transparent[1]{%
    \errmessage{(Inkscape) Transparency is used (non-zero) for the text in Inkscape, but the package 'transparent.sty' is not loaded}%
    \renewcommand\transparent[1]{}%
  }%
  \providecommand\rotatebox[2]{#2}%
  \newcommand*\fsize{\dimexpr\f@size pt\relax}%
  \newcommand*\lineheight[1]{\fontsize{\fsize}{#1\fsize}\selectfont}%
  \ifx\svgwidth\undefined%
    \setlength{\unitlength}{147.99999619bp}%
    \ifx\svgscale\undefined%
      \relax%
    \else%
      \setlength{\unitlength}{\unitlength * \real{\svgscale}}%
    \fi%
  \else%
    \setlength{\unitlength}{\svgwidth}%
  \fi%
  \global\let\svgwidth\undefined%
  \global\let\svgscale\undefined%
  \makeatother%
  \begin{picture}(1,0.75000002)%
    \lineheight{1}%
    \setlength\tabcolsep{0pt}%
    \put(0,0){\includegraphics[width=\unitlength]{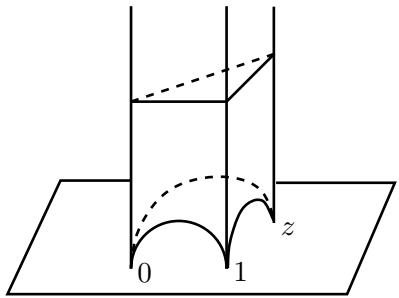}}%
    \put(0.33576053,0.03559865){\color[rgb]{0,0,0}\makebox(0,0)[lt]{\lineheight{1.25}\smash{\begin{tabular}[t]{l}$0$\end{tabular}}}}%
    \put(0.5800971,0.03721681){\color[rgb]{0,0,0}\makebox(0,0)[lt]{\lineheight{1.25}\smash{\begin{tabular}[t]{l}$1$\end{tabular}}}}%
    \put(0.70145623,0.1569579){\color[rgb]{0,0,0}\makebox(0,0)[lt]{\lineheight{1.25}\smash{\begin{tabular}[t]{l}$z$\end{tabular}}}}%
  \end{picture}%
\endgroup%

\caption{Ideal tetrahedron}
\label{Fig:ideal-tet}
\end{figure}

Since there exists a M\"obius transformation\index{M\"obius transformation} taking any three points to $0$, $1$, and $\infty$ in $\CC \cup \{\infty\}$, we may assume our tetrahedron has vertices at $0$, $1$, and $\infty$, and at some point $z \in \CC \setminus \{0,1\}$. So any ideal tetrahedron is parameterized by $z$. See \reffig{ideal-tet}.

The value of $z$ tells us about the geometry of the ideal tetrahedron. For example, the argument of $z$ is the dihedral angle between the vertical planes through $0,1,\infty$ and through $0, z,\infty$. 

The modulus of $z$ also has geometric meaning. Consider the hyperbolic geodesic through $z \in \CC$ that meets the vertical line from $0$ to $\infty$ in a right angle at a point $p_1$. Consider also the geodesic through $1 \in \CC$ that meets the vertical line from $0$ to $\infty$ at a right angle at point $p_2$. The hyperbolic distance between $p_1$ and $p_2$ is exactly $|\log|z||$
(\refex{CrossRatio}).  Hence
\[
\log z = (\mbox{signed dist between altitudes}) + i(\mbox{dihedral angle}).
\] 

\begin{definition}\label{Def:horosphere}
A \emph{horosphere}\index{horosphere} about $\infty$ in $\bdy\HH^3$ is a plane parallel to $\CC$, consisting of points $\{(x+iy, c) \in \CC \times \RR\}$ where $c>0$ is constant. Note for any $c>0$, this plane is perpendicular to all geodesics through $\infty$. When we apply an isometry that takes $\infty$ to some $p \in \CC$, note a horosphere is taken to a Euclidean sphere tangent to $p$.  By definition, this is a horosphere about $p$. A \emph{horoball}\index{horoball} is the region interior to a horosphere.
\end{definition}

\begin{figure}
\includegraphics{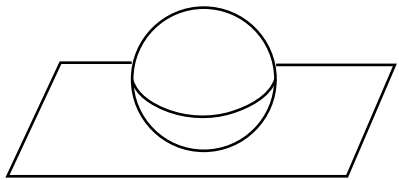}
\caption{Horosphere}
\label{Fig:horosphere}
\end{figure}

The metric on $\HH^3$ induces a metric on a horosphere. For a horosphere $\{x+iy,c)\in\CC\times\RR\}$ about $\infty$, the metric is just the Euclidean metric, rescaled by $1/c$. We may apply an isometry to any horosphere, taking it to one about $\infty$. Thus the induced metric on any horosphere will always be Euclidean. Hence when we intersect horospheres about $0$, $1$, $\infty$ and $z$ with an ideal tetrahedron through those points, we obtain four Euclidean triangles. These four triangles are similar (\refex{TetLabels1}).

\section{Exercises}

\begin{exercise}[Requires geometry]
  Prove \refthm{GeodesicsH2}, that is, show that vertical lines and semi-circles are geodesics, without using isometries of $\HH^2$. One way to solve this problem is to use Riemannian geometry, such as calculations in coordinates on $\HH^2$. Break the problem into two steps. 
  \begin{enumerate}
  \item Prove that vertical lines $L(t) = (x, t)$, $t>0$, are geodesics in $\HH^2$.
  \item Prove that semi-circles $C(t) = (x+ r\,\cos(t), r\,\sin(t))$, $t\in (0, \pi)$ are geodesics in $\HH^2$.
  \end{enumerate}
\end{exercise}

\begin{exercise}[Requires some geometry]\label{Ex:Reflection}
Suppose $C$ is a geodesic in $\HH^2$ that is a Euclidean semi-circle with center $a\in\RR$ and radius $R$. Then the \emph{reflection through $C$} \index{reflection through a geodesic} takes $z$ to $R^2/(\overline{z}-\overline{a}) +a$, where $\overline{z}$ denotes complex conjugation.

Prove the reflection through $C$ is an isometry of $\HH^2$ that fixes $C$ pointwise. Note this is an orientation reversing isometry. 

A similar result holds for reflection through a vertical line. Find a description for the reflection through a vertical line, and prove it is an isometry.
\end{exercise}

\begin{exercise}[Requires geometry]
Prove any isometry of $\HH^2$ is the product of reflections in hyperbolic geodesics. 
\end{exercise}

\begin{exercise}\label{Ex:IsomH2}
Work through the classification of isometries of $\HH^2$ as elliptic,\index{elliptic} parabolic,\index{parabolic} or loxodromic.\index{loxodromic} (E.g.\ Thurston \cite[page 67]{thurston}).
\end{exercise}

\begin{exercise}\label{Ex:IsomH3}
\Reflem{3points} shows there exists an orientation preserving isometry of $\HH^2$ taking any three points of $\bdy\HH^2$ to any other three points, provided we are careful with orientation. Prove a similar statement for $\HH^3$: Given distinct $b$, $c$ and $d$ in $\CC \cup \{\infty\}$, prove there exists an orientation preserving isometry of $\HH^3$ taking $b$ to $1$, $c$ to $0$, and $d$ to $\infty$. Write it down as a matrix in $\PSL(2,\CC)$. Note in $\HH^3$ we no longer have to worry about orientation.
\end{exercise}

\begin{exercise}
Prove the following analogue of \reflem{IntGeodesics} in $\HH^3$. Show two distinct geodesics $\ell_1$ and $\ell_2$ either intersect in a single point in the interior of $\HH^3$, intersect in a single point on $\bdy\HH^3$, or are completely disjoint in $\HH^3\cup\bdy\HH^3$. In the third case, show there exists a unique geodesic that is perpendicular to both $\ell_1$ and $\ell_2$.
\end{exercise}

\begin{exercise}[Cross ratios] \label{Ex:CrossRatio}
Given $a\in \CC$, the image of $a$ under the isometry of \refex{IsomH3} is said to be the \emph{cross ratio}\index{cross ratio} of $a,b,c,d$, and is denoted $\lambda(a,b;c,d)$.

Let $x$ be the point on the geodesic in $\HH^3$ between $c$ and $d$ such that the geodesic from $a$ to $x$ is perpendicular to that between $c$ and $d$.  Let $y$ be the point on the geodesic between $c$ and $d$ such that the geodesic from $b$ to $y$ is perpendicular to that between $c$ and $d$.  Prove the hyperbolic distance between $x$ and $y$ is equal to $|\log|\lambda(a,b;c,d)||$.
\end{exercise}

\begin{exercise}[Areas of ideal triangles]\label{Ex:IdealTriangleArea}
Prove that the area of an ideal hyperbolic triangle is $\pi$. (E.g.\ use calculus.)\index{ideal triangle}
\end{exercise}

\begin{exercise}[Areas of 2/3-ideal triangles]\label{Ex:2/3IdealTriangle}
  A \emph{$2/3$-ideal triangle}\index{$2/3$-ideal triangle} is a triangle with two vertices on the boundary at infinity\index{boundary at infinity} $\bdy\HH^2$, and the third in the interior of $\HH^2$ such that the interior angle at the third vertex is $\theta$.
  
\begin{figure}
%% Accompanies image file 'F2-06-23Tri.eps' (pdf, eps, ps)
\begingroup%
  \makeatletter%
  \providecommand\color[2][]{%
    \errmessage{(Inkscape) Color is used for the text in Inkscape, but the package 'color.sty' is not loaded}%
    \renewcommand\color[2][]{}%
  }%
  \providecommand\transparent[1]{%
    \errmessage{(Inkscape) Transparency is used (non-zero) for the text in Inkscape, but the package 'transparent.sty' is not loaded}%
    \renewcommand\transparent[1]{}%
  }%
  \providecommand\rotatebox[2]{#2}%
  \newcommand*\fsize{\dimexpr\f@size pt\relax}%
  \newcommand*\lineheight[1]{\fontsize{\fsize}{#1\fsize}\selectfont}%
  \ifx\svgwidth\undefined%
    \setlength{\unitlength}{108bp}%
    \ifx\svgscale\undefined%
      \relax%
    \else%
      \setlength{\unitlength}{\unitlength * \real{\svgscale}}%
    \fi%
  \else%
    \setlength{\unitlength}{\svgwidth}%
  \fi%
  \global\let\svgwidth\undefined%
  \global\let\svgscale\undefined%
  \makeatother%
  \begin{picture}(1,0.80000003)%
    \lineheight{1}%
    \setlength\tabcolsep{0pt}%
    \put(0,0){\includegraphics[width=\unitlength]{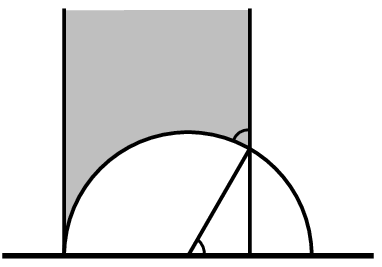}}%
    \put(0.55927577,0.16790706){\color[rgb]{0,0,0}\makebox(0,0)[lt]{\lineheight{1.25}\smash{\begin{tabular}[t]{l}$\theta$\end{tabular}}}}%
    \put(0.59151782,0.48784742){\color[rgb]{0,0,0}\makebox(0,0)[lt]{\lineheight{1.25}\smash{\begin{tabular}[t]{l}$\theta$\end{tabular}}}}%
    \put(0.08804563,0.04042769){\color[rgb]{0,0,0}\makebox(0,0)[lt]{\lineheight{1.25}\smash{\begin{tabular}[t]{l}$-1$\end{tabular}}}}%
    \put(0.79489085,0.04042769){\color[rgb]{0,0,0}\makebox(0,0)[lt]{\lineheight{1.25}\smash{\begin{tabular}[t]{l}$1$\end{tabular}}}}%
  \end{picture}%
\endgroup%

  \caption{$2/3$-ideal triangle.\index{$2/3$-ideal triangle}}
  \label{Fig:two-thirds}
\end{figure}

\begin{enumerate}
\item[(a)] Show that all $2/3$-ideal triangles of angle $\theta$ are congruent to the triangle shown in \reffig{two-thirds}, with one ideal vertex at infinity, one at $-1 \in \bdy\HH^2 = \RR\cup\{\infty\}$, and the third in the interior of $\HH^2$ with edges making angle $\theta$.

\item[(b)] Define a function $A\from (0,\pi) \to \RR$ by: $A(\theta)$ is the area of the $2/3$-ideal triangle with interior angle $\pi-\theta$. Show that
  \[ A(\theta_1 + \theta_2) = A(\theta_1) + A(\theta_2), \]
when this is defined. (Hint: Figure \ref{Fig:area} may be useful.)

\begin{figure}[h!]
%% Accompanies image file 'F2-07-Area23.eps' (pdf, eps, ps)
\begingroup%
  \makeatletter%
  \providecommand\color[2][]{%
    \errmessage{(Inkscape) Color is used for the text in Inkscape, but the package 'color.sty' is not loaded}%
    \renewcommand\color[2][]{}%
  }%
  \providecommand\transparent[1]{%
    \errmessage{(Inkscape) Transparency is used (non-zero) for the text in Inkscape, but the package 'transparent.sty' is not loaded}%
    \renewcommand\transparent[1]{}%
  }%
  \providecommand\rotatebox[2]{#2}%
  \newcommand*\fsize{\dimexpr\f@size pt\relax}%
  \newcommand*\lineheight[1]{\fontsize{\fsize}{#1\fsize}\selectfont}%
  \ifx\svgwidth\undefined%
    \setlength{\unitlength}{289.99999237bp}%
    \ifx\svgscale\undefined%
      \relax%
    \else%
      \setlength{\unitlength}{\unitlength * \real{\svgscale}}%
    \fi%
  \else%
    \setlength{\unitlength}{\svgwidth}%
  \fi%
  \global\let\svgwidth\undefined%
  \global\let\svgscale\undefined%
  \makeatother%
  \begin{picture}(1,0.21034484)%
    \lineheight{1}%
    \setlength\tabcolsep{0pt}%
    \put(0,0){\includegraphics[width=\unitlength]{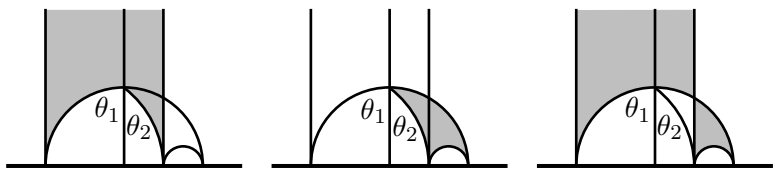}}%
    \put(0.1171777,0.06604334){\color[rgb]{0,0,0}\makebox(0,0)[lt]{\lineheight{1.25}\smash{\begin{tabular}[t]{l}$\theta_1$\end{tabular}}}}%
    \put(0.15859395,0.04280984){\color[rgb]{0,0,0}\makebox(0,0)[lt]{\lineheight{1.25}\smash{\begin{tabular}[t]{l}$\theta_2$\end{tabular}}}}%
    \put(0.46305453,0.0640812){\color[rgb]{0,0,0}\makebox(0,0)[lt]{\lineheight{1.25}\smash{\begin{tabular}[t]{l}$\theta_1$\end{tabular}}}}%
    \put(0.50447077,0.0408477){\color[rgb]{0,0,0}\makebox(0,0)[lt]{\lineheight{1.25}\smash{\begin{tabular}[t]{l}$\theta_2$\end{tabular}}}}%
    \put(0.8065064,0.0661015){\color[rgb]{0,0,0}\makebox(0,0)[lt]{\lineheight{1.25}\smash{\begin{tabular}[t]{l}$\theta_1$\end{tabular}}}}%
    \put(0.84792265,0.04286801){\color[rgb]{0,0,0}\makebox(0,0)[lt]{\lineheight{1.25}\smash{\begin{tabular}[t]{l}$\theta_2$\end{tabular}}}}%
  \end{picture}%
\endgroup%

  \caption{Areas of triangles.}
  \label{Fig:area}
\end{figure}

\item[(c)] It follows that $A$ is $\QQ$-linear.  Since $A$ is continuous, it must be $\RR$-linear.  Show $A(\theta) = \theta$.
\end{enumerate}
\end{exercise}

\begin{exercise}\label{Ex:TriangleAreas}
(Areas of general triangles.)
Using the previous two problems, show that the area of a triangle with interior angles $\alpha$, $\beta$, and $\gamma$ is equal to $\pi -\alpha -\beta -\gamma$.  Note an ideal vertex has interior angle $0$.
\end{exercise}

\begin{figure}
%% Accompanies image file 'F2-08-LabelT.eps' (pdf, eps, ps)
\begingroup%
  \makeatletter%
  \providecommand\color[2][]{%
    \errmessage{(Inkscape) Color is used for the text in Inkscape, but the package 'color.sty' is not loaded}%
    \renewcommand\color[2][]{}%
  }%
  \providecommand\transparent[1]{%
    \errmessage{(Inkscape) Transparency is used (non-zero) for the text in Inkscape, but the package 'transparent.sty' is not loaded}%
    \renewcommand\transparent[1]{}%
  }%
  \providecommand\rotatebox[2]{#2}%
  \newcommand*\fsize{\dimexpr\f@size pt\relax}%
  \newcommand*\lineheight[1]{\fontsize{\fsize}{#1\fsize}\selectfont}%
  \ifx\svgwidth\undefined%
    \setlength{\unitlength}{147.99999619bp}%
    \ifx\svgscale\undefined%
      \relax%
    \else%
      \setlength{\unitlength}{\unitlength * \real{\svgscale}}%
    \fi%
  \else%
    \setlength{\unitlength}{\svgwidth}%
  \fi%
  \global\let\svgwidth\undefined%
  \global\let\svgscale\undefined%
  \makeatother%
  \begin{picture}(1,0.75000002)%
    \lineheight{1}%
    \setlength\tabcolsep{0pt}%
    \put(0,0){\includegraphics[width=\unitlength]{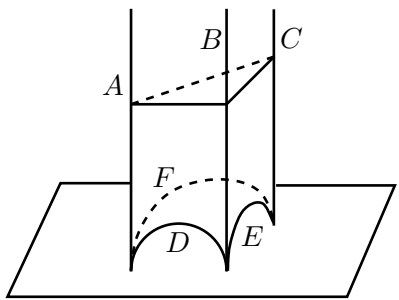}}%
    \put(0.24606419,0.51909426){\color[rgb]{0,0,0}\makebox(0,0)[lt]{\lineheight{1.25}\smash{\begin{tabular}[t]{l}$A$\end{tabular}}}}%
    \put(0.49325196,0.6348231){\color[rgb]{0,0,0}\makebox(0,0)[lt]{\lineheight{1.25}\smash{\begin{tabular}[t]{l}$B$\end{tabular}}}}%
    \put(0.69886718,0.63819385){\color[rgb]{0,0,0}\makebox(0,0)[lt]{\lineheight{1.25}\smash{\begin{tabular}[t]{l}$C$\end{tabular}}}}%
    \put(0.40673622,0.11909954){\color[rgb]{0,0,0}\makebox(0,0)[lt]{\lineheight{1.25}\smash{\begin{tabular}[t]{l}$D$\end{tabular}}}}%
    \put(0.59886855,0.13482968){\color[rgb]{0,0,0}\makebox(0,0)[lt]{\lineheight{1.25}\smash{\begin{tabular}[t]{l}$E$\end{tabular}}}}%
    \put(0.37415236,0.28426591){\color[rgb]{0,0,0}\makebox(0,0)[lt]{\lineheight{1.25}\smash{\begin{tabular}[t]{l}$F$\end{tabular}}}}%
  \end{picture}%
\endgroup%

  \caption{Dihedral angles of an ideal tetrahedron.}
  \label{Fig:labeled-tet}
\end{figure}

\begin{exercise}\label{Ex:TetLabels1}
(Ideal tetrahedra and dihedral angles.)
The dihedral angles on a tetrahedron are labeled $A$, $B$, $C$, $D$, $E$, and $F$ in \reffig{labeled-tet}. Using linear algebra, prove that opposite dihedral angles agree. That is, show $A=E$, $B=F$, and $C=D$.
\end{exercise}

\begin{exercise}\label{Ex:tet-labels}
(Ideal tetrahedra and cross ratios.)
Orient an ideal tetrahedron with vertices $a,b,c,d$. When we apply a M\"obius transformation\index{M\"obius transformation} taking $b,c,d$ to $1,0,\infty$, respectively, the point $a$ goes to the cross ratio\index{cross ratio} $\lambda(a,b;c,d)$. Label the edge from $c$ to $d$ by the complex number $\lambda = \lambda(a,b;c,d)$. We may do this for each edge of the tetrahedron, labeling by a different cross ratio. (Notice you need to keep track of orientation.) Find all labels on the edges of the tetrahedra in terms of $\lambda$.
\end{exercise}

\begin{exercise}\label{Ex:CuspVolume}(Volume of a region in a horoball)
  Let $R$ be the region in $\HH^3$ given by $A\times[1,\infty)$, where $A$ is some region contained in the horosphere about $\infty$ of height $1$, i.e.\ $A\subset \{(x+i\,y, 1)\}$. Prove that $\vol(R) = \area(A)/2$. 
\end{exercise}

%% Ch03_Geometric.tex

\chapter{Geometric Structures on Manifolds}\label{Chap:Geometric}
\blfootnote{Jessica S. Purcell, Hyperbolic Knot Theory}

In this chapter, we give our first examples of hyperbolic manifolds, combining ideas from the previous two chapters.

\section{Geometric structures}

\subsection{Introductory example: The torus}
\begin{figure}[h]
\includegraphics{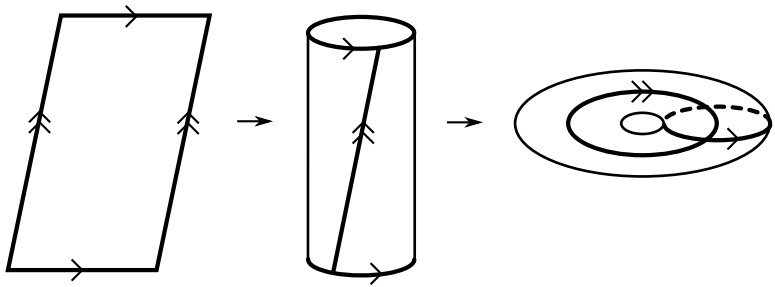}
\caption{A parallelogram glued to a torus}
\label{Fig:TorusGluing}
\end{figure}

A geometric structure you are likely familiar with is a 2-dimensional Euclidean structure\index{Euclidean structure} on a torus. Given any parallelogram, we obtain a torus by gluing the top and bottom sides of the parallelogram, and the right and left sides, as shown in \reffig{TorusGluing}.

\begin{figure}
  \includegraphics{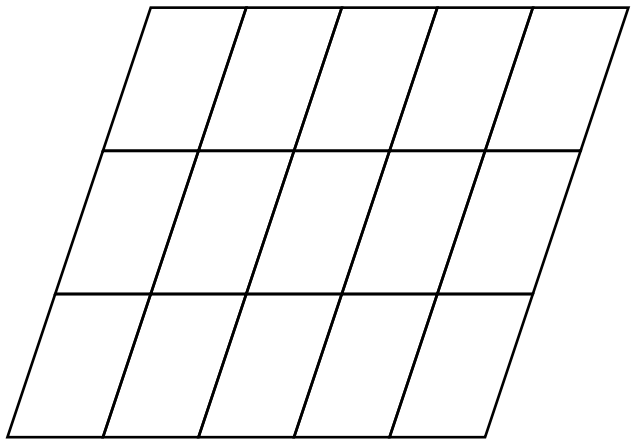}
  \caption{The universal cover of a Euclidean torus.}
  \label{Fig:EuclideanTorus}
\end{figure}

The universal cover of the torus is obtained by gluing copies of the parallelogram to itself in $\RR^2$. We may glue infinitely many copies in two directions, and we obtain a tiling of the plane $\RR^2$ by parallelograms, as in \reffig{EuclideanTorus}. These parallelograms define a lattice in $\RR^2$, and covering transformations of the universal cover $\RR^2$ of the torus are given by Euclidean translations by points of the lattice. That is, if the parallelogram is determined by vectors $\overrightarrow{v}$ and $\overrightarrow{w}$ along its sides, then any covering transformation is of the form $a\overrightarrow{v}+b\overrightarrow{w}$ for $a, b \in \ZZ$. This construction works for any choice of parallelogram.

\begin{figure}
\includegraphics{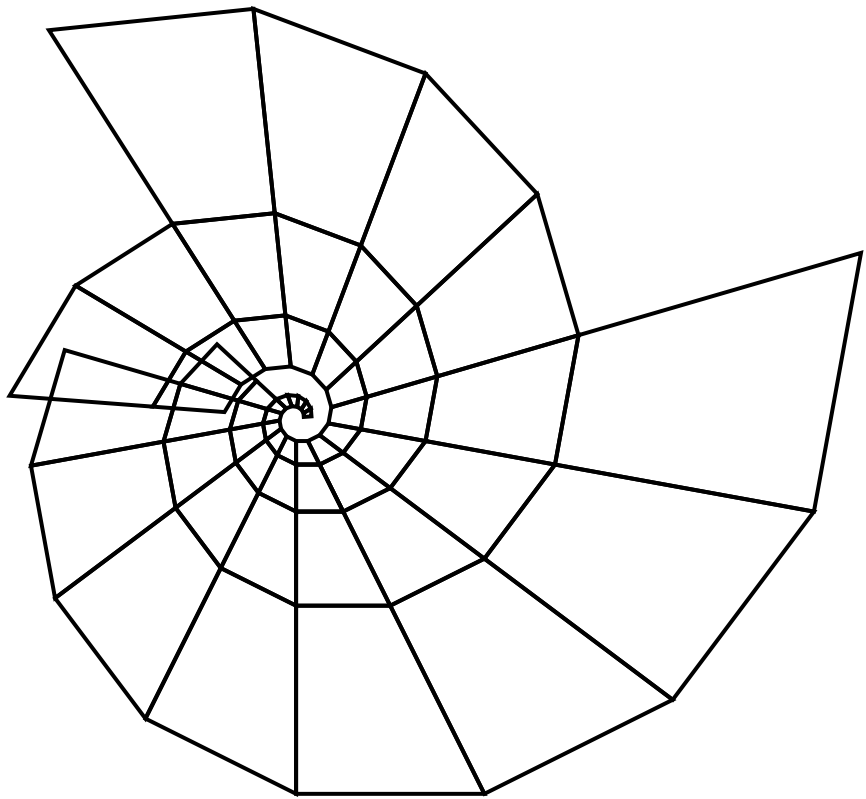}
  \caption{When we construct a torus from a quadrilateral that is not a parallelogram, generally a single point is omitted from the plane.}
  \label{Fig:AffineTorus}
\end{figure}

Now modify this construction by choosing a more general quadrilateral instead of a parallelogram. We can still identify opposite sides in an orientation preserving manner, so when we glue we still get an object homeomorphic to a torus. However, the quadrilateral no longer determines a tiling of $\RR^2$, nor a lattice. Indeed, when we glue copies of the quadrilateral to itself, as we did when constructing the universal cover above, we have to shrink, expand, and rotate the quadrilateral to glue copies, and the result is not a tiling of the plane. See \reffig{AffineTorus}.

These examples of the torus can be generalized to different surfaces and manifolds. The torus was created by gluing quadrilaterals. More generally, we
will glue different types of polygons, including ideal polygons, and
in 3-dimensions, polyhedra.

\begin{definition}\label{Def:TopPolygon}
Let $M$ be a 2-manifold.  A \emph{topological polygonal decomposition}\index{topological polygonal decomposition} of $M$ is a combinatorial way of gluing polygons so that the result is homeomorphic to $M$. 

We allow ideal polygons, i.e.\ those with one or more ideal vertex. Additionally, by \emph{gluing}\index{gluing} we mean an identification that takes faces to faces, edges to edges, and vertices to vertices.
\end{definition}

Both constructions of the torus above give examples of topological polygonal decompositions of the torus. 

\begin{definition}\label{Def:GeomPoly}
A \emph{geometric polygonal decomposition}\index{geometric polygonal decomposition} of $M$ is a topological polygonal decomposition along with a metric on each polygon such that gluing is by isometry and the result of the gluing is a smooth manifold with a complete metric.
\end{definition}

Recall that a metric space is \emph{complete}\index{complete metric space} if every Cauchy sequence converges; and recall that a \emph{Cauchy sequence}\index{Cauchy sequence} is a sequence $\{x_i\}_{i=1}^\infty$ such that for each $\epsilon>0$, there exists a positive integer $N$ such that $d(x_i,x_j)<\epsilon$ if $i,j\geq N$. 

The first construction of the torus gives a complete Euclidean metric on the torus, by pulling back the Euclidean metric on the parallelogram. Because gluings of the sides of the parallelogram are by Euclidean isometries, this will be well-defined. The second construction of the torus does not give a complete Euclidean metric, or any Euclidean metric: gluings of the quadrilaterals are by affine transformations (rotation, translation, scale),\index{affine transformation} not isometries of the Euclidean plane, so we cannot pull back a well-defined metric. Note also that toward the center of \reffig{AffineTorus}, the quadrilaterals are becoming arbitrarily small. In fact, there is a point in the figure that is disjoint from all quadrilaterals (see \refex{DevelopingImageMissesPoint}). 

We will also be studying polygonal decompositions of manifolds and their generalization to three dimensions: polyhedral decompositions.  More generally, we can discuss geometric structures on manifolds.

\subsection{Geometric structures on manifolds}

\begin{definition}\label{Def:GXStructure}
Let $X$ be a manifold, and $G$ a group acting on $X$.  We say a manifold $M$ has a $(G,X)$-structure\index{$(G,X)$-structure} if for every point $x \in M$, there exists a \emph{chart}\index{chart} $(U,\phi)$, that is, a neighborhood $U \subset M$ of $x$ and a homeomorphism $\phi\from U \to \phi(U)\subset X$. We also sometimes refer to the map $\phi$ as a chart when $U$ is understood. Charts satisfy the following: if two charts $(U, \phi)$ and $(V, \psi)$ overlap, then the \emph{transition map}\index{transition map} or \emph{coordinate change map}\index{coordinate change map}
\[
\gamma = \phi \circ \psi^{-1}\from \psi(U\cap V) \to \phi(U\cap V)
\]
is an element of $G$.
\end{definition}

In the examples we encounter here, $X$ will be simply connected, and $G$ a group of \emph{real analytic diffeomorphisms}\index{real analytic diffeomorphism}\index{diffeomorphism!real analytic} acting transitively on $X$.
The reason we need real analytic diffeomorphisms is that they are uniquely 
determined by their restriction to any open set. This is true, for example, of isometries of Euclidean space, and isometries of hyperbolic space. While we present the results in this full generality, the reader who is unfamiliar with real analytic diffeomorphisms can read with Euclidean or hyperbolic isometries in mind. 

Our manifold $X$ will typically admit a known metric as well, and $G$ will be the group of isometries of $X$. It will follow that $M$ inherits a metric from $X$ (\refex{InheritMetric}). 
We will say that $M$ has a \emph{geometric structure}\index{geometric structure}. 

\begin{example}[Euclidean torus]
\label{Example:EuclTorus}

Let $X$ be 2-dimensional Euclidean space, $\EE^2$.  Let $G$ be isometries of Euclidean space $\Isom(\EE^2)$.  The torus admits an $(\Isom(\EE^2),\EE^2)$-structure, also called a \emph{Euclidean structure}.\index{Euclidean structure}

To help us understand the definition, let's look at some charts\index{chart} and
transition maps for this example.

We know the universal cover of the torus is given by tiling the plane $\RR^2$ with parallelograms. For simplicity, we will work with the example in which each parallelogram is a square, and one square has vertices $(0,0)$, $(1,0)$, $(1,1)$, and $(0,1)$ in $\RR^2$. Call this square the \emph{basic} square.

Now pick any point $p$ on the torus.  This will lift to a collection of points on $\RR^2$, one for each copy of the unit square.  Take a disk of radius $1/4$ around each lift. These all project under the covering map to an open neighborhood $U$ of $p$ in the torus. Therefore we have the following charts: $(U,\phi)$ is a chart,\index{chart} where $\phi$ maps $U$ into the disk of radius $1/4$ centered around the lift $\widehat{p}_0$ of $p$ in the basic square.  Another chart is $(U, \psi)$, where $\psi$ maps $U$ into the disk of radius $1/4$ about the lift $\widehat{p}_1$ of $p$ in some other square. Such a lift is given by a translation of $\widehat{p}_0$ by a vector $(m,n)\in\ZZ\times\ZZ$, in the lattice determined by the basic square. Thus $\phi \circ \psi^{-1}$ will be a Euclidean translation by integral values in the $x$ and $y$ direction.  These are Euclidean isometries.

More generally, let $q$ be a point such that a lift $\widehat{q}_0$ of $q$ has distance less than $1/2$ to $\widehat{p}_0$ in the basic square. Thus a disk of radius $1/4$ about $\widehat{p}_0$ overlaps with a disk of radius $1/4$ about $\widehat{q}_0$. These disks project to give open neighborhoods $U$ and $V$ of $p$ and $q$ respectively in the torus.  Since these neighborhoods overlap, we need to ensure that any corresponding charts\index{chart} differ by a Euclidean isometry in the region of overlap. Obtain charts by mapping $U$ to your favorite disk of radius $1/4$ about a lift of $p$ in $\RR^2$.  Map $V$ to your favorite disk of radius $1/4$ about a lift of $q$ in $\RR^2$; see \reffig{GXTorus} for an example. Again, regardless of the choice of $\phi$ and $\psi$, the overlap
\[
\phi \circ \psi^{-1}\from \psi(U\cap V)\to \phi(U\cap V)
\]
will be a Euclidean translation of the intersection of the two disks by some $(n,m)\in\ZZ\times\ZZ$ corresponding to the choice of lifts. Again see \reffig{GXTorus}. 

\begin{figure}
  \begin{center}
    \includegraphics{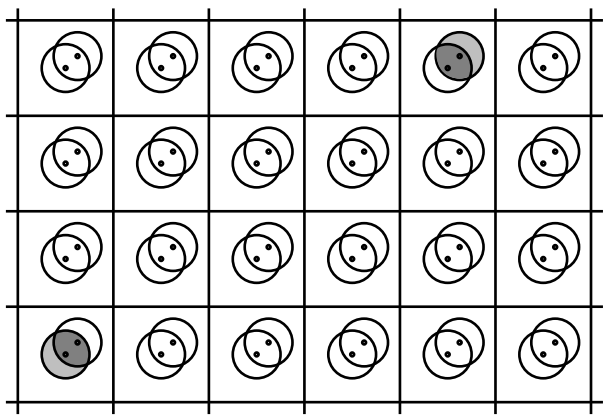}
  \end{center}
  \caption{Euclidean structure\index{Euclidean structure} on a torus: Transition maps are Euclidean translations.}
  \label{Fig:GXTorus}
\end{figure}

This idea extends to arbitrary neighborhoods $U$ and $V$: transition maps will always be translations by $(n,m)\in \ZZ\times\ZZ$. 
Therefore, we conclude that the torus obtained by gluing sides of the square with vertices $(0,0)$, $(0,1)$, $(1,0)$ and $(1,1)$ admits an $(\Isom(\EE^2),\EE^2)$-structure, where $\EE^2$ denotes $\RR^2$ with the standard Euclidean metric. 
\end{example}

\begin{example}[The affine torus]
\label{Example:AffineTorus}\index{affine torus}
Again let $X=\RR^2$, but this time let $G$ be the affine group acting on $\RR^2$.  That is, $G$ consists of invertible affine transformations,\index{affine transformation} i.e.\ linear transformations followed by a translation:
\[
x\mapsto Ax+b.
\]

The torus of \reffig{AffineTorus} admits a $(G, \RR^2)$-structure.  This can be seen in a manner similar to that in the previous example.  Charts\index{chart} will differ by a scaling, rotation, then translation.  
\end{example}

\smallskip

In practice, we rarely use charts\index{chart} to show manifolds have a particular $(G,X)$-structure.\index{$(G,X)$-structure}  Instead, as in the two previous examples, we build manifolds by starting with an existing manifold $X$ and taking the quotient by the action of a group, or by gluing together polygons.

\subsection{Hyperbolic surfaces}

Let $X = \HH^2$, and let $G = \Isom(\HH^2)$, the group of isometries of $\HH^2$.  When a 2-manifold admits an $(\Isom(\HH^2), \HH^2)$-structure, we say the manifold admits a \emph{hyperbolic structure}\index{hyperbolic structure, definition}, or is hyperbolic. More generally, an $n$-manifold that admits an $(\Isom(\HH^n), \HH^n)$-structure admits a \emph{hyperbolic structure}, or is hyperbolic. 

We will look at some examples of hyperbolic 2-manifolds obtained from geometric polygonal decompositions.  To do so, we start with a collection of hyperbolic polygons in $\HH^2$, for example, a collection of triangles. We allow vertices to either be finite or ideal, i.e.\ in the interior of $\HH^2$ or on $\bdy_\infty\HH^2$, respectively. In any case, we will always assume each polygon is convex, and edges are segments of geodesics in $\HH^2$.  Now, to each edge, associate exactly one other edge.  Just as in the case of the torus, glue polygons along associated edges by an isometry of $\HH^2$. 

When does the result of this gluing give a manifold that admits a hyperbolic structure?  We obtain a hyperbolic structure exactly when each point in the result has a neighborhood $U$ and a homeomorphism into $\HH^2$ so that transition maps are in $\Isom(\HH^2)$. The following lemma gives a condition that will guarantee this.

\begin{lemma}\label{Lem:IsometricNbhd}
A gluing of hyperbolic polygons yields a 2-manifold with a hyperbolic structure, with structure agreeing with that in the interior of the polygons, if and only if each point in the gluing has a neighborhood (in the quotient topology) isometric to a disk in $\HH^2$.

More generally, a gluing of $n$-dimensional hyperbolic polyhedra yields a hyperbolic $n$-manifold, with hyperbolic structure agreeing with that in the interior of the polyhedra, if and only if each point has a neighborhood (in the quotient topology) isometric to a ball in $\HH^n$, with the isometry the identity in the interior of polyhedra.
\end{lemma}

Here by a \emph{gluing}\index{gluing} of hyperbolic polyhedra, we mean a collection of geodesic polyhedra embedded in $\HH^n$, along with identifications on faces, called gluing maps or face-pairings,\index{face-pairing isometry} which are given by an isometry on each face. The quotient space of the polyhedra with identifications given by the gluing maps is the gluing.

Additionally, we say that the hyperbolic structure agrees with the structure in the interior of the polyhedra if, for any point in the interior of the polyhedron, a ball $U$ containing that point, lying in the interior of the polyhedron, along with the identity map from $U$ to $U\subset\HH^3$, provides a chart\index{chart} in the hyperbolic structure.

\begin{proof}[Proof of \reflem{IsometricNbhd}]
We will prove the more general statement. Suppose first that a gluing of hyperbolic polyhedra $M$ yields an $n$-manifold with hyperbolic structure, agreeing with the hyperbolic structure in the interior of the polyhedra. Then every point $x$ in $M$ has a neighborhood $U$ and a chart\index{chart} $\phi\from U\to \phi(U)\subset \HH^n$ such that transition maps are isometries of $\HH^n$. By restricting $\phi$ to a subset of $U$, we may assume $\phi(U)$ is a ball in $\HH^n$. The neighborhood $U$ is open in the quotient topology on the gluing. Thus it is made up of portions of open neighborhoods meeting the polyhedra in $\HH^n$, identified by gluing isometries. In the interior of a polyhedron $P$, $\phi$ composed with the identity map on $U\cap P$ is an isometry of $\HH^n$. Thus we may view $U\cap P$ as the intersection of a hyperbolic ball with $P$. Since gluing maps are isometries, they identify faces of $U\cap P$ into a hyperbolic ball, and $\phi$ must be an isometry of $U$ into a ball in $\HH^n$.

Now suppose that under the quotient topology, every point of $M$ has a neighborhood isometric to a ball of $\HH^n$, with isometry the identity for points in the interior of a polyhedron. Then this isometry gives a chart\index{chart} $\phi\from U\to\phi(U)\subset\HH^n$. If $(U,\phi)$ and $(V,\psi)$ are charts\index{chart} and $U\cap V\neq\emptyset$, then $\phi\circ\psi^{-1}\from \psi(U\cap V)\to \phi(U\cap V)$ is the composition of isometries, hence an isometry, so $M$ has an $(\Isom(\HH^n),\HH^n)$-structure. Because charts\index{chart} in the interior of polyhedra are identity maps, the hyperbolic structure agrees with that on the polyhedra.
\end{proof}

When does each point in a gluing of hyperbolic polygons have a neighborhood isometric to a disk in the hyperbolic plane?  Let $x$ be a point in the gluing, and consider its lifts to the polygons. There are three cases. 
\begin{enumerate}
\item If $x$ lifts to a point $\widehat{x}$ in the interior of one of the polygons, then that lift is unique. In this case, for small enough $\epsilon>0$, there is a disk about $\widehat{x}$ of radius $\epsilon$ embedded in the interior of the polygon in $\HH^2$. This projects under the quotient map to a disk about $x$ isometric to a disk in $\HH^2$. 

\item If $x$ lifts to a point on an edge of a polygon, then it has two lifts, $\widehat{x}_0$ and $\widehat{x}_1$, on two different edges that are glued to each other by the gluing map. A neighborhood of $x$ in the quotient topology lifts to give a ``half-neighborhood'' of $\widehat{x}_0$ glued to a corresponding ``half-neighborhood'' of $\widehat{x}_1$. Each contains a half-disk in $\HH^2$, and we may scale the disks so that they glue to a disk under the gluing map. Thus in this case as well, $x$ has a neighborhood isometric to a disk in $\HH^2$. 

\item If $x$ lifts to a finite vertex of a polygon, then it may have several lifts, possibly including several vertices of the collection of polygons. In this case, we need to be more careful. The following lemma gives a condition that will guarantee we have an isometry to a hyperbolic disk in this case as well. 
\end{enumerate}

\begin{lemma}\label{Lem:VertexSum}
A gluing of hyperbolic polygons gives a 2-manifold with a hyperbolic structure if and only if for each finite vertex $v$ of the polygons, the sum of interior angles at each vertex glued to $v$ is $2\pi$.
\end{lemma}

\begin{proof}
  This is an immediate consequence of \reflem{IsometricNbhd} and the observation that around a vertex, portions of the polygons meet in a cycle, with total angle around the finite vertex equal to the sum of interior angles of the polyhedron at that vertex. We need to check that each finite vertex has a neighborhood isometric to a neighborhood in $\HH^2$. This will hold if and only if the sum of interior angles is $2\pi$. 
\end{proof}

\section{Complete structures}

Given a gluing of hyperbolic polygons, suppose the angle sum at each finite vertex is $2\pi$, so that we have a hyperbolic structure by \reflem{VertexSum}. Does it necessarily follow that we have a geometric polygonal decomposition?

Recall from \refdef{GeomPoly} that for a geometric polygonal decomposition, we need a geometric structure on each polygon so that the result of the gluing is a smooth manifold with a complete metric.\index{complete metric space} Our hyperbolic structure gives a smooth manifold with a metric. However, in the presence of ideal vertices, the metric may not be complete.

It will be easier to discuss criteria for completeness using the language of \emph{developing maps} and \emph{holonomy}. Our exposition of these terms is based on that of Thurston \cite{thurston:book}.

\subsection{Developing map and holonomy}

The developing map, which we define in this subsection, encodes information on the $(G,X)$-structure of a manifold.\index{$(G,X)$-structure} It is a local homeomorphism into $X$. When a manifold has a polygonal decomposition, say by polygons in $X=\RR^2$ or $\HH^2$, the developing map ``develops'' the gluing information on the polygons by attaching copies of the polygons along edges in the space $X$, as we did for the torus in Figures~\ref{Fig:EuclideanTorus} and~\ref{Fig:AffineTorus}.

More generally, a developing map can be defined for any manifold $M$ with a $(G,X)$-structure,\index{$(G,X)$-structure} assuming as before that $X$ is a manifold and $G$ is a group of real analytic diffeomorphisms acting transitively on $X$. Any chart\index{chart} $(U,\phi)$ gives a homeomorphism of $U$ onto $\phi(U)\subset X$. To define the developing map, we wish to extend this map.

Suppose $(V,\psi)$ is another chart,\index{chart} and $y\in U\cap V$. Then
\[ \gamma = \phi\circ\psi^{-1}\from \psi(U\cap V)\to \phi(U\cap V) \]
is an element of $G$ acting on $\psi(U\cap V)$. By setting $y\mapsto \gamma$, we obtain a map from $U\cap V$ to $G$. 
Because $G$ is a group of real analytic diffeomorphisms, the element $\gamma$ is uniquely determined in a neighborhood of $\psi(y)$. This implies that the map $y\mapsto \gamma$ is locally constant: we obtain the same element $\gamma$ for all $x$ in a neighborhood of $y$ in $U\cap V$. We let $\gamma(y)$ denote this element of $G$. Then we may define a map $\Phi\from U\cup V \to X$ by
\[
\Phi(x) = \begin{cases} \phi(x) & \mbox{ if } x\in U \\
\gamma(y)\cdot\psi(x) \quad & \mbox{ if } x\in V
\end{cases}
\]
Note that if $U\cap V$ is connected, then $\Phi$ is a well-defined homeomorphism, since for $x\in U\cap V$, we have $\phi(x)= \gamma(y)\cdot\psi(x)$. Thus in this case, $\Phi$ is an extension of $\phi$. However, note that we may run into trouble when $U\cap V$ is not connected, as follows. If $x$ is in a component disjoint from that containing $y$, then $\phi(x)$ may not equal $\gamma(y)\cdot\psi(x)$. This is illustrated in the following example. 

\begin{example}\label{Example:TorusDeveloping}
Consider the Euclidean torus obtained by gluing sides of a square with vertices $(0,0)$, $(1,0)$, $(0,1)$, and $(1,1)$ in $\RR^2$. Suppose the union of two simply connected neighborhoods $U$ and $V$ forms a neighborhood of a longitude for the torus, as in \reffig{TorusDeveloping}, such that $U\cap V$ has two components. Suppose $y$ lies in one component of $U\cap V$. There exist charts\index{chart} $(U,\phi)$ and $(V,\psi)$ sending $y$ to the interior of the basic square in $\RR^2$. Then the transition map $\gamma(y)$ is the identity element of $G$ in this case, since $\phi(U)$ and $\psi(V)$ overlap in the component of $U\cap V$ containing $y$. However, the map $\Phi$ defined above is not well-defined, for if $x$ lies in the other component of $U\cap V$, $\phi(x)$ lies in the basic square, but $\gamma(y)\cdot\psi(x)= \psi(x)$ lies in the square with vertices $(1,0)$, $(2,0)$, $(1,1)$, and $(1,2)$. 

\begin{figure}
%% Accompanies image file 'F3-05-TorDev.eps' (pdf, eps, ps)
\begingroup%
  \makeatletter%
  \providecommand\color[2][]{%
    \errmessage{(Inkscape) Color is used for the text in Inkscape, but the package 'color.sty' is not loaded}%
    \renewcommand\color[2][]{}%
  }%
  \providecommand\transparent[1]{%
    \errmessage{(Inkscape) Transparency is used (non-zero) for the text in Inkscape, but the package 'transparent.sty' is not loaded}%
    \renewcommand\transparent[1]{}%
  }%
  \providecommand\rotatebox[2]{#2}%
  \newcommand*\fsize{\dimexpr\f@size pt\relax}%
  \newcommand*\lineheight[1]{\fontsize{\fsize}{#1\fsize}\selectfont}%
  \ifx\svgwidth\undefined%
    \setlength{\unitlength}{331.61695862bp}%
    \ifx\svgscale\undefined%
      \relax%
    \else%
      \setlength{\unitlength}{\unitlength * \real{\svgscale}}%
    \fi%
  \else%
    \setlength{\unitlength}{\svgwidth}%
  \fi%
  \global\let\svgwidth\undefined%
  \global\let\svgscale\undefined%
  \makeatother%
  \begin{picture}(1,0.37704391)%
    \lineheight{1}%
    \setlength\tabcolsep{0pt}%
    \put(0,0){\includegraphics[width=\unitlength]{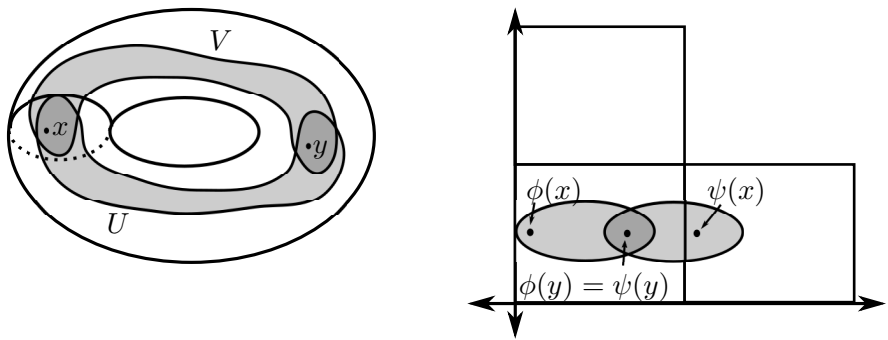}}%
    \put(0.22951926,0.3313478){\color[rgb]{0,0,0}\makebox(0,0)[lt]{\lineheight{0}\smash{\begin{tabular}[t]{l}$V$\end{tabular}}}}%
    \put(0.11340461,0.12321737){\color[rgb]{0,0,0}\makebox(0,0)[lt]{\lineheight{0}\smash{\begin{tabular}[t]{l}$U$\end{tabular}}}}%
    \put(0.0514575,0.23497564){\color[rgb]{0,0,0}\makebox(0,0)[lt]{\lineheight{0}\smash{\begin{tabular}[t]{l}$x$\end{tabular}}}}%
    \put(0.34750489,0.21463658){\color[rgb]{0,0,0}\makebox(0,0)[lt]{\lineheight{0}\smash{\begin{tabular}[t]{l}$y$\end{tabular}}}}%
    \put(0.58190519,0.05333989){\color[rgb]{0,0,0}\makebox(0,0)[lt]{\lineheight{0}\smash{\begin{tabular}[t]{l}$\phi(y)=\psi(y)$\end{tabular}}}}%
    \put(0.59021567,0.16224622){\color[rgb]{0,0,0}\makebox(0,0)[lt]{\lineheight{0}\smash{\begin{tabular}[t]{l}$\phi(x)$\end{tabular}}}}%
    \put(0.79576659,0.16132359){\color[rgb]{0,0,0}\makebox(0,0)[lt]{\lineheight{0}\smash{\begin{tabular}[t]{l}$\psi(x)$\end{tabular}}}}%
  \end{picture}%
\endgroup%

\caption{Neighborhoods $U$ and $V$ on the torus have two components of intersection, one containing $x$ and one containing $y$. The map $\phi$ cannot be extended over $V$ because it will not be well-defined at $x$}
\label{Fig:TorusDeveloping}
\end{figure}
\end{example}

Similarly, as we attempt to extend $\Phi$ by considering other coordinate neighborhoods overlapping $U$ and $V$, the natural extensions using transition maps such as $\gamma(y)$ again may not be well-defined.

To overcome this problem, we use the universal cover of $M$. Recall from algebraic topology that the universal cover $\widetilde{M}$ of $M$ can be defined to be the space of homotopy classes of paths in $M$ that start at a fixed basepoint $x_0$. See, for example \cite[Theorem~82.1]{Munkres:Topology} or \cite[page~64]{Hatcher:AlgebraicTopology}. Let $\alpha\from[0,1]\to M$ be a path representing a point $[\alpha]\in\widetilde{M}$, and let the chart\index{chart} $(U_0,\phi_0)$ contain the basepoint $x_0$.

Now find $0=t_0<t_1<\dots<t_n=1$ and charts\index{chart} $(U_i,\phi_i)$ such that $\alpha([t_i,t_{i+1}])$ is contained in $U_i$ for $i=0, 1, \dots, n-1$. Denote the points $\alpha(t_i)$ by $x_i\in M$. We extend $\phi_0$ to all of $\alpha$ as follows. First, note that each $x_i$, for $i=1, \dots, n-1$, is contained in a connected component of the intersection of two charts,\index{chart} $x_i \in U_{i-1}\cap U_i$. Then the transition map $\gamma_{i-1,i} = \phi_{i-1}\circ\phi_i^{-1}$ gives an element $\gamma_{i-1,i}(x_i)$ in $G$ that is well-defined on the entire connected component. Thus at the first step, we may extend $\phi_0$ to a function from $[0,t_2]$ to $X$ by defining $\Phi_1\from [0,t_2] \to X$ to be the function:
\[
\Phi_1(t) = \begin{cases}
\phi_0(\alpha(t)) & \mbox{ if } t\in[0,t_1] \\
\gamma_{0,1}(x_1)\cdot\phi_1(\alpha(t)) & \mbox{ if } t\in [t_1,t_2]
\end{cases}
\]
This will be well-defined on all of $[0,t_2]$, since $\phi_0(\alpha(t_1)) = \gamma_{0,1}(x_1)\cdot\phi_1(\alpha(t_1))$.

Extend inductively to $\Phi_i\from[0,t_{i+1}]\to X$ by setting:
\[
\Phi_i(t) = \begin{cases}
\Phi_{i-1}(t) & \mbox{ if } t\in [0,t_i] \\
\gamma_{0,1}(x_1)\gamma_{1,2}(x_2)\dots\gamma_{(i-1),i}(x_i)\cdot\phi_i(\alpha(t)) & \mbox{ if } t\in[t_i,t_{i+1}]
\end{cases}
\]
Again this is well-defined, for we know $\phi_{i-1}(\alpha(t_i)) = \gamma_{(i-1),i}(x_i)\cdot\phi_i(\alpha(t_i))$. Thus by induction, at the point $t_i$, 
\begin{align*}
  \Phi_{i-1}(t_i) = & \:\gamma_{0,1}(x_1)\gamma_{1,2}(x_2)\dots\gamma_{(i-2),(i-1)}(x_{i-1})\cdot\phi_{i-1}(\alpha(t_i)) \\
  = & \:\gamma_{0,1}(x_1)\gamma_{1,2}(x_2)\dots\gamma_{(i-1),i}(x_i)\cdot\phi_i(\alpha(t_i)).
\end{align*}

After the $(n-1)$-st step, we have a map $\Phi_{n-1} \from [0,1]\to X$. In fact, note that the definition of $\Phi_{n-1}$ actually provides a map $\Phi_{[\alpha]}\from U\to X$, for some small neighborhood $U$ of $\alpha(1)$, defined by
\[
\Phi_{[\alpha]}(x) = \gamma_{0,1}(x_1)\gamma_{1,2}(x_2)\dots\gamma_{(n-2),(n-1)}(x_n)\cdot\phi_{n-1}(x).
\]
The function $\Phi_{[\alpha]}$ defined in this manner, with fixed initial chart\index{chart} $(U_0,\phi_0)$ and fixed basepoint $x_0$, is an example of a function defined by \emph{analytic continuation}\index{analytic continuation}. It is well known that analytic continuation gives a well-defined function, independent of choice of the charts\index{chart} $(U_1,\phi_1)$, $\dots$, $(U_{n-1}, \phi_{n-1})$, independent of the choice of points $t_1, \dots, t_{n-1}$, and independent of the choice of path $\alpha$ in the homotopy class $[\alpha]\in\widetilde{M}$. For our particular application, we will leave this as an exercise (\refex{AnalyticContinuation}). 

\begin{definition}\label{Def:DevelopingMap}
The \emph{developing map}\index{developing map} $D\from\widetilde{M}\to X$ is the map
\[
D([\alpha]) = \Phi_n(1) =
\gamma_{0,1}(x_1)\gamma_{1,2}(x_2)\dots\gamma_{(n-2),(n-1)}(x_{n-1})\cdot\phi_{n-1}(\alpha(1)),
\]
with notation given above. 
\end{definition}

\begin{proposition}\label{Prop:DevelopingMapProperties}
The developing map $D\from \widetilde{M}\to X$ satisfies the following properties.
\begin{enumerate}
\item\label{Itm:DevMapWellDefined} For fixed basepoint $x_0$ and initial chart\index{chart} $(U_0,\phi_0)$, with $x_0\in U_0$, the map $D$ is well-defined, independent of all other choices used to define it, including charts,\index{chart} points in the intersection of chart neighborhoods, and independent of choice of $\alpha$ in the homotopy class of $[\alpha]$.
\item\label{Itm:DevMapLocalHomeo} $D$ is a local diffeomorphism.
\item\label{Itm:DevMapBaseptChange} If we define a new map in the same way as $D$, except beginning with a new choice of basepoint and initial chart,\index{chart} the resulting map is equal to the composition of $D$ with an element of $G$. 
\end{enumerate}
\end{proposition}

\begin{proof}
Showing the map is well-defined, part \eqref{Itm:DevMapWellDefined}, is a standard exercise in analytic continuation,\index{analytic continuation} and uses heavily the fact that $G$ is analytic. We leave it as \refex{AnalyticContinuation}. We also leave part  \eqref{Itm:DevMapBaseptChange} as an exercise. Part \eqref{Itm:DevMapLocalHomeo} follows from part \eqref{Itm:DevMapWellDefined}, the fact that each $\gamma_{i,(i+1)}(x_i)$ is a diffeomorphism and $\phi_n$ is a local diffeomorphism on $M$, and the topology on $\widetilde{M}$.
\end{proof}

Now consider the case that $[\alpha]\in\widetilde{M}$ is an element of the fundamental group of $M$. That is, $[\alpha]$ is a homotopy class of loops starting and ending at $x_0$. Analytic continuation\index{analytic continuation} along a loop gives a function $\Phi_{[\alpha]}$ whose domain is a neighborhood of the basepoint of the loop; this is a new chart\index{chart} defined in a neighborhood of the basepoint. Since $\phi_0$ and $\Phi_{[\alpha]}$ are both charts\index{chart} defined in a neighborhood of the basepoint, these maps must differ by an element of $G$.  Let $g_{[\alpha]} \in G$ be the element such that $\Phi_{[\alpha]} = g_{[\alpha]} \phi_0$.

Let $T_{[\alpha]}$ denote the covering transformation of $\widetilde{M}$ that corresponds to $[\alpha]$. It follows that
\[
D\circ T_{[\alpha]} = g_{[\alpha]} \circ D.\]
Note also that for $[\alpha], [\beta] \in \widetilde{M}$,
\[
D\circ T_{[\alpha]}\circ T_{[\beta]} = (g_{[\alpha]}\circ D) \circ T_{[\beta]} = g_{[\alpha]}\circ g_{[\beta]} \circ D. 
\]
It follows that the map $\rho\from \pi_1(M) \to G$ defined by $\rho([\alpha]) = g_{[\alpha]}$ is a group homomorphism.

\begin{definition}\label{Def:holonomy}
The element $g_{[\alpha]}$ is the \emph{holonomy}\index{holonomy} of $[\alpha]$.  The group homomorphism $\rho$ is called the \emph{holonomy} of $M$.  Its image is the \emph{holonomy group}\index{holonomy group} of $M$.  
\end{definition}

Note that $\rho$ depends on the choices from the construction of $D$. When $D$ changes, $\rho$ changes by conjugation in $G$ (exercise).

\begin{example}\label{Example:DevelopTorus}
Pick a point $x$ on the torus, say $x$ lies at the intersection of a choice of meridian and longitude curves for the torus, and consider a nontrivial curve $\gamma$ based at $x$.  An example of a nontrivial curve $\gamma$ on the torus is shown in \reffig{TorusCurve}.

\begin{figure}
  \includegraphics{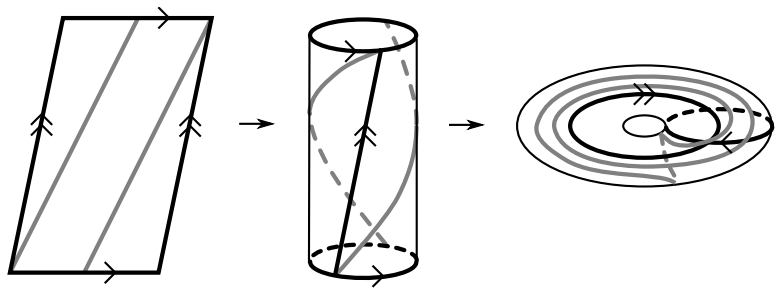}
  \caption{A nontrivial curve $\gamma$ (gray) on the torus. Meridian
    and longitude curves are shown in black.}
  \label{Fig:TorusCurve}
\end{figure}

Now consider a Euclidean structure\index{Euclidean structure} on the torus.  There exists a chart\index{chart} mapping $x$ onto the Euclidean plane.  We can take our chart to be an open parallelogram about $x$, where boundaries of the parallelogram glue in the usual way to form the torus.  As the curve $\gamma$ passes over a meridian or longitude, in the image of the developing map we must glue a new parallelogram to the appropriate side of the parallelogram we just left.  See \reffig{EuclidTorusDevel}, left, for an example.  The tiling of the plane by parallelograms is the image of the developing map, or the developing image of the Euclidean torus.

\begin{figure}
  \includegraphics{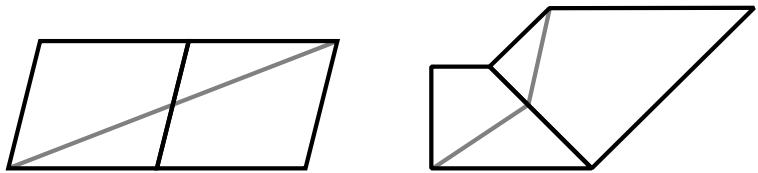}
  \caption{Left to right: developing a Euclidean torus, developing an
    affine torus.}
  \label{Fig:EuclidTorusDevel}
\end{figure}

As for the affine torus, \refexamp{AffineTorus},\index{affine torus} each time a curve crosses a meridian or longitude we attach a rescaled, rotated, translated copy of our quadrilateral to the appropriate edge. \Reffig{EuclidTorusDevel} right shows an example. \Reffig{AffineTorus} shows (part of) the developing image of the affine torus.
\end{example}

\subsection{Completeness of polygonal gluings}

Now we return to the question of determining when a gluing of hyperbolic polygons gives a complete hyperbolic structure.\index{complete metric space} We know there will be a hyperbolic structure provided the angle sum around finite vertices is $2\pi$ (\reflem{VertexSum}). The question of whether the structure is complete or not depends on what happens near ideal vertices. 

Let $M$ be an oriented hyperbolic surface obtained by gluing \emph{ideal} hyperbolic polygons. An \emph{ideal vertex}\index{ideal vertex} of $M$ is
an equivalence class of ideal vertices of the polygons, identified by the gluing.

\begin{figure}
  \begin{center}
%% Accompanies image file 'F3-08-DistD.eps' (pdf, eps, ps)
\begingroup%
  \makeatletter%
  \providecommand\color[2][]{%
    \errmessage{(Inkscape) Color is used for the text in Inkscape, but the package 'color.sty' is not loaded}%
    \renewcommand\color[2][]{}%
  }%
  \providecommand\transparent[1]{%
    \errmessage{(Inkscape) Transparency is used (non-zero) for the text in Inkscape, but the package 'transparent.sty' is not loaded}%
    \renewcommand\transparent[1]{}%
  }%
  \providecommand\rotatebox[2]{#2}%
  \newcommand*\fsize{\dimexpr\f@size pt\relax}%
  \newcommand*\lineheight[1]{\fontsize{\fsize}{#1\fsize}\selectfont}%
  \ifx\svgwidth\undefined%
    \setlength{\unitlength}{138.09152699bp}%
    \ifx\svgscale\undefined%
      \relax%
    \else%
      \setlength{\unitlength}{\unitlength * \real{\svgscale}}%
    \fi%
  \else%
    \setlength{\unitlength}{\svgwidth}%
  \fi%
  \global\let\svgwidth\undefined%
  \global\let\svgscale\undefined%
  \makeatother%
  \begin{picture}(1,0.81696145)%
    \lineheight{1}%
    \setlength\tabcolsep{0pt}%
    \put(0,0){\includegraphics[width=\unitlength]{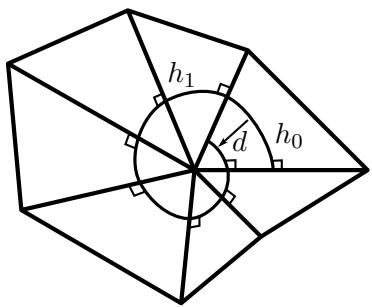}}%
    \put(0.61831948,0.42284852){\color[rgb]{0,0,0}\makebox(0,0)[lt]{\lineheight{0}\smash{\begin{tabular}[t]{l}$d$\end{tabular}}}}%
    \put(0.73450254,0.44277028){\color[rgb]{0,0,0}\makebox(0,0)[lt]{\lineheight{0}\smash{\begin{tabular}[t]{l}$h_0$\end{tabular}}}}%
    \put(0.44173597,0.60932648){\color[rgb]{0,0,0}\makebox(0,0)[lt]{\lineheight{0}\smash{\begin{tabular}[t]{l}$h_1$\end{tabular}}}}%
  \end{picture}%
\endgroup%

  \end{center}
\caption{Extending a horocycle: view inside the manifold.}
\label{Fig:distd}
\end{figure}

Let $v$ be an ideal vertex of $M$. Then $v$ is identified to some ideal vertex $v_0$ of a polygon $P_0$. Let $h_0$ be a horocycle centered at $v_0$ on $P_0$, and extend $h_0$ counterclockwise around $v_0$. The horocycle $h_0$ will meet an edge $e_0$ of $P_0$, which is glued to an edge of some polygon $P_1$ meeting ideal vertex $v_1$ identified to $v$. Note $h_0$ meets $e_0$ at a right angle. It extends to a unique horocycle $h_1$ about $v_1$ in $P_1$. Continue extending the horocycle in this manner, obtaining horocycles $h_2, h_3,\dots$. Since we only have a finite number of polygons with a finite number of vertices, eventually we return to the vertex $v_0$ of $P_0$, obtaining a horocycle $h_n$ about that ideal vertex. Note $h_n$ may not agree with the initial horocycle $h_0$. See \reffig{distd}.

\begin{definition}\label{Def:SignedDistHorocycle}
Let $d(v)$\index{$d(v)$} denote the signed hyperbolic distance between $h_0$ and $h_n$ on $P_0$. See \reffig{distd}. The sign is taken such that if $h_n$ is closer to $v_0$ than $h_0$, then $d(v)$ is positive. This is the direction shown in the figure. 
\end{definition}

\begin{lemma}\label{Lem:dvInd-h0}
The value $d(v)$\index{$d(v)$} does not depend on the initial choice of horocycle $h_0$, nor on the initial choice of $v_0$ in the equivalence class of $v$.
\end{lemma}

\begin{proof}
Exercise.
\end{proof}

It may be easier to compute $d(v)$\index{$d(v)$} if we look at polygons in $\HH^2$, using terminology of developing map and holonomy.\index{holonomy}

\begin{figure}
\begin{center}
%% Accompanies image file 'F3-09-ExtH.eps' (pdf, eps, ps)
\begingroup%
  \makeatletter%
  \providecommand\color[2][]{%
    \errmessage{(Inkscape) Color is used for the text in Inkscape, but the package 'color.sty' is not loaded}%
    \renewcommand\color[2][]{}%
  }%
  \providecommand\transparent[1]{%
    \errmessage{(Inkscape) Transparency is used (non-zero) for the text in Inkscape, but the package 'transparent.sty' is not loaded}%
    \renewcommand\transparent[1]{}%
  }%
  \providecommand\rotatebox[2]{#2}%
  \newcommand*\fsize{\dimexpr\f@size pt\relax}%
  \newcommand*\lineheight[1]{\fontsize{\fsize}{#1\fsize}\selectfont}%
  \ifx\svgwidth\undefined%
    \setlength{\unitlength}{216bp}%
    \ifx\svgscale\undefined%
      \relax%
    \else%
      \setlength{\unitlength}{\unitlength * \real{\svgscale}}%
    \fi%
  \else%
    \setlength{\unitlength}{\svgwidth}%
  \fi%
  \global\let\svgwidth\undefined%
  \global\let\svgscale\undefined%
  \makeatother%
  \begin{picture}(1,0.44819923)%
    \lineheight{1}%
    \setlength\tabcolsep{0pt}%
    \put(0,0){\includegraphics[width=\unitlength]{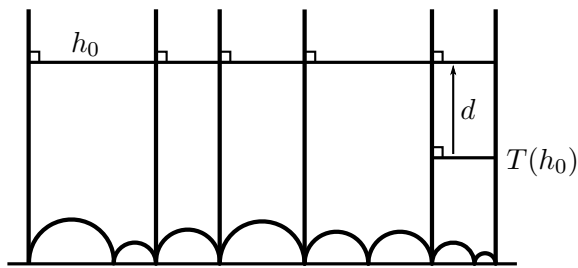}}%
    \put(0.79034388,0.25838449){\color[rgb]{0,0,0}\makebox(0,0)[lt]{\lineheight{0}\smash{\begin{tabular}[t]{l}$d$\end{tabular}}}}%
    \put(0.11111111,0.37412523){\color[rgb]{0,0,0}\makebox(0,0)[lt]{\lineheight{0}\smash{\begin{tabular}[t]{l}$h_0$\end{tabular}}}}%
    \put(0.87037035,0.17042153){\color[rgb]{0,0,0}\makebox(0,0)[lt]{\lineheight{0}\smash{\begin{tabular}[t]{l}$T(h_0)$\end{tabular}}}}%
  \end{picture}%
\endgroup%

\end{center}
\caption{Extending a horocycle.}
\label{Fig:Extend-h}
\end{figure}

Fix an ideal vertex on one of the polygons $P$.  Put $P$ in $\HH^2$ with $v$ at infinity.  Now take $h_0$ to be a horocycle centered at infinity intersected with $P$.  Follow $h_0$ to the right.  When it meets the edge of $P$, a new polygon is glued.  The developing map instructs us how to embed that new polygon as a polygon in $\HH^2$, with one edge the vertical geodesic which is the edge of $P$. Continue along this horocycle, placing polygons in $\HH^2$ according to their developing image.  Eventually the horocycle will meet $P$ again with $v$ at infinity. When this happens, the developing map will instruct us to glue a copy of $P$ to the given edge. This copy of $P$ will be isometric to the original copy of $P$, where the isometry is the holonomy\index{holonomy} of the closed path which encircles the ideal vertex $v$ once in the counterclockwise direction. This holonomy isometry, call it $T$, takes the horocycle $h_0$ on our original copy of $P$ to a horocycle $T(h_0)$, and $T(h_0)$ will be of distance $d(v)$\index{$d(v)$} from the extended horocycle that began with $h_0$.  See \reffig{Extend-h}.

\begin{proposition}\label{Prop:Complete}
Let $S$ be a surface with hyperbolic structure obtained by gluing
hyperbolic polygons.  Then the metric on $S$ is complete\index{complete metric space} if and only
if $d(v)=0$ for each ideal vertex $v$.
\end{proposition}

Before we prove this proposition, let's look at an example.

\begin{example}[Complete 3-punctured sphere]
\label{Example:3-punct-sphere}

A topological polygonal decomposition for the 3-punctured sphere consists of two ideal triangles.\index{ideal triangle} See \reffig{3-punct1}.

\begin{figure}
%% Accompanies image file 'F3-10-3PTri.eps' (pdf, eps, ps)
\begingroup%
  \makeatletter%
  \providecommand\color[2][]{%
    \errmessage{(Inkscape) Color is used for the text in Inkscape, but the package 'color.sty' is not loaded}%
    \renewcommand\color[2][]{}%
  }%
  \providecommand\transparent[1]{%
    \errmessage{(Inkscape) Transparency is used (non-zero) for the text in Inkscape, but the package 'transparent.sty' is not loaded}%
    \renewcommand\transparent[1]{}%
  }%
  \providecommand\rotatebox[2]{#2}%
  \newcommand*\fsize{\dimexpr\f@size pt\relax}%
  \newcommand*\lineheight[1]{\fontsize{\fsize}{#1\fsize}\selectfont}%
  \ifx\svgwidth\undefined%
    \setlength{\unitlength}{235.0524044bp}%
    \ifx\svgscale\undefined%
      \relax%
    \else%
      \setlength{\unitlength}{\unitlength * \real{\svgscale}}%
    \fi%
  \else%
    \setlength{\unitlength}{\svgwidth}%
  \fi%
  \global\let\svgwidth\undefined%
  \global\let\svgscale\undefined%
  \makeatother%
  \begin{picture}(1,0.32071754)%
    \lineheight{1}%
    \setlength\tabcolsep{0pt}%
    \put(0,0){\includegraphics[width=\unitlength]{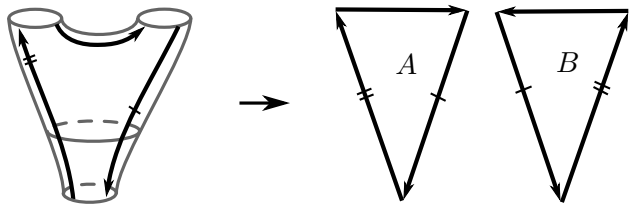}}%
    \put(0.62243477,0.21155023){\color[rgb]{0,0,0}\makebox(0,0)[lt]{\lineheight{0}\smash{\begin{tabular}[t]{l}$A$\end{tabular}}}}%
    \put(0.88028857,0.21498841){\color[rgb]{0,0,0}\makebox(0,0)[lt]{\lineheight{0}\smash{\begin{tabular}[t]{l}$B$\end{tabular}}}}%
  \end{picture}%
\endgroup%

\caption{Topological polygonal decomposition for the 3-punctured sphere.}
\label{Fig:3-punct1}
\end{figure}

Let's try to construct a geometric polygonal decomposition by building the developing image.  We can put one of the ideal triangles in $\HH^2$ as the triangle with vertices at $0$, $1$, $\infty$.  If we glue the other triangle immediately to the right, we have two vertices at $1$ and at $\infty$, but the third can go to any point $x$, where $x>1$.  See \reffig{3punct2}.  These two triangles on the left, labeled $A$ and $B$, give a fundamental region for the 3-punctured sphere.  The developing image will be created by gluing additional copies of these two triangles to edges in the figure by holonomy\index{holonomy} isometries.

\begin{figure}
%% Accompanies image file 'F3-11-3PDev.eps' (pdf, eps, ps)
\begingroup%
  \makeatletter%
  \providecommand\color[2][]{%
    \errmessage{(Inkscape) Color is used for the text in Inkscape, but the package 'color.sty' is not loaded}%
    \renewcommand\color[2][]{}%
  }%
  \providecommand\transparent[1]{%
    \errmessage{(Inkscape) Transparency is used (non-zero) for the text in Inkscape, but the package 'transparent.sty' is not loaded}%
    \renewcommand\transparent[1]{}%
  }%
  \providecommand\rotatebox[2]{#2}%
  \newcommand*\fsize{\dimexpr\f@size pt\relax}%
  \newcommand*\lineheight[1]{\fontsize{\fsize}{#1\fsize}\selectfont}%
  \ifx\svgwidth\undefined%
    \setlength{\unitlength}{126.09207916bp}%
    \ifx\svgscale\undefined%
      \relax%
    \else%
      \setlength{\unitlength}{\unitlength * \real{\svgscale}}%
    \fi%
  \else%
    \setlength{\unitlength}{\svgwidth}%
  \fi%
  \global\let\svgwidth\undefined%
  \global\let\svgscale\undefined%
  \makeatother%
  \begin{picture}(1,0.91361807)%
    \lineheight{1}%
    \setlength\tabcolsep{0pt}%
    \put(0,0){\includegraphics[width=\unitlength]{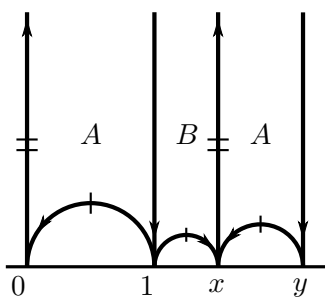}}%
    \put(0.40071341,0.04825189){\color[rgb]{0,0,0}\makebox(0,0)[lt]{\lineheight{0}\smash{\begin{tabular}[t]{l}$1$\end{tabular}}}}%
    \put(0.60521247,0.05901489){\color[rgb]{0,0,0}\makebox(0,0)[lt]{\lineheight{0}\smash{\begin{tabular}[t]{l}$x$\end{tabular}}}}%
    \put(0.85786228,0.06128103){\color[rgb]{0,0,0}\makebox(0,0)[lt]{\lineheight{0}\smash{\begin{tabular}[t]{l}$y$\end{tabular}}}}%
    \put(0.22205994,0.50058555){\color[rgb]{0,0,0}\makebox(0,0)[lt]{\lineheight{0}\smash{\begin{tabular}[t]{l}$A$\end{tabular}}}}%
    \put(0.50756558,0.50058555){\color[rgb]{0,0,0}\makebox(0,0)[lt]{\lineheight{0}\smash{\begin{tabular}[t]{l}$B$\end{tabular}}}}%
    \put(0.72962552,0.50058555){\color[rgb]{0,0,0}\makebox(0,0)[lt]{\lineheight{0}\smash{\begin{tabular}[t]{l}$A$\end{tabular}}}}%
    \put(0.01659108,0.05040108){\color[rgb]{0,0,0}\makebox(0,0)[lt]{\lineheight{0}\smash{\begin{tabular}[t]{l}$0$\end{tabular}}}}%
  \end{picture}%
\endgroup%

\caption{We may choose any $x>1$, $y>x$ when finding a hyperbolic
structure.}
\label{Fig:3punct2}
\end{figure}

We may choose the position of the next copy of the triangle $A$ glued to the right, putting its vertex at the point $y$ as in \reffig{3punct2}.  After this choice, notice we cannot choose where the next vertex of $B$ to the right will go.  This is because the choice $y$ determines an isometry of $\HH^2$ taking the triangle $A$ on the left to the triangle labeled $A$ on the right.  This isometry is exactly the holonomy\index{holonomy} element corresponding to the closed curve running once around the vertex at infinity. The same isometry, which has been determined with the choice of $y$, must take $B$ in the middle to the next triangle glued to the right in our figure. In fact, now that we know this holonomy\index{holonomy} element, we may apply it and its inverse successively to the triangles of \reffig{3punct2}, and we obtain the entire developing image of all triangles adjacent to infinity.

Recall that we want our hyperbolic structure to be complete.\index{complete metric space}  By \refprop{Complete}, we need to look at horocycles.  Pick a collection of horocycles about the vertices $0$, $1$, and $\infty$.  Each of these horocycles extends to give a new horocycle about another copy of $A$.  Each copy of $A$ is obtained by applying a holonomy isometry to the original triangle with vertices at $0$, $1$, and $\infty$.  We want the horocycles obtained under these holonomy
isometries to agree with the horocycles obtained by extending the original horocycles.  This is the condition for completeness.

Here is one way to determine complete structures.\index{complete metric space} Let $\ell_1$ denote the distance in $\HH^2$ between the horocycle at infinity and the horocycle at $0$. See \reffig{3punct3}. The holonomy\index{holonomy} element $\psi$, corresponding to the group element fixing the ideal vertex at infinity, is an isometry of $\HH^2$, hence it preserves distances. Thus, under this isometry, the distance between the image of the horocycle at infinity and the horocycle at $\psi(0)=x$ must also be $\ell_1$. If the structure is complete, then the horocycle about infinity is preserved by $\psi$. Thus the horocycle at $x$ must have the same (Euclidean) diameter as the horocycle at $0$.

Now consider the length of the edge between horocycles at $0$ and $1$, labeled $\ell_3$ in \reffig{3punct3}. There is another holonomy\index{holonomy} isometry $\phi$ mapping the geodesic edge between $0$ and $1$ to one between $x$ and $1$, corresponding to the group element encircling the ideal vertex at $1$.
Again completeness implies that the horocycle at $1$ is fixed by $\phi$. The horocycle at $0$ maps to the horocycle centered at $x$, and again because $\phi$ is an isometry, the distance between horocycles centered at $1$ and $x$ must still be $\ell_3$. We already determined the fact that the horocycle at $x$ has the same (Euclidean) diameter as the one at $0$. The only possible way that the distance $\ell_3$ will also be preserved is if $x=2$, and the picture is symmetric across the edge from $1$ to infinity. 

\begin{figure}
%% Accompanies image file 'F3-12-3PLen.eps' (pdf, eps, ps)
\begingroup%
  \makeatletter%
  \providecommand\color[2][]{%
    \errmessage{(Inkscape) Color is used for the text in Inkscape, but the package 'color.sty' is not loaded}%
    \renewcommand\color[2][]{}%
  }%
  \providecommand\transparent[1]{%
    \errmessage{(Inkscape) Transparency is used (non-zero) for the text in Inkscape, but the package 'transparent.sty' is not loaded}%
    \renewcommand\transparent[1]{}%
  }%
  \providecommand\rotatebox[2]{#2}%
  \newcommand*\fsize{\dimexpr\f@size pt\relax}%
  \newcommand*\lineheight[1]{\fontsize{\fsize}{#1\fsize}\selectfont}%
  \ifx\svgwidth\undefined%
    \setlength{\unitlength}{110.63493347bp}%
    \ifx\svgscale\undefined%
      \relax%
    \else%
      \setlength{\unitlength}{\unitlength * \real{\svgscale}}%
    \fi%
  \else%
    \setlength{\unitlength}{\svgwidth}%
  \fi%
  \global\let\svgwidth\undefined%
  \global\let\svgscale\undefined%
  \makeatother%
  \begin{picture}(1,0.97145619)%
    \lineheight{1}%
    \setlength\tabcolsep{0pt}%
    \put(0,0){\includegraphics[width=\unitlength]{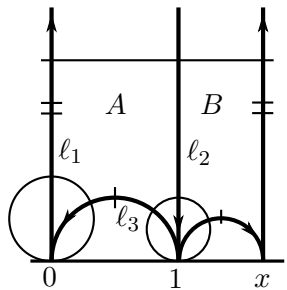}}%
    \put(0.54593968,0.01593257){\color[rgb]{0,0,0}\makebox(0,0)[lt]{\lineheight{0}\smash{\begin{tabular}[t]{l}$1$\end{tabular}}}}%
    \put(0.83776177,0.01980581){\color[rgb]{0,0,0}\makebox(0,0)[lt]{\lineheight{0}\smash{\begin{tabular}[t]{l}$x$\end{tabular}}}}%
    \put(0.32901,0.60990687){\color[rgb]{0,0,0}\makebox(0,0)[lt]{\lineheight{0}\smash{\begin{tabular}[t]{l}$A$\end{tabular}}}}%
    \put(0.65440451,0.60990687){\color[rgb]{0,0,0}\makebox(0,0)[lt]{\lineheight{0}\smash{\begin{tabular}[t]{l}$B$\end{tabular}}}}%
    \put(0.1172453,0.02324212){\color[rgb]{0,0,0}\makebox(0,0)[lt]{\lineheight{0}\smash{\begin{tabular}[t]{l}$0$\end{tabular}}}}%
    \put(0.36645621,0.22640626){\color[rgb]{0,0,0}\makebox(0,0)[lt]{\lineheight{0}\smash{\begin{tabular}[t]{l}$\ell_3$\end{tabular}}}}%
    \put(0.17147775,0.45624835){\color[rgb]{0,0,0}\makebox(0,0)[lt]{\lineheight{0}\smash{\begin{tabular}[t]{l}$\ell_1$\end{tabular}}}}%
    \put(0.60791957,0.45237467){\color[rgb]{0,0,0}\makebox(0,0)[lt]{\lineheight{0}\smash{\begin{tabular}[t]{l}$\ell_2$\end{tabular}}}}%
  \end{picture}%
\endgroup%

  \caption{Lengths between horocycles}
  \label{Fig:3punct3}
\end{figure}

Note at this point that the holonomy\index{holonomy} $\psi$ is completely determined: It fixes $\infty$, takes $0$ to $2$, and maps a point $i\,h$ on a horocycle about infinity to the point $2+i\,h$. This is the translation $\psi(z)=z+2$. Similarly, the holonomy $\phi$ is also completely determined, as it fixes $1$, maps $0$ to $2$ and takes a point on a horocycle on the edge of length $\ell_3$ to a determined point on the edge from $x$ to $1$. 
Because the fundamental group of the 3-punctured sphere is generated by the two loops corresponding to $\psi$ and $\phi$, this determines the complete structure. We have therefore shown:

\begin{proposition}\label{Prop:CompleteStruct3punct}
There is a unique complete hyperbolic structure\index{complete metric space} on the 3-punctured sphere. A fundamental region for the structure is given by two ideal triangles\index{ideal triangle} with vertices $0$, $1$, and $\infty$ and $1$, $2$, and $\infty$, respectively.\qed
\end{proposition}

\end{example}

\begin{example}[Incomplete structure on 3-punctured sphere]
\label{Example:Incomplete3PunctSphere}\index{incomplete 3-punctured sphere}
  What if we choose a different value for $x$ besides $x=2$? Say we let $x=3/2$. To simplify things, let's keep the length of the edge between horocycles at $0$ and $1$ constant as we extend horocycles. Choose horocycles at $0$ and $1$ of (Euclidean) radius $1/2$, so that these horocycles are tangent along the edge between $0$ and $1$, hence the distance between horocycles is $0$. This distance will remain equal to $0$ under each holonomy\index{holonomy} element, so there will be a horocycle at $x=3/2$ tangent to the horocycle about $1$, to preserve distance $0$. This determines where the image of the triangle $A$ must go under the holonomy fixing infinity: its third vertex (called $y$ in \reffig{3punct2}) must have a horocycle about it of the same (Euclidean) size as the horocycle at $3/2$. This determines the holonomy isometry about the vertex at infinity. Apply this holonomy isometry successively, and we obtain a pattern of triangles as in \reffig{3punct-incomplete}.

\begin{figure}
\begin{center}
%% Accompanies image file 'F3-13-3PInc.eps' (pdf, eps, ps)
\begingroup%
  \makeatletter%
  \providecommand\color[2][]{%
    \errmessage{(Inkscape) Color is used for the text in Inkscape, but the package 'color.sty' is not loaded}%
    \renewcommand\color[2][]{}%
  }%
  \providecommand\transparent[1]{%
    \errmessage{(Inkscape) Transparency is used (non-zero) for the text in Inkscape, but the package 'transparent.sty' is not loaded}%
    \renewcommand\transparent[1]{}%
  }%
  \providecommand\rotatebox[2]{#2}%
  \newcommand*\fsize{\dimexpr\f@size pt\relax}%
  \newcommand*\lineheight[1]{\fontsize{\fsize}{#1\fsize}\selectfont}%
  \ifx\svgwidth\undefined%
    \setlength{\unitlength}{194.40000343bp}%
    \ifx\svgscale\undefined%
      \relax%
    \else%
      \setlength{\unitlength}{\unitlength * \real{\svgscale}}%
    \fi%
  \else%
    \setlength{\unitlength}{\svgwidth}%
  \fi%
  \global\let\svgwidth\undefined%
  \global\let\svgscale\undefined%
  \makeatother%
  \begin{picture}(1,0.56907626)%
    \lineheight{1}%
    \setlength\tabcolsep{0pt}%
    \put(0,0){\includegraphics[width=\unitlength]{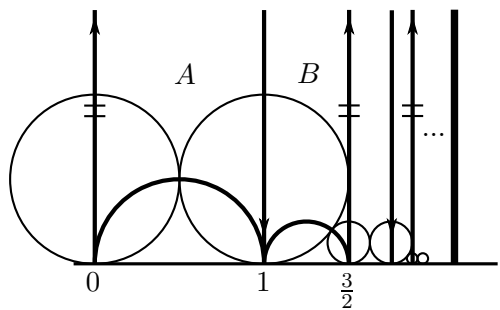}}%
    \put(0.47898299,0.02307376){\color[rgb]{0,0,0}\makebox(0,0)[lt]{\lineheight{0}\smash{\begin{tabular}[t]{l}$1$\end{tabular}}}}%
    \put(0.31951796,0.42504353){\color[rgb]{0,0,0}\makebox(0,0)[lt]{\lineheight{0}\smash{\begin{tabular}[t]{l}$A$\end{tabular}}}}%
    \put(0.55834807,0.42798301){\color[rgb]{0,0,0}\makebox(0,0)[lt]{\lineheight{0}\smash{\begin{tabular}[t]{l}$B$\end{tabular}}}}%
    \put(0.14902998,0.02282459){\color[rgb]{0,0,0}\makebox(0,0)[lt]{\lineheight{0}\smash{\begin{tabular}[t]{l}$0$\end{tabular}}}}%
    \put(0.63844805,0.00984596){\color[rgb]{0,0,0}\makebox(0,0)[lt]{\lineheight{0}\smash{\begin{tabular}[t]{l}$\frac{3}{2}$\end{tabular}}}}%
    \put(0.80011762,0.32363248){\color[rgb]{0,0,0}\makebox(0,0)[lt]{\lineheight{0}\smash{\begin{tabular}[t]{l}...\end{tabular}}}}%
  \end{picture}%
\endgroup%

\end{center}
\caption{Part of developing image of an incomplete structure on a
	3-punctured sphere.}
\label{Fig:3punct-incomplete}
\end{figure}

Notice that the edges of the triangles approach a limit --- the thick line shown on the far right of the figure.  Notice also that this line is not part of the developing image of the 3-punctured sphere.

This hyperbolic structure is incomplete: for any horocycle about infinity in $\HH^2$, the sequence of points at the intersection of the horocycle and the edges of the developing images of ideal triangles\index{ideal triangle} projects to a Cauchy sequence that does not converge. Alternately, the value $d(v)$\index{$d(v)$} is nonzero for $v$ the ideal vertex lifting to the point at infinity.

An incomplete metric space\index{complete metric space} may be completed by adjoining points corresponding to limits of Cauchy sequences, and giving the resulting space the metric topology. In our case, the completion of this incomplete 3-punctured sphere is obtained by attaching a geodesic segment --- the projection of the thick line in \reffig{3punct-incomplete}. Each point of the thick geodesic on the right of \reffig{3punct-incomplete} corresponds to the limiting point of the Cauchy sequence given by a horocycle about infinity at the appropriate height. Note that in the quotient, however, we attach a closed curve of length $d(v)$,\index{$d(v)$} since points on that thick geodesic lying on horospheres of distance $d(v)$ apart will be identified.

Note that horocycles about infinity run straight into this thick geodesic, meeting it at right angles. On the other hand, these horocycles meet infinitely many edges of ideal triangles\index{ideal triangle} on their way into the geodesic, and none of these ideal edges meets the geodesic. It follows that the ideal edges of the two triangles $A$ and $B$ become arbitrarily close to the geodesic attached in the completion, without ever meeting it. Geometrically, it appears that the edges of the ideal triangles spin around the geodesic infinitely many times, while horocycles run directly into it. See \reffig{3punctCompletion}.

\begin{figure}
\begin{center}
\includegraphics{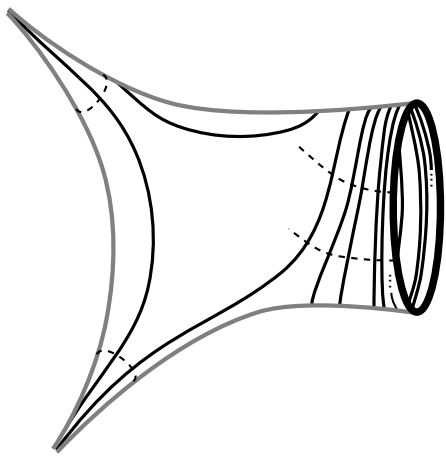}
\end{center}
\caption{The completion of an incomplete structure on a 3-punctured sphere.\index{complete metric space} Attach a geodesic of length $d(v)$.\index{$d(v)$} Ideal edges spin arbitrarily close to the attached geodesic without meeting it. Horocycles (dashed) run directly into the geodesic. (This example has two complete cusped ends and one incomplete end.)}
\label{Fig:3punctCompletion}
\end{figure}
\end{example}

\begin{proof}[Proof of \refprop{Complete}]
Let $S$ be a surface obtained by gluing hyperbolic polygons. 

Suppose first that $d(v)$\index{$d(v)$} is nonzero. Then take a sequence of points on a horocycle about $v$, one point for each intersection of the horocycle with an ideal edge. This gives a Cauchy sequence that does not converge. Therefore, the metric is not complete.\index{complete metric space}

Now suppose $d(v)=0$ for each ideal vertex $v$.  Then some horocycle closes up around each ideal vertex, so we may remove the interior horoball from each polygon. After this removal, the remainder is a compact manifold with boundary.  For any $t>0$, let $S_t$ be the compact manifold obtained by removing interiors of horocycles of distance $t$ from our original choice of horocycle. Then the compact subsets $S_t$ of $S$ satisfy $\bigcup_{t\in\RR^+} S_t= S$ and $S_{t+a}$ contains a neighborhood of radius $a$ about $S_t$. Any Cauchy sequence must be contained in some $S_t$ for sufficiently large $t$.  Hence by compactness of $S_t$, the Cauchy sequence must converge.
\end{proof}

\section{Developing map and completeness}

Here is a better condition for completeness that works in all dimensions and all geometries.  

\begin{theorem}\label{Thm:Developing}
Let $M$ be an $n$-manifold with a $(G,X)$-structure,\index{$(G,X)$-structure} where $G$ acts transitively on $X$, and $X$ admits a complete $G$-invariant metric. Then the metric on $M$ inherited from $X$ is complete\index{complete metric space} if and only if the developing map $D\from \widetilde{M} \to X$ is a covering map.
\end{theorem}

\begin{proof}
Suppose first that the developing map $D\from \widetilde{M}\to X$ is a covering map. Let $\{x_n\}_{n=1}^\infty$ be a Cauchy sequence in $M$. For $n$ large enough, $x_n$ will be contained in an $\epsilon$-ball in $M$ that is evenly covered in $\widetilde{M}$. Thus the sequence lifts to a Cauchy sequence $\{\widetilde{x}_n\}$ in $\widetilde{M}$. Since $D$ is a local isometry, $\{D(\widetilde{x}_n)\}$ is a Cauchy sequence in $X$. Finally since $X$ is complete,\index{complete metric space} $\{D(\widetilde{x}_n)\}$ converges to $y\in X$. Now, because $D$ is a covering map, there is a neighborhood $U$ of $y$ that is evenly covered by $D$. Lift this to a neighborhood $\widetilde{U}$ of $\widetilde{M}$ containing infinitely many points of the sequence $\{\widetilde{x}_n\}$. The lift of $y$ in this neighborhood, call it $\widetilde{y}$, must be a limit point of $\{\widetilde{x}_n\}$. Then the projection of $\widetilde{y}$ to $M$ is a limit point of the sequence $\{x_n\}$, so $M$ is complete.\index{complete metric space}

For the converse, we appeal to a proof by Thurston~\cite[Proposition 3.4.15]{thurston:book}. Suppose $M$ is complete.\index{complete metric space} To show $D\from \widetilde{M} \to X$ is a covering map, we show that any path $\alpha_t$ in $X$ lifts to a path $\widetilde{\alpha_t}$ in $\widetilde{M}$. Since $D$ is a local homeomorphism, this implies that $D$ is a covering map.

First, if $M$ is complete,\index{complete metric space} then $\widetilde{M}$ must also be complete, where the metric $\widetilde{M}$ is the lift of the metric on $M$, as follows. The projection to $M$ of any Cauchy sequence gives a Cauchy sequence in $M$, with limit point $x$. Then $x$ has a compact neighborhood which is evenly covered in $\widetilde{M}$, hence there is a compact neighborhood in $\widetilde{M}$ containing all but finitely many points of the Cauchy sequence and also containing a lift of $x$. Thus the sequence converges in $\widetilde{M}$.

Let $\alpha_t$ be a path in $X$. Because $D$ is a local homeomorphism, we may lift $\alpha_t$ to a path $\widetilde{\alpha_t}$ in $\widetilde{M}$ for $t \in [0, t_0)$, some $t_0>0$. By completeness of $\widetilde{M}$, the lifting extends to $[0, t_0]$. But because $D$ is a local homeomorphism, a lifting to $[0, t_0]$ extends to $[0, t_0+\epsilon)$. Hence the lifting extends to all of $\alpha_t$ and $D$ is a covering map.
\end{proof}

\begin{corollary}\label{Cor:DevelopingCover}
If $X$ is simply connected, and $M$ is a manifold with a $(G,X)$-structure\index{$(G,X)$-structure} as in \refthm{Developing}, then $M$ is complete\index{complete metric space} if and only if the developing map is an isometry of $X$.
\end{corollary}

\begin{proof}
The developing map is a local isometry by construction. \Refthm{Developing} shows that $M$ is complete if and only if the developing map is a covering map. Since $X$ and $\widetilde{M}$ are simply connected, the developing map is a covering map if and only if it is a covering isomorphism. A covering isomorphism that is a local isometry must be an isometry. 
\end{proof}

%%%%%%%%%%%%%%%%%%%%%%%%%%%%%%%%%%%%%%%%%%%%%%%%%%%%%%%%%%%%%%%%%
\section{Exercises}

\begin{exercise} We have seen that Euclidean structures\index{Euclidean structure} on a torus are determined by a parallelogram.
\begin{enumerate}
\item[(a)] Show that by applying translation and rotation isometries of $\EE^2$, we may assume that the parallelogram has vertices $(0,0)$, $(x_1,0)$, $(x_2, y)$, and $(x_1+x_2,y)$ where $x_1>0$ and $y>0$.
\item[(b)] Show that up to rescaling, a parallelogram has vertices $(0,0)$, $(1,0)$, $(x,y)$, and $(x+1,y)$ for some $(x,y)\in \RR^2$ with $y>0$.
\end{enumerate}
\end{exercise}

\begin{exercise}\label{Ex:InheritMetric}
If $X$ is a metric space, and $G$ is a group of isometries acting transitively on $X$, and $M$ is a manifold admitting a $(G,X)$-structure,\index{$(G,X)$-structure} show that $M$ inherits a metric from $X$. That is, explain how to define a metric on $M$ from that on $X$, and show that the metric is well-defined. 
\end{exercise}

\begin{exercise} (Induced structures \cite[Exercise 3.1.5]{thurston:book}).
Let $N$ be a topological space and $M$ a manifold with a $(G,X)$-structure,\index{$(G,X)$-structure} and suppose $\pi \from N \to M$ is a local homeomorphism.  Prove $N$ has a $(G,X)$-structure that is preserved by $\pi$. As a corollary, show that any covering space of $M$ admits a $(G,X)$-structure. 
\end{exercise}

\begin{exercise}[Analytic continuation\index{analytic continuation}]\label{Ex:AnalyticContinuation}
  Prove item \eqref{Itm:DevMapWellDefined} of \refprop{DevelopingMapProperties}. That is, prove the following.
  \begin{enumerate}
  \item[(a)] Suppose $\alpha\from [0,1]\to M$ is a path. Let $(U_i,\phi_i)$ and $(V_j,\psi_j)$ be two choices of charts\index{chart} that cover $\alpha([0,1])$. Let
    \[ 0=t_0<t_1<\dots<t_n=1 \mbox{ and }
    0=s_0<s_1<\dots<s_m=1\]
    be points in $[0,1]$ such that $\alpha([t_i,t_{i+1}])\subset U_i$ and $\alpha([s_j,s_{j+1}])\subset V_j$. Define inductively extensions $\Phi_i(t)$ and $\Psi_j(t)$ as in the definition of the developing map. Finally, let $X\subset [0,1]$ be the set of points on which $\Phi_i(t) = \Psi_j(t)$. Prove that $X=[0,1]$.

    (Hint: You will need to use the fact that $G$ is analytic. One reference for analytic continuation\index{analytic continuation} is \cite[Chapter~IX]{Conway:Complex}.)

  \item[(b)] Suppose $\alpha$ and $\beta$ are homotopic paths with the same endpoints. By part (a), there are well-defined functions $D(\alpha)$ and $D(\beta)$, defined separately on $\alpha$ and on $\beta$ as in \refdef{DevelopingMap}. Prove that $D(\alpha)=D(\beta)$. 
    This proves that the definition of $D$ is independent of choice $\alpha$ in the homotopy class of $[\alpha]\in\widetilde{M}$.

    (Again you will use the fact that $G$ is analytic; see for example \cite[Chapter~IX]{Conway:Complex}.)
  \end{enumerate}
\end{exercise}

\begin{exercise} Prove item \eqref{Itm:DevMapBaseptChange} of \refprop{DevelopingMapProperties}. Show that if we define a new map in the same way as $D$, except we change the basepoint $x_0$ or the initial chart\index{chart} $(U_0,\phi_0)$, then the resulting map is equal to the composition of $D$ with an element of $G$.
\end{exercise}

\begin{exercise}\label{Ex:AffineTorusExplicit}
Let $T$ be the affine torus\index{affine torus} obtained by identifying the sides of the trapezoid with vertices $(0,0)$, $(1,0)$, $(0,1)$, and $(1/2, 1)$.
\begin{enumerate}
\item[(a)] Compute the holonomy\index{holonomy} elements of $T$ corresponding to meridian and longitude (i.e.\ the loop running along the horizontal edge of the trapezoid and the loop running along the vertical edge of the trapezoid). What is the holonomy group\index{holonomy group} of $T$?
\item[(b)] For basepoint $(0,0)$ and initial chart\index{chart} chosen so that the trapezoid is mapped by the identity into $\RR^2$, compute explicitly the developing images of various curves, including the following:
\begin{itemize}
\item The curve running twice along the meridian (based at $(0,0)$).
\item The curve running twice along the longitude.
\item The curve running twice along the meridian and three times along the longitude.
\end{itemize}
\end{enumerate}
\end{exercise}

\begin{exercise}\label{Ex:DevelopingImageMissesPointExplicit}
Let $T$ be the affine torus\index{affine torus} of \refex{AffineTorusExplicit}, obtained by identifying sides of the trapezoid with vertices $(0,0)$, $(1,0)$, $(0,1)$, and $(1/2,1)$, and let $\widetilde{T}$ denote its universal cover. Prove that the developing image $D(\widetilde{T})\subset \RR^2$ misses exactly one point. 
\end{exercise}

\begin{exercise}\label{Ex:DevelopingImageMissesPoint}
Generalize \refex{DevelopingImageMissesPointExplicit}: Let $T$ be any affine torus.\index{affine torus} Prove that either the developing map $D\from \widetilde{T}\to \RR^2$ is a covering map, and $T$ is a Euclidean torus, or the image of the developing map misses a single point in $\RR^2$. 
\end{exercise}

\begin{exercise}
Fix an example of your favorite quadrilateral that is not a parallelogram, and let $T$ be the torus obtained by identifying sides. Use a computer to create a picture such as \reffig{AffineTorus} for your quadrilateral. 
\end{exercise}

\begin{exercise} Prove \reflem{dvInd-h0}, that $d(v)$\index{$d(v)$} is independent of initial choice of horocycle, and independent of choice of $v_0$ in the equivalence class of $v$. 
\end{exercise}

\begin{exercise}
  Prove the holonomy\index{holonomy group} group of the complete structure\index{complete metric space} on a 3-punctured sphere is generated by
  \[ \mat{1&2\\0&1} \quad \mbox{and} \quad \mat{1&0\\2&1}. \]
\end{exercise}

\begin{exercise} How many incomplete hyperbolic structures are there on a 3-punctured sphere?  How can they be parameterized?  Give a geometric interpretation of this parameterization. That is, relate the parameterization to the developing image of the associated hyperbolic structure.
\end{exercise}

\begin{exercise}\label{Ex:1punctTorus} A torus with 1 puncture has a topological polygonal decomposition consisting of two triangles.
\begin{enumerate}
\item[(a)] Find a complete hyperbolic structure on the 1-punctured torus and prove your structure is complete.
\item[(b)] Find all complete hyperbolic structures on the 1-punctured torus.  How are they parameterized?
\end{enumerate}
\end{exercise}

\begin{exercise}\label{Ex:4PunctSphere} A sphere with 4 punctures has a topological polygonal decomposition consisting of four triangles. Repeat \refex{1punctTorus} for the 4-punctured sphere. 
\end{exercise}

%% Ch04_Gluing.tex 

\chapter{Hyperbolic Structures and Triangulations}\label{Chap:GluingCompleteness}
\blfootnote{Jessica S. Purcell, Hyperbolic Knot Theory}

In \refchap{Geometric}, we learned that hyperbolic structures lead to developing maps and holonomy,\index{holonomy} and that the developing map is a covering map if and only if the hyperbolic structure is complete.\index{complete metric space}

In this chapter, we wish to compute explicit complete hyperbolic structures on 3-manifolds, again with our primary examples being knot complements. One of the most straightforward ways to find a hyperbolic structure is to first triangulate the manifold, or subdivide it into tetrahedra, and then to put a hyperbolic structure on each tetrahedron, ensuring the tetrahedra glue to give a $(\PSL(2,\CC), \HH^3)$-structure whose developing map is complete. This method of computing hyperbolic structures has been studied by many, and in particular was implemented on the computer by J.~Weeks as part of his 1985 PhD thesis \cite{Weeks:Thesis}. Here we will describe the conditions required to obtain a complete hyperbolic structure via triangulations, and as usual, work through examples.

\section{Geometric triangulations}

In \refchap{Geometric}, we defined topological and geometric polygonal decompositions of 2-manifolds. We can extend these notions to 3-manifolds by considering decompositions into ideal polyhedra. In \refchap{Fig8Decomp}, we obtained topological ideal polyhedral decompositions for knot complements. For many applications, including those later in this chapter, it simplifies matters greatly to consider decompositions into ideal tetrahedra. 

\begin{definition}\label{Def:TopIdealTriang}
Let $M$ be a 3-manifold. A \emph{topological ideal triangulation}\index{topological ideal triangulation} of $M$ is a combinatorial way of gluing truncated tetrahedra (ideal tetrahedra) so that the result is homeomorphic to $M$. Truncated parts will correspond to the boundary of $M$. As before, a gluing should take faces to faces, edges to edges, etc. 
\end{definition}

\begin{example}\label{Example:Fig8TopTriang}
The figure-8 knot has a topological ideal triangulation consisting of two ideal tetrahedra, as we saw in \refex{CollapseBigons} in \refchap{Fig8Decomp}.
\end{example}

For a given knot complement, it is relatively easy to find topological ideal triangulations. For example, starting with any polyhedral decomposition, choose an ideal vertex $v$ and cone to that vertex: i.e.\ add edges between $v$ and all other ideal vertices, between any two edges meeting $v$ add an ideal triangle\index{ideal triangle} (adding an additional edge opposite $v$ if necessary), and between three triangles meeting $v$ add an ideal tetrahedron. Split off the resulting tetrahedra. This reduces the collection of polyhedra to a collection with at least one fewer ideal vertex. Hence after repeating a finite number of times, we are left with a collection of topological tetrahedra. 

\subsection{An extended example: the $6_1$ knot}\label{Example:61TopTriang}

We work out an example for the $6_1$ knot carefully.  We will see how to decompose the complement into five tetrahedra.  (In fact, the complement of the $6_1$ knot can be decomposed into four tetrahedra, but we won't bother simplifying further here.)

\begin{figure}[h]
  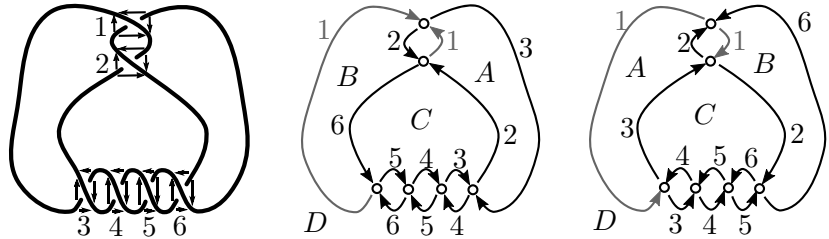
\caption{Left to right: The $6_1$ knot, the top polyhedron, the bottom polyhedron}
  \label{Fig:61Poly}
\end{figure}

We start with a polyhedral decomposition of the $6_1$ knot.  We use the decomposition obtained using the methods of \refchap{Fig8Decomp}. The result is shown in \reffig{61Poly}, with the knot on the left, the top polyhedron in the center, and the bottom polyhedron on the right. Recall all polyhedra are viewed from the outside; that is the ball of the polyhedron is behind the projection plane in each figure. In this example, oriented edges are labeled $1$ through $6$. 

\begin{figure}[h]
  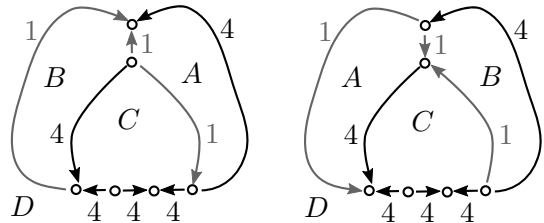
  \caption{Polyhedra for $6_1$ knot with bigons\index{bigon} collapsed}
  \label{Fig:61NoBigons}
\end{figure}

Collapse all bigons,\index{bigon} identifying edges $1$ and $2$, and $3$ through $6$. New edges and orientations are shown in \reffig{61NoBigons}.

\begin{figure}[h]
  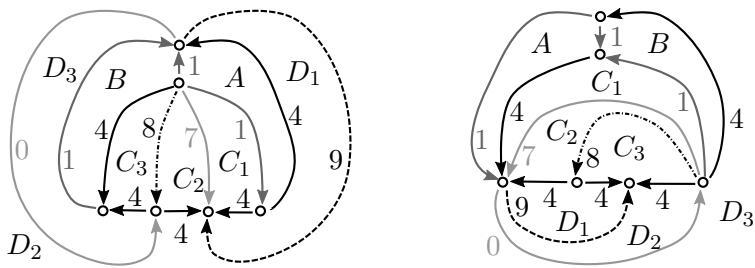
  \caption{A subdivision of faces $C$ and $D$ in the top
    polyhedron (left) leads to a subdivision of the bottom (right)}
  \label{Fig:61Subdivide}
\end{figure}

We cone the top polyhedron to the vertex in the center. This subdivides faces $C$ and $D$ into triangles, shown in \reffig{61Subdivide} in both top and bottom polyhedra. 

\begin{figure}[h]
  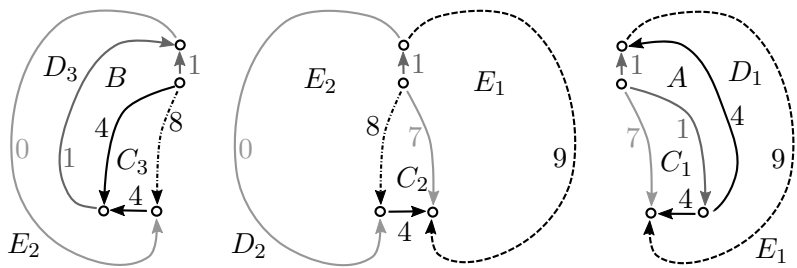
  \caption{The top polyhedron splits into the three tetrahedra shown}
  \label{Fig:61TopTetr}
\end{figure}

Continuing the subdivision in the top polyhedron, two edges meeting in the center vertex bound an ideal triangle;\index{ideal triangle} three triangles bound a tetrahedron. Thus edges labeled $1$, $7$, $9$ bound an ideal triangle $E_1$; edges labeled $1$, $8$, $0$ bound an ideal triangle $E_2$. Triangles $A$, $C_1$, $D_1$, and $E_1$ bound an ideal tetrahedron, as do triangles $B$, $C_3$, $D_3$, and $E_2$. When we split off these tetrahedra a single tetrahedron remains. All tetrahedra making up the top polyhedron are shown in \reffig{61TopTetr}. 

Now we split the bottom polyhedron into tetrahedra. However, first, observe in \reffig{61Subdivide} that edges labeled $7$ and $0$ in the bottom polyhedron run between the same two ideal vertices. Thus these two edges should be flattened and identified in the bottom polyhedron. While we could do that now in one step, we believe it is more geometrically clear how to flatten and identify if we first cut off ideal tetrahedra from the bottom polyhedron.

\begin{figure}[h]
  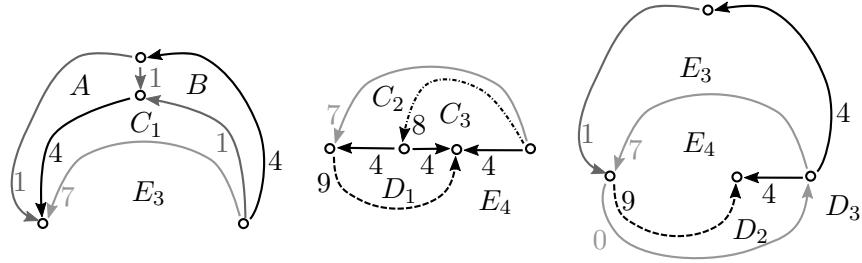
  \caption{Splitting off two tetrahedra in the bottom polyhedron}
  \label{Fig:61BottomSubdivide}
\end{figure}

So first, note there will be an ideal triangle\index{ideal triangle} $E_3$ with edges labeled $4$, $7$, and $1$, and this cuts off an ideal tetrahedron with sides $A$, $B$, $C_1$, $E_3$. Similarly there is an ideal triangle $E_4$ with edges $7$, $9$, and $4$, cutting off an ideal tetrahedron with sides $C_2$, $C_3$, $E_4$, and $D_1$. These two tetrahedra, as well as the remnant of the bottom polyhedron, are shown in \reffig{61BottomSubdivide}. 

Notice that the object on the right of \reffig{61BottomSubdivide} is not a tetrahedron: edges labeled $7$ and $0$ in that polyhedron form a bigon,\index{bigon} which collapses to a single edge which we label $7$.  When we do the collapse, the faces $E_4$ and $D_2$ collapse to a single triangle, which we will label $D_2$.  The faces $E_3$ and $D_3$ also collapse to a single triangle, which we
will label $D_3$. 

\begin{figure}
  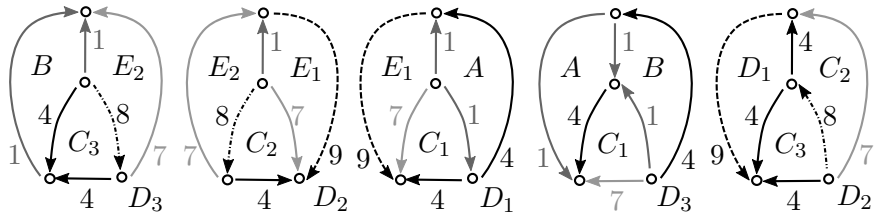
  \caption{Five tetrahedra which glue to give the complement of the
    $6_1$ knot}
  \label{Fig:61Tetr}
\end{figure}

When we have finished, we have five tetrahedra that glue to give the complement of the $6_1$ knot.  All five tetrahedra with their edges and faces labeled are shown in \reffig{61Tetr}. 

\subsection{Geometric ideal triangulations}

\begin{definition}\label{Def:GeomTriang}
A \emph{geometric ideal triangulation}\index{geometric ideal triangulation} of $M$ is a topological ideal triangulation such that each tetrahedron has a (positively oriented)\index{positively oriented tetrahedron}\index{tetrahedron!positively oriented} hyperbolic structure, and the result of gluing is a smooth manifold with a complete metric.\index{complete metric space} We also call such a triangulation a \emph{geometric triangulation}\index{geometric triangulation} for short. 
\end{definition}

As of the writing of this book, it is still an open question as to whether every 3-manifold that admits a complete hyperbolic structure actually admits a geometric ideal triangulation. 
It is known that every cusped hyperbolic 3-manifold can be decomposed into convex ideal polyhedra \cite{EpsteinPenner}: we will go through this in \refchap{Canonical}.
However, subdividing this decomposition into tetrahedra may create degenerate tetrahedra --- actual topological tetrahedra (as opposed to the object on the right of \reffig{61BottomSubdivide}), but tetrahedra that are flat in the hyperbolic structure on $M$.
There are known examples of generalized spaces with singularities that do not admit geometric triangulations\index{geometric triangulation} \cite{Choi:Triangulations}.

\section{Edge gluing equations}

In \refchap{Geometric}, we saw that a gluing of hyperbolic polygons has a hyperbolic structure if and only if the angle sum around each finite vertex is $2\pi$ (\reflem{VertexSum}).
There are similar conditions for a gluing of hyperbolic tetrahedra. We now need to consider gluing around an edge. 

Let $T$ be an ideal tetrahedron embedded in $\HH^3$. Any ideal tetrahedron has six edges. If we select any one, say $e$, we may choose an isometry of $\HH^3$ taking the endpoints of $e$ to $0$ and $\infty$, and sending a third vertex to $1 \in \CC\subset \bdy_\infty\HH^3$. This choice uniquely determines the isometry. The fourth vertex of $T$ will be mapped to some $z'\in\CC$. We may assume that $z'$ has positive imaginary part, for if not, apply an isometry of $\HH^3$ rotating around the geodesic from $0$ to $\infty$ and rescaling so that $z'$ maps to $1$. In this case, the image of $1$ under this isometry will be a complex number with positive imaginary part.

\begin{definition}\label{Def:EdgeInvariant}
For an ideal tetrahedron $T$ embedded in $\HH^3$, and edge $e$ of that tetrahedron, define the number $z(e)$ in $\CC$ to be the complex number with positive imaginary part obtained by applying the unique isometry of $\HH^3$ that takes the vertices of $e$ to $0$ and $\infty$, takes another vertex to $1$, and takes the final vertex of $T$ to $z(e)$. This is called the \emph{edge invariant}\index{edge invariant} of $e$. 
\end{definition}

\begin{remark}\label{Rem:Degenerate}
  Note that it is possible to map an ideal tetrahedron to $\HH^3$ so that three vertices map to $0$, $\infty$, and $1$, and the fourth maps to a point on the real line. In this case, the tetrahedron produced does not have a hyperbolic structure. If the fourth vertex is not $0$ or $1$, it is said to be \emph{flat}.\index{flat tetrahedron}\index{tetrahedron!flat} If the fourth vertex is $0$ or $1$, it is \emph{degenerate}.\index{degenerate tetrahedron}\index{tetrahedron!degenerate} Similarly, a fourth vertex mapped to infinity is a degenerate tetrahedron. An ideal triangulation of a hyperbolic 3-manifold with flat or degenerate tetrahedra is not a geometric ideal triangulation. When looking for geometric triangulations,\index{geometric triangulation}  we must rule out such tetrahedra. Similarly, for geometric triangulations, all edge invariants of all tetrahedra must have positive imaginary part. This ensures the tetrahedra are \emph{positively oriented}.\index{positively oriented tetrahedron}\index{tetrahedron!positively oriented} Finally, the procedure above always chooses an edge invariant with positive imaginary part. However, when we glue many tetrahedra together, at times it is impossible to simultaneously choose all edge invariants to have positive imaginary part; some may have negative imaginary part. Such a tetrahedron is a \emph{negatively oriented} tetrahedron.\index{negatively oriented tetrahedron}\index{tetrahedron!negatively oriented}
\end{remark}

Edge invariants of an ideal tetrahedron determine each other, in the following way.

\begin{lemma}\label{Lem:EdgeInvariants}
Let $T$ be an ideal tetrahedron with edge $e_1$, mapped so that vertices of $T$ lie at $\infty$, $0$, $1$, and $z(e_1)$ (so endpoints of $e_1$ lie at $0$ and $\infty$). Then $T$ has the following additional edge invariants.
\begin{itemize}
\item The edge $e_1'$ opposite $e_1$, with vertices $1$ and $z(e_1)$, has edge invariant $z(e_1')=z(e_1)$.
\item The edge $e_2$ with vertices $\infty$ and $1$ has edge invariant
  \[ z(e_2) = \frac{1}{1-z(e_1)}. \]
\item The edge $e_3$ with vertices $\infty$ and $z(e_1)$ has edge invariant
  \[ z(e_3) = \frac{z(e_1)-1}{z(e_1)}. \]
\end{itemize}
Thus we have the following relationships for these edge invariants.
\[ z(e_1)z(e_2)z(e_3) = -1, \quad \mbox{ and } \quad
1 - z(e_1) + z(e_1)z(e_3) = 0 \]
\end{lemma}

\begin{proof}
The proof is obtained by considering isometries of $\HH^3$ that move the different edges of $T$ onto the geodesic from $0$ to $\infty$. For ease of notation, we set $z=z(e_1)$.

For the first part, we label one more edge. Let $e_3'$ be the edge of $T$ opposite $e_3$. So $e_3'$ has endpoints $0$ and $1$. Note there is a geodesic $\gamma$ in $\HH^3$ that meets the edges $e_3$ and $e_3'$ orthogonally. An elliptic\index{elliptic} isometry rotating about $\gamma$ by angle $\pi$ maps $0$ to $1$ and $1$ to $0$, and maps $\infty$ to $z$ and $z$ to $\infty$, thus it preserves $T$. It takes the edge $e_1'$ with endpoints $1$ and $z$ to an edge with endpoints $0$ and $\infty$. Hence $z(e_1') = z$. 

To determine $z(e_2)$, we apply a M\"obius transformation\index{M\"obius transformation} fixing $\infty$, taking $1$ to $0$, and taking $z$ to $1$.  This transformation is given by
\[ w \mapsto \frac{w-1}{z-1}.\]
It sends $0$ to $-1/(z-1)$.  Thus $z(e_2) = 1/(1-z)$.

As for the edge $e_3$ running from $z$ to $\infty$, to determine its edge invariant we apply a M\"obius transformation\index{M\"obius transformation} fixing $\infty$, sending $z$ to $0$, and sending $0$ to $1$. This is given by
\[ w \mapsto \frac{w-z}{-z}.\]
It sends $1$ to $(1-z)/(-z)$.  Thus $z(e_3) = (z-1)/z$.
\end{proof}

\begin{figure}
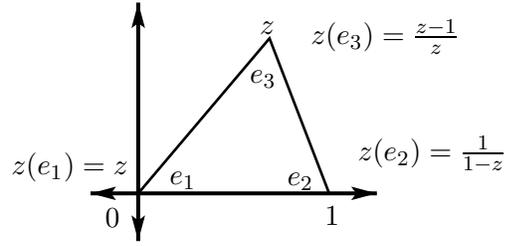
\caption{Edge invariants}
\label{Fig:EdgeInvt}
\end{figure}

The three edge invariants of a tetrahedron are shown in \reffig{EdgeInvt}.

Now consider a gluing of ideal tetrahedra. Fix an edge $e$ of the gluing, and let $T_1$ be a tetrahedron which has edge $e_1$ glued to $e$.  Put $T_1$ in $\HH^3$ with the edge $e_1$ running from $0$ to $\infty$, with a third vertex at $1$, and the fourth vertex at $z(e_1)$, where $z(e_1)$ has positive imaginary part. The gluing identifies each face of $T_1$ with another face.  Let $F_1$ denote the face of $T_1$ with vertices $0$, $z(e_1)$, and $\infty$.  This is glued to a face $F_1'$ in some tetrahedron $T_2$, where the edge $e_2$
in $T_2$ glues to $e$.

Now, we could put $T_2$ in $\HH^3$ with vertices at $0$, $\infty$, $1$, and $z(e_2)$, but since we're gluing to $T_1$, we want the face $F_1'$ to have vertices $0$, $\infty$, and $z(e_1)$ rather than vertices $0$, $\infty$, and $1$. Thus to do the gluing, we apply an isometry of $\HH^3$ fixing $0$ and $\infty$, mapping $1$ to $z(e_1)$. This takes the fourth vertex of $T_2$ to $z(e_1)z(e_2)$.

\begin{figure}
%% Accompanies image file 'F4-08-TriGlu.eps' (pdf, eps, ps)
\begingroup%
  \makeatletter%
  \providecommand\color[2][]{%
    \errmessage{(Inkscape) Color is used for the text in Inkscape, but the package 'color.sty' is not loaded}%
    \renewcommand\color[2][]{}%
  }%
  \providecommand\transparent[1]{%
    \errmessage{(Inkscape) Transparency is used (non-zero) for the text in Inkscape, but the package 'transparent.sty' is not loaded}%
    \renewcommand\transparent[1]{}%
  }%
  \providecommand\rotatebox[2]{#2}%
  \newcommand*\fsize{\dimexpr\f@size pt\relax}%
  \newcommand*\lineheight[1]{\fontsize{\fsize}{#1\fsize}\selectfont}%
  \ifx\svgwidth\undefined%
    \setlength{\unitlength}{188.15379524bp}%
    \ifx\svgscale\undefined%
      \relax%
    \else%
      \setlength{\unitlength}{\unitlength * \real{\svgscale}}%
    \fi%
  \else%
    \setlength{\unitlength}{\svgwidth}%
  \fi%
  \global\let\svgwidth\undefined%
  \global\let\svgscale\undefined%
  \makeatother%
  \begin{picture}(1,0.50470068)%
    \lineheight{1}%
    \setlength\tabcolsep{0pt}%
    \put(0,0){\includegraphics[width=\unitlength]{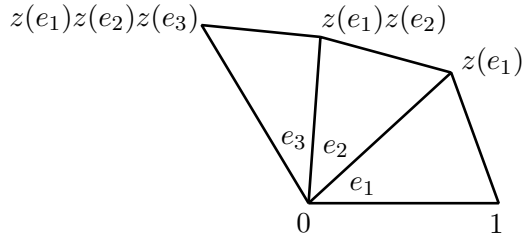}}%
    \put(0.00429205,0.44809977){\color[rgb]{0,0,0}\makebox(0,0)[lt]{\lineheight{1.25}\smash{\begin{tabular}[t]{l}$z(e_1)z(e_2)z(e_3)$\end{tabular}}}}%
    \put(0.62872433,0.44687601){\color[rgb]{0,0,0}\makebox(0,0)[lt]{\lineheight{1.25}\smash{\begin{tabular}[t]{l}$z(e_1)z(e_2)$\end{tabular}}}}%
    \put(0.90886595,0.36125789){\color[rgb]{0,0,0}\makebox(0,0)[lt]{\lineheight{1.25}\smash{\begin{tabular}[t]{l}$z(e_1)$\end{tabular}}}}%
    \put(0.5788675,0.0361634){\color[rgb]{0,0,0}\makebox(0,0)[lt]{\lineheight{1.25}\smash{\begin{tabular}[t]{l}0\end{tabular}}}}%
    \put(0.96526134,0.03207677){\color[rgb]{0,0,0}\makebox(0,0)[lt]{\lineheight{1.25}\smash{\begin{tabular}[t]{l}1\end{tabular}}}}%
    \put(0.68573603,0.11892743){\color[rgb]{0,0,0}\makebox(0,0)[lt]{\lineheight{1.25}\smash{\begin{tabular}[t]{l}$e_1$\end{tabular}}}}%
    \put(0.63179261,0.19370614){\color[rgb]{0,0,0}\makebox(0,0)[lt]{\lineheight{1.25}\smash{\begin{tabular}[t]{l}$e_2$\end{tabular}}}}%
    \put(0.54760805,0.21086999){\color[rgb]{0,0,0}\makebox(0,0)[lt]{\lineheight{1.25}\smash{\begin{tabular}[t]{l}$e_3$\end{tabular}}}}%
  \end{picture}%
\endgroup%

\caption{Vertices of attached triangles.}
\label{Fig:TriGlue}
\end{figure}

Continue attaching tetrahedra counterclockwise around $e$. The next tetrahedron attached will have vertices $0$, $\infty$, $z(e_1) z(e_2)$, and $z(e_1) z(e_2) z(e_3)$ in $\CC$. See \reffig{TriGlue}. Eventually one of the tetrahedra will be glued to $T_1$ again. The fourth vertex of the final tetrahedron will be at the point $z(e_1) z(e_2) \cdots z(e_n)$.

\begin{theorem}[Edge gluing equations]
\label{Thm:Gluing}\index{edge gluing equations}\index{gluing equations!edge equations}
Let $M^3$ admit a topological ideal triangulation such that each ideal tetrahedron has a hyperbolic structure. The hyperbolic structures on the ideal tetrahedra induce a hyperbolic structure on the gluing, $M$, if and only if for each edge $e$,
\[ \prod z(e_i) = 1 \quad \mbox{ and } \quad \sum {\rm arg}(z(e_i)) = 2\pi,\]
where the product and sum are over all edges that glue to $e$.
\end{theorem}

\begin{figure}
\includegraphics{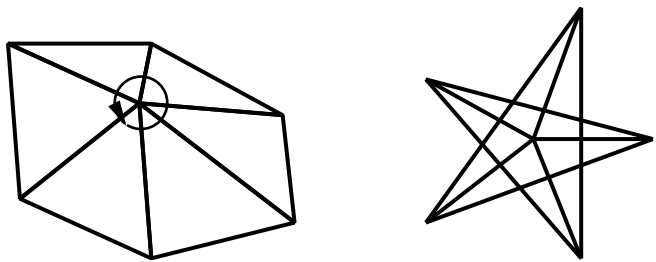}
\caption{Left: Angle sum must be $2\pi$. Right: An example of why this condition is important.}
\label{Fig:AngleSum2pi}
\end{figure}

\begin{proof}
The hyperbolic structure on the tetrahedra induces a hyperbolic structure on $M$ if and only if every point in $M$ has a neighborhood isometric to a ball in $\HH^3$, by \reflem{IsometricNbhd}. Consider a point on an edge. If it has a neighborhood isometric to a ball in $\HH^3$ then the sum of the dihedral angles around the edge must be $2\pi$.  See \reffig{AngleSum2pi}.  This sum of dihedral angles is $\sum\arg(z(e_i))$.  Moreover there must be no nontrivial translation as we move around the edge. Since the last face of the last triangle glues to the triangle with vertices $0$, $1$, and $\infty$, this condition requires that $\prod z(e_i)=1$.

Conversely, if we have $\prod z(e_i)=1$ and $\sum \arg(z(e_i))=2\pi$, then any point on the edge under the gluing has a ball neighborhood isometric to a ball in $\HH^3$.   
\end{proof}

The equations $\prod z(e_i) =1$ (and restrictions $\sum \arg(z(e_i))=2\pi$) are called the \emph{edge gluing equations}\index{edge gluing equations}.  We have one for each edge. However, since by \reflem{EdgeInvariants} the three edge invariants of a tetrahedron are all determined by a single edge invariant, one ideal tetrahedron contributes at most one unknown to the gluing equations. 

\begin{example}[Edge gluing equations for the figure-8 knot]
The figure-8 knot decomposes into two ideal tetrahedra. Choose the two tetrahedra to be \emph{regular}. That is, all dihedral angles are $\pi/3$.  We claim that this gives a hyperbolic structure on the figure-8 knot complement.

We wish to find all such structures.

Thurston worked through this example in detail in his notes; we recall his work here \cite[pages~50--52]{thurston}. 

\Reffig{Fig8tet} shows the two tetrahedra in the decomposition of the figure-8 knot complement, which we obtained in \refchap{Fig8Decomp}. 
These tetrahedra come from the two ideal polyhedra that glue to give the figure-8 knot complement that we discussed in detail in \refchap{Fig8Decomp}; see \reffig{8TopPoly} and \reffig{8BotPoly}. The tetrahedra differ from those in \refchap{Fig8Decomp} in the following ways. First, we have collapsed the bigons.\index{bigon} This gives two remaining edge classes, which we label with one tick mark and with two tick marks. Second, in \reffig{8TopPoly}, we viewed the top ideal polyhedron from the inside; that is, the ball of the polyhedron lay above the plane of projection. To be more consistent in viewing both top and bottom polyhedron, we have rotated our perspective such that now both tetrahedra are viewed from the outside. 

\begin{figure}
  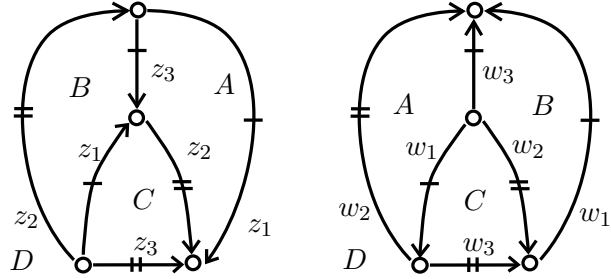
  \caption{The ideal tetrahedra of the figure-8 knot complement.}
  \label{Fig:Fig8tet}
\end{figure}

For each tetrahedron, we label each edge with a complex number $z_i$ or $w_i$, to denote the edge invariant associated with that edge. Note that opposite edges in a tetrahedron have the same edge invariant. We also have relationships between $z_1$, $z_2$, and $z_3$ as in \reflem{EdgeInvariants}, and similarly for $w_1$, $w_2$, and $w_3$. 

There are two edge classes in the tetrahedra in \reffig{Fig8tet}, labeled with one or two tick marks on the edge. We obtain the edge gluing equations by taking the product of edge invariants for all edges identified with each edge class.

For the edge with one tick mark, we obtain the edge gluing equation
\[ z_1^2\,z_3\,w_1^2\,w_3=1. \]
For the edge with two tick marks,
\[ z_2^2\,z_3\,w_2^2\,w_3=1. \]

We set $z_1 = z$ and $w_1 = w$. From \reflem{EdgeInvariants}, the first edge gluing equation gives
\[ z^2\left(\frac{z-1}{z}\right)w^2\left(\frac{w-1}{w}\right)=1, \]
or
\begin{equation}\label{Eqn:Fig8Gluing}
  z\,(z-1)\,w\,(w-1) =1.
\end{equation}
Solving for $z$ in terms of $w$:
\[
z = \frac{1 \pm \sqrt{1+4/(w(w-1))}}{2}.
\]

We need the imaginary parts of $z$ and $w$ to be strictly greater than $0$. For each value of $w$, there is at most one solution for $z$ with positive imaginary part. The solution exists provided that the discriminant $1+4/(w(w-1))$ is not positive real. Thus solutions are parameterized by the region of $\CC$ shown in \reffig{ThurstonRegion} (see also \refex{ThurstonRegion}).

\begin{figure}
%% Accompanies image file 'F4-11-ThuReg.eps' (pdf, eps, ps)
\begingroup%
  \makeatletter%
  \providecommand\color[2][]{%
    \errmessage{(Inkscape) Color is used for the text in Inkscape, but the package 'color.sty' is not loaded}%
    \renewcommand\color[2][]{}%
  }%
  \providecommand\transparent[1]{%
    \errmessage{(Inkscape) Transparency is used (non-zero) for the text in Inkscape, but the package 'transparent.sty' is not loaded}%
    \renewcommand\transparent[1]{}%
  }%
  \providecommand\rotatebox[2]{#2}%
  \newcommand*\fsize{\dimexpr\f@size pt\relax}%
  \newcommand*\lineheight[1]{\fontsize{\fsize}{#1\fsize}\selectfont}%
  \ifx\svgwidth\undefined%
    \setlength{\unitlength}{147bp}%
    \ifx\svgscale\undefined%
      \relax%
    \else%
      \setlength{\unitlength}{\unitlength * \real{\svgscale}}%
    \fi%
  \else%
    \setlength{\unitlength}{\svgwidth}%
  \fi%
  \global\let\svgwidth\undefined%
  \global\let\svgscale\undefined%
  \makeatother%
  \begin{picture}(1,0.87755102)%
    \lineheight{1}%
    \setlength\tabcolsep{0pt}%
    \put(0,0){\includegraphics[width=\unitlength]{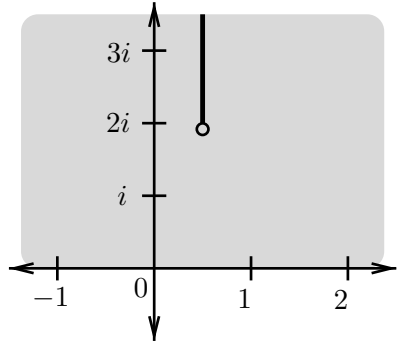}}%
    \put(0.32257437,0.11620111){\color[rgb]{0,0,0}\makebox(0,0)[lt]{\lineheight{1.25}\smash{\begin{tabular}[t]{l}0\end{tabular}}}}%
    \put(0.58848586,0.08740487){\color[rgb]{0,0,0}\makebox(0,0)[lt]{\lineheight{1.25}\smash{\begin{tabular}[t]{l}1\end{tabular}}}}%
    \put(0.8357618,0.08461604){\color[rgb]{0,0,0}\makebox(0,0)[lt]{\lineheight{1.25}\smash{\begin{tabular}[t]{l}2\end{tabular}}}}%
    \put(0.06237029,0.09112329){\color[rgb]{0,0,0}\makebox(0,0)[lt]{\lineheight{1.25}\smash{\begin{tabular}[t]{l}$-1$\end{tabular}}}}%
    \put(0.28067709,0.35050585){\color[rgb]{0,0,0}\makebox(0,0)[lt]{\lineheight{1.25}\smash{\begin{tabular}[t]{l}$i$\end{tabular}}}}%
    \put(0.25190246,0.53824365){\color[rgb]{0,0,0}\makebox(0,0)[lt]{\lineheight{1.25}\smash{\begin{tabular}[t]{l}$2i$\end{tabular}}}}%
    \put(0.25374006,0.7176582){\color[rgb]{0,0,0}\makebox(0,0)[lt]{\lineheight{1.25}\smash{\begin{tabular}[t]{l}$3i$\end{tabular}}}}%
  \end{picture}%
\endgroup%

  \caption{Solutions to edge gluing equations for the figure-8 knot complement are parameterized by the above region.}
  \label{Fig:ThurstonRegion}
\end{figure}

Notice that
\[ z = w = \sqrt[3]{-1} = \frac{1}{2} + \frac{\sqrt{3}}{2}i \]
is one solution to the equations.  We will see that this gives a complete hyperbolic structure on the complement of the figure-8 knot.\index{complete metric space}
\end{example}

%%%%%%%%%%%%%%%%%%%%%%%%%%%%%%%%%%%%%%%%%%%%%%%%%%%%%%%%%%%%%%%%%
\section{Completeness equations}

Suppose now that $M$ is a 3-manifold with torus boundary. In much of this section, we will assume that $M$ admits a topological ideal triangulation, and moreover we have a solution to the edge gluing equations for this triangulation, thus $M$ admits a hyperbolic structure. We need to consider cusps of the manifold to determine whether this is a complete structure or not. 

\begin{definition}\label{Def:CuspNbhd}
Let $M$ be a 3-manifold with torus boundary. Define a \emph{cusp}\index{cusp},
or \emph{cusp neighborhood}\index{cusp neighborhood} of $M$ to be a neighborhood of $\bdy M$ homeomorphic to the product of a torus and an interval, $T^2\times I$. Define a \emph{cusp torus}\index{cusp torus} to be a torus component of $\bdy M$, or the boundary of a cusp. 
\end{definition}

A hyperbolic structure on $M$ induces an affine structure on the boundary of any cusp of $M$.  

\begin{theorem}\label{Thm:EuclidCusp}
Let $M$ be a 3-manifold with torus boundary and hyperbolic structure, i.e.\ with $(\Isom(\HH^3), \HH^3)$-structure. Then the structure on $M$ is complete\index{complete metric space} if and only if for each cusp of $M$, the induced structure on the boundary of the cusp is a Euclidean structure\index{Euclidean structure} on the torus.
\end{theorem}

\begin{proof}
\Refex{EuclidCusp}. Hint: the proof is very similar to that of the analogous result in two dimensions, \refprop{Complete}. 
\end{proof}

\begin{definition}
\label{Def:CuspTriangulation}
Let $M$ have a topological ideal triangulation. If we truncate the vertices of each ideal tetrahedron, we obtain a collection of triangles, each of which lies on the boundary of a cusp. Edges of each triangle inherit a gluing from the gluing of faces of the ideal tetrahedra. This gives a triangulation of each boundary torus, which we call a \emph{cusp triangulation}\index{cusp triangulation}. 
\end{definition}

\begin{figure}
  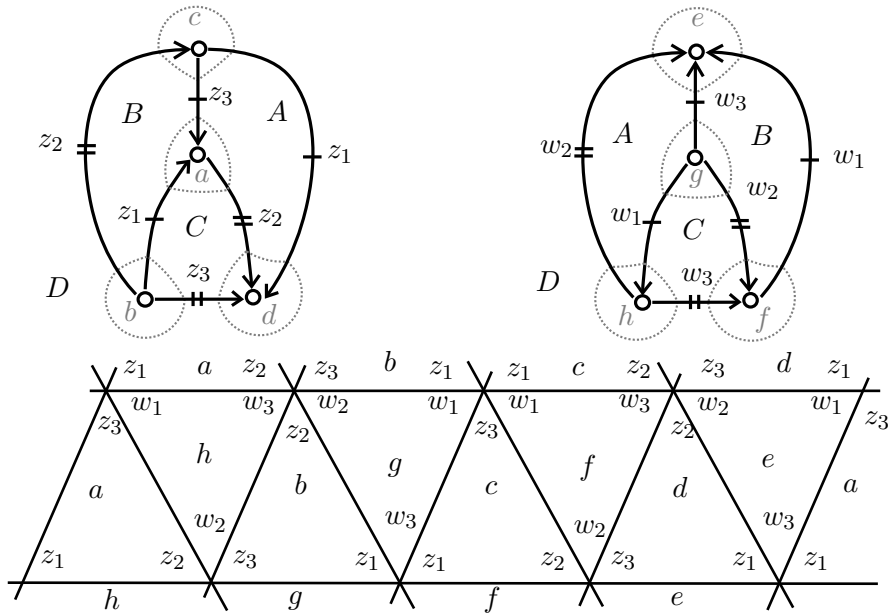
  \caption{Finding the cusp triangulation of the figure-8 knot complement}
  \label{Fig:Fig8CuspTriang}
\end{figure}

An example for the figure-8 knot is shown in \reffig{Fig8CuspTriang}. The truncated ideal vertices give eight triangles, with labels $a$ through $h$. These glue together on the boundary of the cusp to give a triangulation of the torus as shown. Note that the corner of each triangle is labeled with the edge invariant of the tetrahedron corresponding to the edge meeting that corner. 

\Reffig{Fig8CuspTriang} shows a fundamental region of the cusp triangulation. By tracing through gluings of cusp triangles, we may obtain the full developing image of the cusp torus. \Refthm{EuclidCusp} states that the original manifold is complete\index{complete metric space} if and only if the cusp tori are Euclidean, which will hold if and only if the holonomy\index{holonomy} maps for each element of $\pi_1(T)$ on each cusp torus $T$ are pure Euclidean translations, without rotation or scale.

We can determine if holonomy\index{holonomy} maps are Euclidean translations directly from the cusp triangulation. Start with a triangle $\Delta$ whose vertices we may assume lie at $0$, $1$, and $z(e_1)$ in the complex plane $\CC$. Let $\alpha\in\pi_1(T)$. Then the holonomy $\rho(\alpha)$ takes $\Delta$ to a new triangle, which appears in the developing image. The holonomy $\rho(\alpha)$ will be a Euclidean translation if and only if the triangle side from $0$ to $1$ of $\Delta$ is mapped to the side of a triangle of length $1$ pointing in the same direction, without rotation (or scale). To determine whether this holds, we may follow the side of the triangle in the developing image, and obtain exactly its rotation and scale by considering the edge invariants that adjust its length and direction as it is adjusted in the cusp triangulation, as in \reffig{TriGlue}. This can be described efficiently in the following way.

\begin{definition}\label{Def:CompletenessEquations}
Suppose $M$ has a topological ideal triangulation, and let $T$ be the boundary torus of a cusp of $M$. Let $[\alpha]\in\pi_1(T)$, so $\alpha$ is a loop on $T$ in the homotopy class of $[\alpha]$. We associate a complex number $H(\alpha)$ to $\alpha$ as follows.

First, orient the loop $\alpha$ on $T$. The loop $\alpha$ can be homotoped to run through any triangle of the cusp triangulation of $T$ monotonically, i.e.\ in such a way that it cuts off a single corner of each triangle it enters. Denote the edge invariants of the corners cut off by $\alpha$ by $z_1$, $z_2$, \dots, $z_n$. Further associate to each corner a value $\epsilon_i=\pm 1$: if the $i$-th corner cut off by $\alpha$ lies to the left of $\alpha$, set $\epsilon_i=+1$. If the corner lies to the right of $\alpha$, set $\epsilon_i=-1$. Finally, set the value of $H(\alpha)$ to be
\begin{equation}\label{Eqn:CompletenessEquations}
  H(\alpha) = \prod_{i=1}^n z_i^{\epsilon_i}
\end{equation}
\end{definition}

\begin{figure}
%% Accompanies image file 'F4-13-CuspEx.eps' (pdf, eps, ps)
\begingroup%
  \makeatletter%
  \providecommand\color[2][]{%
    \errmessage{(Inkscape) Color is used for the text in Inkscape, but the package 'color.sty' is not loaded}%
    \renewcommand\color[2][]{}%
  }%
  \providecommand\transparent[1]{%
    \errmessage{(Inkscape) Transparency is used (non-zero) for the text in Inkscape, but the package 'transparent.sty' is not loaded}%
    \renewcommand\transparent[1]{}%
  }%
  \providecommand\rotatebox[2]{#2}%
  \newcommand*\fsize{\dimexpr\f@size pt\relax}%
  \newcommand*\lineheight[1]{\fontsize{\fsize}{#1\fsize}\selectfont}%
  \ifx\svgwidth\undefined%
    \setlength{\unitlength}{245.59993172bp}%
    \ifx\svgscale\undefined%
      \relax%
    \else%
      \setlength{\unitlength}{\unitlength * \real{\svgscale}}%
    \fi%
  \else%
    \setlength{\unitlength}{\svgwidth}%
  \fi%
  \global\let\svgwidth\undefined%
  \global\let\svgscale\undefined%
  \makeatother%
  \begin{picture}(1,0.41368051)%
    \lineheight{1}%
    \setlength\tabcolsep{0pt}%
    \put(0,0){\includegraphics[width=\unitlength]{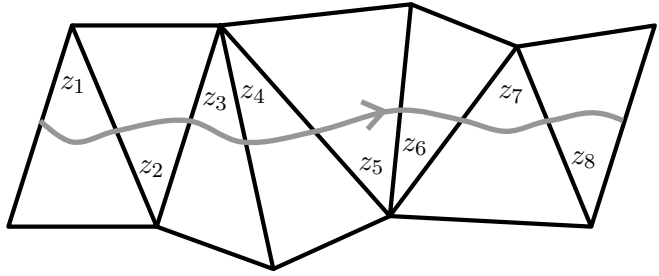}}%
    \put(0.0828799,0.28127925){\color[rgb]{0,0,0}\makebox(0,0)[lt]{\lineheight{0}\smash{\begin{tabular}[t]{l}$z_1$\end{tabular}}}}%
    \put(0.2029712,0.15214354){\color[rgb]{0,0,0}\makebox(0,0)[lt]{\lineheight{0}\smash{\begin{tabular}[t]{l}$z_2$\end{tabular}}}}%
    \put(0.30069109,0.25660136){\color[rgb]{0,0,0}\makebox(0,0)[lt]{\lineheight{0}\smash{\begin{tabular}[t]{l}$z_3$\end{tabular}}}}%
    \put(0.3589262,0.26663777){\color[rgb]{0,0,0}\makebox(0,0)[lt]{\lineheight{0}\smash{\begin{tabular}[t]{l}$z_4$\end{tabular}}}}%
    \put(0.5403859,0.15691383){\color[rgb]{0,0,0}\makebox(0,0)[lt]{\lineheight{0}\smash{\begin{tabular}[t]{l}$z_5$\end{tabular}}}}%
    \put(0.60701235,0.18603138){\color[rgb]{0,0,0}\makebox(0,0)[lt]{\lineheight{0}\smash{\begin{tabular}[t]{l}$z_6$\end{tabular}}}}%
    \put(0.75490668,0.26647255){\color[rgb]{0,0,0}\makebox(0,0)[lt]{\lineheight{0}\smash{\begin{tabular}[t]{l}$z_7$\end{tabular}}}}%
    \put(0.8664496,0.16612358){\color[rgb]{0,0,0}\makebox(0,0)[lt]{\lineheight{0}\smash{\begin{tabular}[t]{l}$z_8$\end{tabular}}}}%
  \end{picture}%
\endgroup%

  \caption{Example for determining $H([\alpha])$}
  \label{Fig:ExampleCusp}
\end{figure}

\begin{example}\label{Example:ExampleCusp}
An example cusp is shown in \reffig{ExampleCusp}. For this example, the value of $H(\alpha)$ is given by
\[ H(\alpha) = z_1\,z_2^{-1}\,z_3\,z_4\,z_5^{-1}\,z_6^{-1}\,z_7\,z_8^{-1}.\]
\end{example}

We will see that $H$ is independent of homotopy class of $\alpha$ (\refex{HomotopyInvarianceH}). For this reason, we sometimes denote the complex number by $H([\alpha])$, or evaluate it on a homotopy class rather than a curve.

\begin{example}\label{Example:Fig8EdgeSequence}
For the figure-8 knot, there is a closed curve on the cusp torus running from the left side of the triangle labeled $a$ on the left of \reffig{Fig8CuspTriang} to the left side of the triangle labeled $a$ on the right of that figure. Call this curve $\alpha$. Then we can compute:
\[ H(\alpha) = z_3\,w_2^{-1}\,z_2\,w_3^{-1}\,z_3\,w_2^{-1}\,z_2\,w_3^{-1} =
\left(\frac{z_2\,z_3}{w_2\,w_3}\right)^2 \]

Another closed curve runs from the base of the triangle labeled $a$ on the left of \reffig{Fig8CuspTriang} to the top of the triangle labeled $h$, also on the left of that figure. Call this curve $\beta$. Then we have:
\[ H(\beta) = z_2^{-1}\,w_1 = \frac{w_1}{z_2}. \]
\end{example}

\begin{proposition}[Completeness equations]\label{Prop:CompletenessEqns}
Let $T$ be the torus boundary of a cusp neighborhood of $M$, where $M$ admits a topological ideal triangulation, and the ideal tetrahedra admit hyperbolic structures that satisfy the edge gluing equations (\refthm{Gluing}). Let $\alpha$ and $\beta$ generate $\pi_1(T)$. If $H(\alpha) = H(\beta)=1$, then the ideal triangulation is a geometric ideal triangulation, i.e.\ the hyperbolic structure on $M$ induced by the hyperbolic structure on the tetrahedra will be a complete structure.\index{complete metric space}
\end{proposition}

The equations $H(\alpha)=1$ and $H(\beta)=1$ are called the \emph{completeness equations}.\index{completeness equations}

\begin{proof}[Proof of \refprop{CompletenessEqns}]
By \refthm{EuclidCusp}, it suffices to show that the induced structure on $T$ is Euclidean. To do so, it suffices to show that the holonomy\index{holonomy} elements $\rho(\alpha)$ and $\rho(\beta)$ are pure translations, with no rotation and scale. Thus we will show $\rho(\alpha)$ and $\rho(\beta)$ do not rotate or scale. 

To show this, let $\Delta$ be a triangle met by the curve $\alpha$ used in defining the complex number $H(\alpha)$, and suppose $\alpha$ meets a side $e_1$ of $\Delta$. Let $v$ be a vector with length equal to the length of $e_1$, pointing in the direction of $e_1$ such that the oriented curve $\alpha$ and the vector $v$ are oriented according to the right hand rule. This is true of the vector $v$ shown on the far left of \reffig{CompletenessEqnsProof}.

\begin{figure}
  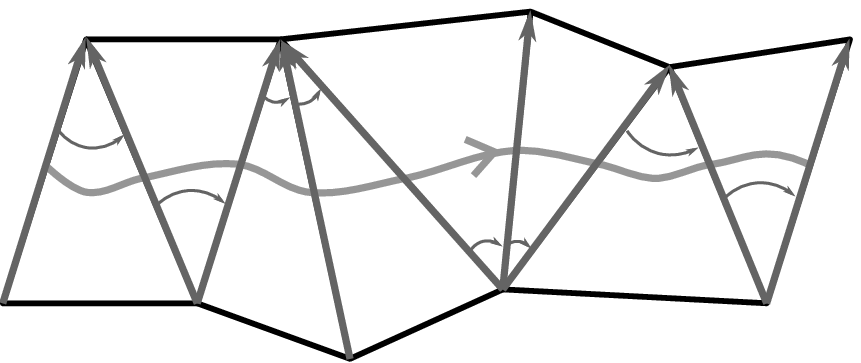
  \caption{A path of vectors in the proof of \refprop{CompletenessEqns}}
  \label{Fig:CompletenessEqnsProof}
\end{figure}

The holonomy\index{holonomy} $\rho(\alpha)$ is Euclidean if and only if the image of $v$ under $\rho(\alpha)$ still has length $v$, and points in the same direction as $v$. We determine the effect of holonomy by considering what happens to $v$ in each triangle of the cusp triangulation. 

We may rotate $v$ around a vertex of the triangle $\Delta$ meeting $e_1$, and scale, so that the result lines up with a second edge $e_2$ of the triangle, having the same length and direction as $e_2$. We know exactly how the rotation and scale is determined when the vertex of the triangle is labeled with edge invariant $z_1$: if we rotate in a counterclockwise direction, $v$ is adjusted by multiplication by $z_1$, as in \reffig{TriGlue}. If we rotate in a clockwise direction, $v$ is adjusted by multiplication by $1/z_1$.

Now, our path $\alpha$ cuts off exactly one corner of each triangle it meets. This defines a path of edges of triangles, namely, starting with $v$, at each step we have a vector lying on the side of a triangle where $\alpha$ enters that triangle. In this triangle, rotate through the corner cut off by $\alpha$ to produce a new vector pointing in the direction of the side where $\alpha$ exits. An example path of such vectors is shown in \reffig{CompletenessEqnsProof}. When $\alpha$ returns to the initial triangle $\Delta$, the final vector of this path will be parallel to the image of $v$ under $\rho(\alpha)$. Then $\rho(\alpha)$ will be a Euclidean transformation if and only if the final vector in the path has length and direction identical to that of $v$.

On the other hand, the final length and direction of the vector $\rho(\alpha)v$ is given by the product of edge invariants at the corners of each triangle in the path of edges, with edge invariant either multiplied or divided depending on whether the rotation is in the counterclockwise or clockwise direction, respectively. This is exactly the complex number $H(\alpha)$. Thus $\rho(\alpha)$ is Euclidean if and only if $H(\alpha)=1$.

The same argument applies to $H(\beta)$ and $\rho(\beta)$. Since the holonomy group\index{holonomy group} of the cusp is generated by $\rho(\alpha)$ and $\rho(\beta)$, the cusp will be Euclidean if and only if $H(\alpha)=H(\beta)=1$.
\end{proof}

\begin{example}\label{Example:Fig8Completeness}
Returning to the example of the figure-8 knot, in \refexamp{Fig8EdgeSequence}, we found that completeness equations\index{completeness equations} are given by
\[ H(\alpha) = \left(\frac{z_2\,z_3}{w_2\,w_3}\right)^2 \quad \mbox{ and } \quad H(\beta) = \frac{w_1}{z_2}. \]
\Reflem{EdgeInvariants} implies that these can be rewritten in terms of variables $z$ and $w$ alone, as
\[
H(\alpha) = \left( \frac{1}{1-z} \cdot \frac{z-1}{z} \cdot \frac{1-w}{1} \cdot \frac{w}{w-1}\right)^2 = \left( \frac{w}{z} \right)^2 \]
and
\begin{equation}\label{Eqn:2ndFig8Complete}
H(\beta) = w\,(1-z)
\end{equation}

If the hyperbolic structure is complete,\index{complete metric space} then by
\refprop{CompletenessEqns}, $H(\alpha)=H(\beta)=1$, so $z=w$.

From \refeqn{2ndFig8Complete}, $z(z-1)=-1$.  Hence the only possibility is $z=w= \frac{1}{2}+i\frac{\sqrt{3}}{2}$.
\end{example}

\section{Computing hyperbolic structures}

Given a triangulation of a 3-manifold $M$ with torus boundary, we may determine a complete hyperbolic structure on $M$ by solving the edge gluing and completeness equations.\index{completeness equations}\index{complete metric space} However, note this amounts to solving a complicated system of nonlinear equations. Consequently, it is difficult to use these to find exact hyperbolic structures on infinite families of manifolds. 

However in practice, topological triangulations, edge gluing equations, and completeness equations\index{completeness equations} can be found very efficiently by computer for specific, finite examples. The resulting nonlinear system of equations can then be solved numerically. The first software to find hyperbolic structures on knots and 3-manifolds was the program SnapPea,\index{SnapPea} written by Weeks \cite{Weeks:Thesis} (see also \cite{weeks:computation}). This program has allowed researchers to run experiments on large classes of hyperbolic 3-manifolds, making observations and testing conjectures, and has been influential in a great deal of results on hyperbolic structures on knots and 3-manifolds.
The SnapPea kernel is now part of a program maintained by Culler, Dunfield, Goerner, and others, reincarnated as SnapPy, and available for free download \cite{SnapPy}.\index{SnapPy} This new program includes much additional functionality, and still remains an excellent tool for research in hyperbolic knot theory.

One issue in the past with finding a hyperbolic structure via SnapPea (SnapPy) is that it would only give a numerical approximation to a hyperbolic structure, and there was no guarantee that the manifold would be actually provably hyperbolic. This has been addressed in a few ways. The program Snap \cite{Snap} deduces exact solutions from the numerical approximations, which can be used to prove hyperbolicity. In another direction, Moser used analytic techniques to prove that a solution to edge gluing and completeness equations\index{completeness equations} exists in a small neighborhood of an approximate solution \cite{Moser:ProvingHyp}. In \cite{HIKMOT}, interval arithmetic is used to prove hyperbolic structures exist when a structure is computed numerically. Thus using these tools, we can often prove that if SnapPy computes a hyperbolic structure on a knot complement, then the knot is indeed hyperbolic.

%%%%%%%%%%%%%%%%%%%%%%%%%%%%%%%%%%%%%%%%%%%%%%%%%%%%%%%%%%%%%%%%%

%%%%%%%%%%%%%%%%%%%%%%%%%%%%%%%%%%%%%%%%%%%%%%%%%%%%%%%%%%%%%%%%%
\section{Exercises}

\begin{exercise}
Write down the edge gluing equations (not completeness equations) for the $6_1$ knot, using the ideal tetrahedra of \refexamp{61TopTriang}.  Make appropriate substitutions such that your equations contain exactly one variable per tetrahedron.
\end{exercise}

\begin{exercise}\label{Ex:edge=tet}
Notice that for both the figure-8 knot complement and for the $6_1$ knot, we had exactly the same number of edges as tetrahedra in the ideal triangulation.
\begin{enumerate}
\item[(a)] Prove that this will always be true. That is, prove that if $M$ is any 3-manifold with (possibly empty) boundary consisting of tori, then for any topological ideal triangulation of $M$, the number of edges of the triangulation will always equal the number of tetrahedra.
\item[(b)] Since we have one unknown per ideal tetrahedra, part (a) implies that the number of gluing equations will equal the number of unknowns.  However, in fact the gluing equations are always redundant.  Prove this fact.
\end{enumerate}
\end{exercise}

\begin{exercise}\label{Ex:52KnotTriang}
In \refchap{Fig8Decomp}, we found a polyhedral decomposition of the $5_2$ knot complement (without bigons).\index{bigon}
Split this into a topological ideal triangulation of the knot complement.
\end{exercise}

\begin{exercise} Using the ideal tetrahedra of \refex{52KnotTriang}, or otherwise, 
write down all edge invariants and all edge gluing equations, one variable per tetrahedron.
\end{exercise}

\begin{exercise}
Find a topological ideal triangulation of the $6_3$ knot, edge invariants, and edge gluing equations.
\end{exercise}

\begin{exercise}\label{Ex:ThurstonRegion}
Check that \reffig{ThurstonRegion} does indeed parameterize the space of hyperbolic structures on the figure-8 knot complement. What is the equation of the vertical ray shown in that picture? 
\end{exercise}

\begin{exercise}\label{Ex:EuclidCusp}
Prove \refthm{EuclidCusp}: the hyperbolic structure on $M$ is complete\index{complete metric space} if and only if for each cusp of $M$, the induced structure on the boundary of the cusp is a Euclidean structure\index{Euclidean structure} on the torus.
\end{exercise}

\begin{exercise}
For the topological triangulation of the $5_2$ knot of \refex{52KnotTriang}:
\begin{enumerate}
\item[(a)] Find the triangulation of the cusp.  Label a fundamental domain, and meridian and longitude.
\item[(b)] Write down completeness equations.\index{completeness equations}
\end{enumerate}
\label{Ex:52CuspTri}
\end{exercise}

\begin{exercise}
Find the cusp triangulation for the complement of the $6_1$ knot from \refexamp{61TopTriang}.  
\end{exercise}

\begin{exercise}
Find completeness equations\index{completeness equations} for the $6_1$ or $6_3$ knot.
\end{exercise}

\begin{exercise}\label{Ex:HomotopyInvarianceH}
Suppose $M$ admits an ideal triangulation that satisfies the edge gluing equations.
\begin{enumerate}
\item [(a)] In \refdef{CompletenessEquations}, we claimed that for any closed curve $\alpha$ in a torus boundary component of $\bdy M$, we could homotope $\alpha$ in such a way that it cuts off a single corner of each triangle that it meets. Prove this.

\item[(b)] Show that $H([\alpha])$ is independent of the choice of $\alpha$ in the homotopy class of $[\alpha]$. In particular, if $\alpha$ is homotoped to run through different triangles, the value of $H([\alpha])$ is unchanged.
\end{enumerate}
\end{exercise}

\begin{figure}[h]
  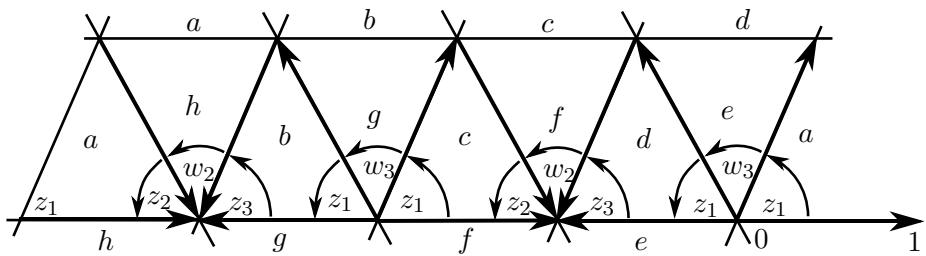
  \caption{Path of vectors going from $e_1$, the oriented edge from $0$ to $1$, to $\rho(-\alpha)(e_1)$}
  \label{Fig:ThurstonFig8}
\end{figure}

\begin{exercise}
In Thurston's 1979 notes \cite{thurston}, he computed completeness equations\index{completeness equations} for the figure-8 knot using a method similar to our proof of \refprop{CompletenessEqns}. Namely, he found a path of vectors from an edge on a triangle $\Delta$ to the same edge on $\rho(-\alpha)(\Delta)$ and $\rho(\beta)(\Delta)$. His path of vectors for $\rho(\beta)$ agrees with ours. His path of vectors for $\rho(-\alpha)$ is different from our path for $\rho(\alpha)$, and is shown in \reffig{ThurstonFig8}.
\begin{enumerate}
\item[(a)] Prove that the completeness equation obtained from Thurston's path of vectors is equivalent to our completeness equation.
\item[(b)] More generally, prove that if we replace our path of vectors used to construct the complex number $H([\alpha])$ by any other path of vectors obtained by rotating around vertices of the cusp triangulation, with same starting and ending vectors, then the equation we obtain from multiplying (and dividing) by edge invariants corresponding to the path of vectors gives a completeness equation that is equivalent to $H([\alpha])=1$. 
\end{enumerate}
\end{exercise}

\begin{exercise}
What breaks down when you try to find triangulations and edge gluing equations for non-hyperbolic knots and links, such as the trefoil or the $(2,4)$-torus link?
\end{exercise}

\begin{exercise}
Use the computer program SnapPy to determine which of the knots with seven or fewer crossings admit a hyperbolic structure \cite{SnapPy}. For those that do admit a hyperbolic structure, use SnapPy to find the cusp triangulation of the knot. Obtain a screen shot of this information, which should include cusp triangles as well as a fundamental parallelogram for the cusp. 
\end{exercise}

%% Ch05_Margulis.tex

\chapter{Discrete Groups and the Thick--Thin Decomposition}\label{Chap:Margulis}
\blfootnote{Jessica S. Purcell, Hyperbolic Knot Theory}

Suppose we have a complete hyperbolic structure on an orientable 3-manifold $M$. Then the developing map $D\from \widetilde{M}\to\HH^3$ is a covering map, by \refthm{Developing}. Since $\widetilde{M}$ and $\HH^3$ are both simply connected, it follows that the developing map is an isometry. Thus we may view $\HH^3$ as the universal cover of $M$. The covering transformations are then the elements of the holonomy group\index{holonomy group} $\rho(\pi_1(M)) = \Gamma \leq \PSL(2,\CC)$. Hence $M$ is homeomorphic to the quotient $M\cong \HH^3/\Gamma$.

Subgroups $\Gamma$ of $\PSL(2,\CC)$ can have very nice properties, and have been investigated for many decades. In this chapter, we discuss some classical results in the area and their consequences for hyperbolic 3-manifolds. Some of our discussion follows closely work of J{\o}rgensen and Marden; we recommend the book~\cite{marden} for more details, generalizations, and consequences. 

\section{Discrete subgroups of hyperbolic isometries}

\subsection{Isometries and subgroups}
In \refthm{ClassifyPSL(2,C)} we classified elements of $\PSL(2,\CC)$ as elliptic,\index{elliptic} parabolic,\index{parabolic} or loxodromic\index{loxodromic} depending on their fixed points. One of the first things we need is an extension of that theorem.

Before we give the extension, recall that we can view an element of $\PSL(2,\CC)$ as a matrix
\[ A = \mat{a&b\\ c&d}, \mbox{ with } a,b,c,d\in \CC \mbox{ and } ad-bc=1,\]
and the matrix is well-defined up to multiplication by $\pm \Id$. In this chapter, we will frequently write an isometry of $\HH^3$ as a 2 by 2 matrix with determinant $1$, omitting and ignoring the $\pm$ sign. The sign very rarely affects our arguments, but the reader should be aware that we are suppressing it, for example in the following definition. 

\begin{definition}\label{Def:ConjugateTrace}
We say $A\in\PSL(2,\CC)$ is \emph{conjugate}\index{conjugate} to $B\in\PSL(2,\CC)$ if there exists $U\in\PSL(2,\CC)$ such that $A=UBU^{-1}$. The \emph{trace}\index{trace} of $A$ is the trace of its normalized matrix:
\[ \tr\mat{a&b\\c&d} = a+d. \]
\end{definition}
Note there is a sign ambiguity in our definition of trace; again this will not affect our arguments. Note also that conjugate elements have the same trace. 

\begin{lemma}\label{Lem:MoreClassifyPSL}
For $A\in\PSL(2,\CC)$,
\begin{itemize}
\item $A$ is \emph{parabolic}\index{parabolic} if and only if $\tr(A)=\pm 2$, and if and only if $A$ is conjugate to \[ z\mapsto z+1.\]
\item $A$ is \emph{elliptic}\index{elliptic} if and only if $\tr(A)\in(-2,2)\subset\RR\subset\CC$, and if and only if $A$ is conjugate to
\[ z\mapsto e^{2i\theta}z, \quad \mbox{ with } 2\theta\neq 2\pi n \mbox{ for any } n\in\ZZ. \]
\item $A$ is \emph{loxodromic}\index{loxodromic} if and only if $\tr(A)\in\CC-[-2,2]$, and if and only if $A$ is conjugate to
\[ z\mapsto \zeta^2z, \quad \mbox{with } |\zeta|>1.\]
\end{itemize}
\end{lemma}

\begin{proof}
\Refex{MoreClassifyPSL}
\end{proof}

\begin{definition}\label{Def:DiscreteGroup}
A subgroup of $\PSL(2,\CC)$ is said to be \emph{discrete}\index{discrete group} if it contains no sequence of distinct elements converging to the identity element.
A discrete subgroup of $\PSL(2,\CC)$ is often called a \emph{Kleinian group}\index{Kleinian group}. 
\end{definition}

An example of a discrete group is a subgroup generated by a single loxodromic\index{loxodromic} element, or a single parabolic\index{parabolic} element. These are the simplest such groups. They are so simple that they are examples of what are called \emph{elementary groups}; see \refdef{Elementary}. Examples of discrete groups in general can be quite complicated. In \refprop{FreePropDisc}, we will prove that the holonomy group\index{holonomy group} of any complete\index{complete metric space} hyperbolic 3-manifold is always a discrete group. Meanwhile, consider the example of the figure-8 knot complement. 

\begin{example}\label{Example:Fig8GroupDiscrete}
Let $K$ be the figure-8 knot, and give $S^3-K$ its complete hyperbolic structure by gluing two regular ideal tetrahedra, with face-pairings\index{face-pairing isometry} as in \reffig{Fig8tet}. We will find generators of the holonomy group,\index{holonomy group} which is a discrete subgroup of $\PSL(2,\CC)$, as we will see in \refprop{FreePropDisc}. These are obtained by face-pairing isometries,\index{face-pairing isometry} as follows. 

Place the two ideal tetrahedra in $\HH^3$, putting ideal vertices for one tetrahedron at $0$, $1$, $\omega$, and $\infty$, where $\omega = \half + i\frac{\sqrt{3}}{2}$, and putting the ideal vertices of the other tetrahedron at $1$, $\omega$, $\omega+1$, and $\infty$. This glues the faces labeled $A$ along the ideal triangle\index{ideal triangle} with vertices $1$, $\omega$, and $\infty$, to obtain one connected fundamental region for the knot complement, shown in \reffig{Fig8FundRegion}.

\begin{figure}
%% Accompanies image file 'F5-01-FunReg.eps' (pdf, eps, ps)
\begingroup%
  \makeatletter%
  \providecommand\color[2][]{%
    \errmessage{(Inkscape) Color is used for the text in Inkscape, but the package 'color.sty' is not loaded}%
    \renewcommand\color[2][]{}%
  }%
  \providecommand\transparent[1]{%
    \errmessage{(Inkscape) Transparency is used (non-zero) for the text in Inkscape, but the package 'transparent.sty' is not loaded}%
    \renewcommand\transparent[1]{}%
  }%
  \providecommand\rotatebox[2]{#2}%
  \newcommand*\fsize{\dimexpr\f@size pt\relax}%
  \newcommand*\lineheight[1]{\fontsize{\fsize}{#1\fsize}\selectfont}%
  \ifx\svgwidth\undefined%
    \setlength{\unitlength}{108bp}%
    \ifx\svgscale\undefined%
      \relax%
    \else%
      \setlength{\unitlength}{\unitlength * \real{\svgscale}}%
    \fi%
  \else%
    \setlength{\unitlength}{\svgwidth}%
  \fi%
  \global\let\svgwidth\undefined%
  \global\let\svgscale\undefined%
  \makeatother%
  \begin{picture}(1,0.88750553)%
    \lineheight{1}%
    \setlength\tabcolsep{0pt}%
    \put(0,0){\includegraphics[width=\unitlength]{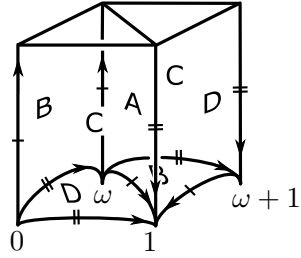}}%
    \put(0.28714522,0.18425914){\color[rgb]{0,0,0}\makebox(0,0)[lt]{\lineheight{0}\smash{\begin{tabular}[t]{l}$\omega$\end{tabular}}}}%
    \put(-0.00700774,0.02006241){\color[rgb]{0,0,0}\makebox(0,0)[lt]{\lineheight{0}\smash{\begin{tabular}[t]{l}$0$\end{tabular}}}}%
    \put(0.45784954,0.01632118){\color[rgb]{0,0,0}\makebox(0,0)[lt]{\lineheight{0}\smash{\begin{tabular}[t]{l}$1$\end{tabular}}}}%
    \put(0.76931309,0.1706501){\color[rgb]{0,0,0}\makebox(0,0)[lt]{\lineheight{0}\smash{\begin{tabular}[t]{l}$\omega+1$\end{tabular}}}}%
  \end{picture}%
\endgroup%

\caption{A connected fundamental region for the figure-8 knot complement}
\label{Fig:Fig8FundRegion}
\end{figure}

The manifold $S^3-K$ is obtained by gluing the remaining faces labeled $B$, $C$, and $D$. These gluings, or face-pairings,\index{face-pairing isometry} correspond to holonomy\index{holonomy} isometries, which we will denote by $T_B$, $T_C$, and $T_D$, respectively.
A calculation (\refex{Fig8Gluings}) shows that the gluing isometries are given by:
\begin{equation}\label{Eqn:Fig8Gluings}
T_B = \frac{i}{\sqrt{\omega}} \mat{1&1\\1&-\omega^2}, \quad T_C=\mat{1&\omega\\0&1}, \quad T_D = \mat{2&-1\\1&0}.
\end{equation}
These three gluing isometries generate the holonomy group\index{holonomy group} for $S^3-K$. In fact, $T_B$ can be written as a (somewhat complicated) product involving $T_C$ and $T_D$ and their inverses.

Riley was the first to prove that $S^3-K$ has a hyperbolic structure \cite{Riley:Fig8}. He did so by taking a presentation of the fundamental group of $S^3-K$ with two generators, and finding an explicit representation of the fundamental group into $\PSL(2,\CC)$. Exercises~\ref{Ex:Riley},
\ref{Ex:Riley2},
and~\ref{Ex:Riley3} explore this work a little further. 
\end{example}

We finish this section with one quick condition equivalent to a group being discrete.

\begin{lemma}\label{Lem:EquivDiscrete}
A subgroup $G\leq \PSL(2,\CC)$ is discrete if and only if it does not contain an infinite sequence of distinct elements that converges to some element $A \in \PSL(2,\CC)$.
\end{lemma}

\begin{proof}
One implication is trivial: If $G$ is not discrete, by definition it contains an infinite sequence of distinct elements converging to the identity in $\PSL(2,\CC)$.

For the other direction, suppose $\{A_n\} \subset G$ is an infinite sequence of distinct elements of $G$ converging to $A\in\PSL(2,\CC)$. Consider $\{A_{n+1}A_n^{-1}\}\subset G$. Note the sequence converges to the identity. To show $G$ is not discrete, it remains to show that $\{A_{n+1}A_n^{-1}\}$ contains infinitely many distinct elements. Suppose not. Then $A_{n+1}=CA_n$ for some fixed $C\in G$ and some subsequence. Since $A_{n+1}A_n^{-1}\to\Id$, we must have $C=\Id$, and thus $A_{n+1}=A_n$. This contradicts the fact that $\{A_n\}$ is a sequence of distinct elements. Thus $G$ is not discrete.
\end{proof}

\subsection{Sequences of isometries}
We can learn a lot about subgroups of $\PSL(2,\CC)$ by considering sequences of group elements. For example, note that the definition of a discrete group involves sequences. We also have the following result, which will be used later in the chapter.

\begin{lemma}\label{Lem:Sequences}
Let $\{A_n\}$ be a sequence of elements of $\PSL(2,\CC)$. Then either a subsequence of $\{A_n\}$ converges to some $A\in\PSL(2,\CC)$, or there exists a point $q\in\bdy\HH^3$ such that for all $x\in\HH^3$, the sequence $\{A_n(x)\}$ has a subsequence converging to $q$. 
\end{lemma}

\begin{proof}
Let $p_n$, $q_n$ denote the fixed points of $A_n$; note we could have $p_n=q_n$. Then $\{p_n\}$ and $\{q_n\}$ are sequences in $\bdy\HH^3 \cong S^2$, which is compact, so they have convergent subsequences. Replace $A_n$, $p_n$, $q_n$ by a subsequence such that $p_n\to p$ and $q_n\to q$. Again note that $p$ could equal $q$.

\smallskip

\textbf{Case 1.} Suppose $p\neq q$. Then for large enough $n$, $p_n\neq q_n$.
Consider an isometry $R_n$ of $\HH^3$ mapping $p_n$ to $0$ and $q_n$ to $\infty$.
Furthermore, for concreteness, fix a point $y\in \bdy\HH^3$, independent of $n$ that is disjoint from the sequences $\{p_n\}$, $\{q_n\}$ and from $p$ and $q$. We may take $R_n$ to map $y$ to $1$.

If we view $R_n$ as a sequence of matrices for example, we see that $R_n$ converges to the hyperbolic isometry $R\in\PSL(2,\CC)$ taking $p$ to $0$, $q$ to $\infty$, and $y$ to $1$. Consider $B_n = R_nA_nR_n^{-1}$. This is an isometry in $\PSL(2,\CC)$ fixing $0$ and $\infty$. Hence it has the form $B_n(z) = a_nz$ for $a_n\in\CC$. If $\{|a_n|\}$ has a bounded subsequence, then some subsequence $a_n$ converges to $a\in\CC$. Hence there is a subsequence $B_n$ with $B_n \to B$, where $B$ is the hyperbolic isometry $B(z)=az$. This is an element of $\PSL(2,\CC)$. It follows that $A_n = R_n^{-1}B_nR_n$ converges to $A=R^{-1}B R \in \PSL(2,\CC)$.

If $|a_n|\to\infty$, then for any $z\in\bdy\HH^3$, $B_n(z)\to\infty$. Thus for any point $x\in\HH^3$, $B_n(x)\to\infty$. It follows that for all $x\in\HH^3$, $A_n(x) = R_n^{-1}B_nR_n(x)$ converges to $q\in\bdy\HH^3$.

\smallskip

\textbf{Case 2.} Now suppose $p=q$. Then again we will conjugate $A_n$ by an isometry $R_n$ taking $q_n$ to infinity. For concreteness, choose $y_1$ and $y_2$ disjoint from $\{p_n\}$, $\{q_n\}$, and $q$. Let $R_n$ be the isometry taking $y_1$ to $1$, $y_2$ to $0$, and $q_n$ to $\infty$. Then $R_n$ converges to the isometry $R$ taking $y_1$, $y_2$, and $q$ to $1$, $0$, and $\infty$, respectively. Finally let $B_n = R_nA_nR_n^{-1}$. Note $B_n$ fixes $\infty$, hence it is of the form $B_n=a_nz+b_n$ for $a_n$, $b_n\in\CC$. If $a_n=1$, $B_n$ is parabolic\index{parabolic} and has unique fixed point $\infty$. Otherwise, the other fixed point of $B_n$ is $b_n/(1-a_n)$.

If $\{|b_n|\}$ has a bounded subsequence, then some subsequence $b_n\to b$. In that case, either $a_n=1$ for large $n$, or since $p_n, q_n$ converge to $p=q$, the fixed point $b_n/(1-a_n)$ converges to $\infty$. Thus $a_n$ converges to $1$. In any case, $B_n(z)$ converges to $B(z)=z+b$. This is an element of $\PSL(2,\CC)$. It follows that $A_n=R_n^{-1}B_nR_n$ converges to $R^{-1}BR\in\PSL(2,\CC)$.

If $\{|b_n|\}$ has no bounded subsequence, then $b_n\to\infty$. We know that the fixed point $b_n/(1-a_n)$ converges to $\infty$ because it is a fixed point of $B_n$, so $(1-a_n)/b_n\to 0$. Rewrite $B_n$ to have the form
\[ B_n(z) = b_n\left(\frac{(a_n-1)z}{b_n}+1\right)+z.\]
Then as $n\to\infty$, $B_n(z) \to \infty$ for all $z\in \bdy\HH^3$. Thus $B_n(x)\to\infty$ for all $x\in\HH^3$. It follows that $A_n(x) = R_n^{-1}B_nR_n(x)$ converges to $q$ for all $x\in\HH^3$. 
\end{proof}

\subsection{Action of groups of isometries}

We return to the problem of showing that holonomy groups\index{holonomy group} of complete hyperbolic 3-manifolds are discrete. We will show this by considering the action of these groups on $\HH^3$. 

\begin{definition}\label{Def:ProperlyDiscont}
The action of a group $G\leq \PSL(2,\CC)$ on $\HH^3$ is \emph{properly discontinuous}\index{properly discontinuous}\index{group action!properly discontinuous} if for every closed ball $B\subset \HH^3$, the set $\{\gamma\in G \mid \gamma(B)\cap B\neq \emptyset\}$ is a finite set.
\end{definition}

\begin{definition}\label{Def:FreeAction}
The action of a group $G\leq\PSL(2,\CC)$ is \emph{free}\index{free group action}\index{group action!free} if the identity element of $G$ is the only element to have a fixed point in $\HH^3$. 
\end{definition}

Note that parabolics\index{parabolic} and loxodromics\index{loxodromic} have fixed points on $\bdy \HH^3$, but not in the interior of $\HH^3$. However, elliptics\index{elliptic} have fixed points in the interior of $\HH^3$. Thus the action of $G$ is free\index{free group action}\index{group action!free} if and only if $G$ contains no elliptics.

\begin{lemma}\label{Lem:DiscretePropDisc}
A subgroup of $\PSL(2,\CC)$ is discrete if and only if its action on $\HH^3$ is properly discontinuous.\index{properly discontinuous}\index{group action!properly discontinuous}
\end{lemma}

\begin{proof}
Suppose $G$ is a subgroup of $\PSL(2,\CC)$ that is not discrete, so there exists a sequence $\{A_n\}$ in $G$ with $A_n\to\Id$. Then for all $x\in \HH^3$, the hyperbolic distance $d(x,A_nx)\to 0$. Let $B$ be any closed ball about $x$ with radius $R>0$. For $n$ such that $d(x,A_nx)<R$, the set
\[ \{A\in G \mid A(B)\cap B\neq \emptyset\} \] contains $A_n$. Since this is true for infinitely many $A_n$, the action is not properly discontinuous.\index{properly discontinuous}\index{group action!properly discontinuous}

Now suppose that for $G\leq\PSL(2,\CC)$, there exists a closed ball $B$ of radius $R$ such that the set $\{A\in G\mid A(B)\cap B\neq\emptyset\}$ is infinite. Let $\{A_n\}$ be a sequence of distinct elements in this set. Note that for $x\in B$, the hyperbolic distance $d(x,A_nx)$ is bounded by $4R$, for all $n$. Thus $\{A_nx\}$ has no subsequence converging to a point on $\bdy\HH^3$. \Reflem{Sequences} implies that $\{A_n\}$ has a subsequence converging to $A\in\PSL(2,\CC)$. Then \reflem{EquivDiscrete} implies $G$ is not discrete.
\end{proof}

We are now ready to prove the main result in this section, namely that a complete hyperbolic 3-manifold has a discrete holonomy group,\index{holonomy group} and conversely a discrete subgroup of $\PSL(2,\CC)$ that acts freely\index{free group action}\index{group action!free} gives rise to a complete hyperbolic 3-manifold. 

\begin{proposition}\label{Prop:FreePropDisc}
The action of a group $G\leq\PSL(2,\CC)$ on $\HH^3$ is free\index{free group action}\index{group action!free} and properly discontinuous\index{properly discontinuous}\index{group action!properly discontinuous} if and only if $\HH^3/G$ is a 3-manifold with a complete hyperbolic structure and with covering projection $\HH^3\to\HH^3/G$.
\end{proposition}

\begin{proof}
Suppose the action of $G$ on $\HH^3$ is free\index{free group action}\index{group action!free} and properly discontinuous.\index{properly discontinuous}\index{group action!properly discontinuous} Let $x\in\HH^3/G$, and let $\widetilde{x}\in\HH^3$ be a point that projects to $x$ under the map $\HH^3\to\HH^3/G$. Because the action of $G$ is properly discontinuous, there is a closed ball $B_x$ that intersects only finitely many of its translates. Because the action is free,\index{free group action}\index{group action!free} we may shrink $B_x$ until all its translates are disjoint. Then the interior of $B_x$ maps isometrically to a neighborhood of $x$ in $\HH^3/G$, so $\HH^3/G$ is a hyperbolic manifold. Moreover, this neighborhood is evenly covered (by translates of the interior of $B_x$), and so the quotient map is a covering projection.

Conversely, suppose $\HH^3/G$ is a hyperbolic manifold and $p\from \HH^3\to \HH^3/G$ is a covering projection. For any $x\in\HH^3$, the action of $G$ permutes the preimages $\{p^{-1}p(x)\}$. Only the identity of $G$ fixes $x$, so the action is free.\index{free group action}\index{group action!free}

Let $B\subset\HH^3$ be a closed ball. Consider the compact set $B\times B$. For any $(x,y)\in B\times B$, we claim there exist neighborhoods $U_{xy}$ of $x$ and $V_{xy}$ of $y$ such that $g(U_{xy})\cap V_{xy}\neq \emptyset$ for at most one $g\in G$. To see this, if $y$ is not in the orbit of $x$, then $p(x)$ and $p(y)$ have disjoint neighborhoods in $\HH^3/G$. Shrink these neighborhoods to be evenly covered, and let $U_{xy}$ and $V_{xy}$ be neighborhoods of $x$ and $y$ respectively homeomorphic to the disjoint neighborhoods of $p(x)$ and $p(y)$. For any $g\in G$, $g(U_{xy})\cap V_{xy} = \emptyset$ in this case. On the other hand, if $y=g_1(x)$ for some $g_1\in G$, then take $U_{xy}$ to be homeomorphic to an evenly covered neighborhood of $p(x)=p(y)$ in $\HH^3/G$, and let $V_{xy}=g_1(U_{xy})$. Then $g(U_{xy})\cap V_{xy}\neq\emptyset$ only when $g=g_1$.

Now $B\times B$ is compact, and the set $\{ U_{xy}\times V_{xy}\}_{(x,y)\in B\times B}$ forms an open cover. Thus there is a finite subcover $\{U_1\times V_1, \dots, U_n\times V_n\}$, where $U_i\times V_i$ has the property that $g(U_i)\cap V_i\neq \emptyset$ only when $g=g_i\in G$.

If $\gamma\in G$ is a group element such that there exists $x \in \gamma(B)\cap B$, then consider $(\gamma^{-1}(x),x) \in B\times B$. There must be some $U_i\times V_i$ containing $(\gamma^{-1}(x),x)$. Since $x\in \gamma(U_i)\cap V_i$, it follows that $\gamma=g_i$. Thus $\gamma$ must be one of the elements $g_1, \dots, g_n$ associated to the finite covering. It follows that the action is properly discontinuous.\index{properly discontinuous}\index{group action!properly discontinuous}
\end{proof}

\Refprop{FreePropDisc} implies that if $\HH^3/G$ is a hyperbolic 3-manifold, then $G$ contains no elliptics\index{elliptic}. For this reason, we will exclude elliptic elements from discrete groups $G$ whenever possible to simplify our proofs in the rest of the chapter. In fact, many results below also hold for discrete groups that contain elliptics. Details can be found, for example, in Marden \cite{marden}.

%%%%%%%%%%%%%%%%%%%%%%%%%%%%%%%%%%%%%%%%%%%%%%%%%%%%%%%%%%%%%%%%%
\section{Elementary groups}

\begin{definition}\label{Def:Elementary}
A subgroup $G\leq\PSL(2,\CC)$ is \emph{elementary}\index{elementary group} if one of the following holds.
\begin{enumerate}
\item The union of all fixed points on $\bdy\HH^3$ of all nontrivial elements of $G$ is a single point on $\bdy\HH^3$.
\item The union of all fixed points on $\bdy\HH^3$ of all nontrivial elements of $G$ consists of exactly two points on $\bdy\HH^3$. 
\item\label{Itm:FixedPtInt} There exists $x\in \HH^3$ such that for all $g\in G$, $g(x)=x$.
\end{enumerate}
The group is \emph{nonelementary}\index{nonelementary group} if it is not elementary. 
\end{definition}

Elementary groups will be important subgroups of the discrete groups we study. Because of that, we will need to know more about their form.

\begin{proposition}\label{Prop:ClassInfElementary}
Let $G$ be a discrete nontrivial elementary subgroup\index{elementary group} of $\PSL(2,\CC)$ without elliptics. Then either
\begin{enumerate}
\item the union of fixed points of nontrivial elements of $G$ is a single point on $\bdy\HH^3$, $G$ is isomorphic to $\ZZ$ or $\ZZ\times\ZZ$, and $G$ is generated by parabolics\index{parabolic} (fixing the same point on $\bdy\HH^3$), or
\item the union of fixed points of nontrivial elements of $G$ consists of two points on $\bdy \HH^3$, $G$ is isomorphic to $\ZZ$, and $G$ is generated by a single loxodromic\index{loxodromic} leaving invariant the line between the fixed points.
\end{enumerate}
\end{proposition}

\begin{proof}
If the union of all fixed points of nontrivial elements of $G$ consists of a single point on $\bdy\HH^3$, then $G$ must contain only parabolics\index{parabolic} fixing that point. Conjugate so that the fixed point is $\infty$ in $\HH^3$. Then $G$ fixes a horosphere about $\infty$, which is isometric to the Euclidean plane $P$. The group $G$ acts on $P$ by Euclidean translations. Since $G$ is discrete, $G$ must be generated by either one translation, in which case $G\cong\ZZ$, or two linearly independent translations, in which case $G\cong\ZZ\times\ZZ$.

If the union of all fixed points of nontrivial elements of $G$ consists of two points, then $G$ contains only loxodromics\index{loxodromic} fixing the axis between them. The group $G$ acts on the axis; the fact that the group is discrete means that there is some finite minimal translation distance $\tau$ under this group action. Let $A\in G$ realize the minimal translation distance, i.e.\ $d(x,Ax)=\tau$ for $x$ on the axis. We claim $G = \langle A \rangle$.
First, we show all $C\in G$ translate by distance $n\tau$ for some $n\in\ZZ$, for if some $C\in G$ has translation distance that is not a multiple of $\tau$, then $C(x)$ lies between $A^n(x)$ and $A^{n+1}(x)$ for any $x$ on the axis. But then $CA^{-n}\in G$ translates $A^n(x)$ a distance strictly less than $\tau$, which is a contradiction. Thus all $C\in G$ translate along the axis a distance equal to a multiple of $\tau$. Now suppose $C\in G$ translates by $n\tau$ for some integer $n$. Then $CA^{-n}$ fixes the axis pointwise. Because $G$ contains no elliptics, $C=A^n$. So $G$ is cyclic generated by $A$. 
\end{proof}

Consider the first case of \refprop{ClassInfElementary}.

\begin{figure}
  \includegraphics{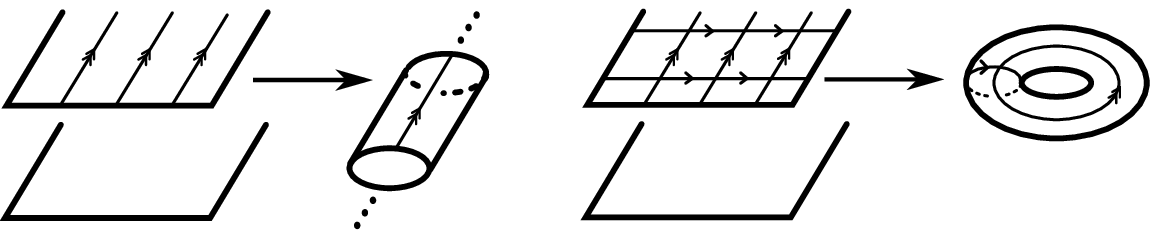}
  \caption{Left: The quotient of horosphere $\bdy C$ under the group $\ZZ$ generated by a single parabolic\index{parabolic} gives a cylinder, or annulus. Right: The quotient of $\bdy C$ under $\ZZ\times\ZZ$ is a torus}
  \label{Fig:Rank1Rank2Cusp}
\end{figure}

\begin{definition}\label{Def:RankOneRankTwoCusp}
Suppose $G$ is an infinite elementary\index{elementary group} discrete group in $\PSL(2,\CC)$ fixing a single point on $\bdy\HH^3$. We may conjugate $G$ so that fixed point is the point at infinity. Let $H$ be the closed horoball of height $1$:
\[ H = \{ (x,y,z) \mid z\geq 1\}. \]
\Refprop{ClassInfElementary} tells us that $G$ is isomorphic to $\ZZ$ or $\ZZ\times\ZZ$.

If $G\cong\ZZ$, the quotient of the horoball $H/G$ is homeomorphic to the space $A\times [1,\infty)$, where $A$ is an annulus, or cylinder; see \reffig{Rank1Rank2Cusp}. We say that $H/G$ is a \emph{rank-1 cusp}\index{rank-1 cusp}\index{cusp!rank-1}.

If $G\cong\ZZ\times\ZZ$, the quotient of the horoball $H/G$ is homeomorphic to $T\times[1,\infty)$, where $T$ is a Euclidean torus; see \reffig{Rank1Rank2Cusp}, right. We say that $H/G$ is a \emph{rank-2 cusp}\index{rank-2 cusp}\index{cusp!rank-2}.
\end{definition}

\Refprop{ClassInfElementary} has an immediate corollary giving information about $\ZZ\times\ZZ$ subgroups of discrete groups, which will be used in later chapters. 

\begin{corollary}[$\ZZ\times\ZZ$ subgroups]\label{Cor:ZxZSubgroup}
Suppose a discrete group $G$ without elliptics has a subgroup isomorphic to $\ZZ\times\ZZ$. Then the subgroup is generated by two parabolic\index{parabolic} elements fixing the same point on the boundary at infinity\index{boundary at infinity} $\bdy\HH^3$ of $\HH^3$.
\end{corollary}

\begin{proof}
Let $A$ and $B$ denote the generators of the subgroup of $G$ isomorphic to $\ZZ\times\ZZ$. Since $A$ and $B$ commute, they must have the same fixed points on the boundary at infinity\index{boundary at infinity} $\bdy\HH^3$ of $\HH^3$ (\refex{PSL(2,C)Commute}). Thus $H\cong \langle A,B\rangle$ is an elementary discrete group\index{elementary group} isomorphic to $\ZZ\times\ZZ$. By \refprop{ClassInfElementary}, $H$ must be generated by parabolics\index{parabolic} fixing the same point on $\bdy \HH^3$. 
\end{proof}

%%%%%%%%%%%%%%%%%%%%%%%%%%%%%%%%%%%%%%%%%%%%%%%%%%%%%%%%%%%%%%%%%
Discrete elementary groups\index{elementary group} are often defined in terms of the set of accumulation points of the group on $\bdy \HH^3$; for example this is the definition in \cite{thurston}. We review that definition here as well. 

\begin{definition}\label{Def:LimitSet}
Let $G\leq\PSL(2,\CC)$ be a discrete group, and let $x\in \HH^3$ be any point. The \emph{limit set $\Lambda(G)$}\index{limit set} is defined to be the set of accumulation points on $\bdy\HH^3$ of the orbit $G(x)$.
\end{definition}

\begin{lemma}\label{Lem:LimitSetIndX}
The limit set $\Lambda(G)$ is well-defined, independent of choice of $x$ in \refdef{LimitSet}.
\end{lemma}

\begin{proof}
Suppose $\{A_n\}\subset G$ is a sequence such that $A_n(x)$ converges to a point $p\in \Lambda(G)\subset \bdy \HH^3$. Let $y\in \HH^3$. Then the distance between $x$ and $y$ is a constant, equal to the distance between $A_n(x)$ and $A_n(y)$ for all $n$. Thus as $n\to\infty$, $A_n(x)$ and $A_n(y)$ lie a bounded distance apart, but $A_n(x)$ approaches $p$. This is possible only if $A_n(y)$ approaches the same point $p$ on $\bdy\HH^3$. 
\end{proof}

Consider a few examples of groups $G$ and limit sets $\Lambda(G)$. If $G$ is generated by a single loxodromic\index{loxodromic} element $g$, then its limit set $\Lambda(G)$ consists of the two fixed points of $g$ on $\bdy\HH^3$: one is an accumulation point for $g^n(x)$, and the other for $g^{-n}(x)$. If $G$ is generated by a single parabolic\index{parabolic} element, then $\Lambda(G)$ consists of a single point. If $G$ contains both a loxodromic element $g$ and a parabolic element $h$, then $\Lambda(G)$ contains the fixed points of $g$ on $\bdy\HH^3$, as well as the fixed points of $h^n\circ g$ for all $n$; this is a countably infinite set. Finally, if $G$ is the identity group, consisting only of the identity element, then $\Lambda(G)$ is empty. 

The following is often given as the definition of an elementary discrete subgroup of $\PSL(2,\CC)$. 

\begin{lemma}\label{Lem:Elementary2}
A discrete subgroup $G\leq\PSL(2,\CC)$ with no elliptics is elementary\index{elementary group} if and only if $\Lambda(G)$ consists of $0$, $1$, or $2$ points. \qed
\end{lemma}

%%%%%%%%%%%%%%%%%%%%%%%%%%%%%%%%%%%%%%%%%%%%%%%%%%%%%%%%%%%%%%%%%
%%%%%%%%%%%%%%%%%%%%%%%%%%%%%%%%%%%%%%%%%%%%%%%%%%%%%%%%%%%%%%%%%
\section{Thick and thin parts}

We are now ready to put together facts about elementary\index{elementary group} and nonelementary\index{nonelementary group} discrete groups to prove a remarkable result on the geometry and topology of hyperbolic 3-manifolds, namely that any such manifold decomposes into a thick part and completely classified thin parts. To state the result precisely, we give a few definitions. 

\begin{definition}\label{Def:InjRad}
Suppose $M$ is a complete hyperbolic 3-manifold and $x\in M$. The \emph{injectivity radius}\index{injectivity radius} of $x$, denoted $\injrad(x)$, is defined to be the supremal radius $r$ such that a metric $r$-ball around $x$ is embedded.
\end{definition}

\begin{definition}\label{Def:ThickThin}
Let $M$ be a complete hyperbolic 3-manifold, and let $\epsilon>0$. Define the \emph{$\epsilon$-thin part}\index{thin part} of $M$, denoted $M^{<\epsilon}$ to be
\[ M^{<\epsilon} = \{ x\in M \mid \injrad(x)<\epsilon/2 \}. \]
Similarly, the \emph{$\epsilon$-thick part}\index{thick part}, denoted $M^{>\epsilon}$ is defined to be
\[ M^{>\epsilon} = \{ x\in M \mid \injrad(x)>\epsilon/2 \}. \]
We also have closed versions $M^{\geq\epsilon}$ and $M^{\leq\epsilon}$ defined in the obvious way.
\end{definition}

\begin{theorem}[Structure of thin part]\label{Thm:ThinPart}\index{thin part!structure of thin part}
There exists a universal constant $\epsilon_3>0$ such that for $0<\epsilon \leq \epsilon_3$, the $\epsilon$-thin part of any complete, orientable, hyperbolic 3-manifold $M$ consists of tubes around short geodesics, rank-1 cusps, and/or rank-2 cusps. 
\end{theorem}

\begin{figure}
%% Accompanies image file 'F5-03-Thin.eps' (pdf, eps, ps)
\begingroup%
  \makeatletter%
  \providecommand\color[2][]{%
    \errmessage{(Inkscape) Color is used for the text in Inkscape, but the package 'color.sty' is not loaded}%
    \renewcommand\color[2][]{}%
  }%
  \providecommand\transparent[1]{%
    \errmessage{(Inkscape) Transparency is used (non-zero) for the text in Inkscape, but the package 'transparent.sty' is not loaded}%
    \renewcommand\transparent[1]{}%
  }%
  \providecommand\rotatebox[2]{#2}%
  \newcommand*\fsize{\dimexpr\f@size pt\relax}%
  \newcommand*\lineheight[1]{\fontsize{\fsize}{#1\fsize}\selectfont}%
  \ifx\svgwidth\undefined%
    \setlength{\unitlength}{173.42262268bp}%
    \ifx\svgscale\undefined%
      \relax%
    \else%
      \setlength{\unitlength}{\unitlength * \real{\svgscale}}%
    \fi%
  \else%
    \setlength{\unitlength}{\svgwidth}%
  \fi%
  \global\let\svgwidth\undefined%
  \global\let\svgscale\undefined%
  \makeatother%
  \begin{picture}(1,0.81608809)%
    \lineheight{1}%
    \setlength\tabcolsep{0pt}%
    \put(0,0){\includegraphics[width=\unitlength]{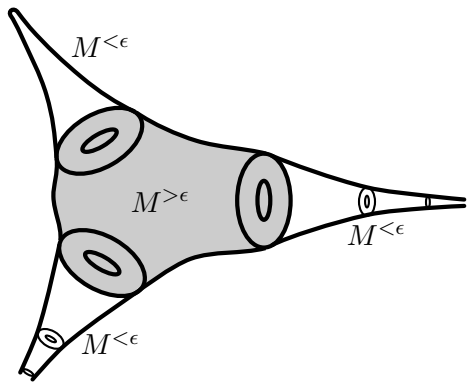}}%
    \put(0.26677595,0.3716222){\color[rgb]{0,0,0}\makebox(0,0)[lt]{\lineheight{0}\smash{\begin{tabular}[t]{l}$M^{>\epsilon}$\end{tabular}}}}%
    \put(0.15843459,0.06290793){\color[rgb]{0,0,0}\makebox(0,0)[lt]{\lineheight{0}\smash{\begin{tabular}[t]{l}$M^{<\epsilon}$\end{tabular}}}}%
    \put(0.13639371,0.70821892){\color[rgb]{0,0,0}\makebox(0,0)[lt]{\lineheight{0}\smash{\begin{tabular}[t]{l}$M^{<\epsilon}$\end{tabular}}}}%
    \put(0.7363484,0.30048265){\color[rgb]{0,0,0}\makebox(0,0)[lt]{\lineheight{0}\smash{\begin{tabular}[t]{l}$M^{<\epsilon}$\end{tabular}}}}%
  \end{picture}%
\endgroup%

\caption{A schematic picture of a hyperbolic 3-manifold $M$, with $M^{<\epsilon}$ a collection of cusps and tubes}
\label{Fig:SchematicMargulis}
\end{figure}

A cartoon illustrating \refthm{ThinPart} is given in \reffig{SchematicMargulis}. 

\begin{definition}\label{Def:MargulisConstant}
The supremum of all constants $\epsilon_3$ satisfying \refthm{ThinPart} is called the \emph{Margulis constant}\index{Margulis constant}. More generally, given a complete hyperbolic 3-manifold $M$, a number $\epsilon>0$ is said to be a \emph{Margulis number for $M$}\index{Margulis number for $M$} if $M^{<\epsilon}$ satisfies the conclusions of \refthm{ThinPart}, i.e.\ $M^{<\epsilon}$ consists of tubes around short geodesics, rank-1, and/or rank-2 cusps. The Margulis constant is therefore the infimum over all complete hyperbolic 3-manifolds $M$ of the supremum of all Margulis numbers for $M$.
\end{definition}

As of the writing of this book, the optimal Margulis constant is still unknown, although there are bounds on its value.
R.~Meyerhoff gave what is currently the best lower bound on $\epsilon_3$ in \cite[Section~9]{meyerhoff}, that it is at least $0.104$. As for an upper bound, M.~Culler has discovered a closed hyperbolic 3-manifold with Margulis number less than $0.616$ using SnapPea \cite{weeks:computation}.

We save the proof of \refthm{ThinPart} until the end of this section. We will see that it is a consequence of a well-known theorem concerning the structure of discrete groups of isometries, commonly called the Margulis lemma,\index{Margulis lemma} which appears in a paper of Ka{\v{z}}dan and Margulis \cite{KazhdanMargulis}. The actual Margulis lemma is very general, concerning discrete groups acting on symmetric spaces. We restrict to the case of discrete subgroups of $\PSL(2,\CC)$ acting freely\index{free group action}\index{group action!free} on hyperbolic space. The consequence we will need is the following. 

\begin{theorem}[Universal Elementary Neighborhoods]\label{Thm:Margulis}
There is a universal constant $\epsilon_3>0$ such that for all $x\in\HH^3$, and for any discrete group $G\leq \PSL(2,\CC)$ without elliptics, if $H$ denotes the subgroup of $G$ generated by all elements of $G$ that translate $x$ distance less than $\epsilon_3$, then $H$ is elementary.\index{elementary group}
\end{theorem}

We will give a proof of \refthm{Margulis} in \refsec{Margulis}. Before that, a few remarks are in order. First, the Margulis lemma holds when we allow elliptics;\index{elliptic} this appears in Wang \cite{Wang1969} in the full generality of the theorem of Ka{\v{z}}dan and Margulis. Second, the form of \refthm{Margulis} above, concerning discrete subgroups of $\PSL(2,\CC)$, is due to J{\o}rgensen and Marden, only their result is more general in that it also includes elliptics. Their proof appears in \cite{marden}, and is the basis for the proof that we include below in \refsubsec{Technical}.

However, before we discuss the proof, we show how \refthm{Margulis} implies \refthm{ThinPart} (Structure of thin part). 
First, we need to relate translation distance to injectivity radius.\index{injectivity radius} 

\begin{lemma}\label{Lem:InjRad}
Let $M$ be a complete, orientable, hyperbolic 3-manifold with $M\cong \HH^3/\Gamma$ for a discrete group $\Gamma\leq\PSL(2,\CC)$. For any $x\in M$ with lift $\tilde{x}\in \HH^3$,
\[ \injrad(x) = \half\inf_{A\neq\Id \in \Gamma} \{ d(\tilde{x}, A\tilde{x})\}. \]
Moreover, this is realized. That is, there exists nontrivial $A\in \Gamma$ such that $2\,\injrad(x) = d(\tilde{x},A\tilde{x})$. 
\end{lemma}

\begin{proof}
A metric $r$-ball is embedded at $x$ if and only if for all $A\neq \Id\in \Gamma$, the metric $r$-ball $B(r,\tilde{x})$ is disjoint from the metric $r$-ball $A(B(r,\tilde{x})) = B(r,A\tilde{x})$. This holds if and only if the translation distance $d(\tilde{x}, A\tilde{x})$ is at least $2r$ for all $A$.

Now suppose $\injrad(x)=b$. Then a metric $b$-ball is embedded, but for any $\epsilon>0$, a metric $b+\epsilon$-ball is not embedded. Thus for each $\epsilon>0$, there is $A_\epsilon\in\Gamma$ such that $d(\tilde{x},A_\epsilon(\tilde{x})<2(b+\epsilon)$. If the set $\{A_\epsilon\}$ contains infinitely many distinct elements, then we obtain a sequence $\{A_n\}$ such that $A_n(\tilde{x})$ is of bounded distance from $\tilde{x}$. By \reflem{Sequences}, $A_n\to A\in \PSL(2,\CC)$, implying $\Gamma$ is not discrete by \reflem{EquivDiscrete}. This is a contradiction. Thus $\{A_\epsilon\}$ is a finite set. Let $A\in\Gamma$ be such that $d(\tilde{x},A\tilde{x})$ is minimal. This $A$ satisfies the conclusion of the lemma. 
\end{proof}

We are now ready to complete the proof of \refthm{ThinPart}, assuming \refthm{Margulis}.

\begin{proof}[Proof of \refthm{ThinPart} (Structure of thin part)]
Take $\epsilon_3>0$ as in \refthm{Margulis}. Let $M \cong \HH^3/\Gamma$ be a complete, orientable, hyperbolic 3-manifold, so $\Gamma\leq \PSL(2,\CC)$ is a discrete subgroup with no elliptics.

For $\epsilon\leq\epsilon_3$, if $x\in M^{<\epsilon}$, then by definition $\injrad(x)<\epsilon/2$.\index{injectivity radius} By \reflem{InjRad}, it follows that there exists $A\neq\Id\in \Gamma$ such that $d(\tilde{x},A\tilde{x}) < \epsilon$ for any lift $\tilde{x}$ of $x$. But \refthm{Margulis} implies that the subgroup $\Gamma_\epsilon$ of $\Gamma$ generated by all $A \in \Gamma$ such that $d(\tilde{x}, A\tilde{x})<\epsilon$ is elementary.\index{elementary group} Since $\Gamma_\epsilon$ contains $A\neq\Id$, \refprop{ClassInfElementary} implies that $\Gamma_\epsilon$ either fixes a single point $\zeta\in\bdy\HH^3$ and is generated by parabolics\index{parabolic} fixing $\zeta$, or $\Gamma_\epsilon$ is generated by a single loxodromic\index{loxodromic} preserving an axis $\ell \subset \HH^3$.

Suppose first that $\Gamma_\epsilon$ fixes a single point $\zeta\in\bdy\HH^3$. Then $\Gamma_\epsilon$ is generated by one or two parabolics\index{parabolic} (\refprop{ClassInfElementary}), and $\tilde{x}$ lies on a horosphere $H$ about $\zeta$ that is fixed by $\Gamma_\epsilon$. 
Suppose $\tilde{y}$ lies in the horoball bounded by $H$. Then the height of $\tilde{y}$ is at least $C$: $\tilde{y}$ has coordinates $(a+b\,i,t)$ with $t\geq C$. A generator $A$ of $\Gamma_\epsilon$ takes $\tilde{y}$ to a point with the same height $t$. A calculation in this case (\refex{CalculateHypDist}) shows that
\[ \epsilon > d(\tilde{x},A\tilde{x}) \geq d(\tilde{y},A\tilde{y}),\]
and it follows that in the quotient $\HH^3/\Gamma$, the point $\tilde{y}$ maps to $M^{<\epsilon}$. Since this is true for every point in the horoball bounded by $H$, $M^{<\epsilon}$ contains the quotient of a horoball under the elementary group\index{elementary group} $\Gamma_\epsilon$; this is a rank-1 or rank-2 cusp. 

Now suppose that $H$ is generated by a single loxodromic\index{loxodromic} $A$ preserving the axis $\ell$. Let $R$ denote the distance from $\tilde{x}$ to the axis $\ell$, and let $T_R$ denote the set of points in $\HH^3$ of distance $R$ from the axis $\ell$. Then $T_R$ bounds a tube consisting of all points in $\HH^3$ of distance at most $R$ from $\ell$. If $\tilde{y}$ is any point within this tube, then one can calculate (\refex{CalculateDistTube}) that $d(\tilde{y},A\tilde{y}) \leq d(\tilde{x},A\tilde{x})<\epsilon$, so $M^{<\epsilon}$ contains the quotient of a tube about $\ell$ under the elementary group\index{elementary group} $\langle A\rangle$. This is a tube around a short geodesic.
\end{proof}

%%%%%%%%%%%%%%%%%%%%%%%%%%%%%%%%%%%%%%%%%%%%%%%%%%%%%%%%%%%%%%%%%
%%%%%%%%%%%%%%%%%%%%%%%%%%%%%%%%%%%%%%%%%%%%%%%%%%%%%%%%%%%%%%%%%
\section{Hyperbolic manifolds with finite volume}

In \refchap{GluingCompleteness} we gave a method that will allow us to compute (complete) hyperbolic structures on many 3-manifolds, including many knot complements. Once we have a hyperbolic structure on a 3-manifold, we have equipped the manifold with a Riemannian metric with very nice properties, for example the metric can be described in local coordinates by \refeqn{3dHypMetric}.

One of the simplest invariants we can compute from a hyperbolic metric is the volume of the underlying manifold. This gives a good measure of the ``size'' of the manifold. In \refchap{Volume}, we will discuss volumes in some detail, including how to compute volumes of hyperbolic 3-manifolds including knot and link complements. Meanwhile, we give an application of the thick--thin decomposition of hyperbolic 3-manifolds to classifying those with finite volume.

\begin{theorem}\label{Thm:FteVolIffTorusBdy}
A hyperbolic 3-manifold $M$ has finite volume if and only if $M$ is closed (compact without boundary), or $M$ is homeomorphic to the interior of a compact manifold $\overline{M}$ with torus boundary components.
\end{theorem}

\begin{proof}
If $M$ is closed then a fundamental domain for $M$ in its universal cover $\HH^3$ is a compact set, hence has finite volume. If $M$ is the interior of a manifold with torus boundary, then each such boundary component will be realized as a cusp in the complete hyperbolic structure on $M$. The complement of the cusps of $M$ in $M$ is compact, hence has finite volume. We now show that each cusp has finite volume.

Consider the universal cover $\HH^3$. For any cusp $C$, we may apply an isometry to $\HH^3$ so that the point at infinity projects to that cusp, and a horoball of height $1$ projects to an embedded horoball neighborhood of the cusp. On the horosphere of height $1$, some parallelogram $A$ will be a fundamental region for the torus of the cusp, since the structure is complete (\refthm{EuclidCusp}). Then the volume of the cusp is given by
\[ \int_C d\vol = \int_{t=1}^{\infty}\int_A d\vol = \int_{t=1}^\infty\int_A \frac{dx\,dy\,dt}{t^3} = \half\area(A). \]
(See \refex{CuspVolume}.) Thus every cusp has finite volume. Since the volume of $M$ is the sum of the volumes of the compact region with cusps removed, as well as a finite number of finite-volume cusps, the manifold $M$ has finite volume.

To prove the converse, we use \refthm{ThinPart}. Suppose $M$ is a complete hyperbolic manifold with finite volume. Fix $\epsilon>0$ less than the universal constant $\epsilon_3$ of \refthm{ThinPart}, and consider $M^{<\epsilon}$ and $M^{\geq\epsilon}$. By \refthm{ThinPart}, $M^{<\epsilon}$ consists of cusps and tubes. Note that a rank-1 cusp has infinite volume, hence since $M$ has finite volume, $M^{<\epsilon}$ consists of rank-2 cusps and tubes, each of which has finite volume. On the other hand, $M^{\geq\epsilon}$ has finite volume. Moreover, any point in $M^{\geq\epsilon}$ is contained in an embedded ball of radius at least $\half\epsilon$. If two points in $M^{\geq\epsilon}$ have distance at least $\epsilon$, then the balls of radius $\half\epsilon$ about each are disjointly embedded in $M^{\geq\epsilon}$. Thus a collection of points with pairwise distance at least $\epsilon$ in $M^{\geq\epsilon}$ leads to a pairwise disjoint collection of $\epsilon/2$-balls. Because $M$ has finite volume there can only be finitely many of these. Starting with any such collection of points, we may complete the collection to a maximal collection of points of $M^{\geq\epsilon}$ of distance at least $\epsilon$; there are finitely many of these and the $\epsilon/2$-balls around each are embedded. Then the closed $\epsilon$-balls about the collection must contain $M^{\geq\epsilon}$. The union of these balls is a compact set, and $M^{\geq\epsilon}$ is a closed subset. Hence $M^{\geq\epsilon}$ is compact.

Now the union of $M^{\geq\epsilon}$ and any tubes of $M^{<\epsilon}$ is the union of compact sets, hence compact. This is a manifold with boundary homeomorphic to a finite collection of tori corresponding to the finite number of cusps of $M^{<\epsilon}$. Attach a closed collar neighborhood of each torus boundary component, and call the result $N$; each collar neighborhood is homeomorphic to $T^2\times[0,1]$, where $T^2$ is a torus. Then by construction, the manifold $M$ is homeomorphic to the interior of $N$. 
\end{proof}

By \refthm{FteVolIffTorusBdy}, the complement of any knot or link in $S^3$ with a hyperbolic structure must have finite hyperbolic volume.

%%%%%%%%%%%%%%%%%%%%%%%%%%%%%%%%%%%%%%%%%%%%%%%%%%%%%%%%%%%%%%%%%
%%%%%%%%%%%%%%%%%%%%%%%%%%%%%%%%%%%%%%%%%%%%%%%%%%%%%%%%%%%%%%%%%
\section{Universal elementary neighborhoods}\label{Sec:Margulis}

In this section, we give a proof of \refthm{Margulis}, on the existence of universal elementary neighborhoods.

In fact, we split this section into two subsections. The first gives a proof of \refthm{Margulis} that is elementary, in the sense that it uses only the machinery of subgroups of $\PSL(2,\CC)$ developed in this chapter. However, it is also quite technical, requiring calculations that, upon first glance, may seem mysterious and arbitrary. Nevertheless, by the end of \refsubsec{Technical}, the proof of \refthm{Margulis} is complete. 

The second subsection is an attempt to put \refthm{Margulis} into a wider mathematical context. Although we have presented a proof that uses only the tools of $\PSL(2,\CC)$, related theorems hold for much more general Lie groups. The technical calculations of \refsubsec{Technical} can be seen as instances of more general, and in some sense simpler, mathematical phenomena, put into a broader context. 

%%%%%%%%%%%%%%%%%%%%%%%%%%%%%%%%%%%%%%%%%%%%%%%%%%%%%%%%%%%%%%%%%
\subsection{A technical proof in $\PSL(2,\CC)$}\label{Subsec:Technical}

We now give a complete proof of \refthm{Margulis}, restricting to the setting of $\PSL(2,\CC)$.

We need a few more tools before we begin. Namely, \refprop{ClassInfElementary} classifies elementary discrete groups without elliptics. We also need the following result giving more information on \emph{nonelementary} discrete groups without elliptics.

\begin{lemma}\label{Lem:InfiniteElementary}
If $G$ is a nonelementary discrete subgroup\index{nonelementary group} of $\PSL(2,\CC)$ that contains no elliptics, then the following hold.
\begin{enumerate}
\item\label{Itm:Gfinite} $G$ is infinite.
\item\label{Itm:LoxNoFixed} For any nontrivial $A\in G$, there exists a loxodromic\index{loxodromic} $B\in G$ that has no common fixed points with $A$.
\item\label{Itm:LoxOneCommonIndiscrete} If $B\in G$ is loxodromic, then there is no nontrivial $C\in G$ that has exactly one fixed point in common with $B$. 
\item\label{Itm:2Lox} $G$ contains two loxodromic elements with no fixed points in common. 
\end{enumerate}
\end{lemma}

\begin{proof}
The group $G$ must be nontrivial; since it contains no elliptics it must contain a loxodromic\index{loxodromic} or parabolic.\index{parabolic} Such an element has infinite order, so $G$ is infinite, proving \refitm{Gfinite}.

Next we show \refitm{LoxOneCommonIndiscrete}. Suppose $B$ is loxodromic,\index{loxodromic} and $C$ has exactly one fixed point in common with $B$; we will show that the group generated by $B$ and $C$ is indiscrete, contradicting the fact that $G$ is discrete. Conjugate the group. \Reflem{MoreClassifyPSL} implies we may assume $B=\mat{\rho&0\\0&1/\rho}$, and since $C$ has exactly one fixed point in common with $B$, it has the form $C=\mat{a&b\\0&1/a}$ where $b\neq 0$. Then
\[ B^nCB^{-n}C^{-1} = \mat{1& ab\;(\rho^{2n}-1) \\ 0& 1}. \]
If $|\rho|<1$, let $n\to\infty$. If $|\rho|>1$, let $n\to-\infty$. In either case, $B^nCB^{-n}C^{-1}$ approaches the parabolic\index{parabolic} $\mat{1&-ab\\0&1}$. \Reflem{EquivDiscrete} now implies that the subgroup generated by $B$ and $C$ is not discrete, therefore $G$ is not discrete. 

Now we show \refitm{LoxNoFixed}. There are two cases depending on whether $A$ is parabolic\index{parabolic} or loxodromic.\index{loxodromic} Note that if we can show the result for a conjugate group $UGU^{-1}$ for $U\in\PSL(2,\CC)$, then the result holds for $G$, so in both cases we will replace $G$ by a conjugate group at the first step.

\textbf{Case 1.} Suppose $A$ is parabolic.\index{parabolic} Then by \reflem{MoreClassifyPSL}, $A$ is conjugate to $z\mapsto z+1$, so we may assume $A=\mat{1&1\\0&1}$ and $A$ fixes $\infty$. Because $G$ is nonelementary, there exists $C\in G$ that does not fix $\infty$. If $C$ is loxodromic,\index{loxodromic} we are done. If not, $C$ must be parabolic,\index{parabolic} and $C=\mat{a&b\\c&d}$ with $c\neq 0$. Note that $A^nC$ cannot fix $\infty$ for any integer $n$, and $\tr(A^nC) = a+nc+d = nc\pm 2$. For $|n|$ sufficiently large, this cannot be in $[-2,2]$, so $A^nC$ is the desired loxodromic by \reflem{MoreClassifyPSL}.

\textbf{Case 2.} Suppose $A$ is loxodromic.\index{loxodromic} Then after conjugating, \reflem{MoreClassifyPSL} implies we may assume $A=\mat{\rho&0\\0&\rho^{-1}}$ with $|\rho|>1$, so $A$ fixes $0$ and $\infty$. Because $G$ is nonelementary and discrete, \refitm{LoxOneCommonIndiscrete} implies there is $C=\mat{a&b\\c&d} \in G$ that does not fix either $0$ or $\infty$ (so $b, c\neq 0$). If $C$ happens to be loxodromic,\index{loxodromic} we are done. If not, $C$ is parabolic,\index{parabolic} so $a+d=\pm 2$. Then $A^nC$ also has distinct fixed points from those of $A$ for any integer $n$, and $\tr(A^nC)=a\rho^n+d\rho^{-n}$. For $|n|$ large, this lies outside $[-2,2]$, hence $A^nC$ is loxodromic\index{loxodromic} by \reflem{MoreClassifyPSL}. This concludes the proof of \refitm{LoxNoFixed}. 

Finally, to prove part \refitm{2Lox}, we use part \refitm{LoxNoFixed}. Suppose $A\in G$ is not the identity. Then \refitm{LoxNoFixed} implies there is a loxodromic\index{loxodromic} $B\in G$ with distinct fixed points from $A$. If $A$ is also loxodromic, we are done. Otherwise, apply \refitm{LoxNoFixed} to $B$, to obtain a loxodromic $C$ with no fixed points in common with $B$. Then $B$ and $C$ are the desired loxodromics. 
\end{proof}

The following theorem, on convergence of nonelementary discrete groups, is due to J{\o}rgensen and Klein \cite{jorgensen-klein}, using previous work of J{\o}rgensen \cite{jorgensen}.

\begin{theorem}[J{\o}rgensen and Klein, 1982]\label{Thm:JorgensenKlein}
  Let
  \[ G_n= \langle A_{1,n}, A_{2,n}, \dots, A_{r,n}\rangle\] be a sequence of $r$-generator, nonelementary,\index{nonelementary group} discrete subgroups of $\PSL(2,\CC)$ such that $A_k = \lim_{n\to\infty} A_{k,n}$ exists and is an element of $\PSL(2,\CC)$ for each $k$. Then $G=\langle A_1,A_2, \dots, A_r\rangle$ is also nonelementary and discrete. Moreover, for sufficiently large $n$, the map $A_k\to A_{k,n}$ for each $k$ extends to a homomorphism from $G$ to $G_n$.
\end{theorem}

The proof of \refthm{JorgensenKlein} follows from an analysis of various properties of elements of $\PSL(2,\CC)$ and discrete subgroups. Its proof is not unlike many of the other results proved in this chapter. However, its proof would lead us a little further afield than we wish to go, into technicalities of $\PSL(2,\CC)$. The full proof can be found in the original papers; Marden also gives an exposition closely following the original proof in \cite{marden}. We will refer the interested reader to those references.

Meanwhile, we don't actually need the full strength of \refthm{JorgensenKlein}; we only need the following immediate consequence. 

\begin{corollary}\label{Cor:LimitNonelementary}
Suppose $\{\langle A_n,B_n\rangle\}$ is a sequence of nonelementary\index{nonelementary group} discrete subgroups of $\PSL(2,\CC)$ such that $\lim A_n=A$ and $\lim B_n=B$ in $\PSL(2,\CC)$. Then $\langle A, B\rangle$ is a nonelementary discrete subgroup of $\PSL(2,\CC)$.\qed
\end{corollary}

\begin{proof}[Proof of \refthm{Margulis}]
First we establish some notation. For fixed $x\in\HH^3$ and $A\in\PSL(2,\CC)$, let $d(x,Ax)$ denote the distance in $\HH^3$ between $x$ and $Ax$. For fixed $r>0$, let $G(r,x)$ denote the set
\[ G(r,x) = \{A\in G \mid d(x, Ax) < r\}.
\]
The group generated by $G(r,x)$ will be denoted by $\langle G(r,x) \rangle$.

Our goal is to show that there exists $r>0$ such that for all discrete $G$ and for all $x$, the group $\langle G(r,x)\rangle$ is elementary.\index{elementary}

As a first step, we show that if we fix a discrete group $G$ with no elliptics and fix $x$, then there exists $r>0$ such that the group $\langle G(r,x)\rangle$ is elementary. 
For suppose this is not the case. Then for a sequence $r_n\to 0$, each $\langle G(r_n,x)\rangle$ is nonelementary. It follows that there exists a sequence of distinct $A_n\in G(r_n,x)$ with $d(x,A_nx)<r_n$. But then \reflem{Sequences} implies that $A_n$ must converge to some $A\in\PSL(2,\CC)$. Using \reflem{EquivDiscrete}, we see that this contradicts the fact that $G$ is a discrete group. So for $r>0$ sufficiently small, $\langle G(r,x)\rangle$ is elementary, and it follows that $G(r,x)$ contains finitely many elements. By choosing $r>0$ smaller than the translation distance of each of these elements, we find that $G(r,x)$ contains only the identity element. Note that the identity group is elementary. 

Now we will prove the more general result, that there is a universal $r>0$, independent of $G$ and $x$, such that $\langle G(r,x)\rangle$ is always elementary. Again suppose not. Then there is a sequence $r_n\to 0$, a sequence of discrete groups $G_n\leq\PSL(2,\CC)$ without elliptics, and a sequence of points $x_n\in\HH^3$ such that $\langle G_n(r_n,x_n)\rangle$ is not elementary.

We will simplify the argument by replacing $x_n$ with a fixed $x$ for all $n$: choose any $x\in \HH^3$, and let $R_n \in \PSL(2,\CC)$ be an isometry mapping $x_n$ to $x$. Consider the group $R_n G_n R_n^{-1}$. Note that $A\in G_n(r_n,x_n)$ if and only if $R_n A R_n^{-1}$ is in $R_nG_nR_n^{-1}(r_n,x)$, and so $\langle R_nG_nR_n^{-1}(r_n,x)\rangle$ is nonelementary. Thus if we replace $G_n$ by $R_nG_nR_n^{-1}$, we may work with a single fixed value of $x$. So we assume there is a fixed $x$ and sequences $r_n\to 0$ and $G_n$ so that $\langle G_n(r_n,x)\rangle$ is nonelementary.

Now fix $n$. Our next goal is to find $A_n$ and $B_n$ in $G_n(r_n,x)$ such that $\langle A_n,B_n\rangle$ is nonelementary. Since $\langle G_n(r_n,x)\rangle$ is nonelementary, \reflem{InfiniteElementary} implies that there exist loxodromics\index{loxodromic} $S_n$ and $T_n$ with no common fixed points in $\langle G_n(r_n,x)\rangle$, and certainly they generate a nonelementary group. However, we need to take some care to ensure that $A_n$ and $B_n$ are actually in $G_n(r_n,x)$. To do this, we use the first part of this proof: consider the groups $\langle G_n(\rho,x)\rangle$ as $\rho$ ranges between $0$ and $r_n$. We have observed that for some $\rho_n<r_n$, the group $\langle G_n(\rho_n,x)\rangle$ will consist only of the identity element. As $\rho$ increases, the sets $G_n(\rho,x)$ will be nested. There will be some value $0<\mu_n\leq r_n$ such that $\langle G_n(\rho,x)\rangle$ is elementary for $\rho<\mu_n$ but $\langle G_n(\mu_n,x)\rangle$ is nonelementary. We may assume $\mu_n=r_n$. 

Moreover, there is some $\tau_n<r_n$ such that for $\tau_n\leq \rho < r_n$, the groups $\langle G_n(\rho,x)\rangle$ are all elementary and isomorphic, equal to the group $\langle G_n(\tau_n,x)\rangle$. 

Suppose that the elementary group $\langle G_n(\tau_n,x)\rangle$ is infinite with two fixed points on $\bdy \HH^3$. Then \refprop{ClassInfElementary} implies that it contains a loxodromic\index{loxodromic} $A_n\in G_n(\tau_n,x)$ fixing a line $\ell$. Since $\langle G_n(r_n,x)\rangle$ is not elementary, $G_n(r_n,x)$ must contain a loxodromic $B_n$ that does not fix $\ell$. Then $A_n$ and $B_n$ are loxodromics in $G_n(r_n,x)$ with no common fixed points. So $\langle A_n, B_n\rangle$ is not elementary.

Now suppose that the elementary group $\langle G_n(\tau_n,x)\rangle$ fixes a single point $\zeta \in \bdy\HH^3$. Then $G_n(\tau_n,x)$ contains a parabolic\index{parabolic} $A_n$. Since $\langle G_n(r_n,x)\rangle$ is not elementary, $G_n(r_n,x)$ contains some $B_n$ that does not fix $\zeta$. So again $A_n$ and $B_n$ are elements of $G_n(r_n,x)$ with no common fixed points, and $\langle A_n, B_n\rangle$ is not elementary.

Finally suppose that the elementary group $\langle G_n(\tau_n,x)\rangle$ consists only of the identity element. Since $\langle G_n(r_n,x)\rangle$ is nonelementary with no elliptics, the generating set $G_n(r_n,x)$ must contain two elements $A_n$ and $B_n$ with no common fixed point. Thus $\langle A_n, B_n \rangle$ is not elementary.

In all cases, we have a nonelementary subgroup with two generators, $\langle A_n,B_n\rangle$, and $A_n$, $B_n \in G_n(r_n,x)$. Note that $A_n(x)\to x$ and $B_n(x)\to x$, so \reflem{Sequences} implies there are subsequences of $\{A_n\}$ and $\{B_n\}$ converging to $A\in\PSL(2,\CC)$ and $B\in\PSL(2,\CC)$, respectively. Then \refcor{LimitNonelementary} implies that $\langle A,B \rangle$ is nonelementary.

On the other hand, $A_n, B_n \in G_n(r_n,x)$, so as $n\to\infty$, $A_n$ and $B_n$ must converge to elements of $\PSL(2,\CC)$ fixing $x$. Thus $\langle A, B\rangle$ fixes $x$, hence it is elementary by definition. This contradiction finishes the proof. 
\end{proof}

%%%%%%%%%%%%%%%%%%%%%%%%%%%%%%%%%%%%%%%%%%%%%%%%%%%%%%%%%%%%%%%%%
\subsection{A sketch of a broader result in Lie groups}

The Universal Elementary Neighborhoods theorem, \refthm{Margulis}, which we proved using properties of $\PSL(2,\CC)$ in the previous subsection, actually follows quickly from a broader result in Lie groups due to Ka{\v{z}}dan and Margulis~\cite{KazhdanMargulis}. We will not go into many details on Lie groups here, but we do include a sketch of some of the ideas. 

\begin{definition}\label{Def:Commutator}
  Let $G$ be a group with subgroups $H$ and $K$. The group $[H,K]$ is defined to be the subgroup of $G$ generated by elements $[h,k] = hkh^{-1}k^{-1}$ for all $h\in H$ and $k\in K$.

  The $m$-th \emph{commutator}\index{commutator} $G^m$ of $G$ is defined recursively by $G^1 = [G,G]$, and $G^{m+1} = [G,G^m]$, for $m\geq 1$.

  A group is \emph{nilpotent}\index{nilpotent group} if for some integer $m$, $G^m=\{1\}$. 
\end{definition}

The following is due to Zassenhaus, proved in 1937~\cite{Zassenhaus}. 

\begin{theorem}[Zassenhaus Theorem]\label{Thm:Zassenhaus}\index{Zassenhaus Theorem}
  Let $G$ be a Lie group. Then there is a neighborhood of the identity $U_Z \subset G$ such that for each discrete subgroup $\Gamma\leq G$, the group generated by $\Gamma \cap U_Z$ is nilpotent.
\end{theorem}

\begin{proof}[Proof sketch]
The derivative of the commutator map $[\cdot,\cdot]\from G\times G\to G$ at $(1,1)$ can be shown to be identically $0$, so $[\cdot,\cdot]$ is a strict contraction in a neighborhood $U$ of the identity, in both variables. Thus for $\gamma_1, \dots, \gamma_m\in \Gamma\cap U$, the iterated commutator
\[ y_m = [\gamma_1, [\gamma_2, [\dots [\gamma_{m-1},\gamma_m]\dots]]] \]
must lie in $U$ and must satisfy $\lim_{m\to\infty} y_m = 1$. Because $\Gamma$ is a discrete group, there exists an integer $N$ such that for $n\geq N$, $y_N =1$. Then the group is nilpotent. 
\end{proof}

This in turn implies a more general result on Lie groups, proved in \cite{KazhdanMargulis}.
Before we state the theorem, we say a few words about the general setting in which the theorem applies.

Let $G$ be a Lie group, and $K$ a maximal compact subgroup of $G$. We may give $G$ a left-invariant Riemannian metric that is also right-invariant under $K$. Then the space $G/K$ becomes a Riemannian manifold with $G$ acting on $X$ on the left by isometries of $X$. We say $X = G/K$ is the \emph{homogeneous space associated with $G$}.\index{homogeneous space}

For example, in the setting of $\HH^3$, we may take $G$ to be the group of isometries of $\HH^3$, and $K$ the subgroup fixing a point $x\in\HH^3$. This is isomorphic to the compact Lie group $O(2)$. Then the quotient $G/K$ is $\HH^3$, with its usual metric and action of $G$ by isometries. We will apply the theorem in this setting. 

\begin{theorem}[Kazhdan--Margulis Theorem]\label{Thm:MargulisLemmaGeneral}\index{Ka{\v{z}}dan--Margulis Theorem}
  Let $X$ be the homogeneous space associated with a Lie group $G$. There exists a constant $\eta = \eta(X)$ satisfying the following.
  Let $x\in X$, and let $\Gamma$ be any discrete group generated by elements $\{ g_1, \dots, g_\ell\}\subset G$ such that $d(x,g_j(x))\leq \eta$ for all $j$. Then there exists a subgroup $\Gamma'$ of $\Gamma$ of finite index such that $\Gamma'$ is nilpotent.
\end{theorem}

A proof of this version of the Margulis Lemma can be found in~\cite{kapovich}. See also~\cite{BallmannGromovSchroeder} for a version that applies to Riemannian manifolds with negative sectional curvature, or~\cite{benedetti-petronio} for another proof when $X=\HH^n$.

\begin{proof}[Proof sketch]
The proof begins by taking a Zassenhaus neighborhood $U_Z$ of $1\in G$ from \refthm{Zassenhaus}.
There exists $\epsilon>0$ depending only on $X$ such that the ball of radius $\epsilon$ around $1\in G$ is contained in $U_Z$: $B_\epsilon(1)\subset U_Z$.

Next, because $X$ is homogeneous, we may assume $x$ is the projection of $1\in G$ to $X=G/K$, removing the dependence of the argument upon $x$. 

The value of $\eta$ is determined from $\epsilon$ as follows. Because $K$ is compact, there is an $\epsilon/10$-dense subset of $K$ consisting of a finite number of elements; say $N$ elements. Choose $\eta$ such that whenever $\{g_1, \dots, g_\ell\}$ satisfy $d(x,g_j(x))\leq \eta$, any word $w=w(g_1, \dots, g_\ell)$ in the $g_j$ of length at most $N$ satisfies $d(x,wx)\leq \epsilon/5$. 

For this value of $\eta$, whenever such $\{g_1, \dots, g_\ell\}$ generate a discrete group $\Gamma$, the group $\Gamma\cap B_\epsilon(1)=\Gamma'$ is nilpotent, by \refthm{Zassenhaus}. The choice of $\eta$ allows one to show that $\Gamma'$ also has finite index in $\Gamma$.
\end{proof}

Assuming the Ka{\v{z}}dan--Margulis theorem, \refthm{MargulisLemmaGeneral}, we obtain a quick proof of \refthm{Margulis}, the Universal Elementary Neighborhoods theorem, which we now explain.

Recall that the \emph{center}\index{center} of a group is the subgroup of all elements that commute with every other element.

\begin{lemma}\label{Lem:NilpotentCenter}
  A non-trivial nilpotent group has non-trivial center.
\end{lemma}

\begin{proof}
Suppose $G$ is nilpotent, with $G^n=[G,G^{n-1}] =1$ but $G^{n-1}\neq 1$. Then
$[G,G^{n-1}]=1$ if and only if for every $x\in G^{n-1}$ and every $g\in G$, the product $x^{-1}g^{-1}xg=1$, which holds if and only if $xg=gx$. Thus $G^{n-1}$ lies in the center of $G$, and is nontrivial.
\end{proof}

\begin{corollary}\label{Cor:Nilpotent}
  A nilpotent subgroup $G$ of $\PSL(2,\CC)$ without elliptics must satisfy one of the following:
  \begin{itemize}
  \item $G = \{1\}$
  \item $G-\{1\}$ consists of loxodromic\index{loxodromic} elements with the same fixed points at infinity.
  \item $G-\{1\}$ consists of parabolic\index{parabolic} elements with the same fixed point at infinity.
  \end{itemize}
\end{corollary}

\begin{proof}
By \reflem{NilpotentCenter}, if $G$ is nontrivial then there is a nontrivial element $g\in G$ that commutes with every other element of $G$. By \refex{PSL(2,C)Commute}, every element in $G$ must have the same fixed points as $g$. The cases follow depending on whether $g$ is loxodromic\index{loxodromic} or parabolic.\index{parabolic}
\end{proof}

\begin{proof}[Proof of \refthm{Margulis} assuming \refthm{MargulisLemmaGeneral}]
Let $G$ be a discrete subgroup of $\PSL(2,\CC)$ without elliptics, let $x\in\HH^3$, and let $\eta$ be the constant from the Ka{\v{z}}dan--Margulis Theorem, \refthm{MargulisLemmaGeneral}. If $H$ denotes the subgroup of $G$ generated by elements of $G$ that translate $x$ distance less than $\eta$, then there exists a nilpotent subgroup $H'$ of $H$ such that $H/H'$ is finite. By \refcor{Nilpotent}, $H'$ has one of three forms. 

If $H'$ is trivial, then $H$ is a finite group. Since there are no elliptics in $G$, $H$ must also be trivial, and so $H$ is elementary.\index{elementary group}

Since $H'$ is a finite index subgroup of $H$, for any $h\in H$ there exists an integer $m$ such that $h^m\in H'$. Then $h^m$ has the same fixed points as $H'$, and hence $h$ has the same fixed points as $H'$. Thus either $H'-\{1\}$ consists of loxodromic\index{loxodromic} elements with two fixed points at infinity, and all elements of $H$ have the same fixed points at infinity, or $H'-\{1\}$ consists of parabolics\index{parabolic} with one point at infinity, and all elements of $H$ have the same fixed point at infinity. In either case, $H$ is elementary.\index{elementary group}
\end{proof}

%%%%%%%%%%%%%%%%%%%%%%%%%%%%%%%%%%%%%%%%%%%%%%%%%%%%%%%%%%%%%%%%%
\section{Exercises}

\begin{exercise} Is a subgroup of $\PSL(2,\CC)$ generated by a single elliptic\index{elliptic} element always discrete? Prove it is discrete, or give a counterexample.
\end{exercise}

\begin{exercise}\label{Ex:MoreClassifyPSL}
Prove \reflem{MoreClassifyPSL}, giving more properties of parabolic,\index{parabolic} elliptic,\index{elliptic} and loxodromics\index{loxodromic} in $\PSL(2,\CC)$. 
\end{exercise}

\begin{exercise}\label{Ex:Fig8Gluings}
Prove that the gluing isometries for the figure-8 knot complement are the elements of $\PSL(2,\CC)$ given in \refeqn{Fig8Gluings}.
\end{exercise}

\begin{exercise}\label{Ex:Riley}
  R.~Riley gave a presentation of the fundamental group of the figure-8 knot complement in \cite{Riley:Fig8}:
\[ \pi_1(S^3-K) = \langle a,b \mid yay^{-1}=b\rangle, \]
where $y=a^{-1}bab^{-1}$. He let
\[ A=\mat{1&1\\0&1}, \quad B=\mat{1&0\\-\sigma&1},\]
where $\sigma$ is a primitive cube root of unity, 
and let
\[\rho\from \pi_1(S^3-K)\to \langle A,B\rangle \leq \PSL(2,\CC)\]
be the representation $\rho(a)=A$, $\rho(b)=B$.
Prove the representation $\rho$ gives an isomorphism of groups.
\end{exercise}

\begin{exercise}\label{Ex:Riley2}
Let $A$ and $B$ in $\PSL(2,\CC)$ be as in \refex{Riley}. 
Find an explicit element $U=\mat{a&b\\c&d}$ of $\PSL(2,\CC)$ such that Riley's $A$ and $B$ are conjugate via $U$ to our isometries $T_C$ and $T_D^{-1}$,  respectively. That is, find $U$ such that
\[ A = UT_CU^{-1}, \quad B=UT_D^{-1}U^{-1}. \]
Even better: $U$ can be written as a composition of a parabolic\index{parabolic} fixing infinity $T$, followed by a rotation $R$: $U=RT$. Find $T$ and $R$. 
%% Answer: $U=\mat{e^{-i\pi/6}&0\\0&e^{i\pi/6}}\mat{1&-1\\0&1}$
\end{exercise}

\begin{exercise}\label{Ex:Riley3}
Note that Riley's isometries $A$ and $B$ of \refex{Riley} do not give face-pairings\index{face-pairing isometry} of the fundamental domain in \reffig{Fig8FundRegion}. Find a fundamental domain for the figure-8 knot such that $A$ and $B$ are face-pairing isometries.

Hint: \refex{Riley2} might be helpful. 
\end{exercise}

\begin{exercise} If a group $G$ acts on Euclidean space $\RR^n$ or hyperbolic space $\HH^n$, extend the definitions of properly discontinuous\index{properly discontinuous}\index{group action!properly discontinuous} and free actions\index{free group action}\index{group action!free} in the obvious way.

Show directly by definitions that each of the following groups $G$ acts freely\index{free group action}\index{group action!free} and properly discontinuously\index{properly discontinuous}\index{group action!properly discontinuous} on the given space $X$.
\begin{enumerate}
\item $X=\RR^2$, $G$ is generated by two translations $\phi\from \RR^2\to\RR^2$ and $\psi\from \RR^2\to\RR^2$ given by $\phi(x,y)=(x+t,y)$ and $\psi(x,y)=(x,y+s)$ for $s,t\in\RR$.
\item $X=\HH^2$, $G$ is the holonomy group\index{holonomy group} of the (complete) 3-punctured sphere.
\item $X=\HH^3$, $G$ is generated by face-pairing isometries\index{face-pairing isometry} of an ideal polyhedron such that the face identifications give a complete hyperbolic 3-manifold.
\end{enumerate}
\end{exercise}

\begin{exercise}\label{Ex:FiniteElementary}
Show that the following give finite elementary groups.\index{elementary group}
\begin{enumerate}
\item Cyclic groups fixing an axis in $\HH^3$.
\item Orientation preserving symmetries of an ideal platonic solid (tetrahedron, octahedron/cube, icosahedron/dodecahedron). 
\item Dihedral groups preserving an ideal polygon with $n$ sides inscribed in a plane in $\HH^3$.
\end{enumerate}
\end{exercise}

\begin{exercise}
Show that the finite groups in \refex{FiniteElementary} are the only finite elementary groups.\index{elementary group}
\end{exercise}

\begin{exercise}
  Let $G$ be a subgroup of $\PSL(2,\CC)$. Show that the following are equivalent.
  \begin{enumerate}
  \item $G$ is discrete.
  %% \item No sequence $\{A_n\}$ of distinct elements of $G$ converges to an element $A$ of $\PSL(2,\CC)$.
  \item $G$ has no limit points in the interior of $\HH^3$. That is, for any $x\in\HH^3$, there is no $y\in\HH^3$ and no sequence of distinct elements $\{A_n\}$ in $G$ such that $A_n(y)=x$.
  \end{enumerate}
\end{exercise}

\begin{exercise}\label{Ex:PSL(2,C)Commute}
Let $A$ and $B$ in $\PSL(2,\CC)$ be distinct from the identity.  Prove that the following are equivalent.
\begin{enumerate}
\item[(a)] $A$ and $B$ commute.
\item[(b)] Either $A$ and $B$ have the same fixed points, or $A$ and $B$ have order $2$ and each interchanges the fixed points of the other.
\item[(c)] Either $A$ and $B$ are parabolic\index{parabolic} with the same fixed point at infinity, or the axes of $A$ and $B$ coincide, or $A$ and $B$ have order $2$ and their axes intersect orthogonally in $\HH^3$.
\end{enumerate}
\end{exercise}

\begin{exercise}
  Suppose that $A$ and $B$ in $\PSL(2,\CC)$ are loxodromics\index{loxodromic} with exactly one fixed point in common. Show that $\langle A, B\rangle$ is not discrete.
\end{exercise}

\begin{exercise}
  State and prove a version of \refthm{ThinPart}, the structure of the thin part, for hyperbolic 2-manifolds. 
\end{exercise}

\begin{exercise}\label{Ex:CalculateHypDist}
Suppose $A$ is a parabolic\index{parabolic} fixing the point $\zeta$ and $p$ is a point in $\HH^3$ such that $d(p, A(p)) < \epsilon$. After applying an isometry, we may assume that $\zeta=\infty$, that $A=\mat{1&\alpha\\0&1}$ for some $\alpha\in\CC$, and $p$ lies on a horosphere $H_C$ that is a Euclidean plane of constant height $t=C$ for some $C>0$:
\[ H_C = \{(x+y\,i,C) \mid C > 0\}. \]
\begin{enumerate}
\item[(a)] Prove that if a point $q$ lies inside the horoball bounded by $H_C$ on a horosphere $H_t$ of height $t\geq C$, then the Euclidean distance from $q$ to $A(q)$ measured along $H_t$ is at most the Euclidean distance from $p$ to $A(p)$ measured along $H$.
\item[(b)] Prove that the hyperbolic distances, measured in $\HH^3$, satisfy
\[ \epsilon > d(p, A(p)) \geq d(q, A(q)). \]
\end{enumerate}
\end{exercise}

\begin{exercise}\label{Ex:CalculateDistTube}
Suppose $A$ is a loxodromic\index{loxodromic} fixing an axis $\ell$, and $p$ is a point in $\HH^3$ such that $d(p,A(p))<\epsilon$.
\begin{enumerate}
\item Prove that the distance from any $q\in\HH^3$ to $\ell$ is the same as the distance from $A(q)$ to $\ell$.

\item We can use cylindrical coordinates in $\HH^3$ about the geodesic $\ell$. Let $r$ denote the distance from $\ell$, $\theta$ the rotation about $\ell$ (measured modulo $2\pi$), and $\zeta$ the translation distance along $\ell$. Finally, let $\widehat{\HH}^3$ denote the cover of $\HH^3$ in which $\theta$ is no longer measured modulo $2\pi$, but is a real number.

Using these coordinates, it can be shown that the distance $d$ between points $p_1$ and $p_2$ in $\widehat{\HH}^3$ with cylindrical coordinates $(r_1,\theta_1, \zeta_1)$ and $(r_2,\theta_2,\zeta_2)$ with $|\theta_1-\theta_2|<\pi$ is given by
\[\cosh d = \cosh(\zeta_1-\zeta_2)\cosh r_1 \cosh r_2 - \cos(\theta_1-\theta_2)\sinh r_1 \sinh r_2. \]
(See \cite[Lemma~2.1]{GabaiMeyerMilley:Tubes})

Using this formula, prove that if $x,y\in\HH^3$ are points such that $d(y,\ell)\leq d(x,\ell)$, then
\[ d(y,A(y)) \leq d(x, A(x)). \]
\end{enumerate} 
\end{exercise}

%% Ch06_Completion.tex

\chapter{Completion and Dehn Filling}\label{Chap:CompletionDehnFilling}
\blfootnote{Jessica S. Purcell, Hyperbolic Knot Theory}

In \refchap{Geometric} we considered some incomplete structures on hyperbolic 2-manifolds, particularly the 3-punctured sphere, \refexamp{Incomplete3PunctSphere}. 
In this chapter, we examine incomplete hyperbolic structures on 3-manifolds with torus boundary, and their completions.

\section{Mostow--Prasad rigidity}

We begin by stating a few important results on complete hyperbolic structures on manifolds to set up some context for the rest of the chapter. 

Many surfaces admit infinitely many complete hyperbolic structures. For example, in exercises~\ref{Ex:1punctTorus} and~\ref{Ex:4PunctSphere} you found 2-parameter families of complete hyperbolic structures on the 1-punctured torus and 4-punctured sphere. This flexibility is only possible in two dimensions.  In higher dimensions, there is only \emph{one} complete structure\index{complete metric space} on a finite volume hyperbolic manifold, up to isometry. This result was proved in the case $M$ is a closed manifold by Mostow \cite{mostow}, and extended to the case of open manifolds with finite volume by Prasad \cite{prasad}. Recall that by \refthm{FteVolIffTorusBdy}, an open hyperbolic 3-manifold has finite volume if and only if it is the interior of a manifold with torus boundary components. 

\begin{theorem}[Mostow--Prasad rigidity]
\label{Thm:MostowGeom}\index{Mostow--Prasad rigidity}
If $M_1^n$ and $M_2^n$ are complete\index{complete metric space} hyperbolic $n$-manifolds with finite volume and $n\geq 3$, then any isomorphism of fundamental groups $\phi\from \pi_1(M_1) \to \pi_1(M_2)$ is realized by a unique isometry. 
\end{theorem}

We will not include the proof in this book, as it leads us a little further away from knots and links than we wish to stray. However, the proof of the theorem can be found in the original papers, or in books on hyperbolic geometry including \cite{benedetti-petronio} and \cite{ratcliffe}.

Recall also Gordon and Luecke's knot complement theorem, \refthm{GordonLuecke} from \refchap{KnotIntro}, which states that knots with homeomorphic complement are equivalent.\index{Gordon--Luecke knot complement theorem}

Knots with homeomorphic complements have isomorphic fundamental group. By \refthm{MostowGeom}, Mostow--Prasad rigidity,\index{Mostow--Prasad rigidity} any complete hyperbolic structure on the knot complement is the \emph{only} complete hyperbolic structure.\index{complete metric space} So the complete hyperbolic structure on a knot complement distinguishes any two knots. This is one reason hyperbolic geometry gives many very nice knot invariants!

%%%%%%%%%%%%%%%%%%%%%%%%%%%%%%%%%%%%%%%%%%%%%%%%%%%%%%%%%%%%%%%%%
\section{Completion of incomplete structures}

What about incomplete structures on a manifold $M$ with torus boundary? There are many of these. For the figure-8 knot complement, for example, we found a 1-complex parameter family of incomplete structures, parameterized by $w \in \CC$ as in \reffig{ThurstonRegion}. If we take the \emph{completion} of a hyperbolic structure on a 3-manifold, we obtain surprising topological results.

As a warm up, recall completions of incomplete structure on 2-manifolds.  In \refchap{Geometric}, we saw an example of an incomplete structure on a hyperbolic 3-punctured sphere. Recall that in the developing map for an incomplete structure, ideal polygons approached a limiting line.  By selecting a point on a horocycle about infinity, approaching this line, we obtained a Cauchy sequence that did not converge. See \reffig{3punct-incomplete}. Adjoining a point where each horocycle met the limiting line, we obtained the completion. The completion was given by attaching a geodesic of length $d(v)$,\index{$d(v)$} as in \reffig{3punctCompletion}.

Now consider an incomplete structure on a 3-manifold $M$ such that $M$ is the interior of a compact manifold with torus boundary. Let $C$ be a cusp torus of $M$. Then the torus $C$ inherits an affine structure from the hyperbolic structure on $M$, and because the structure on $M$ is not complete, the affine structure is not Euclidean (\refthm{EuclidCusp}).\index{Euclidean structure}

Let $\alpha$ and $\beta$ generate $\pi_1(C) \cong \ZZ \times\ZZ$. Corresponding to $\alpha$ and $\beta$ are two holonomy\index{holonomy} isometries $\rho(\alpha)$ and $\rho(\beta)$. To simplify notation, we will drop the $\rho$, abusing notation slightly, and simply refer to these isometries as $\alpha$ and $\beta$.  Assume the action of $\alpha$ and $\beta$ does not induce a Euclidean structure\index{Euclidean structure} on $C$, so the hyperbolic structure on $M$ is not complete. To form its completion, we remove a small neighborhood $N(C)$ of $C$, take the completion of $N(C)$, and then reattach this neighborhood to $M$. Thus to analyze the completion of $M$, we analyze the completion of neighborhoods of cusp tori.

\begin{proposition}\label{Prop:CompletionAttachGeodesic}
The completion of $N(C)$ is obtained by adjoining some portion of a geodesic to $N(C)$.
\end{proposition}

\begin{proof}
Consider the developing map for the affine torus $C$.\index{affine torus} The image will miss a single point (\refex{DevelopingImageMissesPoint}), for example as in \reffig{AffineTorus}. This image is obtained by considering the action of $\alpha$ and $\beta$ restricted to a horosphere. More precisely, if $C$ has a fundamental domain that is a quadrilateral, then we build its developing image by starting with a copy of that quadrilateral on $\CC$, which we identify with a horosphere about infinity, and attaching copies of the quadrilateral according to instructions given by the holonomy\index{holonomy} isometries corresponding to $\alpha$ and $\beta$, acting on the fixed horosphere. 

If we shift the original choice of horosphere up, we will see the same image of the developing map.  In particular, the developing map will still miss a single point, with the same complex value for each choice of horosphere. These missed points form a vertical geodesic in $\HH^3$. We may apply an isometry so that this vertical geodesic runs from $0$ to $\infty$ in $\HH^3$. Notice that the developing image of the neighborhood $N(C)$ is obtained by taking developing images of $C$ on all horospheres about $\infty$ above some fixed initial height. Thus the developing image $N(C)$ misses the single geodesic from $0$ to $\infty$ in $\HH^3$. Hence the completion of $N(C)$ is obtained by adjoining some portion of this geodesic to $N(C)$.
\end{proof}

As in the case of incomplete 2-manifolds, the length of the portion of adjoined geodesic of \refprop{CompletionAttachGeodesic} will be determined by considering the action of the holonomy.\index{holonomy} Considering this action leads to the following result on the topology of the completion.

\begin{proposition}\label{Prop:CompletionTop}
Let $N(C)$ be the neighborhood of a cusp torus $C$ of an incomplete hyperbolic manifold, so $N(C)$ is homeomorphic to $C\times(0,1)$. Then the completion of $N(C)$ is either homeomorphic to the 1-point compactification of $N(C)$ obtained by crushing $C\times\{1\}$ to a point, or it is homeomorphic to the solid torus obtained by attaching a solid torus to $C\times\{1\}$. 
\end{proposition}

\begin{proof}
As in the proof of \refprop{CompletionAttachGeodesic}, consider the developing image of $N(C)$ and assume it misses the geodesic from $0$ to $\infty$. Note the group $\langle\alpha,\beta\rangle$ acts on the geodesic from $0$ to $\infty$. Since points in our completion should be identified to their images under the holonomy\index{holonomy} action, we should identify each point $z$ on the geodesic from $0$ to $\infty$ with $\langle \alpha, \beta\rangle \cdot z$. There are two cases. 

\smallskip

\textbf{Case 1.}  The image of $z$ under the action of $\alpha$ and $\beta$ is dense in the line from $0$ to $\infty$.  In this case, the completion is the 1-point compactification.  It is not a manifold (\refex{1PtCmpt}).

\smallskip

\textbf{Case 2.} The image of $z$ is a discrete set of points on the line, each of some distance $d(C)$ apart. In this case the completion is obtained by adjoining a geodesic circle of length $d(C)$ to $N(C)$. Denote the completion by $\overline{N(C)}$. We wish to understand the topology of $\overline{N(C)}$.

We may obtain a manifold homeomorphic to $N(C)$ by removing a small, closed tubular neighborhood\index{tubular neighborhood} of the geodesic circle adjoined to form $\overline{N(C)}$. Notice that a tubular neighborhood of a circle is a solid torus, with the geodesic at its core. Thus we obtain a manifold homeomorphic to $\overline{N(C)}$ by attaching a solid torus to the torus $C\times\{1\}$ of $N(C)$
\end{proof}

\begin{definition}\label{Def:DehnFilling}
Let $M$ be a manifold with torus boundary component $T$. Let $s$ be an isotopy class of simple closed curves on $T$; $s$ is called a \emph{slope}.\index{slope} The manifold obtained from $M$ by attaching a solid torus to $T$ so that $s$ bounds a disk in the resulting manifold is called the \emph{Dehn filling of $M$ along $s$}\index{Dehn filling} and is denoted $M(s)$. 
\end{definition}

\begin{figure}
%% Accompanies image file 'F6-01-Dehnfl.eps' (pdf, eps, ps)
\begingroup%
  \makeatletter%
  \providecommand\color[2][]{%
    \errmessage{(Inkscape) Color is used for the text in Inkscape, but the package 'color.sty' is not loaded}%
    \renewcommand\color[2][]{}%
  }%
  \providecommand\transparent[1]{%
    \errmessage{(Inkscape) Transparency is used (non-zero) for the text in Inkscape, but the package 'transparent.sty' is not loaded}%
    \renewcommand\transparent[1]{}%
  }%
  \providecommand\rotatebox[2]{#2}%
  \newcommand*\fsize{\dimexpr\f@size pt\relax}%
  \newcommand*\lineheight[1]{\fontsize{\fsize}{#1\fsize}\selectfont}%
  \ifx\svgwidth\undefined%
    \setlength{\unitlength}{178.11917496bp}%
    \ifx\svgscale\undefined%
      \relax%
    \else%
      \setlength{\unitlength}{\unitlength * \real{\svgscale}}%
    \fi%
  \else%
    \setlength{\unitlength}{\svgwidth}%
  \fi%
  \global\let\svgwidth\undefined%
  \global\let\svgscale\undefined%
  \makeatother%
  \begin{picture}(1,0.39937447)%
    \lineheight{1}%
    \setlength\tabcolsep{0pt}%
    \put(0,0){\includegraphics[width=\unitlength]{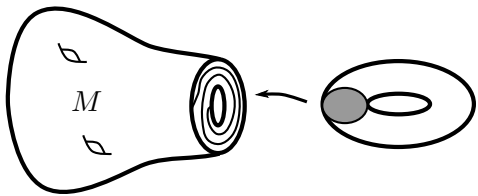}}%
    \put(0.13908975,0.18142582){\color[rgb]{0,0,0}\makebox(0,0)[lt]{\lineheight{0}\smash{\begin{tabular}[t]{l}$M$\end{tabular}}}}%
  \end{picture}%
\endgroup%

\caption{A cartoon describing Dehn filling.\index{Dehn filling} After filling, the curve shown on the torus boundary component of $M$ will bound a disk.}
\label{Fig:DehnFilling}
\end{figure}

A cartoon describing Dehn filling is shown in \reffig{DehnFilling}.

By \refprop{CompletionTop}, the space obtained by taking the completion of an incomplete hyperbolic structure on $M$ either fails to be a manifold, or is homeomorphic to a Dehn filling\index{Dehn filling} of $M$.

Dehn filling\index{Dehn filling} is a very important topological procedure in 3-manifold topology, due to work of Wallace and Lickorish in the 1960s. Independently, they showed the following theorem \cite{wallace, lickorish}. A nice, highly readable proof can be found in the book \cite{rolfsen}.

\begin{theorem}[Fundamental theorem of Wallace and Lickorish]\label{Thm:WallaceLickorish}\index{Dehn filling!fundamental theorem of Wallace and Lickorish}\index{Wallace, fundamental Dehn filling theorem}\index{Lickorish, fundamental Dehn filling theorem}
Let $M$ be a closed, orientable 3-manifold.  Then $M$ is obtained by Dehn filling\index{Dehn filling} the complement of a link in $S^3$. \qed
\end{theorem}

\Refthm{WallaceLickorish} gives a topological result on manifolds. By considering completions of hyperbolic 3-manifolds, we can make Dehn filling a geometric procedure. 

\begin{definition}\label{Def:ConeManifold}
Consider the geodesic running from $0$ to $\infty$ in $\HH^3$. We may write points in $\HH^3$ in cylindrical coordinates $(r, \theta,\zeta)$ where $r$ is the distance from this geodesic, $\theta$ is a rotation angle around the geodesic, measured modulo $2\pi$, and $\zeta$, the height, is translation distance in the direction of the geodesic. In these coordinates, the metric is given by
\[ dr^2 + \sinh^2 r\, d\theta^2 + \cosh^2 r\, d\zeta^2, \]
with $\theta$ measured modulo $2\pi$.

Now fix $\alpha>0$. Adjust the metric so $\theta$ is measured modulo $\alpha$. Then a neighborhood of a point on the geodesic from $0$ to $\infty$ is called a \emph{hyperbolic cone}\index{hyperbolic cone} with \emph{cone angle}\index{cone angle} $\alpha$. Note the definition makes sense when $\alpha>2\pi$. A cross section perpendicular to the geodesic is a 2-dimensional cone with cone angle $\alpha$. 

A 3-dimensional \emph{hyperbolic cone manifold}\index{hyperbolic cone manifold}\index{cone manifold} is a manifold $M$ in which each point $x$ either has a neighborhood isometric to a ball in $\HH^3$, or has a neighborhood isometric to a hyperbolic cone.

In a hyperbolic cone manifold, the set of points that only have neighborhoods of the second kind form a geodesic link in $M$ called the \emph{singular locus}\index{singular locus}. The hyperbolic metric on $M$ is smooth everywhere except at points on the singular locus.
\end{definition}

\begin{proposition}\label{Prop:CompletionConeMfld}
When the completion $\overline{M}$ of $M$ is topologically equivalent to attaching a solid torus, obtained by Dehn filling,\index{Dehn filling} it has the structure of a cone manifold.\index{hyperbolic cone manifold}\index{cone manifold} The singular locus $\Sigma$ is the geodesic (link) attached in the completion.
\end{proposition}

\begin{proof}
As before, let $C$ be a cusp torus of $M$ with neighborhood $N(C)$, whose developing image misses the geodesic from $0$ to $\infty$ in $\HH^3$. Let $\zeta \in \pi_1(C)$ generate the kernel of the action of $\pi_1(C) \cong \langle \alpha, \beta \rangle$ on the line from $0$ to $\infty$. The isometry $\zeta$ will be a rotation about this line by some angle $\alpha$. Then a perpendicular cross section of the circle added to $\overline{N(C)}$ to form the completion will be a 2-dimensional hyperbolic cone, of cone angle\index{cone angle} $\alpha$. Thus a neighborhood of a point on the completion is isometric to a hyperbolic cone.

Thus when we attach $\overline{N(C)}$ to $M$, the result $\overline{M}$ is a hyperbolic cone manifold\index{hyperbolic cone manifold}\index{cone manifold} with singular locus along the attached geodesic.
\end{proof}

There is one very important case of \refprop{CompletionConeMfld}. When the cone angle\index{cone angle} at the singular locus of $\overline{M}$ is actually $2\pi$, then the hyperbolic structure on $\overline{M}$ is smooth everywhere. Thus $\overline{M}$ is a hyperbolic manifold. We conclude:

\begin{corollary}\label{Cor:HypDehnFilling}
  When the holonomy\index{holonomy} $\rho(\pi_1(C))$ acts on the geodesic omitted from the developing image of $N(C)$ by a fixed translation, and when the generator $\zeta \in \pi_1(C)$ of the kernel has holonomy\index{holonomy} a rotation by $2\pi$, then the completion of $M$ is a complete hyperbolic manifold, homeomorphic to the Dehn filled manifold $M(\zeta)$.\qed
\end{corollary}

%%%%%%%%%%%%%%%%%%%%%%%%%%%%%%%%%%%%%%%%%%%%%%%%%%%%%%%%%%%%%%%%%
\section{Hyperbolic Dehn filling space}

We re-interpret the above section in the language of complex lengths of isometries of $\HH^3$.

Anytime $M$ admits a hyperbolic structure, consider a cusp torus $C$ for $M$. 
The fundamental group of the torus is isomorphic to $\ZZ\times \ZZ$, generated by some $\alpha$ and $\beta$.

\begin{remark}\label{Rem:MeridLongitude}
When $M$ is a knot complement, $M\cong S^3\setminus N(K)$, we often choose $\alpha$ to be the \emph{meridian},\index{meridian} i.e.\ the curve on $\bdy N(K)$ bounding a disk in $N(K)\subset S^3$, and $\beta$ to be the \emph{standard longitude},\index{longitude!standard} i.e.\ the curve on $\bdy N(K)$ that is homologous to $0$ in $S^3\setminus N(K)$.
\end{remark}

Consider the holonomy\index{holonomy} elements of $\alpha$ and $\beta$. These are some isometries of $\HH^3$. As above, we will continue to abuse notation and denote the holonomy isometries corresponding to $\alpha$ and $\beta$ by $\alpha$ and $\beta$.

Recall the classification of isometries of $\HH^3$, from \reflem{MoreClassifyPSL}. Any isometry is one of three types: parabolic,\index{parabolic} elliptic,\index{elliptic} or loxodromic.\index{loxodromic} Since $\alpha$ and $\beta$ generate $\ZZ\times \ZZ$, they must commute.  This is possible only if $\alpha$ and $\beta$ are parabolic,\index{parabolic} fixing the same point on the boundary at infinity,\index{boundary at infinity} or if $\alpha$ and $\beta$ share the same axis (\refex{PSL(2,C)Commute}).

If $\alpha$ and $\beta$ are parabolic,\index{parabolic} fixing a point at infinity, then they must fix an entire horosphere about infinity. Conjugating to put their fixed point at $\infty$ in $\bdy\HH^3$, they are of the form $\alpha(z) = z+a$, $\beta(z) = z+b$. Hence they restrict to Euclidean isometries on the horosphere, and the hyperbolic structure is complete.

Now suppose $\alpha$ and $\beta$ are not parabolic.\index{parabolic}
In this case, because $\alpha$ and $\beta$ commute but are not parabolic, they share an axis, and are both given by rotation and/or dilation along this axis. The hyperbolic structure is not complete, and the axis must be exactly the geodesic whose points are omitted from the developing image of $C$ for each horosphere.

\begin{definition}
Suppose the interior of $M$ has a hyperbolic structure, and $C$ is a cusp torus of $M$, with $N(C)$, homeomorphic to $T^2\times I$, a neighborhood of $C$. Let $\alpha, \beta \in \pi_1(C)$ be generators. Suppose the interior of $M$ has a hyperbolic structure, and the holonomy\index{holonomy} elements corresponding to $\alpha$ and $\beta$ are not parabolic, so they share an axis. 
Fix a direction on the axis of $\alpha$ and $\beta$.  Any element $\gamma$ of $\pi_1(C)$ translates some signed distance $d$ along the axis, and rotates by total angle $\theta \in \RR$, where the sign of $\theta$ is given by the right hand rule.  Let $\calL(\gamma) = d + i\theta$.  The value $\calL(\gamma)$ is called the \emph{complex length}\index{complex length} of $\gamma$.  This defines a function $\calL$ from $\pi_1(C) = H_1(C; \ZZ)$ to $\CC$.
\label{Def:ComplexLength}
\end{definition}

Notice that if $\gamma = p\alpha + q\beta$, then $\calL(\gamma) = p\calL(\alpha) + q\calL(\beta)$, so $\calL$ is a linear map. We may extend it canonically to a linear map $\calL\from H_1(C; \RR) \to \CC$. The value $\calL(c)$ for any $c\in H_1(C; \RR) \cong \RR^2$ will be called the complex length of $c$.

Suppose that the complex length of a simple closed curve $\gamma$ on $C$ equals $2\pi i$. Then in the completion of $M$, $\gamma$ will bound a smooth hyperbolic disk. This implies that the completion of $M$ is a manifold homeomorphic to the Dehn filled manifold $M(\gamma)$, and that $M(\gamma)$ admits a complete hyperbolic structure.

Suppose instead that the complex length of a closed curve $\gamma$ on $C$ equals $\theta i \neq 2\pi i$.  Then in the completion of $M$, $\gamma$ will bound a hyperbolic cone, with cone angle\index{cone angle} $\theta$.  The completion of $M$ is still homeomorphic to the Dehn filled manifold $M(\gamma)$. However, the metric on $M(\gamma)$ inherited from the completion of $M$ is not smooth. The core of the added solid torus is the singular locus, with cone angle $\theta$.

For an incomplete structure, there will be a unique element $c \in H_1(C;\RR)$ so that $\calL(c) = 2\pi i$.

\begin{definition}\label{Def:DehnFillingCoeff}
We say $c \in H_1(C;\RR)$ such that $\calL(c) = 2\pi i$ is the \emph{Dehn filling coefficient}\index{Dehn filling coefficient} of the boundary component $C$.
\end{definition}

When $c$ is of the form $(p,q)$, with $p$ and $q$ relatively prime integers, it corresponds to a simple closed curve and the completion is smooth.

We have been looking at a fixed incomplete hyperbolic structure on $M$, and examining possible completions for this fixed structure. Now we turn our attention to a topological manifold $X$, homeomorphic to $M$, and consider all possible hyperbolic structures on $X$.

\begin{definition}\label{Def:DehnFillingSpace}
Let $X$ be a 3-manifold with cusp torus $C$. The subset of $H_1(C;\RR)$ consisting of Dehn filling coefficients of hyperbolic structures on $X$ is called the \emph{hyperbolic Dehn filling space for $X$}.\index{hyperbolic Dehn filling space}\index{Dehn filling space}

If $X$ admits a complete hyperbolic structure, then we let $\infty$ correspond to the complete hyperbolic structure on $X$.
\end{definition}

\begin{theorem}[Thurston's hyperbolic Dehn filling theorem]
  \label{Thm:HypDehnSurgery}\index{Thurston's hyperbolic Dehn filling theorem}
  \index{hyperbolic Dehn filling theorem}
Let $X$ be a 3-manifold homeomorphic to the interior of a compact manifold with boundary a single torus $T$, such that $X$ admits a complete hyperbolic structure. Then hyperbolic Dehn filling space\index{hyperbolic Dehn filling space}\index{Dehn filling space} for $X$ always contains an open neighborhood of $\infty$ in $\RR^2\cup \{\infty\} \cong H_1(T;\RR)\cup\{\infty\}$.

More generally, if $X$ is the interior of a compact manifold with torus boundary components $T_1, \dots, T_n$, and $X$ admits a complete hyperbolic structure, then the hyperbolic Dehn filling space for $X$ contains an open neighborhood of $\infty$ for each $T_i$.
\end{theorem}

\Refthm{HypDehnSurgery} is an important result, and the result, its proofs, and its extensions continue to have useful consequences. 
The first proof of \refthm{HypDehnSurgery} was sketched in Thurston's 1979 notes \cite{thurston}, and uses results on holonomy\index{holonomy} representations. A proof in the case that $X$ admits a geometric triangulation was given in \cite{NeumannZagier}, presented with expanded details in  \cite{benedetti-petronio}. This proof was extended to the case of more general hyperbolic 3-manifolds by Petronio and Porti \cite{PetronioPorti}. Martelli puts these proofs together to give a complete exposition in his recent book~\cite{Martelli}. Additionally, precise universal bounds on the size of the open neighborhood of infinity provided by the theorem were given by Hodgson and Kerckhoff \cite{hk:univ}, about 25 years after \refthm{HypDehnSurgery} was proved.
All the proofs require work.

In chapters~\ref{Chap:Essential} and~\ref{Chap:Volume} we will give full proofs of related results that are weaker than what is claimed in \refthm{HypDehnSurgery}. Here, we provide only a short sketch of the argument that goes into the proof of \refthm{HypDehnSurgery}, and then focus on applications.

\begin{proof}[Proof sketch of \refthm{HypDehnSurgery}]
Suppose first that $X$ is homeomorphic to the interior of a compact manifold with a single torus boundary component $T$, and $X$ admits a complete hyperbolic structure.

Because $X$ is hyperbolic, there is a holonomy\index{holonomy} representation
\[ \rho\from \pi_1(X) \to \PSL(2,\CC) \]
whose image is a discrete group. The fundamental group of the cusp torus $\pi_1(T)$ has image generated by two parabolics\index{parabolic} $\rho(\alpha)$ and $\rho(\beta)$, which we may assume fix the point at infinity in $\HH^3$.

Now Thurston shows that there exists a one-complex parameter family of deformations of the holonomy\index{holonomy} representation \cite[Theorem~5.6]{thurston}.

Each small deformation of the complete hyperbolic structure taking $\rho(\alpha)$ to a loxodromic\index{loxodromic} must take $\rho(\beta)$ to a loxodromic with the same fixed points. As in the discussion above, this extends to an incomplete hyperbolic structure, with Dehn filling coefficient some complex number $d+i\theta = z$, and where $z=\infty$ corresponds to the complete hyperbolic structure on $X$.

To complete the proof, one shows that $z$ varies continuously over a neighborhood of infinity.

When there are $k>1$ cusps, the proof is similar. In this case, there is a $k$-complex parameter family of deformations, with completions giving Dehn filling coefficients $d_1+i\theta_1, \dots, d_k+i\theta_k$. Again one shows that these vary in a neighborhood of $(\infty, \dots, \infty)$. 
\end{proof}

\begin{corollary}\label{Cor:FiniteFillings}
Let $X$ be a manifold with a single torus boundary component such that the interior of $X$ admits a complete hyperbolic metric. Then there are at most finitely many Dehn fillings\index{Dehn filling} of $X$ which do not admit a complete hyperbolic metric. \qed
\end{corollary}

\begin{corollary}\label{Cor:FiniteFillings2}
Let $X$ be a manifold with $n$ torus boundary components $T_1, \dots, T_n$.  For each $T_i$, exclude finitely many Dehn fillings.\index{Dehn filling} The remaining Dehn fillings yield a manifold with a complete hyperbolic structure. \qed
\end{corollary}

Corollaries~\ref{Cor:FiniteFillings} and~\ref{Cor:FiniteFillings2} follow immediately from \refthm{HypDehnSurgery}.

Notice that \refcor{FiniteFillings2} does not rule out the fact that a manifold with more than one torus boundary component may have infinitely many non-hyperbolic Dehn fillings, as in the following example.

\begin{example}\label{Example:WLInfExFillings}
The Whitehead link is the link shown in \reffig{WhiteheadDehn}.\index{Whitehead link} We will see that it admits a complete hyperbolic structure (\refprop{WhiteheadGeom}). 

\begin{figure}
  \includegraphics{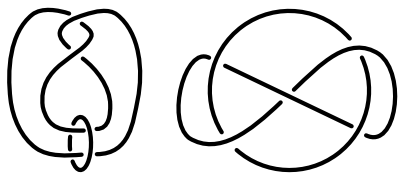}
  \caption{Two diagrams of the Whitehead link.\index{Whitehead link}}
  \label{Fig:WhiteheadDehn}
\end{figure}

If we erase one of the link components, that action can be seen as attaching a solid torus to the link complement in a trivial way. This is called \emph{trivial Dehn filling}.\index{trivial Dehn filling}\index{Dehn filling!trivial} For this example, perform trivial Dehn filling on the component that clasps itself in \reffig{WhiteheadDehn}, leaving a single unknotted component, a trivial knot in $S^3$. Its complement is a solid torus. 

\begin{definition}\label{Def:LensSpace}
  A \emph{lens space}\index{lens space} is the 3-manifold obtained by gluing together two solid tori along their common torus boundary components.
\end{definition}
  
Thus any Dehn filling\index{Dehn filling} of a trivial knot in $S^3$ is a lens space.\index{lens space}

\begin{theorem}\label{Thm:LensSpace}
  A lens space\index{lens space} cannot admit a hyperbolic structure.
\end{theorem}
\begin{proof}
  \Refex{LensSpaceNotHyp}.
\end{proof}

There are infinitely many Dehn fillings\index{Dehn filling} on the trivial knot in $S^3$ that produce lens spaces.\index{lens space} Thus there are infinitely many non-hyperbolic Dehn fillings of the Whitehead link complement.\index{Whitehead link}
\end{example}

The fundamental theorem of Wallace and Lickorish, \refthm{WallaceLickorish}, implies that any closed orientable 3-manifold is obtained by Dehn filling a link complement in $S^3$. In fact we may take that link complement to be hyperbolic, due to work of Myers \cite{Myers}. Thus the hyperbolic Dehn filling theorem implies that in some sense, ``almost all'' 3-manifolds are hyperbolic.

There are still many unanswered questions about hyperbolic Dehn filling space.\index{hyperbolic Dehn filling space}\index{Dehn filling space} As of the writing of this book, the following questions are all unknown. 

\begin{question}\label{Question:HypDehnSurgery}
What is the topology of hyperbolic Dehn filling space?\index{hyperbolic Dehn filling space}\index{Dehn filling space} For example, is it connected? Is it path connected? That is, if a finite volume manifold $M(s)$ admits a complete hyperbolic structure, and if $M$ also admits a complete hyperbolic structure, is there necessarily a deformation of the hyperbolic structure running from the complete structure on $M$ to the complete structure on $M(s)$?

Stronger: If $M(s)$ admits a complete hyperbolic structure, and $M$ admits a complete hyperbolic structure, can we deform the hyperbolic structure on $M$ through cone manifolds\index{hyperbolic cone manifold}\index{cone manifold} with cone angles\index{cone angle} increasing monotonically from $0$ (at the complete structure on $M$) to $2\pi$ (at the complete structure on $M(s)$)?
\end{question}

As of the writing of this book, we do not even know if hyperbolic Dehn filling space is connected for the simplest of examples --- the figure-8 knot complement. The following example is discussed in \cite{Cooper-HK:orbifolds}.

\begin{example}[Dehn filling space for the figure-8 knot]\label{Example:HypDehnFillingFig8}

Thurston identified part of the boundary of the neighborhood about infinity separating hyperbolic Dehn fillings from non-hyperbolic ones. This is done on pages 58 through 61 of his notes \cite{thurston}. To determine these boundaries, he considers what is happening to the two hyperbolic structures on the tetrahedra as the values of their edge invariants approach the boundaries given by the gluing equations (the boundaries of the region in \reffig{ThurstonRegion}).  When both tetrahedra degenerate, the hyperbolic structure collapses and the limiting manifold is not hyperbolic.

However, when only one tetrahedron degenerates, we still have a hyperbolic structure for a little while. In this case, we will be gluing a positively oriented tetrahedron\index{positively oriented tetrahedron}\index{tetrahedron!positively oriented} to a negatively oriented one.\index{negatively oriented tetrahedron}\index{tetrahedron!negatively oriented} We can make sense of this by cutting the negatively oriented tetrahedron into pieces and subtracting them from the positively oriented one, leaving a polyhedron $P$. Faces of $P$ may then be identified to give a hyperbolic structure. No one knows exactly where
this stops working, although Hodgson's 1986 PhD thesis \cite{hodgson:thesis} gives evidence that the boundary should be as shown in \reffig{Fig8Boundary}.

\begin{figure}
\begin{center}
%% Accompanies image file 'F6-03-Fig8DF.eps' (pdf, eps, ps)
\begingroup%
  \makeatletter%
  \providecommand\color[2][]{%
    \errmessage{(Inkscape) Color is used for the text in Inkscape, but the package 'color.sty' is not loaded}%
    \renewcommand\color[2][]{}%
  }%
  \providecommand\transparent[1]{%
    \errmessage{(Inkscape) Transparency is used (non-zero) for the text in Inkscape, but the package 'transparent.sty' is not loaded}%
    \renewcommand\transparent[1]{}%
  }%
  \providecommand\rotatebox[2]{#2}%
  \newcommand*\fsize{\dimexpr\f@size pt\relax}%
  \newcommand*\lineheight[1]{\fontsize{\fsize}{#1\fsize}\selectfont}%
  \ifx\svgwidth\undefined%
    \setlength{\unitlength}{294.73746872bp}%
    \ifx\svgscale\undefined%
      \relax%
    \else%
      \setlength{\unitlength}{\unitlength * \real{\svgscale}}%
    \fi%
  \else%
    \setlength{\unitlength}{\svgwidth}%
  \fi%
  \global\let\svgwidth\undefined%
  \global\let\svgscale\undefined%
  \makeatother%
  \begin{picture}(1,0.44357278)%
    \lineheight{1}%
    \setlength\tabcolsep{0pt}%
    \put(0,0){\includegraphics[width=\unitlength]{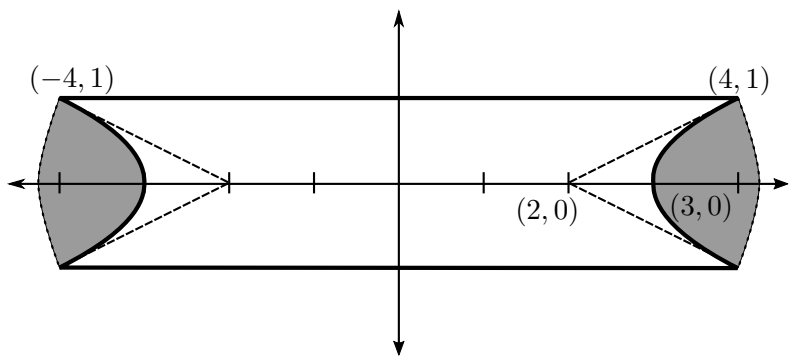}}%
    \put(0.02821436,0.34475965){\color[rgb]{0,0,0}\makebox(0,0)[lt]{\lineheight{0}\smash{\begin{tabular}[t]{l}$(-4,1)$\end{tabular}}}}%
    \put(0.89514089,0.34338876){\color[rgb]{0,0,0}\makebox(0,0)[lt]{\lineheight{0}\smash{\begin{tabular}[t]{l}$(4,1)$\end{tabular}}}}%
    \put(0.65017025,0.17573377){\color[rgb]{0,0,0}\makebox(0,0)[lt]{\lineheight{0}\smash{\begin{tabular}[t]{l}$(2,0)$\end{tabular}}}}%
    \put(0.84688663,0.18053197){\color[rgb]{0,0,0}\makebox(0,0)[lt]{\lineheight{0}\smash{\begin{tabular}[t]{l}$(3,0)$\end{tabular}}}}%
  \end{picture}%
\endgroup%

\end{center}
\caption{Hyperbolic Dehn filling space\index{hyperbolic Dehn filling space}\index{Dehn filling space} for the figure-8 knot complement is known to include the unshaded region exterior to the dark curve shown, is conjectured to contain the two shaded regions, and is conjectured to contain no other points. Figure modified from \cite{Cooper-HK:orbifolds}}
\label{Fig:Fig8Boundary}
\end{figure}

In \refex{fig8}, you are asked to study how tetrahedra
degenerate in the figure-8 knot complement.
\end{example}

\begin{question}
What is the hyperbolic Dehn filling space\index{hyperbolic Dehn filling space}\index{Dehn filling space} for the figure-8 knot complement?
\end{question}

\begin{definition}\label{Def:Exceptional}
  Dehn fillings that do not yield a hyperbolic manifold are called \emph{exceptional}.\index{exceptional Dehn filling}\index{Dehn filling!exceptional}
\end{definition}

There are many interesting problems on exceptional Dehn fillings. We include an example, that as of the writing of this book is open.

It is known (and you can prove as an exercise) that no hyperbolic manifold can contain an embedded 2-sphere that does not bound a 3-ball. A manifold that contains such a 2-sphere is called \emph{reducible}\index{reducible 3-manifold}. If you start with a hyperbolic 3-manifold, perform Dehn filling,\index{Dehn filling} and obtain a reducible manifold, the Dehn filling is called \emph{reducible}.\index{reducible Dehn filling}

\begin{conjecture}[The cabling conjecture]\label{Conj:Cabling}
No hyperbolic knot complement admits a reducible\index{reducible 3-manifold} Dehn filling.\index{Dehn filling}\index{Cabling conjecture}
\end{conjecture}

The original wording of the cabling conjecture is that only cables of knots admit reducible\index{reducible 3-manifold} Dehn fillings.\index{Dehn filling} The conjecture listed as \refconj{Cabling} is the remaining case to prove.

%%%%%%%%%%%%%%%%%%%%%%%%%%%%%%%%%%%%%%%%%%%%%%%%%%%%%%%%%%%%%%%%%
\subsection{Triangulations and Dehn filling}

When $X$ is a hyperbolic 3-manifold that admits an ideal triangulation, then Dehn filling\index{Dehn filling} of $X$ can frequently be performed by adjusting the edge invariants, as in \refdef{EdgeInvariant}, of the ideal tetrahedra making up $X$. That is, given a triangulation of a 3-manifold $X$ with torus boundary, we may solve a non-linear system of equations in the tetrahedra's edge parameters to find a hyperbolic structure on a Dehn filling of $X$. To do so, use the edge gluing equations of \refchap{GluingCompleteness}, but not the completeness equations.

Carefully, let $\mu$ and $\lambda$ be generators of $H_1(\bdy X)$, with associated completeness equations $H([\mu]) = H([\lambda])=1$, with $H([\mu]) = \prod_j z_{i_j}$ as in \refdef{CompletenessEquations}, and similarly for $H([\lambda])$.

Let $s=p\mu + q\lambda \in H_1(\bdy X)$ be the slope of the Dehn filling.\index{Dehn filling} To find a complete hyperbolic structure on $X(s)$, we solve the system of equations consisting of edge gluing equations and the \emph{Dehn filling equation}\index{Dehn filling equation}
\begin{equation}\label{Eqn:DehnFillingEquation}
p \log H([\mu]) + q \log H([\lambda]) = 2\pi i.
\end{equation}
Note a solution to these equations will produce an incomplete hyperbolic structure on $X$, with Dehn filling coefficient $(p,q) \in H_1(\bdy X;\RR)$. In fact, this process is valid for any $(p,q)\in \RR\oplus\RR \cong H_1(\bdy X;\RR)$, not just relatively prime integers.

The system of edge gluing equations along with \refeqn{DehnFillingEquation} may not have a solution. If it does have a solution, it may not be the case that all tetrahedron parameters have positive imaginary part. For such a solution, the corresponding tetrahedra are not all positively oriented; some are negatively oriented as well.\index{negatively oriented tetrahedron}\index{tetrahedron!negatively oriented}

In practice, it is possible to implement this process by computer, to find solutions to gluing and Dehn filling equations numerically, and this has been implemented in SnapPy \cite{SnapPy}. Indeed, in 2009, Schleimer and Segerman investigated hyperbolic Dehn filling space\index{hyperbolic Dehn filling space}\index{Dehn filling space} for thousands of manifolds, using SnapPy \cite{SchleimerSegerman}. Graphically, they identified regions of hyperbolic Dehn filling space for which SnapPy computed positively oriented tetrahedra,\index{positively oriented tetrahedron}\index{tetrahedron!positively oriented} negatively oriented tetrahedra,\index{negatively oriented tetrahedron}\index{tetrahedron!negatively oriented} and degenerate tetrahedra, as well as regions for which no solution was found.

\begin{figure}
  \includegraphics{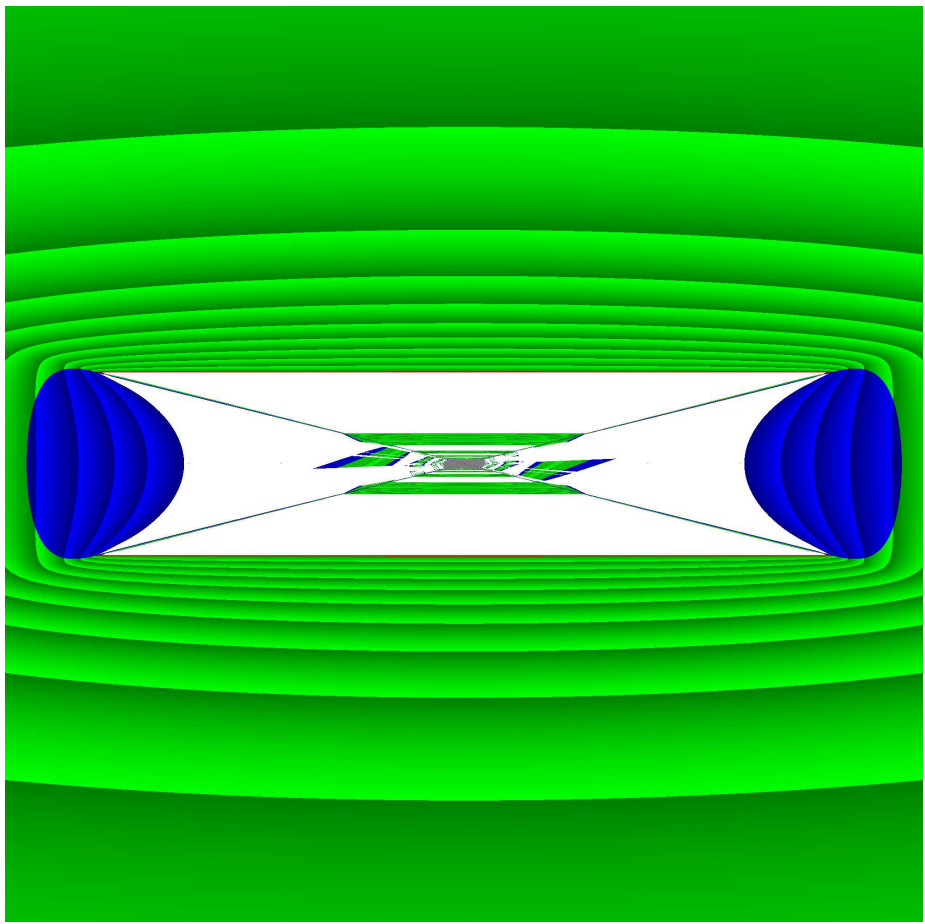}
\caption{Computer generated picture of hyperbolic Dehn filling space\index{hyperbolic Dehn filling space}\index{Dehn filling space} for the figure-8 knot complement, generated by Schleimer and Segerman. Compare with the conjectural picture of the space, \reffig{Fig8Boundary}}
\label{Fig:SchleimerSegerman}
\end{figure}

Schleimer and Segerman's computed space for the figure-8 knot is shown in \reffig{SchleimerSegerman}.
Green regions are those for which the computer found a solution with all positively oriented tetrahedra.\index{positively oriented tetrahedron}\index{tetrahedron!positively oriented} Blue regions are those for which the computer only found solutions with some negatively oriented tetrahedra.\index{negatively oriented tetrahedron}\index{tetrahedron!negatively oriented} In the white regions, the computer failed to recognize a solution. The gray region around the origin is where no solution was found. The shading in green and blue regions corresponds to volume; lines in those regions are level sets of volume. It is difficult to see in the printed version, but there is also a thin red line between green and white regions. Red indicates that all tetrahedra are flat\index{flat tetrahedron}\index{tetrahedron!flat} and non-degenerate,\index{degenerate tetrahedron}\index{tetrahedron!degenerate} i.e.\ the cross-ratio of the four ideal points for each tetrahedron is real, but bounded away from $0$, $1$, and $\infty$. Any point where at least one tetrahedron is degenerate, i.e.\ its cross-ratio is near $0$, $1$, or $\infty$, would be shaded purple. 

Note that the green and blue regions away from the origin in \reffig{SchleimerSegerman} match the conjectured picture for Dehn filling space in \reffig{Fig8Boundary}. The blue and green regions in the interior are conjectured to be noise, and not to correspond to actual hyperbolic structures.

\begin{figure}
  \includegraphics{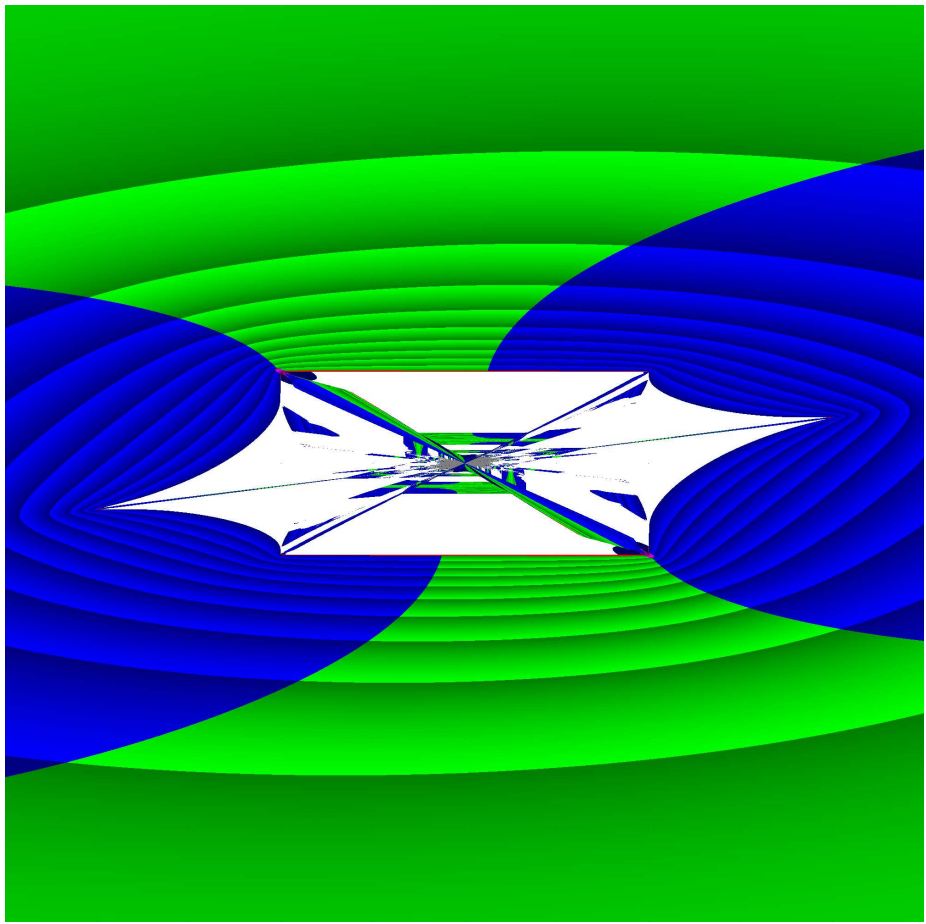}
\caption{Computer generated picture of hyperbolic Dehn filling space\index{hyperbolic Dehn filling space}\index{Dehn filling space} for the $5_2$ knot complement, generated by Schleimer and Segerman}
\label{Fig:52SchleimerSegerman}
\end{figure}

\begin{figure}
  \includegraphics{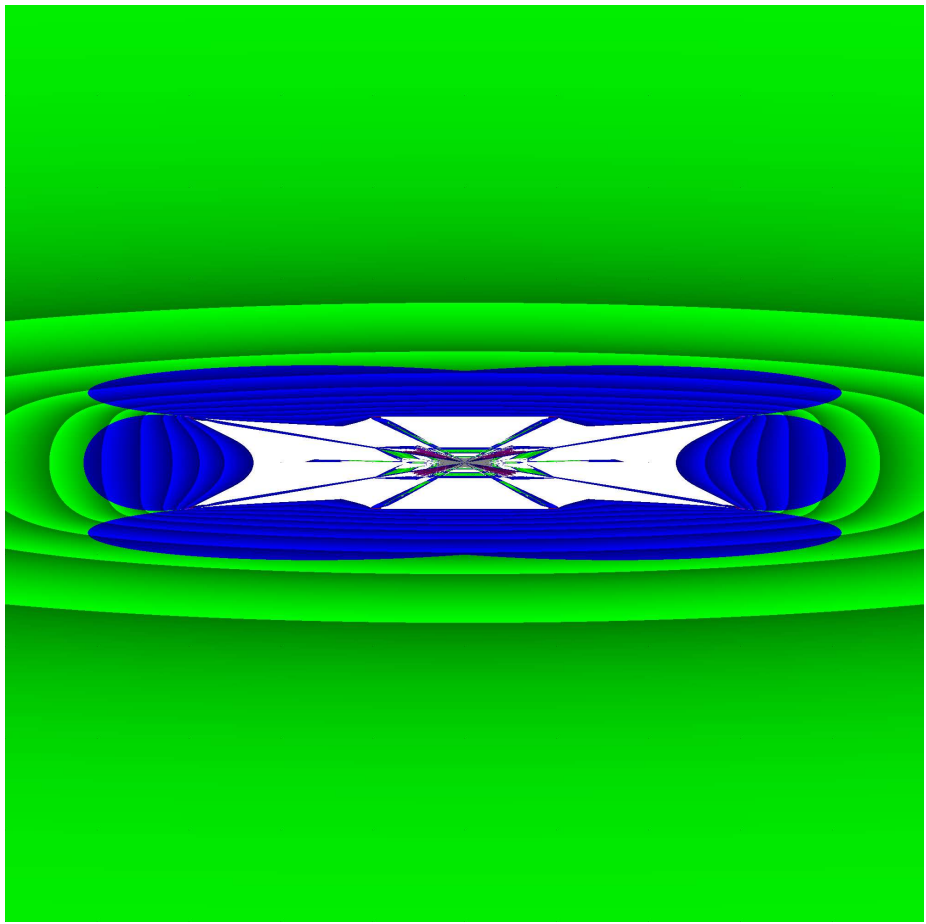}
\caption{Computer generated picture of hyperbolic Dehn filling space\index{hyperbolic Dehn filling space}\index{Dehn filling space} for the $6_3$ knot complement, generated by Schleimer and Segerman}
\label{Fig:63SchleimerSegerman}
\end{figure}

In addition, we include Schleimer and Segerman's images of hyperbolic Dehn filling space\index{hyperbolic Dehn filling space}\index{Dehn filling space} for the $5_2$ knot, and for the $6_3$ knot, in \reffig{52SchleimerSegerman} and \reffig{63SchleimerSegerman}, respectively.
The color scheme is the same as for the figure-8 knot, above. The triangulation used in these cases is known as the \emph{canonical triangulation}\index{canonical triangulation}, which will be defined in \refchap{TwoBridge}. Using different triangulations can lead to different regions of negatively oriented triangles. However, the boundary between hyperbolic and non-hyperbolic structures (green or blue versus white regions) seems to be independent of choice of triangulation. 
These figures and many more can be found on Segerman's website \cite{SchleimerSegerman}.

%%%%%%%%%%%%%%%%%%%%%%%%%%%%%%%%%%%%%%%%%%%%%%%%%%%%%%%%%%%%%%%%%

\section{A brief summary of geometric convergence}

\Refthm{HypDehnSurgery}, the hyperbolic Dehn filling theorem, actually gives information on convergence of geometry of spaces. The 3-manifolds obtained by hyperbolic Dehn filling on a complete hyperbolic manifold $M$ are ``close'' geometrically to $M$. This statement can be made precise, and often explicit, which is very useful: if we can bound geometric quantities for $M$, then the fact that (certain) Dehn fillings are geometrically close often translates into a bound on the same geometric quantities for Dehn fillings.

In this section, we define convergence of spaces and state a stronger version of the hyperbolic Dehn filling theorem that includes such convergence. We also survey briefly a few results and consequences of these results.

\subsection{Convergence of spaces}
Given two abstract metric spaces, we need a way to measure distance between them, and to describe when a sequence of spaces converges to another space. 
Convergence of metric spaces has been studied by Gromov \cite{Gromov}. In the case of hyperbolic spaces, \cite{CanaryEpsteinGreen:Notes}, \cite[Chapter~E]{benedetti-petronio}, \cite[Chapter~8]{kapovich}, and \cite[Chapter~6]{Cooper-HK:orbifolds} give further details and examples. 

There are actually several different definitions of geometric convergence of metric spaces in the literature on hyperbolic 3-manifolds and cone manifolds,\index{hyperbolic cone manifold}\index{cone manifold} which can be confusing. However, many are equivalent; see for example~\cite[Theorem~3.2.9]{CanaryEpsteinGreen:Notes} and~\cite[Theorem~8.11]{kapovich}. We give one definition here, as well as some examples that motivate other equivalent definitions.

One way of measuring ``distance'' between spaces is via quasi-isometries.

\begin{definition}\label{Def:QuasiIsometry}
Let $X$ and $Y$ be metric spaces with distance functions $d_X$ and $d_Y$, respectively. For $K>1$ and $c>0$, a bijection $f\from X \to Y$ is a \emph{$(K,c)$-quasi-isometric embedding}\index{quasi-isometric embedding} if for all distinct points $x,y\in X$,
\[ \frac{1}{K} d_X(x,y) - c \leq d_Y(f(x),f(y)) \leq K d_X(x,y) + c.\]

Let $f\from X\to Y$ be a $(K,c)$-quasi-isometric embedding. We say $f$ is a \emph{$(K,c)$-quasi-isometry}\index{quasi-isometry} if there also exists a map $\overline{f}\from Y \to X$ that is a $(K,c)$-quasi-isometric embedding as well as an \emph{approximate inverse}:\index{approximate inverse}
\[ \mbox{for all $x\in X$ and $y\in Y$, } d_X(\overline{f}\circ f(x),x)\leq c \mbox{ and } d_Y(f\circ\overline{f}(y),y)\leq c. \]
\end{definition}

\begin{definition}\label{Def:GeomConvergence}
Let $\{X_n\}$ be a sequence of metric spaces, each with a basepoint $x_n\in X_n$. Let $X$ be a metric space with basepoint $x\in X$. Let $B_r(x_n)$ denote the set of points in $X_n$ with distance at most $r$ from $x_n$. Similarly, let $B_r(x)$ denote the set of points in $X$ with distance at most $r$ from $x$.

We say that the sequence $(X_n,x_n)$ converges to the metric space $(X,x)$ in the \emph{quasi-isometric topology}\index{quasi-isometric topology} (or \emph{Gromov--Hausdorff topology}\index{Gromov--Hausdorff topology} or \emph{geometric topology}\index{geometric topology}) if the following holds.

Suppose that for all $\epsilon>0$ and all $r>0$ there exists an integer $N$ such that if $n>N$, then there exists a $(1+1/n, \epsilon)$-quasi-isometry \[f_n\from B_r(x_n) \to B_r(x).\]

In this case, we also say that the space $X$ is a \emph{geometric limit}\index{geometric limit} of the sequence $X_n$. The spaces $X_n$ \emph{converge geometrically}\index{geometric convergence} to $X$.
\end{definition}

In other words, the spaces $X_n$ with basepoints $x_n$ converge in the quasi-isometric topology, or converge geometrically, if there are better and better quasi-isometries between larger and larger closed and bounded sets about the basepoints. The maps $f_n$ are becoming closer and closer to actual isometries on larger and larger compact sets. 

Notice that a basepoint, and compact sets around basepoints, feature prominently in the definition. The choice of basepoint does affect geometric limits, as the following example shows.

\begin{example}[A 2-dimensional geometric limit]
Suppose $S$ is a compact surface with genus three. Let $\gamma \subset S$ be a simple closed curve such that cutting $S$ along $\gamma$ yields two components: a genus one surface with one boundary component and a genus two surface with one boundary component. 
Let $X_n$ be the metric space obtained by giving $S$ a hyperbolic metric in which $\gamma$ has length $1/n$. Let $x_n$ be a basepoint that lies on the genus one side of $\gamma$ in $X_n$ and let $y_n$ be a basepoint that lies on the genus two side. See \reffig{GeomLimExample}.

\begin{figure}
  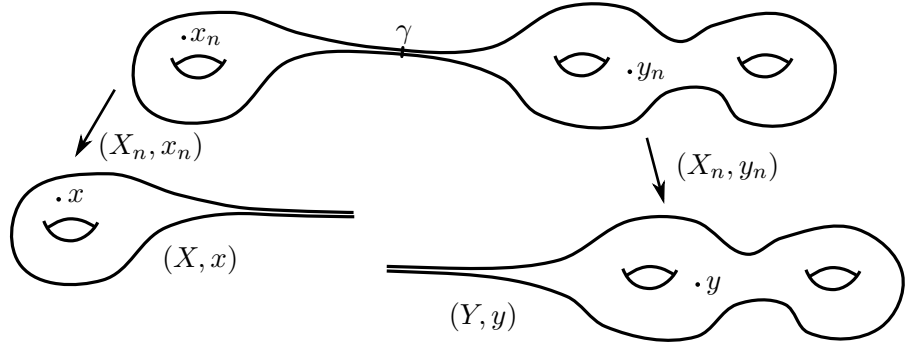
  \caption{Changing the basepoint from $x_n$ to $y_n$ can change the homeomorphism type of a geometric limit.}
  \label{Fig:GeomLimExample}
\end{figure}

Then $(X_n, x_n)$ has a geometric limit $(X,x)$ such that $X$ is homeomorphic to a torus with one cusp, as shown on the left of \reffig{GeomLimExample}. However, $(X_n,y_n)$ has a geometric limit $(Y,y)$ where $Y$ is homeomorphic to a genus two surface with one cusp, as shown on the right of \reffig{GeomLimExample}. Thus changing the basepoint can change the homeomorphism type of a geometric limit.
\end{example}

The following gives a 3-dimensional example of geometric convergence. 

\begin{example}[Geometric convergence of ideal tetrahedra]\label{Example:TetrahedraConvergence}
Consider a sequence of ideal tetrahedra $T_n$ in $\HH^3$ with vertices at $0$, $1$, $\infty$, and $z_n$, where $z_n$ is converging to some $z_\infty\in\CC$ with $\Im(z_\infty)>0$. For $n$ sufficiently large, there will be a point $p\in\HH^3$ in the interior of all tetrahedra $T_n$ and in the ideal tetrahedron $T_\infty$ with vertices at $0$, $1$, $\infty$, and $z_\infty$. For fixed $R>0$, consider the compact set given by taking the intersection of $T_n$ with a closed ball $B_R(p) \subset \HH^3$. 
For large $n$, there will be better and better quasi-isometries from the balls $B_R(p)\cap T_n$ to $B_R(p)\cap T_\infty$. It follows that the ideal tetrahedra $(T_n,p)$ converge geometrically to the ideal tetrahedron $(T_\infty,p)$.
\end{example}

\begin{example}[Polyhedral convergence]\label{Example:PolyhedralConvergence}
More generally, suppose $M$ is a complete hyperbolic 3-manifold obtained by gluing faces of an ideal polyhedron (possibly with infinite volume) in $\HH^3$, and suppose $M_n$ is another complete hyperbolic 3-manifold obtained by face-pairings\index{face-pairing isometry} of slightly deformed polyhedra, with face-pairings of $M_n$ pairing the same (combinatorial) faces as those of $M$. Suppose also that the polyhedra making up $M_n$ converge geometrically to those making up $M$ as $n\to\infty$, with appropriate basepoints. Then with the same basepoints, $(M_n, x_n)$ converges to $(M,x)$ in the quasi-isometric topology. This is made precise in \cite{marden}; convergence of spaces in this manner is called \emph{polyhedral convergence}.\index{polyhedral convergence}

Finally, taking the example one step farther, note that face-pairings\index{face-pairing isometry} are isometries in $\PSL(2,\CC)$, generating a discrete group $G_n$ when $M_n$ is hyperbolic. If $M$ is also hyperbolic, with face-pairings generating the discrete group $G$, we say that the sequence of groups $G_n$ converges to $G$ \emph{geometrically}.\index{geometric convergence of groups} These notions of convergence can be shown to be equivalent to convergence in the quasi-isometric topology; for example see \cite[theorem~8.11]{kapovich}.
\end{example}

Our main reason for defining geometric convergence is that it allows us to restate a much stronger version of Thurston's hyperbolic Dehn filling theorem, \refthm{HypDehnSurgery}, as follows.

\begin{theorem}[Hyperbolic Dehn filling with geomeric convergence]\label{Thm:GeomConvDehnFilling}\index{Thurston's hyperbolic Dehn filling theorem!geometric convergence}\index{hyperbolic Dehn filling theorem!geometric convergence}
  Let $M$ admit a complete hyperbolic structure with fixed horoball neighborhood of a cusp $C$. Let $s_n$ be a sequence of slopes on $\bdy C$ such that the length of a geodesic representative of $s_n$, measured in the induced Euclidean metric on $\bdy C$, approaches infinity. Then for large enough $n$, the Dehn filled manifolds $M(s_n)$ are hyperbolic and approach $M$ as a geometric limit.

  Similarly, if $M$ has multiple cusps, then Dehn filled manifolds along slopes with lengths approaching infinity approach $M$ as a geometric limit. 
\end{theorem}

\begin{proof}[Proof idea]
In the proof sketch of Thurston's hyperbolic Dehn filling theorem, \refthm{HypDehnSurgery}, we noted that incomplete hyperbolic structures on $M$ can be obtained by one-complex parameter families of deformations of the complete hyperbolic structure. These deformations are continuous maps in the quasi-isometric topology. 
\end{proof}

\subsection{Some consequences of geometric convergence}

The fact that $M$ is a geometric limit in \refthm{GeomConvDehnFilling} implies that geometric properties of Dehn fillings of $M$ converge to those of $M$. For example, the thick parts of a cusped finite-volume manifold and its high Dehn fillings will be quasi-isometric. A geodesic in the cusped manifold will map to curves that will eventually be isotopic to geodesics in the filled manifolds, with lengths approaching the length of the original. Unfortunately, \refthm{GeomConvDehnFilling} does not give any information on how high the Dehn fillings need to be in order to guarantee concrete bounds on geometry change. However, since the theorem appeared, there has been progress in making it more concrete.

For example, the volume of a finite volume hyperbolic 3-manifold $M$ is one of its most useful geometric properties. If 3-manifolds $M_n$ converge to $M$ as a geometric limit, then their volumes converge:
\[ \lim_{n\to\infty}\vol(M_n)\to \vol(M). \]

Much more can be said on volumes and Dehn filling.\index{Dehn filling} As a first step, the following theorem is also due to Thurston, and appears in the same notes in which he outlined the proof of \refthm{GeomConvDehnFilling}, as \cite[theorem~6.5.6]{thurston}.

\begin{theorem}[Volume under Dehn filling]\label{Thm:VolumeDF}\index{volume bound!upper bound on Dehn filling}
If $M$ is hyperbolic with cusp $C$, and $s$ is a slope on $\bdy C$ such that $M(s)$ is hyperbolic, then
\[  \vol(M) > \vol(M(s)). \]
Similarly if $M$ has multiple cusps $C_1,\dots, C_n$ and slopes $s_1, \dots, s_n$, one on each $C_j$, such that $M(s_1, \dots, s_n)$ is hyperbolic, then
\[ \vol(M)>\vol(M(s_1, \dots, s_n)). \]
\end{theorem}

While \refthm{GeomConvDehnFilling} implies that volumes of $M(s)$ approach volumes of $M$, \refthm{VolumeDF} implies that volume strictly decreases under Dehn filling for any slope giving a hyperbolic manifold. The slope need not be in the neighborhood of infinity provided by \refthm{HypDehnSurgery}; the volume decreases regardless. 

The full proof of \refthm{VolumeDF} can be found in \cite{thurston}. In this book, we will give a full proof of a slightly weaker result in \refchap{AngleStruct}, so we will delay the discussion of the proof ideas until then.

For volumes, even more can be said, and there have been concrete results bounding the change in volume under Dehn filling by Neumann and Zagier~\cite{NeumannZagier} and by Hodgson and Kerckhoff~\cite{hk:univ}, among others. We state one additional result along these lines here.

Note that if $M$ has a complete hyperbolic structure with cusp $C$, then $\bdy C$ has a Euclidean structure,\index{Euclidean structure} and any slope $s\subset \bdy C$ is isotopic to a geodesic with well-defined Euclidean length $\ell_{\bdy C}(s)$. Provided the length of $s$ is at least $2\pi$, a lower bound on volume under Dehn filling\index{Dehn filling} can also be obtained. We will prove the following theorem in \refchap{Volume}. 

\begin{theorem}[\cite{fkp:dfvjp}]\label{Thm:FKP}\index{volume bound!lower bound on Dehn filling}
Suppose $M$ is a hyperbolic manifold with cusps $C_1, \dots, C_n$ and slopes $s_1, \dots, s_n$, one on each $\bdy C_i$, such that the minimal length slope $\ell_{\min} = \min\{\ell_{\bdy C_j}(s_j)\}$ has length at least $2\pi$. Then the Dehn filled manifold $M(s_1, \dots, s_n)$ is hyperbolic with volume satisfying
  \[ \vol(M(s_1, \dots, s_n)) \geq \left( 1 - \left(\frac{2\pi}{\ell_{\min}}\right)^2 \right)^{3/2} \vol(M). \]
\end{theorem}

\section{Exercises}

\begin{exercise}\label{Ex:fig8} (Incomplete structures on the figure-8 knot)
Thurston's notes contain a figure showing all parameterizations of hyperbolic structures on the figure-8 knot \cite[page 52]{thurston}.  For any $w$ in this region, formula 4.3.2 in the notes gives us a corresponding $z$ so that if two tetrahedra with edge invariants $z$ and $w$ are glued, we obtain a (possibly incomplete) hyperbolic structure on the figure-8 knot.

Analyze what happens to the tetrahedra corresponding to $z$ and to $w$ as $w$ approaches a point on the boundary of this region.

More specifically, if $w$ approaches certain points on the boundary of this region, tetrahedra corresponding to both $z$ and $w$ start to become degenerate.  Which points are these?  Prove that the two tetrahedra are becoming degenerate in this case.

As $w$ approaches other values on the boundary, only one of the tetrahedra degenerates.  Which points are these?  Prove that only one tetrahedron is degenerating in this case.
\end{exercise}

\begin{exercise}
We have seen that the completion of an incomplete hyperbolic 3-manifold is no longer homeomorphic to the original hyperbolic 3-manifold.  Is this true for completions of incomplete structures on the 3-punctured sphere?  What surface do we obtain when we complete an incomplete hyperbolic structure on a 3-punctured sphere?  Prove it.
\end{exercise}

\begin{exercise}\label{Ex:zcrossz}
Suppose $M$ is a closed manifold with a complete hyperbolic structure. Prove that $\pi_1(M)$ cannot contain a $\ZZ \times \ZZ$ subgroup. Conclude that $M$ cannot contain an embedded torus $T$ such that $\pi_1(T)$ injects into $\pi_1(M)$.  [Such a torus is called \emph{incompressible}\index{incompressible torus}.  A Dehn filling\index{Dehn filling} resulting in a closed manifold with an embedded incompressible torus is another example of an exceptional filling.]\index{exceptional Dehn filling}\index{Dehn filling!exceptional}
\end{exercise}

\begin{exercise}
Let $M$ be an orientable 3-manifold with a decomposition into ideal polyhedra, each with a hyperbolic structure, such that the polyhedra induce a hyperbolic structure on $M$.  Let $v$ be an ideal vertex of $M$, i.e.\ an equivalence class of ideal vertices of the polyhedra, where vertices are equivalent if and only if they are identified under the gluing of the polyhedra.

Recall that ${\rm link}(v)$ is defined to be the boundary of a neighborhood of $v$ in $M$.
\begin{enumerate}
\item[(a)] Prove $\link(v)$ always inherits a similarity structure from the hyperbolic structure on $M$.  Here a similarity structure is a $({\rm Sim}(\EE^2), \EE^2)$-structure, where ${\rm Sim}(\EE^2)$ is a subgroup of the group of affine transformations consisting of elements of the form $x\mapsto Ax+b$, where $A$ is a linear map that rotates and/or scales only.  Thus ${\rm Sim}(\EE^2)$ is formed by rotations, scalings, and translations.
\item[(b)] Prove that the only closed, orientable surface which admits a similarity structure is a torus.  It follows that $\link(v)$ is always homeomorphic to a torus when $M$ is an orientable manifold with hyperbolic structure (even incomplete).
\end{enumerate}
\label{Ex:link}
\end{exercise}

\begin{exercise}\label{Ex:1PtCmpt}
Let $M$ be an orientable 3-manifold that admits an incomplete hyperbolic structure with completion given by attaching the one-point compactification of a cusp neighborhood $N(C)$. Prove that the completion is not a manifold.
\end{exercise}

\begin{exercise}
Prove that a reducible\index{reducible 3-manifold} manifold cannot be hyperbolic. That is, it admits no complete hyperbolic structure.
\end{exercise}

\begin{exercise}\label{Ex:LensSpaceNotHyp}
Prove that a lens space\index{lens space} cannot admit a complete hyperbolic structure.
\end{exercise}

\begin{exercise}
(On complex length of $A$ in $\PSL(2,\CC)$)

\begin{enumerate}
\item[(a)] Suppose $A \in \PSL(2,\CC)$ has axis the geodesic from $0$ to $\infty$.  Then the matrix of $A$ may be parameterized by a single complex number $\lambda$.  What is the form of this matrix?
\item[(b)] Denote the trace of a matrix $A$ by ${\rm tr}(A)$, and its complex length\index{complex length} by $\calL(A)$.  Prove that ${\rm tr}(A) = 2\cosh(\calL(A)/2)$.
\end{enumerate}
\end{exercise}

\begin{exercise}
By computer, investigate the hyperbolic Dehn filling space\index{hyperbolic Dehn filling space}\index{Dehn filling space} obtained by filling a single component of the Whitehead link.\index{Whitehead link} Identify regions for which the result has a decomposition into positively oriented tetrahedra, negatively oriented tetrahedra, etc. 
\end{exercise}

\begin{exercise}
  (Algebraic versus geometric convergence) The purpose of this exercise is to work through a basic example of Thurston \cite{thurston} showing a difference between algebraic and geometric convergence of discrete groups. Let $A_n \in \PSL(2,\CC)$ have matrix representation
  \[ A_n = \mat{\exp(w_n) & n\sinh(w_n)\\ 0& \exp(-w_n)} \quad \mbox{where} \quad
  w_n = \frac{1}{n^2} + i\frac{\pi}{n}.\]
  \begin{enumerate}
  \item Show that the matrices $A_n$ converge to the matrix $A= \mat{1&i\pi\\0&1}$; thus the group $\langle A_n \rangle$ converges algebraically to the group $\langle A \rangle$.
  \item Show that $\langle A_n \rangle$ does not converge geometrically to the group $\langle A \rangle$, by finding a subsequence $A_{n_j}$ converging to an element of $\PSL(2,\CC)$ that does not lie in $\langle A \rangle$.
  \end{enumerate}
\end{exercise}

%    Part II
\part{Tools, Techniques, and Families of Examples}\label{Part:Examples}

\chapter{Twist Knots and Augmented Links}\label{Chap:TwistKnots}

In this chapter, \blfootnote{Jessica S. Purcell, Hyperbolic Knot Theory}
we study a class of knots that have some of the simplest hyperbolic geometry, namely twist knots. This class includes the figure-8 knot, the $5_2$ knot, and the $6_1$ knot that we have encountered so far. We also generalize to give examples of knots and links whose geometry is relatively explicit. This will equip us with many examples.

From now on, we say a knot or link in $S^3$ is \emph{hyperbolic}\index{hyperbolic knot or link} if its complement $S^3-K$ admits a complete hyperbolic structure. Similarly, a \emph{hyperbolic 3-manifold}\index{hyperbolic 3-manifold} is a 3-manifold that admits a complete hyperbolic structure. Note that the completeness of the hyperbolic structure is implied in this terminology.

\section{Twist knots and Dehn fillings}

Recall the definition of twist knots from \refchap{KnotIntro}. 

\begin{figure}[h]
  \includegraphics{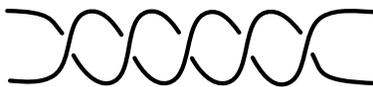}
  \caption{A twist region of a diagram}
  \label{Fig:TwistRegion2}
\end{figure}

A twist region\index{twist region} is a string of bigon\index{bigon} regions in the diagram graph of a knot diagram, with the bigons arranged end-to-end at their vertices, as in \reffig{TwistRegion2}. Recall also that a twist region is maximal in the sense that there are no additional bigon regions meeting the vertices on either end. A single crossing adjacent to no bigons is also a twist region. Recall also that twist regions are required to be alternating.

The \emph{twist knot} $J(2,n)$\index{twist knot}, defined in \refdef{TwistKnot}, is the knot with a diagram consisting of exactly two twist regions, one of which contains two crossings, and the other containing $n\in\ZZ$ crossings. The direction of crossing depends on the sign of $n$.

\begin{figure}
  \includegraphics{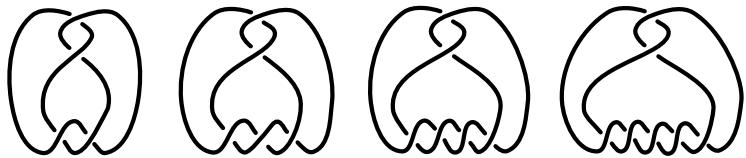}
  \caption{Twist knots $J(2,2)$ (the figure-8 knot), $J(2,3)$ (the $5_2$ knot), $J(2,4)$ (the $6_1$ or Stevedore knot), and $J(2,5)$}
  \label{Fig:TwistKnots2}
\end{figure}

Twist knots $J(2,2)$, $J(2,3)$, $J(2,4)$, and $J(2,5)$ are shown again in \reffig{TwistKnots2}.

\begin{figure}
  \includegraphics{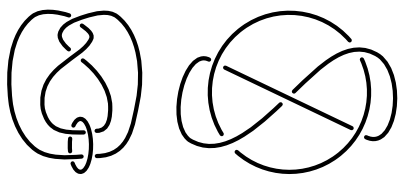}
  \caption{Two diagrams of the Whitehead link.\index{Whitehead link}}
  \label{Fig:Whitehead}
\end{figure}

\begin{definition}\label{Def:WhiteheadLink}
  The \emph{Whitehead link}\index{Whitehead link} is the link shown in \reffig{Whitehead}. Note the two links shown are isotopic. 
\end{definition}

We will show in \refprop{WhiteheadGeom} that the complement of the Whitehead link\index{Whitehead link} is hyperbolic. 

\begin{proposition}\label{Prop:TwistKnotWhitehead}
The complement of the twist knot $J(2,n)$ is obtained by Dehn filling\index{Dehn filling} the hyperbolic manifold isometric to the complement of the Whitehead link.\index{Whitehead link}
\end{proposition}

\begin{proof}
The proof uses topological properties of the sphere $S^3$ and the solid torus. Recall first that the sphere $S^3$ is the union of two solid tori whose cores are linked exactly once, but each core alone is unknotted.

The diagram of the Whitehead link\index{Whitehead link} on the left of \reffig{Whitehead} has a component at the bottom that is unknotted and does not cross itself. The complement of this component in $S^3$ is a solid torus. Note then that the other component is a knot in a solid torus, as shown on the left of \reffig{TKWProof}.

\begin{figure}
%% Accompanies image file 'F7-04-TKWpf.eps' (pdf, eps, ps)
\begingroup%
  \makeatletter%
  \providecommand\color[2][]{%
    \errmessage{(Inkscape) Color is used for the text in Inkscape, but the package 'color.sty' is not loaded}%
    \renewcommand\color[2][]{}%
  }%
  \providecommand\transparent[1]{%
    \errmessage{(Inkscape) Transparency is used (non-zero) for the text in Inkscape, but the package 'transparent.sty' is not loaded}%
    \renewcommand\transparent[1]{}%
  }%
  \providecommand\rotatebox[2]{#2}%
  \newcommand*\fsize{\dimexpr\f@size pt\relax}%
  \newcommand*\lineheight[1]{\fontsize{\fsize}{#1\fsize}\selectfont}%
  \ifx\svgwidth\undefined%
    \setlength{\unitlength}{311.85595322bp}%
    \ifx\svgscale\undefined%
      \relax%
    \else%
      \setlength{\unitlength}{\unitlength * \real{\svgscale}}%
    \fi%
  \else%
    \setlength{\unitlength}{\svgwidth}%
  \fi%
  \global\let\svgwidth\undefined%
  \global\let\svgscale\undefined%
  \makeatother%
  \begin{picture}(1,0.41817061)%
    \lineheight{1}%
    \setlength\tabcolsep{0pt}%
    \put(0,0){\includegraphics[width=\unitlength]{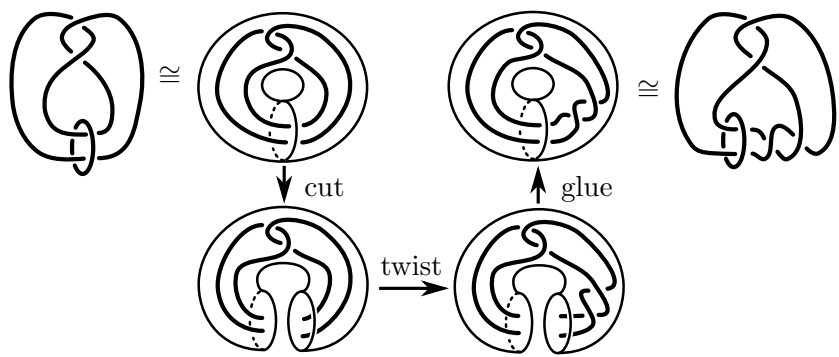}}%
    \put(0.18031699,0.32969486){\color[rgb]{0,0,0}\makebox(0,0)[lt]{\lineheight{0}\smash{\begin{tabular}[t]{l}$\cong$\end{tabular}}}}%
    \put(0.35694935,0.19551739){\color[rgb]{0,0,0}\makebox(0,0)[lt]{\lineheight{0}\smash{\begin{tabular}[t]{l}cut\end{tabular}}}}%
    \put(0.44894158,0.09963823){\color[rgb]{0,0,0}\makebox(0,0)[lt]{\lineheight{0}\smash{\begin{tabular}[t]{l}twist\end{tabular}}}}%
    \put(0.66707038,0.19560006){\color[rgb]{0,0,0}\makebox(0,0)[lt]{\lineheight{0}\smash{\begin{tabular}[t]{l}glue\end{tabular}}}}%
    \put(0.75892262,0.31362925){\color[rgb]{0,0,0}\makebox(0,0)[lt]{\lineheight{0}\smash{\begin{tabular}[t]{l}$\cong$\end{tabular}}}}%
  \end{picture}%
\endgroup%

  \caption{The Whitehead link\index{Whitehead link} complement is homeomorphic to a knot in a solid torus, which we cut, twist, and reglue. The result is homeomorphic to the complement of $J(2,2)\cup U$}
  \label{Fig:TKWProof}
\end{figure}

Now we apply a homeomorphism to the solid torus, which we view as $S^1\times D^2$. There is a homeomorphism given by slicing along a disk $\{x\}\times D^2$ of the solid torus, rotating one full time, then gluing back together. This homeomorphism is shown in the center of \reffig{TKWProof}.

The homeomorphism replaces the original link in the solid torus by a link with two additional crossings. By applying the homeomorphism repeatedly, we see that the complement of the Whitehead link\index{Whitehead link} is homeomorphic to the complement of the link with any even number of crossings encircled by the unknotted component. In particular, it is homeomorphic to the complement of the link $J(2,2k)\cup U$, where $U$ is a single unknotted component. By the Mostow--Prasad rigidity\index{Mostow--Prasad rigidity} theorem (\refthm{MostowGeom}), these link complements have isometric hyperbolic structures.

To obtain the knot $J(2,2k)$, attach a solid torus to $S^3-(J(2,2k)\cup U)$, filling in $U$ in a trivial way to give $S^3-J(2,2k)$. Thus $J(2,2k)$ is obtained from a manifold isometric to the complement of the Whitehead link\index{Whitehead link} by Dehn filling.\index{Dehn filling}

So far our proof only works for $J(2,n)$ with $n$ even. Now we consider the case of the knot $J(2,2k+1)$, with odd second component. We may isotope the Whitehead link,\index{Whitehead link} starting with the diagram on the left of \reffig{Whitehead}, to reverse the two crossings at the top, and insert a crossing encircled by the unknotted component at the bottom. This is shown in \reffig{TKWOdd}, left.

\begin{figure}
%% Accompanies image file 'F7-05-TKWOdd.eps' (pdf, eps, ps)
\begingroup%
  \makeatletter%
  \providecommand\color[2][]{%
    \errmessage{(Inkscape) Color is used for the text in Inkscape, but the package 'color.sty' is not loaded}%
    \renewcommand\color[2][]{}%
  }%
  \providecommand\transparent[1]{%
    \errmessage{(Inkscape) Transparency is used (non-zero) for the text in Inkscape, but the package 'transparent.sty' is not loaded}%
    \renewcommand\transparent[1]{}%
  }%
  \providecommand\rotatebox[2]{#2}%
  \newcommand*\fsize{\dimexpr\f@size pt\relax}%
  \newcommand*\lineheight[1]{\fontsize{\fsize}{#1\fsize}\selectfont}%
  \ifx\svgwidth\undefined%
    \setlength{\unitlength}{322.14482117bp}%
    \ifx\svgscale\undefined%
      \relax%
    \else%
      \setlength{\unitlength}{\unitlength * \real{\svgscale}}%
    \fi%
  \else%
    \setlength{\unitlength}{\svgwidth}%
  \fi%
  \global\let\svgwidth\undefined%
  \global\let\svgscale\undefined%
  \makeatother%
  \begin{picture}(1,0.21486269)%
    \lineheight{1}%
    \setlength\tabcolsep{0pt}%
    \put(0,0){\includegraphics[width=\unitlength]{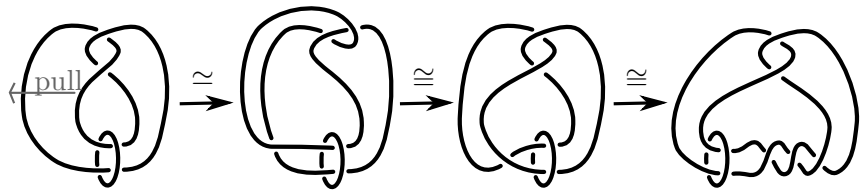}}%
    \put(0.0325797,0.11879199){\color[rgb]{0.39215686,0.39215686,0.39215686}\makebox(0,0)[lt]{\lineheight{0}\smash{\begin{tabular}[t]{l}pull\end{tabular}}}}%
    \put(0.21604048,0.12387844){\color[rgb]{0,0,0}\makebox(0,0)[lt]{\lineheight{0}\smash{\begin{tabular}[t]{l}$\simeq$\end{tabular}}}}%
    \put(0.47503741,0.12199698){\color[rgb]{0,0,0}\makebox(0,0)[lt]{\lineheight{0}\smash{\begin{tabular}[t]{l}$\cong$\end{tabular}}}}%
    \put(0.72387716,0.12074272){\color[rgb]{0,0,0}\makebox(0,0)[lt]{\lineheight{0}\smash{\begin{tabular}[t]{l}$\cong$\end{tabular}}}}%
  \end{picture}%
\endgroup%

  \caption{A sequence of homeomorphisms of the Whitehead link\index{Whitehead link} complement}
  \label{Fig:TKWOdd}
\end{figure}

Following that figure, we may then reflect the diagram in the plane of projection, reversing all the crossings. This is a homeomorphism of the knot complement, hence an isometry. Now just as in the even case, we may insert any even number of crossings into the two strands encircled by the unknotted component. To obtain $J(2,2k+1)$, simply Dehn fill\index{Dehn filling} the unknotted component in the obvious way. 
\end{proof}

\begin{corollary}\label{Cor:TwistKnotWhitehead}
  The complement of the Whitehead link\index{Whitehead link} is a geometric limit of $S^3-J(2,n)$.
\end{corollary}

\begin{proof}
Because they are obtained by Dehn filling\index{Dehn filling} the complement of the Whitehead link,\index{Whitehead link} all but finitely many link complements $S^3-J(2,n)$ lie in any given neighborhood of infinity in the Dehn surgery space for a cusp of the complement of the Whitehead link. \Refthm{GeomConvDehnFilling} implies that the Whitehead link is therefore a geometric limit of these manifolds. 
\end{proof}

%%%%%%%%%%%%%%%%%%%%%%%%%%%%%%%%%%%%%%%%%%%%%%%%%%%%%%%%%%%%%%%%%

In order to study the geometry of twist knots, we study the geometry of the geometric limit, the Whitehead link complement.\index{Whitehead link}

\begin{figure}
%% Accompanies image file 'F7-06-Oct.eps' (pdf, eps, ps)
\begingroup%
  \makeatletter%
  \providecommand\color[2][]{%
    \errmessage{(Inkscape) Color is used for the text in Inkscape, but the package 'color.sty' is not loaded}%
    \renewcommand\color[2][]{}%
  }%
  \providecommand\transparent[1]{%
    \errmessage{(Inkscape) Transparency is used (non-zero) for the text in Inkscape, but the package 'transparent.sty' is not loaded}%
    \renewcommand\transparent[1]{}%
  }%
  \providecommand\rotatebox[2]{#2}%
  \newcommand*\fsize{\dimexpr\f@size pt\relax}%
  \newcommand*\lineheight[1]{\fontsize{\fsize}{#1\fsize}\selectfont}%
  \ifx\svgwidth\undefined%
    \setlength{\unitlength}{96.56466579bp}%
    \ifx\svgscale\undefined%
      \relax%
    \else%
      \setlength{\unitlength}{\unitlength * \real{\svgscale}}%
    \fi%
  \else%
    \setlength{\unitlength}{\svgwidth}%
  \fi%
  \global\let\svgwidth\undefined%
  \global\let\svgscale\undefined%
  \makeatother%
  \begin{picture}(1,0.99666533)%
    \lineheight{1}%
    \setlength\tabcolsep{0pt}%
    \put(0,0){\includegraphics[width=\unitlength]{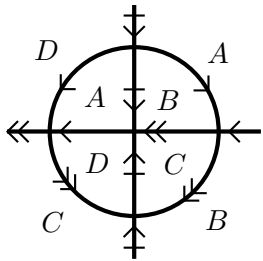}}%
    \put(0.30628015,0.6012646){\color[rgb]{0,0,0}\makebox(0,0)[lt]{\lineheight{0}\smash{\begin{tabular}[t]{l}$A$\end{tabular}}}}%
    \put(0.58787258,0.58871128){\color[rgb]{0,0,0}\makebox(0,0)[lt]{\lineheight{0}\smash{\begin{tabular}[t]{l}$B$\end{tabular}}}}%
    \put(0.30837229,0.34268627){\color[rgb]{0,0,0}\makebox(0,0)[lt]{\lineheight{0}\smash{\begin{tabular}[t]{l}$D$\end{tabular}}}}%
    \put(0.61716304,0.34059388){\color[rgb]{0,0,0}\makebox(0,0)[lt]{\lineheight{0}\smash{\begin{tabular}[t]{l}$C$\end{tabular}}}}%
    \put(0.78619836,0.77323333){\color[rgb]{0,0,0}\makebox(0,0)[lt]{\lineheight{0}\smash{\begin{tabular}[t]{l}$A$\end{tabular}}}}%
    \put(0.10962593,0.78536595){\color[rgb]{0,0,0}\makebox(0,0)[lt]{\lineheight{0}\smash{\begin{tabular}[t]{l}$D$\end{tabular}}}}%
    \put(0.13849588,0.11422912){\color[rgb]{0,0,0}\makebox(0,0)[lt]{\lineheight{0}\smash{\begin{tabular}[t]{l}$C$\end{tabular}}}}%
    \put(0.77740905,0.12050527){\color[rgb]{0,0,0}\makebox(0,0)[lt]{\lineheight{0}\smash{\begin{tabular}[t]{l}$B$\end{tabular}}}}%
  \end{picture}%
\endgroup%

  \caption{Shown is the boundary of an ideal octahedron (one vertex at infinity). Pairing faces as shown gives the complement of the Whitehead link.\index{Whitehead link}}
  \label{Fig:Octahedron}
\end{figure}

\begin{proposition}\label{Prop:WhiteheadGeom}
The complete hyperbolic structure on the complement of the Whitehead link\index{Whitehead link} is obtained by gluing faces of a regular ideal octahedron,\index{regular ideal octahedron} with the face pairings as shown in \reffig{Octahedron}.
\end{proposition}

A \emph{regular ideal octahedron}\index{ideal octahedron, regular}\index{regular ideal octahedron} is the ideal octahedron in $\HH^3$ with all dihedral angles equal to $\pi/2$. 

\begin{proof}
The fact that the Whitehead link\index{Whitehead link} complement is obtained by face pairings of an ideal octahedron can be readily seen by applying the methods of \refchap{Fig8Decomp} to the diagram of the Whitehead link on the right of \reffig{Whitehead}. After collapsing bigons,\index{bigon} we obtain two ideal polyhedra with four triangular faces and one quadrilateral face. Glue the quadrilaterals to obtain an ideal octahedron.\index{ideal octahedron, regular}\index{regular ideal octahedron} The form is shown in \reffig{Octahedron}. We leave the details for \refex{WhiteheadOctahedron}.

In a regular ideal octahedron, all dihedral angles are $\pi/2$, so horospheres intersect a neighborhood of each ideal vertex in a square. We need to check that the face pairings give a hyperbolic structure in this case. Note first that every point in the interior of an octahedron and in the interior of a face of the octahedron has a neighborhood isometric to a ball in $\HH^3$. We need to show that each point on an edge also has such a neighborhood, and then \reflem{IsometricNbhd} will imply that the gluing is a manifold with a (possibly incomplete) hyperbolic structure.

Note first that each of the edges (there are three) is glued four times. Thus the total angle around each edge will be $4\pi/2 =2\pi$. This is not quite enough to show that each point on an edge has a neighborhood isometric to a ball in $\HH^3$, because composing the gluings around an edge may introduce nontrivial translation or scale. To show that this does not happen, consider each end of an ideal edge within a cusp. Any horosphere intersects a neighborhood of an ideal vertex of the regular ideal octahedron\index{ideal octahedron, regular}\index{regular ideal octahedron} in a Euclidean square. Under the developing map, squares can only patch together in squares to give a tiling of the universal cover of each cusp by Euclidean squares. There are four squares meeting around a vertex in the cusp corresponding to one of our ideal edges. Note that the squares cannot be scaled or sheared. It follows that edges glue up without shearing singularities, and the structure is hyperbolic.

To show that the structure is complete, we use \refthm{EuclidCusp}: the structure is complete if and only if for each cusp, the induced structure on the boundary is Euclidean. But as already noted, each cusp is tiled by Euclidean squares corresponding to intersections of a horosphere with an ideal vertex of the regular ideal octahedron.\index{ideal octahedron, regular}\index{regular ideal octahedron} Under the developing map, squares can only patch together to give a Euclidean structure:\index{Euclidean structure} there will be no rotation or scale. Thus the hyperbolic structure must be complete. 
\end{proof}

In \refchap{AngleStruct}, we will obtain a formula to calculate the volume of a regular hyperbolic ideal octahedron. For now, we state that the volume is a constant $\voct = 3.66...$.\index{ideal octahedron, regular}\index{regular ideal octahedron}\index{regular ideal octahedron! volume}

\begin{corollary}\label{Cor:TwistVolume}
The volume of a hyperbolic twist knot is universally bounded
  \[ \vol(S^3-J(2,n)) < \voct, \]
  and as $n\to \infty$, $\vol(S^3-J(2,n))\to \voct$. 
\end{corollary}

\begin{proof}
The Dehn filling\index{Dehn filling} bound follows immediately from Thurston's theorem on volume change under Dehn filling, \refthm{VolumeDF}. The convergence follows from \refthm{GeomConvDehnFilling}. 
\end{proof}

We have not yet discussed which twist knots are hyperbolic. We have seen that the figure-8 knot is hyperbolic, and similar methods can be used to show each of the knots in \reffig{TwistKnots} are hyperbolic. More generally, we will see in \refchap{Alternating} (or by other methods in \refchap{TwoBridge}) that all twist knots $J(2,n)$ with $n\geq 2$ or $n\leq -3$ are hyperbolic. When $n=1$ or $-2$, the standard diagram of $J(2,n)$ can be easily reduced to a diagram with only a single twist region, which is not hyperbolic, and when $n=-1$ its diagram can be easily reduced to that of the unknot, which is also not hyperbolic. All other twist knots are hyperbolic.

%%%%%%%%%%%%%%%%%%%%%%%%%%%%%%%%%%%%%%%%%%%%%%%%%%%%%%%%%%%%%%%%%
\section{Double twist knots and the Borromean rings}

The results of the previous section generalize immediately to knots and links with exactly two twist regions, but with any number of crossings in either twist region.

\begin{figure}
%% Accompanies image file 'F7-07-DoubTw.eps' (pdf, eps, ps)
\begingroup%
  \makeatletter%
  \providecommand\color[2][]{%
    \errmessage{(Inkscape) Color is used for the text in Inkscape, but the package 'color.sty' is not loaded}%
    \renewcommand\color[2][]{}%
  }%
  \providecommand\transparent[1]{%
    \errmessage{(Inkscape) Transparency is used (non-zero) for the text in Inkscape, but the package 'transparent.sty' is not loaded}%
    \renewcommand\transparent[1]{}%
  }%
  \providecommand\rotatebox[2]{#2}%
  \newcommand*\fsize{\dimexpr\f@size pt\relax}%
  \newcommand*\lineheight[1]{\fontsize{\fsize}{#1\fsize}\selectfont}%
  \ifx\svgwidth\undefined%
    \setlength{\unitlength}{70.34094858bp}%
    \ifx\svgscale\undefined%
      \relax%
    \else%
      \setlength{\unitlength}{\unitlength * \real{\svgscale}}%
    \fi%
  \else%
    \setlength{\unitlength}{\svgwidth}%
  \fi%
  \global\let\svgwidth\undefined%
  \global\let\svgscale\undefined%
  \makeatother%
  \begin{picture}(1,1.06616641)%
    \lineheight{1}%
    \setlength\tabcolsep{0pt}%
    \put(0,0){\includegraphics[width=\unitlength]{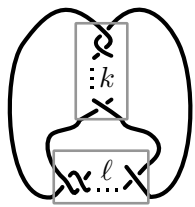}}%
    \put(0.48983503,0.6333513){\color[rgb]{0,0,0}\makebox(0,0)[lt]{\lineheight{0}\smash{\begin{tabular}[t]{l}$k$\end{tabular}}}}%
    \put(0.49929811,0.13093165){\color[rgb]{0,0,0}\makebox(0,0)[lt]{\lineheight{0}\smash{\begin{tabular}[t]{l}$\ell$\end{tabular}}}}%
  \end{picture}%
\endgroup%

  \caption{A double twist knot or link has two twist regions, one with $k$ crossings and one with $\ell$ crossings}
  \label{Fig:DoubleTwist}
\end{figure}

\begin{definition}\label{Def:DoubleTwistKnot}
  The \emph{double twist knot or link} $J(k,\ell)$ is the knot or link with a diagram consisting of exactly two twist regions, one of which contains $k$ crossings, and the other contains $\ell$ crossings, for $k, \ell\in\ZZ$. See \reffig{DoubleTwist}. Note that $J(k,\ell)$ is a knot if and only if at least one of $k, \ell$ is even; otherwise it is a link with two components.
\end{definition}

Just as for twist knots, double twist knots are obtained by Dehn filling\index{Dehn filling} a simple link complement. 

\begin{figure}
  \includegraphics{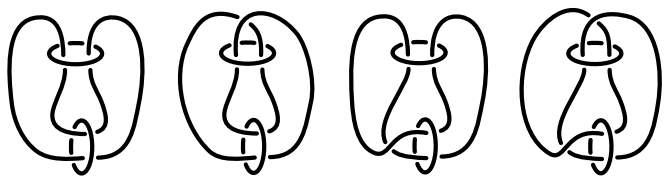}
  \caption{Complements of $J(k,\ell)$ are obtained by Dehn filling\index{Dehn filling} one of these four links. The link on the left is known as the Borromean rings.\index{Borromean rings}}
  \label{Fig:Borromean}
\end{figure}

\begin{proposition}\label{Prop:JklBorromean}
  The complement of the link $J(k,\ell)$ is obtained by Dehn filling\index{Dehn filling} the complement of one of the four links shown in \reffig{Borromean}, depending on the parity of $k$ and $\ell$.
\end{proposition}

\begin{proof}
The proof is nearly identical to that of \refprop{TwistKnotWhitehead}, except now it is done in two steps, since there are two unknotted components. Apply a homeomorphism of a solid torus as in \reffig{TKWOdd} two times. The details are left to the reader. 
\end{proof}

The link on the left of \reffig{Borromean} is equivalent to a link more famously known as the \emph{Borromean rings}\index{Borromean rings}; its more common diagram is shown in \reffig{2BorromeanRings}.
We will call the other links of \reffig{Borromean} the \emph{Borromean twisted sisters}\index{Borromean twisted sisters}, and say the links are in the \emph{Borromean family}\index{Borromean family}. In fact, the middle two links are equivalent.

\begin{proposition}\label{Prop:GeomBorromean}
The complements of the Borromean rings\index{Borromean rings} and the Borromean twisted sisters\index{Borromean twisted sisters} all admit complete hyperbolic structures obtained by gluing two regular ideal octahedra. 
\end{proposition}

\begin{proof}
Because the Borromean rings has a diagram that is alternating,\index{alternating knot or link} its complement can be split into ideal polyhedra using the methods of \refchap{Fig8Decomp}. However, we present a new way to decompose the complements of links of the Borromean family that we will generalize below.

View the diagrams of \reffig{Borromean} in three dimensions. The two link components in each diagram that will be Dehn filled\index{Dehn filling} to produce $J(k,\ell)$ should be viewed as lying perpendicular to the plane of the paper, which is the plane of projection $S^2\subset S^3$. The other link component(s) should be viewed as lying in the plane of projection except at crossings; when the component crosses itself it dips briefly above or below the plane of projection, then returns to the plane.

The components lying perpendicular to the plane of projection are unknotted, and each bounds a 2-punctured disk, shown as shaded in \reffig{Shaded2PunctDisk}.

\begin{figure}
  \includegraphics{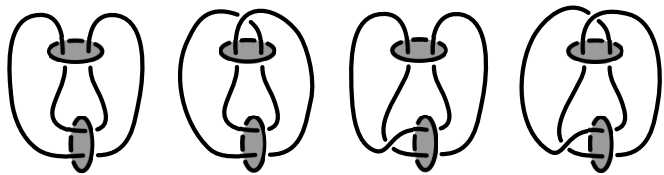}
  \caption{Shaded 2-punctured disks.}
  \label{Fig:Shaded2PunctDisk}
\end{figure}

As the first step of the decomposition, slice each of these disks up the middle, replacing a single 2-punctured disk with two parallel copies of the 2-punctured disk. This move is shown on the left of \reffig{BorromeanDecomp}.

\begin{figure}
  \includegraphics{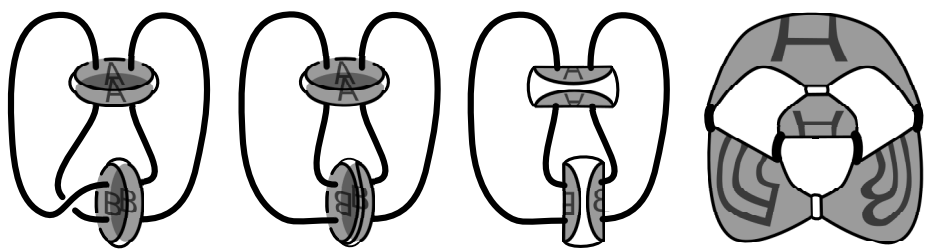}
  \caption{Left: slice 2-punctured disks up the middle (obtain parallel 2-punctured disks, shown here pulled apart). Middle left: Untwist single crossings. Middle right: Cut along plane of projection. Right: collapse remnants of the link to ideal vertices.}
  \label{Fig:BorromeanDecomp}
\end{figure}

Now if a 2-punctured disk is adjacent to a crossing in the plane of projection, the next step is to rotate that 2-punctured disk $180^\circ$ to unwind the crossing, as in the middle left of \reffig{BorromeanDecomp}. Note this rotation pulls the diagram along with it on one side, but the rotation is only performed on the 2-punctured disk adjacent to the crossing, not on the parallel 2-punctured disk. After this step, all crossings in the plane of projection have been removed.

Next, cut along the plane of projection, splitting the complement into two identical pieces as in the middle right of \reffig{BorromeanDecomp}. 

Finally, for each piece, collapse remnants of the link to ideal vertices, as on the right of \reffig{BorromeanDecomp}. We claim the result in that figure is topologically an octahedron. To see this, note it has two ideal vertices colored white, coming from crossing circles, and four ideal vertices colored black, coming from the component of the link on the plane of projection. There are four shaded faces that all have three edges, hence all shaded faces are triangles. There are four white faces, including the one running through the point at infinity in the plane of projection, and each of these white faces also has three edges, so each is a triangle. Thus the result is an ideal octahedron. Recall that there is actually another octahedron coming from our decomposition: the other octahedron comes from the region below the plane of projection after slicing along that plane. So the link complements in the Borromean family\index{Borromean family} all decompose into two ideal octahedra. 

Note that the face pairings of the two ideal octahedra that give back the original link complement will be different for the different links; they can be found by tracing backwards through the decomposition process above. To undo the step of cutting along the plane of projection, we glue matching white faces of the opposite octahedra together in pairs. To undo the step of slicing along 2-punctured disks, we glue remaining shaded triangles in pairs; however there are two options depending on whether or not we untwisted a crossing. If both parallel 2-punctured disks were adjacent to no crossings, then corresponding shaded triangles on the same ideal octahedron are glued across an ideal vertex (one of the white vertices of \reffig{BorromeanDecomp}). If there was an adjacent crossing, then a shaded triangle on one octahedron is glued to the opposite shaded triangle on the other octahedron across the (white) ideal vertex.

Finally, to see that these link complements all admit a complete hyperbolic structure, we give each of the two octahedra the geometry of a hyperbolic regular ideal octahedron,\index{ideal octahedron, regular}\index{regular ideal octahedron} then argue as in the proof of \refprop{WhiteheadGeom}. We check: each edge of the decomposition comes from the intersection of a 2-punctured disk with the plane of projection, and each edge class in the manifold is obtained by gluing four such edges. Thus the total angle around each edge will be $4(\pi/2) = 2\pi$. Again horospheres meet ideal vertices in Euclidean squares, and so the developing image cannot scale, shear, or rotate these squares. Thus edges glue without shearing singularities, and cusps are Euclidean. Hence the result is a complete hyperbolic structure.
\end{proof}

\begin{corollary}\label{Cor:VolJnm}
  The volume of a double twist knot satisfies
  \[ \vol(J(k,\ell)) < 2\,\voct, \]
  where $\voct = 3.66\dots$ is the volume of a regular ideal octahedron.\index{ideal octahedron, regular}\index{regular ideal octahedron}
\end{corollary}

%%%%%%%%%%%%%%%%%%%%%%%%%%%%%%%%%%%%%%%%%%%%%%%%%%%%%%%%%%%%%%%%%
\section{Augmenting and highly twisted knots}

The above procedure can be generalized. 

\begin{definition}\label{Def:Augmenting}
For any twist region of any link diagram, a new link is obtained by adding a single unknotted link component to the diagram, encircling the two strands of the link component. The link is said to be \emph{augmented}.\index{augmenting link diagram}\index{augmented link}\index{augmented link!definition} The added link component is called a \emph{crossing circle}\index{crossing circle}. We will refer to the original link components as \emph{knot strands}.\index{knot strand}

When a crossing circle is added to each twist region of the diagram, the link is said to be \emph{fully augmented}.\index{fully augmented link}\index{augmented link!fully augmented}
\end{definition}

The complement of an augmented link is homeomorphic to the complement of the link with any even number of crossings added to or removed from the twist region, by the same argument illustrated in \reffig{TKWProof}. Thus the complement of a fully augmented link is homeomorphic to the complement of the fully augmented link\index{augmented link} with one or zero crossings adjacent to each crossing circle. When there is one crossing adjacent to a crossing circle, we say the crossing circle is adjacent to a \emph{half-twist}.\index{half-twist}

The four links of \reffig{Borromean} are examples of fully augmented links.\index{augmented link} The decomposition of \refprop{GeomBorromean} goes through more generally for all fully augmented links.
This decomposition appears in the appendix to \cite{lackenby:alt-volume} by Agol and D.~Thurston. These links have very beautiful geometric properties, explored further in \cite{futer-purcell}, \cite{purcell:volume}, \cite{purcell:cusps}, and in the survey article \cite{Purcell:FullyAugmented}. Some of these results are included below, modeled off the exposition in \cite{Purcell:FullyAugmented}.

\begin{theorem}\label{Thm:FullyAugmentedDecomp}
Let $L$ be any fully augmented link,\index{fully augmented link!polyhedral decomposition} and assume we have applied a homeomorphism to $S^3-L$ so that $L$ has one or zero crossings adjacent to each crossing circle. Then the link complement $S^3-L$ decomposes into two identical ideal polyhedra with the following properties.
\begin{enumerate}
\item Faces of the polyhedra can be checkerboard colored. White faces correspond to regions of the plane of projection. Shaded faces are all triangles, and come from 2-punctured disks bounded by crossing circles, which we call \emph{crossing disks}.\index{crossing disk}
\item Ideal vertices are all 4-valent (before gluing).
\item Gluing the polyhedra identifies exactly four edges to a single edge class in the link complement.
\end{enumerate}
\end{theorem}

\begin{proof}
The decomposition is obtained very similarly to that in the proof of \refprop{GeomBorromean}. First, each crossing circle bounds a 2-punctured disk, which we shade. After applying a homeomorphism removing all pairs of crossings in each twist region, we may assume that the shaded disk is either adjacent to no crossings, or adjacent to a single crossing (half-twist).

Slice along the 2-punctured disks, splitting each into two parallel 2-punctured disks. Next, apply a $180^\circ$ rotation to those 2-punctured disks adjacent to a crossing, unwinding the crossing. Then slice along the plane of projection, splitting the complement into two identical pieces. Finally, shrink remnants of the link to ideal vertices. We check that each item of the theorem holds.

First, note faces are already checkerboard colored, with shaded faces coming from 2-punctured disks and white faces coming from the plane of projection. Note that edges of the decomposition come from intersections of white and shaded faces. There are exactly three edges bordering each shaded face, so each shaded face is a triangle.

Ideal vertices of the polyhedra come from remnants of the link. For those ideal vertices coming from a component of the link in the plane of projection, the ideal vertex will be adjacent to two edges coming from the 2-punctured disk on one of its ends, and two edges coming from the 2-punctured disk on its other end. Thus it is 4-valent. An ideal vertex coming from a crossing circle is also adjacent to four edges: two from each point where the link component meets the plane of projection.

Finally, note that each edge class contains four edges: two in each polyhedron lying on the parallel copies of the 2-punctured disk. 
\end{proof}

Just as with the family of Borromean rings,\index{Borromean rings}\index{Borromean family} we can show that many of these links are hyperbolic, but not all. For example, if a fully augmented link has only one crossing circle, then the polyhedral decomposition of \refthm{FullyAugmentedDecomp} will have white bigon\index{bigon} regions, and collapsing these will collapse the entire polyhedron to a triangle (exercise). Similarly, if there are parallel crossing circles then there will be white bigon faces. We wish to rule these out.

\begin{figure}
%% Accompanies image file 'F7-11-UnredA.eps' (pdf, eps, ps)
\begingroup%
  \makeatletter%
  \providecommand\color[2][]{%
    \errmessage{(Inkscape) Color is used for the text in Inkscape, but the package 'color.sty' is not loaded}%
    \renewcommand\color[2][]{}%
  }%
  \providecommand\transparent[1]{%
    \errmessage{(Inkscape) Transparency is used (non-zero) for the text in Inkscape, but the package 'transparent.sty' is not loaded}%
    \renewcommand\transparent[1]{}%
  }%
  \providecommand\rotatebox[2]{#2}%
  \newcommand*\fsize{\dimexpr\f@size pt\relax}%
  \newcommand*\lineheight[1]{\fontsize{\fsize}{#1\fsize}\selectfont}%
  \ifx\svgwidth\undefined%
    \setlength{\unitlength}{330.57112885bp}%
    \ifx\svgscale\undefined%
      \relax%
    \else%
      \setlength{\unitlength}{\unitlength * \real{\svgscale}}%
    \fi%
  \else%
    \setlength{\unitlength}{\svgwidth}%
  \fi%
  \global\let\svgwidth\undefined%
  \global\let\svgscale\undefined%
  \makeatother%
  \begin{picture}(1,0.2050273)%
    \lineheight{1}%
    \setlength\tabcolsep{0pt}%
    \put(0,0){\includegraphics[width=\unitlength]{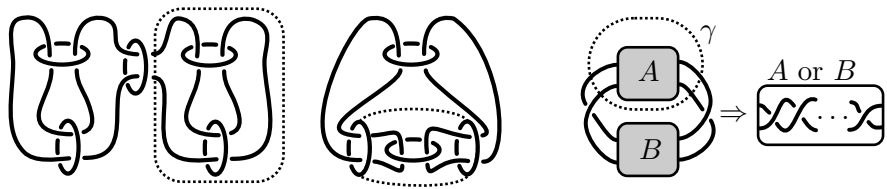}}%
    \put(0.71750726,0.12108906){\color[rgb]{0,0,0}\makebox(0,0)[lt]{\lineheight{1.25}\smash{\begin{tabular}[t]{l}$A$\end{tabular}}}}%
    \put(0.81041871,0.07582471){\color[rgb]{0,0,0}\makebox(0,0)[lt]{\lineheight{1.25}\smash{\begin{tabular}[t]{l}$\Rightarrow$\end{tabular}}}}%
    \put(0.71922615,0.03228058){\color[rgb]{0,0,0}\makebox(0,0)[lt]{\lineheight{1.25}\smash{\begin{tabular}[t]{l}$B$\end{tabular}}}}%
    \put(0.86441506,0.12223498){\color[rgb]{0,0,0}\makebox(0,0)[lt]{\lineheight{1.25}\smash{\begin{tabular}[t]{l}$A$ or $B$\end{tabular}}}}%
    \put(0.78855418,0.1709365){\color[rgb]{0,0,0}\makebox(0,0)[lt]{\lineheight{1.25}\smash{\begin{tabular}[t]{l}$\gamma$\end{tabular}}}}%
    \put(0.92082578,0.07784211){\color[rgb]{0,0,0}\makebox(0,0)[lt]{\lineheight{1.25}\smash{\begin{tabular}[t]{l}$\dots$\end{tabular}}}}%
  \end{picture}%
\endgroup%

  \caption{Left: a fully augmented link\index{fully augmented link} with a diagram that is not prime, with dotted lines indicating the closed curve contradicting the definition of prime. Middle: a fully augmented link that is not reduced, with parallel crossing circles indicated by the dotted lines. Removing one of the parallel crossing circles will give a reduced link. Right: The diagram is twist-reduced if one of the regions $A$ or $B$ consists only of bigons\index{bigon} in a twist region.}
  \label{Fig:UnreducedAug}
\end{figure}

\begin{definition}
  A fully augmented link is called \emph{reduced}\index{reduced fully augmented link}\index{fully augmented link!reduced}\index{fully augmented link} if the following hold.
  \begin{enumerate}
  \item Its diagram is connected. 
  \item Its diagram is \emph{prime},\index{prime!diagram} i.e.\ any closed curve meeting the diagram twice bounds a region on one side with no crossings. 
  \item None of its crossing circles are parallel. That is, there are no closed curves in the diagram running over exactly two crossing circles and meeting exactly two white faces on either side of the two crossing circles. See \reffig{UnreducedAug}.
  \end{enumerate}
\end{definition}

Reduced fully augmented links come from adding crossing circles to links with reduced diagrams as in the sense of the following definition.

\begin{definition}\label{Def:TwReduced}
  A diagram is \emph{twist-reduced}\index{twist-reduced} if whenever a simple closed curve $\gamma$ meets the diagram exactly twice in two crossings, running from one side of each crossing to the opposite side, then the curve $\gamma$ bounds a portion of the diagram containing a string of bigons\index{bigon} arranged end-to-end. See \reffig{UnreducedAug}, right.
\end{definition}

Note that if a diagram is not twist-reduced, then there exists a curve $\gamma$ meeting the diagram in exactly two crossings, with those crossings not separated by a string of bigons.\index{bigon} Because the two crossings are not separated by bigons, they lie in different twist regions. Augmenting the two twist regions results in two parallel crossing circles. Thus a diagram that is not twist-reduced has an associated fully augmented link that is not reduced. 

We will encounter twist-reduced diagrams again, for example in \refdef{TwistReduced}. Meanwhile, the following gives a way of building large numbers of reduced fully augmented links. 

\begin{lemma}\label{Lem:TwReducedGivesReducedAug}
Let $K$ be a link with a connected, prime, twist-reduced diagram. Then the fully augmented link obtained from $K$ by adding crossing circles to each twist region gives a reduced fully augmented link.\index{fully augmented link}\index{fully augmented link!reduced}\index{reduced fully augmented link}
\end{lemma}

\begin{proof}
Adding crossing circles to twist regions of a diagram does not change whether it is prime or connected. If the resulting fully augmented link is not reduced, there must be two parallel crossing circles. Thus in the original $K$, there are two distinct twist regions in the diagram with the property that when crossing circles are added around them, the crossing circles are parallel. Then an isotopy of one of the crossing circles to the other traces out two arcs on the plane of projection disjoint from $K$. Straightening these, and drawing arcs across $K$ over the crossing circle defines a closed curve in the diagram whose boundary meets the diagram exactly twice, once in each of the two distinct twist regions. Isotope slightly to give a curve in $K$ contradicting the definition of a twist-reduced diagram. 
\end{proof}

\begin{lemma}\label{Lem:ExistsRightAngled}
Suppose a fully augmented link is reduced and contains at least two crossing circles. Then the polyhedra in the decomposition of \refthm{FullyAugmentedDecomp} admit a hyperbolic structure in which all dihedral angles are $\pi/2$. \index{fully augmented link}\index{fully augmented link}\index{augmented link!fully}\index{fully augmented link!polyhedral decomposition}
\end{lemma}

The proof of the lemma uses circle packings. 

\begin{definition}\label{Def:CirclePacking}
A \emph{circle packing}\index{circle packing} is a connected collection of circles with disjoint interiors. The \emph{intersection graph}\index{intersection graph}\index{circle packing!intersection graph} of a circle packing is the graph with a vertex at the center of each circle, and an edge between vertices whenever the corresponding circles are tangent. 
\end{definition}

\begin{figure}
  \includegraphics{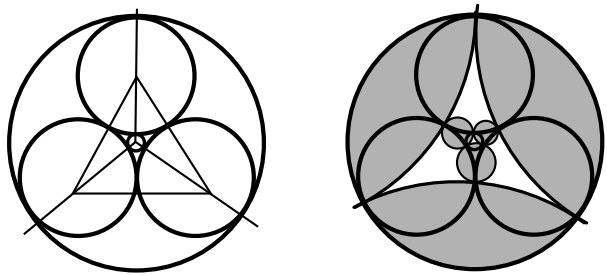}
  \caption{Left: A circle packing\index{circle packing} and its intersection graph.\index{circle packing!intersection graph} Right: Gray circles meeting white circles of a circle packing}
  \label{Fig:CirclePacking}
\end{figure}

\Reffig{CirclePacking} shows an example of a circle packing and most of its intersection graph on the left --- the vertex of the intersection graph in the unbounded region has been omitted. 

\begin{theorem}[Circle packing theorem]\label{Thm:Andreev}\index{circle packing theorem}\index{Koebe--Andreev--Thurston theorem}
Let $G$ be a finite planar graph that is simple, meaning $G$ has no loops and no multiple edges between a pair of vertices. Then $G$ is (isotopic to) the intersection graph of a circle packing\index{circle packing} on $S^2$. If $G$ is a triangulation of $S^2$, then the circle packing is unique up to M\"obius transformation.\index{M\"obius transformation}
\end{theorem}

\Refthm{Andreev} is also known as the Koebe--Andreev--Thurston theorem. It was first proved by Koebe \cite{Koebe}. We will use it here without giving its proof, as the proof is somewhat unrelated to the topic at hand. 

\begin{proof}[Proof of \reflem{ExistsRightAngled}]
Consider a polyhedron $P$ from \refthm{FullyAugmentedDecomp} for a reduced fully augmented link with at least two crossing circles. Edges and vertices of the polyhedron form a graph $\Gamma$ on $S^2$. Form a new graph $G$ on $S^2$ by taking a vertex for each white face of $\Gamma$, and an edge between vertices of $G$ if two white faces are adjacent across an ideal vertex of $\Gamma$. 

If we superimpose $G$ on $P$, then notice that each region of $G$ will contain exactly one shaded triangular face of $P$. Thus $G$ is a triangulation of $S^2$. We show that $G$ has no loops and no multiple edges.

Suppose first that $G$ has a loop. Then the edge of $G$ forming the loop can be superimposed on $P$ to run from a white face, through an ideal vertex of $P$, then back to the same white face. White faces correspond to regions of the diagram, and ideal vertices correspond to remnants of the link. Thus there is a closed curve $\gamma$ on the link diagram that runs from a region back to itself crossing over a single component of the link diagram. Because the link diagram consists of closed curves, this is possible only if the curve $\gamma$ runs along a crossing circle from one white region back to the same white region. Pushing off the crossing circle slightly, this contradicts the fact that the diagram is prime.

Now suppose that the graph $G$ has a multi-edge. Then there is a pair of white faces $W_1$ and $W_2$ of $P$ and a pair of ideal vertices $v_1$ and $v_2$ such that $v_1$ and $v_2$ are both adjacent to $W_1$ and $W_2$. Form a loop in $P$ running from $W_1$ through $v_1$ to $W_2$, then back through $v_2$ to $W_1$. This loop corresponds to a loop $\gamma$ in the diagram meeting the regions on the plane of projection corresponding to $W_1$ and $W_2$ and meeting two distinct link components between those regions. If the link components came from components on the plane of projection, then this contradicts the fact that the diagram is prime. If the link components came from crossing circles, then it contradicts the fact that the diagram is reduced. If one link component lies in the plane of projection and the other is a crossing circle, then we may slide slightly off the crossing circle to obtain a loop meeting exactly three components in the plane of projection. This is impossible for a closed curve and closed link components. 

It follows that $G$ is a finite, simple, planar graph that is a triangulation of $S^2$. The circle packing\index{circle packing} theorem, \refthm{Andreev},\index{circle packing theorem}\index{Koebe--Andreev--Thurston theorem} implies that there is a unique circle packing of $S^2$ with $G$ as its intersection graph. View $S^2$ as the boundary at infinity\index{boundary at infinity} of $\HH^3$. The circle packing of $G$ is then a circle packing on $\bdy \HH^3$. Each Euclidean circle on $\bdy\HH^3$ is the boundary of a plane in $\HH^3$. Color these planes \emph{white}. 

Because the intersection graph of the circle packing\index{circle packing!intersection graph} is a triangulation, regions complementary to the circle packing meet exactly three circles from the packing. There is a unique Euclidean circle running through the three points of tangency of the circle packing. Again this defines a geodesic plane in $\HH^3$. This plane will intersect the white planes at right angles. Color this plane \emph{gray}. See \reffig{CirclePacking}.

For each white plane, remove from $\HH^3$ the region bounded by that plane that is disjoint from the other white planes. Similarly for each gray plane. The result is a right-angled hyperbolic ideal polyhedron that is isomorphic to $P$, proving the lemma. 
\end{proof}

\begin{theorem}\label{Thm:FullyAugCompleteHyp}
The complement of a reduced fully augmented link\index{fully augmented link!reduced}\index{fully augmented link!hyperbolic} with at least two crossing circles admits a complete hyperbolic structure, which is obtained by putting a right-angled structure on each of the polyhedra of \refthm{FullyAugmentedDecomp}. 
\end{theorem}

\begin{proof}
By \reflem{ExistsRightAngled}, there exists a right-angled ideal hyperbolic polyhedron with the combinatorics of one of the polyhedron of \refthm{FullyAugmentedDecomp}. We give each of the polyhedra of \refthm{FullyAugmentedDecomp} the hyperbolic structure of this right-angled hyperbolic polyhedron, and glue by corresponding face-pairing isometries\index{face-pairing isometry} to obtain the fully augmented link.

To show this admits a complete hyperbolic structure, we need to show the angle around each edge is $2\pi$, that there is no shearing around edges, and that the cusps are all Euclidean. Because each edge class contains four edges, and each edge has dihedral angle $\pi/2$, the angle sum around each edge is $2\pi$.

Now consider cusps. Any horosphere meets an ideal vertex of the right-angled polyhedron in a rectangle. The developing image of a cusp is obtained by gluing these rectangles according to the gluing isometries on the faces. Note that a white face is glued by a reflection to the identical white face on the opposite polyhedron, so gluing across white sides of a rectangle does not scale or rotate. But then the gluing across shaded faces cannot scale or rotate either. Hence the developing image of each cusp is a tiling of the plane by Euclidean rectangles. Thus around each vertex there cannot be shearing, and the structure on the cusp must be Euclidean. So this gives the complete hyperbolic structure on the fully augmented link. 
\end{proof}

\begin{corollary}\label{Cor:AugmentedGeodSfces}
In a reduced fully augmented link\index{fully augmented link} with at least two crossing circles, each shaded 2-punctured disk bounded by a crossing circle is a totally geodesic surface embedded in the link complement. The \emph{white surface}, obtained by gluing together regions corresponding to regions on the plane of projection (white faces), is also a totally geodesic surface embedded in the hyperbolic link complement. Moreover, these shaded 2-punctured disks and white surfaces meet at right angles whenever they intersect.
\end{corollary}

\begin{proof}
In the polyhedral decomposition, these surfaces become white and shaded faces, which are straightened to portions of geodesic planes to obtain the hyperbolic structure. Thus we know that these surfaces are \emph{pleated},\index{pleated surface} i.e.\ they decompose into ideal polygons, each of which is totally geodesic. In general pleated surfaces are bent along the edges bounding each polygon, so they are not necessarily totally geodesic. However, in this case, white faces meet shaded faces at angle $\pi/2$, thus in the gluing, white faces glue to white faces with angle $\pi$, i.e.\ no bending, and similarly for shaded faces. If follows that these surfaces are totally geodesic. 
\end{proof}

%%%%%%%%%%%%%%%%%%%%%%%%%%%%%%%%%%%%%%%%%%%%%%%%%%%%%%%%%%%%%%%%%
\section{Cusps of fully augmented links}

For many applications in later chapters, it will be useful to know more explicit information about the geometry of fully augmented links,\index{fully augmented link} particularly the geometry of their \emph{cusps}.\index{cusp} Recall from \refthm{EuclidCusp} that each cusp admits a Euclidean structure.\index{Euclidean structure} In this section, we will determine properties of that Euclidean structure for fully augmented links. The exposition is similar to that in \cite{futer-purcell}.

Consider the universal cover of the complement of a hyperbolic fully augmented link. By \refcor{AugmentedGeodSfces}, the universal cover will contain the lift of embedded totally geodesic white surfaces, which will be a collection of disjoint totally geodesic planes that we color white in $\HH^3$. It will also contain the lifts of embedded totally geodesic shaded 2-punctured disks bounded by crossing circles. These will also be totally geodesic planes in $\HH^3$ and we call them \emph{shaded}. The white planes and shaded planes meet at right angles in $\HH^3$. They cut out all the translates of the two ideal polyhedra of \refthm{FullyAugmentedDecomp} under the developing map.

Apply an isometry so that the boundary $\widetilde{T}$ of a neighborhood of the point at infinity in $\HH^3$ projects under the covering map to a cusp torus $T$ of the fully augmented link. Because each link component meets both white and shaded surfaces, in the universal cover we will see vertical planes corresponding to white and shaded surfaces running into the point at infinity, meeting $\widetilde{T}$ in a rectangular lattice.

If we forget the fact that the edges of the lattice have lengths, but consider each rectangle on $\widetilde{T}$ as a topological object with two opposite shaded sides and two opposite white sides, then we obtain the following.

\begin{lemma}\label{Lem:FullyAugCusps}
Let $T$ be a cusp torus of a fully augmented link,\index{fully augmented link}\index{fully augmented link!cusp} with universal cover $\widetilde{T}$ tiled by rectangles coming from white and shaded surfaces. Let $s$ denote a step along a shaded surface between two white surfaces, and let $w$ denote a step along a white surface between two shaded ones. Then a fundamental domain for $T$ is given as follows.
\begin{itemize}
\item If $T$ comes from a crossing circle without a half-twist, then it has meridian $w$ and longitude $2s$.
\item If $T$ comes from a crossing circle with a half-twist, it has meridian $w\pm s$ (depending on the direction of the twist) and longitude $2s$.
\item If $T$ comes from a knot strand, i.e.\ a component that is not a crossing circle, then it has meridian $2s$ and longitude $nw+ks$, where $n$ is the number of twist regions met by the strand, with multiplicity, and $k$ is some integer. 
\end{itemize}
\end{lemma}

\begin{proof}
From the construction of the polyhedral decomposition of $S^3-L$, each crossing circle gives rise to an ideal vertex of each polyhedron. Thus a fundamental domain for a crossing circle consists of two rectangles, given by neighborhoods of the corresponding 4-valent ideal vertices.

In the case that there are no half-twists, the shaded faces adjacent to the ideal vertex are glued to each other. Thus an arc running along a white face has its endpoints glued into a meridian, and thus the meridian in this case is $w$. As for the longitude, a white face on one polyhedron is glued to a white face on the other. Thus a longitude steps along two shaded sides, one on one polyhedron and one on the other, before closing up. See \reffig{CrossingCircleCusp}.

\begin{figure}
%% Accompanies image file 'F7-13-CCCusp.eps' (pdf, eps, ps)
\begingroup%
  \makeatletter%
  \providecommand\color[2][]{%
    \errmessage{(Inkscape) Color is used for the text in Inkscape, but the package 'color.sty' is not loaded}%
    \renewcommand\color[2][]{}%
  }%
  \providecommand\transparent[1]{%
    \errmessage{(Inkscape) Transparency is used (non-zero) for the text in Inkscape, but the package 'transparent.sty' is not loaded}%
    \renewcommand\transparent[1]{}%
  }%
  \providecommand\rotatebox[2]{#2}%
  \newcommand*\fsize{\dimexpr\f@size pt\relax}%
  \newcommand*\lineheight[1]{\fontsize{\fsize}{#1\fsize}\selectfont}%
  \ifx\svgwidth\undefined%
    \setlength{\unitlength}{191.34446526bp}%
    \ifx\svgscale\undefined%
      \relax%
    \else%
      \setlength{\unitlength}{\unitlength * \real{\svgscale}}%
    \fi%
  \else%
    \setlength{\unitlength}{\svgwidth}%
  \fi%
  \global\let\svgwidth\undefined%
  \global\let\svgscale\undefined%
  \makeatother%
  \begin{picture}(1,0.47610877)%
    \lineheight{1}%
    \setlength\tabcolsep{0pt}%
    \put(0,0){\includegraphics[width=\unitlength]{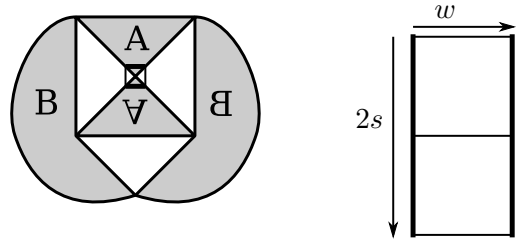}}%
    \put(0.83549029,0.43888375){\color[rgb]{0,0,0}\makebox(0,0)[lt]{\lineheight{1.25}\smash{\begin{tabular}[t]{l}$w$\end{tabular}}}}%
    \put(0.67969488,0.2242939){\color[rgb]{0,0,0}\makebox(0,0)[lt]{\lineheight{1.25}\smash{\begin{tabular}[t]{l}$2s$\end{tabular}}}}%
  \end{picture}%
\endgroup%

  \caption{A fundamental region for a crossing circle.}
  \label{Fig:CrossingCircleCusp}
\end{figure}

For a knot strand $K$ meeting no half-twists, there will be one ideal vertex of one polyhedron, hence one rectangular vertex neighborhood, for each portion of $K$ between adjacent crossing circles. These rectangles are glued end to end along shaded faces coming from the crossing disks to complete a longitude. Thus there will be $n$ such rectangles, and a longitude is given by $n$ steps along white faces, or $nw$. There will be $n$ identical rectangles glued end to end in the other polyhedron. These two blocks of $n$ rectangles will be glued along their white faces to form a $2\times n$ block, making up the fundamental domain of $K$. A meridian is given by two steps along shaded faces. See \reffig{KnotStrandCusp}.

\begin{figure}
  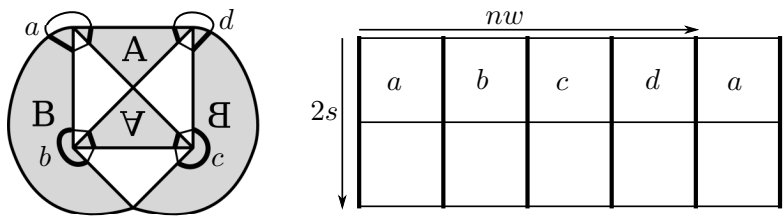
  \caption{A fundamental region for a knot strand with no half twists.}
  \label{Fig:KnotStrandCusp}
\end{figure}

If there are half-twists, then the gluing changes along shaded faces at half-twists. A shaded triangle on one polyhedron will be glued to the opposite shaded triangle on the other polyhedron. This introduces shearing into the fundamental domain, as in \reffig{AugShearing}.
\begin{figure}
  \includegraphics{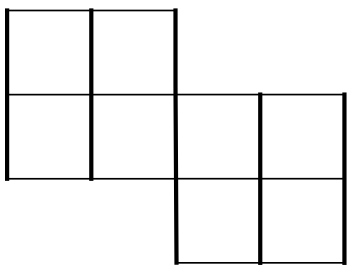}
  \caption{Adding a half twist shifts the gluing along the shaded faces, shearing the fundamental domain.}
  \label{Fig:AugShearing}
\end{figure}
Since the shearing only occurs as shaded faces are glued, it does not affect the longitude of a crossing circle or the meridian of a knot strand: these are both $2s$. However, it will adjust a meridian of a crossing circle by adding $\pm s$, and it will adjust the longitude of a knot strand by adding $\pm s$ for each half-twist. Thus the longitude of a knot strand becomes $nw + ks$ for some integer $k$.
\end{proof}

\Reflem{FullyAugCusps} is purely topological. We now wish to give geometric information on the rectangles forming the cusps. To do so, we need to find more explicit embedded cusp neighborhoods of the cusps of a fully augmented link. An embedded cusp neighborhood lifts to a disjoint collection of horoballs in the universal cover $\HH^3$, one for each ideal vertex of each translate of the ideal polyhedra under the developing map. We will find an embedded cusp neighborhood by finding a collection of embedded horoballs about ideal vertices of the polyhedra forming a fully augmented link.

\begin{figure}
%% Accompanies image file 'F7-16-Midpt.eps' (pdf, eps, ps)
\begingroup%
  \makeatletter%
  \providecommand\color[2][]{%
    \errmessage{(Inkscape) Color is used for the text in Inkscape, but the package 'color.sty' is not loaded}%
    \renewcommand\color[2][]{}%
  }%
  \providecommand\transparent[1]{%
    \errmessage{(Inkscape) Transparency is used (non-zero) for the text in Inkscape, but the package 'transparent.sty' is not loaded}%
    \renewcommand\transparent[1]{}%
  }%
  \providecommand\rotatebox[2]{#2}%
  \newcommand*\fsize{\dimexpr\f@size pt\relax}%
  \newcommand*\lineheight[1]{\fontsize{\fsize}{#1\fsize}\selectfont}%
  \ifx\svgwidth\undefined%
    \setlength{\unitlength}{116.25000286bp}%
    \ifx\svgscale\undefined%
      \relax%
    \else%
      \setlength{\unitlength}{\unitlength * \real{\svgscale}}%
    \fi%
  \else%
    \setlength{\unitlength}{\svgwidth}%
  \fi%
  \global\let\svgwidth\undefined%
  \global\let\svgscale\undefined%
  \makeatother%
  \begin{picture}(1,0.71290281)%
    \lineheight{1}%
    \setlength\tabcolsep{0pt}%
    \put(0,0){\includegraphics[width=\unitlength]{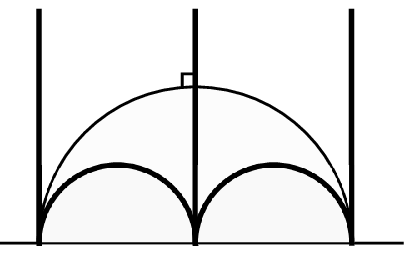}}%
    \put(0.51612904,0.55161299){\color[rgb]{0,0,0}\makebox(0,0)[lt]{\lineheight{1.25}\smash{\begin{tabular}[t]{l}$m$\end{tabular}}}}%
    \put(0.46431958,0.02597036){\color[rgb]{0,0,0}\makebox(0,0)[lt]{\lineheight{1.25}\smash{\begin{tabular}[t]{l}$0$\end{tabular}}}}%
    \put(0.84164069,0.02909799){\color[rgb]{0,0,0}\makebox(0,0)[lt]{\lineheight{1.25}\smash{\begin{tabular}[t]{l}$1$\end{tabular}}}}%
  \end{picture}%
\endgroup%

  \caption{When an ideal triangle\index{ideal triangle} in $\HH^2$ has vertices at $0$, $1$, and $\infty$, one of its midpoints\index{midpoint} will lie in $\HH^2$ at height $1$.}
  \label{Fig:Midpoint}
\end{figure}

\begin{definition}\label{Def:MidpointIdealTriangle}
Let $T\subset \HH^3$ be an ideal triangle.\index{ideal triangle} For each edge $e$ of $T$, define the \emph{midpoint}\index{midpoint} $m$ of $e$ to be the point such that the geodesic from $m$ to the opposite ideal vertex is perpendicular to $e$. Note this point is unique. See \reffig{Midpoint}.

For each edge $e$ of the ideal polyhedral decomposition of a fully augmented link, define its \emph{midpoint} to be the midpoint of that edge on one of the two ideal triangles\index{ideal triangle} adjacent to the edge. Note that since the two polyhedra are symmetric by a reflection in the white faces, both triangles adjacent to $e$ have the same midpoint, so the midpoint of each edge is well-defined. 
\end{definition}

\begin{lemma}\label{Lem:EmbeddedHoroballs}
Let $L$ be a hyperbolic fully augmented link,\index{fully augmented link}\index{fully augmented link!cusp} with decomposition into ideal polyhedra $P_1$ and $P_2$. For each ideal vertex of $P_i$, there is a unique horoball meeting the midpoint\index{midpoint} of each edge through that ideal vertex. The collection of all such horoballs, intersected with $P_i$ and $P_j$, glue to give an embedded cusp neighborhood of all the cusps of $S^3-L$. 
\end{lemma}

\begin{proof}
Place $P_i$ in $\HH^3$ so that the ideal vertex of interest lies at infinity, and so that one of the two shaded faces meeting the ideal vertex has its ideal vertices at $0$, $1$, and $\infty$ in $\HH^3$. Note that the edges of that shaded face have midpoints\index{midpoint} at height $1$, as in \reffig{Midpoint}. Because the polyhedron is right-angled, the other shaded face will have ideal points at some points $ci$, $1+ci$, and $\infty$ for some $c\in \RR$. Thus again the midpoints of these edges lie at height $1$. Then the horoball of height $1$ about infinity meets the midpoint of each edge through the ideal vertex.

The above discussion applies to any vertex of any polyhedron, and so this proves the first statement of the lemma. However, to show that these horoballs glue to give an embedded cusp neighborhood, we need to show that under the developing map, these horoballs have disjoint embedded interiors in $\HH^3$.

Develop in a neighborhood of infinity. Since white faces are glued by reflection, the developing map takes $P_i$ to a reflected copy of $P_i$, where the reflection is through the vertical plane determined by a white face meeting infinity. Note that the reflection isometry takes points at height $1$ to points at height $1$. Similarly, because the polyhedra are right angled, developing by gluing shaded faces produces shaded faces of the same width, and thus midpoints\index{midpoint} are height $1$. Thus a horoball of height $1$ through infinity will meet the midpoints of all edges through infinity under the developing image.

We claim that this horoball cannot meet any white faces besides those that have an ideal vertex at infinity. Consider a white face that does not meet infinity. It lies in a hemisphere with boundary a circle $C$ on $\CC$. Because the white surface is embedded in the fully augmented link complement, the lifts of this surface are disjointly embedded in $\HH^3$. Thus the boundary circles of all lifts of white faces meet only at points of tangency corresponding to ideal vertices. Thus the circle $C$ meets the boundaries of the vertical planes containing white faces meeting $\infty$ only in points of tangency. The vertical planes have boundary on $\CC$ a collection of parallel vertical lines, and these lines must be exactly distance $1$ apart. Then the diameter of $C$ can be at most $1$. It follows that the height of the hemisphere containing a white face that does not run through infinity must be at most $1/2$; therefore the horoball at height $1/2$ cannot meet it.

By an isometry, the previous argument applies to any ideal vertex. Thus we have proved that the horoballs through the midpoints\index{midpoint} of ideal vertices of $P_i$ only meet white faces that run through the center of the horoball at infinity.

Suppose, by way of contradiction, that under the developing map one of these horoballs $H$ centered at a point $p$ in $\CC$ has diameter strictly greater than $1$, so that the collection of interiors of horoballs will not be embedded. Then $p$ must lie on the boundary of one of the vertical white planes, else $H$ intersects a vertical white plane in a compact region, giving a contradiction. 

So $p$ is an ideal vertex of a polyhedron meeting a white face on a vertical plane $V$, and some other white plane $W$. The boundary $\bdy W$ is a circle on $\CC$ of diameter at most $1$, as we have seen above. The vertex $p$ also meets two shaded faces, and at least one of these, call it $S$, is not a vertical plane. Then $S\cap V$ is an ideal edge of a polyhedron, and it must have a midpoint.\index{midpoint} The midpoint is obtained by taking a perpendicular from a point on $\bdy W$ to the semicircle $S\cap V$ on the vertical plane $V$. The set of all points obtained by dropping a perpendicular from $\bdy W$ to $V$ is a circle of diameter equal to the diameter of $\bdy W$ on the plane $V$; see \reffig{Diameter}. But $H$ has diameter greater than $1$, so this entire circle lies inside of $H$. This contradicts the fact that $H$ does not contain any of the midpoints of edges through $p$.

\begin{figure}
  \includegraphics{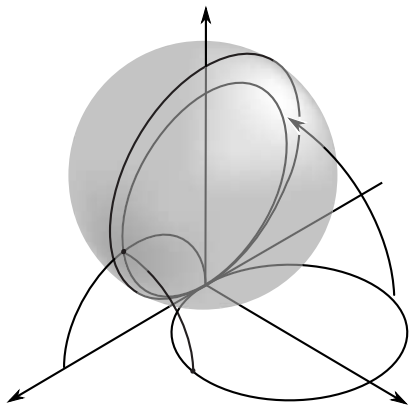}
  \caption{If $H$ is centered at a point $p$ on a vertical white plane, and has diameter greater than $1$, then it must contain the midpoint\index{midpoint} of an edge through $p$.}
  \label{Fig:Diameter}
\end{figure}

Thus when we expand all horoballs to the midpoints\index{midpoint} of their adjacent edges, all those centered at points on $\CC$ have diameter at most $1$, while that at infinity has height exactly $1$, so their interiors are embedded. Since the above discussion applies to any ideal vertex of any polyhedron, we conclude that under the developing map, interiors of all such horoballs are embedded, and thus the quotient under the covering map gives an embedded horoball neighborhood of each cusp of $S^3-L$. 
\end{proof}

\begin{corollary}\label{Cor:AugSideLength}
Let $L$ be a hyperbolic fully augmented link.\index{fully augmented link}\index{fully augmented link!cusp} There exists an embedded horoball neighborhood of the cusps of $S^3-L$ such that, when measured in the induced Euclidean metric on the boundary of each cusp, the sides of the steps $s$ and $w$ (of \reflem{FullyAugCusps}) have lengths $\ell(s)=1$ and $\ell(w)\geq 1$. 
\end{corollary}

\begin{proof}
If we place the ideal vertices of a shaded triangle at $0$, $1$, and $\infty$, then the midpoints\index{midpoint} of the edges from $0$ to $\infty$ and from $1$ to $\infty$ are of height $1$. Thus the horoball neighborhood through these points is at height $1$, and distance along the boundary of this horoball is just Euclidean distance. Since the shaded triangle meets this plane in a line segment from $0$ to $1$, the length of the step $s$ is $\ell(s)=1$.

To find $w$, we note that there will be horoballs of diameter $1$ centered at all the corners of the rectangle containing the step $w$. Two of these will be centered at $0$ and $1$, the other two at some $ci$ and $1+ci$ in $\CC$, for some $c=\ell(w)$. Because the four horoballs are disjoint, we must have $\ell(w)\geq 1$. 
\end{proof}

The above results lead to consequences on slope lengths of Dehn fillings.\index{Dehn filling}

\begin{theorem}\label{Thm:AugSlopeLengths}
Let $L$ be a hyperbolic fully augmented link.\index{fully augmented link}\index{fully augmented link!cusp} Let $C_1, \dots, C_k$ be crossing circles of $L$. Let $s_j$ be a slope on $N(C_j)$ such that Dehn filling\index{Dehn filling} along $s_j$ replaces the crossing circle $C_j$ by a twist region with $n_j$ crossings (with $n_j$ even if and only if $C_j$ is not adjacent to a half-twist). Then there is an embedded horoball neighborhood of all cusps of $S^3-L$ such that on the boundary of each cusp, the length of $s_j$ is at least $\ell(s_j) \geq \sqrt{n_j^2+1}$. 
\end{theorem}

\begin{proof}
The slope of the Dehn filling that replaces a crossing circle with $2a_j$ crossings runs over one meridian and $a_j$ longitudes. 

If $n_j=2a_j$ is even, then $C_j$ meets no half-twist. Then \reflem{FullyAugCusps} implies that the slope $s_j$ will have the form $w+2a_j s$ or $w-2a_j s$. Because $w$ and $s$ run in orthogonal directions, and each has length at least one by \refcor{AugSideLength}, the length of $w\pm 2a_js$ is at least $\sqrt{1+(2a_j)^2} =\sqrt{1+n_j^2}$, as claimed.

If $n_j=2a_j+1$ is odd, then $C_j$ meets a half-twist, and \reflem{FullyAugCusps} implies that $s_j$ has the form $w\pm s +2a_j s$ or $w\pm s-2a_j s$, with the signs the same: $s_j = w\pm (2a_j+1)s$. Again because $w$ and $s$ are orthogonal, \refcor{AugSideLength} implies the length of $s_j$ is at least $\sqrt{1 + (2a_j+1)^2} = \sqrt{1+n_j^2}$. 
\end{proof}

%%%%%%%%%%%%%%%%%%%%%%%%%%%%%%%%%%%%%%%%%%%%%%%%%%%%%%%%%%%%%%%%%
\section{Exercises}

\begin{exercise}\label{Ex:Whitehead}
  In this exercise, you investigate the two diagrams of the Whitehead link\index{Whitehead link} shown in \reffig{Whitehead}. 
  \begin{enumerate}
  \item Show by a sequence of diagrams that the two links in that figure are isotopic.
  \item Use SnapPy \cite{SnapPy} to show that the two link complements are isometric. The check using SnapPy is not mathematically rigorous, but in this case the link has a special property: it is \emph{arithmetic}\index{arithmetic link}. We will not define an arithmetic link here (we won't use the definition elsewhere), but a consequence of arithmeticity is that the program Snap \cite{Snap} can be used to give a mathematically rigorous certification that the two links shown are isometric.
  \end{enumerate}
\end{exercise}

\begin{exercise}\label{Ex:WhiteheadOctahedron}
  Use the methods of \refchap{Fig8Decomp} to prove that the complement of a Whitehead link\index{Whitehead link} can be decomposed into two ideal pyramids with a square base, which in turn can be glued to an ideal octahedron. 
\end{exercise}

\begin{exercise}
Shown on the left of \reffig{2BorromeanRings} is the diagram of a link which we claim is the Borromean rings.\index{Borromean rings} Shown on the right is a more familiar diagram of the Borromean rings. There are several different ways to prove the complements of these hyperbolic manifolds are isometric.
  \begin{enumerate}
  \item Show by a sequence of diagram moves that the links are isotopic. Why does this suffice to show the complements are isometric?
  \item Find a hyperbolic structure on each by hand, and show by hand that the manifolds are isometric. This will take some work, and sounds tedious. The exercise here is to think about why this will be tedious: list the steps involved. 
  \item Use computational tools. Use SnapPy \cite{SnapPy} to show they are isometric. As in \refex{Whitehead}, these links are arithmetic,\index{arithmetic link} so you can check using Snap \cite{Snap} that the link complements are isometric, which gives a mathematically rigorous certification. 
  \end{enumerate}

  \begin{figure}
    \includegraphics{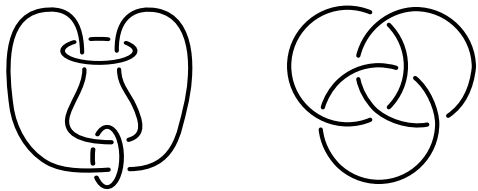}
    \caption{Two different diagrams of the Borromean rings.\index{Borromean rings}}
    \label{Fig:2BorromeanRings}
  \end{figure}
\end{exercise}

\begin{exercise}
  The $(p,q,r)$-pretzel link. As $p$, $q$, $r$ go to infinity, find geometric limits of pretzel links. Find a universal upper bound on their volumes.
\end{exercise}

\begin{exercise} (Topology of the solid torus)
  A solid torus $V$ is homeomorphic to $S^1\times D^2$, where a specified homeomorphism $h\from S^1\times D^2 \to V$ is called a \emph{framing}\index{framing of solid torus}.
  \begin{enumerate}
  \item[(a)] A non-trivial simple closed curve in $\bdy V$ is called a \emph{meridian}\index{meridian} if it bounds a disk in $V$. Prove that if $\mu$ is a meridian, then for some framing $h\from S^1\times D^2\to V$, $\mu = h(\{1\}\times \bdy D^2)$.
  \item[(b)] A non-trivial simple closed curve $\lambda$ in $\bdy V$ is called a \emph{longitude}\index{longitude} if it represents a generator of $\pi_1(V)\cong \ZZ$. Prove that if $\lambda$ is a longitude, then for some framing $h\from S^1\times D^2\to V$, $\lambda = h(S^1\times\{1\})$.
  \item[(c)] Prove that there are infinitely many ambient isotopy\index{ambient isotopy} classes of longitudes in a solid torus. 
  \end{enumerate}
\end{exercise}

\begin{exercise}
For the Whitehead link,\index{Whitehead link} find slopes of Dehn filling\index{Dehn filling} giving the twist knot $J(2,n)$ for $n$ even. Write them as $p\mu + q\lambda$, for relatively prime integers $p$ and $q$, where $\mu$ is a meridian and $\lambda$ is the longitude that bounds a disk in $S^3$. This is called the \emph{standard longitude}.\index{standard longitude!Whitehead link}

  Repeat for $n$ odd, using the isometric link. 
\end{exercise}

\begin{exercise}
  Using the meridian and standard longitude as a basis for two boundary components of the exterior of the Borromean rings\index{Borromean rings} (i.e.\ take a longitude on each component that bounds a disk in $S^3$), find the slopes of the Dehn fillings\index{Dehn filling} of the Borromean rings that give $J(2k,2\ell)$.

  Repeat for $J(2k,2\ell+1)$ and $J(2k+1,2\ell+1)$. 
\end{exercise}

\begin{exercise}
The simplest fully augmented link\index{fully augmented link} has a single crossing circle; it comes from augmenting a knot with only one twist region. Show that when we apply the decomposition of this chapter to the fully augmented link with only one crossing circle, the result is not a decomposition into two ideal polyhedra. What does the decomposition give?
\end{exercise}

\begin{exercise}
  Prove a result analogous to \refthm{AugSlopeLengths} for knot strand cusps. If $K_i$ is a knot strand cusp of a hyperbolic fully augmented link,\index{fully augmented link} and $s_i$ is a slope on $K_i$ that represents a nontrivial filling (i.e.\ $s_i$ is not a meridian), then the length of $s_i$ is at least $m_i$, where $m_i$ denotes the number of crossing disks that $K_i$ intersects, counted with multiplicity.
\end{exercise}

\begin{exercise}\label{Ex:FullyAug2Bridge}
In \refchap{TwoBridge} we will consider a class of links called \emph{two-bridge links}\index{two-bridge links} which have twist regions arranged in two rows, illustrated in \reffig{2BridgeDiagram}. Show that the complement of the fully augmented link\index{fully augmented link} coming from a 2-bridge\index{2-bridge knot or link} link can be obtained by gluing a collection of regular ideal octahedra. How many regular ideal octahedra?
\end{exercise}

%% Ch08_Essential.tex

\chapter{Essential Surfaces}\label{Chap:Essential}
\blfootnote{Jessica S. Purcell, Hyperbolic Knot Theory}

We have already encountered hyperbolic surfaces embedded in hyperbolic 3-manifolds, for example the 3-punctured spheres that bound ``shaded surfaces'' in fully augmented links.\index{fully augmented link} In this chapter, we explore surfaces more carefully. We will see that many results can be deduced about the geometry of 3-manifolds from the topology of the surfaces they contain. 

%%%%%%%%%%%%%%%%%%%%%%%%%%%%%%%%%%%%%%%%%%%%%%%%%%%%%%%%%%%%%%%%%
\section{Incompressible surfaces}

In section gives many of the definitions needed to describe surfaces in 3-manifolds. These are topological in nature, and are standard in 3-manifold topology. Good references are \cite{hempel:3mflds}, \cite{Hatcher:3MfldNotes}, and \cite{Schultens:3Mfld}.

\begin{remark}\label{Rem:SmoothCategory}
  Whenever we step from geometric arguments to topological ones, we typically need to take some care to discuss whether our objects will be merely continuous, or piecewise linear, or smooth, because the topology of manifolds, maps, etc.\ can behave very differently under different assumptions. We will assume throughout that our manifolds and maps are smooth, i.e.\ in the $C^\infty$ category,\index{smooth category} unless otherwise stated. This allows us to assume basic results on differentiable manifolds:
  \begin{itemize}
  \item Submanifolds have \emph{tubular neighborhoods}.\index{tubular neighborhood} That is, for any submanifold $S$ embedded in a manifold $M$ of codimension $k$, there exists an open neighborhood of $S$ embedded in $M$ diffeomorphic to $S\times D^k$. Thus a link $L$ embedded in $S^3$ lies in a tube $L\times D^2$ embedded in $S^3$ with $L$ at its core. A surface $S$ lies in a thickened surface $S\times D^1 = S\times I$. A tubular neighborhood is also sometimes called a \emph{regular neighborhood}\index{regular neighborhood}. 
  \item An isotopy of a submanifold can be extended to an isotopy of the ambient manifold, i.e.\ to an ambient isotopy.\index{ambient isotopy}
  \item Submanifolds can be perturbed to intersect transversely.\index{transverse intersection}
  \end{itemize}
  More information on these results can be found in a standard text on differential topology. 
\end{remark}

Our submanifolds will typically be surfaces, and throughout, these will almost always be \emph{properly} embedded in the ambient 3-manifold,\index{proper embedding, definition} where a \emph{proper} embedding is one in which the boundary of the surface is mapped by the embedding to the boundary of the 3-manifold: $S\cap \bdy M=\bdy S$ where the intersection is transverse. Additionally, we will assume throughout that the ambient 3-manifold is orientable; this is the case for knot complements in $S^3$, for example. 

\begin{definition}\label{Def:Compressible}
Let $F$ be a connected surface properly embedded in a 3-manifold. An embedded disk $D\subset M$ with $\partial D \subset F$ is said to be a \emph{compression disk}\index{compression disk} for $F$ if $\partial D$ does not bound a disk on $F$. A surface that admits a compression disk is \emph{compressible}.\index{compressible} If the surface contains no compression disk, and is not the sphere $S^2$, projective plane $P^2$, or disk $D^2$, then we say it is \emph{incompressible}.\index{incompressible}
\end{definition}

\begin{figure}
\includegraphics{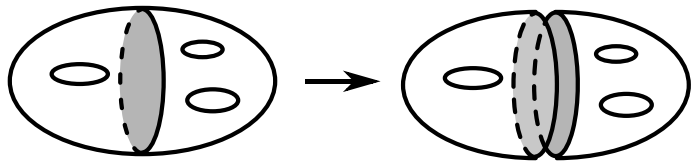}
  \caption{A compressible\index{compressible} surface can surgered along a disk, replaced with a simpler surface or surfaces}
  \label{Fig:Surgery}
\end{figure}

Compressible surfaces can be simplified, as follows. Suppose $S$ is a surface properly embedded in a 3-manifold $M$, and $D$ is a disk embedded in $M$ with $\bdy D$ contained in $S$. Let $\nu(D)$ be a thickened $D$: i.e.\ $\nu(D)$ is homeomorphic to $D\times I$, with $(\bdy D)\times I$ a regular neighborhood of $\bdy D$ in $S$, and $\nu(D)$ embedded in a tubular neighborhood\index{tubular neighborhood} of $D$ in $M$. 
We may form a new (possibly disconnected) properly embedded surface from $S$ and $D$ by the following procedure: remove $\bdy \nu(D) \cap S$ from $S$ and attach the two parallel disks $\bdy \nu(D) - (\bdy \nu(D)\cap S)$ to $S$. See \reffig{Surgery}. (Technically, as a final last step we need to smooth corners, to remain in the $C^\infty$ category. Such a smoothing is easily done, and we will assume it is done without comment for related constructions.)

\begin{definition}\label{Def:Surgery}
The process of replacing $S$ by attaching the two disks $\bdy\nu(D)-(\bdy\nu(D)\cap S)$ to curves $S-(\bdy\nu(D)\cap S)$ is called \emph{surgery} of $S$ along $D$. Usually, we use the verb to describe the procedure: \emph{surger}\index{surger} $S$ along $D$. 
\end{definition}

An incompressible surface\index{incompressible} admits no compression disk. That means that if $D$ is a disk embedded in the ambient 3-manifold with $\bdy D \subset S$, then $\bdy D$ also bounds a disk $E\subset S$. Thus in this case, if we surger $S$ along $D$ we obtain two surfaces: one diffeomorphic to $S$, and one diffeomorphic to a 2-sphere ($D\cup E$). Often in our applications the 2-sphere bounds a ball and contracts to a point. Thus surgery on an incompressible surface does nothing to simplify the surface. 

\begin{example}\label{Example:Satellite}
  As an example of an incompressible\index{incompressible} surface, consider the torus $T$ marked by dashed lines embedded in the knot complement shown in \reffig{Satellite}. On its outside, this torus bounds a manifold homeomorphic to the figure-8 knot complement.

  We claim there cannot be a compression disk\index{compression disk} on the outside, in the complement of the figure-8 knot. For suppose $D$ is such a disk. Surger the torus $T$ along $D$. The result is a sphere $S$ embedded in $S^3$. Any sphere in $S^3$ bounds two balls, one on either side. One ball bounded by $S$ must contain the figure-8 knot and the compression disk\index{compression disk} $D$. The other is a ball contained in the figure-8 knot complement, disjoint from $D\times I$. If we undo the surgery along $D$, we glue this ball along two disks, $D\times\{0\}$ and $D\times \{1\}$, yielding a solid torus. This implies that the figure-8 knot complement is homeomorphic to a solid torus, which is the unknot complement. But this is a contradiction, for example because the figure-8 knot is hyperbolic.

\begin{figure}
\includegraphics{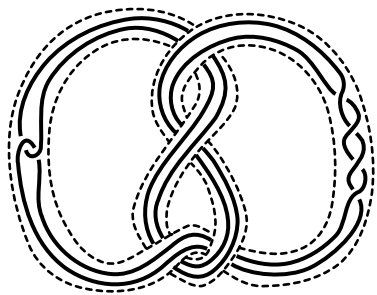}
  \caption{The torus shown (dotted lines) bounds the figure-8 knot complement on one side, the Whitehead link on the other, hence is incompressible.\index{incompressible}}
  \label{Fig:Satellite}
\end{figure}

So if the torus $T$ of \reffig{Satellite} is compressible,\index{compressible} then a compression disk\index{compression disk} must lie on the inside of the torus. But the inside of the torus is homeomorphic to a solid torus in $S^3$ containing a knot complement. In fact, the inside is homeomorphic to the complement of the Whitehead link, the knot in the solid torus shown in \reffig{TKWProof}. Again if there were an embedded compression disk\index{compression disk} for the torus on this side, surgering would give a sphere embedded in the Whitehead link complement, bounding two balls in $S^3$. A similar argument to that above would imply that one of the components of the Whitehead link lies in a ball in a solid torus in the link complement. But then the two link components are unlinked, which is a contradiction: the Whitehead link is nontrivially linked.

Thus the torus in \reffig{Satellite} is incompressible.\index{incompressible}
\end{example}

\begin{definition}\label{Def:BoundaryParallel}
  An embedded surface $F$ in a 3-manifold $M$ is said to be \emph{boundary parallel}\index{boundary parallel} if it can be isotoped into the boundary of $M$.
\end{definition}

\begin{definition}\label{Def:Satellite}
A \emph{satellite knot}\index{satellite knot} is a knot whose complement contains an incompressible\index{incompressible} torus\index{incompressible torus} that is not boundary parallel.\index{boundary parallel}

Equivalently, a satellite knot can be formed as follows. Start with a knot $K'$ in a solid torus $V$, with $K'$ chosen so that it is not contained in a ball in $V$, and $K'$ not isotopic to the core of the solid torus. Let $K''$ be a nontrivial knot in $S^3$. Form the satellite knot (complement) by removing a tubular neighborhood\index{tubular neighborhood} of $K''$ and replacing it with $V-K'$ in a trivial way (that is, attach $V$ so that the meridian curve of $K''$ still bounds a disk in $V$). The knot $K''$ is called the \emph{companion knot}.\index{companion knot} The satellite knot lies in a regular neighborhood of the companion. 
\end{definition}

For orientable surfaces in orientable 3-manifolds, incompressibility is equivalent to the fundamental group injecting in the 3-manifold.

\begin{lemma}\label{Lem:IncompressVsPi1}
  An orientable surface in an orientable 3-manifold is incompressible\index{incompressible} if and only if it is $\pi_1$-injective,\index{$\pi_1$-injective} i.e.\ if the fundamental group of the surface injects into the fundamental group of the 3-manifold under the homomorphism induced by inclusion.

  A nonorientable surface $S$ is $\pi_1$-injective if and only if the boundary of a regular neighborhood of $S$ is an orientable incompressible surface.
\end{lemma}

\begin{proof}
\Refex{ProveIncompressVsPi1} and \refex{NonorIncompressVsPi1}.
\end{proof}

There is an additional notion of incompressibility for properly embedded surfaces with boundary. 

\begin{definition}\label{Def:BoundaryIncompressible}
  Let $F$ be a surface with boundary properly embedded in a 3-manifold $M$. A \emph{boundary compression disk}\index{boundary compression disk} for $F$ is a disk $D$ with $\bdy D$ consisting of two arcs, $\bdy D=\alpha\cup\beta$, such that $\alpha=D\cap F\subset F$ and $\beta=D\cap \bdy M\subset \bdy M$, and such that there is no arc $\gamma$ of $\bdy F$ such that $\gamma\cup\alpha$ bounds a disk on $F$.

  If $F$ admits a boundary compression disk it is \emph{boundary compressible}\index{boundary compressible}, otherwise it is \emph{boundary incompressible}\index{boundary incompressible}. 
\end{definition}

We will give an example of a class of knots that always contains a boundary incompressible surface. First we define the class of knots. 

\begin{figure}
\includegraphics{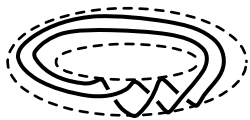}
  \caption{The torus knot $T(2,3)$}
  \label{Fig:TorusKnot}
\end{figure}

\begin{definition}\label{Def:TorusKnot}
  Let $(p,q)$ be relatively prime integers. Then the pair $(p,q)\in\ZZ\times\ZZ\cong H_1(T^2;\ZZ)$ defines a nontrivial simple closed curve on a torus. View the torus $T$ as the boundary of a neighborhood of an unknot in $S^3$. Notice that there is one compression disk\index{compression disk} to the inside of $T$; its boundary is a meridian $m$ for the unknot. There is also a compression disk to the outside of $T$; its boundary is a longitude $\ell$ of the unknot. We choose a basis so that the curve $(p,q)$  has minimal representative intersecting $m$ a total of $|p|$ times, and $\ell$ a total of $|q|$ times, with the signs of $p,q$ determining the direction (right or left handed screw motion).

  A \emph{torus knot}\index{torus knot}, or \emph{$(p,q)$-torus knot}, is the knot in $S^3$ given by the $(p,q)$ curve on the unknotted torus $T$ in $S^3$. It is frequently denoted by $T(p,q)$. \Reffig{TorusKnot} shows an example. 
\end{definition}

\begin{example}\label{Example:TorusKnot}
Suppose $T(p,q)=K$ is a torus knot with $|p|,|q| \geq 2$. The surface $T-K$ is an annulus; we will denote the annulus by $F$. We show that $F$ is boundary incompressible\index{boundary incompressible} in $S^3-K$.

For suppose there is a boundary compression disk\index{boundary compression disk} $D$ for $F$. Since $F\cup K$ is the torus $T$, the disk $D$ lies to one side of $T$. Because $T$ is an unknotted torus in $S^3$, the disk $D$ is either (freely) homotopic to a disk with boundary $m$ or to one with boundary $\ell$. By definition, $T(p,q)$ has intersection number $|p|$ with $m$ and $|q|$ with $\ell$. On the other hand, a boundary compression disk intersects $F$ in exactly one arc, and $K$ in exactly one arc. This is impossible when $|p|,|q|\geq 2$. 
\end{example}

We now fold all our definitions into one.

\begin{definition}\label{Def:EssentialSurface}
A surface $F$ properly embedded in a 3-manifold $M$ is \emph{essential}\index{essential} if one of the following holds.
\begin{enumerate}
\item $F$ is a 2-sphere that does not bound a 3-ball.
\item $F$ is a disk and either $\bdy F \subset \bdy M$ does not bound a disk on $\bdy M$, or $\bdy F\subset \bdy M$ does bound a disk $E$ on $\bdy M$, but $E\cup F$ does not bound a 3-ball. 
\item $F$ is not a disk or sphere, and is incompressible,\index{incompressible} boundary incompressible,\index{boundary incompressible} and not boundary parallel.\index{boundary parallel}
\end{enumerate}
\end{definition}

\begin{definition}\label{Def:NoBadEssential}
  A 3-manifold is said to be:
  \begin{itemize}
  \item \emph{irreducible}\index{irreducible} if it contains no essential\index{essential} 2-sphere, 
  \item \emph{boundary irreducible}\index{boundary irreducible} if it contains no essential disk, 
  \item \emph{atoroidal}\index{atoroidal} if it contains no essential torus, and
  \item \emph{anannular}\index{anannular} if it contains no essential annulus.
  \end{itemize}
\end{definition}

\begin{theorem}\label{Thm:Atoroidal}
  A manifold that contains an embedded essential\index{essential} torus cannot be hyperbolic.
\end{theorem}

The proof of this theorem was part of an exercise in \refchap{CompletionDehnFilling}, but we will go through the argument here.

\begin{proof}
Suppose $M$ contains an essential torus. By \reflem{IncompressVsPi1}, the fundamental group of $M$ contains a $\ZZ\times\ZZ$ subgroup. \Refcor{ZxZSubgroup} implies that the subgroup is generated by two parabolic elements\index{parabolic} fixing the same point on the boundary of $\HH^3$ at infinity. But then the thick--thin decomposition (\refthm{ThinPart} structure of thin part) implies that the torus is parallel to a cusp torus, hence it is boundary parallel.\index{boundary parallel} This contradicts the definition of essential. 
\end{proof}

\begin{corollary}\label{Cor:Satellite}
A satellite knot\index{satellite knot} complement does not admit a hyperbolic structure. \qed
\end{corollary}

\begin{theorem}\label{Thm:Anannular}
  Suppose $M$ is a 3-manifold with torus boundary components whose interior has a complete finite volume hyperbolic metric. Then $M$ cannot contain an essential\index{essential} annulus.\index{anannular}
\end{theorem}

\begin{proof}
  Suppose not. Suppose $M$ is hyperbolic, and $A$ is an essential annulus in $M$. Consider the core curve $\gamma$ of $A$. Note $\gamma$ is isotopic to $\bdy A$, hence $\gamma$ is isotopic into $\bdy M$. Under the hyperbolic structure on the interior of $M$, torus boundary components become cusps. Thus in the hyperbolic structure, $\gamma$ is isotopic into a cusp, so is isotopic to closed curves of arbitrarily small length. Then the loop $\gamma$ corresponds to a covering transformation of $\HH^3\to M$ that is parabolic.\index{parabolic} But then both boundary components of $\bdy A$ correspond to the same parabolic element. Thus the annulus has boundary components given by the same curve in the same cusp. This is possible only if the annulus is boundary parallel.\index{boundary parallel}
\end{proof}

\begin{corollary}\label{Cor:TorusKnot}
A torus knot\index{torus knot} complement does not admit a hyperbolic structure.
\end{corollary}

\begin{proof}
Consider the curve $T(p,q)$ embedded on an unknotted torus $T\subset S^3$. If $(p,q)\in\{(\pm 1,q)\mid q\in\ZZ\}\cup \{(p,\pm 1)\mid p\in\ZZ\}$ then $T(p,q)$ is the unknot in $S^3$, which does not have hyperbolic complement. Otherwise 
the annulus $T-T(p,q)$ is incompressible\index{incompressible} (\refex{TorusKnot}). We showed in \refexamp{TorusKnot} that it is boundary incompressible.\index{boundary incompressible} A similar argument shows it cannot be boundary parallel.\index{boundary parallel} So it is essential. Thus a torus knot is not hyperbolic. 
\end{proof}

Theorems~\ref{Thm:Atoroidal} and~\ref{Thm:Anannular} give surfaces which preclude a knot from being hyperbolic. In fact, an even stronger result is known.

\begin{theorem}[Thurston, Hyperbolization]\label{Thm:SfcesHyperbolic}
  A knot complement admits a complete hyperbolic structure if and only if it is not a satellite knot or a torus knot.

  More generally, a compact 3-manifold with nonempty torus boundary has interior admitting a complete hyperbolic structure if and only if it is irreducible,\index{irreducible} boundary irreducible,\index{boundary irreducible} atoroidal,\index{atoroidal} and anannular.\index{anannular} \qed
\end{theorem}

\Refthm{SfcesHyperbolic} follows from the geometrization theorem for Haken manifolds, which is a very deep result; see \cite{thurston:bulletin}. The proof requires a book of its own (e.g.\ \cite{kapovich}), and we will not include it here. However we will use the theorem to show many knot and link complements are hyperbolic. 

Note that \refthm{SfcesHyperbolic} turns the geometric problem of determining whether a manifold admits a hyperbolic structure into a topological problem of finding surfaces in 3-manifolds, or proving such surfaces cannot exist. It has been applied to show many knots and 3-manifolds admit a hyperbolic structure, although unfortunately it does not give much information on such a structure, beyond the fact that it exists.  We will see some results along these lines in the rest of this chapter.

%%%%%%%%%%%%%%%%%%%%%%%%%%%%%%%%%%%%%%%%%%%%%%%%%%%%%%%%%%%%%%%%%
\section{Torus decomposition, Seifert fibering, and geometrization}

While we are most interested in hyperbolic spaces and hyperbolic geometry, we will also need to identify non-hyperbolic spaces as we encounter them in knot theory. \Refthm{SfcesHyperbolic} gives us a characterization of hyperbolic knots, but we also have the tools now to study non-hyperbolic knots. This section gives a brief overview of the terminology and results that we will need. 

\begin{definition}\label{Def:SeifertFiberedSolidTorus}
  Let $p$ and $q$ be relatively prime integers.
  A \emph{Seifert fibered solid torus}\index{Seifert fibered solid torus}\index{solid torus!Seifert fibered} of type $(p,q)$ is a solid torus $S^1\times D^2$ constructed as the union of disjoint circles, as follows. Begin with a solid cylinder $D^2\times [0,1]$, fibered by intervals $\{x\}\times [0,1]$. Glue the disk $D^2\times\{0\}$ to $D^2\times\{1\}$ by a $2\pi p/q$ rotation. The fiber $\{0\}\times [0,1]$ in $D^2\times[0,1]$ becomes a circle; this is called the \emph{exceptional fiber}\index{exceptional fiber}\index{Seifert fibered solid torus!exceptional fiber}. Every other fiber $\{x\}\times [0,1]$ is glued to $q$ segments to form a circle. These are called \emph{normal fibers}\index{normal fiber}\index{Seifert fibered solid torus!normal fiber}. If $q=1$, the Seifert fibered solid torus is called a \emph{regularly fibered solid torus}.\index{regularly fibered solid torus}\index{Seifert fibered solid torus!regularly fibered solid torus}
\end{definition}

\begin{definition}\label{Def:SeifertFiberedSpace}
  A \emph{Seifert fibered space}\index{Seifert fibered space} is an orientable 3-manifold $M$ that is the union of pairwise disjoint circles, called \emph{fibers}\index{fiber}, such that every fiber has neighborhood diffeomorphic to a fibered solid torus, preserving fibers. 
\end{definition}

\begin{example}\label{Example:S3SeifertFibered}
The 3-sphere $S^3$ is the union of two solid tori $V$ and $W$. For relatively prime integers $(p,q)$, give $V$ the fibering of a Seifert fibered solid torus of type $(p,q)$, and give $W$ the fibering of a Seifert fibered solid torus of type $(q,p)$. Then when we glue $\bdy V$ to $\bdy W$ to form $S^3$, the fibers on the boundaries are identified. Thus $S^3$ is Seifert fibered. 
\end{example}

\begin{example}[Torus knot complements]\label{Example:TorusKnotSeifertFibered}
  The complement of a $(p,q)$-torus knot (\refdef{TorusKnot}) is Seifert fibered, as follows. Take the Seifert fibering of $S^3$ of the previous example. A regular fiber on $\bdy V$ is a $(p,q)$-torus knot. Thus when we remove a fibered solid torus neighborhood of this regular fiber, the result is a Seifert fibered space homeomorphic to the exterior of a torus knot, $S^3-N(T(p,q))$. 
\end{example}

A Seifert fibered space is never hyperbolic. The following theorem follows from work of many people, including work of Casson and Jungreis \cite{CassonJungreis} and Gabai \cite{Gabai:Convergence}. 

\begin{theorem}[Characterization of Seifert fibered spaces]\label{Thm:SeifertFibered}
  A compact orientable irreducible\index{irreducible} 3-manifold $M$ with infinite fundamental group is a Seifert fibered space if and only if $\pi_1(M)$ contains a normal infinite cyclic subgroup. \qed
\end{theorem}

For further information on Seifert fibered spaces, see \cite{Schultens:3Mfld}, \cite{Hatcher:3MfldNotes}, or \cite{scott:geometries}.

The following theorem applies to manifolds that admit an embedded essential torus. It was proved by Jaco, Shalen \cite{jaco-shalen}, and Johannson \cite{johannson}. 

\begin{theorem}[JSJ decomposition]\label{Thm:JSJ}
For any compact irreducible,\index{irreducible} boundary irreducible\index{boundary irreducible} 3-manifold $M$, there exists a (possibly empty) finite collection $\calT$ of disjoint essential\index{essential} tori such that each component of the 3-manifold obtained by cutting $M$ along $\calT$ is either atoroidal\index{atoroidal} or Seifert fibered. Moreover, a minimal such collection $\calT$ is unique up to isotopy. \qed
\end{theorem}

\begin{definition}\label{Def:JSJ}
  The minimal collection of tori $\calT$ as in \refthm{JSJ} is called the \emph{JSJ-decomposition}\index{JSJ-decomposition} of $M$, or sometimes the \emph{torus decomposition}\index{torus decomposition} of $M$. The union of the Seifert fibered pieces of $M$ cut along $\calT$ is called the \emph{characteristic submanifold}\index{characteristic submanifold} of $M$. 
\end{definition}

In \refexamp{Satellite}, the torus decomposition consists of the single essential\index{essential} torus shown in \reffig{Satellite}. Cutting along it splits the 3-manifold into two hyperbolic pieces, hence the characteristic submanifold is empty.

Thurston's hyperbolization theorem implies if the JSJ-decomposition of $M$ is nontrivial, then atoroidal components of $M$ cut along $\calT$ are hyperbolic. More generally, even in the closed case we now know the following theorem. 

\begin{theorem}[Geometrization of closed 3-manifolds] \label{Thm:Geometrization}\index{Geometrization theorem for closed 3-manifolds}
Let $M$ be a closed, orientable, irreducible\index{irreducible} 3-manifold.
\begin{enumerate}
\item If $\pi_1(M)$ is finite, then $M$ is spherical; i.e.\ $M$ is homeomorphic to $S^3/\Gamma$ where $\Gamma$ is a finite subgroup of $O(4)$ acting on $S^3$ without fixed points.
\item If $\pi_1(M)$ is infinite and contains a $\ZZ\times\ZZ$ subgroup, then $M$ is either Seifert fibered or contains an incompressible\index{incompressible} torus\index{incompressible torus} (so is not hyperbolic).
\item If $\pi_1(M)$ is infinite and contains no $\ZZ\times\ZZ$ subgroup, then $M$ is hyperbolic. \qed
\end{enumerate}
\end{theorem}

The second part of the theorem follows from work of Casson and Jungreis \cite{CassonJungreis} and Gabai \cite{Gabai:Convergence}. The first and third parts were proved by Perelman \cite{perelman02}, \cite{perelman03}.

%%%%%%%%%%%%%%%%%%%%%%%%%%%%%%%%%%%%%%%%%%%%%%%%%%%%%%%%%%%%%%%%%
\section{Normal surfaces, angled polyhedra, and hyperbolicity}

In this section, we will use \refthm{SfcesHyperbolic} to prove that many manifolds are hyperbolic. We will be considering 3-manifolds
that admit an ideal polyhedral decomposition, for example a decomposition into ideal tetrahedra, but also more general ideal polyhedra as in \refchap{Fig8Decomp} and \refchap{TwistKnots}. We will see that we need only consider surfaces that intersect the polyhedra in simple ways: in disks with well-behaved boundaries.

\subsection{Normal surfaces}

To describe nice positions of embedded surfaces in polyhedra, we give the following definition. 

\begin{figure}
%% Accompanies image file 'F8-04-Normal.eps' (pdf, eps, ps)
\begingroup%
  \makeatletter%
  \providecommand\color[2][]{%
    \errmessage{(Inkscape) Color is used for the text in Inkscape, but the package 'color.sty' is not loaded}%
    \renewcommand\color[2][]{}%
  }%
  \providecommand\transparent[1]{%
    \errmessage{(Inkscape) Transparency is used (non-zero) for the text in Inkscape, but the package 'transparent.sty' is not loaded}%
    \renewcommand\transparent[1]{}%
  }%
  \providecommand\rotatebox[2]{#2}%
  \newcommand*\fsize{\dimexpr\f@size pt\relax}%
  \newcommand*\lineheight[1]{\fontsize{\fsize}{#1\fsize}\selectfont}%
  \ifx\svgwidth\undefined%
    \setlength{\unitlength}{311.24349976bp}%
    \ifx\svgscale\undefined%
      \relax%
    \else%
      \setlength{\unitlength}{\unitlength * \real{\svgscale}}%
    \fi%
  \else%
    \setlength{\unitlength}{\svgwidth}%
  \fi%
  \global\let\svgwidth\undefined%
  \global\let\svgscale\undefined%
  \makeatother%
  \begin{picture}(1,0.26146868)%
    \lineheight{1}%
    \setlength\tabcolsep{0pt}%
    \put(0,0){\includegraphics[width=\unitlength]{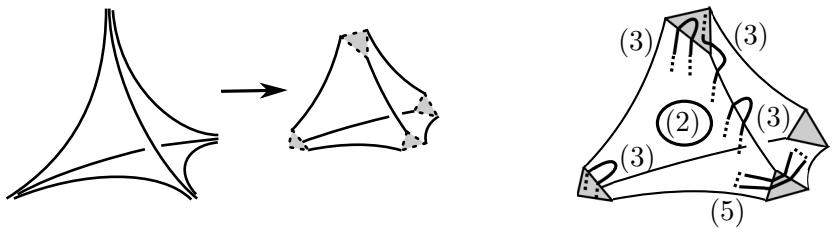}}%
    \put(0.79727682,0.11133845){\color[rgb]{0,0,0}\makebox(0,0)[lt]{\lineheight{0}\smash{\begin{tabular}[t]{l}(2)\end{tabular}}}}%
    \put(0.90697602,0.11847872){\color[rgb]{0,0,0}\makebox(0,0)[lt]{\lineheight{0}\smash{\begin{tabular}[t]{l}(3)\end{tabular}}}}%
    \put(0.74210271,0.07044465){\color[rgb]{0,0,0}\makebox(0,0)[lt]{\lineheight{0}\smash{\begin{tabular}[t]{l}(3)\end{tabular}}}}%
    \put(0.74080445,0.21000297){\color[rgb]{0,0,0}\makebox(0,0)[lt]{\lineheight{0}\smash{\begin{tabular}[t]{l}(3)\end{tabular}}}}%
    \put(0.87954542,0.21725698){\color[rgb]{0,0,0}\makebox(0,0)[lt]{\lineheight{0}\smash{\begin{tabular}[t]{l}(3)\end{tabular}}}}%
    \put(0.85115274,0.00423578){\color[rgb]{0,0,0}\makebox(0,0)[lt]{\lineheight{0}\smash{\begin{tabular}[t]{l}(5)\end{tabular}}}}%
  \end{picture}%
\endgroup%

  \caption{Left: truncating an ideal polyhedron. Boundary faces are shaded, boundary edges dashed, (regular) faces are white, and (regular) edges are solid black. Right: examples of curves that cannot be the boundary of a normal disk, along with the number of the property of \refdef{NormalDisk} that they violate}
  \label{Fig:NormalDisk}
\end{figure}

\begin{definition}\label{Def:NormalDisk}
Let $P$ be an ideal polyhedron. Truncate the ideal vertices of $P$, so they become \emph{boundary faces}, and denote the truncated polyhedron by $\overline{P}$. The edges between (regular) faces and boundary faces are called \emph{boundary edges}. See \reffig{NormalDisk}, left. 
  
Let $D$ be a disk embedded in $\overline{P}$ with $\bdy D \subset \bdy \overline{P}$. We say that $D$ is \emph{normal}\index{normal!disk} if it satisfies the following conditions.
  \begin{enumerate}
  \item $\bdy D$ meets the faces, boundary faces, edges, and boundary edges of $\overline{P}$ transversely.
  \item $\bdy D$ does not lie entirely on a single face or boundary face of $\overline{P}$.
  \item Any arc of intersection of $\bdy D$ with a face of $\overline{P}$ does not have both endpoints on the same edge, or on the same boundary edge, or on an adjacent edge and boundary edge. Similarly, any arc of intersection of $\bdy D$ with a boundary face does not have both endpoints on the same boundary edge. 
  \item $\bdy D$ meets any edge at most once.
  \item $\bdy D$ meets any boundary face at most once.
  \end{enumerate}
\end{definition}

\Reffig{NormalDisk} illustrates some of these conditions. 

\begin{definition}\label{Def:NormalForm}
  A surface is in \emph{normal form},\index{normal form}\index{normal} with respect to a polyhedral decomposition, or is \emph{normal}, for short, if it intersects the (truncated) polyhedra in a collection of normal disks.
\end{definition}

Normal surfaces in 3-manifolds are well-studied objects, and the following theorem is classical, dating back to work of Kneser in the late 1920s \cite{kneser}, and Haken and Schubert in the 1960s \cite{Haken,Schubert:NormalSfces}, and many others since then.

\begin{theorem}\label{Thm:NormalForm}
Suppose $M$ admits an ideal polyhedral decomposition.
  
If $M$ contains an essential\index{essential} 2-sphere, then it contains one in normal form.\index{normal}

If $M$ is irreducible\index{irreducible} and $M$ contains an essential disk, then it contains one in normal form.
  
If $M$ is irreducible\index{irreducible} and boundary irreducible,\index{boundary irreducible} and contains an essential surface, then that surface can be isotoped in $M$ to meet the polyhedra in normal form.\index{normal} 
\end{theorem}

\begin{proof}
Let $S$ be an essential surface in $M$. We may isotope $S$ so that it intersects faces, boundary faces, edges, and boundary edges of the truncated polyhedra of $M$ transversely. Let $f$ denote the number of times $S$ meets a face or boundary face, and let $e$ denote the number of times $S$ meets an edge or boundary edge. The pair $(f,e)$ is called the \emph{complexity}\index{complexity of essential surface} of $S$ in $M$, and we order it lexicographically. We will adjust $S$ to remove intersections with the polyhedra that violate normality while reducing its complexity. Since the complexity is finite, it follows that a finite number of adjustments give the result. 

First we claim we can adjust $S$ so that it meets the polyhedra only in disks, lowering the complexity. If $S$ is a sphere, this is done by replacing $S$. If $M$ is irreducible,\index{irreducible} then this is done by isotopy of $S$. The argument is similar to arguments below and so we leave it as \refex{EssentialSfcesMeetPolyinDisks}.

We now assume that the components of intersection of $S$ with a truncated polyhedron $P$ are all disks. Suppose that $\bdy(S\cap P)$ contains a simple closed curve of intersection contained entirely in a face or boundary face. Then there must be an innermost such curve $\gamma$, and $\gamma$ bounds a disk $E$ inside that face disjoint from $S$.

If $S$ is a 2-sphere, then surger\index{surger} $S$ along $E$, obtaining two 2-spheres, $S'$ and $S''$, which we push slightly to be disjoint from a neighborhood of $E$. Thus $S'$ and $S''$ have strictly fewer intersections with $\bdy P$. Either $S'$ or $S''$ must still be essential, else $S$ could not be essential, so say $S''$ is essential. Replace $S$ with $S''$. Then $S''$ has strictly smaller complexity than $S$. Repeating a finite number of times, we may assume $S$ does not meet faces or boundary faces of $P$ in closed curves in this case. 

If $M$ is irreducible\index{irreducible} and $S$ is a disk, then as before surger along $E$. This gives two surfaces, a disk $S'$ and a sphere $S''$, both of strictly smaller complexity than $S$. Because $M$ is irreducible, the sphere $S''$ bounds a ball, and we may isotope $S$ through that ball to the disk $S'$. Repeat for each closed curve of intersection of $S$ with faces, removing all such intersections. 
  
If $M$ is irreducible\index{irreducible} and boundary irreducible,\index{boundary irreducible} then $S$ is not a sphere or disk. Because $S$ is essential, $\gamma$ bounds a disk $E'$ in $S$. Then the sphere $E\cup E'$ must bound a ball by irreducibility. Isotope $S$ through this ball, removing the intersection $\gamma$ and reducing complexity. Repeating finitely many times eliminates all closed curves of intersection of $S$ with faces of $P$.

Now suppose an arc of intersection of $S$ with a face or boundary face of $P$ has both its endpoints on the same edge or boundary edge, or on an edge and adjacent boundary edge. Then there must be an outermost such arc $\alpha$, bounding a disk $E$ on that face or boundary face, with $E$ disjoint from $S$. In the case that the face is a boundary face, or the endpoints of $S$ do not both lie on a boundary edge, we may slide $S$ through a neighborhood of $E$ to isotope across the edge, decreasing complexity.

When the arc $\alpha$ lies on a regular face, and both endpoints of $\alpha$ lie on boundary edges, then we have to take more care since we cannot isotope the surface $S$ past the boundary edge without changing its topology. In this case, we know $S$ is a surface with boundary, so not a sphere, so we are in the case that $M$ is irreducible.\index{irreducible} If $S$ is a disk, surger along $E$ and push off the face slightly, obtaining two disks with lower complexity. One of them must be essential, since $S$ is essential. Replace $S$ with this essential disk. 

If $S$ is not a disk, then $M$ is irreducible\index{irreducible} and boundary irreducible.\index{boundary irreducible} Since $S$ is essential, $E$ is not a boundary compression disk\index{boundary compression disk} for $S$, thus $\alpha$ bounds a disk $E'$ on $S$. Then $E\cup E'$ is a disk with boundary on $\bdy M$. Because $M$ is boundary irreducible, it must be parallel to a disk $E''$ on $\bdy M$. Then $E\cup E'\cup E''$ is a sphere, so bounds a ball, and we may isotope $S$ through this ball to remove the intersection $\alpha$, strictly decreasing complexity. 

Finally, we need to show for each polyhedron $P$, any disk of $S\cap P$ has boundary meeting each edge of $P$ at most once, and meeting each boundary face of $P$ at most once. We leave these as exercises~\ref{Ex:NormalDisk} and~\ref{Ex:NormalBdry}.
\end{proof}

\subsection{Angle structures and combinatorial area}

We are interested in 3-manifolds that are hyperbolic. If a 3-manifold admits a decomposition into ideal tetrahedra,
recall from \refthm{Gluing} (edge gluing equations) and \refdef{CompletenessEquations} (completeness equations) that a complete hyperbolic structure satisfies a system of nonlinear equations. If we take the log of the edge gluing equations, the complex product becomes a sum of real and imaginary parts:
\[ \log\left(\prod z(e_j)\right) = \sum ( \log|z(e_j)| + i\,\Arg\, z(e_j) ) = 2\pi\,i. \]
The imaginary parts encode relations on dihedral angles of the tetrahedra. If we ignore the real part, then finding solutions to the imaginary parts involves finding dihedral angles that satisfy a system of linear equations. Thus by considering dihedral angles alone, we reduce a complicated nonlinear problem to a linear problem. This can significantly simplify computations.

\begin{definition}\label{Def:AngleStructures}
An \emph{angle structure}\index{angle structure} on an ideal triangulation $T$ of a manifold $M$ is a collection of (interior) dihedral angles, one for each edge of each tetrahedron, satisfying the following conditions.
\begin{enumerate}
\item[(0)]\label{Itm:OppEdges} Opposite edges of the tetrahedron have the same angle.
\item[(1)]\label{Itm:Range} Dihedral angles lie in $(0,\pi)$.
\item[(2)]\label{Itm:VertexSum} The sum of angles around any ideal vertex of any tetrahedron is $\pi$.
\item[(3)]\label{Itm:EdgeSum} The sum of angles around any edge class of $M$ is $2\pi$.
\end{enumerate}
The set of all angle structures for triangulation $T$ is denoted by $\mathcal{A}(T)$. 
\end{definition}

Conditions~(0) and~(1) are required for nonsingular tetrahedra. Condition~(2) ensures that a triangular cross-section of any ideal vertex of any tetrahedron is actually a Euclidean triangle. Finally, condition~(3) is the imaginary part of the edge gluing equations. Note we have not included the completeness equations. For many results, we don't need them!

An angle structure\index{angle structure} on an ideal tetrahedron uniquely determines the shape of that tetrahedron. If it has assigned dihedral angles $\alpha$, $\beta$, $\gamma$ in clockwise order, then there is a unique hyperbolic tetrahedron with those dihedral angles, and its edge invariant corresponding to the edge with angle $\alpha$ can be shown to be
\begin{equation}\label{Eqn:AngleEdgeInvariant}
  z(\alpha) = \frac{\sin\gamma}{\sin\beta}e^{i\alpha}.
\end{equation}
Thus if a triangulation has an angle structure (not all do), we can think of the manifold as being built of hyperbolic tetrahedra.

\begin{figure}
  \includegraphics{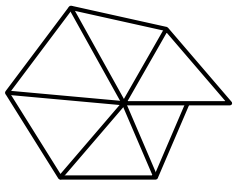}
  \caption{Angle structures typically have shearing singularities.\index{angle structure}}
  \label{Fig:Shearing}
\end{figure}

Note that because we have discarded the completeness equations, an angle structure\index{angle structure} typically will not give a complete hyperbolic structure. In fact, because we have discarded the nonlinear part of the edge gluing equations, an angle structure typically won't even give a hyperbolic structure. There will likely be shearing singularities\index{shearing} around each edge, as in \reffig{Shearing}.

Even so, much useful information can be extracted from angle structures,\index{angle structure} which we will see in later chapters. In this chapter, we will show that if an angle structure exists on a triangulation of a manifold, then the manifold admits (some) hyperbolic structure.

The idea of an angled triangulation can be generalized to ideal polyhedra as well. First, we need to assign to each edge a dihedral angle, which is a number lying in the range $(0,\pi)$. Once that has been done, we can measure a combinatorial area of normal disks embedded in the polyhedra, as follows.

\begin{definition}\label{Def:CombinatorialArea}
Let $D$ be a normal disk\index{normal} in a (truncated) ideal polyhedral decomposition of $M$, such that each ideal edge of $M$ has been assigned an interior dihedral angle in the range $(0,\pi)$. Let $\alpha_1, \dots, \alpha_n$ be the angles assigned to the ideal edges met by $\bdy D$. Then the \emph{combinatorial area} of $D$\index{combinatorial area}\index{combinatorial area!disk} is defined as:
\[ a(D) = \sum_{i=1}^n (\pi-\alpha_i) - 2\pi + \pi|\bdy D \cap \bdy M|. \]
Here $|\bdy D \cap \bdy M|$ indicates the number of components of intersection of $\bdy D$ with boundary faces.

If $S$ is a surface in normal form,\index{normal} the \emph{combinatorial area}\index{combinatorial area}\index{combinatorial area!surface} of $S$ is defined to be the sum of combinatorial areas of the normal disks making up $S$.
\end{definition}

Note in the case $D$ is contained in a hyperbolic plane, meeting each edge of the polyhedron orthogonally, the combinatorial area\index{combinatorial area} of $D$ agrees with the actual hyperbolic area (\refex{AreaPolygon}). 

We can now generalize the idea of an angle structure\index{angle structure} on a triangulation to an angle structure on an ideal polyhedral decomposition.

\begin{definition}\label{Def:AnglePolyhedra}
An \emph{angled polyhedral structure}\index{angled polyhedral structure} on a 3-manifold $M$ is a decomposition of $M$ into ideal polyhedra, along with a collection of (interior) dihedral angles, one for each edge of each polyhedra, that satisfy the following conditions.
  \begin{enumerate}
  \item Each dihedral angle lies in the range $(0,\pi)$. 
  \item Every normal\index{normal} disk has non-negative combinatorial area.\index{combinatorial area}
  \item Interior angles around an edge sum to $2\pi$. 
  \end{enumerate}
\end{definition}

\begin{example}
An angle structure\index{angle structure} on an ideal triangulation of $M$ is an example of an angled polyhedral structure.\index{angled polyhedral structure} To show this, suppose we have an angle structure on a triangulation of $M$. Then by \refdef{AngleStructures}, the dihedral angles are in the correct range, and interior angles around edges must sum to $2\pi$ to satisfy the definition of an angle structure. So we need only consider normal disks,\index{normal} and show that each has non-negative combinatorial area.\index{combinatorial area} Two examples of normal disks are shown in \reffig{BdyTriangleBdyBigon}.

\begin{figure}
\includegraphics{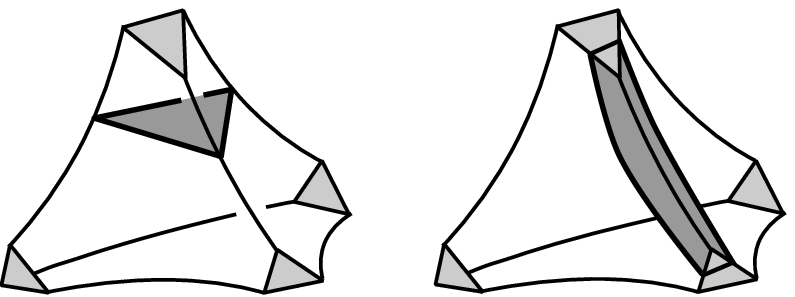}
  \caption{Normal disks:\index{normal} a vertex triangle\index{vertex triangle} and a boundary bigon\index{boundary bigon}}
  \label{Fig:BdyTriangleBdyBigon}
\end{figure}

The triangle in the figure has combinatorial area\index{combinatorial area}
\[ a(D) = (\pi-\alpha) + (\pi-\beta) + (\pi-\gamma) -2\pi = \pi-(\alpha+\beta+\gamma) = 0,\]
since $\alpha$, $\beta$, $\gamma$ encircle an ideal vertex of a tetrahedron. We call this normal disk a \emph{vertex triangle}\index{vertex triangle}. 

The other normal disk shown also has zero combinatorial area:\index{combinatorial area}
\[ a(D) = 0 - 2\pi + \pi\cdot 2.\]
This disk is called a \emph{boundary bigon}\index{boundary bigon}.

\begin{lemma}\label{Lem:CombAreaTriangulation}
Let $M$ be a triangulated 3-manifold with an angle structure.\index{angle structure} Then the combinatorial area\index{combinatorial area} of any normal disk\index{normal} $D$ in an ideal tetrahedron of $M$ is non-negative. It is zero if and only if $D$ is a vertex triangle or a boundary bigon.
\end{lemma}

\begin{proof}
As we have seen above, the combinatorial areas of boundary bigons and vertex triangles are zero.

If $D$ is a normal disk\index{normal} that meets at least two boundary faces, its combinatorial area is at least the sum $\sum (\pi-\alpha_i)$, and any term $(\pi-\alpha_i)$ is positive, so $a(D)\geq 0$ in this case.

If $D$ meets exactly one boundary face, then because $\bdy D$ cannot meet edges adjacent to a boundary edge, it must meet an opposite edge in each of the triangles on either side of the boundary face. These cannot be opposite edges in the tetrahedron. Hence the combinatorial area is
\[ a(D) \geq \pi-\alpha + \pi-\beta -2\pi + \pi = \pi-\alpha-\beta = \gamma>0,\]
where here we let $\alpha$, $\beta$, and $\gamma$ denote the angles of an ideal tetrahedron with $\alpha+\beta+\gamma=\pi$. Thus if $D$ meets just one boundary face, $a(D)$ is strictly positive. 

If $D$ meets no boundary faces, then it is either a vertex triangle or a quad separating two opposite edges. In the first case, the combinatorial area is zero. In the second, the combinatorial area is
\[ a(D) = 2(\pi-\alpha)+2(\pi-\beta) -2\pi = 2(\pi-\alpha-\beta) = 2\gamma>0.\]
In all cases, the combinatorial area is non-negative.
\end{proof}

Thus we have shown:
\begin{theorem}\label{Thm:AngleStructIsAngledPoly}
  An angle structure\index{angle structure} on an ideal triangulation of $M$ is an angled polyhedral structure.\index{angled polyhedral structure}\qed
\end{theorem}
\end{example}

\begin{lemma}[Gauss--Bonnet]\label{Lem:GaussBonnet}
A normal surface\index{normal} $S$ in an angled polyhedral structure\index{angled polyhedral structure} satisfies
\[ a(S) = -2\pi\chi(S).\]
\end{lemma}

\begin{proof}
Recall that $\chi(S)$, the Euler characteristic of $S$,\index{Euler characteristic} is given by $\chi(S)=v-e+f$, where $v$ is the number of vertices in a polygonal decomposition of $S$, $e$ is the number of edges, and $f$ is the number of faces. In our case, the intersection of $S$ with the polyhedra determines a polygonal decomposition. Then $f$ is the number of normal disks,\index{normal} or intersections of $S$ with interiors of the polyhedra. The value $e$ is the number of intersections of $S$ with faces of the polyhedra, and $v$ is the number of intersections of $S$ with ideal edges of the polyhedra. Intersections of $S$ with boundary edges and boundary faces do not affect Euler characteristic at all.

By definition,
\begin{align*}
  a(S) & = \sum_{D} a(D) = \sum_{D} \left(\sum_i(\pi-\alpha_i)+\pi|\bdy D\cap\bdy M| -2\pi\right) \\
  & = \pi \sum_{D} \left( \left(\sum_i 1\right)+|\bdy D\cap\bdy M|\right) - \sum_{D}\sum_i\alpha_i - \sum_{D} 2\pi,
\end{align*}
where the sum is over normal disks\index{normal} $D\subset S$. 
Note that the last term in the sum is $-2\pi f$, since we add $-2\pi$ for each normal disk of intersection in $S\cap P$.

The term $\sum_D\sum_i \alpha_i$ gives the sum of all interior angles met by the surface $S$. This is $2\pi v$.

Finally, we claim that $(\sum_i 1+|\bdy D\cap\bdy M|)$ counts the number of edges in faces (not boundary faces) in the normal disk\index{normal} $D$. To see this, orient $\bdy D$ and give each edge the corresponding direction. Its initial endpoint is either on an ideal edge of the polyhedron or a boundary edge. The sum counts all the initial endpoints of edges of $\bdy D$ on faces, without counting initial endpoints of edges on boundary faces. Denote this by
\[ \left(\sum_i 1 + |\bdy D\cap\bdy M|\right) = e(D).\]

Now take the sum $\sum_D e(D)$. The sum over all normal disks\index{normal} counts each edge exactly twice, so its value is $2 e$. Thus $\pi \sum_D (\sum_i 1 + |\bdy D \cap \bdy M|) = 2\pi e$. 

Putting it together, we find
\[ a(S) = 2\pi e - 2\pi v - 2\pi f = -2\pi\chi(S).\qedhere\]
\end{proof}

\subsection{Hyperbolicity}

\begin{theorem}\label{Thm:HypAngleStruct}
Let $M$ be a manifold admitting an angled polyhedral structure.\index{angled polyhedral structure} Then $M$ is irreducible\index{irreducible} and boundary irreducible,\index{boundary irreducible} and its boundary consists of tori.

Moreover, if the angled polyhedral structure is actually an angle structure on a triangulation of $M$, then $M$ is atoroidal\index{atoroidal} and anannular.\index{anannular} Hence any manifold admitting an angle structure\index{angle structure} is hyperbolic.
\end{theorem}

\begin{proof}
Suppose $S$ is an essential\index{essential} sphere in $M$. By \refthm{NormalForm}, we can put $S$ into normal form\index{normal form} with respect to the polyhedral decomposition of $M$, and obtain a combinatorial area\index{combinatorial area} for $S$. By definition of an angled polyhedral structure,\index{angled polyhedral structure} each normal disk\index{normal} of $S$ has non-negative combinatorial area, so $a(S)$ is non-negative. But \reflem{GaussBonnet} implies $a(S) = -4\pi$. This contradiction proves that $M$ is irreducible.\index{irreducible} A similar argument shows that $S$ cannot be an essential\index{essential} disk, so $M$ is boundary irreducible.\index{boundary irreducible}

Now consider the boundary components of $M$. These are obtained by gluing boundary faces. Pushing in slightly, we find that $\bdy M$ is parallel to a normal surface\index{normal} made up of boundary parallel\index{boundary parallel} disks. By the definition of an angle structure,\index{angle structure} \refdef{AngleStructures}, and the definition of combinatorial area, \refdef{CombinatorialArea}, it follows that each such disk has combinatorial area zero. So $\bdy M$ has combinatorial area\index{combinatorial area} zero. Because $\bdy M$ consists of closed surfaces in an orientable manifold, it must be a disjoint union of tori. 

Now suppose the angled polyhedral structure\index{angled polyhedral structure} is an angle structure\index{angle structure} on a triangulation of $M$, and suppose $S$ is an essential torus. Then $S$ can be put into normal form\index{normal} by \refthm{NormalForm}, and \reflem{GaussBonnet} implies that $a(S)=0$, so each normal disk of $S$ has zero combinatorial area.\index{combinatorial area} Then \reflem{CombAreaTriangulation} implies that each normal disk is a vertex triangle\index{vertex triangle} or a boundary bigon.\index{boundary bigon} Since $S$ is a closed surface embedded in $M$, it does not meet boundary faces of $M$, hence each normal disk\index{normal} is a vertex triangle. But vertex triangles join to form the boundary of $M$, hence a component of $\bdy M$ is a torus, and $S$ is parallel to $\bdy M$. This contradicts the fact that $S$ is essential.

Finally, suppose $S$ is an essential annulus in the manifold $M$ with a triangulation and angle structure.\index{angle structure} Then again $a(S)=0$, so $S$ is made up of vertex triangles\index{vertex triangle} and boundary bigons.\index{boundary bigon} There must be at least one boundary bigon. This must be glued to another boundary bigon, since the edges on the face of the tetrahedron run between boundary edges. Then $S$ is made up entirely of boundary bigons. The only possibility is that the boundary bigons encircle a single edge of the triangulation of $M$. This is not incompressible.\index{incompressible} So $S$ is anannular.\index{anannular}

The fact that $M$ is hyperbolic now follows from \refthm{SfcesHyperbolic}. 
\end{proof}

%%%%%%%%%%%%%%%%%%%%%%%%%%%%%%%%%%%%%%%%%%%%%%%%%%%%%%%%%%%%%%%%%
\section{Pleated surfaces and a 6-theorem}

When a 3-manifold admits a hyperbolic structure, then that structure can often be used to induce a hyperbolic structure on a properly embedded essential\index{essential} surface with punctures. If the surface is totally geodesic inside the 3-manifold, such as for the white and shaded surfaces in a fully augmented link,\index{fully augmented link} \refcor{AugmentedGeodSfces}, then the induced hyperbolic structure on the surface is unique. But usually a properly embedded surface is not totally geodesic. In this case, frequently we may still straighten the surface.

\begin{definition}\label{Def:HomotopicBdyIncompr}
  An embedded surface $S$ in a 3-manifold $M$ is \emph{homotopically boundary incompressible},\index{homotopically boundary incompressible} or \emph{homotopically $\bdy$-incompressible}\index{homotopically $\bdy$-incompressible}, if for any properly embedded arc $\alpha$ in $S$ that is not homotopic rel endpoints into $\bdy S$, the arc $\alpha$ in $M$ is not homotopic rel endpoints into $\bdy M$. That is, a nontrivial arc in $S$ remains nontrivial in $M$. (This is also sometimes called \emph{algebraically $\bdy$-incompressible}.)\index{algebraically $\bdy$-incompressible}
\end{definition}

Notice that a homotopically $\bdy$-incompressible surface is boundary incompressible.\index{boundary incompressible} However, now arcs may be immersed and a homotopy gives a singular disk. 

\begin{lemma}\label{Lem:PrePleating}
  Let $S$ be a surface with non-empty boundary properly embedded in a 3-manifold $M$ whose interior admits a complete 
  hyperbolic structure. Suppose $S$ is homotopically $\bdy$-incompressible. Then the ideal edges of any ideal triangulation of $S$ can be homotoped to be geodesics in $M$. Similarly, each ideal triangle\index{ideal triangle} can be homotoped to be totally geodesic in $M$. 
\end{lemma}

\begin{proof}
Any edge of an ideal triangulation on $S$ is homotopically non-trivial on $S$. Because $S$ is homotopically $\bdy$-incompressible, each edge must also be homotopically non-trivial in $M$. Thus it lifts to an arc with distinct endpoints in the universal cover $\HH^3$ of $M$. Such an arc is homotopic to a unique geodesic in $\HH^3$. The image of this geodesic (and the homotopy) under the covering map gives the desired geodesic (and homotopy) in $M$.

Now for any ideal triangle of $S$, the edges of the triangle are homotopic to geodesics in $M$. Lift one edge to be a geodesic in $\HH^3$. Because the interior of the triangle is homotopically trivial in $S$ and in $M$, it lifts to the interior of a triangle in $\HH^3$, and thus we may choose lifts of the other two edges of the triangle such that the three bound a unique totally geodesic triangle, homotopic to a lift of $S$, in $\HH^3$. The image of this triangle (and homotopy) under the covering map gives the desired totally geodesic triangle in $M$. 
\end{proof}

We often refer to the homotopy of \reflem{PrePleating} as \emph{straightening}\index{straightening edges and triangles}. 
After straightening edges and triangles in \reflem{PrePleating}, note that the surface will typically be bent along the geodesic ideal edges. In addition, note that the homotopy is not at all guaranteed to leave the surface embedded. Thus after such a straightening the surface is frequently only immersed, not embedded. However, the process still can give significant geometric information.

\begin{definition}\label{Def:PleatedSurface}
  A \emph{pleated surface}\index{pleated surface} in a hyperbolic 3-manifold $M$ is a pair $(S, \varphi)$ consisting of a surface $S$ with complete hyperbolic structure, and a local isometry $\varphi\from S\to \varphi(S)\subset M$ such that each point in $S$ lies in a geodesic mapped by $\varphi$ to a geodesic. When $(S,\varphi)$ is a pleated surfaces, will also sometimes say that the image $\varphi(S) \subset M$ is pleated. 
\end{definition}

\begin{proposition}\label{Prop:Pleating}
A homotopically $\bdy$-incompressible surface (with non-empty boundary) properly embedded in a hyperbolic 3-manifold can be pleated. That is, it is homotopic to the image of a local isometry $\varphi\from S\to\varphi(S)$ coming from a pleated surface.\index{pleated surface}
\end{proposition}

\begin{proof}
This follows almost immediately from \reflem{PrePleating}; we only need to describe the hyperbolic structure on $S$. The straightening process of \reflem{PrePleating} maps each ideal triangle\index{ideal triangle} of $S$ to a hyperbolic ideal triangle. We define a hyperbolic structure on $S$ by taking isometric ideal triangles in $S$ and attaching them along edges such that the map from $S$ into the homotopic surface in $M$ is an isometry. Thus the hyperbolic structure on $S$ can be viewed as pulling the triangles of $S$ out of $M$ and lining them up, without bending, in $\HH^2$. This gives a fundamental domain for $S$. The surface is obtained by applying isometric gluing maps on edges. 
\end{proof}

%%%%%%%%%%%%%%%%%%%%%%%%%%%%%%%%%%%%%%%%%%%%%%%%%%%%%%%%%%%%%%%%%
Let $M$ be a compact 3-manifold with a torus boundary component $T$ such that the interior of $M$ admits a complete hyperbolic structure. Then recall that in the complete hyperbolic structure, the boundary of any embedded horoball neighborhood of the cusp corresponding to $T$ inherits a Euclidean structure\index{Euclidean structure} (\refthm{EuclidCusp}).

\begin{definition}\label{Def:SlopeLength}
  Recall that an isotopy class of simple closed curves on the torus $T$ is a \emph{slope}.\index{slope} In a Euclidean metric on $T$, any slope can be isotoped to a geodesic. The \emph{slope length} of $s$\index{slope length} is defined to be the length of such a geodesic, denoted $\ell(s)$.

  For $M$ a compact 3-manifold with torus boundary components and hyperbolic interior, a fixed choice of cusp neighborhoods gives a fixed Euclidean structure\index{Euclidean structure} on each torus boundary component. The slope length of $s$ is measured in this fixed Euclidean structure. 
\end{definition}

%%%%%%%%%%%%%%%%%%%%%%%%%%%%%%%%%%%%%%%%%%%%%%%%%%%%%%%%%%%%%%%%%

Hyperbolic geometry and pleated surfaces\index{pleated surface} can be used to give a proof of the following theorem.

\begin{theorem}[A 6-theorem]\label{Thm:6Theorem}\index{6-theorem, weaker}
  Suppose $M$ is a compact manifold with torus boundary components, such that the interior of $M$ admits a complete hyperbolic structure. Let $s_1, \dots, s_n$ be slopes on distinct boundary components of $M$ such that each slope length\index{slope length} $\ell(s_i)$ is strictly larger than $6$ on a collection of disjoint embedded horospherical tori for $M$. Then the manifold $M(s_1, \dots, s_n)$ obtained by Dehn filling\index{Dehn filling} $M$ along slopes $s_1, \dots,s_n$ is irreducible,\index{irreducible} boundary irreducible,\index{boundary irreducible} anannular,\index{anannular} and atoroidal.\index{atoroidal}
\end{theorem}

\Refthm{6Theorem} is \emph{a} 6-theorem, but not exactly \emph{the} 6-theorem.\index{6-theorem} The 6-theorem, proved independently and simultaneously by Agol \cite{agol:bounds} and Lackenby \cite{lackenby:word}, is stronger. Their theorem states that if slope lengths are at least six, then the Dehn filled manifold cannot be reducible,\index{reducible 3-manifold} toroidal, Seifert fibered or have finite fundamental group. The geometrization theorem implies such a manifold must be hyperbolic. When the Dehn filling gives a closed manifold, our 6-theorem, in \refthm{6Theorem}, does not rule out Seifert fibered or finite fillings. However, in the case that the manifold we obtain after Dehn filling\index{Dehn filling} still has boundary, it will be hyperbolic by Thurston's hyperbolization theorem, \refthm{SfcesHyperbolic}. Also, our proof uses a little less machinery, while still giving a nice introduction to the geometric arguments involved. We highly recommend reading the original papers \cite{agol:bounds} and \cite{lackenby:word}.

Our proof of \refthm{6Theorem} will follow three simple steps. First, assuming that $M(s)$ is reducible,\index{reducible 3-manifold} boundary reducible, annular, or toroidal, we show that there is a punctured 2-sphere or punctured torus $S$ embedded in $M$ that is essential,\index{essential} with boundary components on $\bdy M$ tracing out slopes $s_i$. Second, because $S$ is essential, it can be pleated,\index{pleated surface} and it inherits a hyperbolic metric and cusp neighborhoods from the metric on $M$ and its embedded horocusps. Third, arguments in hyperbolic geometry show that the slope lengths\index{slope length} are at most six.

\begin{lemma}\label{Lem:EssentialPunctSurface}
Let $M$, $s_1, \dots, s_n$ be as in the statement of \refthm{6Theorem}. 
Suppose $M(s_1, \dots, s_n)$ contains an embedded essential\index{essential} sphere, disk, annulus, or torus. Then $M$ contains an essential, homotopically $\bdy$-incompressible punctured sphere or torus $S$, with $\bdy S$ some subset of the slopes $s_1, \dots, s_n$. Moreover, if $S$ is a punctured sphere then it has at least three punctures. 
\end{lemma}

\begin{proof}
Let $T$ be an embedded essential sphere, disk, annulus, or torus in $M(s_1, \dots, s_n)$. Note that $M \subset M(s_1, \dots, s_n)$. If $T$ is embedded in $M$, then it cannot be essential in $M$ because $M$ is hyperbolic. But if $T$ is compressible\index{compressible} or boundary compressible\index{boundary compressible} in $M$, then a compression disk\index{compression disk} for $T$ is embedded in $M\subset M(s_1, \dots, s_n)$, hence is a compression disk for $T$ in $M(s_1, \dots, s_n)$, contradicting the fact that $T$ is essential. Similarly, if $T$ is boundary parallel\index{boundary parallel} in $M$ then it is compressible\index{compressible} in $M(s_1, \dots, s_n)$. So $T$ cannot be embedded in $M$; it must meet the solid tori attached to $M$ to form $M(s_1, \dots, s_n)$. We may assume it meets the cores of the added solid tori transversely, else we could isotope the surface to lie in $M$, which would be a contradiction. Thus when we drill these cores from $M(s_1, \dots, s_n)$, the surface $S$ given by removing neighborhoods of the cores from $T$ is a surface with boundary whose boundary components come from the set of slopes $\{s_1, \dots, s_n\}$. Note $S$ is a punctured sphere or punctured torus.

Now, $S$ cannot be compressible\index{compressible} in $M$, else a compression disk\index{compression disk} is a compression disk for $T$ in $M(s_1, \dots, s_n)$. Any boundary compression disk\index{boundary compression disk} $D$ must have boundary consisting of an arc $\alpha$ in $S$ and an arc $\beta$ running along a boundary component of $M$. Using $D$, we may isotope $T$ through $D$ to the core of a solid torus in $M(s_1, \dots, s_n)$, and then slightly past, removing two intersections of $T$ with cores of filled solid tori. So after repeating this move finitely many times, we may assume $S$ is boundary incompressible.\index{boundary incompressible} Note also that $S$ cannot be boundary parallel,\index{boundary parallel} or $T$ is compressible.\index{compressible} So $S$ is essential.

To show $S$ is homotopically $\bdy$-incompressible, apply a proof similar to that of \reflem{IncompressVsPi1} to show that if $S$ is not homotopically $\bdy$-incompressible, then the boundary of a regular neighborhood of $S$ is boundary compressible.\index{boundary compressible} The same argument as above implies that intersections of $S$ with cores of solid tori can be removed in this case.

Finally, if $S$ is a punctured sphere, then it must have at least three punctures else it will be an essential disk or annulus in $M$, but the hyperbolicity of $M$ rules out such surfaces. 
\end{proof}

The following lemma is from \cite{FuterSchleimer}, and uses arguments of \cite{agol:bounds}. 
\begin{lemma}\label{Lem:PleatedSurfaceLength}
Suppose $M$ is an orientable hyperbolic 3-manifold with a cusp, and horoball neighborhood $C$ about the cusp. Suppose $f\from S\to M$ is a pleating of a punctured surface $S$,\index{pleated surface} with $n$ punctures of $S$ mapping to $C$. Suppose finally that for each puncture of $S$,a loop about the puncture is represented by a geodesic of length $\lambda$ on $\bdy C$ in $M$. Then in the hyperbolic metric on $S$ given by the pleating, the preimage $f^{-1}(C)\subset S$ contains horospherical cusp neighborhoods $R_1, \dots, R_n$ of the $n$ punctures of $S$, with disjoint interiors, such that
  \[ \ell(\bdy R_i) = \area(R_i) \geq \lambda \quad \mbox{for each } i.\]
\end{lemma}

\begin{proof}
The pleating of $S$ gives $S$ an ideal triangulation. Start with a cusp neighborhood $C_0\subset C$ such that $f$ maps all ideal edges of the triangulation to geodesic rays running into the cusp. That is, $C_0$ does not intersect any edge of the triangulation in a compact arc. See \reffig{Pleating}, left. Then $f(S)\cap C_0$ consists of tips of triangles, and $f^{-1}(C_0)$ is a collection of embedded cusps $R_1^0, \dots, R_n^0$ in $S$.

\begin{figure}
%% Accompanies image file 'F8-07-Pleat.eps' (pdf, eps, ps)
\begingroup%
  \makeatletter%
  \providecommand\color[2][]{%
    \errmessage{(Inkscape) Color is used for the text in Inkscape, but the package 'color.sty' is not loaded}%
    \renewcommand\color[2][]{}%
  }%
  \providecommand\transparent[1]{%
    \errmessage{(Inkscape) Transparency is used (non-zero) for the text in Inkscape, but the package 'transparent.sty' is not loaded}%
    \renewcommand\transparent[1]{}%
  }%
  \providecommand\rotatebox[2]{#2}%
  \newcommand*\fsize{\dimexpr\f@size pt\relax}%
  \newcommand*\lineheight[1]{\fontsize{\fsize}{#1\fsize}\selectfont}%
  \ifx\svgwidth\undefined%
    \setlength{\unitlength}{259.62864876bp}%
    \ifx\svgscale\undefined%
      \relax%
    \else%
      \setlength{\unitlength}{\unitlength * \real{\svgscale}}%
    \fi%
  \else%
    \setlength{\unitlength}{\svgwidth}%
  \fi%
  \global\let\svgwidth\undefined%
  \global\let\svgscale\undefined%
  \makeatother%
  \begin{picture}(1,0.28486561)%
    \lineheight{1}%
    \setlength\tabcolsep{0pt}%
    \put(0,0){\includegraphics[width=\unitlength]{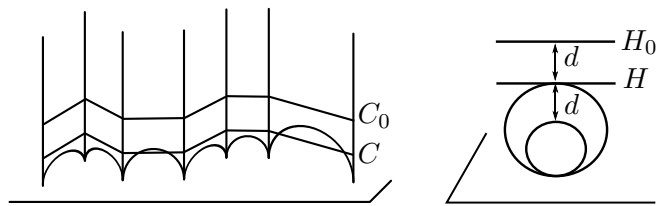}}%
    \put(0.50709233,0.11699405){\color[rgb]{0,0,0}\makebox(0,0)[lt]{\lineheight{1.25}\smash{\begin{tabular}[t]{l}$C_0$\end{tabular}}}}%
    \put(0.50709233,0.05921923){\color[rgb]{0,0,0}\makebox(0,0)[lt]{\lineheight{1.25}\smash{\begin{tabular}[t]{l}$C$\end{tabular}}}}%
    \put(0.80575553,0.19980658){\color[rgb]{0,0,0}\makebox(0,0)[lt]{\lineheight{1.25}\smash{\begin{tabular}[t]{l}$d$\end{tabular}}}}%
    \put(0.80648515,0.13019641){\color[rgb]{0,0,0}\makebox(0,0)[lt]{\lineheight{1.25}\smash{\begin{tabular}[t]{l}$d$\end{tabular}}}}%
    \put(0.88883355,0.22359709){\color[rgb]{0,0,0}\makebox(0,0)[lt]{\lineheight{1.25}\smash{\begin{tabular}[t]{l}$H_0$\end{tabular}}}}%
    \put(0.89138687,0.1688833){\color[rgb]{0,0,0}\makebox(0,0)[lt]{\lineheight{1.25}\smash{\begin{tabular}[t]{l}$H$\end{tabular}}}}%
  \end{picture}%
\endgroup%

  \caption{Left: Choose $C_0$ to meet only tips of triangles. Right: Distances between horoballs}
  \label{Fig:Pleating}
\end{figure}

Lift $M$ to the universal cover $\HH^3$. The cusps $C$ and $C_0$ both lift to collections of disjoint horoballs. Because $C_0$ is contained in $C$, for each horoball lift $H$ of $C$, there is a horoball lift $H_0$ of $C_0$ contained in $H$. Let $d$ denote the hyperbolic distance between $H_0$ and $H$. Note that since $C$ is embedded, the distance from $H_0$ to any other lift of $C_0$ must be at least $2d$. See \reffig{Pleating}, right. 
Projecting back to $M$, a geodesic between $H_0$ and any other lift of $C_0$ projects to a geodesic from $C_0$ to $C_0$ of length at least $2d$. Pulling back to $S$, the distance from $R_i^0$ to any other $R_j^0$ is at least $2d$. Let $R_1, \dots, R_n$ be cusps in $S$ of distance $d$ from $R_1^0, \dots, R_n^0$. They must be embedded.

We now show that the lengths $\ell(\bdy R_i)$ are at least $\lambda$ for all $i$. Let $\gamma_0$ be a Euclidean geodesic on $\bdy C_0$ representing $f(\bdy R_i^0)$. Since pleating may decrease distance, $\ell(\gamma_0)\leq \ell(\bdy R_i^0)$. Moreover, letting $\gamma$ be the loop on $\bdy C$ homotopic to $\gamma_0$, we have $\lambda = \ell(\gamma) = e^{-d}\ell(\gamma_0)$, because $\gamma$ and $\gamma_0$ lie on cusp boundaries of hyperbolic distance $d$ apart. Moreover, $\ell(\bdy R_i) = e^{-d}\ell(\bdy R_i^0)$. Putting this together,
\[ \lambda \leq \ell(\gamma) = e^{-d}\ell(\gamma_0) \leq e^{-d}\ell(\bdy R_i^0) = e^{-d}\cdot e^d\ell(\bdy R_i) = \ell(\bdy R_i). \]

Finally, we need to show that $f(R_i)$ is contained in $C$. We know $f(R_i^0)$ is contained in $C_0$. By construction, $f(R_i)$ is contained in a $d$-neighborhood of $C_0$. But a $d$-neighborhood of $C_0$ is the cusp $C$. Thus $f(R_i)$ lies in $C$.
\end{proof}

Now we present a result from \cite{boroczky}.

\begin{theorem}[B{\"o}r{\"o}czky cusp density theorem]\label{Thm:Boroczky}\index{B{\"o}r{\"o}czky cusp density theorem!2-dimensional}\index{cusp density theorem!2-dimensional}
  Let $S$ be a hyperbolic surface with cusps, and let $H$ be an embedded horoball neighborhood for the cusps of $S$. Then
  \[ \area(H) \leq \frac{3}{\pi}\area(S). \]
\end{theorem}

\begin{proof}
  Given $S$ and $H$, we claim there exists an ideal triangulation of $S$ such that for $T$ any triangle, $H$ meets $T$ only in connected neighborhoods of its ideal vertices in noncompact sets. For example this will hold for a subdivision of the canonical decomposition\index{canonical decomposition} of $S$ with respect to $H$, which is defined in \refchap{Canonical}. We will assume such a decomposition exists.

Map $T$ isometrically to the triangle $T' \subset \HH^2$ with ideal vertices at $0$, $1$, and $\infty$. The image of $H$ determines horoballs $H_0$, $H_1$, and $H_\infty$ about $0$, $1$, $\infty$, respectively. The area of $H\cap T$ in $S$ is given by the sum
\begin{equation}\label{Eqn:AreaSum}
  \area(H\cap T) = \area(H_0\cap T') + \area(H_1\cap T') + \area(H_\infty\cap T'),
\end{equation}
and the area of $H$ is given by the sum of all such areas over all triangles $T$. Since the area of $S$ is just $\pi$ times the number of triangles, to maximize the cusp density $\area(H)/\area(S)$ we need to maximize the cusp density within ideal triangles,\index{ideal triangle} or maximize the sum of \refeqn{AreaSum} within each triangle.

So consider the triangle $T'\subset \HH^2$ with vertices at $0$, $1$, and $\infty$, and with horoballs $H_\infty$ about $\infty$ of height $h_\infty$, $H_0$ about $0$ of diameter $h_0$, and $H_1$ about $1$ of diameter $h_1$. By the observation that (open) horoballs do not meet edges of the triangulation in intervals with compact closure, we know that $h_0$ and $h_1$ are at most $2$, and $h_\infty$ is at least $1/2$. These give constraints on $h_0$, $h_1$, and $h_\infty$. We also have constraints coming from the fact that $H_0$, $H_1$, and $H_\infty$ are disjoint.

If one of $H_0$, $H_1$, or $H_\infty$ is not tangent to one of the other two, then we may expand it, increasing cusp area, until either it is as large as possible given our constraints, or it is tangent to one of the other horoballs. If one of the $H_i$ is as large as possible, but still not tangent to the other horoballs, then the other two horoballs are much smaller than our constraints, and we may expand them until they are tangent to $H_i$. In any case, we may assume that the cusp area is maximized when each horoball is tangent to one other horoball; because there are three, it follows that one horoball, without loss of generality $H_\infty$, is tangent to the other two.

Now we may compute the cusp area directly. Since $H_0$ and $H_1$ are tangent to $H_\infty$, they have diameters $h_0=h_1=h_\infty$, and $H_\infty$ has height $H_\infty$. Then the areas of the cusps satisfy (\refex{CuspAreaTriangle}):
\[ \area(H_\infty\cap T') = \frac{1}{h_\infty}, \quad
\area(H_0\cap T') = \area(H_1\cap T') = h_\infty, \]
so
\[ \area(H\cap T') = \frac{1}{h_\infty} + 2 h_\infty. \]
We maximize this equation for $h_\infty$ subject to constraints: $h_\infty$ is at least $1/2$, and $H_1$ and $H_0$ are disjoint, so $h_\infty$ is at most $1$.
We find that the function has a critical point at $\sqrt{2}$, but that it reaches its maximum value when $h_\infty=1/2$ and when $h_\infty=1$, and the maximum is $3$.

Now let $n$ be the number of triangles in $S$. Then
\[
\frac{\area(H)}{\area(S)} \leq \frac{ 3\cdot n}{\pi \cdot n} = \frac{3}{\pi}.\qedhere
\]
\end{proof}

\begin{proof}[Proof of \refthm{6Theorem}, a 6-theorem]\index{6-theorem, weaker}
Suppose by way of contradiction that $M(s_1, \dots, s_n)$ is reducible,\index{reducible 3-manifold} boundary reducible, annular, or toroidal. Then \reflem{EssentialPunctSurface} implies $M$ contains an embedded essential\index{essential} punctured 2-sphere or torus, whose boundary components on $\bdy M$ are parallel to slopes $s_1, \dots, s_n$.

By \refprop{Pleating}, $S$ may be pleated.\index{pleated surface} By \reflem{PleatedSurfaceLength}, the pleating induces horoball neighborhoods $R_1, \dots, R_m$ of cusps of $S$ for which
\[ \ell(\bdy R_i) = \area(R_i) \geq \ell(s_{j}),\]
where $f(\bdy R_i)$ is the slope $s_{j}$. Let $H$ denote the union of the horoball neighborhoods $R_i$. 

Now \refthm{Boroczky} and the Gauss--Bonnet theorem imply
\begin{equation}\label{Eqn:6Thm}
\sum_i \ell(s_{j_i}) \leq \sum_i \ell(\bdy R_i) = \area(H) \leq \frac{3}{\pi}\area(S) = \frac{3}{\pi}\cdot 2\pi|\chi(S)| = 6|\chi(S)|. 
\end{equation}

On the other hand, each $\ell(s_{j_i})>6$, and there are $m$ of these, where $m$ is the number of boundary components of $S$. If $S$ is a punctured sphere, $|\chi(S)| = m-2$ and \refeqn{6Thm} implies $6m < 6(m-2)$, which is a contradiction. If $S$ is a punctured torus, \refeqn{6Thm} implies $6m < 6m$, again a contradiction. 
\end{proof}

Our 6-theorem\index{6-theorem, weaker} has immediate consequences to determining when certain knots and links are hyperbolic.

\begin{definition}\label{Def:HighlyTwisted}
  For an integer $c>0$, we say a knot or link is \emph{$c$-highly twisted} if it admits a diagram in which every twist region has at least $c$ crossings. If $c$ is understood from the context, we also say such a knot or link is \emph{highly twisted}.\index{highly twisted}
\end{definition}

The following theorem was first proved in \cite{futer-purcell}. 

\begin{theorem}[Hyperbolicity of highly twisted links]\label{Thm:FuterPurcellFilling}\index{highly twisted!hyperbolicity}
Let $K\subset S^3$ be a link with a prime, twist-reduced diagram, as in \refdef{TwReduced}. Assume that $K$ has at least two twist regions. If every twist region of the diagram contains at least six crossings, then the complement of $K$ is hyperbolic. 
\end{theorem}

\begin{proof}
By \reflem{TwReducedGivesReducedAug}, when we add crossing circles to every twist region of $K$, we obtain a fully augmented link\index{fully augmented link} $L$; $K$ now forms the knot strands of this link. Remove all crossings of the knot strands except possibly single crossings in twist regions. Then $K$ is obtained from $L$ by performing Dehn filling along slopes $s_j$ on crossing circles $C_j$ that replace the crossing circle with a twist region with at least six crossings.

By \refthm{AugSlopeLengths}, each slope $s_j$ has length at least $\sqrt{(6)^2+1} =\sqrt{37}>6$. Thus by our 6-Theorem, \refthm{6Theorem},\index{6-theorem, weaker} the link complement $L$ obtained by the Dehn filling is irreducible,\index{irreducible} boundary irreducible,\index{boundary irreducible} anannular,\index{anannular} and atoroidal.\index{atoroidal} By Thurston's hyperbolization theorem, \refthm{SfcesHyperbolic}, it is hyperbolic. 
\end{proof}

The full 6-theorem can be used to identify knots in $S^3$.

\begin{corollary}\label{Cor:KnotsInS3}
  Suppose $M$ is a hyperbolic 3-manifold with a single cusp. Then $M$ is the complement of a knot in $S^3$ if and only if there exists a slope $s$ on the cusp of $M$ of length at most six such that $M(s)$ is homeomorphic to $S^3$.
\end{corollary}

\begin{proof}
  A manifold with a torus boundary component is the complement of a knot in $S^3$ if and only if a Dehn filling gives $S^3$. The full 6-theorem\index{6-theorem} of Agol and Lackenby implies such a slope on a horocusp of a hyperbolic 3-manifold must have length at most six. 
\end{proof}

Since only finitely many slopes on a fixed hyperbolic 3-manifold have length at most six, to determine whether a hyperbolic 3-manifold is a knot in $S^3$, it suffices to check finitely many Dehn fillings,\index{Dehn filling} and then to identify whether or not the filled manifold is $S^3$. 

%%%%%%%%%%%%%%%%%%%%%%%%%%%%%%%%%%%%%%%%%%%%%%%%%%%%%%%%%%%%%%%%%
\section{Exercises}

\begin{exercise}
  Prove that the unknot is the only knot $K$ in $S^3$ such that $S^3-N(K)$ admits a properly embedded compression disk\index{compression disk} for $\bdy N(K)$. 
\end{exercise}

\begin{exercise}
Prove that the paragraph starting with ``Equivalently'' in \refdef{Satellite} is indeed an equivalent definition of a satellite knot. You may assume Alexander's theorem from 3-manifold topology that states that an embedded torus in $S^3$ bounds a solid torus on at least one side.
\end{exercise}

\begin{exercise}\label{Ex:ProveIncompressVsPi1}
  Prove \reflem{IncompressVsPi1}. For one direction, you may use the \emph{loop theorem}, which is a classical result in 3-manifold topology:

  \begin{theorem}[Loop theorem \cite{papa}]\label{Thm:LoopTheorem}\index{loop theorem}
    If $N$ is a 3-manifold with boundary, and there is a map $f\from D^2\to N$ such that the loop $f(\bdy D^2)\subset \bdy N$ is homotopically nontrivial in $\bdy N$, then there is an embedding with the same property.
  \end{theorem}
\end{exercise}

\begin{exercise}\label{Ex:NonorIncompressVsPi1}
  Prove that a nonorientable surface $S$ properly embedded in a 3-manifold $M$ is $\pi_1$-injective\index{$\pi_1$-injective} if and only if $\widetilde{S}$, the boundary of a regular neighborhood of $S$ in $M$, is an orientable incompressible\index{incompressible} surface.
\end{exercise}
    
\begin{exercise}\label{Ex:TorusKnot}
  \begin{enumerate}
  \item Prove that a $(p,q)$-torus knot $T(p,q)$ is nontrivial if $|p|, |q|\geq 2$.
  \item Let $T$ denote the torus in $S^3$ on which the torus knot $T(p,q)$ lies, and let $A$ denote the annulus $T-T(p,q)$. Prove $A$ is incompressible\index{incompressible} if $|p|, |q|\geq 2$.
  \end{enumerate}
  Hint for both parts: Seifert--Van Kampen theorem.
\end{exercise}

\begin{exercise}\label{Ex:EssentialSfcesMeetPolyinDisks}
Suppose $M$ is a 3-manifold with an ideal polyhedral decomposition and $S$ a properly embedded essential\index{essential} surface in $M$ with complexity as in the proof of \refthm{NormalForm}. 
\begin{enumerate}
\item Prove that if $S$ is an embedded essential sphere, then we may replace $S$ with an embedded essential sphere $S'$ meeting each polyhedron in disks such that the complexity of $S'$ is at most that of $S$.
\item Prove that if $M$ is irreducible\index{irreducible} and $S$ is an essential surface properly embedded in $M$, then $S$ can be isotoped to meet polyhedra only in disks, reducing (or at worst fixing) complexity. 
\end{enumerate}
\end{exercise}

\begin{exercise}\label{Ex:NormalDisk}
  Suppose $M$ is a 3-manifold with an ideal polyhedral decomposition, and $S$ is an essential\index{essential} surface properly embedded in $M$. Suppose there is a disk $D$ of intersection of $S$ with a polyhedron $P$ such that $\bdy D$ meets an edge of $P$ more than once. Prove $S$ can be isotoped to remove at least two intersections with that edge. 
\end{exercise}

\begin{exercise}\label{Ex:NormalBdry}
  Suppose $M$ is a 3-manifold with an ideal polyhedral decomposition, and $S$ is an essential\index{essential} surface with boundary properly embedded in $M$. Suppose there is a disk $D$ of intersection of $S$ with a polyhedron such that $\bdy D$ meets a boundary face more than once. Then if $S$ is a disk, prove it can be replaced by an essential disk meeting the boundary face fewer times. If $S$ is not a disk and $M$ is irreducible\index{irreducible} and boundary irreducible,\index{boundary irreducible} prove $S$ can be isotoped to meet the boundary face fewer times.
\end{exercise}

\begin{exercise}\label{Ex:AreaPolygon}
  Prove that if $D$ is a hyperbolic polygon, then its hyperbolic area is
  \[ \area(D) = \sum(\pi-\alpha_i) - 2\pi + \pi \, v,\]
  where $\alpha_i$ is the angle of the $i$-th finite vertex, and $v$ is the number of ideal vertices of $D$.
\end{exercise}

\begin{exercise}\label{Ex:CuspAreaTriangle}
  Consider the ideal triangle\index{ideal triangle} $\Delta$ with vertices at $0$, $1$, and $\infty$, and horoballs $H_\infty$ about $\infty$ of height $h_\infty$, $H_0$ about $0$ of diameter $h_0$, and $H_1$ about $1$ of diameter $h_1$. Prove that the areas of $H_i\cap \Delta$ satisfy
  \[ \area(H_\infty\cap \Delta) = \frac{1}{h_\infty}, \quad \area(H_0\cap\Delta)=h_0, \quad \mbox{and } \area(H_1\cap\Delta)=h_1, \]
  so the total cusp area of $\Delta$ is $1/h_\infty + h_0 + h_1$.
\end{exercise}

\begin{exercise}
  Let $S$ be a hyperbolic surface with a cusp, with horoball cusp neighborhood $C$. Show that the length of the boundary of $C$ is equal to the area of $C$.
\end{exercise}

\begin{exercise}
  Let $T'$ be the ideal triangle\index{ideal triangle} in $\HH^2$ with vertices at $0$, $1$, and $\infty$, with horoballs $H_\infty$ about $\infty$ of height $h_\infty$, $H_0$ about $0$ of diameter $h_0$, and $H_1$ about $1$ of diameter $h_1$. Show that
  \[ \area(H_\infty\cap T') = \frac{1}{h_\infty}, \quad \area(H_0\cap T') = h_0, \quad \area(H_1\cap T') = h_1.\]
\end{exercise}

%% Ch09_AngleStruct.tex

\chapter{Volume and Angle Structures}\label{Chap:AngleStruct}
\blfootnote{Jessica S. Purcell, Hyperbolic Knot Theory}

Those hyperbolic 3-manifolds that admit a triangulation by positively oriented geometric tetrahedra\index{positively oriented tetrahedron}\index{tetrahedron!positively oriented} exhibit many additional nice properties. The existence of such a triangulation often gives a simpler way to prove many results in hyperbolic geometry. We present some of the techniques and consequences in this chapter.

In the theory of knots and links, these tools have been applied to great effect to an infinite class of knots and links called 2-bridge links,\index{2-bridge knot or link} which we will describe (and triangulate) in the next chapter.

%%%%%%%%%%%%%%%%%%%%%%%%%%%%%%%%%%%%%%%%%%%%%%%%%%%%%%%%%%%%%%%%%
\section{Hyperbolic volume of ideal tetrahedra}

Ideal tetrahedra are building blocks of many complete hyperbolic manifolds. In this section, we will calculate volumes of ideal tetrahedra. 

Recall that a hyperbolic ideal tetrahedron is completely determined by $z\in\CC$ with positive imaginary part, as in \refdef{EdgeInvariant}. It is also determined by three dihedral angles, as the following lemma shows. 

\begin{lemma}\label{Lem:AnglesDetermineTet}
Let $\alpha$, $\beta$, $\gamma$ be angles in $(0,\pi)$ such that $\alpha+\beta+\gamma=\pi$. Then $\alpha$, $\beta$, and $\gamma$ determine a unique hyperbolic ideal tetrahedron up to isometry of $\HH^3$. Conversely, any hyperbolic tetrahedron determines unique $\{\alpha, \beta, \gamma\} \subset (0,\pi)$ with $\alpha+\beta+\gamma=\pi$. 
\end{lemma}

\begin{proof}
First we prove the converse. Given an ideal tetrahedron with ideal vertices on $\bdy \HH^3$ at $0$, $1$, $\infty$, and $z$, note that a horosphere about $\infty$ intersects the tetrahedron in a Euclidean triangle. Let $\alpha$, $\beta$, $\gamma$ denote the interior angles of the triangle; these are dihedral angles of the tetrahedron. Each angle $\alpha$, $\beta$, $\gamma$ lies in $(0,\pi)$, and the sum $\alpha+\beta+\gamma=\pi$, as desired. \Refex{TetLabels1} shows that taking a different collection of vertices to $0$, $1$, and $\infty$ will give the same dihedral angles $\alpha$, $\beta$, $\gamma$, so these three angles are uniquely determined by the tetrahedron. 

Now, suppose $\alpha$, $\beta$, and $\gamma$ in $(0,\pi)$ are given, with $\alpha+\beta+\gamma=\pi$. Then these three numbers determine a Euclidean triangle, uniquely up to scale, with interior angles $\alpha$, $\beta$, $\gamma$. View the triangle as lying in $\CC$; we may adjust such a triangle so that it has vertices at $0$, $1$, and some $z\in\CC$ with positive imaginary part. This determines a tetrahedron with edge parameter $z$. If we rotate and scale the triangle so that different vertices map to $0$ and $1$, this corresponds to mapping different ideal vertices of the tetrahedron to $0$ and $1$. The parameter $z$ will be adjusted as in \reflem{EdgeInvariants}, but the tetrahedron will be the same up to isometry. 
\end{proof}

Lemmas~\ref{Lem:AnglesDetermineTet} and~\ref{Lem:EdgeInvariants} give two different ways of uniquely describing an ideal tetrahedron, either by a single complex number $z$ or by a triple of angles $\alpha$, $\beta$, $\gamma$ with $\alpha+\beta+\gamma=\pi$. We will compute volumes of an ideal tetrahedron, and we choose to compute volumes using a parameterization by angles rather than edge parameter, although computations can be done either way. (See exercises.)

\begin{definition}\label{Def:LobachevskyFunction}
The \emph{Lobachevsky function}\index{Lobachevsky function} $\Lambda(\theta)$ is the function defined by
\[
\Lambda(\theta) = - \int_0^\theta \log|2\sin u|\, du.
\]
\end{definition}

\begin{theorem}\label{Thm:VolTet}
Suppose $\alpha$, $\beta$, and $\gamma$ are angle measures strictly between $0$ and $\pi$, and suppose $\alpha+\beta+\gamma=\pi$, so they determine a hyperbolic ideal tetrahedron $\Delta(\alpha, \beta, \gamma)$. Then the volume $\vol(\Delta(\alpha, \beta, \gamma))$ is equal to
\[
\vol(\Delta(\alpha, \beta, \gamma)) = \Lambda(\alpha)+\Lambda(\beta)+\Lambda(\gamma),
\]
where $\Lambda$ is the Lobachevsky function of \refdef{LobachevskyFunction}. 
\end{theorem}

\begin{example}\label{Example:VolFig8}
The figure-8 knot complement has complete hyperbolic structure built of two regular ideal tetrahedra. Therefore the volume of the figure-8 knot complement is $6\Lambda(\pi/3)$, which can be numerically calculated to be approximately $2.0299$.
\end{example}

Our proof of \refthm{VolTet} follows that given by Milnor in \cite{Milnor:HypGeom} and also in \cite[Chapter~7]{thurston}. Milnor, in turn, credits Lobachevsky for several of his calculations.

First, we need a lemma concerning the Lobachevsky function.

\begin{lemma}\label{Lem:Lobachevsky}
  The Lobachevsky function $\Lambda(u)$ satisfies:
  \begin{enumerate}
  \item It is well-defined and continuous on $\RR$ (even though the defining integral is improper).
  \item $\Lambda(-\theta)=-\Lambda(\theta)$, i.e.\ $\Lambda(\theta)$ is odd.
  \item $\Lambda(\theta)$ is periodic of period $\pi$. 
  \item\label{Itm:Kubert2} It satisfies the expression $\Lambda(2\theta) = 2\Lambda(\theta) + 2\Lambda(\theta+\pi/2).$
  \end{enumerate}
\end{lemma}

\begin{proof}
To prove the lemma, we will relate the Lobachevsky function to the well-known  \emph{dilogarithm function}\index{dilogarithm function}
\begin{equation}\label{Eqn:Dilogarithm}
\psi(z) = \sum_{n=1}^\infty z^n/n^2 \quad \mbox{for } |z|\leq 1.
\end{equation}
For more information on the dilogarithm, see for example \cite{Zagier:Dilogarithm}. Note that for $|z|<1$, the derivative of $\psi(z)$ satisfies
\[
\psi'(z) = \sum_{n=1}^\infty \frac{z^{n-1}}{n} = \frac{1}{z}\left( \sum_{n=1}^\infty \frac{z^n}{n} \right).
\]
The sum on the right hand side is a well-known Taylor series:
\[ -\log(1-z) = \sum_{n=1}^\infty \frac{z^n}{n} \quad \mbox{for } |z|<1. \]
Thus the analytic continuation\index{analytic continuation} of $\psi(z)$ is given by
\begin{equation}\label{Eqn:DilogIntegral}
\psi(z) = -\int_0^z \frac{\log(1-u)}{u}\,du \quad \mbox{ for } z\in \CC-[1,\infty).
\end{equation}

For $0 < u < \pi$, consider $\psi(e^{2iu}) - \psi(1)$. Although the integral formula \refeqn{DilogIntegral} above is not defined at $z=1$, the summation of \refeqn{Dilogarithm} is defined and continuous at $z=1$ (in fact, $\psi(1)=\pi^2/6$), so we may write
\[
\psi(e^{2iu})-\psi(1) = -\int_1^{\displaystyle{e^{2iu}}} \frac{\log(1-w)}{w}\,dw.
\]
Substitute $w=e^{2i\theta}$ into this expression to obtain
\begin{align*}
  \psi(e^{2iu})-\psi(1) &= - \int_{\theta=0}^u \log(1-e^{2i\theta})\,(2i)\,d\theta \\
  &= - \int_0^u\log\left(-2ie^{i\theta}\left(\frac{e^{i\theta}-e^{-i\theta}}{2i}\right) \right) (2i)\,d\theta \\
  &= - \int_0^u 2i(\log(-i)+\log(e^{i\theta}) + \log(2\sin\theta))\,d\theta \\
  &= - \int_0^u (\pi - 2\theta + 2i\log(2\sin\theta))\,d\theta.\\
\end{align*}
Take the imaginary parts of both sides of the above equation. Note $\psi(1)$ is real, hence
\[
\Im(\psi(e^{2iu})-\psi(1))=\Im(\psi(e^{2iu})) = \Im\left(\sum_{n=1}^\infty \frac{e^{2inu}}{n^2}\right) = \sum_{n=1}^\infty \frac{\sin(2nu)}{n^2}.\]
On the other side, this equals
\[ \Im(\psi(e^{2iu})-\psi(1)) = 2\int_0^u-\log(2\sin\theta)\,d\theta = 2\Lambda(u).\]
Thus for $0\leq u \leq \pi$, we have the uniformly convergent Fourier series for $\Lambda(u)$ given by
\begin{equation}\label{Eqn:LobachevskiDilog}
  \Lambda(u) = \half\sum_{n=1}^\infty \frac{\sin(2nu)}{n^2} \quad \mbox{for } 0\leq u\leq\pi.
\end{equation}
This shows $\Lambda(u)$ is well-defined and continuous for $0\leq u\leq \pi$. It also shows that $\Lambda(u)$ can be defined on $-\pi\leq u\leq 0$, and it is an odd function on this range. Finally, it shows that $\Lambda(0)=\Lambda(\pi) =0$. 

Notice now that the derivative $d\Lambda(\theta)/d\theta = -2\log|2\sin\theta|$ is periodic of period $\pi$. Then for $\theta>\pi$,
\begin{align*}
  \Lambda(\theta) & = \int_0^\theta \Lambda'(u)\,du = \int_0^\pi\Lambda'(u)\,du + \int_\pi^\theta \Lambda'(u)\,du \\
  & = \Lambda(\pi) + \int_0^{\theta-\pi}\Lambda'(u)\,du = \Lambda(\theta-\pi),
\end{align*}
by the periodicity of $\Lambda'$, and the fact that $\Lambda(\pi)=0$. This shows that $\Lambda$ is well-defined and continuous for $\theta\geq 0$; a similar result implies it is well-defined and continuous for $\theta\leq 0$, and it will be odd everywhere. 

It only remains to show the last item of the lemma. To do so, begin with the identity
\[
2\sin(2\theta) = 4\sin\theta\cos\theta = (2\sin\theta)(2\sin(\theta+\pi/2)).
\]
Then note that
\begin{align*}
  \Lambda(2\theta) &= \int_0^{2\theta} - \log|2\sin u| \, du \\
  &= 2\int_0^\theta - \log|2\sin(2w)|\,dw \quad (\mbox{letting } w=u/2) \\
  &= 2\int_0^\theta -\log|2\sin w|\, dw + 2\int_0^\theta-\log|2\sin(w+\pi/2)|\,dw \\
  &= 2\Lambda(\theta) + 2\int_{\pi/2}^{\theta+\pi/2} -\log|2\sin v|\, dv \\
  &= 2\Lambda(\theta) +2\Lambda(\theta+\pi/2) - 2\Lambda(\pi/2).
\end{align*}
Finally, note that if we substitute $u=\pi/2$ into \refeqn{LobachevskiDilog}, we obtain $\Lambda(\pi/2)=0$. This finishes the proof of the lemma. 
\end{proof}

\begin{remark}\label{Rem:Kubert}
  Item~\eqref{Itm:Kubert2} of \reflem{Lobachevsky} is a special case of more general identities known as the \emph{Kubert identities}\index{Kubert identities}, which have the following form. For any nonzero integer $n$, 
  \[ \Lambda(n\theta) = \sum_{k=0}^{n-1} n\Lambda(\theta + k\pi/n). \]
  You are asked to prove these identities in the exercises. 
\end{remark}

To prove \refthm{VolTet}, we will subdivide our ideal tetrahedron into six 3-dimensional simplices, each simplex with some finite and some infinite vertices. Such a simplex will be described by a region in $\HH^3$. To obtain the volume, we integrate the hyperbolic volume form $d\vol = dx\,dy\,dz/z^3$ over the region describing the simplex, and then sum the six results.

More carefully, given an ideal tetrahedron in $\HH^3$, we have been viewing the tetrahedron as having vertices $0$, $1$, $\infty$, and $z$. The three points $0$, $1$, and $z$ determine a Euclidean circle on $\CC$, which is the boundary of a Euclidean hemisphere, giving a hyperbolic plane in $\HH^3$. To this picture, apply a hyperbolic isometry that takes the circle on $\CC$ through $0$, $1$, $z$ to the unit circle in $\CC$, taking $0$, $1$, $z$ to some points $p$, $q$, $r$ on $S^1\subset \CC$. 

\begin{figure}
  \includegraphics{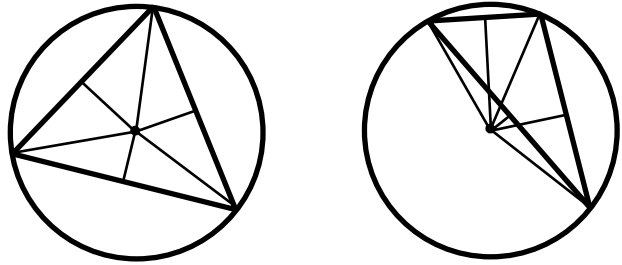}
  \caption{Left is a tetrahedron for which the point $(0,0)$ lies in the interior of the triangle on $\CC$, right is one for which it is exterior. Both show subdivisions into six triangles. }
  \label{Fig:TrianglesVolume}
\end{figure}

Now, drop a perpendicular from $\infty$ to the hemisphere; this will be a vertical ray from $(0,0,1)\in \HH^3$ to $\infty$. There will be two cases to consider: the case that the point $(0,0)\in \CC$ is interior to the triangle determined by $p$, $q$, $r$, and the case that the point $(0,0)$ is exterior to that triangle. The cases are shown in \reffig{TrianglesVolume}.

Consider first the case that the point $(0,0)$ is interior to the triangle determined by $p$, $q$, and $r$. Then the ray from $(0,0,1)$ to $\infty$ lies interior to the tetrahedron. Now, on the hemisphere whose boundary is the unit circle, draw perpendicular arcs from $(0,0,1)$ to each edge of the tetrahedron lying on that hemisphere. Also draw arcs from $(0,0,1)$ to the vertices of the tetrahedron, as shown in \reffig{TrianglesVolume}. Now cone to $\infty$. This divides the original tetrahedron up into six simplices. Similarly, if $(0,0)$ is not interior to the triangle determined by $p$, $q$, and $r$, it still makes sense to draw the same arcs and rays, as in \reffig{TrianglesVolume}, right. However, in this case the six simplices obtained overlap each other. In either case, we have the following result.

\begin{figure}
%% Accompanies image file 'F9-02-Simplx.eps' (pdf, eps, ps)
\begingroup%
  \makeatletter%
  \providecommand\color[2][]{%
    \errmessage{(Inkscape) Color is used for the text in Inkscape, but the package 'color.sty' is not loaded}%
    \renewcommand\color[2][]{}%
  }%
  \providecommand\transparent[1]{%
    \errmessage{(Inkscape) Transparency is used (non-zero) for the text in Inkscape, but the package 'transparent.sty' is not loaded}%
    \renewcommand\transparent[1]{}%
  }%
  \providecommand\rotatebox[2]{#2}%
  \newcommand*\fsize{\dimexpr\f@size pt\relax}%
  \newcommand*\lineheight[1]{\fontsize{\fsize}{#1\fsize}\selectfont}%
  \ifx\svgwidth\undefined%
    \setlength{\unitlength}{106.70078945bp}%
    \ifx\svgscale\undefined%
      \relax%
    \else%
      \setlength{\unitlength}{\unitlength * \real{\svgscale}}%
    \fi%
  \else%
    \setlength{\unitlength}{\svgwidth}%
  \fi%
  \global\let\svgwidth\undefined%
  \global\let\svgscale\undefined%
  \makeatother%
  \begin{picture}(1,1.26709001)%
    \lineheight{1}%
    \setlength\tabcolsep{0pt}%
    \put(0,0){\includegraphics[width=\unitlength]{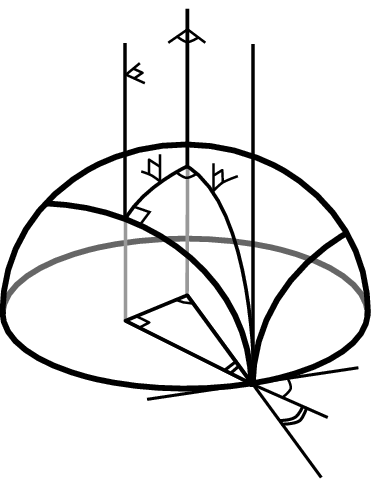}}%
    \put(0.83890048,0.06910366){\color[rgb]{0,0,0}\makebox(0,0)[lt]{\lineheight{0}\smash{\begin{tabular}[t]{l}$\frac{\pi}{2}-\zeta$\end{tabular}}}}%
    \put(0.52104817,1.07418171){\color[rgb]{0,0,0}\makebox(0,0)[lt]{\lineheight{0}\smash{\begin{tabular}[t]{l}$\zeta$\end{tabular}}}}%
    \put(0.88867339,0.20479862){\color[rgb]{0,0,0}\makebox(0,0)[lt]{\lineheight{0}\smash{\begin{tabular}[t]{l}$\zeta$\end{tabular}}}}%
  \end{picture}%
\endgroup%

  \caption{One of the six simplices obtained from subdividing an ideal tetrahedron}
  \label{Fig:SixSimplices}
\end{figure}

\begin{lemma}\label{Lem:SixSimplices}
Each of the six simplices obtained as above has the following properties, illustrated in \reffig{SixSimplices}.
\begin{enumerate}
\item It has two finite vertices and two ideal vertices.
\item Three of its dihedral angles are $\pi/2$, the other dihedral angles are $\zeta$, $\zeta$, and $\pi/2-\zeta$ for some $\zeta\in(0,\pi/2)$.
\end{enumerate}
\end{lemma}

\begin{proof}
Note that by construction, the two ideal vertices are at $\infty$ and one of $p$, $q$, $r$, i.e.\ one of the vertices of the original ideal tetrahedron. The other vertices are at $(0,0,1)$, and some point on the unit hemisphere where an arc from $(0,0,1)$ meets an edge of the original tetrahedron in a right angle.
  
Consider the dihedral angles of the faces meeting infinity. Each of these is a cone (to $\infty$) over an edge on the unit hemisphere. The dihedral angles agree with the dihedral angles of the vertical projection of the simplex to $\CC$, which is the triangle $T$ shown in \reffig{TrianglesVolume}; these angles are $\pi/2$, $\zeta$, and $\pi/2-\zeta$ for some $\zeta\in (0,\pi/2)$. The fourth face of the tetrahedron lies on the hemisphere. It meets both vertical faces through $(0,0,1)$ in right angles. The final face is a subset of a vertical plane whose boundary on $\CC$ is a line $L$ containing a side of the projection triangle $T$. The angle this vertical plane meets with the unit hemisphere is obtained by measuring the angles between the line $L$ and a tangent to the unit circle at the points where these intersect. Notice this angle is complementary to $\pi/2-\zeta$, hence is $\zeta$. 
\end{proof}

A simplex with the form of \reflem{SixSimplices} is called an \emph{orthoscheme},\index{orthoscheme} named by Scl\"afli in the 1950s \cite{SchlafliI, SchlafliII}. Around that time, he computed volumes of orthoschemes.

\begin{lemma}\label{Lem:VolSubSimplex}
Let $S(\zeta)$ denote a simplex obtained as above, with properties of \reflem{SixSimplices}. That is, $S(\zeta)$ has two finite vertices and two ideal vertices, three dihedral angles of $\pi/2$, and other dihedral angles $\zeta$, $\zeta$, and $\pi/2-\zeta$ for $\zeta\in(0, \pi/2)$. Then the volume of $S(\zeta)$ is
\[ \vol(S(\zeta)) = \frac{1}{2}\Lambda(\zeta). \]
\end{lemma}

\begin{proof}
The proof is a computation.

Apply an isometry to $\HH^3$ so that one ideal vertex of $S(\zeta)$ lies at $\infty$, the other on the unit circle, with one of the finite vertices at $(0,0,1)$; this is the same position of the simplex in the proof of \reflem{SixSimplices} above. When we project vertically to $\CC$, we obtain a triangle $T$ with one vertex at $0$, one on the unit circle, and the last some $v\in\CC$. The angle at $v$ is $\pi/2$, and the other two angles are $\zeta$ and $\pi/2-\zeta$. By applying a M\"obius transformation\index{M\"obius transformation} that rotates and reflects (but does not affect volume), we may assume $v$ is the point $\cos(\zeta) \in \RR\subset\CC$, and the third point, on the unit circle, is the point $\cos(\zeta) + i\, \sin(\zeta)$.

Now the triangle $T$ is described by the region
\[ 0\leq x \leq \cos(\zeta) \quad \mbox{and} \quad 0\leq y\leq x\tan(\zeta).\]
Then $\vol(S(\zeta))$ is given by
\[  \vol(S(\zeta)) = \int_T\int_{z\geq\sqrt{1-x^2-y^2}} d\vol
= \int_{0}^{\cos(\zeta)}\int_0^{x\tan(\zeta)}\int_{\sqrt{1-x^2-y^2}}^{\infty} \frac{dz\,dy\,dx}{z^3} \\
\]
Integrating with respect to $z$, we obtain
\[
\vol(S(\zeta))  = \int_0^{\cos(\zeta)}\int_0^{x\tan(\zeta)}\frac{dx\,dy}{2(1-x^2-y^2)},
\]
which we rewrite 
\[
\vol(S(\zeta))  = \int_0^{\cos(\zeta)} \int_0^{x\tan(\zeta)} \frac{dx\,dy}{2((\sqrt{1-x^2})^2-y^2)}, 
\]
and integrate with respect to $y$:
\begin{align*}
  \vol(S(\zeta)) 
  & = \int_0^{\cos(\zeta)} \frac{1}{4\sqrt{1-x^2}} \log \left(\frac{\sqrt{1-x^2}+x\tan\zeta}{\sqrt{1-x^2}-x\tan\zeta}\right)\,dx \\
  & = \int_0^{\cos(\zeta)} \frac{1}{4\sqrt{1-x^2}}\log \left(\frac{\sqrt{1-x^2}\cos(\zeta) + x\sin(\zeta)}{\sqrt{1-x^2}\cos(\zeta)-x\sin(\zeta)}\right) \, dx.
\end{align*}
Using the substitution $x=\cos(\theta)$, the integral becomes
\begin{align*}
  \vol(S(\zeta)) &= \int_{\pi/2}^{\zeta} \frac{1}{4}\log\left(\frac{\sin\theta\cos\zeta+\cos\theta\sin\zeta}{\sin\theta\cos\zeta-\cos\theta\sin\zeta}\right)\,(-d\theta) \\
  &= -\frac{1}{4}\left( \int_{\pi/2}^{\zeta}\log\left(\frac{2\sin(\theta+\zeta)}{2\sin(\theta-\zeta)}\right)\,d\theta \right)\\
  &= \frac{1}{4}\left( \int_{\pi/2}^{\zeta} -\log(2\sin(\theta+\zeta))\,d\theta - \int_{\pi/2}^{\zeta} -\log(2\sin(\theta-\zeta))\,d\theta \right) \\
  &= \frac{1}{4}\left( \int_{\pi/2+\zeta}^{2\zeta} -\log(2\sin(u))\,du - \int_{\pi/2-\zeta}^0-\log(2\sin(u))\,du \right) \\
  &= \frac{1}{4}(\Lambda(2\zeta)-\Lambda(\pi/2+\zeta)+\Lambda(\pi/2-\zeta)).
\end{align*}

To finish, we use \reflem{Lobachevsky}. Since $\Lambda(\theta)$ is periodic of period $\pi$, note that $\Lambda(\pi/2-\zeta) = \Lambda(-\pi/2-\zeta)$. Since $\Lambda$ is an odd function, $\Lambda(-\pi/2-\zeta)=-\Lambda(\pi/2+\zeta)$. Finally, since $\Lambda(2\zeta)= 2\Lambda(\zeta)+2\Lambda(\zeta+\pi/2)$, the above becomes
\[ \vol(S(\zeta)) = \frac{1}{4}( 2\Lambda(\zeta) +2\Lambda(\zeta+\pi/2) - 2\Lambda(\pi/2+\zeta)) = \frac{1}{2}\Lambda(\zeta).\qedhere
\]
\end{proof}

\begin{proof}[Proof of \refthm{VolTet}]
For an ideal tetrahedron with dihedral angles $\alpha$, $\beta$, and $\gamma$, place the tetrahedron in $\HH^3$ with vertices at $\infty$, and at $p$, $q$, $r$ all on the unit circle in $\CC$. As above, drop a perpendicular ray to the unit hemisphere. 

\underline{Case 1.} Suppose first that the ray lies in the interior of the ideal tetrahedron. Then subdivide the tetrahedron into six simplices as before. Each of the simplices has the properties of \reflem{SixSimplices}, and is determined by some $\zeta\in(0,\pi/2)$. By \reflem{VolSubSimplex}, its volume is determined by $\zeta$ as well, so it remains to calculate $\zeta$ for each of the six simplices making up the ideal tetrahedron. Project vertically to the complex plane $\CC$; the angles determining the simplex can then be easily computed using Euclidean geometry. In particular, there are two with angle $\alpha$, two with angle $\beta$, and two with angle $\gamma$. See the left of \reffig{SubSimplices}.

\begin{figure}
  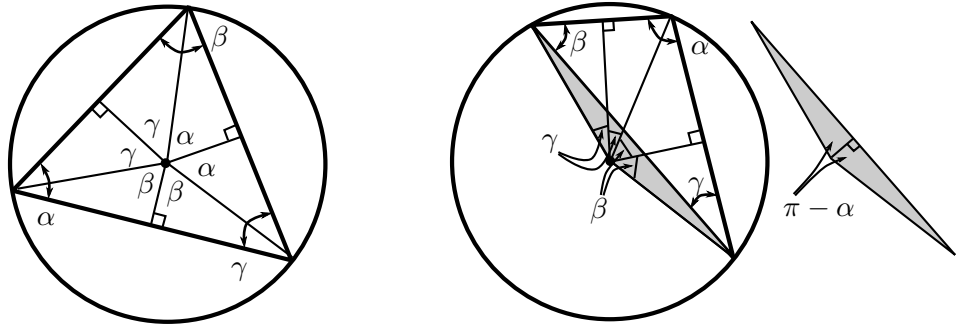
  \caption{Left: angles of subsimplicies when perpendicular ray lies interior to the tetrahedron. Right: angles when it is exterior}
  \label{Fig:SubSimplices}
\end{figure}

Then the volume of the tetrahedron $\Delta(\alpha, \beta,\gamma)$ is
\begin{align*}
  \vol(\Delta(\alpha, \beta, \gamma)) &= 2\vol(S(\alpha)) + 2\vol(S(\beta)) + 2\vol(S(\gamma)) \\
  &= \Lambda(\alpha)+\Lambda(\beta)+\Lambda(\gamma)
\end{align*}

\underline{Case 2.} Now suppose that the ray from $\infty$ to the point $(0,0,1)$ lies outside of the ideal tetrahedron. We may still draw perpendicular lines from $(0,0,1)$ to the edges of the ideal tetrahedron on the unit hemisphere, and lines from $(0,0,1)$ to vertices of the ideal tetrahedron; the right of \reffig{SubSimplices} shows the projection to $\CC$ and the corresponding angles. Note that we may still cone to $\infty$, obtaining six simplices with the properties of \reflem{SixSimplices}, only now they overlap. However, by adding and subtracting volumes of overlapping simplices, we still will obtain the volume of the ideal tetrahedron. In particular, we have the following.
\begin{align*}
  \vol(\Delta(\alpha, \beta, \gamma)) & = 2\vol(S(\gamma)) + 2\vol(S(\beta)) - 2\vol(S(\pi-\alpha)) \\
  &= \Lambda(\gamma) + \Lambda(\beta) - \Lambda(\pi-\alpha)
\end{align*}
Since $\Lambda$ is an odd function and has period $\pi$, $-\Lambda(\pi-\alpha) = \Lambda(\alpha)$. Hence $\vol(\Delta(\alpha,\beta,\gamma)) = \Lambda(\alpha)+\Lambda(\beta)+\Lambda(\gamma)$ in this case as well. 
\end{proof}

The formula for volume of a tetrahedron has the following useful consequences.

\begin{theorem}\label{Thm:VolConcaveDown}
Let $\mathcal{A}$ be the set of possible angles on a tetrahedron:
\[ \mathcal{A} = \{(\alpha, \beta, \gamma)\in (0,\pi)^3 \mid \alpha+\beta+\gamma=\pi\}.\]
Then the function $\vol\from \mathcal{A}\to \RR$ given by \[ \vol(\alpha,\beta,\gamma)=\Lambda(\alpha)+\Lambda(\beta) +\Lambda(\gamma)\]
is strictly concave down on $\mathcal{A}$. Moreover, we can compute its first two derivatives. For $a=(a_1, a_2, a_3)\in\mathcal{A}$ a point and $w=(w_1,w_2,w_3) \in T_a\mathcal{A}$ a nonzero tangent vector, the first two derivatives of $\vol$ in the direction of $w$ satisfy
\[ \frac{\partial\vol}{\partial w} = \sum_{i=1}^3 -w_i\log\sin a_i, \quad
\frac{\partial^2\vol}{\partial w^2} < 0.\]
\end{theorem}

\begin{proof}
First, note that since $w$ is a tangent vector to $\mathcal{A}$, and the sum of the three coordinates of each point in $\mathcal{A}$ is $\pi$, it follows that $w_1+w_2+w_3=0$.

Next, by \refthm{VolTet}, the directional derivative of $\vol$ at $a$ in the direction of $w$ is given by
\begin{align*}
  \frac{\partial \vol}{\partial w} &= \sum_{i=1}^3 -w_i\log|2\sin a_i| \\
  & = \sum_{i=1}^3 w_i(-\log 2) + \sum_{i=1}^3 -w_i\log|\sin a_i| \\
  & = 0 + \sum_{i=1}^3 -w_i\log \sin a_i.
\end{align*}
The last line holds since $w_1+w_2+w_3=0$ and since $a_i\in(0,\pi)$, hence $\sin a_i >0$.

For the second derivative, we know $a_1+a_2+a_3=\pi$, so at least two of $a_1, a_2, a_3$ are strictly less than $\pi/2$. Without loss of generality, say $a_1$ and $a_2$ are less than $\pi/2$.

Then the second derivative is
\[ \frac{\partial^2\vol}{\partial w^2} = \sum_{i=1}^3 -w_i^2\cot a_i. \]
Since $a_3 = \pi-a_1-a_2$ and $w_3=-w_1-w_2$, we may write
\[ w_3^2\cot a_3 = (w_1+w_2)^2\cot(\pi-a_1-a_2) = -(w_1+w_2)^2\frac{\cot a_1\cot a_2 -1}{\cot a_1 + \cot a_2}, \]
where the last equality is an exercise in trig identities.

Then we obtain
\begin{align*}
  -\frac{\partial^2\vol}{\partial w^2} & = w_1^2\cot a_1 +w_2^2\cot a_2 -(w_1+w_2)^2\frac{\cot a_1 \cot a_2 -1}{\cot a_1 + \cot a_2} \\
  & = \frac{(w_1+w_2)^2 + (w_1\cot a_1 - w_2\cot a_2)^2}{\cot a_1 + \cot a_2}.
\end{align*}
The denominator of the last fraction is positive, because $a_1, a_2\in(0,\pi/2)$. The numerator is the sum of squares, hence at least zero. In fact, if it equals zero, then we have $w_1=-w_2$ and $\cot a_1 = -\cot a_2$. But $a_1,a_2\in(0,\pi/2)$, so this is impossible. Thus numerator and denominator are strictly positive, and so $\partial^2\vol/\partial w^2$ is strictly negative, hence strictly concave down. 
\end{proof}

\begin{theorem}\label{Thm:MaxVolTet}
  The regular ideal tetrahedron, with dihedral angles $\alpha=\beta=\gamma=\pi/3$, maximizes volume over all ideal tetrahedra.
\end{theorem}

\begin{proof}
Because $\vol$ is continuous, we know it obtains a maximum on the cube $[0,\pi]^3$. First we consider the boundary of that cube, and we show the maximum cannot occur there. If any angle is $\pi$, then $\alpha+\beta+\gamma=\pi$ implies the other two angles are $0$. Thus to show the maximum does not occur on the boundary of the cube, it suffices to show the maximum does not occur when one of the angles is zero. So suppose $\alpha=0$. 
Since $\Lambda(0)=\Lambda(\pi)=0$ by \refeqn{LobachevskiDilog}, and since $\Lambda(\beta)+\Lambda(\pi-\beta) = \Lambda(\beta)+\Lambda(-\beta)=0$ by \reflem{Lobachevsky}, the volume in this case will be $0$. So the maximum does not occur on the boundary.

Thus we seek a maximum in the interior. We maximize
$\vol(\alpha,\beta, \gamma) = \Lambda(\alpha)+\Lambda(\beta)+\Lambda(\gamma)$ subject to the constraint $\pi=\alpha+\beta+\gamma =: f(\alpha,\beta,\gamma)$. The theory of Lagrange multipliers tells us that at the maximum, there is a scalar $\lambda$ such that
\[ \nabla \vol = \lambda \nabla f, \quad \mbox{or}\]
\[ \log\sin\alpha = \log\sin\beta = \log\sin\gamma = \lambda. \]
This will be satisfied when $\sin\alpha=\sin\beta=\sin\gamma$. Since $\alpha, \beta, \gamma \in (0,\pi)$, and $\alpha+\beta+\gamma=\pi$, it follows that $\alpha=\beta=\gamma=\pi/3$, and the tetrahedron is regular.
\end{proof}

%%%%%%%%%%%%%%%%%%%%%%%%%%%%%%%%%%%%%%%%%%%%%%%%%%%%%%%%%%%%%%%%%
\section{Angle structures and the volume functional}

Note that in \refthm{VolTet}, we showed that the volume of an ideal tetrahedron can be computed given only its dihedral angles. A dihedral angle can be obtained by taking the imaginary part of the log of a tetrahedron's edge invariant. Thus the imaginary parts alone of the edge invariants allow us to assign a volume to the structure. These are exactly the angles of an angle structure.\index{angle structure}

Recall from \refdef{AngleStructures} that we defined an angle structure\index{angle structure} on an ideal triangulation $\mathcal{T}$ of a manifold $M$ to be a collection of (interior) dihedral angles satisfying:
\begin{enumerate}
\item[(0)] Opposite edges of a tetrahedron have the same angle.
\item[(1)] Dihedral angles lie in $(0,\pi)$.
\item[(2)] The sum of angles around any ideal vertex of any tetrahedron is $\pi$.
\item[(3)] The sum of angles around any edge class of $M$ is $2\pi$.
\end{enumerate}

The set of all angle structures\index{angle structure} for a triangulation $\mathcal{T}$ is denoted by $\mathcal{A}(\mathcal{T})$. For $M$ is an orientable 3-manifold with boundary consisting of tori, and $\mathcal{T}$ a triangulation of $M$, we will study the set of angle structures $\mathcal{A}(\mathcal{T})$. 

\begin{proposition}\label{Prop:AngleStructsPoly}
Let $\mathcal{T}$ be an ideal triangulation of a 3-manifold $M$ consisting of $n$ tetrahedra, and as usual denote the set of angle structures\index{angle structure} by $\mathcal{A}(\mathcal{T})$. If $\mathcal{A}(\mathcal{T})$ is nonempty, then it is a convex, finite-sided, bounded polytope in $(0,\pi)^{3n}\subset \RR^{3n}$.
\end{proposition}

\begin{proof}
For each tetrahedron of $\mathcal{T}$, an angle structure\index{angle structure} selects three dihedral angles lying in $(0,\pi)$. Thus $\mathcal{A}(\mathcal{T})$ is a subset of $(0,\pi)^{3n}$. The equations coming from conditions~(2) and~(3) are linear equations whose solution set is an affine subspace of $\RR^{3n}$. When we intersect the solution space with the cube $(0,\pi)^{3n}$, we obtain a bounded, convex, finite-sided polytope. 
\end{proof}

There is no guarantee that $\mathcal{A}(\mathcal{T})$ is nonempty. However, \refprop{AngleStructsPoly} implies that if it is nonempty, then we may view a point of $\mathcal{A}(\mathcal{T})$ as a point in $(0,\pi)^{3n}$. We write $a\in\mathcal{A}(\mathcal{T})$ as $a=(a_1, \dots, a_{3n})$. 

\begin{definition}\label{Def:VolFunctional}
The \emph{volume functional}\index{volume functional}
$\mathcal{V}\from \mathcal{A}(\mathcal{T})\to\RR$ is defined by
\[ \mathcal{V}(a_1, \dots, a_{3n}) = \sum_{i=1}^{3n} \Lambda(a_i). \]
\end{definition}

Thus $\mathcal{V}(a)$ is the sum of volumes of hyperbolic tetrahedra associated with the angle structure\index{angle structure} $a$.

A reason angle structures are so useful comes from the following two theorems.

\begin{theorem}[Volume and angle structures]\label{Thm:VolAngleStructs}
  Let $M$ be an orientable 3-manifold with boundary consisting of tori, with ideal triangulation $\mathcal{T}$. If a point $A\in\mathcal{A}(\mathcal{T})$ is a critical point for the volume functional\index{volume functional} $\mathcal{V}$ then the ideal hyperbolic tetrahedra obtained from the angle structure\index{angle structure} $A$ give $M$ a complete hyperbolic structure. 
\end{theorem}

The converse is also true:

\begin{theorem}\label{Thm:VolAngleStructsConverse}
If $M$ is finite volume hyperbolic 3-manifold with boundary consisting of tori such that $M$ admits a positively oriented hyperbolic ideal triangulation\index{positively oriented tetrahedron}\index{tetrahedron!positively oriented} $\mathcal{T}$, then the angle structure\index{angle structure} $A\in\mathcal{A}(\mathcal{T})$ giving the angles of $\mathcal{T}$ for the complete hyperbolic structure is the unique global maximum of the volume functional\index{volume functional} $\mathcal{V}$ on $\mathcal{A}(\mathcal{T})$.
\end{theorem}

The two theorems are attributed to Casson and Rivin, and follow from proofs in \cite{Rivin:EuclidStructs}. The first direct proof of the results are written in Chan's honors thesis \cite{Chan:Honours}, using work of Neumann and Zagier \cite{NeumannZagier}. A very nice self-contained exposition and proof of both theorems is given in \cite{FuterGueritaud:Survey}. We will follow the ideas of Futer and Gu{\'e}ritaud to show \refthm{VolAngleStructs} in \refsec{LeadingTrailing}.

To prove the converse, we will follow a simple proof of Chan using the Schl\"afli formula for the variation of volumes of ideal tetrahedra. Chan credits his proof to unpublished ideas of Schlenker.

%%%%%%%%%%%%%%%%%%%%%%%%%%%%%%%%%%%%%%%%%%%%%%%%%%%%%%%%%%%%%%%%%
\section{Leading--trailing deformations}\label{Sec:LeadingTrailing}

\begin{lemma}\label{Lem:DerivativesVol}
Let $M$ be an orientable 3-manifold with boundary consisting of tori, with ideal triangulation $\mathcal{T}$ consisting of $n$ tetrahedra. Then the volume functional\index{volume functional} $\mathcal{V}\from \mathcal{A}(\mathcal{T})\to \RR$ is strictly concave down on $\mathcal{A}(\mathcal{T})$. For $a=(a_1, \dots, a_{3n})\in\mathcal{A}(\mathcal{T})$ and $w=(w_1, \dots, w_{3n})\in T_a\mathcal{A}(\mathcal{T})$ a non-zero tangent vector, the first two directional derivatives of $\mathcal{V}$ satisfy
\[ \frac{\partial \mathcal{V}}{\partial w} = \sum_{i=1}^{3n} -w_i\log\sin a_i \quad \mbox{and} \quad \frac{\partial^2 \mathcal{V}}{\partial w^2} < 0. \]
\end{lemma}

\begin{proof}
Because the volume functional\index{volume functional} $\mathcal{V}$ is the sum of volumes of ideal tetrahedra, the formulas for derivatives follow by linearity from \refthm{VolConcaveDown}. Because the second derivative is strictly negative, the volume functional is strictly concave down. 
\end{proof}

We will need to take derivatives in carefully specified directions. To that end, we now define a vector $w=(w_1, \dots, w_{3n})\in \RR^{3n}$ and show that $w$ lies in $T_a\mathcal{A}(\mathcal{T})$. Again the ideas follow from \cite{FuterGueritaud:Survey}.

\begin{definition}\label{Def:LeadingTrailing}
Let $C$ be a cusp of $M$ with a cusp triangulation corresponding to the ideal tetrahedra of $\mathcal{T}$. Let $\zeta$ be an oriented closed curve on $C$, isotoped to run monotonically through the cusp triangulation, as in \refdef{CompletenessEquations}. Let $\zeta_1, \dots, \zeta_k$ be the oriented segments of $\zeta$ in distinct triangles. For the segment $\zeta_i$ in triangle $t_i$, define the \emph{leading corner} of $t_i$ to be the corner of the triangle that is opposite the edge where $\zeta_i$ enters $t_i$, and define the \emph{trailing corner} to be the corner opposite the edge where $\zeta_i$ exits.

Each corner of the triangle $t_i$ is given a dihedral angle $a_j$ in an angle structure,\index{angle structure} thus corresponds to a coordinate of $\mathcal{A}(\mathcal{T}) \subset \RR^{3n}$. Similarly for any $a\in \mathcal{A}(\mathcal{T})$, each corner of $t_i$ corresponds to a coordinate of the tangent space $T_a\mathcal{A}(\mathcal{T}) \subset \RR^{3n}$. 

We define a vector $w(\zeta_i) \in \RR^{3n}$ by setting the coordinate corresponding to the leading corner of $t_i$ equal to $+1$, and the coordinate corresponding to the trailing corner of $t_i$ equal to $-1$. Set all other coordinates equal to zero. The \emph{leading--trailing deformation}\index{leading--trailing deformation} corresponding to $\zeta$ is defined to be the vector $w(\zeta)=\sum_i w(\zeta_i)$.
\end{definition}

\begin{figure}
  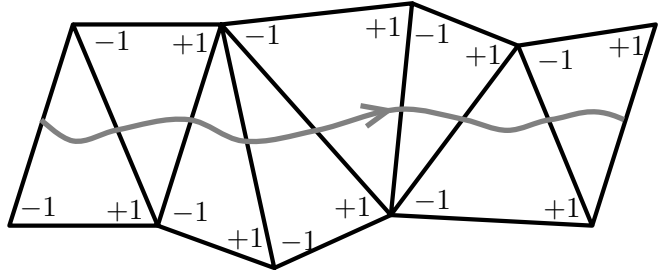
  \caption{For the oriented curve $\zeta$ shown, leading corners are marked with $+1$ and trailing corners with $-1$}
  \label{Fig:LeadingTrailing}
\end{figure}

An example is shown in \reffig{LeadingTrailing}.

\begin{lemma}\label{Lem:TantoAngleStruct}
Let $\sigma$ be a curve encircling a vertex of the cusp triangulation on cusp $C$. Let $\mu$ be an embedded curve isotopic to a generator of the holonomy group\index{holonomy group} of the cusp torus. Then the corresponding leading--trailing deformation vectors $w(\sigma)$ and $w(\mu)$ both lie in the tangent space $T_a\mathcal{A}(\mathcal{T})$, for any $a\in\mathcal{A}(\mathcal{T})$.
\end{lemma}

\begin{proof}
The space $\mathcal{A}(\mathcal{T})$ is a submanifold of $\RR^3$ cut out by linear equations corresponding to (2) and (3) of \refdef{AngleStructures}, namely that angles at each ideal vertex of a tetrahedron sum to $\pi$, and angles about an edge of $M$ sum to $2\pi$. Let $f_i(a) = a_i+a_{i+1}+a_{i+2}$ be the sum of angles of the $i$-th tetrahedron, and let $g_e(a) = \sum a_{e_i}$ be the sum of angles about the edge $e$. So $\mathcal{A}(\mathcal{T}) \subset (0,\pi)^{3n}$ is the space cut out by all equations $f_i=\pi$ and $g_e = 2\pi$. Thus to see that $w(\zeta)$ is a tangent vector to $\mathcal{A}(\mathcal{T})$ at a point $a$, we need to show that the vector is orthogonal to the gradient vectors $\nabla f_i$ and $\nabla g_e$ at $a$, for all $i$ and all $e$. 

Note $\nabla f_i$ is the vector $(0,\dots,0,1,1,1,0,\dots,0)$, with $0$s away from the $i^{\rm{th}}$ tetrahedron $t_i$ and $1$s in the three positions corresponding to the angles of $t_i$. There are four cusp triangles coming from this tetrahedron, corresponding to its four ideal vertices. Suppose $\zeta$ is a curve in the cusp triangulation of $C$. If no segment of $\zeta$ runs through a triangle of $t_i$, then $w(\zeta)$ has only $0$s in the position corresponding to the $1$s of $\nabla f_i$, hence $\nabla f_i \cdot w(\zeta)=0$ in this case.
So suppose that some segment $\zeta_j$ of $\zeta$ runs through a triangle of $t_i$. Then one corner of the triangle is a leading corner for $\zeta_j$, and one is a trailing corner, so $w(\zeta_j)$ has one $0$, one $+1$, and one $-1$ in the three positions corresponding to angles of $t_i$. Hence $\nabla f_i \cdot w(\zeta_j)=0$. By linearity, $\nabla f_i\cdot w(\zeta)=0$.

So it remains to show that for each edge $e$, $\nabla g_e \cdot w(\zeta) = 0$, where $\zeta$ is one of the curves $\sigma$ or $\mu$ in the hypothesis of the lemma. Note that $\nabla g_e$ is a vector $(\epsilon_1, \dots, \epsilon_{3n})$, where $\epsilon_j$ is one of the integers $0$, $1$, or $2$, counting the number of times a dihedral angle of a tetrahedron occurs in the gluing equation $g_e$. We will consider the segments of $\zeta$ one at a time. Note that any segment $\zeta_j$ of $\zeta$ contributes $0$, $+1$, and $-1$ to opposite edges of exactly one tetrahedron $t_j$, as illustrated in \reffig{LeadingTrailingVertex}, left. 

\begin{figure}
  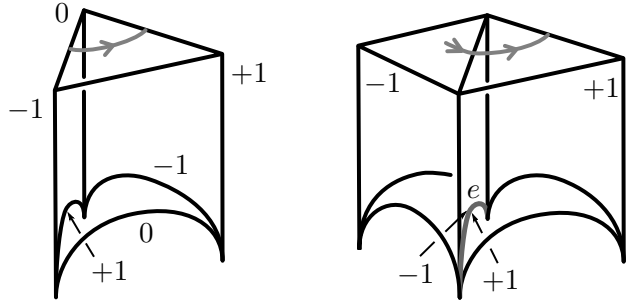
\caption{Left: Effect of $w(\sigma_j)$ on the edges of a tetrahedron. Right: 
  If the lower edge is identified to $e$, then the contribution of $+1$ from $w(\zeta_j)$ to $e$ cancels with the $-1$ contribution from $w(\zeta_{j-1})$.}
  \label{Fig:LeadingTrailingVertex}
\end{figure}

If the edge $e$ is not identified to any of the edges of the tetrahedron $t_j$, then $\nabla g_e\cdot w(\zeta_j)=0$. Similarly, if $e$ is identified only to one or both of the edges for which $w(\zeta_j)$ contributes a $0$, then although the corresponding coordinate of $\nabla g_e$ will be $1$ or $2$, the dot product $\nabla g_e \cdot w(\zeta_j)$ will still be $0$.

If $e$ is identified to one or both of the edges labeled with a $+1$ by $w(\zeta_j)$, then there will be a contribution of $+1$ or $+2$ (respectively) to $\nabla g_e \cdot w(\zeta_j)$ coming from these labels. We will show that in this case, there exists one or two (respectively) segments of $\zeta$ each contributing $-1$, so that the positive contributions cancel.

Suppose first that $e$ is identified to the non-vertical edge of $t_j$ labeled $+1$. Then consider the segment $\zeta_{j-1}$. This lies in a tetrahedron $t_{j-1}$ glued to $t_j$ along a face containing the edge identified to $e$. In the cusp triangulation, $\zeta_{j-1}$ exits its cusp triangle at this face. Thus the opposite corner of the cusp triangle is a trailing corner, and is assigned a $-1$. This trailing corner corresponds to an edge opposite $e$. So $e$ picks up a $-1$ from $w(\zeta_{j-1})$. See \reffig{LeadingTrailingVertex}, right. Hence the $+1$ contribution of $w(\zeta_j)$ is canceled in this case with this $-1$ from $w(\zeta_{j-1})$.

Now suppose $e$ is identified to the vertical edge of $t_j$ labeled $+1$. Let $\zeta_j, \zeta_{j+1}, \dots, \zeta_{j+r}$ be a maximal collection of segments in cusp triangles adjacent to $e$. Note if $\zeta=\sigma$ encircles a vertex, that vertex will not correspond to the endpoint of $e$. Then $r=1$, i.e.\ there are just two segments of $\zeta$ adjacent to $e$. If $\zeta=\mu$ is a generator of cusp homology, then $r\geq 1$. Because we are assuming $\zeta$ is embedded and meets each edge of the cusp triangulation at most once, we know $\zeta_j, \dots, \zeta_{j+r}$ do not encircle $e$ completely in this case. See \reffig{LeadingTrailingTan}.

\begin{figure}
%% Accompanies image file 'F9-06-LTMu.eps' (pdf, eps, ps)
\begingroup%
  \makeatletter%
  \providecommand\color[2][]{%
    \errmessage{(Inkscape) Color is used for the text in Inkscape, but the package 'color.sty' is not loaded}%
    \renewcommand\color[2][]{}%
  }%
  \providecommand\transparent[1]{%
    \errmessage{(Inkscape) Transparency is used (non-zero) for the text in Inkscape, but the package 'transparent.sty' is not loaded}%
    \renewcommand\transparent[1]{}%
  }%
  \providecommand\rotatebox[2]{#2}%
  \newcommand*\fsize{\dimexpr\f@size pt\relax}%
  \newcommand*\lineheight[1]{\fontsize{\fsize}{#1\fsize}\selectfont}%
  \ifx\svgwidth\undefined%
    \setlength{\unitlength}{136.92424679bp}%
    \ifx\svgscale\undefined%
      \relax%
    \else%
      \setlength{\unitlength}{\unitlength * \real{\svgscale}}%
    \fi%
  \else%
    \setlength{\unitlength}{\svgwidth}%
  \fi%
  \global\let\svgwidth\undefined%
  \global\let\svgscale\undefined%
  \makeatother%
  \begin{picture}(1,0.99041756)%
    \lineheight{1}%
    \setlength\tabcolsep{0pt}%
    \put(0,0){\includegraphics[width=\unitlength]{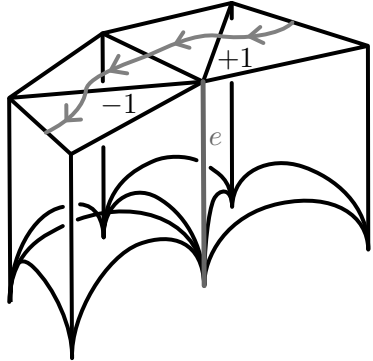}}%
    \put(0.55463674,0.59688313){\color[rgb]{0.4,0.4,0.4}\makebox(0,0)[lt]{\lineheight{1.25}\smash{\begin{tabular}[t]{l}$e$\end{tabular}}}}%
    \put(0.5767691,0.81129051){\color[rgb]{0,0,0}\makebox(0,0)[lt]{\lineheight{1.25}\smash{\begin{tabular}[t]{l}$+1$\end{tabular}}}}%
    \put(0.25630698,0.68308535){\color[rgb]{0,0,0}\makebox(0,0)[lt]{\lineheight{1.25}\smash{\begin{tabular}[t]{l}$-1$\end{tabular}}}}%
  \end{picture}%
\endgroup%

  \caption{If $w(\zeta_j)$ contributes $+1$ to a vertical edge meeting $e$, there is a maximal collection of segments running through cusp triangles adjacent to $e$.}
  \label{Fig:LeadingTrailingTan}
\end{figure}

In both cases $\zeta=\sigma$ and $\zeta=\mu$, for segments $\zeta_{j+k}$ with $0<k<r$, note $w(\zeta_{j+k})$ contributes only $0$s to the edge $e$. Since the segment after $\zeta_{j+r}$ is no longer adjacent to the vertical edge $e$, it follows that $w(\zeta_{j+r})$ contributes $-1$ to $e$. Then the $+1$ contribution from $w(\zeta_j)$ cancels with the $-1$ contribution from $w(\zeta_{j+r})$.

Finally, it could be the case that $e$ is identified to both edges labeled $+1$ by $w(\zeta_j)$, so that $\nabla g_e$ has a $2$ in that coordinate and $\nabla g_e\cdot w(\zeta_j)$ picks up a $+2$ from these two edges. But in this case, combining both arguments above implies that one of the $+1$ contributions is canceled by a $-1$ coming from $w(\zeta_{j-1})$ and one by a $-1$ coming from $w(\zeta_{j+r})$ for appropriate $r$. Thus both are canceled.

We have shown that for each $j$, each $+1$ contribution of $w(\zeta_j)$ to $\nabla g_e \cdot w(\zeta_j)$ is canceled by a $-1$ contribution from some $w(\zeta_k)$. Provided none of the $-1$ contributions from $w(\zeta_k)$ are repeated for distinct $j$, this shows that $\nabla g_e \cdot w(\zeta) \leq 0$. The fact that these contributions are not repeated follows from the uniqueness of the choice of $\zeta_{j-1}$ and $\zeta_{j+r}$.

A similar argument implies $\nabla g_e \cdot w(\zeta) \geq 0$. Thus $\nabla g_e\cdot w(\zeta)=0$, as desired. 
\end{proof}

\begin{lemma}\label{Lem:ReH}
Let $\zeta$ be one of the curves $\sigma$ or $\mu$ of \reflem{TantoAngleStruct}, and let  $w(\zeta)\in T_a\mathcal{A}(\mathcal{T})$ be the corresponding leading--trailing deformation vector. Let $H(\zeta)$ be the complex number associated to the curve $\zeta$ given in \refdef{CompletenessEquations} (completeness equations). Then
\[ \frac{\partial \mathcal{V}}{\partial w(\zeta)} = \Re(\log H(\zeta)). \]
\end{lemma}

\begin{proof}
Let $\zeta_1, \dots, \zeta_k$ denote segments of $\zeta$ in cusp triangles $t_1, \dots, t_k$, respectively. Label the dihedral angles of triangle $t_i$ by $\alpha_i$, $\beta_i$, $\gamma_i$, in clockwise order, so that $\alpha_i$ is the angle cut off by $t_i$. By \refdef{CompletenessEquations},
\[ \Re(\log H(\zeta)) = \sum_i \epsilon_i \Re (\log |z(\alpha_i)|), \]
where $z(\alpha_i)$ is the edge invariant associated with the edge labeled $\alpha_i$, and $\epsilon_i =+1$ if $\alpha_i$ is to the left of $\zeta_i$ and $\epsilon_i=-1$ if $\alpha_i$ is to the right of $\zeta_i$.

On the other hand, comparing \reffig{ExampleCusp} and \reffig{LeadingTrailing}, we see that when $\alpha_i$ is to the left of $\zeta_i$, the vector $w(\zeta_i)$ has a $+1$ in the position corresponding to $\beta_i$ and a $-1$ in the position corresponding to $\gamma_i$, and when $\alpha_i$ is to the right of $\zeta_i$, the vector $w(\zeta_i)$ has a $-1$ in the position corresponding to $\beta_i$ and a $+1$ in the position corresponding to $\gamma_i$. Then \reflem{DerivativesVol} implies that
\begin{align*}
\frac{\partial \mathcal{V}}{\partial w} &= \sum_{j=1}^{3n} -w_j\log\sin a_j \\
&= \sum_i -\epsilon_i \log\sin\beta_i + \epsilon_i\log\sin\gamma_i \\
&= \sum_i \epsilon_i \log \left( \frac{\sin\gamma_i}{\sin\beta_i} \right)\\
& = \sum_i \epsilon_i \Re(\log |z(\alpha_i)|), \quad \mbox{by \refeqn{AngleEdgeInvariant}}
\end{align*}
This is what we needed to show.
\end{proof}

We now have the tools we need to prove \refthm{VolAngleStructs}, to show that a critical point $a\in\mathcal{A}(\mathcal{T})$ of the volume functional\index{volume functional} corresponds to a complete hyperbolic structure on the manifold $M$.

\begin{proof}[Proof of \refthm{VolAngleStructs}]
Suppose $a\in\mathcal{A}(\mathcal{T})$ is a critical point of the volume functional\index{volume functional} $\mathcal{V}$. Then $a$ assigns a dihedral angle to each tetrahedron of $\mathcal{T}$, giving each ideal tetrahedron a unique hyperbolic structure. By \refthm{Gluing}, gluing these tetrahedra will give a hyperbolic structure on $M$ if and only if the edge gluing equations are satisfied for each edge. By \refthm{EuclidCusp}, the hyperbolic structure will be complete if and only if the induced geometric structure on each cusp torus is a Euclidean structure,\index{Euclidean structure} and we obtain a Euclidean structure when the completeness equations are satisfied by \refprop{CompletenessEqns}.

Consider first the edge gluing equations. Notice that any angle structure\index{angle structure} gives hyperbolic ideal tetrahedra satisfying the imaginary part of the gluing equations, so we need to show that our tetrahedra satisfy the real part. Fix an edge of the triangulation, and let $\sigma$ be a curve on a cusp torus encircling an endpoint of that edge. The real part of the gluing equation corresponding to this edge will be satisfied if and only if $\Re (\log H(\sigma)) = 0$. But \reflem{ReH} implies that $\Re(\log H(\sigma)) = \frac{\partial\mathcal{V}}{\partial w(\sigma)}$, and this is zero because our angle structure\index{angle structure} is a critical point of the volume functional.\index{volume functional} So the gluing equations hold.

As for the completeness equations, for any cusp torus $C$, and $\mu_1$ and $\mu_2$ generators of the first homology group of $C$, the completeness equations require that $H(\mu_1)=H(\mu_2)=1$. By \reflem{ReH}, we know
\[ \Re\log(H(\mu_i)) = \frac{\partial \mathcal{V}}{\partial w(\mu_i)} = 0,\]
since $a$ is a critical point. Thus the real part of each of these completeness equations is satisfied.

Consider the developing image of a fundamental domain for the cusp torus $C$. Because we know the angles given by $a$ satisfy the gluing equations, the structure on $C$ is at least an affine structure on the torus. Therefore the developing image of the fundamental domain is a quadrilateral in $\CC$. Since the real parts of the completeness equations for $\mu_1$ and $\mu_2$ are both satisfied, it follows that the holonomy\index{holonomy} elements corresponding to $\mu_1$ and $\mu_2$ do not scale either side of the fundamental domain. But then the holonomy elements cannot effect a non-trivial rotation either; a quadrilateral in $\CC$ whose opposite sides are the same length is a parallelogram. Thus the developing image of a fundamental domain is a parallelogram, the holonomy elements corresponding to $\mu_1$ and $\mu_2$ must be pure translations, and the cusp torus admits a Euclidean structure.\index{Euclidean structure} So the completeness equations hold. 
\end{proof}

\Refthm{VolAngleStructs} gives us a way of proving not only that a 3-manifold $M$ is hyperbolic, but also that it admits a positively oriented\index{positively oriented tetrahedron}\index{tetrahedron!positively oriented} geometric triangulation.\index{geometric triangulation} To use the theorem, first, fix a triangulation $\mathcal{T}$. Then show the space of angle structures\index{angle structure} $\mathcal{A}(\mathcal{T})$ is nonempty. Finally, show that the volume functional\index{volume functional} achieves its maximum in the interior of $\overline{\mathcal{A}(\mathcal{T})}$. The last step can often be accomplished by considering angle structures\index{angle structure} on the boundary $\overline{\mathcal{A}(\mathcal{T})}-\mathcal{A}(\mathcal{T})$ and proving such structures cannot maximize volume. We will follow exactly this procedure for 2-bridge knots\index{2-bridge knot or link} in \refchap{TwoBridge}.

The following proposition is a useful tool for examining the maximum of the volume functional\index{volume functional} on the boundary $\overline{\mathcal{A}(\mathcal{T})}-\mathcal{A}(\mathcal{T})$.

\begin{proposition}\label{Prop:NoSingleDegenerate}
  Suppose an angle structure\index{angle structure} $a \in \overline{\mathcal{A}(\mathcal{T})}$ maximizes the volume functional\index{volume functional} $\mathcal{V}$. Suppose that for some tetrahedron $\Delta_i$, one of the three angles of $\Delta_i$ in the angle structure $a$ is $0$. Then two of the angles are $0$ and the third is $\pi$. 
\end{proposition}

\begin{proof}
Suppose instead that one angle, say $a_{i}$ is $0$, but the other two angles of $\Delta_i$ are nonzero: $a_{i+1} \neq 0$ and $a_{i+2}\neq 0$. We will find a path through $\mathcal{A}(\mathcal{T})$ with endpoint the angle structure\index{angle structure} $a$, and we will show that the derivative of this path is positive, and in fact unbounded, as it approaches the endpoint corresponding to $a$. It will follow that $a$ cannot be a maximum, which contradicts our assumption on $a$.

The space of angle structures\index{angle structure} $\mathcal{A}(\mathcal{T})$ is a bounded open convex subset of $\RR^{3n}$. Its tangent space can be extended to its boundary. We may choose a tangent vector $w$ in $T_a\overline{\mathcal{A}(\mathcal{T})}$ pointing into the interior of $\mathcal{A}(\mathcal{T})$, and take the path corresponding to geodesic flow in the direction of this tangent vector. \Refthm{VolConcaveDown} implies that the derivative of the volume functional\index{volume functional} along this path is the sum of terms of the form $\sum_{i=1}^{3n} -w_{i}\log\sin(a_{i})$.

Consider the contribution from $\Delta_i$. The terms
\[ w_{i+1}\log\sin(a_{i+1}) \mbox{ and } w_{i+2}\log\sin(a_{i+2}) \]
are bounded, since $a_{i+1}$ and $a_{i+2}$ are bounded away from zero. But as the path approaches the angle structure\index{angle structure} $a$, the term coming from $-w_{i}\log\sin(a_{i})$ approaches positive infinity. Thus such a point cannot be a maximum. 
\end{proof}

%%%%%%%%%%%%%%%%%%%%%%%%%%%%%%%%%%%%%%%%%%%%%%%%%%%%%%%%%%%%%%%%%
\section{The Schl{\"a}fli formula}\label{Sec:NeumannZagier}

In this short section, we prove \refthm{VolAngleStructsConverse}, the converse to \refthm{VolAngleStructs}. Our proof uses the Schl{\"a}fli formula for ideal tetrahedra, which can be stated as follows.

\begin{theorem}[Schl{\"a}fli's formula for ideal tetrahedra]\label{Thm:Schlafli}\index{Schl{\"a}fli's formula}
Let $P$ be an ideal tetrahedron. Let $H_1, \dots, H_n$ be a collection of horospheres centered on the ideal vertices of $P$. For each edge $e_{ij}$, running between the $i$-th to the $j$-th ideal vertices of $P$, let $\ell(e_{ij})$ denote the signed distance between $H_i$ and $H_j$ (that is, $\ell(e_{i,j})$ is defined to be negative if $H_i\cap H_j \neq \emptyset$). Finally, let $\theta_{ij}$ denote the dihedral angle along edge $e_{i,j}$. Then the variation in the volume of $P$ satisfies
\begin{equation}\label{Eqn:Schlafli}
  d\mathcal{V}(P) = -\frac{1}{2}\sum_{i,j} \ell(e_{ij}) d\theta_{i,j}.
\end{equation}
\end{theorem}

Schl{\"a}fli's formula was originally proved for finite spherical simplices by Schl{\"a}fli in the 1850s. It has been extended in many directions, including to finite and ideal polyhedra in spaces of constant curvature. A proof of a formula that contains the result in \refthm{Schlafli} can be found in \cite{Milnor:Schlafli}; see also \cite{Rivin:EuclidStructs}. These sources note that the right hand side of \refeqn{Schlafli} is independent of the choice of horospheres. 

Using this, we can finish the proof. 

\begin{proof}[Proof of \refthm{VolAngleStructsConverse}]
Choose a horosphere about each cusp in the complete hyperbolic structure on $M$. Because the hyperbolic structure is complete, this choice gives a well-defined horosphere about each ideal vertex of each ideal tetrahedron in the positively oriented\index{positively oriented tetrahedron}\index{tetrahedron!positively oriented} hyperbolic ideal triangulation $\mathcal{T}$. Thus for each tetrahedron, we may use this choice to define the edge lengths  $\ell(e_{ij})$ of \refthm{Schlafli}.

Now note that because the total angle around each edge is a constant $2\pi$, the contributions to the variation of the volume coming from each simplex add to zero for each edge. Thus the right hand side of \refeqn{Schlafli} is zero. It follows that the complete structure is a critical point for the volume functional.\index{volume functional}

On the other hand, since the volume functional\index{volume functional} is strictly concave down on $\mathcal{A}$, \refthm{VolConcaveDown}, it must follow that the complete structure is the unique global maximum. 
\end{proof}

%%%%%%%%%%%%%%%%%%%%%%%%%%%%%%%%%%%%%%%%%%%%%%%%%%%%%%%%%%%%%%%%%
\section{Consequences}

\Refthm{VolAngleStructs} and its converse have a number of important immediate consequences. We leave many proofs as exercises. 

\begin{corollary}[Lower volume bounds, angle structures]\label{Cor:VolumeBound}\index{volume bound!lower bound from angle structures}
Suppose $M$ has ideal triangulation $\mathcal{T}$ such that the volume functional\index{volume functional} $\mathcal{V}\from \mathcal{A}(\mathcal{T})\to\RR$ has a critical point $p \in \mathcal{A}(\mathcal{T})$.\index{angle structure} Then for any other point $q\in\overline{\mathcal{A}(\mathcal{T})}$, the volume functional satisfies
\[ \mathcal{V}(q) \leq \vol(M), \]
with equality if and only if $q=p$, i.e.\ $q$ also gives the complete hyperbolic metric on $M$.
\end{corollary}

\begin{proof}
By \reflem{DerivativesVol}, the volume functional\index{volume functional} is strictly concave down on $\overline{\mathcal{A}(\mathcal{T})}$, and so for any point $q\in \overline{\mathcal{A}(\mathcal{T})}$, $\mathcal{V}(q)$ is at most the maximum value of the volume functional, which is the value $\mathcal{V}(p)$ by hypothesis, and with equality if and only if $q=p$. By \refthm{VolAngleStructs}, $\vol(M)=\mathcal{V}(p)$. 
\end{proof}

More is conjectured to be true. \Refcor{VolumeBound} only gives a bound when there is a known critical point of the volume functional\index{volume functional} in the interior of the space of angle structures.\index{angle structure} If the maximum of the volume functional occurs on the boundary, it still seems to be the case in practice that the maximum is bounded by the volume of the complete hyperbolic structure. However, the following conjecture is currently still open. 

\begin{conjecture}[Casson's conjecture]\label{Conj:Casson}\index{Casson's conjecture}
Let $M$ be a cusped hyperbolic 3-manifold, and let $\mathcal{T}$ be any ideal triangulation of $M$. If the space of angle structures\index{angle structure} $\mathcal{A}(\mathcal{T})$ is nonempty, then the maximum value for the volume functional\index{volume functional} on $\overline{\mathcal{A}(\mathcal{T})}$ is at most the volume of the complete hyperbolic structure on $M$. 
\end{conjecture}

We may use angle structures to find hyperbolic Dehn fillings\index{Dehn filling} of triangulated 3-manifolds as well. Recall from \refchap{CompletionDehnFilling} that the $(p,q)$ Dehn filling on a triangulated manifold satisfies equation \refeqn{DehnFillingEquation}:
\[
  p \log H(\mu) + q\log H(\lambda) = 2\pi i.
\]

\begin{theorem}[Angle structures and Dehn fillings]\label{Thm:DehnFillingTriangulated}
  Let $M$ be a manifold with torus boundary components $T_1, \dots, T_n$, with generators $\mu_j, \lambda_j$ of $\pi_1(T_j)$ for each $j$. For each $j$, let $(p_j, q_j)$ denote a pair of relatively prime integers. Let $\mathcal{A}_{(p_1, q_1),\dots,(p_n,q_n)} \subset \mathcal{A}$ be the set of all angle structures\index{angle structure} that satisfy the imaginary part of the Dehn filling\index{Dehn filling} equations:
  \[ \Im( p_j\log H(\mu_j) + q_j\log H(\lambda_j) ) = 2\pi. \]
Then a critical point of the volume functional\index{volume functional} $\mathcal{V}$ on $\mathcal{A}_{(p_1, q_1),\dots,(p_n,q_n)}$ gives the complete hyperbolic structure on the Dehn filling $M((p_1, q_1), \dots, (p_n,q_n))$ of $M$.
\end{theorem}

\begin{proof}
As in the proof of \refthm{VolAngleStructs}, we will have a complete hyperbolic structure on the Dehn filling\index{Dehn filling} if and only if each edge gluing equation is satisfied and additionally each Dehn filling equation
\[ p_j \log H(\mu_j) + q_j\log H(\lambda_j) = 2\pi i \]
is satisfied.

The proof that edge gluing equations are satisfied follows exactly as in the proof of \refthm{VolAngleStructs}. As for the Dehn filling equations, the imaginary part of each equation is satisfied by the given constraint on the space of angle structures.\index{angle structure} By \reflem{ReH}, the real part satisfies
\[ \Re ( p_j\log(H(\mu_j)) + q_j\log(H(\lambda_j))) = p_j \frac{\partial \mathcal{V}}{\partial w(\mu_j)} + q_j\frac{\partial \mathcal{V}}{\partial w(\lambda_j)} = 0,\]
because this is a critical point for the volume functional.\index{volume functional} Thus each Dehn filling equation is satisfied.
\end{proof}

\Refcor{VolumeBound} and \refthm{HypDehnSurgery} imply the following (weaker) version of Thurston's theorem on volume change under Dehn filling,\index{Dehn filling} \refthm{VolumeDF}.

\begin{corollary}\label{Cor:VolDecreaseDehnFilling}
Let $M$ be as in \refthm{DehnFillingTriangulated}. 
If $s$ is a slope in the neighborhood of $\infty$ provided by Thurston's hyperbolic Dehn filling theorem, \refthm{HypDehnSurgery}, then the volume of $M(s)$ is strictly smaller than the volume of $M$. 
\end{corollary}

\begin{proof}
  \Refex{VolDecreaseDehnFilling}.
\end{proof}

We also have the tools to prove a rigidity theorem originally due to Weil. The following follows from Mostow--Prasad rigidity,\index{Mostow--Prasad rigidity} \refthm{MostowGeom}, but was first proved over a decade before that theorem, and can now be proved easily using angle structures.\index{angle structure} 

\begin{corollary}[Weil rigidity theorem]\label{Cor:WeilRigidity}\index{Weil rigidity theorem}
  Suppose $M$ is a 3-manifold with boundary consisting of tori, and suppose the interior of $M$ admits a complete hyperbolic metric. Then the metric is locally rigid, i.e.\ there is no local deformation of the metric through complete hyperbolic structures. 
\end{corollary}

\begin{proof}
  \Refex{WeilRigidity}. 
\end{proof}

%%%%%%%%%%%%%%%%%%%%%%%%%%%%%%%%%%%%%%%%%%%%%%%%%%%%%%%%%%%%%%%%%%
\section{Exercises}

\begin{exercise}
Find a formula for volume of a tetrahedron with ideal vertices $0$, $1$, $\infty$ and $z$ in terms of $z$ alone. 
\end{exercise}

\begin{exercise}
  Give a proof that the figures of \reffig{SubSimplices} are correct. That is, given a Euclidean triangle with vertices on the unit circle, and angles $\alpha$, $\beta$, and $\gamma$, prove that the angles around the origin are given as shown in the figure.
\end{exercise}

\begin{exercise}
 Prove the Kubert identities for the Lobachevsky function:\index{Kubert identities}
  \[ \Lambda(n\theta) = \sum_{k=0}^{n-1} n\Lambda(\theta + k\pi/n). \]
  Hint: cyclotomic identities:
  \[ 2\sin(n\theta) = \prod_{k=0}^{n-1} 2\sin\left( \theta + \frac{k\pi}{n}\right). \]
\end{exercise}  

\begin{exercise}
Find an explicit convex polytope describing the set of angle structures\index{angle structure} on the complement of the figure-8 knot.
\end{exercise}

\begin{exercise}
Find an explicit convex polytope describing the set of angle structures\index{angle structure} on the complement of the $5_2$ knot.
\end{exercise}

\begin{exercise}
  Suppose an ideal tetrahedron has dihedral angles $\alpha$, $\beta$, $\gamma$ in clockwise order. Prove equation \eqref{Eqn:AngleEdgeInvariant}: that the edge invariant of the tetrahedron assigned to the edge with angle $\alpha$ is
  \[ z(\alpha) = \frac{\sin(\gamma)}{\sin(\beta)}\,e^{i\alpha}. \]
\end{exercise}

\begin{exercise}\label{Ex:RegIdealOctahedron}
  An ideal octahedron can be obtained by gluing two identical ideal pyramids over an ideal quadrilateral base along the ideal quadrilateral. We triangulate this by running an edge from the ideal point opposite the quadrilateral on one pyramid, through the quadrilateral, to the ideal point opposite the quadrilateral on the other pyramid, and then we stellar subdivide. Using angle structures\index{angle structure} on this collection of ideal tetrahedra, prove that the maximal volume ideal hyperbolic octahedron is the regular one: the one for which the quadrilateral base is a square.\index{ideal octahedron, regular}\index{regular ideal octahedron}
\end{exercise}

\begin{exercise}\label{Ex:DoubleIdealPyramids}
  Generalize \refex{RegIdealOctahedron} to the ideal object obtained by gluing two ideal pyramids over an ideal $n$-gon base. Using stellar subdivision and angle structures,\index{angle structure} prove that the volume of the ideal double pyramid with base an $n$-gon is maximized when the $n$-gon is regular.
\end{exercise}

\begin{exercise}\label{Ex:VolDecreaseDehnFilling}
  Prove that volume decreases locally under Dehn filling,\index{Dehn filling} \refcor{VolDecreaseDehnFilling}.
\end{exercise}

\begin{exercise}\label{Ex:WeilRigidity}
  Prove the Weil rigidity theorem, \refcor{WeilRigidity}. First prove it for manifolds admitting an angle structure.\index{angle structure} Extend to all manifolds using the following theorem of Luo, Schleimer, and Tillmann, which can be found in  \cite{LuoSchleimerTillmann}. (You may assume the theorem for the exercise.)
  
  \begin{theorem}[Geometric triangulations exist virtually]\label{Thm:VirtualTriangulation}\index{geometric triangulation!exist virtually}
    Let $M$ be a 3-manifold with boundary consisting of tori, such that the interior of $M$ admits a complete hyperbolic structure. Then $M$ has a finite cover $N$ such that $N$ decomposes into positively oriented ideal tetrahedra.\index{positively oriented tetrahedron}\index{tetrahedron!positively oriented}
  \end{theorem}
\end{exercise}

%% Ch10_TwoBridge.tex

\chapter{Two-Bridge Knots and Links}\label{Chap:TwoBridge}
\blfootnote{Jessica S. Purcell, Hyperbolic Knot Theory}

In this chapter we will study in detail a class of knots and links that has particularly nice geometry, namely the class of 2-bridge knots and links.\index{2-bridge knot or link} The key property of these links that we will explore is the fact that they admit a geometric triangulation\index{geometric triangulation} that can be read off of a diagram, or an algebraic description of the link. In this chapter, we will define 2-bridge knots and links, describe their triangulations, and mention some of their geometric properties and consequences. 

\section{Rational tangles and 2-bridge links}

In 1956, H.~Schubert showed that the class of 2-bridge knots and links are classified by a rational number, and any continued fraction expansion of this number gives a diagram of the knot \cite{Schubert}. In this section, we work through the description of 2-bridge knots and links via rational numbers. Additional references are \cite{BurdeZieschang} and \cite{murasugi}. 

First, we define tangles.

\begin{definition}\label{Def:Tangle}
  A \emph{tangle}\index{tangle} is a 1-manifold properly embedded in a 3-ball $B$. That is, it is a collection of arcs with endpoints on $\bdy B$ and interiors disjointly embedded in the interior of $B$, possibly along with a collection of simple closed curves. For our purposes, we will consider only tangles consisting of two arcs, thus with four endpoints embedded on the boundary of a ball.
\end{definition}

The simplest tangle is a rational tangle.

\begin{definition}\label{Def:RationalTangle}
A \emph{rational tangle}\index{rational tangle}\index{tangle!rational tangle} is a tangle obtained by embedding two disjoint arcs on the surface of a 4-punctured sphere (pillowcase), and then pushing the interiors slightly into the $3$-ball bounded by the 4-punctured sphere. 
\end{definition}

The tangles are called rational because they can be defined by a rational number, as follows. Recall that a rational number can be described by a continued fraction:\index{continued fraction}
\[
\frac{p}{q} = [a_n, a_{n-1}, \hdots, a_1] =  a_n + \cfrac{1}{a_{n-1} + \cfrac{1}{\ddots \, + \cfrac{1}{a_1}}}.
\]

\begin{figure}[h]
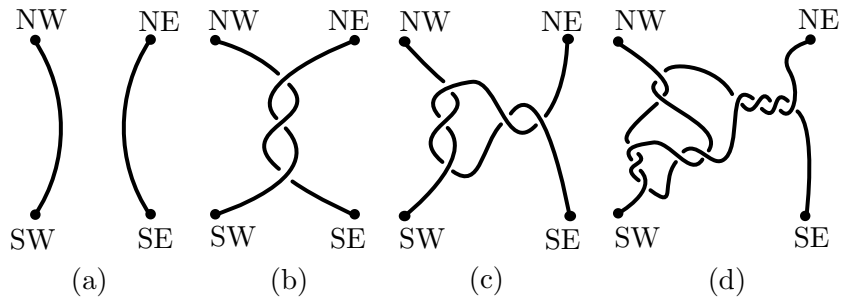
\caption{Building a rational tangle from the continued fraction $[4, -2, -2, 3]$.}
 \label{Fig:RatTangle1}
\end{figure}

Now, given a continued fraction $[a_n, \dots, a_1]$, we will form a rational tangle. Start by labeling the four points on the pillowcase NW, NE, SW, and SE.
If $n$ is even, connect NE to SE and NW to SW by attaching two arcs as in \reffig{RatTangle1}(a). 
Perform a homeomorphism of $B^3$ that rotates the points NW and NE $|a_1|$ times, creating a vertical band of $|a_1|$ crossings in the two arcs. If $a_1>0$, rotate in a counterclockwise direction, so that the overcrossings of the result have positive slope. This is called a \emph{positive crossing}.\index{positive crossing}\index{crossing!positive} If $a_1<0$ rotate in a clockwise direction, so that overcrossings have negative slope, forming a \emph{negative crossing}.\index{negative crossing}\index{crossing!negative} In \reffig{RatTangle1}(b), three positive crossings have been added. After twisting, relabel the points NW, NE, SW, and SE to match their original orientation. Next, apply a homeomorphism of $B^3$ that rotates NE and SE $|a_2|$ times, adding crossings in a horizontal band. Again these crossings will be positive if $a_2>0$, and negative if $a_2<0$. Repeat this process for each $a_i$. When finished, we obtain a rational tangle. An example is shown in \reffig{RatTangle1}. 

If $n$ is odd, start with two arcs connecting NW to NE and SW to SE. In this case we add a horizontal band of crossings first, and then continue as before, alternating between horizontal and vertical bands for each $a_i$.

Any rational tangle may be built by this process. As a convention, we require that the left-most term $a_n$ in the continued fraction expansion corresponds to a horizontal band of crossings. If we build a rational tangle ending with a vertical band, as in \reffig{RatTangle1}(b), then we insert a $0$ into the corresponding continued fraction, representing a horizontal band of $0$ crossings. For example, the continued fraction corresponding to the tangle in \reffig{RatTangle1}(b) is $[0, 3]$. This convention ensures that any continued fraction completely specifies a single rational tangle. There are two trivial rational tangles, namely $0 = [0]$, with untwisted strands connecting NW to NE and SW to SE, and $\infty = [0,0] = 0+\frac{1}{0}$, with untwisted strands connecting NW to SW and NE to SE. The tangle $\infty$ is shown in \reffig{RatTangle1}(a). 

\begin{proposition}[\cite{Conway}]\label{Prop:Conway}
Equivalence classes of rational tangles are in one-to-one correspondence with the set $\QQ \cup \infty$. In particular, tangles $T(a_n, \dots, a_1)$ and $T(b_m,\dots, b_1)$ are equivalent if and only if the continued fractions $[a_n, \dots, a_1]$ and $[b_m, \dots, b_1]$ are equal. \qed
\end{proposition}

\Refprop{Conway} allows us to put all our tangles into nice form. 

\begin{corollary}\label{Cor:Conway}
For any rational tangle, there exists an equivalent tangle $T(a_n, \dots, a_1)$ for which $a_i\neq 0$ for $1\leq i<n$, and either all $a_i\geq 0$ or all $a_i\leq 0$.
\end{corollary}

\begin{proof}
This follows immediately from \refex{ContinuedFractions}.
\end{proof}

Thus we will assume that positive tangles have only positive crossings, and negative tangles have only negative crossings. In either case, this will make the tangle diagram alternating.\index{alternating diagram!tangle}

\begin{definition}
The \emph{numerator closure}\index{numerator closure} ${\rm{num}}(T)$ of a rational tangle $T$ is formed by connecting NW to NE and SW to SE by simple arcs with no crossings. The \emph{denominator closure}\index{denominator closure} ${\rm{denom}}(T)$ is formed by connecting NW to SW and NE to SE by simple arcs with no crossings.
\end{definition}

\begin{definition}\label{Def:2BridgeKnot}
A \emph{2-bridge knot or link}\index{2-bridge knot or link}\index{2-bridge knot or link!definition} is the denominator closure of a rational tangle.
\end{definition}

Notice that the denominator closure of the tangle $T(a_n, a_{n-1},\dots, a_1)$ is always equivalent to the denominator closure of the tangle $T(0,a_{n-1}, \dots, a_1)$, since $a_n$ corresponds to horizontal crossings that can simply be unwound after forming the denominator closure. Thus when we consider 2-bridge knots, we may assume that in our rational tangle, $a_n=0$.

\begin{definition}\label{Def:2BridgeNotation}
  The 2-bridge knot or link\index{2-bridge knot or link} that is the denominator closure of the tangle $T(a_n, a_{n-1}, \dots, a_1)$ (and $T(0,a_{n-1}, \dots, a_1)$) is denoted by $K[a_{n-1}, \dots, a_1]$.
\end{definition}

The above discussion of twisting and taking denominator closure gives a nice correspondence between diagrams of 2-bridge knots and continued fraction expansions of rational numbers $p/q$ with $|p/q|\leq 1$. This is summarized in the following lemma. 

\begin{figure}[h]
  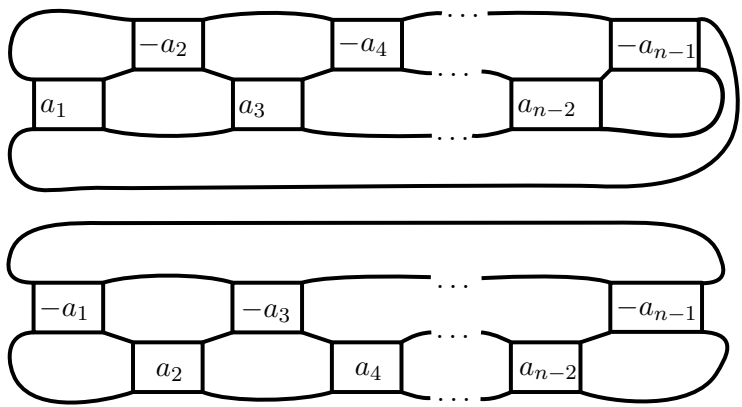
  \caption{The diagram of $K[a_{n-1}, \dots, a_1]$. Top: $n$ odd; bottom: $n$ even. Box labeled $\pm a_i$ denotes a (horizontal) twist region with $|a_i|$ crossings, with sign of the crossings equal to that of $\pm a_i$.}
  \label{Fig:2BridgeDiagram}
\end{figure}

\begin{lemma}\label{Lem:DiagramCF}
  Suppose $[0,a_{n-1},\dots, a_1]$ is a continued fraction with either $a_i>0$ for all $i$ or $a_i<0$ for all $i$. Then the diagram of $K[a_{n-1}, \dots, a_1]$ contains $n-1$ twist regions, arranged left to right. The twist region on the far left contains $|a_1|$ crossings, with sign opposite that of $a_1$ when $n$ is even, the next twist region to the right contains $|a_2|$ crossings, with sign opposite that of $a_2$ when $n$ is odd, and so on, with the $i$-th twist region from the left containing $|a_i|$ crossings, with sign the same as $a_i$ if $i$ and $n$ are both even or both odd, and opposite sign of $a_i$ if one of $i$ and $n$ is even and the other odd. Twist regions connect as illustrated in \reffig{2BridgeDiagram}.
\end{lemma}

\begin{proof}
  The tangle $T(0,a_{n-1}, \dots, a_1)$ is obtained by forming horizontal or vertical bands of $|a_i|$ crossings, for $i=1,\dots, n-1$. Thus the diagram of the denominator closure has $n-1$ twist regions with the numbers of crossings as claimed. To put it into the form of \reffig{2BridgeDiagram}, isotope the diagram by rotating vertical twist regions to be horizontal. Note that the rotation changes the sign of the crossing. Thus when $n$ is even, this rotates twist regions with odd index $i$ to be horizontal, and thus sign becomes opposite that of $a_i$, for even $n$ and odd $i$. When $n$ is odd, this rotates twist regions with even index to be horizontal, again with sign opposite that of $a_i$, for odd $n$ and even $i$. 
\end{proof}

\begin{lemma}\label{Lem:a1an}
  For a 2-bridge knot or link\index{2-bridge knot or link} $K[a_{n-1}, \dots, a_1]$, we may always assume $|a_1|\geq 2$ and $|a_{n-1}|\geq 2$.
\end{lemma}
\begin{proof}
  Exercise. One way to see this is to consider the form of a 2-bridge knot or link with $|a_1|=1$ or $|a_{n-1}|=1$, and show that the corresponding twist region can be subsumed into another twist region. 
\end{proof}

For the rest of this chapter, we will assume the conclusions of \refcor{Conway} and \reflem{a1an}, namely that if $K[a_{n-1}, \dots, a_1]$ is a 2-bridge knot or link, then either $a_i>0$ for all $i$ or $a_i<0$ for all $i$, and $|a_{n-1}|\geq 2$ and $|a_1|\geq 2$.

%%%%%%%%%%%%%%%%%%%%%%%%%%%%%%%%%%%%%%%%%%%%%%%%%%%%%%%%%%%%%%%%%
\section{Triangulations of 2-bridge links}

We now describe a way to triangulate 2-bridge link complements that was first observed by Sakuma and Weeks \cite{SakumaWeeks}. A description was also given by Futer in the appendix of \cite{GueritaudFuter:2bridge}; we base our exposition here off of the latter paper.

\begin{figure}[h]
  \includegraphics{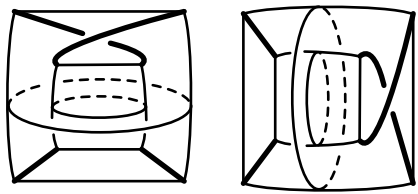}
  \caption{Vertical (left) and horizontal (right) blocks of the form $S\times I$. The 4-punctured spheres on the outside and inside correspond to $S\times\{1\}$ and $S\times\{0\}$, respectively}
  \label{Fig:4PunctSphereBlocks}
\end{figure}

Consider again our construction of a rational tangle. We started with two strands in a 4-punctured sphere, or pillowcase. To form each crossing, we either rotate the points NE and NW or the points NE and SE. In the former case, we call the crossing a \emph{vertical} crossing,\index{vertical crossing}\index{crossing!vertical} and in the latter a \emph{horizontal} crossing.\index{horizontal crossing}\index{crossing!horizontal} For all but the first crossing in a tangle, adding a crossing can be seen as stacking a region $S^2\times I$ to the outside of a pre-existing tangle, where $S^2\times I$ contains four strands, two of them forming a crossing. Positive vertical and horizontal crossings in $S^2\times I$ are shown in \reffig{4PunctSphereBlocks}, negative ones will be in the opposite direction. If we drill the four strands from $S^2\times I$, the region becomes $S\times I$, where $S$ is a 4-punctured sphere. 

\begin{lemma}\label{Lem:StackingSpheres}
  Let $K:=K[a_{n-1}, \dots, a_1]$ be a 2-bridge link\index{2-bridge knot or link} and let $C$ denote the number of crossings of $K$; so $C=|a_1|+\dots+|a_{n-1}|$. Assume either $a_i<0$ for all $i$ or $a_i>0$ for all $i$, and $|a_1|\geq 2$ and $|a_{n-1}|\geq 2$. Let $N$ be the manifold obtained from the complement $S^3-K$ by removing a ball neighborhood of the first and last crossings, and let $S$ denote the 4-punctured sphere. Then $N$ is homeomorphic to $S\times [a,b]$, obtained from stacking $C-2$ copies of $S\times I$ end to end, with each copy of $S\times I$ corresponding to either a horizontal or vertical crossing.
\begin{itemize}
  \item If $n$ is even, the first $a_1-1$ copies of $S\times I$ are vertical, followed by $a_2$ horizontal copies, $a_3$ vertical, etc, finishing with $a_{n-1}-1$ vertical copies of $S\times I$.
  \item If $n$ is odd, the first $a_1-1$ copies of $S\times I$ are horizontal, followed by $a_2$ vertical copies, $a_3$ horizontal, etc, finishing with $a_{n-1}-1$ vertical copies.
\end{itemize}
The $i$-th copy of $S\times I$ is glued along $S\times\{1\}$ to $S\times\{0\}$ on the $(i+1)$-st copy, $i=2, 3, \dots, C-1$. \qed
\end{lemma}

\begin{figure}[h]
  \includegraphics{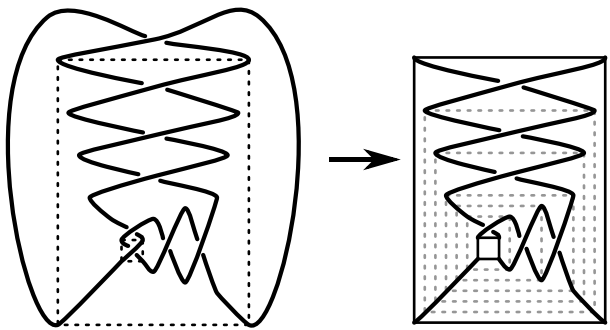}
  \caption{On the left is $K[4,2,2]$. On the right, remove neighborhoods of inside and outside crossings to obtain a manifold homeomorphic to $S\times[a,b]$. Each crossing is contained in a block $S\times I$ of the form of \reffig{4PunctSphereBlocks}}
  \label{Fig:TangleExample}
\end{figure}

An example of \reflem{StackingSpheres} is shown in \reffig{TangleExample}.

We will obtain a triangulation of a 2-bridge link complement by first finding a triangulation of the manifold $N$ in \reflem{StackingSpheres}. To do so, we will consider each of the blocks $S\times I$ separately, and then consider how they fit together.

Denote the blocks of \reflem{StackingSpheres} by $S_2\times I, S_3\times I, \dots, S_{C-1}\times I$, where $S_i\times I$ corresponds to the $i$-th crossing of the tangle. In the description below, we will consider the case that all crossings are positive, i.e.\ $a_j>0$ for all $j$, so that if the $i$-th crossing is vertical, then $S_i\times I$ has the form of the left of \reffig{4PunctSphereBlocks}, and if it is horizontal, then $S_i\times I$ has the form of the right. The case of all negative crossings will be similar. 

In $S_i\times I$, the 4-punctured sphere $S_i$ is embedded at any level $S_i\times\{t\}$. We will focus in particular on $S_i\times\{0\}$ and $S_i\times \{1\}$. 

\begin{lemma}\label{Lem:TriangulateSi}
There is an ideal triangulation of $S_i$ such that when we isotope the triangulation to $S_i\times\{1\}$, edges are horizontal (from NE to NW and from SE to SW), vertical (from SW to NW and from SE to NE), and diagonal, and when we isotope the triangulation to $S_i\times\{0\}$, edges are still horizontal, vertical, and diagonal, but the diagonals are opposite those of $S_i\times\{1\}$.
\end{lemma}

\begin{proof}
First consider the outside, $S_i\times\{1\}$. 
Draw vertical and horizontal ideal edges on $S_i\times\{1\}$; that is, draw horizontal edges from NE to NW and from SE to SW, and draw vertical edges from SE to NE and from SW to NW. Now isotope from $S_i\times\{1\}$ through $S_i\times\{t\}$ inside to $S_i\times\{0\}$, and track these ideal edges through the isotopy.

\begin{figure}[h]
  \includegraphics{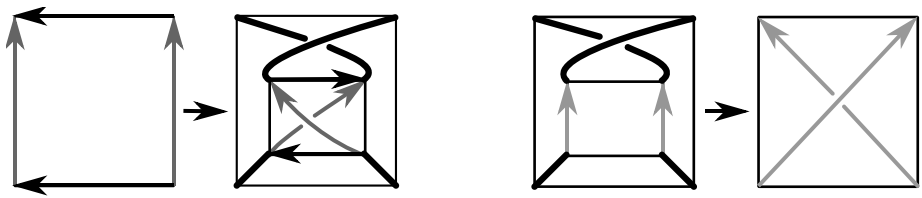}
  \caption{Effect on ideal edges of isotopies between $S_i\times\{1\}$ and $S_i\times\{0\}$.}
  \label{Fig:4PunctTriangles}
\end{figure}

In the case that $S_i\times I$ has a vertical crossing, as on the left of \reffig{4PunctSphereBlocks}, notice that the isotopy takes the horizontal edges in $S_i\times\{1\}$ to horizontal edges in $S_i\times\{0\}$, but it takes vertical edges in $S_i\times\{1\}$ to diagonal edges in $S_i\times\{0\}$, as shown on the left of \reffig{4PunctTriangles}. Now consider the vertical edges in $S_i\times\{0\}$, i.e.\ the ideal edges running from SE to NE and SW to NW on the inside of the block. When we isotope $S_i\times\{0\}$ to $S_i\times\{1\}$, notice that these edges become diagonal edges on $S_i\times\{1\}$, as shown on the right of \reffig{4PunctTriangles}. Notice that these diagonals are exactly opposite the diagonals on the inside on the left of the figure; see also \reffig{4PunctFullTriang}. 

When $S_i\times I$ has a horizontal crossing, as on the right of \reffig{4PunctSphereBlocks}, the vertical ideal edges in $S_i\times\{1\}$ are isotopic to vertical edges in $S_i\times\{0\}$. Horizontal edges on $S_i\times\{1\}$ isotope to diagonal edges on $S_i\times\{0\}$, and horizontal edges on $S_i\times\{0\}$ isotope to diagonal edges on $S_i\times\{1\}$. Again the diagonal edges have opposite slopes on the inside and outside. 

In either case, add all horizontal and vertical edges to $S_i$ on $S_i\times\{1\}$ and add all horizontal and vertical edges to $S_i$ on $S_i\times\{0\}$. Since either horizontal or vertical edges are duplicated, we add six ideal edges total. This is the triangulation claimed in the lemma.
\end{proof}

\begin{figure}[h]
%% Accompanies image file 'F10-06-4PFull.eps' (pdf, eps, ps)
\begingroup%
  \makeatletter%
  \providecommand\color[2][]{%
    \errmessage{(Inkscape) Color is used for the text in Inkscape, but the package 'color.sty' is not loaded}%
    \renewcommand\color[2][]{}%
  }%
  \providecommand\transparent[1]{%
    \errmessage{(Inkscape) Transparency is used (non-zero) for the text in Inkscape, but the package 'transparent.sty' is not loaded}%
    \renewcommand\transparent[1]{}%
  }%
  \providecommand\rotatebox[2]{#2}%
  \newcommand*\fsize{\dimexpr\f@size pt\relax}%
  \newcommand*\lineheight[1]{\fontsize{\fsize}{#1\fsize}\selectfont}%
  \ifx\svgwidth\undefined%
    \setlength{\unitlength}{345.66291046bp}%
    \ifx\svgscale\undefined%
      \relax%
    \else%
      \setlength{\unitlength}{\unitlength * \real{\svgscale}}%
    \fi%
  \else%
    \setlength{\unitlength}{\svgwidth}%
  \fi%
  \global\let\svgwidth\undefined%
  \global\let\svgscale\undefined%
  \makeatother%
  \begin{picture}(1,0.23293368)%
    \lineheight{1}%
    \setlength\tabcolsep{0pt}%
    \put(0,0){\includegraphics[width=\unitlength]{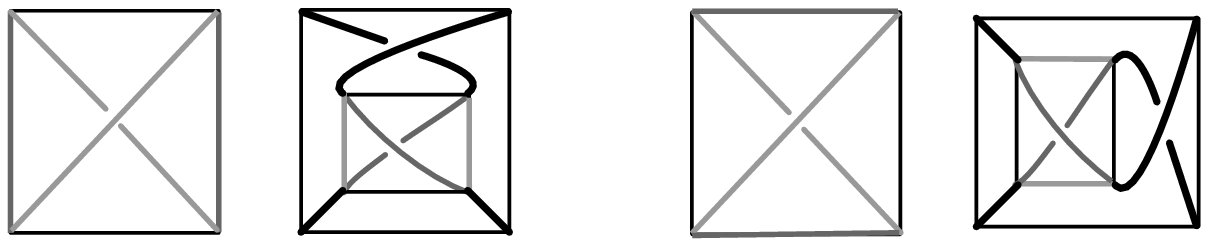}}%
    \put(0.04329919,0.00553736){\color[rgb]{0,0,0}\makebox(0,0)[lt]{\lineheight{0}\smash{\begin{tabular}[t]{l}$S\times\{1\}$\end{tabular}}}}%
    \put(0.29788247,0.00553736){\color[rgb]{0,0,0}\makebox(0,0)[lt]{\lineheight{0}\smash{\begin{tabular}[t]{l}$S\times\{0\}$\end{tabular}}}}%
    \put(0.59957604,0.00553736){\color[rgb]{0,0,0}\makebox(0,0)[lt]{\lineheight{0}\smash{\begin{tabular}[t]{l}$S\times\{1\}$\end{tabular}}}}%
    \put(0.86619534,0.0078752){\color[rgb]{0,0,0}\makebox(0,0)[lt]{\lineheight{0}\smash{\begin{tabular}[t]{l}$S\times\{0\}$\end{tabular}}}}%
  \end{picture}%
\endgroup%

  \caption{Triangulation of $S_i\times\{1\}$ and $S_i\times\{0\}$, shown for both horizontal and vertical (positive) crossings.}
  \label{Fig:4PunctFullTriang}
\end{figure}

The triangulation of $S_i\times\{1\}$ and $S_i\times\{0\}$ for both vertical and horizontal crossings is shown in \reffig{4PunctFullTriang}. (We have removed the arrows from the edges as we will not need to work with directed edges, and keeping track of direction will unnecessarily complicate the discussion.) Note the ideal edges cut $S_i\times\{1\}$ into four ideal triangles:\index{ideal triangle} two on the front and two on the back. Similarly, these ideal triangles can be isotoped to the inside $S_i\times\{0\}$, giving two triangles on the front and two on the back, although notice that the isotopy does not take both triangles in the front of $S_i\times\{1\}$ to triangles in the front of $S_i\times\{0\}$.

So far we only have ideal triangles on surfaces $S_i$, and no ideal tetrahedra. The ideal tetrahedra are obtained when we put blocks $S_{i-1}\times I$ and $S_i\times I$ together, and we now describe how this works.

\begin{lemma}\label{Lem:StackingTetrahedra}
  With $S_{i-1}$ and $S_i$ triangulated as in \reflem{TriangulateSi}, gluing $S_{i-1}\times I$ to $S_i\times I$ by identifying $S_{i-1}\times\{1\}$ and $S_i\times\{0\}$ gives rise to two ideal tetrahedra, each with two faces on $S_{i-1}$ and two on $S_i$. 
\end{lemma}

\begin{proof}
Consider the triangulations of $S_{i-1}\times\{1\}$ and $S_i\times\{0\}$, shown in \reffig{4PunctFullTriang}. Notice that the diagonal edges of $S_{i-1}\times\{1\}$ are exactly opposite the diagonal edges of $S_i\times\{0\}$, and so these edges do not match up. The horizontal and vertical edges on $S_{i-1}\times\{1\}$ and $S_i\times\{0\}$ can be identified, but the diagonal edges cannot. To keep these edges embedded, we view the diagonals of $S_{i-1}\times\{1\}$ as inside of the diagonals of $S_i\times\{0\}$, as shown in \reffig{Pillowcases}. 

\begin{figure}[h]
  \includegraphics{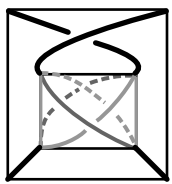}
  \caption{When blocks are glued, diagonals of the triangulated surfaces $S_i$ and $S_{i-1}$ are as shown.}
  \label{Fig:Pillowcases}
\end{figure}

With horizontal and vertical edges identified, notice that the interior of the region between $S_{i-1}\times\{1\}$ and $S_i\times\{0\}$ lies in two components: one on the front of the figure, and one on the back. Each of these components is bounded by four triangular faces, six ideal edges, and four ideal vertices; each is an ideal tetrahedron as desired.
\end{proof}

By \reflem{StackingTetrahedra}, when we glue $S_i\times I$ to $S_{i+1}\times I$, we obtain two additional tetrahedra, and these will be attached along $S_i$ to the two tetrahedra from $S_i\times I$ and $S_{i-1}\times I$. Thus as we run from $S_2\times I$ out to $S_{C-1}\times I$, we obtain pairs of tetrahedra for each gluing of blocks, and these are glued inside to outside to form a triangulation of $N\cong S\times[a,b]$.

Now, at this stage, we have a triangulation of $N$, but there will be four triangular faces on the very inside corresponding to $S_2$ that are unglued, and four triangular faces on the very outside corresponding to $S_{C-1}$ that are unglued. To complete the description of the triangulation of the 2-bridge link complement, we need to describe what happens at the outermost and innermost crossings, e.g.\ on the left of \reffig{TangleExample}.

\begin{proposition}\label{Prop:2BridgeTriang}
  Let $K:=K[a_{n-1}, \dots, a_1]$ be a 2-bridge link\index{2-bridge knot or link!triangulation}\index{2-bridge knot or link} with at least two twist regions, with either $a_i>0$ for all $1\leq i\leq n-1$, or $a_i < 0$ for all $i$. Assume $|a_1|\geq 2$ and $|a_{n-1}|\geq 2$. Let $C=|a_1|+\dots+|a_{n-1}|$ denote the number of crossings of $K$. 
  Then $S^3-K$ has a decomposition into $2(C-3)$ ideal tetrahedra denoted by $T_i^1, T_i^2$, for $i=2, \dots, C-2$.
  \begin{itemize}
  \item For $2\leq i \leq C-2$, the tetrahedra $T_i^1$ and $T_i^2$ each have two faces on $S_i$ and two on $S_{i+1}$.
  \item The two faces of $T_2^1$ on $S_2$ glue to the two faces of $T_2^2$ on $S_2$.
  \item Similarly, the two faces of $T_{C-2}^1$ on $S_{C-1}$ glue to the two faces of $T_{C-2}^2$ on $S_{C-1}$. 
  \end{itemize}
\end{proposition}

\begin{proof}
The tetrahedra come from the triangulation of $N$. By the previous lemmas, there are two tetrahedra for each pair of adjacent crossings, omitting the first and last, thus $2(C-3)$ tetrahedra. By \reflem{StackingTetrahedra}, each tetrahedron in each pair has two faces on $S_i$ and two on $S_{i+1}$, where $S_i\times I$ and $S_{i+1}\times I$ are blocks corresponding to the two adjacent crossings. We label the tetrahedra $T_2^1$, $T_2^2$, $T_3^1$, $T_3^2$, $\dots$, $T_{C-2}^1$, $T_{C-2}^2$, so that the first item is satisfied.

It remains to show that the innermost and outermost tetrahedra, $T_2^1$, $T_2^2$ and $T_{C-2}^1$, $T_{C-2}^2$, glue as claimed. 
We will focus on the outermost tetrahedra here. The innermost case is similar, but we leave its description to the reader.

For the outermost crossing, recall that we may assume our 2-bridge knot is the denominator closure of a tangle $T(a_n,a_{n-1}, \dots, a_1)$ with $a_n=0$. Thus the outermost crossing will be vertical, not horizontal. Hence we restrict to pictures with a single vertical crossing on the outside.

The outermost 4-punctured sphere $S_{C-1}$ will be triangulated as shown on the left of \reffig{OutsideTriangle}. Notice that when we add the outside crossing as shown, vertical edges and horizontal edges all become isotopic and hence are identified, by isotopies swinging the endpoints around the strand of the crossing. One such isotopy is indicated by the small arrow on the left of \reffig{OutsideTriangle}.

The diagonal edges are not identified to horizontal or vertical edges. When we follow the isotopy of \reffig{OutsideTriangle}, the diagonal edge in the front wraps once around a strand of the knot, as shown on the right of \reffig{OutsideTriangle}. Thus the triangle in the upper left corner of $S_{C-1}$ maps under the isotopy to a triangle with two of its edges identified, looping around a strand of the 2-bridge knot.

\begin{figure}[h]
  \includegraphics{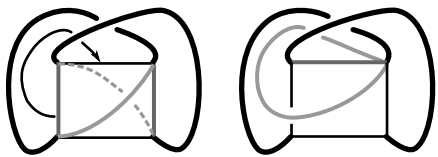}
  \caption{Identifying triangles of the outermost 4-punctured sphere}
  \label{Fig:OutsideTriangle}
\end{figure}

Now consider the triangle in front in the lower right corner. We will isotope the triangle by dragging its vertex on the SE corner around the strand of the knot to the NW corner. If we perform this isotopy while holding the diagonal fixed, note that the lower left triangle flips around backwards to be identified to the upper right triangle in the front. Thus the two triangles on the front of $S_{C-1}\times \{1\}$ will be identified under isotopy. Similarly for the two back triangles. 

Thus inserting the outermost crossing identifies the four outside triangular faces of the outermost tetrahedra in pairs. 

The tetrahedra $T_{C-2}^1$ and $T_{C-2}^2$ have triangular faces on $S_{C-1}$, shown in \reffig{Pillowcases}. One of these, say $T_{C-2}^1$, will lie in front in that figure and one will lie in back. 

However, note that in \reffig{Pillowcases}, we have isotoped the surface $S_{C-1}$ to be in the position of $S_{C-1}\times\{0\}$, while in \reffig{OutsideTriangle}, when we glue faces of $S_{C-1}$, we have isotoped $S_{C-1}$ to be in the position of $S_{C-1}\times\{1\}$. Isotoping from $S_{C-1}\times\{0\}$ to $S_{C-1}\times\{1\}$ will move the faces of $T_{C-2}^1$ and $T_{C-2}^2$. In particular, the face of $T_{C-2}^1$ lying on the upper right of \reffig{Pillowcases} will be moved by isotopy to lie in the back on the upper right, and the face of $T_{C-2}^1$ lying in the lower left will be moved by isotopy to lie in front, in the lower right. See \reffig{TetrFaces}.

\begin{figure}[h]
%% Accompanies image file 'F10-09-TetFac.eps' (pdf, eps, ps)
\begingroup%
  \makeatletter%
  \providecommand\color[2][]{%
    \errmessage{(Inkscape) Color is used for the text in Inkscape, but the package 'color.sty' is not loaded}%
    \renewcommand\color[2][]{}%
  }%
  \providecommand\transparent[1]{%
    \errmessage{(Inkscape) Transparency is used (non-zero) for the text in Inkscape, but the package 'transparent.sty' is not loaded}%
    \renewcommand\transparent[1]{}%
  }%
  \providecommand\rotatebox[2]{#2}%
  \newcommand*\fsize{\dimexpr\f@size pt\relax}%
  \newcommand*\lineheight[1]{\fontsize{\fsize}{#1\fsize}\selectfont}%
  \ifx\svgwidth\undefined%
    \setlength{\unitlength}{146.00548553bp}%
    \ifx\svgscale\undefined%
      \relax%
    \else%
      \setlength{\unitlength}{\unitlength * \real{\svgscale}}%
    \fi%
  \else%
    \setlength{\unitlength}{\svgwidth}%
  \fi%
  \global\let\svgwidth\undefined%
  \global\let\svgscale\undefined%
  \makeatother%
  \begin{picture}(1,0.5374865)%
    \lineheight{1}%
    \setlength\tabcolsep{0pt}%
    \put(0,0){\includegraphics[width=\unitlength]{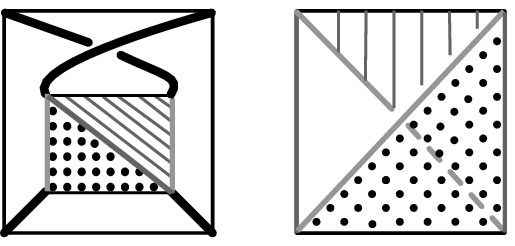}}%
    \put(0.66686454,0.01449338){\color[rgb]{0,0,0}\makebox(0,0)[lt]{\lineheight{0}\smash{\begin{tabular}[t]{l}$S\times\{1\}$\end{tabular}}}}%
    \put(0.11971145,0.01310953){\color[rgb]{0,0,0}\makebox(0,0)[lt]{\lineheight{0}\smash{\begin{tabular}[t]{l}$S\times\{0\}$\end{tabular}}}}%
  \end{picture}%
\endgroup%

  \caption{Locations of faces of $T_{C-2}^1$ under isotopy from $S\times\{0\}$ to $S\times\{1\}$}
  \label{Fig:TetrFaces}
\end{figure}

Thus the identification of triangles on the outside identifies faces of $T_{C-2}^1$ to faces of $T_{C-2}^2$.
\end{proof}

\subsection{The cusp triangulation}

Now we consider the view of the tetrahedra from a cusp. Consider first the manifold $N$ with ball neighborhoods of the first and last crossings removed. The manifold $N$ is homeomorphic to the product of a 4-punctured sphere and a closed interval. Note that in $N$, there are four distinct cusps, corresponding to the product of $I$ and the four distinct punctures of the 4-punctured sphere. Note that each cusp meets each 4-punctured sphere $S_i$, and that a curve on $S_i$ running around the puncture forms a meridian.
Finally, note that between each $S_i$ and $S_{i+1}$ lie two tetrahedra, as in \reffig{Pillowcases}. Each tetrahedron has exactly one ideal vertex on each of the four cusps. Thus the cusp triangulation of $N$ consists of four disjoint cusp neighborhoods. Each cusp neighborhood meets each $S_i$, in the same order, and each cusp neighborhood meets each tetrahedron in the decomposition in the same order. Thus the four cusp triangulations look identical at this stage, at least combinatorially. We will create one of these four cusps. 

In order to see the pattern of tetrahedra in one of these cusps, note that there will be a stack of triangles in each cusp, each triangle corresponding to the tip of a tetrahedron. By \refprop{2BridgeTriang}, the triangles will be sandwiched between 4-punctured spheres $S_i$ and $S_{i+1}$, with the bottom of the stack of triangles bounded by $S_2$ and the top by $S_{C-1}$. (Recall that the 4-punctured sphere $S_i$ actually lies in a block $S_i\times I$, so when we refer to $S_i$ in the following it may be helpful to recall that we are referring to a surface isotopic to $S_i\times\{t\}$ for appropriate $t$.)

\begin{figure}[h]
%% Accompanies image file 'F10-10-CuSphe.eps' (pdf, eps, ps)
\begingroup%
  \makeatletter%
  \providecommand\color[2][]{%
    \errmessage{(Inkscape) Color is used for the text in Inkscape, but the package 'color.sty' is not loaded}%
    \renewcommand\color[2][]{}%
  }%
  \providecommand\transparent[1]{%
    \errmessage{(Inkscape) Transparency is used (non-zero) for the text in Inkscape, but the package 'transparent.sty' is not loaded}%
    \renewcommand\transparent[1]{}%
  }%
  \providecommand\rotatebox[2]{#2}%
  \newcommand*\fsize{\dimexpr\f@size pt\relax}%
  \newcommand*\lineheight[1]{\fontsize{\fsize}{#1\fsize}\selectfont}%
  \ifx\svgwidth\undefined%
    \setlength{\unitlength}{229.87997246bp}%
    \ifx\svgscale\undefined%
      \relax%
    \else%
      \setlength{\unitlength}{\unitlength * \real{\svgscale}}%
    \fi%
  \else%
    \setlength{\unitlength}{\svgwidth}%
  \fi%
  \global\let\svgwidth\undefined%
  \global\let\svgscale\undefined%
  \makeatother%
  \begin{picture}(1,0.2862958)%
    \lineheight{1}%
    \setlength\tabcolsep{0pt}%
    \put(0,0){\includegraphics[width=\unitlength]{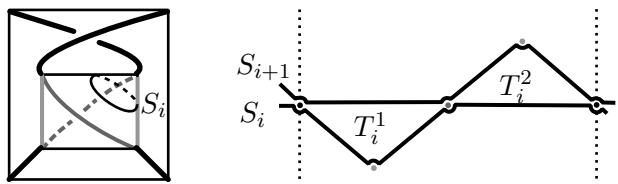}}%
    \put(0.38240439,0.09990473){\color[rgb]{0,0,0}\makebox(0,0)[lt]{\lineheight{0}\smash{\begin{tabular}[t]{l}$S_i$\end{tabular}}}}%
    \put(0.22031917,0.1163563){\color[rgb]{0,0,0}\makebox(0,0)[lt]{\lineheight{0}\smash{\begin{tabular}[t]{l}$S_i$\end{tabular}}}}%
    \put(0.57048403,0.07705908){\color[rgb]{0,0,0}\makebox(0,0)[lt]{\lineheight{0}\smash{\begin{tabular}[t]{l}$T_i^1$\end{tabular}}}}%
    \put(0.80911787,0.14541762){\color[rgb]{0,0,0}\makebox(0,0)[lt]{\lineheight{0}\smash{\begin{tabular}[t]{l}$T_i^2$\end{tabular}}}}%
    \put(0.37743284,0.18007076){\color[rgb]{0,0,0}\makebox(0,0)[lt]{\lineheight{0}\smash{\begin{tabular}[t]{l}$S_{i+1}$\end{tabular}}}}%
  \end{picture}%
\endgroup%

  \caption{Form of two meridians running over $S_i$ and $S_{i+1}$}
  \label{Fig:CuspSpheres}
\end{figure}

When we run along a meridian of the cusp on $S_i$, we stay on edges of the cusp triangulation. Moreover, note that we pass over exactly three ideal edges; see the left of \reffig{CuspSpheres}. Thus in the cusp triangulation, running along such a meridian on the surface $S_i$ will correspond to running over three edges of triangles. This is shown in \reffig{CuspSpheres} for $S_i$. The vertical dotted lines indicate boundaries of a fundamental region for the cusp torus; in this case running from one dotted line to the other corresponds to a meridian.

When the $i$-th, $(i+1)$-st, and $(i+2)$-nd blocks are all vertical crossings, note that the surfaces $S_i$, $S_{i+1}$, and $S_{i+2}$ will all share an edge; a horizontal edge in \reffig{4PunctFullTriang}. Similarly adjacent horizontal crossings also lead to surfaces sharing an edge.

\begin{notation}\label{Not:RLNotation}
  In the cusp triangulation we see 4-punctured spheres $S_2, \dots, S_{C-1}$. Give the label $R$ to each 4-punctured sphere $S_i$ corresponding to a block $S_i\times I$ containing a horizontal crossing. Give the label $L$ to each corresponding to a block containing a vertical crossing. A tetrahedron that lies between layers labeled $R$ and $L$ is called a \emph{hinge tetrahedron}\index{hinge tetrahedron}. 
\end{notation}

The labels $R$ and $L$ are given for historical reasons; they refer to moves to the right and left for a path in the Farey graph given by the rational number of our tangle. We won't delve into the history of this notation here, but we will use this notation for ease of reference with other literature. For more information see \cite{GueritaudFuter:2bridge}.

\begin{example}
A cusp triangulation for an example $N$ is shown in \reffig{CuspExample}. That figure follows some standard conventions. Because we have many surfaces $S_i$, we connect edges of $S_i$ to form a single connected jagged line, identifiable as one surface in the cusp, and put a little space between multiple such surfaces at a vertex they share. We also put vertices of the triangles in two columns (in a fundamental domain). Finally, we shade the hinge layers. 

\begin{figure}[h]
  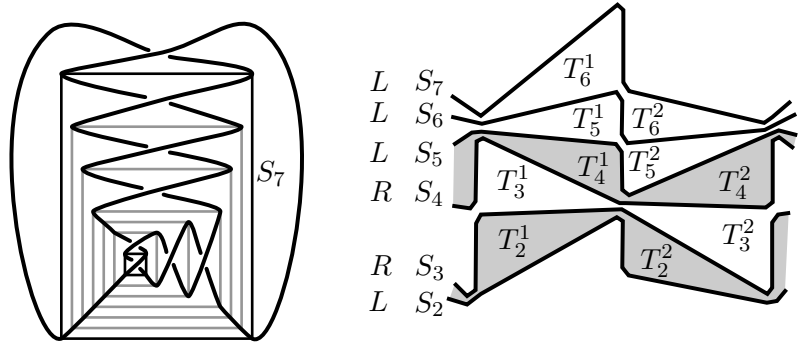
  \caption{On the right is shown the cusp triangulation of one of the four cusps of $N$ on the left}
  \label{Fig:CuspExample}
\end{figure}

Note in the figure as we move inside to out, we move from the bottom of the cusp triangulation to the top. Tetrahedron $T_2^1$ lies between surfaces $S_2$ associated with the second innermost crossing (vertical, $L$) and $S_3$ associated with the third innermost crossing (horizontal, $R$). It is a hinge tetrahedron. Tetrahedron $T_3^1$ lies between $S_3$ and $S_4$, both of which are associated with horizontal crossings, $R$. Note there is an ideal edge shared by all three surfaces $S_2$, $S_3$, and $S_4$, and this corresponds to a shared vertex of $T_2^1$ and $T_3^1$ in the cusp triangulation (in the center). 
\end{example}

\begin{figure}[h]
  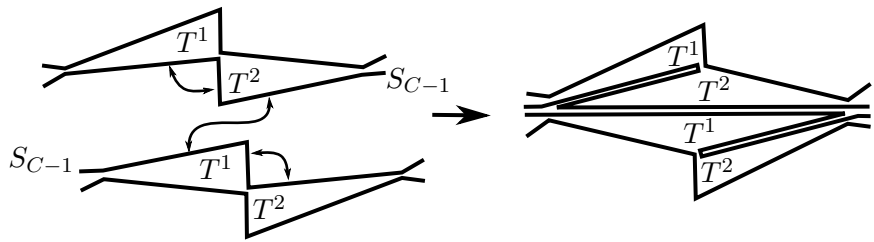
  \caption{Gluing tetrahedra across $S_{C-1}$ yields a hairpin turn}
  \label{Fig:HairpinTurn}
\end{figure}

Now we determine what happens to cusps when we put in the innermost and outermost crossings. At the outermost crossing, note that the cusp corresponding to the vertex SE becomes identified with the cusp corresponding to the vertex NW, and similarly for SW and NE. Thus the four identical cusp triangulations we have obtained so far will be glued. Recall that the gluing is along triangle faces of $S_{C-1}$ in the case of the outermost crossing. The faces of $T_{C-2}^1$ are glued to faces of $T_{C-2}^2$. The result is a ``folding'' of triangles. See \reffig{HairpinTurn}. We call this a \emph{hairpin turn}\index{hairpin turn}. 

If $K$ is a knot, if we follow a longitude of the cusp, starting at one of the corners of $S_2$, we will see 4-punctured spheres $S_2, S_3, \dots, S_{C-2}$, then a hairpin turn on $S_{C-1}$ corresponding to the outside crossing as in \reffig{HairpinTurn}. Continuing, we will pass $S_{C-2}, S_{C-3}, \dots, S_3$, then another hairpin turn on $S_2$ corresponding to the inside crossing, then $S_3, \dots, S_{C-1}$ and a hairpin turn, and finally $S_{C-2}, \dots, S_3$ and the original $S_2$ with a hairpin turn. A hairpin turn appears in the cusp triangulation as a single edge stretching across a meridian, adjacent to two triangles whose third vertex is 3-valent.

We summarize:

\begin{proposition}\label{Prop:2BridgeCuspTriang}
  Let $K:=K[a_{n-1}, \dots, a_1]$ be a 2-bridge knot\index{2-bridge knot or link}\index{2-bridge knot or link!cusp triangulation} with at least two twist regions, such that either $a_i>0$ for all $i$, or $a_i<0$ for all $i$, and $|a_1|\geq 2$ and $|a_{n-1}|\geq 2$. Let $C=|a_1|+\dots+|a_{n-1}|$ denote the number of crossings of $K$. The cusp triangulation of $K$ has the following properties.
  \begin{itemize}
  \item It is made of four pieces, each piece bookended by hairpin turns corresponding to 4-punctured spheres $S_2$ and $S_{C-1}$. Between lies a sequence of 4-punctured spheres $S_3, \dots, S_{C-2}$. We call the 4-punctured spheres \emph{zig-zags}.\index{zig-zag}
  \item The first and third pieces are identical; the second and fourth are also identical and given by rotating the first piece $180^\circ$ about a point in the center of the edge of the final hairpin turn (and swapping some labels $T_i^1$ as $T_i^2$). Thus the second and fourth pieces follow the first in reverse. 
  \item When running in a longitudinal direction, the first piece begins with $|a_{n-1}|-1$ zig-zags labeled $L$; the first of these is $S_{C-1}$, the hairpin turn corresponding to the outside crossing of the knot. These zig-zags are followed by $|a_{n-2}|$ zig-zags labeled $R$, then $|a_{n-3}|$ labeled $L$, and so on. If $n$ is even, finish with $|a_1|-1$ zig-zags labeled $L$, the last of which is the final hairpin turn, corresponding to $S_2$ at the inside crossing. If $n$ is odd, the final $|a_1|-1$ zig-zags are labeled $R$.
  \item A meridian follows a single segment of the zig-zag in a hairpin turn, or three segments of any other zig-zag. 
  \end{itemize}
Note we see each $S_i$ exactly four times, including seeing $S_2$ twice for each of the two hairpin turns in the cusp triangulation corresponding to the inside crossing, and seeing $S_{C-1}$ twice for each hairpin turn corresponding to the outside crossing.
\end{proposition}

\begin{figure}[h]
  \includegraphics{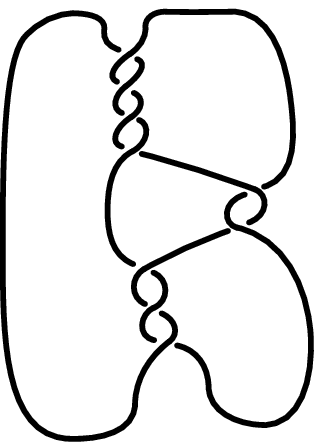}
  \hspace{.1in}
  \includegraphics{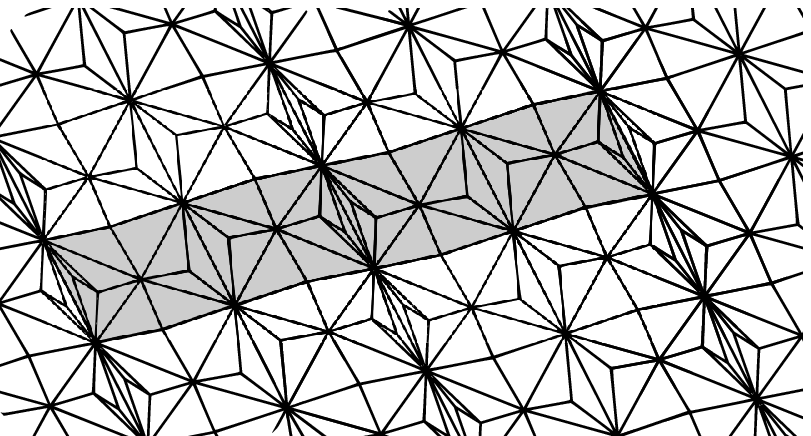}
  \caption{An example of a 2-bridge knot\index{2-bridge knot or link!cusp triangulation} and its cusp triangulation from SnapPy. The shaded region shows a fundamental domain for the cusp torus, stretching from one hairpin turn through three others back to the same hairpin turn.}
  \label{Fig:SnapPyExample}
\end{figure}

An example sketched by SnapPy (\cite{SnapPy}) is shown in \reffig{SnapPyExample}.

%%%%%%%%%%%%%%%%%%%%%%%%%%%%%%%%%%%%%%%%%%%%%%%%%%%%%%%%%%%%%%%%%
\section{Positively oriented tetrahedra}

The triangulation described in the last section has nice geometry. In particular, when the 2-bridge link\index{2-bridge knot or link} has at least two twist regions, we can find angle structures\index{angle structure} on the triangulation. These can be used to prove that the 2-bridge link is hyperbolic (\refcor{2BridgeHyperbolic}), and to show that in the complete hyperbolic structure on the link complement, the tetrahedra are all geometric. Thus we will obtain our first infinite class of knots and links with known geometric triangulations.\index{geometric triangulation}

The main theorem of the next two sections is \refthm{2BridgeGeometric}, below. It was originally proved by Futer in the appendix to \cite{GueritaudFuter:2bridge}.

\begin{theorem}\label{Thm:2BridgeGeometric}
Let $K$ be a 2-bridge knot or link\index{2-bridge knot or link}\index{2-bridge knot or link!triangulation} with a reduced alternating diagram\index{alternating knot or link}\index{alternating diagram} with at least two twist regions. Let $\mathcal{T}$ be the triangulation of $S^3-K$ as described above.
Then $S^3-K$ is hyperbolic, and in the complete hyperbolic structure on $S^3-K$, all tetrahedra of $\mathcal{T}$ are positively oriented.\index{positively oriented tetrahedron}\index{tetrahedron!positively oriented}
\end{theorem}

The proof of \refthm{2BridgeGeometric} uses angle structures\index{angle structure} on $\mathcal{T}$, as in \refdef{AngleStructures}, and is done in two steps. First, the space of angle structures\index{angle structure} $\mathcal{A}(\mathcal{T})$ is shown to be non-empty. In \refthm{HypAngleStruct}, we showed that the existence of an angle structure\index{angle structure} is enough to conclude that the manifold admits a hyperbolic structure. We conclude that these 2-bridge link complements are hyperbolic. 

Second, in the next section, we show that the volume functional\index{volume functional} cannot achieve its maximum on the boundary of $\mathcal{A}(\mathcal{T})$. This is all that is needed: because the volume functional is strictly concave down on $\mathcal{A}(\mathcal{T})$ (\reflem{DerivativesVol}), it achieves a maximum in the interior of $\mathcal{A}(\mathcal{T})$. By \refthm{VolAngleStructs}, the maximum corresponds to the complete hyperbolic structure, and at that structure, all angles are strictly positive, meaning all tetrahedra are geometric --- positively oriented.\index{positively oriented tetrahedron}\index{tetrahedron!positively oriented}

\begin{proposition}\label{Prop:2BridgeAngleNonempty}
Let $\mathcal{T}$ be the triangulation of a 2-bridge knot or link\index{2-bridge knot or link!angle structure} complement with at least two twist regions, as described above. Then the space of angle structures\index{angle structure} $\mathcal{A}(\mathcal{T})$ is nonempty.
\end{proposition}

\begin{figure}[h]
  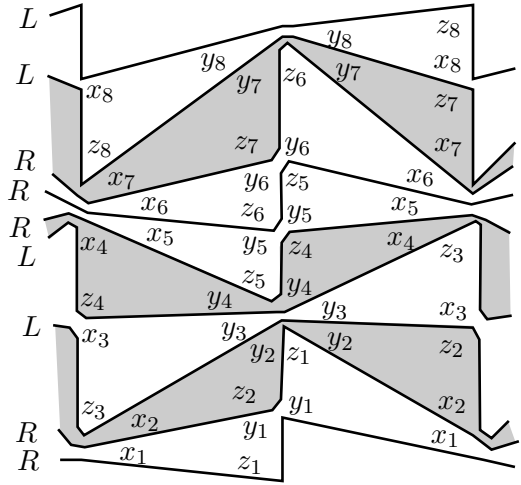
  \caption{Labels on the cusp triangulation of $N = S\times [a,b]$ for an example}
  \label{Fig:CuspLabels}
\end{figure}

The proof of the proposition is not hard, but requires additional notation. First, we need to label the angles of each of the tetrahedra constructed in the previous section. Remember that the tetrahedra were constructed in pairs, and the pairs of tetrahedra lie between two 4-punctured spheres of the manifold $N=S\times[a,b]$, as in \reffig{Pillowcases}. In order to show some angle structure\index{angle structure} exists, we will first assume that the angles on each of these pairs of tetrahedra agree. Let $z_i$ denote the angle on the outside diagonal edges of tetrahedra $T_i^1$ and $T_i^2$. Because opposite edges have the same angle, $z_i$ is also the angle on the inside diagonal edge. Denote the angle at the horizontal edges by $x_i$ and the angle at vertical edges by $y_i$. We may add these angles to the cusp triangulation. The cusp triangulation was obtained by adding layers of zig-zagging 4-punctured spheres. Each 4-punctured sphere shares two edges with the previous 4-punctured sphere, and has one new edge. In the cusp triangulation, this forms a sequence of triangles in which two vertices are shared, but one new vertex is added. The new vertex corresponds to the diagonal edge, so is labeled $z_i$. Note that angles labeled $x_i$ are glued together, as are angles labeled $y_i$. Finally, we have oriented tetrahedra so that angles read $x_i$, $y_i$, $z_i$ in clockwise order around a cusp triangle. This completely determines the labeling on all the cusp triangles of $N=S\times[a,b]$. An example is shown in \reffig{CuspLabels}. 

In the example, a 4-punctured sphere corresponding to a horizontal crossing is labeled $R$, and one corresponding to a vertical crossing is labeled $L$, as in \refnot{RLNotation}.

Now $(x_i,y_i,z_i)$ give us labels for angles of all the tetrahedra on the 2-bridge link complement. We will need $x_i+y_i+z_i=\pi$ for each $i$ to satisfy condition \refitm{VertexSum} of the definition of an angle structure,\index{angle structure} \refdef{AngleStructures}. We also need sums of angles around edge classes to be $2\pi$.

Away from hairpin turns, the edge gluings of $S^3-K$ agree with those of $N=S\times[a,b]$, so we will first consider angle sums around edges of $N$, simplifying these conditions to a system of equations in terms of the $z_i$ alone, then deal with hairpin turns later. 
There will be four cases depending on whether the $i$-th tetrahedron lies between two horizontal crossings, two vertical crossings, a horizontal followed by a vertical crossing, or a vertical followed by a horizontal crossing. These cases are denoted by $RR$, $LL$, $RL$, and $LR$, respectively. 

\begin{figure}[h]
  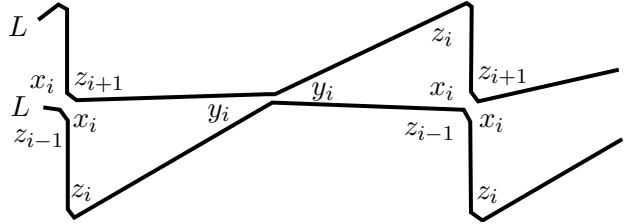
  \caption{Labels in the $LL$ case}
  \label{Fig:LLCusp}
\end{figure}

The labels for two consecutive $LL$ 4-punctured spheres are shown in \reffig{LLCusp}. Note in this case, there is a 4-valent vertex in the cusp triangulation (or 4-valent ideal edge in the decomposition into tetrahedra). In order for the angle sum around this edge to be $2\pi$, we need $2x_i + z_{i+1}+z_{i-1}=2\pi$, or $x_i = \half(2\pi-z_{i-1} - z_{i+1})$. Then in order for $x_i+y_i+z_i=\pi$, we need $y_i= \half(z_{i-1} - 2z_i + z_{i+1})$. 

A similar picture occurs in the $RR$ case. Again there is a 4-valent vertex, and reading the labels around that vertex we find that we need the formulas $y_i = \half(2\pi-z_{i-1}-z_{i+1})$, and $x_i = \half(z_{i-1}-2z_i+z_{i+1})$.

\begin{figure}[h]
  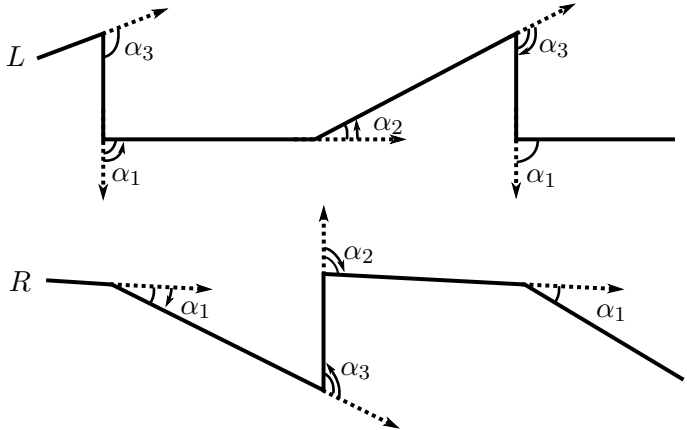
  \caption{Pleating angles for 4-punctured spheres}
  \label{Fig:PleatingAngles}
\end{figure}

In the $LR$ and $RL$ cases, there is not a single edge all of whose labels we can read off the diagram. In these cases, we find restrictions by considering \emph{pleating angles}\index{pleating angle}. Pleating angles $\alpha_1$, $\alpha_2$, and $\alpha_3$ are the angles determining the bending of the pleated\index{pleated surface} 4-punctured sphere. They are shown for 4-punctured spheres labeled $L$ and $R$ in \reffig{PleatingAngles}.

\begin{lemma}\label{Lem:PleatingAngles}
  If the angle structure\index{angle structure} gives a Euclidean structure\index{Euclidean structure} on the cusp, then it will be the case that pleating angles as in \reffig{PleatingAngles} satisfy $\alpha_1+\alpha_2-\alpha_3=0$.
\end{lemma}

\begin{proof}
  Exercise.
\end{proof}

To find an angle structure,\index{angle structure} we will assume this pleating condition holds in the $LR$ and $RL$ case.

\begin{figure}[h]
  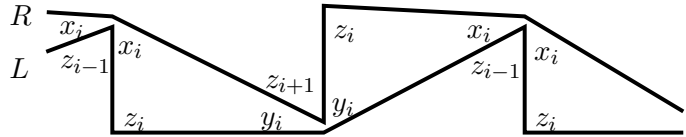
  \caption{Labels in the $LR$ case}
  \label{Fig:LRCusp}
\end{figure}

The $LR$ labels are shown in \reffig{LRCusp}. Note that the pleating angles for the 4-punctured sphere at the bottom of the diagram are $\alpha_1=\pi-z_i$, $\alpha_2 = \pi-(2y_i+z_{i+1})$, and $\alpha_3=\pi-z_{i-1}$. Thus the condition $\alpha_1+\alpha_2-\alpha_3=0$ implies $y_i=\half(\pi+z_{i-1}-z_i-z_{i+1})$. The pleating angles on the 4-punctured sphere on the top of the diagram in \reffig{LRCusp} are $\alpha_1=\pi-(2x_i+z_{i-1})$, $\alpha_2=\pi-z_i$, and $\alpha_3=\pi-z_{i+1}$. Thus the pleating condition for this 4-punctured sphere gives $x_i=\half(\pi-z_{i-1}-z_i+z_{i+1})$. Conditions can be obtained in a similar manner in the $RL$ case. 

In summary, away from hairpin turns, labels must satisfy the conditions given in \reftable{LabelConditions}.

\begin{table}
\begin{tabular}{ccc}
  & $LL$ & $RR$ \\
  \hline\noalign{\smallskip}
  $x_i$ & ${\half(2\pi-z_{i-1} - z_{i+1})}$ & ${\half(z_{i-1}-2z_i+z_{i+1})}$  \vspace{.1in} \\
  $y_i$ & ${\half(z_{i-1} - 2z_i + z_{i+1})}$ & ${\half(2\pi-z_{i-1}-z_{i+1})}$ \\
  $z_i$ & $z_i$ & $z_i$  \\

  \multicolumn{3}{ c } {\vspace{.1in}} \\

  & $LR$ & $RL$ \\
  \hline\noalign{\smallskip}
  $x_i$ & ${\half(\pi-z_{i-1}-z_i+z_{i+1})}$ & ${\half(\pi+z_{i-1}-z_i-z_{i+1})}$ \vspace{.1in} \\
  $y_i$ & ${\half(\pi+z_{i-1}-z_i-z_{i+1})}$ & ${\half(\pi-z_{i-1}-z_i+z_{i+1})}$ \\
  $z_i$ & $z_i$ & $z_i$ \smallskip
\end{tabular}
\caption{Label conditions in terms of the $z_i$}
\label{Table:LabelConditions}
\end{table}

Notice this allows us to express $x_i$ and $y_i$ in terms of $z_{i-1}, z_i,$ and $z_{i+1}$ alone. Note also that the sum of the angles $x_i+y_i+z_i=\pi$ in each case.

Finally, we claim that with the conditions in \reftable{LabelConditions}, the angle sum around each edge in $N$ is $2\pi$. To see this, note first that we have constructed the angles so that the sum is $2\pi$ around 4-valent edges. We now check the remaining edges. The angle sum around one such edge will be
\[ z_{j-1} + 2x_j + \sum_{i=j+1}^{k-1} 2x_i + 2x_k + z_{k+1},\]
where $j$ and $k$ are indices of hinge tetrahedra, with $j$ between $LR$ and $k$ between $RL$, and $j<k$, and all 4-punctured spheres labeled $R$ between them; refer to \reffig{CuspLabels}. 
By the formulas in the tables, this is
\[ z_{j-1} + \pi-z_{j-1}-z_j+z_{j+1} + \sum_{i=j+1}^{k-1} (z_{i-1}-2z_i+z_{i+1}) + \pi+z_{k-1}-z_k-z_{k+1} + z_{k+1}. \]
This is a telescoping sum; all terms cancel except $2\pi$, as desired.

The angle sum around another such edge will be
\[ z_{j-1} + 2y_j + \sum_{i=j+1}^{k-1} 2y_i + 2y_k + z_{k+1},\]
where $j$ and $k$ are hinge indices, with $j$ between $RL$ and $k$ between $LR$, and $j<k$, and all 4-punctured spheres are labeled $L$ between them. Again check that everything cancels except $2\pi$.

\begin{figure}
  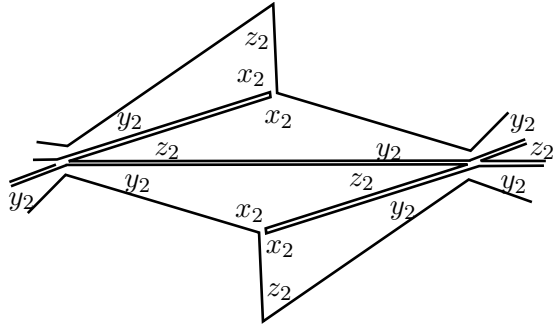
  \caption{Labels from a hairpin turn}
  \label{Fig:HairpinLabels}
\end{figure}

We still need to consider the hairpin turns. With the gluing that comes from a hairpin turn, labels are as shown in \reffig{HairpinLabels}, for the $LL$ case. The cases $RR$, $LR$, $RL$ are similar (exercise).

If we set $z_1=0$, the interior angle in which the 4-punctured sphere $S_2$ is bent at the hairpin turn, then all the equations in \reftable{LabelConditions} hold, depending on whether the hairpin turn occurs in the case $LL$, $RR$, $LR$, or $RL$. It remains only to check the edge equations. For the edge at the sharp bend, the equation will be identical to one of the previous equations, only now with angle $z_1=0$ included. The sum is still $2\pi$. As for the final edges, in the case $S_2$ is $R$, these contribute $2z_2 + 4x_2 + \dots$, where the remainder of terms depends on whether the hairpin turn occurs at a hinge or not. In either case, the sum is $2\pi$. Similarly when $S_2$ is $L$, and similarly for the outside hairpin turn that occurs at the 4-punctured sphere $S_{C-1}$.

We are now ready to show the space of angle structures is nonempty.

\begin{proof}[Proof of \refprop{2BridgeAngleNonempty}]
We show the space of angle structures\index{angle structure} is nonempty by showing there is a choice of $(z_1, z_2, \dots, z_{C-2}, z_{C-1})$ with $z_1=z_{C-1}=0$, all other $z_i\in (0,\pi)$, and $x_i, y_i \in (0,\pi)$. For this to hold, the equations in \reftable{LabelConditions} tell us that:
\begin{equation}\label{Eqn:2BridgeZConditions}
  \begin{cases}
    2z_i < z_{i-1}+z_{i+1} & \mbox{ if } i \mbox{ is not a hinge ($LL$ or $RR$)} \\
    |z_{i+1}-z_{i-1}| < \pi-z_i & \mbox{ if } i \mbox{ is a hinge index ($LR$ or $RL$)} \\
  \end{cases}
\end{equation}

The first equation is called the \emph{convexity equation}.\index{convexity equation} The second is the \emph{hinge equation}.\index{hinge equation}

We find a point with all $z_i\in(0,\pi)$ that satisfies convexity and hinge equations. Namely, let $z_1=z_{C-1}=0$. For each hinge index $i$, let $z_i=\pi/3$. Between hinge indices, choose a sequence to satisfy the convexity equations. For example, if $j,k$ are consecutive hinge indices with $j<k$, then for all $j\leq i\leq k$, take
\[ z_i = \frac{\pi}{3} - \frac{2(i-j)(k-i)}{(k-j)^2}. \]
Then the sequence $(z_1, z_2, \dots, z_{C-2}, z_{C-1})$ satisfies all required conditions. Letting $x_i$ and $y_j$ be as in the tables, this gives an angle structure. 
\end{proof}

\begin{corollary}\label{Cor:2BridgeHyperbolic}
  Let $K[a_{n-1}, \dots, a_1]$ be a 2-bridge knot or link\index{2-bridge knot or link}\index{2-bridge knot or link!hyperbolic} with $a_i>0$ for all $i$, or $a_i<0$ for all $i$, and $|a_1|\geq 2$ and $|a_{n-1}|\geq 2$. Assume also that $n\geq 3$, so there are at least two twist regions in the diagram of $K$ given by the denominator closure of the rational tangle $T(0,a_{n-1}, \dots, a_1)$. Then $S^3-K$ is hyperbolic. 
\end{corollary}

\begin{proof}
The link complement $S^3-K$ admits a triangulation as in \refprop{2BridgeTriang}. Then \refprop{2BridgeAngleNonempty} implies the set of angle structures\index{angle structure} on this triangulation is nonempty. By \refthm{HypAngleStruct} in \refchap{Essential}, any manifold admitting an angle structure must also admit a hyperbolic structure. 
\end{proof}

\Refcor{2BridgeHyperbolic} is a special case of a stronger theorem due to Menasco determining when any alternating knot or link is hyperbolic \cite{menasco:alt}. We will return to that theorem in \refchap{Alternating}. 

\begin{remark}
\Refcor{2BridgeHyperbolic} will also follow from \refthm{2BridgeGeometric}, which we will finish proving in the next section, by an appeal to \refthm{VolAngleStructs} (volume and angle structures).\index{angle structure} While the proof of \refcor{2BridgeHyperbolic} given above appears short, in fact recall that the proof of \refthm{HypAngleStruct} requires the difficult hyperbolization theorem of Thurston, \refthm{SfcesHyperbolic}, whose proof is beyond the scope of this book. By contrast, finishing the proof of \refthm{2BridgeGeometric} requires only calculus and some calculations, and we go through it in the next section. Moreover, when finished, we will additionally know that the hyperbolic structure on 2-bridge links arises from a geometric triangulation\index{geometric triangulation} of the link complements, and that triangulation can be explicitly described. Thus a proof of \refcor{2BridgeHyperbolic} using the calculations in the next section is in many ways a ``better'' proof, worth finishing.
\end{remark}

%%%%%%%%%%%%%%%%%%%%%%%%%%%%%%%%%%%%%%%%%%%%%%%%%%%%%%%%%%%%%%%%%
\section{Maximum in interior}

In this section we conclude the proof of \refthm{2BridgeGeometric}, by proving the following.

\begin{proposition}\label{Prop:MaxInInterior}
  For the 2-bridge links\index{2-bridge knot or link}\index{2-bridge knot or link!angle structure} of \refprop{2BridgeAngleNonempty}, the volume functional\index{volume functional} $\mathcal{V}\from \mathcal{A}(\mathcal{T}) \to \RR$ cannot have a maximum on the boundary of the space of angle structures.\index{angle structure} 
\end{proposition}

\begin{remark}[Summary of proof]
  The proof is given by a series of lemmas and calculations, and is quite technical. However, the idea of the proof is straightforward. First, we show that we can use the conditions on angle structures\index{angle structure} obtained in \reftable{LabelConditions}; this is done in \reflem{VolMaxAtSymmPt}. We then assume the maximum occurs on the boundary. Using the conditions of \reftable{LabelConditions}, we find restrictions on the tetrahedra that arise; this is done in \reflem{FlatTetrConsequences}. Finally, we show that in all cases that remain there is a path from the purported maximum on the boundary of the space of angle structures\index{angle structure} to the interior for which the directional derivative of $\mathcal{V}$ is strictly increasing. This contradicts the fact that the boundary point is a maximum.
\end{remark}

We note that the technical arguments required for the proof of \refprop{MaxInInterior} are used only in this section, and are not required for other chapters of the book. 

Let $K$ be a 2-bridge knot or link as in \refprop{2BridgeAngleNonempty}. To obtain angle structures\index{angle structure} on $S^3-K$, we made some simplifying assumptions in the proof of \refprop{2BridgeAngleNonempty}. Namely, when constructing the triangulation $\mathcal{T}$, we had two tetrahedra $T_i^1$ and $T_i^2$ at each level, and we assumed that the angles on the two tetrahedra agreed. This led to the calculations of the previous section.

\begin{lemma}\label{Lem:VolMaxAtSymmPt}
The maximum of the volume functional\index{volume functional} $\mathcal{V}\from \mathcal{A}(\mathcal{T}) \to \RR$ must occur at a point for which the angles $(x_i^1,y_i^1,z_i^1)$ of $T_i^1$ agree with those $(x_i^2,y_i^2,z_i^2)$ of $T_i^2$, for all $i$, where $T_i^1$ and $T_i^2$ are the two tetrahedra constructed at the $i$-th level.
\end{lemma}

\begin{proof}
Suppose the volume is maximized at an angle structure\index{angle structure} $A$ for which angles of $T_i^1$ and $T_i^2$ do not agree. Because of the symmetry of the construction of $\mathcal{T}$, note that we obtain a new angle structure $A'$ by swapping angles of $T_i^1$ with the corresponding angles of $T_i^2$, for all tetrahedra of $A$. Note that since $A$ and $A'$ contain isometric ideal tetrahedra, $\mathcal{V}(A)=\mathcal{V}(A')$. Then $A$ and $A'$ are distinct angle structures, and the volume is maximized on both.

By \refthm{VolConcaveDown}, the volume functional\index{volume functional} is strictly concave down on $\mathcal{A}(\mathcal{T})$. Thus if the volume obtains its maximum in the interior, then that maximum is unique, and the fact that $\mathcal{V}(A)=\mathcal{V}(A')$ gives an immediate contradiction in this case. If $A$ lies on the boundary, then $A'$ also lies on the boundary. Because $\mathcal{A}(\mathcal{T})$ is convex (\refprop{AngleStructsPoly}), the line between $A$ and $A'$ lies in $\mathcal{A}(\mathcal{T})$. But then along this line, the second directional derivative in the direction of the line is strictly negative, which implies the maximum cannot occur at the endpoints. This is a contradiction.
\end{proof}

By \reflem{VolMaxAtSymmPt}, we may assume angles of $T_i^1$ and $T_i^2$ agree. Thus we may use the conditions on angles in \reftable{LabelConditions} that we calculated in the previous section to prove \refprop{MaxInInterior}.

Now assume that the maximum of $\mathcal{V}$ does occur on the boundary of the space of angle structures\index{angle structure} $\mathcal{A}(\mathcal{T})$. Then there will be a flat tetrahedron (or more accurately, a pair of flat tetrahedra). We will slowly narrow in on what type of tetrahedron it is, and where it occurs in the triangulation.

Note we will switch notation slightly. Rather than referring to the two tetrahedra between $S_i$ and $S_{i+1}$ as $T_i^1$ and $T_i^2$, we will simply refer to such a tetrahedron by $\Delta_i$. Because the angles of $T_i^1$ and $T_i^2$ can be assumed to agree by \reflem{VolMaxAtSymmPt}, this will simplify our notation. 

\begin{lemma}\label{Lem:FlatTetrConsequences}
Suppose the maximum of the volume functional\index{volume functional} occurs on the boundary of $\mathcal{A}(\mathcal{T})$. 
  \begin{enumerate}
  \item[(1)] Then there exists a flat tetrahedron in the triangulation of the 2-bridge link.
  \item[(2)] The flat tetrahedron is not adjacent to any other flat tetrahedra.
  \item[(3)] The flat tetrahedron is not adjacent to a hairpin turn.
  \item[(4)] The flat tetrahedron occurs at a hinge, and satisfies $z_i=\pi$, and for the two adjacent tetrahedra, $z_{i-1}=z_{i+1}$.
  \item[(5)] If some tetrahedron $\Delta_i$ is type LL or RR, then $\Delta_{i-1}$ and $\Delta_{i+1}$ cannot both be flat.
  \end{enumerate}
\end{lemma}

\begin{proof}
By \refprop{NoSingleDegenerate}, if the volume takes its maximum at an angle structure\index{angle structure} for which a tetrahedron has an angle equal to $0$, then it must have two angles equal to $0$ and one equal to $\pi$. This is a flat tetrahedron. Because we are assuming the maximum is on the boundary, there must be a flat tetrahedron in the triangulation, say tetrahedron $\Delta_i$ is flat, where $2\leq i\leq C-2$. This proves item (1). 

There are three cases for the angles, namely $(x_i, y_i, z_i)$ can equal $(0,0,\pi)$,  $(0,\pi,0)$, or $(\pi, 0,0)$. There are also four possibilities for the tetrahedron: type $LL$, $RR$, $LR$, or $RL$. The equations of \reftable{LabelConditions} give us angles of adjacent tetrahedra in all cases, and an analysis of these will lead to the conclusions of the lemma.

\bigskip

\subsubsection*{Case $(x_i,y_i,z_i) = (0,0,\pi)$}
\begin{enumerate}
\item[$LL$,] $RR$: Equations of \reftable{LabelConditions} imply
  \[0  = \half(2\pi-z_{i-1}-z_{i+1}), \mbox{ which implies } z_{i-1}=z_{i+1}=\pi.\]
  In this case, both adjacent tetrahedra must be flat.
\item[$LR$,] $RL$: Equations of \reftable{LabelConditions} imply
  \[0 = \half(z_{i+1}-z_{i-1}), \mbox{ or } z_{i-1}=z_{i+1}. \]
  Note in this case, it is not necessarily true that both adjacent tetrahedra are flat, but if one is flat then so is the other.
\end{enumerate}

\subsubsection*{Case $(x_i,y_i,z_i)=(0,\pi,0)$}.
\begin{enumerate}
\item[$LL$:] $0=\half(2\pi-z_{i-1}-z_{i+1})$ implies $z_{i-1}=z_{i+1}=\pi$. 
\item[$RR$:] $0=\half(z_{i-1}+z_{i+1})$ implies $z_{i-1}=z_{i+1}=0$.
\item[$LR$:] $0=\half(\pi-z_{i-1}+z_{i+1})$ implies $\pi+z_{i+1}=z_{i-1}$. Since angles lie in $[0,\pi]$, it follows that $z_{i+1}=0$ and $z_{i-1}=\pi$.
\item[$RL$:] Similar to the last case, $z_{i+1}=\pi$ and $z_{i-1}=0$.\\
For all types of tetrahedra in this case, the two tetrahedra adjacent to $\Delta_i$ are flat. 
\end{enumerate}

\subsubsection*{Case $(x_i,y_i,z_i)=(\pi,0,0)$.}
\begin{enumerate}
\item[$LL$:] Equations of \reftable{LabelConditions} imply $z_{i-1}=z_{i+1}=0$.
\item[$RR$:] $z_{i-1}=z_{i+1}=\pi$.
\item[$LR$:] $z_{i-1}=0$, $z_{i+1}=\pi$.
\item[$RL$:] $z_{i+1}=0$, $z_{i-1}=\pi$. \\
Again this shows that the two adjacent tetrahedra are both flat in this case. 
\end{enumerate}

In all cases, if two adjacent tetrahedra are flat, then the next adjacent tetrahedron is also flat. It follows that if there are two adjacent flat tetrahedra, then all tetrahedra are flat, and the structure has zero volume, which cannot be a maximum for the volume. Thus we cannot have two adjacent flat tetrahedra. This proves item (2).

Moreover, the only case that does not immediately imply multiple adjacent flat tetrahedron is the first case, with $z_i=\pi$, for the hinge tetrahedra $RL$ or $LR$, and the calculation above gives the relationship $z_{i-1}=z_{i+1}$, proving item (4).

If the tetrahedron is adjacent to a hairpin turn, then $i=2$ or $i=C-2$, and $z_i=\pi$. We also have $z_1=0$ and $z_{C-1}=0$, hence in either case the equations above imply that a next adjacent tetrahedron, corresponding to $z_3$ or $z_{C-3}$, is flat, and thus all tetrahedra are flat, contradicting item (2). This proves (3).

Now suppose $\Delta_{i-1}$ and $\Delta_{i+1}$ are flat. By the previous work, we know $z_{i-1}=z_{i+1}=\pi$. If $\Delta_i$ is type $LL$, the equations of \reftable{LabelConditions} imply $x_i=\half(2\pi-\pi-\pi)=0$, so $\Delta_i$ is flat. Similarly if $\Delta_i$ is of type $RR$, then $y_i=0$ and $\Delta_i$ is flat. But then we have three adjacent flat tetrahedra, contradicting item (2). This proves item (5). 
\end{proof}

We now know that any flat tetrahedron occurring in a maximum for $\mathcal{V}$ on the boundary has a very particular form. To finish the proof of \refprop{MaxInInterior}, we will show that the maximum cannot occur in the remaining cases. For the argument, we will find a path through the space of angle structures\index{angle structure} starting at the purported maximum for $\mathcal{V}$ on the boundary, and then show that the derivative at time $0$ in the direction of this path is strictly positive. This will contradict the fact that the point is a maximum.

The paths we consider adjust the angles of the flat tetrahedron $\Delta_i$ by
\[ (x_i(\epsilon),y_i(\epsilon),z_i(\epsilon)) = ((1+\lambda)\epsilon, (1-\lambda)\epsilon, \pi -2\epsilon),\]
where $\epsilon\to 0$ and $\lambda$ will be a carefully chosen constant. 
In such a path, we will leave as many angles unchanged away from the $i$-th tetrahedron as possible. However, the equations in \reftable{LabelConditions} imply that many angles of adjacent tetrahedra must change with $\epsilon$ as well. 

Recall from \reflem{DerivativesVol} that the derivative of the volume functional\index{volume functional} in the direction of a vector $w=(w_1, \dots, w_n)$ at a point $a=(a_1, \dots, a_n)$ is 
\[ \frac{\partial \mathcal{V}}{\partial w} = \sum_{i=1}^{3n} -w_i\log\sin a_i. \]
The terms of the sum are grouped into threes, with each group corresponding to a single tetrahedron, with derivative coming from \refthm{VolConcaveDown}. 

\begin{lemma}\label{Lem:FrancoisPath}
  Let $\gamma(t)$ be a path through $\overline{\mathcal{A}(\mathcal{T})}$ with the angles of the $i$-th tetrahedron $\Delta_i$ in $\gamma(t)$ satisfying $(x_i,y_i,z_i) = ((1+\lambda)t, (1-\lambda)t,\pi-2t)$. Then the derivative of the volume of $\Delta_i$ along this path at $t=0$ satisfies
  \[ \left. \frac{d\vol(\Delta_i)}{dt} \right\rvert_{t=0} = \log \left( \frac{4}{1-\lambda^2}\left(\frac{1-\lambda}{1+\lambda}\right)^\lambda\right).\]
\end{lemma}

\begin{proof}
By \refthm{VolConcaveDown}, the derivative in the direction of $\gamma'(0) = w= ((1+\lambda), (1-\lambda), -2)$ is
\begin{align*}
  \frac{\partial \vol}{\partial w} & = \lim_{t\to 0} \big[  -(1+\lambda)\log\sin((1+\lambda)t) -(1-\lambda)\log \sin((1-\lambda)t)\\
    & \hspace{.5in} +2\log\sin(\pi-2t) \big].
\end{align*}
Using the Taylor expansion $\sin(At) = At$ near $t=0$, this becomes
\begin{align*}
  \frac{\partial \vol}{\partial w} &
  \lim_{t\to 0} \big[-(1+\lambda)\log((1+\lambda) t) - (1-\lambda)\log((1-\lambda)t) + 2\log(2t)\big] \\
  & = \log \left( \frac{4}{(1+\lambda)(1-\lambda)}\left(\frac{1-\lambda}{1+\lambda}\right)^\lambda \right)\qedhere
\end{align*}
\end{proof}

We denote the location of a flat tetrahedron by a vertical line:
\[\dots LL|RR\dots.\]
By \reflem{FlatTetrConsequences}, a vertical line can only appear at a hinge: $L|R$ or $R|L$; at least two letters lie between consecutive vertical lines; and patterns $L|RR|L$ and $R|LL|R$ cannot occur. The remaining cases are $LR|LR$ and $RL|RL$, which we deal with simultaneously; $RR|LR$ and $LL|RL$ and their reversals $RL|RR$ and $LR|LL$; and $RR|LL$ and $LL|RR$.

In all cases, we find a path $\gamma(t)$ through $\mathcal{A}(\mathcal{T})$ with $\gamma(0)$ a point on the boundary with the flat tetrahedron specified in the given case. 

\subsection*{Case $LR|LR$ and $RL|RL$:}

Begin with the $LR|LR$ case. 
Let $\Delta_i$ denote the flat tetrahedron, with $(x_i,y_i,z_i)=(0,0,\pi)$. We take a path $\gamma(t)$ to satisfy $(x_i(t), y_i(t), z_i(t)) = (t,t,\pi-2t)$, i.e.\ $\lambda=0$ in \reflem{FrancoisPath}, and we will keep as many other angles constant as possible. The formulas in \reftable{LabelConditions} imply that angles of tetrahedra $\Delta_{i-1}$ and $\Delta_{i+1}$ must also vary, as in the following table. In the table, we let $z_{i-1}=z_{i+1}=w$ (required by \reflem{FlatTetrConsequences}(4)), and we let $u=z_{i-2}$, $v=z_{i+2}$.

\bigskip

\begin{tabular}{c|c|c|c}
  Angle & $\Delta_{i-1}$ & $\Delta_i$ & $\Delta_{i+1}$ \\
  \hline\noalign{\smallskip}
  $x$ & $\half(2\pi - u -w - 2t)$ & $t$ & $\half(2t - w + v)$ \\
  $y$ & $\half(u-w+2t)$ & $t$ & $\half(2\pi - 2t - w - v)$ \\
  $z$ & $w$ & $\pi-2t$ & $w$
\end{tabular}

\bigskip

Thus the derivative vector to the path at time $t=0$ is
\[ \gamma'(t) = (0, \dots, 0, \underbrace{-1, 1, 0,}_{\Delta_{i-1}} \underbrace{1,1,-2,}_{\Delta_i} \underbrace{1,-1,0,}_{\Delta_{i+1}} 0, \dots, 0). \]

Hence the derivative of the volume functional\index{volume functional} in the direction of the path is given by
\begin{align*}
  \left.\frac{d\mathcal{V}}{dt}\right\rvert_{t=0} & =
  \log\sin\left(\half(2\pi-u-w)\right) - \log\sin\left(\half(u-w)\right) + \log 4 \\
  & \hspace{.5in} -\log\sin\left(\half(v-w)\right) + \log\sin\left(\half(2\pi-v-w)\right) \\
  & = \log \left( \frac{ 4 \sin(u/2-w/2)\sin(v/2-w/2)}{\sin(u/2+w/2)\sin(v/2+w/2)} \right) > 0.
\end{align*}
Note this is strictly positive, hence the volume functional\index{volume functional} cannot have a maximum at this boundary point. The calculation is similar for $RL|RL$.

\subsection*{Remaining cases:}

We will first take care of cases $RR|LR$ and $RR|LL$.

As in the previous case, we will take a path such that the flat tetrahedron $\Delta_i$ changes. This time, we will find a fixed $\lambda$ such that angles of $\Delta_i$ satisfy $(x_i,y_i,z_i)=((1-\lambda)t,(1+\lambda)t, \pi-2t)$ for $t\in [0, \epsilon]$, for some $\epsilon>0$. At time $t=0$, we require $z_{i-1}=z_{i+1}=w$, a constant. Set $z_{i-2}=u$ and $z_{i+2}=v$, also constant. Additionally, adjust $z_{i-1}$ so that at time $t$, $z_{i-1}=w-2\lambda t$. In the argument below, we will assume that $i-2\neq 1$, so there is a tetrahedron $\Delta_{i-2}$. We also need to consider the case $i-2=1$; we will do this at the very end of the proof. 
Assuming $i-2\neq 1$, the angles that are modified are shown in the tables below for the cases $RR|LR$ and $RR|LL$.

\bigskip

\noindent\begin{tabular}{c|cccccccc}
  & & $R$ && $R$ & $|$ & $L$ && $R$ \\
  Angle & $\Delta_{i-2}$ &&  $\Delta_{i-1}$ && $\Delta_i$ && $\Delta_{i+1}$ &\\
  \hline\noalign{\smallskip}
  $x$ & $A+\half w-\lambda t$ && $x_{i-1}(t,\lambda)$ && $(1+\lambda)t$ && $x_{i+1}(t,\lambda)$ &\\
  $y$ & $A'-\half w+\lambda t$ && $\half(\pi-2t)$ && $(1-\lambda)t$ && $\half(2t - w + v)$ &\\
  $z$ & $u$ && $w-2\lambda t$ && $\pi-2t$ && $w$ &
\end{tabular}

\bigskip

\noindent\begin{tabular}{c|cccccccc}
  & & $R$ && $R$ & $|$ & $L$ && $L$ \\
  Angle & $\Delta_{i-2}$ &&  $\Delta_{i-1}$ && $\Delta_i$ && $\Delta_{i+1}$ &\\
  \hline\noalign{\smallskip}
  $x$ & $A+\half w-\lambda t$ && $x_{i-1}(t,\lambda)$ && $(1+\lambda)t$ && $\half(\pi+2t-v)$ &\\
  $y$ & $A'-\half w+\lambda t$ && $\half(\pi-2t)$ && $(1-\lambda)t$ && $y'_{i+1}(t,\lambda)$ &\\
  $z$ & $u$ && $w-2\lambda t$ && $\pi-2t$ && $w$ &
\end{tabular}

\bigskip

Here $A$ and $A'$ are constants, $x_{i-1}(t,\lambda) = \half(u-2w+4\lambda t + \pi-2t)$,  $x_{i+1}(t,\lambda) = \half(2\pi-2t - w - v)$, and $y'_{i+1}(t,\lambda) = \half(\pi-2t-2w+v)$.

If $i>3$, we may use the table to compute the derivative in the direction of the path, and find in the case $RR|LR$, $d\mathcal{V}/dt \rvert_{t=0}$ equals:
\begin{equation}\label{Eqn:RRLRDeriv}
\left.\frac{d\mathcal{V}}{dt}\right\rvert_{t=0} =  \log \left( \frac{4}{1-\lambda^2}\frac{\sin(\frac{v}{2}+\frac{w}{2})}{\sin(\frac{v}{2}-\frac{w}{2})}
\frac{\sin x_{i-1}}{\sin y_{i-1}} \left( \frac{1-\lambda}{1+\lambda} \cdot
\frac{\sin x_{i-2}}{\sin{y_{i-2}}} \frac{\sin^2 z_{i-1}}{\sin^2 x_{i-1}} \right)^\lambda\right).
\end{equation}

And in the case $RR|LL$, $d\mathcal{V}/dt \rvert_{t=0}$ equals:
\begin{equation}\label{Eqn:RRLLDeriv}
  \left.\frac{d\mathcal{V}}{dt}\right\rvert_{t=0} =
  \log \left( \frac{4}{1-\lambda^2} \frac{\sin x_{i-1} }{\sin y_{i-1}}
  \frac{\sin y_{i+1}}{\sin x_{i+1}}
  \left( \frac{1-\lambda}{1+\lambda} \cdot
  \frac{\sin x_{i-2}}{\sin y_{i-2}} \frac{\sin^2 z_{i-1}}{\sin^2 x_{i-1}}
  \right)^\lambda\right).
\end{equation}

\begin{lemma}\label{Lem:LambdaFact}
  Let $X$, $Y$ be positive constants, and let
  \[ f(\lambda) = \log \left( \frac{4}{1-\lambda^2}\, X \left(\frac{1-\lambda}{1+\lambda} \, Y \right)^{\lambda} \right). \]
  Then $f$ has a critical point at $\lambda = (Y-1)/(Y+1)$, and $f$ takes the value $\log (X (Y+1)^2/Y)$ at this point. 
\end{lemma}

\begin{proof}
  Calculus.
\end{proof}

Now apply \reflem{LambdaFact} to \refeqn{RRLRDeriv} and \refeqn{RRLLDeriv}, choosing $\lambda$ to be the value given by that lemma at time $t=0$. 
For this value of $\lambda$, we obtain the following:

The derivative $d\mathcal{V}/dt \rvert_{t=0}$ in the case $RR|LR$ equals:
\begin{align}
\nonumber  & \log
  \left( \frac{\sin(\frac{v}{2}+\frac{w}{2})}{\sin(\frac{v}{2}-\frac{w}{2})}
  \frac{\sin x_{i-1}}{\sin y_{i-1}} \left( 1 + 
  \frac{\sin x_{i-2}}{\sin{y_{i-2}}} \frac{\sin^2 z_{i-1}}{\sin^2 x_{i-1}} \right)^2
\frac{\sin y_{i-2}}{\sin{x_{i-2}}} \frac{\sin^2 x_{i-1}}{\sin^2 z_{i-1}}
  \right) \\
\label{Eqn:RRLRDerivSimp}  & \geq \log \left( 
  \frac{\sin x_{i-1}}{\sin y_{i-1}} \left( 1 + 
  \frac{\sin x_{i-2}}{\sin{y_{i-2}}} \frac{\sin^2 z_{i-1}}{\sin^2 x_{i-1}} \right)^2
\frac{\sin y_{i-2}}{\sin{x_{i-2}}} \frac{\sin^2 x_{i-1}}{\sin^2 z_{i-1}}
  \right).
\end{align}

The derivative $d\mathcal{V}/dt \rvert_{t=0}$ in the case $RR|LL$ equals:
\begin{equation}\label{Eqn:RRLLDerivSimp}
  \log \left( \frac{\sin x_{i-1} }{\sin y_{i-1}}
  \frac{\sin y_{i+1}}{\sin x_{i+1}}
  \left( 1 + \frac{\sin x_{i-2}}{\sin y_{i-2}} \frac{\sin^2 z_{i-1}}{\sin^2 x_{i-1}}
  \right)^2
\frac{\sin y_{i-2}}{\sin{x_{i-2}}} \frac{\sin^2 x_{i-1}}{\sin^2 z_{i-1}}
  \right).
\end{equation}

The remaining quantities $\sin a / \sin b$ are geometric: by the law of sines, they give a ratio of lengths of triangles, and the triangles are those from our cusp triangulation, as in \reffig{CuspLabels}.

\begin{figure}[h]
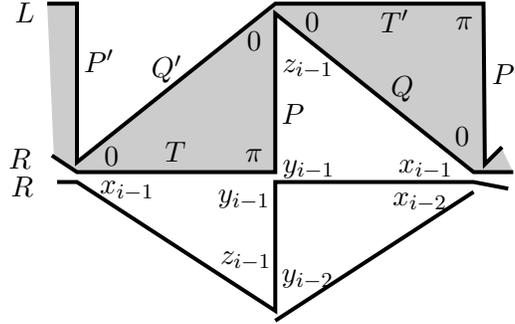
  \caption{Segments of length $P$, $Q$, and $T$ shown on the zigzag corresponding to the first $R$ after a flat hinge tetrahedron. For the $RR|LL$ case, segments of length $P'$, $Q'$ and $T'$ also shown on the first $L$ zigzag after the flat hinge tetrahedron.}
  \label{Fig:Segments}
\end{figure}

\begin{lemma}\label{Lem:SineRelations}
  In the case $RR|L$, let $P$, $Q$, $T$ be the lengths of segments on the middle zigzag $R$, with $P$ opposite the angle $x_{i-1}$, $Q$ opposite the angle $y_{i-1}$ and $T$ opposite the angle $z_{i-1}$, as in \reffig{Segments}. Then the following hold:
\[
\frac{\sin x_{i-2}}{\sin y_{i-2}} \frac{\sin^2 z_{i-1}}{\sin^2 x_{i-1}} = \frac{T}{P}, \qquad
 \frac{\sin x_{i-1}}{\sin y_{i-1}} = \frac{P}{Q}.\]
\end{lemma}

\begin{proof}
  The equations follow from the law of sines.
\end{proof}

\begin{lemma}\label{Lem:Geometrical}
  Suppose there is a subword $L|R^k L$ with $k\geq 2$. Let $Q$, $P$, and $T$ be lengths of segments of the zigzag corresponding to the first $R$, with $P$ and $T$ adjacent to the angle labeled $z=\pi$ on the hinge tetrahedron $L|R$. Then $P+T > Q$.

  Similarly, if there is a subword $R|L^k R$ with $k\geq 2$, and $Q'$, $P'$, and $T'$ denote the lengths of the segments of the zigzag corresponding to the first $L$, with $P'$ and $T'$ adjacent to the angle $z=\pi$ on the hinge tetrahedron $R|L$, then $P'+T'>Q'$.
\end{lemma}

The labels $P$, $Q$, and $T$ are illustrated in \reffig{Segments}. In the case there are at least two $L$'s at the top of the figure, $P'$, $Q'$ and $T'$ will be labeled as shown there as well.

\begin{proof}
  \cite[Lemma~8.2]{GueritaudFuter:2bridge}.
\end{proof}

Now we can show in the $RR|LR$ case the derivative $d\mathcal{V}/dt\rvert_{t=0}$ is positive. From \refeqn{RRLRDerivSimp}, we obtain
\[
\left.\frac{d\mathcal{V}}{dt}\right\rvert_{t=0} \geq
\log \left( \frac{P}{Q} \left(1 + \frac{T}{P} \right)^2 \frac{P}{T} \right) =
\log \left( \frac{P+T}{T} \cdot \frac{P+T}{Q}\right) > \log(1)=0. \]

A similar calculation holds in the $LL|RL$ case. By swapping the indices $i-1$, $i+1$, and $i-2$, $i+2$, the same argument shows the derivative is strictly positive in the $RL|RR$ and $LR|LL$ cases, provided $i+2$ is not the index of a hairpin turn.

We now finish the $RR|LL$ case. 

\begin{lemma}\label{Lem:Geom2}
  In the $RR|LL$ case, with $P'$, $Q'$, and $T'$ as in \reffig{Segments}, $P'/T' = T/P$, and $\sin(y_{i+1})/\sin(x_{i+1}) = P'/Q'$. 
\end{lemma}

\begin{proof}
  By \reflem{FlatTetrConsequences}, item (5), there is no flat tetrahedron either directly before or directly after the sequence $RR|LL$, so the angles of $\Delta_{i-2}$, $\Delta_{i-1}$, $\Delta_{i+1}$, and $\Delta_{i+2}$ are all positive. Thus the parameter $w$ can vary freely in an open interval when $t=0$. Since the volume is maximized, the derivative with respect to $w$ satisfies
  \[ \left.\frac{d\mathcal{V}}{dw}\right\rvert_{t=0} =
  \log \left( \sqrt{\frac{\sin x_{i-2}}{\sin y_{i-2}}} \cdot \frac{\sin z_{i-1}}{\sin x_{i-1}} \cdot \frac{\sin z_{i+1}}{\sin y_{i+1}} \cdot \sqrt{\frac{\sin y_{i+2}}{\sin x_{i+2}}} \right) = 0. \]
  Thus
  \[ \frac{\sin y_{i-2} \sin^2 x_{i-1}}{\sin x_{i-2} \sin^2 z_{i-1}} \cdot \frac{\sin^2 y_{i+1} \sin x_{i+2}}{\sin^2 z_{i+1} \sin y_{i+2}} =1. \]

Using an expanded version of \reffig{Segments}, one can check (exercise) that
  \[ \frac{\sin y_{i-2} \sin^2 x_{i-1}}{\sin x_{i-2} \sin^2 z_{i-1}}  = \frac{P}{T}, \quad \mbox{ and } \quad
  \frac{\sin^2 y_{i+1} \sin x_{i+2}}{\sin^2 z_{i+1} \sin y_{i+2}}  = \frac{P'}{T'}. \]
  This shows $P'/T'=T/P$.

  Similarly using \reffig{Segments}, one can check that $\sin(y_{i+1})/\sin(x_{i+1}) = P'/Q'$.
\end{proof}

By \reflem{Geom2} and \refeqn{RRLLDerivSimp}, we find that in the $RR|LL$ case,
\begin{align*}
  \left.\frac{d\mathcal{V}}{dt}\right\rvert_{t=0} & =
  \log \left( \frac{P}{Q}\frac{P'}{Q'}\left( 1+\frac{T}{P}\right)^2\frac{P}{T}\right) \\
  & = \log \left( \frac{P}{Q}\left(1+\frac{T}{P}\right) \frac{P'}{Q'}\left(1+\frac{P'}{T'}\right) \frac{T'}{P'}\right) \\
  & = \log \left( \frac{P+T}{Q} \cdot \frac{P'+T'}{Q'} \right) \\
  & > \log (1) = 0. 
\end{align*}

A similar calculation takes care of the $LL|RR$ case.

So far, we have argued only for $i>3$. 
It remains to consider what happens when $i=3$. In this case, $i-2=1$ is the index of a hairpin turn, and the terms $\sin y_1/\sin x_1$ disappear from the computations of $d\mathcal{V}/dt$ in \refeqn{RRLRDeriv} and \refeqn{RRLLDeriv}. We have a result similar to \reflem{Geometrical}: $R^aL$ is a tessellated Euclidean triangle, and lengths still behave as in \reflem{Geometrical} to give the same result; see \cite[Lemma~1.5]{GueritaudFuter:2bridge}.

This concludes the proof of \refprop{MaxInInterior}. \qed

\bigskip

We now assemble the pieces to obtain the stronger result, \refthm{2BridgeGeometric}. 

\begin{proof}[Proof of \refthm{2BridgeGeometric}]
Let $K$ be a knot or link with a reduced alternating diagram with at least two twist regions. Let $\mathcal{T}$ be the triangulation of $S^3-K$ described in this chapter. By \refprop{2BridgeAngleNonempty}, the space of angle structures\index{angle structure} $\mathcal{A}(\mathcal{T})$ is nonempty. By \refprop{MaxInInterior}, the volume functional\index{volume functional} $\mathcal{V}\from \mathcal{A}(\mathcal{T}) \to \RR$ cannot have a maximum on the boundary of the space of angle structures. It follows that the maximum of $\mathcal{V}$ is on the interior of the space of angle structures.\index{angle structure} Let $A\in\mathcal{A}(\mathcal{T})$ denote this critical point. Thus by \refthm{VolAngleStructs} (volume and angle structures), the ideal hyperbolic tetrahedra obtained from the angle structure $A$ give $S^3-K$ a complete hyperbolic structure. Note since $A$ lies in the interior, the ideal hyperbolic tetrahedra it determines are all positively oriented,\index{positively oriented tetrahedron}\index{tetrahedron!positively oriented} as claimed. 
\end{proof}

%%%%%%%%%%%%%%%%%%%%%%%%%%%%%%%%%%%%%%%%%%%%%%%%%%%%%%%%%%%%%%%%%
\section{Exercises}

\begin{exercise}
  Sketch rational tangles and diagrams of 2-bridge links\index{2-bridge knot or link} associated to the following continued fractions: $[3,2]$, $[0,3,2]$, $[1,3,2]$. 
\end{exercise}

\begin{exercise}\label{Ex:ContinuedFractions}
  Continued fractions. Show every rational number has a continued fraction expansion $p/q=[a_n, a_{n-1}, \dots, a_1]$ such that if $i<n$, then $a_i\neq 0$, and such that if $p/q>0$, then each $a_i\geq 0$, while if $p/q<0$, then each $a_i\leq 0$.
\end{exercise}

\begin{exercise}\label{Ex:a1an}
  Prove \reflem{a1an}. That is, show that if $K[a_{n-1}, \dots, a_1]$ is a 2-bridge knot or link,\index{2-bridge knot or link} then we may assume that $|a_{n-1}|\geq 2$ and $|a_1|\geq 2$.
\end{exercise}

\begin{exercise}
  Work through the identification of tetrahedra at the innermost crossing. Prove that faces of the innermost tetrahedra are glued in pairs, two triangles of one tetrahedron glued to triangles of the opposite tetrahedron. Why is there no need to consider both horizontal and vertical crossings for the innermost crossing?
\end{exercise}

\begin{exercise}\label{Ex:2TwistRegions}
  In \refprop{2BridgeTriang}, we require at least two twist regions. Show that this requirement is necessary by showing that the construction fails to give a triangulation of a knot or link with just one twist region. What breaks down?
\end{exercise}

\begin{exercise}
  This exercise asks you to consider hairpin turns.
  \begin{enumerate}
  \item Prove if $|a_{n-1}|\geq 3$, there is a 3-valent vertex of the cusp triangulation.
  \item Prove all vertices aside from possibly a single vertex in a hairpin turn must have valence at least four.
  \item If $|a_{n-1}|=2$, prove the vertex corresponding to the outside hairpin turn may have arbitrarily high valence.
  \end{enumerate}
\end{exercise}

\begin{exercise}
  Use the methods of this chapter to find the form of the cusp triangulation for the twist knot $J(2,n)$. How many tetrahedra are in its decomposition?
\end{exercise}

\begin{exercise}
  Find the form of the cusp triangulation for $J(k, \ell)$, where one of $k$, $\ell$ is even. How many tetrahedra are in its decomposition?
\end{exercise}

\begin{exercise}
  Find the form of the cusp triangulation of a 2-bridge knot\index{2-bridge knot or link}\index{2-bridge knot or link!cusp triangulation} with exactly three twist regions.
\end{exercise}

\begin{exercise}
  Prove \reflem{PleatingAngles}: that pleating angles as in \reffig{PleatingAngles} satisfy $\alpha_1+\alpha_2-\alpha_3=0$ when the cusp is Euclidean.
\end{exercise}

\begin{exercise}
  Determine the labels of the hairpin turns of the form $RR$, $LR$, $RL$, similar to \reffig{HairpinLabels}. 
\end{exercise}

\begin{exercise}
  In the cases $RR|LR$ and $RR|LL$, compute the derivative $d\mathcal{V}/dt \rvert_{t=0}$ and check that it agrees with the formulas given in \refeqn{RRLRDeriv} or \refeqn{RRLLDeriv}. 
\end{exercise}

\begin{exercise}
Give the proof of \reflem{Geometrical}. 
\end{exercise}

\begin{exercise}
  Work through the geometric details of \reflem{Geom2}. First, sketch the zigzag labeled $L$ at the top of \reffig{Segments}, along with angles at its corners, and show that:
  \[ \frac{\sin y_{i-2} \sin^2 x_{i-1}}{\sin x_{i-2} \sin^2 z_{i-1}}  = \frac{P}{T}, \quad \mbox{ and } \quad
\frac{\sin^2 y_{i+1} \sin x_{i+2}}{\sin^2 z_{i+1} \sin y_{i+2}}  = \frac{P'}{T'}. \]
  Also show
  $\sin(y_{i+1})/\sin(x_{i+1}) = P'/Q'$. 
\end{exercise}

\begin{exercise}
  Go carefully through the proof of cases $RR|LR$ and $RR|LL$ when the index of the flat tetrahedron is $i=3$. 
\end{exercise}

\chapter{Alternating Knots and Links}\label{Chap:Alternating}
\blfootnote{Jessica S. Purcell, Hyperbolic Knot Theory}

Alternating knots and links\index{alternating knot or link} need their own chapter, because there is a wealth of geometric information coming from them. Of all knots, alternating knots seem to have hyperbolic geometry most closely related to their diagrams. As of the writing of this book, there are many open conjectures concerning how the geometry and diagrams interact. 

An \emph{alternating knot or link}\index{alternating knot or link} has a diagram with an orientation such that when following the knot in the direction of the orientation, the crossings alternate between over and under, all the way along the diagram.\index{alternating diagram} Alternating knots account for large numbers of knots with small crossing numbers, but they are less prevalent among knots with higher crossing numbers. Indeed, the proportion of alternating knots and links among all prime $n$-crossing knots and links is known to approach zero exponentially as $n$ approaches infinity \cite{SundbergThistlethwaite, Thistlethwaite:Tangles}. The first non-alternating knot in the knot tables has eight crossings. Note that it takes some work to prove that a knot is non-alternating: one must show that among all possible knot diagrams, there is no alternating diagram.\index{alternating diagram} In this chapter, we won't consider the question of whether a knot with a non-alternating diagram actually is alternating. Instead we will assume we have an alternating knot diagram, and consider what this implies for the geometry of the knot complement. 

One main result of the chapter is a proof of a theorem originally due to Menasco that identifies when alternating knots and links are hyperbolic \cite{menasco:alt}. We also define checkerboard surfaces of these links, and show they are essential.\index{essential}

%%%%%%%%%%%%%%%%%%%%%%%%%%%%%%%%%%%%%%%%%%%%%%%%%%%%%%%%%%%%%%%%%
\section{Alternating diagrams and hyperbolicity}

Since we are interested in hyperbolic knots and links, we will consider only \emph{connected diagrams} of knots and links throughout; that is, the underlying 4-valent diagram graph is connected.\index{connected diagram} For note that if a diagram is not connected, it contains an obvious essential\index{essential} 2-sphere, namely one separating two diagram components. Since a hyperbolic knot or link complement can contain no essential 2-sphere, we restrict to connected diagrams. 

We also wish to work with diagrams that have been simplified in obvious ways. For example, we wish to untwist all reducible crossings,\index{reducible crossing} like those shown in \reffig{Nugatory} in \refchap{Fig8Decomp}.

We also wish to work with prime diagrams, also defined in \refchap{Fig8Decomp}. We recall the definition again. 

\begin{figure}
  \begin{center}
%% Accompanies image file 'F11-01-Prime.eps' (pdf, eps, ps)
\begingroup%
  \makeatletter%
  \providecommand\color[2][]{%
    \errmessage{(Inkscape) Color is used for the text in Inkscape, but the package 'color.sty' is not loaded}%
    \renewcommand\color[2][]{}%
  }%
  \providecommand\transparent[1]{%
    \errmessage{(Inkscape) Transparency is used (non-zero) for the text in Inkscape, but the package 'transparent.sty' is not loaded}%
    \renewcommand\transparent[1]{}%
  }%
  \providecommand\rotatebox[2]{#2}%
  \newcommand*\fsize{\dimexpr\f@size pt\relax}%
  \newcommand*\lineheight[1]{\fontsize{\fsize}{#1\fsize}\selectfont}%
  \ifx\svgwidth\undefined%
    \setlength{\unitlength}{320.70650482bp}%
    \ifx\svgscale\undefined%
      \relax%
    \else%
      \setlength{\unitlength}{\unitlength * \real{\svgscale}}%
    \fi%
  \else%
    \setlength{\unitlength}{\svgwidth}%
  \fi%
  \global\let\svgwidth\undefined%
  \global\let\svgscale\undefined%
  \makeatother%
  \begin{picture}(1,0.1923997)%
    \lineheight{1}%
    \setlength\tabcolsep{0pt}%
    \put(0,0){\includegraphics[width=\unitlength]{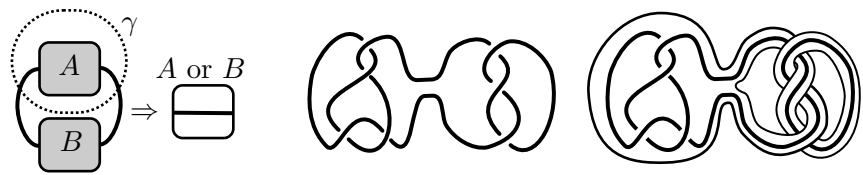}}%
    \put(0.05648321,0.11583667){\color[rgb]{0,0,0}\makebox(0,0)[lt]{\lineheight{1.25}\smash{\begin{tabular}[t]{l}$A$\end{tabular}}}}%
    \put(0.14093648,0.06445575){\color[rgb]{0,0,0}\makebox(0,0)[lt]{\lineheight{1.25}\smash{\begin{tabular}[t]{l}$\Rightarrow$\end{tabular}}}}%
    \put(0.05825497,0.02665881){\color[rgb]{0,0,0}\makebox(0,0)[lt]{\lineheight{1.25}\smash{\begin{tabular}[t]{l}$B$\end{tabular}}}}%
    \put(0.17037064,0.1117024){\color[rgb]{0,0,0}\makebox(0,0)[lt]{\lineheight{1.25}\smash{\begin{tabular}[t]{l}$A$ or $B$\end{tabular}}}}%
    \put(0.130306,0.16839836){\color[rgb]{0,0,0}\makebox(0,0)[lt]{\lineheight{1.25}\smash{\begin{tabular}[t]{l}$\gamma$\end{tabular}}}}%
  \end{picture}%
\endgroup%

  \end{center}
  \caption{Left: a prime diagram. Middle: a diagram that is not prime. Right: a swallow--follow torus}
  \label{Fig:Prime}
\end{figure}

\begin{definition}\label{Def:Prime}
A diagram is \emph{prime}\index{prime}\index{prime!diagram} if, for every simple closed curve $\gamma$ in the plane of projection, if $\gamma$ meets the knot exactly twice transversely away from crossings, then $\gamma$  bounds a region of the diagram with no crossings. See \reffig{Prime}, left.
\end{definition}

\Reffig{Prime}, middle, shows an example of a diagram that is not prime: a curve running through the center of the diagram meets the knot exactly twice, with crossings on both sides. It is constructed of two simpler knots via the following procedure, which we also defined in \refchap{KnotIntro}; see \reffig{KnotSum}.

\begin{definition}\label{Def:ConnectedSum}
Given two knots $K_1$ and $K_2$ in $S^3$, form their \emph{knot sum}\index{knot sum} or \emph{connected sum}\index{connected sum} as follows. For each knot $K_i\subset S^3$, take a ball $B_i$ in $S^3$ such that $B_i\cap K_i$ is a single unknotted arc, with $K_i$ meeting $\bdy B_i$ transversely in two points. That is, $(B_i,K_i)$ is homeomorphic to the product of an interval and a disk with a single marked point.

Now, remove $B_i$ from $S^3-K_i$. The result is homeomorphic to an arc in a ball. Glue $(S^3-K_1)-B_1$ to $(S^3-K_2)-B_2$ via a homeomorphism taking $(\bdy B_1, \bdy B_1\cap K_1)$ to $(\bdy B_2, \bdy B_2\cap K_2)$. (Here the notation means that $\bdy B_1$ is mapped to $\bdy B_2$ in such a way that the two points $\bdy B_1\cap K_1$ are mapped to the two points $\bdy B_2\cap K_2$.)
\end{definition}

\begin{definition}\label{Def:PrimeKnot}
A knot or link is said to be \emph{prime}\index{prime!knot or link} if it cannot be expressed as a connected sum of knots.\index{connected sum}\index{knot sum}
\end{definition}

The knot in the middle of \reffig{Prime} is a connected sum. Again from the point of view of hyperbolic geometry, knots that are not prime are not interesting, since they always contain an incompressible\index{incompressible} torus\index{incompressible torus} called a \emph{swallow--follow torus}\index{swallow--follow torus}. The torus is built by taking the boundary of one of the balls $\bdy B_1-N(K_1)$ in the construction of the connected sum, and then attaching a tube from one component of $N(K_1)$ on $\bdy B_1$ to the other, following $K_2$. This forms a torus which ``swallows'' $K_1$, then ``follows'' $K_2$. The right side of \reffig{Prime} shows a swallow--follow torus for the given example. 

Notice that prime diagrams and prime knots are not the same thing in general. Every knot, whether or not it is prime, admits a diagram that is not prime: simply insert a nugatory crossing. In general, if a knot admits a prime diagram, it still may not be a prime knot.

Recall the definitions of meridian and longitude of a knot or link.

\begin{definition}\label{Def:MeridianLongitude}
A curve on the boundary of a neighborhood of a knot or link in $S^3$ that bounds a disk inside the neighborhood of the link is called a \emph{meridian}\index{meridian}.

For a knot, a \emph{longitude}\index{longitude} is a curve on the boundary of a neighborhood of the knot that intersects a meridian exactly once. The \emph{standard longitude}\index{longitude!standard} is the longitude that is homologous to zero in $H_1(S^3-K)$.

More generally, a \emph{standard longitude} of a component $K_1$ of a link is the longitude that is trivial in $H_1(S^3-K_1)$. 
\end{definition}

\begin{lemma}\label{Lem:PrimeEquiv}
Let $K$ be a knot in $S^3$. Then $K$ is a connected sum\index{connected sum} of nontrivial knots if and only if $S^3-N(K)$ contains an essential\index{essential} annulus that meets the boundary of the neighborhood $N(K)$ in two simple meridians. 
\end{lemma}

We call such an annulus an \emph{essential meridional annulus}\index{essential meridional annulus}. 

\begin{proof}
Suppose first that $S^3-N(K)$ contains an essential annulus with boundary two meridians of $N(K)$. Cut $S^3$ along the sphere obtained from the union of this annulus and two disks bounded by the meridians. This separates $S^3$ into two balls $B_1$ and $B_2$, each containing an arc of $K$. Form new knots $K_1$ and $K_2$ by attaching to each $(\bdy B_i, B_i \cap K)$ a ball with an unknotted arc. By construction, $K$ is the connected sum\index{connected sum} of $K_1$ and $K_2$.

Now suppose that $K$ is a connected sum\index{connected sum} of nontrivial knots. Then $S^3-N(K)$ is obtained from nontrivial knots $K_1$ and $K_2$ by removing 3-balls $B_1$ and $B_2$ from $S^3-N(K_1)$ and $S^3-N(K_2)$, respectively, each meeting $\bdy N(K_i)$ transversely in two simple meridians. The result has boundary an annulus $A\cong B_i-N(K_i)$, and these annuli are glued to form $S^3-N(K)$. We claim $A$ is the essential annulus required. It meets $N(K)$ in meridians, as required. If it is compressible,\index{compressible} then a disk $D$ with boundary isotopic to the essential core curve of $A$ lies inside $(S^3-K_i)-B_i$ for one of $i=1,2$. Slice $\bdy B_i$ along $\bdy D$ to obtain a disk $E_i$ meeting $K_i$ exactly once. Attach to $E_i$ the disk $D$. This is a sphere meeting $K_i$ exactly once. But $K_i$ is a closed curve in $S^3$, hence it meets any sphere an even number of times. This contradiction proves that $A$ is incompressible.\index{incompressible}

Now suppose that there is a boundary compression disk\index{boundary compression disk} $D$ for $A$. An arc of $\bdy D$ must run from one (meridian) boundary component of $A$ to the other along $A$. The other arc of $\bdy D$ must run along $K$, either along $K_1$ or $K_2$, say $K_1$. But then $D$ can be used to isotope $K_1$ through $B_1$ to $\bdy B_1$, contradicting the fact that $K_1$ is nontrivial.

Finally, $A$ cannot be boundary parallel,\index{boundary parallel} else one side $(S^3-K_i)-B_i$ is homeomorphic to an unknotted arc in the ball $B_i$, again contradicting the fact that $K_i$ is nontrivial. This concludes the proof that a connected sum\index{connected sum} of knots contains an essential annulus with boundary two meridians of $N(K)$.
\end{proof}

\subsection{Polyhedral decomposition, revisited}

We know from the above discussion that alternating diagrams\index{alternating diagram} that are not connected and not prime cannot have hyperbolic complement. More generally, we need to determine which alternating diagrams lead to essential spheres, disks, tori, and annuli in the complement to rule out hyperbolicity. Our main tool will be a polyhedral decomposition of the link complement. 

Recall that in \refchap{Fig8Decomp} we worked through a decomposition of the figure-8 knot complement into two ideal polyhedra. This was extended in the exercises. In particular, following the methods of that chapter, the exercises outline a proof of the following theorem.

\begin{theorem}\label{Thm:PolyAltKnot}
  Let $L$ be an alternating link.\index{alternating knot or link}\index{alternating knot or link!polyhedral decomposition} Then the complement of $L$ can be obtained by gluing two ideal polyhedra that satisfying:
  \begin{enumerate}
  \item The polyhedra are obtained by labeling the boundary of two balls with the projection graph of the alternating diagram\index{alternating diagram} of $L$, and declaring each vertex to be ideal. On one ball, the outside boundary is labeled with the diagram, on the other the inside.
  \item Ideal vertices are 4-valent, corresponding to overcrossings in one polyhedron, undercrossings in the other.
  \item Ideal edges correspond to crossing arcs\index{crossing arc} in the diagram, and each edge class contains four ideal edges of the two polyhedra.
  \item Faces correspond to regions of the diagram, and are checkerboard colored, white and shaded.
  \item Each face on one polyhedron is glued to the identical face on the opposite polyhedron. The gluing rotates the face by one edge in the clockwise direction for white faces, and rotates by one edge in the counterclockwise direction for shaded faces. \qed
  \end{enumerate}
\end{theorem}

\begin{figure}
\includegraphics{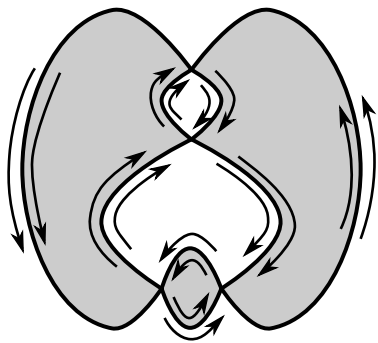}
  \caption{Ideal polyhedral decomposition of an alternating link (the figure-8 knot). Shown is one ideal polyhedron. The other is identical (with head on the opposite side) and gluing of faces is by a rotation in each face as shown}
  \label{Fig:PolyAltKnot}
\end{figure}

The theorem is illustrated in \reffig{PolyAltKnot}.

In the exercises of \refchap{Fig8Decomp}, we collapsed bigons\index{bigon} to a single edge. We will actually keep bigons around in this chapter, as they make certain arguments simpler. 

\subsection{Angled polyhedra and alternating links}

\begin{proposition}\label{Prop:AltAngleStruct}
Let $K$ be a knot or link with a connected, prime, alternating diagram.\index{alternating diagram}\index{alternating knot or link}\index{alternating knot or link!angled polyhedral structure}
If we assign a dihedral angle of $\pi/2$ to each ideal edge of the polyhedral decomposition of $S^3-K$ of \refthm{PolyAltKnot}, then we obtain an angled polyhedral structure,\index{angled polyhedral structure} as in \refdef{AnglePolyhedra}.
\end{proposition}

\begin{proof}
We need to check that the polyhedra with interior dihedral angles $\pi/2$ satisfy the three conditions of \refdef{AnglePolyhedra}. The first condition is immediate: $\pi/2$ lies in $(0,\pi)$. The third condition is also straightforward: each ideal edge of the ideal polyhedral decomposition appears exactly four times in the decomposition, hence interior angles sum to $2\pi$.

The second condition takes the most work. We need to show that every normal disk\index{normal} has non-negative combinatorial area.\index{combinatorial area} Recall the combinatorial area of a normal disk $D$ is defined to be
\[ a(D) = \sum_{i=1}^n (\pi-\alpha_i) - 2\pi + \pi|\bdy D \cap \bdy M|, \]
where $\alpha_1, \dots, \alpha_n$ are the dihedral angles met by $\bdy D$, and $|\bdy D\cap \bdy M|$ is the number of times $\bdy D$ meets a boundary face. In our case, each $\alpha_i=\pi/2$, so the sum is
\[ a(D) = \frac{\pi}{2}|\bdy D \cap e(M)| - 2\pi + \pi|\bdy D \cap \bdy M|,\]
where $|\bdy D\cap e(M)|$ is the number of times $\bdy D$ meets an ideal edge (not a boundary edge).

Notice that if $\bdy D$ meets at least four ideal edges, or at least two boundary faces, then the combinatorial area\index{combinatorial area} of $D$ is non-negative. The only possible ways it could be negative is if $\bdy D$ meets three or fewer edges and no boundary faces, or if it meets one boundary face and at most one edge. We rule these out.

First, suppose $\bdy D$ meets exactly one boundary face. The endpoints of the arc of $\bdy D$ on the boundary face must be on distinct boundary edges, by definition of a normal disk.\index{normal} So $\bdy D$ runs through at least two distinct regular faces of the polyhedron, and so $\bdy D$ must meet an edge of the polyhedron to connect into a closed curve. If $\bdy D$ meets only one edge, then it cannot meet an edge adjacent to the boundary face, by definition of normal.\index{normal} So $\bdy D$ encloses boundary faces on both sides. See \reffig{AltAngleStructProof}.
\begin{figure}[h]
%% Accompanies image file 'F11-03-1Ed1Bd.eps' (pdf, eps, ps)
\begingroup%
  \makeatletter%
  \providecommand\color[2][]{%
    \errmessage{(Inkscape) Color is used for the text in Inkscape, but the package 'color.sty' is not loaded}%
    \renewcommand\color[2][]{}%
  }%
  \providecommand\transparent[1]{%
    \errmessage{(Inkscape) Transparency is used (non-zero) for the text in Inkscape, but the package 'transparent.sty' is not loaded}%
    \renewcommand\transparent[1]{}%
  }%
  \providecommand\rotatebox[2]{#2}%
  \newcommand*\fsize{\dimexpr\f@size pt\relax}%
  \newcommand*\lineheight[1]{\fontsize{\fsize}{#1\fsize}\selectfont}%
  \ifx\svgwidth\undefined%
    \setlength{\unitlength}{214.19189072bp}%
    \ifx\svgscale\undefined%
      \relax%
    \else%
      \setlength{\unitlength}{\unitlength * \real{\svgscale}}%
    \fi%
  \else%
    \setlength{\unitlength}{\svgwidth}%
  \fi%
  \global\let\svgwidth\undefined%
  \global\let\svgscale\undefined%
  \makeatother%
  \begin{picture}(1,0.26341828)%
    \lineheight{1}%
    \setlength\tabcolsep{0pt}%
    \put(0,0){\includegraphics[width=\unitlength]{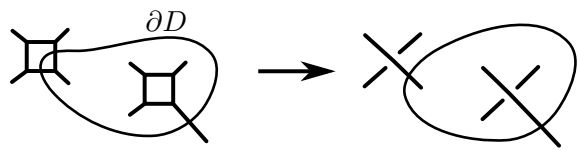}}%
    \put(0.23886396,0.20850665){\color[rgb]{0,0,0}\makebox(0,0)[lt]{\lineheight{1.25}\smash{\begin{tabular}[t]{l}$\bdy D$\end{tabular}}}}%
  \end{picture}%
\endgroup%

\caption{If $\bdy D$ meets just one edge and one boundary face, it determines a simple closed curve in the diagram of $K$ as shown}
\label{Fig:AltAngleStructProof}
\end{figure}
But recall that the graph of the polyhedron is exactly the diagram graph of $K$, and boundary faces correspond to crossings of $K$. Then $\bdy D$ gives a simple closed curve in the diagram of $K$ meeting a single crossing and a single strand of the link. We may slide $\bdy D$ off the crossing slightly so that it meets one more strand of the link near this crossing. Then $\bdy D$ is a closed curve in the link diagram meeting the diagram exactly twice transversely away from crossings, enclosing crossings on either side. This contradicts the fact that the diagram is prime.

Now suppose $\bdy D$ meets no boundary faces, but has negative combinatorial area.\index{combinatorial area} Then $\bdy D$ meets fewer than four edges, but more than zero edges by definition of normal.\index{normal} Because edges correspond to strands of the link, and the link consists of closed curves, it follows that $\bdy D$ meets exactly two ideal edges. Transferring $\bdy D$ to the diagram, it becomes a closed curve meeting the diagram exactly twice. But then because $K$ has a prime diagram, there are no crossings on one side of $\bdy D$. Transferring back to the polyhedron, this means an arc of $\bdy D$ meets the same edge of the polyhedron two times. This contradicts the definition of normal.\index{normal} 
\end{proof}

\begin{corollary}\label{Cor:AltIrreducible}
Let $K$ be a knot or link with a connected, prime, alternating diagram.\index{alternating diagram} Then $S^3-K$ is irreducible\index{irreducible} and boundary irreducible.\index{boundary irreducible}\index{alternating knot or link}
\end{corollary}

\begin{proof}
For such a knot or link, $S^3-K$ admits an angled polyhedral structure\index{angled polyhedral structure} by \refprop{AltAngleStruct}. Then the result follows from the first part of \refthm{HypAngleStruct}. 
\end{proof}

In fact, we may say more. The following is proved in \cite{menasco:alt}, using Thurston's \refthm{SfcesHyperbolic}; we also give a proof here.

\begin{theorem}\label{Thm:AltHyperbolic}
A knot with a connected prime alternating diagram\index{alternating diagram} is either a $(2,q)$-torus knot or it is hyperbolic. \index{alternating knot or link}\index{alternating knot or link!hyperbolic}
\end{theorem}

We have already shown alternating knots have complements that are irreducible\index{irreducible} and boundary irreducible.\index{boundary irreducible} To prove \refthm{AltHyperbolic} we need to consider essential annuli and tori, and we do so in the next subsections. 

%%%%%%%%%%%%%%%%%%%%%%%%%%%%%%%%%%%%%%%%%%%%%%%%%%%%%%%%%%%%%%%%%
\subsection{Alternating knots and essential annuli}

There are alternating knots that contain essential\index{essential} annuli, namely the $(2,q)$-torus knots. However, all other alternating knots are anannular. In this section, we prove that fact. In Menasco's original proof classifying hyperbolic alternating knots, he proves knots are anannular\index{anannular} by appealing to an algebraic result of Simon \cite{Simon:AlgKnot}. We take a more direct approach here, giving a geometric proof of this fact using the angled polyhedral structure\index{angled polyhedral structure} of the previous subsection.

First, we need more terminology to describe $(2,q)$-torus knots. The following definition is \refdef{TwReduced}, repeated here for convenience. 

\begin{figure}[h]
%% Accompanies image file 'F11-04-TwRed.eps' (pdf, eps, ps)
\begingroup%
  \makeatletter%
  \providecommand\color[2][]{%
    \errmessage{(Inkscape) Color is used for the text in Inkscape, but the package 'color.sty' is not loaded}%
    \renewcommand\color[2][]{}%
  }%
  \providecommand\transparent[1]{%
    \errmessage{(Inkscape) Transparency is used (non-zero) for the text in Inkscape, but the package 'transparent.sty' is not loaded}%
    \renewcommand\transparent[1]{}%
  }%
  \providecommand\rotatebox[2]{#2}%
  \newcommand*\fsize{\dimexpr\f@size pt\relax}%
  \newcommand*\lineheight[1]{\fontsize{\fsize}{#1\fsize}\selectfont}%
  \ifx\svgwidth\undefined%
    \setlength{\unitlength}{324bp}%
    \ifx\svgscale\undefined%
      \relax%
    \else%
      \setlength{\unitlength}{\unitlength * \real{\svgscale}}%
    \fi%
  \else%
    \setlength{\unitlength}{\svgwidth}%
  \fi%
  \global\let\svgwidth\undefined%
  \global\let\svgscale\undefined%
  \makeatother%
  \begin{picture}(1,0.18693556)%
    \lineheight{1}%
    \setlength\tabcolsep{0pt}%
    \put(0,0){\includegraphics[width=\unitlength]{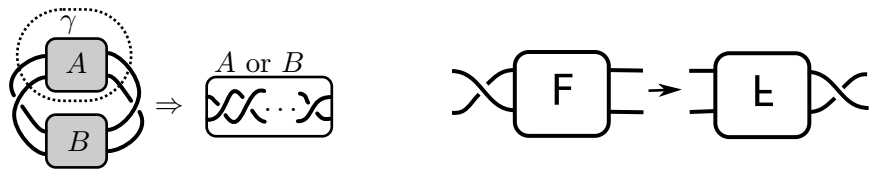}}%
    \put(0.06551547,0.11232076){\color[rgb]{0,0,0}\makebox(0,0)[lt]{\lineheight{1.25}\smash{\begin{tabular}[t]{l}$A$\end{tabular}}}}%
    \put(0.16957052,0.06613856){\color[rgb]{0,0,0}\makebox(0,0)[lt]{\lineheight{1.25}\smash{\begin{tabular}[t]{l}$\Rightarrow$\end{tabular}}}}%
    \put(0.06726921,0.02171105){\color[rgb]{0,0,0}\makebox(0,0)[lt]{\lineheight{1.25}\smash{\begin{tabular}[t]{l}$B$\end{tabular}}}}%
    \put(0.23855086,0.11348993){\color[rgb]{0,0,0}\makebox(0,0)[lt]{\lineheight{1.25}\smash{\begin{tabular}[t]{l}$A$ or $B$\end{tabular}}}}%
    \put(0.0590851,0.16493274){\color[rgb]{0,0,0}\makebox(0,0)[lt]{\lineheight{1.25}\smash{\begin{tabular}[t]{l}$\gamma$\end{tabular}}}}%
    \put(0.29610565,0.06819683){\color[rgb]{0,0,0}\makebox(0,0)[lt]{\lineheight{1.25}\smash{\begin{tabular}[t]{l}$\dots$\end{tabular}}}}%
  \end{picture}%
\endgroup%

  \caption{Left: A twist-reduced diagram. Right: A flype}
  \label{Fig:Flype}
\end{figure}

\begin{definition}\label{Def:TwistReduced}
  A diagram is \emph{twist-reduced}\index{twist-reduced} if, whenever $\gamma$ is a simple closed curve on the plane of projection meeting the diagram exactly twice in two crossings, running from one side of the crossing to the opposite side, the curve $\gamma$ bounds a string of bigons\index{bigon}
on one side. See \reffig{Flype}, left.
\end{definition}

\begin{definition}
  Let $\gamma$ be a simple closed curve meeting the diagram of $K$ transversely exactly four times in knot strands, with two intersections adjacent to a crossing on the outside of $\gamma$. A \emph{flype}\index{flype} is a move on the diagram that rotates the region inside $\gamma$ by $180^\circ$, moving the crossing outside $\gamma$ to lie between the opposite two strands. 
  See \reffig{Flype}, right. 
\end{definition}

\begin{lemma}\label{Lem:TwistReduced}
  Every knot or link $K$ has a twist-reduced diagram.

Moreover, if a diagram of $K$ is connected, prime, and alternating, then there is a twist-reduced diagram of $K$ that is connected, prime, and alternating.\index{alternating knot or link}
\end{lemma}

\begin{proof}
Start with a diagram of a knot or link. It has a finite number of twist regions, and a finite number of crossings in each twist region. Suppose the diagram is not twist reduced. Then there is a curve meeting the diagram exactly four times adjacent to two distinct twist regions. Slide the curve so that all crossings of both twist regions are on the outside of the curve, say with one twist region on the left and one on the right. Perform a sequence of flypes. Each flype will remove a crossing in the twist region on the left, and either add or remove a crossing in the twist region on the right (depending on the direction of crossings on the right and the direction of the flype). 
Continue until there are no crossings on the left. When finished, the diagram has one fewer twist region and at most the same number of crossings as before. Repeat, strictly reducing the number of twist regions. Since the number of twist regions is finite, the process will terminate in a twist-reduced diagram.

Finally, note that the process of flyping takes a connected diagram to a connected diagram. It also takes a prime diagram to a prime diagram and an alternating diagram\index{alternating diagram} to an alternating diagram (\refex{FlypePrimeAlt}). Thus if the original diagram of a link is connected, prime, and alternating, then the twist-reduced diagram, obtained by performing flypes, is also connected, prime, and alternating.
\end{proof}

\begin{definition}\label{Def:TwistNumber}
  The \emph{twist-number}\index{twist-number} of a knot diagram is the number of twist regions in a twist-reduced diagram.
\end{definition}

\begin{example}\label{Example:2qTorusKnot}
  A $(2,q)$-torus knot\index{torus knot!$(2,q)$-torus knot} has twist-number $1$. Any knot with a prime alternating diagram\index{alternating diagram}\index{alternating knot or link}\index{alternating knot or link!twist-number} that is not a $(2,q)$-torus knot has twist number at least $2$. In particular, the figure-8 knot shown in \reffig{Fig8Diagram} has twist number $2$. More generally, any twist knot $J(k,\ell)$ has twist number $2$. 
\end{example}

We are now ready to consider essential annuli. 

\begin{lemma}\label{Lem:Squares}
Suppose $K$ is a knot or link with a connected prime alternating diagram, with $S^3-N(K)$ given its (truncated ideal) polyhedral decomposition. Suppose $S$ is an essential\index{essential} annulus embedded in $S^3-N(K)$. Then when $S$ is isotoped into normal form,\index{normal} it contains at least one normal disk $D$ meeting a boundary face, and $\bdy D$ either meets exactly two boundary faces and no edges, or $\bdy D$ meets exactly one boundary face and exactly two edges. 
\end{lemma}

\begin{proof}
When we put $S$ into normal form,\index{normal form}\index{normal} \reflem{GaussBonnet} implies the combinatorial area\index{combinatorial area} of $S$ is $0$. Because each normal disk of $S$ has non-negative combinatorial area,\index{combinatorial area} in fact each normal disk of $S$ must have combinatorial area\index{combinatorial area} $0$. Because $S$ is a surface with boundary, there is at least one normal disk of $S$ that meets a boundary face; this is $D$. Now, considering the formula for the combinatorial area\index{combinatorial area} of $D$, there are only two possibilities: $\bdy D$ either meets exactly two boundary faces and no edges, or $\bdy D$ meets one boundary face and exactly two edges.
\end{proof}

\begin{lemma}\label{Lem:21Annulus}
  If $K$ has a prime alternating diagram, and $S$ is an embedded normal annulus\index{normal} properly embedded in the truncated polyhedral decomposition of $S^3-N(K)$, containing at least one normal disk $D_2$ whose boundary meets exactly one boundary face and exactly two edges of the polyhedra, then:
  \begin{enumerate}
  \item $S$ contains a subannulus $S'$ for which all normal disks\index{normal} meet exactly one boundary face and two edges.
  \item $K$ is a $(2,q)$-torus link with two components, and there is an annulus $\Sigma$ bounded by the two components of the link that is obtained by gluing bigon\index{bigon} faces of the polyhedral decomposition.
  \item A component of $\bdy S$ and $\bdy S'$ runs along at least one longitude of the link, so $\bdy S$ is not a meridian.
  \item The other component of $\bdy S'$ runs along the core of the annulus $\Sigma$.
  \end{enumerate}
\end{lemma}

\begin{proof}
Let $S$, $K$, and $D_2$ be as in the statement of the lemma. The disk $D_2$ is glued to normal disks\index{normal} $D_1$ and $D_3$ in the opposite polyhedron. The gluing maps a side of $D_2$ in a face to a side of $D_1$, and the gluing map on a face rotates the side either clockwise or counterclockwise. Without loss of generality, say clockwise. Thus a side of $D_2$ running from a boundary face to an edge is glued to a side of $D_1$ running from a boundary face to an edge, although rotated. Similarly for a side of $D_2$. See \reffig{21Curve}, left. 
Since $D_1$ and $D_3$ also have combinatorial area\index{combinatorial area} $0$, they must each meet one boundary face and exactly two edges. Repeat for disks meeting $D_1$ and $D_3$. Eventually this string of disks will glue up. Thus if there is one normal disk\index{normal} meeting a single boundary face and two edges, then there is a cycle of normal disks meet a single boundary face and two edges, gluing to form a subannulus $S'$ of $S$ as claimed. 

\begin{figure}[h]
%% Accompanies image file 'F11-05-21Crv.eps' (pdf, eps, ps)
\begingroup%
  \makeatletter%
  \providecommand\color[2][]{%
    \errmessage{(Inkscape) Color is used for the text in Inkscape, but the package 'color.sty' is not loaded}%
    \renewcommand\color[2][]{}%
  }%
  \providecommand\transparent[1]{%
    \errmessage{(Inkscape) Transparency is used (non-zero) for the text in Inkscape, but the package 'transparent.sty' is not loaded}%
    \renewcommand\transparent[1]{}%
  }%
  \providecommand\rotatebox[2]{#2}%
  \newcommand*\fsize{\dimexpr\f@size pt\relax}%
  \newcommand*\lineheight[1]{\fontsize{\fsize}{#1\fsize}\selectfont}%
  \ifx\svgwidth\undefined%
    \setlength{\unitlength}{308.50374985bp}%
    \ifx\svgscale\undefined%
      \relax%
    \else%
      \setlength{\unitlength}{\unitlength * \real{\svgscale}}%
    \fi%
  \else%
    \setlength{\unitlength}{\svgwidth}%
  \fi%
  \global\let\svgwidth\undefined%
  \global\let\svgscale\undefined%
  \makeatother%
  \begin{picture}(1,0.25476796)%
    \lineheight{1}%
    \setlength\tabcolsep{0pt}%
    \put(0,0){\includegraphics[width=\unitlength]{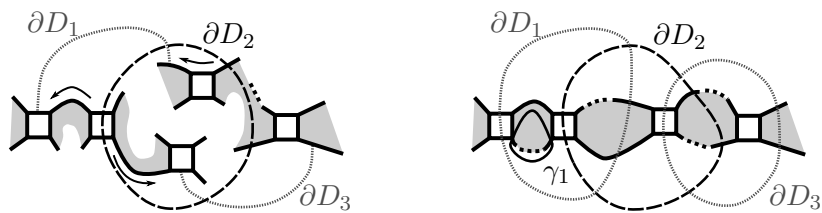}}%
    \put(0.23656555,0.21726209){\color[rgb]{0,0,0}\makebox(0,0)[lt]{\lineheight{1.25}\smash{\begin{tabular}[t]{l}$\bdy D_2$\end{tabular}}}}%
    \put(0.02291352,0.23138296){\color[rgb]{0.30196078,0.30196078,0.30196078}\makebox(0,0)[lt]{\lineheight{1.25}\smash{\begin{tabular}[t]{l}$\bdy D_1$\end{tabular}}}}%
    \put(0.35505652,0.01895875){\color[rgb]{0.30196078,0.30196078,0.30196078}\makebox(0,0)[lt]{\lineheight{1.25}\smash{\begin{tabular}[t]{l}$\bdy D_3$\end{tabular}}}}%
    \put(0.78731394,0.21630058){\color[rgb]{0,0,0}\makebox(0,0)[lt]{\lineheight{1.25}\smash{\begin{tabular}[t]{l}$\bdy D_2$\end{tabular}}}}%
    \put(0.58824845,0.23042158){\color[rgb]{0.30196078,0.30196078,0.30196078}\makebox(0,0)[lt]{\lineheight{1.25}\smash{\begin{tabular}[t]{l}$\bdy D_1$\end{tabular}}}}%
    \put(0.92898676,0.01983923){\color[rgb]{0.30196078,0.30196078,0.30196078}\makebox(0,0)[lt]{\lineheight{1.25}\smash{\begin{tabular}[t]{l}$\bdy D_3$\end{tabular}}}}%
    \put(0.6529279,0.05157117){\color[rgb]{0,0,0}\makebox(0,0)[lt]{\lineheight{1.25}\smash{\begin{tabular}[t]{l}$\gamma_1$\end{tabular}}}}%
  \end{picture}%
\endgroup%

\caption{Curve $\bdy D_2$ must be glued to arcs of $\bdy D_1$ and $\bdy D_3$ as shown on the left. Since $\bdy D_1$ and $\bdy D_3$ are disjoint, the only possibility for $\bdy D_1$ is that shown on the right. Then there is a curve $\gamma$ meeting the diagram exactly twice; it must bound an unknotted strand, forming a bigon.\index{bigon}}
  \label{Fig:21Curve}
\end{figure}

Sketch $\bdy D_2$ onto the boundary of the polyhedron, which has the combinatorics of the diagram graph. We will add to this picture $\bdy D_1$ and $\bdy D_3$ by superimposing, as in \reffig{21Curve}. By what we know of the gluing maps, an arc of $\bdy D_1$ must have its endpoints rotated once clockwise from an arc of $\bdy D_2$, as shown on the left of \reffig{21Curve}, and similarly for an arc of $\bdy D_3$. 

Because $D_1$ and $D_3$ are disjoint, $\bdy D_1$ must lie to one side of the arc of $\bdy D_3$ shown in \reffig{21Curve}, and thus $\bdy D_1$ and $\bdy D_3$ close up as shown on the right of that figure. Now inside of $\bdy D_1$, we may draw a curve $\gamma_1$ running through the shaded face between boundary faces met by $\bdy D_1$ and $\bdy D_2$, and through a single white face as in \reffig{21Curve}. This gives a curve meeting the diagram exactly twice. Because the diagram is prime, $\gamma_1$ bounds an arc of the diagram with no crossings on one side. Thus that shaded face is a simple bigon.\index{bigon} (Similar arguments show that other dotted lines in \reffig{21Curve}, right, are also single edges, but we will not use this.)

Repeating this argument with $D_1$ replacing $D_2$, and so on, we find that $S'$ is made up of disks bounding a closed chain of bigons. Thus the diagram of $K$ contains a single twist region, and $K$ is a $(2,q)$-torus knot or link. To see it is a link, note that disks $D_i$ for $i$ odd must all be disjoint, and disks $D_i$ for $i$ even must also be disjoint. If there are an odd number of bigons in the chain, this will be impossible. So there are an even number of bigons, $K$ is a 2-component link, and the surface $\Sigma$ made up of the bigons lies between the strands of $K$ and is an annulus. 

Now note that $\bdy D_1$ and $\bdy D_2$ together meet both ideal vertices on either side of a bigon face in the polyhedral decomposition. One of $D_1$, $D_2$ lies in one polyhedron, and one in the other. But then some $\bdy D_i$ will meet each ideal vertex in the diagram graph. It follows that $\bdy S$ meets each ideal vertex in each polyhedron at least once. This implies that $\bdy S$ runs along at least one longitude.

Finally, the arc of $\bdy D_i$ lying in a shaded face is a simple arc through the bigon. Thus it runs from one crossing arc\index{crossing arc} bounding the bigon to the other. When we glue all the disks $D_i$, the boundary of $S'$ traces the core of the annulus $\Sigma$.
\end{proof}

\begin{lemma}\label{Lem:22Annulus}
Suppose $K$ is a knot or link with a connected, twist-reduced, prime, alternating diagram. Suppose $S$ is an essential\index{essential} annulus in normal form\index{normal form} in the polyhedral decomposition of $S^3-K$ such that $S$ contains a normal disk\index{normal} whose boundary meets exactly two boundary faces and no ideal edges of the polyhedra. Then all normal disks of $S$ meet exactly two boundary faces, and the diagram of $K$ is that of a $(2,q)$-torus knot or link. Further, $\bdy S$ runs along at least one longitude of the knot or link, so $\bdy S$ is not a meridian.
\end{lemma}

\begin{proof}
Suppose $D_2$ is a normal disk\index{normal} of $S$ such that $\bdy D_2$ meets exactly two boundary faces and no edges. Then the fact that the diagram is prime and twist-reduced implies that $D_2$ is either a boundary bigon\index{boundary bigon} or $\bdy D_2$ encircles a portion of a twist region in the diagram.

Suppose first that all normal disks\index{normal} of $S$ that meet boundary faces are boundary bigons. Then by considering how such disks must glue, note that there can only be four disks, and they must encircle a single edge of the polyhedral decomposition. Thus $S$ is an annulus encircling a crossing circle. This contradicts the fact that $S$ is essential.

Now suppose all normal disks\index{normal} encircle portions of twist regions. Because these match up to form an annulus, the diagram of $K$ must consist only of a single twist region, and the knot is a $(2,q)$-torus knot or link. As in the proof of \reflem{21Annulus}, trace the boundary of $S$ in this case. Superimpose onto a single polyhedron to obtain a string of boundaries of normal squares, with every other normal square coming from the same polyhedron. Because normal squares in a polyhedron are disjoint, this forces each square to bound either a single bigon, or a pair of adjacent bigons;\index{bigon} see \refex{22Bigons}. Then $\bdy S$ meets every ideal vertex at least once, so $\bdy S$ is not a meridian. 

So finally suppose normal disks\index{normal} of $S$ consist both of curves encircling twist regions and boundary bigons.\index{boundary bigon} There must be a disk $D_2$ that is a boundary bigon adjacent to a disk $D_1$ encircling a portion of a twist region. Superimpose $\bdy D_1$ and $\bdy D_2$ on a single polyhedron. By following the gluing maps, we find $\bdy D_1$ bounds a single bigon face of the polyhedron, and $\bdy D_2$ bounds an edge sharing the same boundary face (ideal vertex) with the bigon. See \reffig{TwistBigon}, left.

\begin{figure}[h]
  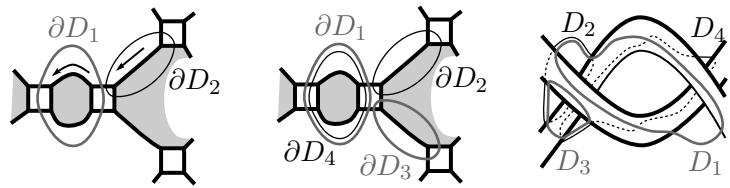
  \caption{Left: Boundary bigon\index{boundary bigon} adjacent to a disk encircling a portion of twist region must have the form shown on the left. Middle: two next disks must have the form shown. Right: sketch the disks in the 3-dimensional knot complement}
  \label{Fig:TwistBigon}
\end{figure}

There is another normal disk\index{normal} $D_3$ attached to $D_2$, in the opposite polyhedron from that containing $D_2$. Because $D_1$ and $D_3$ are disjoint, $D_3$ must be a boundary bigon\index{boundary bigon} of the form shown in \reffig{TwistBigon}. Then $D_4$ cannot be a boundary bigon, for if it is, $D_4$ is not glued to $D_1$ (its side is in the wrong region), thus $D_4$ is glued to another disk $D_5$. Because $D_5$ and $D_1$ are disjoint, $D_5$ must be a boundary bigon, and then some $D_6$ will also be a boundary bigon contained inside $D_2$, and so on, and there will be infinitely many boundary bigons spiraling around the same edge class. This is impossible. So $D_4$ bounds a portion of twist region, and $\bdy D_4$ is parallel to $\bdy D_1$ when superimposed (although recall that the disk $D_1$ lies in the opposite polyhedron from $D_4$). The disks $D_1$ through $D_4$ are shown superimposed on the same polyhedron in \reffig{TwistBigon}, middle, and in the link complement in \reffig{TwistBigon}, right.

We claim the annulus $S$ can be isotoped so that these four disks become two normal boundary bigons, and all other normal disks\index{normal} of $S$ are unchanged. The isotopy is by sliding past a crossing of the twist region where the boundary bigons\index{boundary bigon} cause the annulus to double back on itself. The isotopy is shown in \reffig{RemoveTwist}.

\begin{figure}[h]
\includegraphics{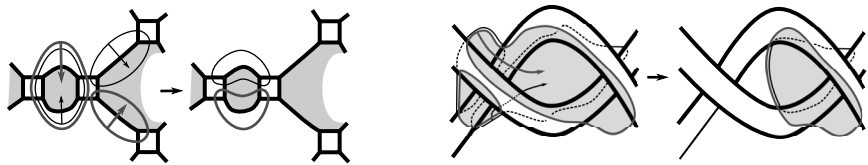}
  \caption{An isotopy of $S$ removes normal disks\index{normal} bounding portions of twist region. Shown on the left is the effect of the isotopy in the diagram graph. Shown on the right is the result of the isotopy in the link complement}
  \label{Fig:RemoveTwist}
\end{figure}

Repeating this move a finite number of times, we remove all disks bounding twist regions, and $S$ is made up only of boundary bigons.\index{boundary bigon} This is a contradiction.
\end{proof}

\begin{proposition}\label{Prop:AltAnannular}
  If $K$ is a knot or link with a connected, twist-reduced, prime, alternating diagram,\index{alternating diagram} and $K$ is not a $(2,q)$-torus knot or link, then $K$ is anannular.\index{anannular}\index{alternating knot or link}\index{alternating knot or link!anannular}

  If $K$ is a $(2,q)$-torus knot or link, then any essential annulus in $S^3-K$ has boundary tracing out at least one longitude. Thus there is no essential\index{essential} meridional annulus.\index{essential meridional annulus}
\end{proposition}

\begin{proof}
Suppose $S^3-N(K)$ contains an embedded essential annulus $S$. We may isotope it into normal form,\index{normal} and by \reflem{Squares} each normal disk making up $S$ either meets one boundary face and two ideal edges, or two boundary faces and no ideal edges. In the former case, \reflem{21Annulus} implies $K$ is a $(2,q)$-torus link and $S$ is not meridional. In the latter case, \reflem{22Annulus} implies $K$ is a $(2,q)$-torus knot or link and $S$ is not meridional. 
\end{proof}

\begin{corollary}\label{Cor:PrimeDiagram}
If $K$ has a connected prime alternating diagram,\index{alternating diagram}\index{alternating knot or link} then $K$ is a prime link. 
\end{corollary}

\begin{proof}
By \reflem{PrimeEquiv}, the link $K$ is not prime if and only if $S^3-N(K)$ contains an essential meridional annulus.\index{essential meridional annulus} By \reflem{TwistReduced}, $K$ has a diagram that is connected, prime, alternating, and twist-reduced. Then \refprop{AltAnannular} implies that $S^3-N(K)$ cannot contain an essential meridional annulus. So $K$ is prime. 
\end{proof}

%%%%%%%%%%%%%%%%%%%%%%%%%%%%%%%%%%%%%%%%%%%%%%%%%%%%%%%%%%%%%%%%%
\subsection{Closed surfaces and alternating knots}

Our goal is still to prove \refthm{AltHyperbolic}, that a knot with a connected, prime, alternating diagram\index{alternating diagram} is either a $(2,q)$-torus knot or is hyperbolic. We now consider closed essential surfaces embedded in $S^3-K$. 

\begin{lemma}\label{Lem:Meridional}
  Suppose $S$ is a closed essential\index{essential} surface embedded in the complement of a knot or link $K$ with a prime, connected, alternating diagram.\index{alternating diagram} Then $S$ contains a closed curve that encircles a meridian of $K$ at a crossing. 
\end{lemma}

\begin{proof}
Put $S$ into normal form\index{normal form} with respect to the polyhedral decomposition of $S^3-K$. Let $D$ be an innermost normal disk\index{normal} in the polyhedron; that is, $D$ cuts off a portion of a polyhedron that contains no other normal disks of $S$. Now, because $S$ is a closed surface, $D$ must meet a regular (i.e.\ not boundary) face $F$ of the polyhedron. Moreover, normality implies an arc of $\bdy D$ meets $F$ on two distinct edges $e_1$ and $e_2$ bordering $F$. These edges correspond to crossing arcs.\index{crossing arc} Recall from the construction of the polyhedra (e.g.\ in \refchap{Fig8Decomp}) that each such edge is identified to an edge on the opposite side of an ideal vertex in the polyhedron. Because the diagram is alternating, the two edges that are identified to $e_1$ and $e_2$ must lie on opposite sides of $D$; see \reffig{OppositeSides}.

\begin{figure}[h]
%% Accompanies image file 'F11-08-OppSid.eps' (pdf, eps, ps)
\begingroup%
  \makeatletter%
  \providecommand\color[2][]{%
    \errmessage{(Inkscape) Color is used for the text in Inkscape, but the package 'color.sty' is not loaded}%
    \renewcommand\color[2][]{}%
  }%
  \providecommand\transparent[1]{%
    \errmessage{(Inkscape) Transparency is used (non-zero) for the text in Inkscape, but the package 'transparent.sty' is not loaded}%
    \renewcommand\transparent[1]{}%
  }%
  \providecommand\rotatebox[2]{#2}%
  \newcommand*\fsize{\dimexpr\f@size pt\relax}%
  \newcommand*\lineheight[1]{\fontsize{\fsize}{#1\fsize}\selectfont}%
  \ifx\svgwidth\undefined%
    \setlength{\unitlength}{220.06720734bp}%
    \ifx\svgscale\undefined%
      \relax%
    \else%
      \setlength{\unitlength}{\unitlength * \real{\svgscale}}%
    \fi%
  \else%
    \setlength{\unitlength}{\svgwidth}%
  \fi%
  \global\let\svgwidth\undefined%
  \global\let\svgscale\undefined%
  \makeatother%
  \begin{picture}(1,0.41433597)%
    \lineheight{1}%
    \setlength\tabcolsep{0pt}%
    \put(0,0){\includegraphics[width=\unitlength]{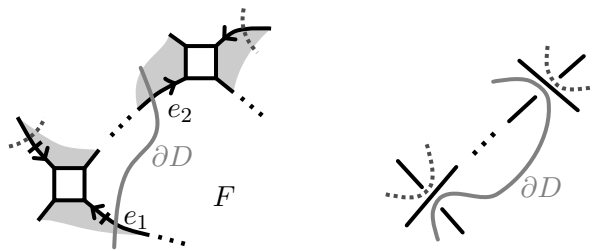}}%
    \put(0.24747163,0.14595957){\color[rgb]{0.50196078,0.50196078,0.50196078}\makebox(0,0)[lt]{\lineheight{1.25000012}\smash{\begin{tabular}[t]{l}$\bdy D$\end{tabular}}}}%
    \put(0.19944964,0.04689465){\color[rgb]{0,0,0}\makebox(0,0)[lt]{\lineheight{1.25000012}\smash{\begin{tabular}[t]{l}$e_1$\end{tabular}}}}%
    \put(0.2750515,0.23019119){\color[rgb]{0,0,0}\makebox(0,0)[lt]{\lineheight{1.25000012}\smash{\begin{tabular}[t]{l}$e_2$\end{tabular}}}}%
    \put(0.35204535,0.07806089){\color[rgb]{0,0,0}\makebox(0,0)[lt]{\lineheight{1.25000012}\smash{\begin{tabular}[t]{l}$F$\end{tabular}}}}%
    \put(0.87944212,0.09336866){\color[rgb]{0.50196078,0.50196078,0.50196078}\makebox(0,0)[lt]{\lineheight{1.25}\smash{\begin{tabular}[t]{l}$\bdy D$\end{tabular}}}}%
  \end{picture}%
\endgroup%

  \caption{Crossing arcs\index{crossing arc} identified to edges meeting $\bdy D$ lie on opposite sides of $D$. This is shown on the left in the polyhedron, and on the right in the diagram of the knot.}
  \label{Fig:OppositeSides}
\end{figure}

Because $D$ meets edges $e_1$ and $e_2$, another normal disk\index{normal} of $S$ in the same polyhedron must meet the opposite edges identified to $e_1$ and $e_2$, and thus there is an arc of a normal disk of $S$ on either side of $D$. In \reffig{OppositeSides}, these are shown as dashed arcs. But $D$ was chosen to be innermost, so one of those arcs must also belong to $D$. Then $D$ contains an arc running from one crossing arc on one side of an ideal vertex back to the identified crossing arc\index{crossing arc} on the other side of the ideal vertex. This arc glues up in $S$ to be a closed curve encircling a meridian at the crossing.
\end{proof}

\begin{corollary}\label{Cor:AltAtoroidal}
  If $K$ is a knot or link with a prime alternating diagram,\index{alternating diagram}\index{alternating knot or link}\index{alternating knot or link!atoroidal} then its complement is atoroidal.\index{atoroidal}
\end{corollary}

\begin{proof}
Suppose $S$ is an essential torus in $S^3-K$. Then $S$ contains a closed curve encircling a meridian at a crossing, by \reflem{Meridional}. This closed curve bounds a disk in $S^3$ that meets the knot exactly once in a meridian at a crossing. Surger along this disk and push both ends away from the crossing. We obtain a sphere $S'$ that meets the knot exactly twice in two meridians, with a crossing on the outside. Then $S'-N(K):=A$ is a meridional annulus. Because the link is prime, by \refcor{PrimeDiagram}, the annulus $A$ cannot be essential. It follows that $A$ is boundary parallel,\index{boundary parallel} and thus the original torus $S$ is boundary parallel,\index{boundary parallel} not essential.
\end{proof}

\begin{proof}[Proof of \refthm{AltHyperbolic}]
\Refcor{AltIrreducible} implies that a knot or link $K$ with a prime alternating diagram\index{alternating diagram}\index{alternating knot or link!hyperbolic} is irreducible\index{irreducible} and boundary irreducible.\index{boundary irreducible} \Refcor{AltAtoroidal} implies that it is atoroidal.\index{atoroidal} By \refprop{AltAnannular}, if it is not a $(2,q)$-torus knot, then it is also anannular.\index{anannular} The fact that $S^3-K$ is hyperbolic then follows from Thurston's \refthm{SfcesHyperbolic}. 
\end{proof}

%%%%%%%%%%%%%%%%%%%%%%%%%%%%%%%%%%%%%%%%%%%%%%%%%%%%%%%%%%%%%%%%%
\section{Checkerboard surfaces}\label{Sec:Checkerboard}

Recall that the diagram graph of a knot or link is a 4-valent graph embedded in the plane of projection. We may checkerboard color the regions of the graph, white and shaded. This checkerboard coloring\index{checkerboard coloring} may be used to define two surfaces embedded in any link complement.

\begin{figure}
\includegraphics{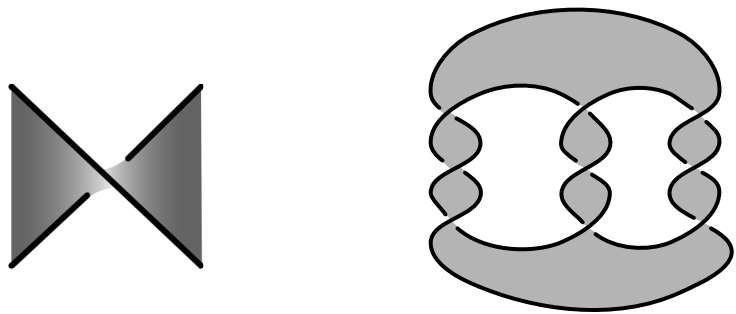}
  \caption{A twisted band is shown on the left, and a checkerboard surface on the right.\index{checkerboard surface}}
  \label{Fig:TwistedBand}
\end{figure}

\begin{definition}\label{Def:CheckerboardSurfaces}\index{checkerboard surface}\index{surface!checkerboard}
Let $K$ be a knot or link. Consider all the shaded regions in the checkerboard coloring\index{checkerboard coloring} of the complement of the diagram graph of $K$. By removing a neighborhood of each vertex, these can be embedded as disks in the link complement, with boundary lying on the knot and along crossing arcs.\index{crossing arc} At each crossing, attach a \emph{twisted band}\index{twisted band} between crossing arcs on opposite sides of the crossing; see \reffig{TwistedBand}. Note the twisting is in the same direction as the crossing. The result is a surface embedded in $S^3-N(K)$, with boundary on $N(K)$. This is called the \emph{shaded checkerboard surface}. The \emph{white checkerboard surface} is obtained similarly, using the opposite regions of the checkerboard coloring.\index{checkerboard coloring} 
\end{definition}

The main result of this section will be to show that if $K$ is alternating, then checkerboard surfaces are essential. Our main tool again will be the checkerboard polyhedral decomposition of an alternating link. 

For an alternating knot or link, the checkerboard surfaces are closely related to the polyhedral decomposition of \refthm{PolyAltKnot}. In that theorem, we obtained two polyhedra with checkerboard colored faces that glue to give the link complement. The shaded checkerboard surface is obtained from the shaded faces of the polyhedra, the white checkerboard surface from the white faces.

\begin{definition}\label{Def:Cut}
  Let $\Sigma$ be a properly embedded surface in a compact manifold $M$ with torus boundary components. Let $N(\Sigma)$ be a regular neighborhood of $\Sigma$. The manifold \emph{cut along $\Sigma$}\index{cut along surface $\Sigma$} is the manifold
  \[ M\cut\Sigma := M - N(\Sigma). \]
  The boundary of $M\cut\Sigma$ is a union of two subsurfaces. One of these is the surface $\bdy(N(\Sigma))\subset\bdy(M\cut\Sigma)$; it is homeomorphic to the double cover $\widetilde{\Sigma}$ of $\Sigma$. The other is the remnant of $\bdy M$, consisting of $\bdy M-(\bdy M\cap N(\Sigma))$, containing annuli and tori. The latter surface is called the \emph{parabolic locus}\index{parabolic locus} of $M\cut\Sigma$. 
\end{definition}

\begin{definition}\label{Def:BoundedPolyDecomp}
A \emph{bounded polyhedral decomposition} of a manifold $M\cut\Sigma$ is a decomposition of $M\cut\Sigma$ into truncated ideal polyhedra with interior and boundary faces, as well as \emph{surface faces},\index{surface faces} which are unglued, and which come from $\widetilde{\Sigma} \subset \bdy(M\cut\Sigma)$. As in \refdef{NormalDisk}, boundary edges are still defined to lie between boundary faces and other faces, \emph{interior edges}\index{interior edge} lie between pairs of interior faces, and \emph{surface edges}\index{surface edge} lie between surface faces and interior faces. We do not allow two surface faces to be adjacent along an edge. Moreover, under the gluing, each edge class either contains no surface edges, or it contains exactly two surface edges. 
\end{definition}

Our main example of a bounded polyhedral decomposition comes from checkerboard surfaces and alternating knots. 

\begin{lemma}\label{Lem:BoundedPolyAlt}
  Suppose $M=S^3-N(K)$ is the exterior of an alternating knot or link $K$,\index{alternating knot or link!checkerboard surface}\index{alternating knot or link!bounded polyhedral decomposition} and suppose $\Sigma$ is the shaded checkerboard surface. Then the cut manifold $M\cut\Sigma$ has a bounded polyhedral decomposition into the two checkerboard colored polyhedra of \refthm{PolyAltKnot}. Surface faces are shaded faces; boundary faces glue to form the parabolic locus of $M\cut\Sigma$. 
  A similar statement holds for the white checkerboard surface.\index{checkerboard surface}
\end{lemma}

\begin{proof}
  The decomposition is just as before, only in the gluing of the two polyhedra, leave the shaded faces unglued.
\end{proof}

There is a theorem for normal surfaces in bounded polyhedral decompositions that is completely analogous to \refthm{NormalForm}, for normal surfaces\index{normal} in ideal polyhedral decompositions. 

\begin{theorem}\label{Thm:BddNormalForm}
  Let $M\cut\Sigma$ have a bounded polyhedral decomposition. 
  \begin{enumerate}
  \item If $M$ is reducible,\index{reducible 3-manifold} then $M$ contains a normal 2-sphere.\index{normal}
  \item If $M$ is irreducible\index{irreducible} and boundary reducible, then $M$ contains a normal disk.\index{normal}
  \item If $M$ is irreducible and boundary irreducible,\index{boundary irreducible} then any essential\index{essential} surface in $M$ can be isotoped into normal form.\index{normal}
  \end{enumerate}
\end{theorem}

\begin{proof}
The proof is nearly identical to that of \refthm{NormalForm}, except we can no longer isotope an essential surface $S$ through surface faces, as they are now part of the boundary of $M\cut\Sigma$. We modify the proof of \refthm{NormalForm} where required to avoid such moves. Note that the proofs of the first two parts of the theorem require surgering, not isotoping, and so their arguments go through unchanged. So suppose $M$ is irreducible\index{irreducible} and boundary irreducible,\index{boundary irreducible} and $S$ is essential. 

First, if a component of $\bdy S$ lies entirely in a surface face and bounds a disk in that face, then since $S$ is incompressible,\index{incompressible} that curve bounds a disk in $S$ as well, hence $S$ has a disk component, parallel into a surface face, contradicting the fact that it is essential.

If an arc of intersection of $S$ with a face has both its endpoints on the same surface edge, and the arc lies in an interior face, then the arc and the edge bound a disc $D$ with one arc of $\bdy D$ on $S$ and one arc on a surface face. Because $S$ is essential, it is boundary incompressible;\index{boundary incompressible} it follows that the arc of intersection can be pushed off. A similar argument implies that an arc of intersection of $S$ with an interior face that has one endpoint on a boundary edge and one on a surface edge can be pushed off. For all other arcs of intersection with endpoints on one edge, or an edge and adjacent boundary edges, the argument follows just as before. 
\end{proof}

As in the case of polyhedral decompositions, we may put angled structures on bounded polyhedral decomposition and assign to normal disks\index{normal} and normal surfaces a combinatorial area,\index{combinatorial area} exactly as in \refdef{CombinatorialArea}.

\begin{definition}\label{Def:BoundedAngledPolyhedra}
  A \emph{bounded angled polyhedral structure}\index{bounded angled polyhedral structure}\index{angled polyhedral structure!bounded} is a decomposition of $M\cut\Sigma$ into ideal polyhedra, glued along interior faces, along with a collection of dihedral angles, one for each (surface or interior) edge, that satisfy the following.
  \begin{enumerate}
  \item Each dihedral angle lies in the range $(0,\pi)$. 
  \item Each normal disk\index{normal} in a polyhedron has nonnegative combinatorial area.\index{combinatorial area}
  \item Under the gluing, dihedral angles sum to $2\pi$ around an edge class meeting no surface edges. They sum to $\pi$ if they meet surface edges. 
  \end{enumerate}
\end{definition}

\begin{proposition}\label{Prop:BddAngledPolyAlt}
Suppose $M=S^3-N(K)$ is the exterior of an alternating knot or link\index{alternating knot or link} $K$, and $\Sigma$ is the shaded checkerboard surface. Then the cut manifold $M\cut\Sigma$ has a bounded angled polyhedral structure.\index{angled polyhedral structure!bounded}\index{bounded angled polyhedral structure}\index{checkerboard surface}
\end{proposition}

\begin{proof}
As in \refprop{AltAngleStruct}, label each ideal edge of the checkerboard polyhedra by $\pi/2\in(0,\pi)$. Then the proof of \refprop{AltAngleStruct} carries through to show that every normal disk\index{normal} has non-negative combinatorial area.\index{combinatorial area} We only need to check that dihedral angles sum to $\pi$ at surface edges. Note that because each ideal edge is adjacent to both white and shaded faces, in fact each ideal edge is a surface edge. Because we no longer glue shaded faces, each edge class contains exactly two surface edges. Thus the sum of dihedral angles at each edge is $\pi/2+\pi/2 =\pi$, as required. 
\end{proof}

\begin{lemma}[Bounded Gauss--Bonnet]\label{Lem:BddGaussBonnet}
Let $S$ be a surface properly embedded in $M\cut\Sigma$, in normal form\index{normal form}\index{normal} with respect to a bounded angled polyhedral structure\index{angled polyhedral structure!bounded}\index{bounded angled polyhedral structure} on $M\cut\Sigma$. Let $p$ denote the number of times $\bdy S$ intersects a boundary edge adjacent to a surface face. Then 
\[ a(S) = -2\pi\chi(S) + \frac{\pi}{2}\,p. \]
\end{lemma}

\begin{proof}
Exercise. 
\end{proof}

\begin{definition}\label{Def:BdryPi1Injective}
  Let $S$ be a surface properly embedded in a compact 3-manifold $M$ with boundary. We say $S$ is \emph{boundary $\pi_1$-injective}\index{boundary $\pi_1$-injective} if whenever $\alpha \subset S$ is an arc properly embedded in $S$ that is not homotopic rel endpoints to $\bdy S$ in $S$, then $\alpha$ is not homotopic rel endpoints to $\bdy M$ inside $M$.

We say the surface $S$ is \emph{$\pi_1$-essential}\index{$\pi_1$-essential} if
it is $\pi_1$-injective,\index{$\pi_1$-injective} boundary $\pi_1$-injective, and not parallel into $\bdy M$. 
\end{definition}

Note that boundary $\pi_1$-injective is stronger than boundary incompressible,\index{boundary incompressible} and $\pi_1$-essential is stronger than essential.\index{essential} For checkerboard surfaces, we have this stronger result.

\begin{theorem}\label{Thm:CheckerboardEssential}
  Let $K$ be a link with a connected, prime, reduced alternating diagram,\index{alternating diagram} and let $\Sigma$ be one of its checkerboard surfaces. Then $\Sigma$ is $\pi_1$-injective\index{$\pi_1$-injective} and boundary $\pi_1$-injective,\index{boundary $\pi_1$-injective} hence it is $\pi_1$-essential.\index{alternating knot or link!checkerboard surface}\index{checkerboard surface!$\pi_1$-essential}
\end{theorem}

\begin{proof}
We claim first that $\Sigma$ is $\pi_1$-injective and boundary $\pi_1$-injective if and only if the surface $\widetilde{\Sigma} = \bdy N(\Sigma)$ is incompressible\index{incompressible} and boundary incompressible.\index{boundary incompressible} The proof uses the loop theorem; we leave it as an exercise.

Now we claim that if $D$ is a compression disk\index{compression disk} for $\widetilde{\Sigma}$ in $S^3-N(K)$, then we may assume $D$ is properly embedded in $(S^3-N(K))\cut\Sigma$. This is because $\widetilde{\Sigma}$ separates $S^3-N(K)$ into $N(\Sigma)$ and $(S^3-N(K))\cut\Sigma$. An innermost disk argument implies that $D$ can be isotoped to be disjoint $\widetilde{\Sigma}$ in its interior, so $D$ either lies in $(S^3-N(K))\cut\Sigma$, as desired, or in the product $N(\Sigma)$. If $D$ is in the product, then an isotopy mapping $N(\Sigma)$ to $\Sigma$ takes the disk $D$ to a disk parallel to $\Sigma$, hence parallel to $\widetilde{\Sigma}$, contradicting the fact that it is a compression disk\index{compression disk} for $\widetilde{\Sigma}$.

Thus $D$ is an essential disk in $(S^3-N(K))\cut\Sigma$ with boundary completely contained in $\widetilde{\Sigma}$. Put $D$ into normal form\index{normal form} with respect to the bounded polyhedral decomposition of $(S^3-N(K))\cut\Sigma$. Because $\bdy D$ meets no boundary faces, \reflem{BddGaussBonnet}, the bounded Gauss--Bonnet lemma, implies that $a(D)=-2\pi$. But each normal disk making up $D$ has nonnegative combinatorial area,\index{combinatorial area} by \refprop{BddAngledPolyAlt}. This is a contradiction. 

Now suppose that $D$ is a boundary compression disk\index{boundary compression disk} for $\widetilde{\Sigma}$. An innermost disk and outermost arc argument implies that $D$ is isotopic to a disk with interior disjoint from $\widetilde{\Sigma}$, and again this disk must lie in $(S^3-N(K))\cut\Sigma$. Put the disk into normal form.\index{normal form} One arc of $\bdy D$ lies on $\widetilde{\Sigma}$ and one arc lies on boundary faces. Note this arc begins and ends on boundary edges adjacent to a surface face, but if it meets any other boundary edges in its interior they must be adjacent to interior edges. Then the bounded Gauss--Bonnet lemma, \reflem{BddGaussBonnet} implies that $a(D) = -2\pi + \pi = -\pi$. Again this contradicts the fact that normal disks\index{normal} have nonnegative combinatorial area.\index{combinatorial area}

So $\Sigma$ is $\pi_1$-injective\index{$\pi_1$-injective} and boundary $\pi_1$-injective.\index{boundary $\pi_1$-injective} Because it has boundary on $N(K)$ it cannot be boundary parallel.\index{boundary parallel} So it is $\pi_1$-essential.\index{$\pi_1$-essential}
\end{proof}

By \refthm{CheckerboardEssential}, every alternating knot contains a pair of $\pi_1$-essential checkerboard surfaces.\index{alternating knot or link!checkerboard surface} The converse is also true: independently, Howie \cite{Howie:Alternating} and Greene \cite{Greene:Alternating} showed that if a 3-manifold contains a pair of essential surfaces satisfying certain conditions required of checkerboard surfaces, then the 3-manifold is the complement of an alternating knot in $S^3$ and the surfaces are isotopic to checkerboard surfaces. \index{checkerboard surface!$\pi_1$-essential}

In addition to being $\pi_1$-essential,\index{$\pi_1$-essential} checkerboard surfaces of a hyperbolic alternating knot also exhibit other nice geometric properties. We discuss these in the next chapter. 

%%%%%%%%%%%%%%%%%%%%%%%%%%%%%%%%%%%%%%%%%%%%%%%%%%%%%%%%%%%%%%%%%
\section{Exercises}

\begin{exercise}
Give a proof that a knot or link with a connected, prime, alternating diagram\index{alternating diagram} is a prime link.

One way to prove this is to note that essential meridional annuli\index{essential meridional annulus} can be put into normal form\index{normal form} while ensuring boundary components of the annuli stay well-behaved, and then analyzing the number and form of normal disks that could possibly arise. 
\end{exercise}

\begin{exercise}\label{Ex:FlypePrimeAlt}
(Flypes and alternating diagrams)
\begin{enumerate}
\item Prove that a flype takes a prime diagram to a prime diagram.
\item Prove that a flype takes an alternating diagram\index{alternating diagram} to an alternating diagram.
\end{enumerate}
\end{exercise}

\begin{exercise}\label{Ex:NoNormalBigons}
Prove the following result, which will be used in \refchap{Quasifuchsian}.

\begin{proposition}\label{Prop:NoNormalBigons}
Let $K$ be a knot or link with a connected, prime, alternating diagram.\index{alternating diagram} Then in the polyhedral decomposition of the link complement, there can be no bigon\index{bigon} in normal form.\index{normal form} That is, there is no normal disk\index{normal} embedded in a polyhedron that meets exactly two interior edges. 
\end{proposition}
\end{exercise}

\begin{exercise}\label{Ex:22Bigons}
Work through the details of the proof of \reflem{22Annulus}: Suppose $K$ is a knot or link with a connected, twist-reduced, prime, alternating diagram. Suppose $S$ is an essential\index{essential} annulus in normal form with respect to the polyhedral decomposition. Suppose some normal disk\index{normal} making up $D_2$ meets exactly two boundary faces, but runs through opposite sides of those boundary faces.
\begin{enumerate}
\item[(a)] Show that $\bdy D_2$ encircles a twist region of the diagram. 
\item[(b)] Assume that $\bdy D_2$ runs through exactly two boundary faces and exactly two white faces. One arc of $\bdy D_2$ in a white face is glued to the side of a normal disk\index{normal} $D_2$, and the other arc of $\bdy D_2$ in a white face is glued to the side of a normal disk $D_3$. Following the example of \reffig{21Curve} left, sketch the images of these arcs of $\bdy D_1$ and $\bdy D_3$ superimposed on the same polyhedron containing $\bdy D_2$.
\item[(c)] Prove that $\bdy D_1$ and $\bdy D_3$ must each encircle a string of adjacent bigons.\index{bigon}
\item[(d)] Prove that in fact, $\bdy D_1$, $\bdy D_2$, and $\bdy D_3$ either all encircle a single bigon each, or they all encircle a pair of bigons\index{bigon} each. Sketch these curves superimposed on the same polyhedron. 
\end{enumerate}
\end{exercise}

\begin{exercise}
  Prove \reflem{BddGaussBonnet}, the bounded Gauss--Bonnet theorem.
\end{exercise}

\begin{exercise}
  Prove that a properly embedded surface $S$ in a compact 3-manifold $M$ is $\pi_1$-injective\index{$\pi_1$-injective} if and only if the surface $\widetilde{S} = \bdy N(S)$ is incompressible.\index{incompressible}
\end{exercise}

\begin{exercise}
  Prove that a properly embedded surface $S$ in a compact 3-manifold $M$ is boundary $\pi_1$-injective\index{boundary $\pi_1$-injective} if and only if the surface $\widetilde{S} = \bdy N(S)$ is boundary incompressible.\index{boundary incompressible}
\end{exercise}

\chapter{The Geometry of Embedded Surfaces}\label{Chap:Quasifuchsian}
\blfootnote{Jessica S. Purcell, Hyperbolic Knot Theory}

In this chapter, we discuss the geometry of essential surfaces\index{essential} embedded in hyperbolic 3-manifolds.

In the first section, we show that specific surfaces embedded in hyperbolic 3-manifolds always admit isometries. This allows us to cut along such surfaces and reglue, obtaining new manifolds whose geometry can be understood from the geometry of the original. The most straightforward instance of this uses the 3-punctured sphere, and was discovered by Adams \cite{adams:3-punct}, building on work of Wielenberg \cite{Wielenberg}. Ruberman discovered a similar result for 4-punctured spheres and related surfaces \cite{Ruberman:mutation}. Both techniques are still frequently used to build examples of hyperbolic knots and links with particular geometric properties (for example volume: \cite{SBurton}, \cite{adamsetal:VolDet}, short geodesics: \cite{Millichap}, cusp shapes: \cite{DangPurcell}).

We then return to more general essential\index{essential} surfaces, and discuss a geometric classification of such surfaces as quasifuchsian (or Fuchsian), accidental, or virtual fibered. We illustrate the behavior of such surfaces using examples from knot complements, especially alternating knots.\index{alternating knot or link} We show that for hyperbolic alternating links, their checkerboard surfaces are always quasifuchsian.

\section{Belted sums and mutations}

This section describes two techniques for building distinct links with related hyperbolic structures; for example they have the same volume. The techniques both arose in the 1980s by cutting along an embedded surface in a hyperbolic 3-manifold and regluing via isometry.

\subsection{3-punctured spheres and belted sums}

Suppose $M$ is a hyperbolic 3-manifold that contains an embedded incompressible\index{incompressible} 3-punctured sphere $S$. We have seen examples of this: in \refchap{TwistKnots}, each crossing circles of a reduced fully augmented link\index{fully augmented link} bounds an embedded essential\index{essential} 3-punctured sphere. In \refcor{AugmentedGeodSfces} we noted that these 3-punctured spheres are always totally geodesic in fully augmented links.\index{fully augmented link}
We now generalize this.

\begin{theorem}\label{Thm:3PunctTotallyGeodesic}
Let $M$ be a 3-manifold admitting a complete, finite volume hyperbolic structure, so $M$ is the interior of a compact manifold $\overline{M}$ with torus boundary. Let $S$ be a $\pi_1$-injective (equivalently, incompressible)\index{incompressible}
3-punctured sphere properly embedded in $\overline{M}$. Then $S$ is isotopic to a properly embedded 3-punctured sphere that is totally geodesic in the hyperbolic structure on $M$. 
\end{theorem}

\begin{proof}
Let $\alpha$, $\beta$, and $\gamma = \alpha\cdot \beta$ be generators of $\pi_1(S)$ that encircle the three punctures of $S$. Because $M$ admits a complete hyperbolic structure, there is a representation $\rho\from \pi_1(M)\to \PSL(2,\CC)$ taking $\alpha$, $\beta$, and $\gamma$ to parabolic elements.\index{parabolic}
We may conjugate to adjust the images of three points at infinity; we conjugate so that the fixed point of $\rho(\alpha)$ is $\infty$, so that $\rho(\alpha)$ translates $0\in\CC$ to $2\in\CC$, and so that the fixed point of $\rho(\beta)$ is $0$. Then the three parabolics have the form
\[ \rho(\alpha) = \mat{1 & 2 \\ 0 & 1}, \quad \rho(\beta) = \mat{1 & 0 \\ z & 1}, \quad \mbox{and } \rho(\alpha\cdot\beta) = \mat{ 1+2z & 2 \\ z & 1}. \]
Because $\rho(\alpha\cdot\beta)$ is parabolic, its trace is $2+2z = \pm 2$, so $z=0$ or $z=-2$. If $z=0$, $\rho(\beta)$ is the identity, contradicting the fact that $S$ is $\pi_1$-injective. Thus $z=-2$.

Now note that both $\rho(\alpha)$ and $\rho(\beta)$ (and hence $\rho(\alpha\cdot\beta)$) preserve the real line $\RR \subset \CC\subset \bdy_\infty \HH^3$. Hence $\rho(\pi_1(S))$ preserves the vertical plane $P$ in $\HH^3$ whose boundary is the real line. Thus under the covering map $p\from \HH^3\to \HH^3/\rho(\pi_1(M)) = M$, the plane $P$ maps to a totally geodesic surface homeomorphic to $S$ in $M$.

It remains to show that $p(P)$ is embedded in $M$, and $S$ is isotopic to the embedded totally geodesic surface $p(P)$. 
To do so, consider $p^{-1}(S)$. This is a disjoint union of embedded, possibly non-geodesic planes in $\HH^3$. One lift $\widetilde{S}$ is fixed by $\rho(\alpha)$, $\rho(\beta)$, and $\rho(\gamma)$ above. It follows that $\widetilde{S}$ and $P$ have the same limit set;\index{limit set} recall limit set is defined in \refdef{LimitSet}. Then for any $\gamma\in\pi_1(M)$, $\rho(\gamma)(P)$ has the same limit set as $\rho(\gamma)(\widetilde{S})$. Because $S$ is embedded, translates of $\widetilde{S}$ are disjoint, and it follows that translates of $P$ are disjoint embedded planes in $\HH^3$. Then $p(P)$ is an embedded surface in $M$, an isotopy from $\widetilde{S}$ to $P$ projects to an isotopy from $S$ to $p(P)$ in $M$, and $S$ is isotopic to a properly embedded totally geodesic 3-punctured sphere. 
\end{proof}

\begin{corollary}\label{Cor:Gluing3PunctSpheres}
  Let $M$ and $M'$ be hyperbolic 3-manifolds containing essential\index{essential} embedded 3-punctured spheres $S$ and $S'$, respectively. Then $M\cut S$ and $M'\cut S'$ are hyperbolic 3-manifolds, each with two totally geodesic 3-punctured sphere boundary components. Moreover:
  \begin{enumerate}
  \item Any manifold $M''$ obtained by identifying 3-punctured sphere boundary components of $M\cut S$ to those of $M'\cut S'$ will be hyperbolic, containing embedded essential 3-punctured spheres, and $M\cut S$ and $M'\cut S'$ embed isometrically in $M''$. In particular, $\vol(M'') = \vol(M)+\vol(M')$.
  \item Any manifold $M'''$ obtained by identifying the 3-punctured sphere boundary components of $M\cut S$ via homeomorphism will be hyperbolic, containing an embedded essential 3-punctured sphere, and $\vol(M''')=\vol(M)$.
  \end{enumerate}
\end{corollary}

\begin{proof}
Because $S$ and $S'$ are isotopic to totally geodesic hyperbolic surfaces, we obtain $M\cut S$ and $M'\cut S'$, respectively, by removing a collection of half spaces from $\HH^3$ corresponding to lifts of $S$ and $S'$, and then taking the quotient. Note that the geometry of $M\cut S$ and $M'\cut S'$ therefore agree with geometry of $M$ and $M'$, respectively, away from $S$ and $S'$. In particular, $\vol(M\cut S) = \vol(M)$ and $\vol(M'\cut S')=\vol(M')$. 
Moreover, $M\cut S$ and $M'\cut S'$ have totally geodesic 3-punctured sphere boundary.

By \refprop{CompleteStruct3punct}, there is a unique hyperbolic structure on a 3-punctured sphere. Therefore, any gluing of 3-punctured spheres can be obtained by isometry. So $M''$ is obtained by gluing $M\cut S$ to $M'\cut S'$ by isometry. Thus $M\cut S$ and $M'\cut S'$ isometrically embed in $M''$, and
the volume of $M''$ is equal to $\vol(M) + \vol(M')$.

Similarly, $M\cut S$ isometrically embeds in $M'''$ and $\vol(M''')=\vol(M)$. 
\end{proof}

\Refthm{3PunctTotallyGeodesic} and \refcor{Gluing3PunctSpheres} were used by Adams to construct explicit examples of links in $S^3$ with additive volumes.

\begin{figure}[h]
%% Accompanies image file 'F12-01-BeltSu.eps' (pdf, eps, ps)
\begingroup%
  \makeatletter%
  \providecommand\color[2][]{%
    \errmessage{(Inkscape) Color is used for the text in Inkscape, but the package 'color.sty' is not loaded}%
    \renewcommand\color[2][]{}%
  }%
  \providecommand\transparent[1]{%
    \errmessage{(Inkscape) Transparency is used (non-zero) for the text in Inkscape, but the package 'transparent.sty' is not loaded}%
    \renewcommand\transparent[1]{}%
  }%
  \providecommand\rotatebox[2]{#2}%
  \newcommand*\fsize{\dimexpr\f@size pt\relax}%
  \newcommand*\lineheight[1]{\fontsize{\fsize}{#1\fsize}\selectfont}%
  \ifx\svgwidth\undefined%
    \setlength{\unitlength}{270.6264782bp}%
    \ifx\svgscale\undefined%
      \relax%
    \else%
      \setlength{\unitlength}{\unitlength * \real{\svgscale}}%
    \fi%
  \else%
    \setlength{\unitlength}{\svgwidth}%
  \fi%
  \global\let\svgwidth\undefined%
  \global\let\svgscale\undefined%
  \makeatother%
  \begin{picture}(1,0.2270753)%
    \lineheight{1}%
    \setlength\tabcolsep{0pt}%
    \put(0,0){\includegraphics[width=\unitlength]{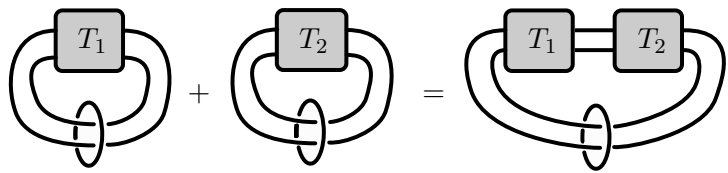}}%
    \put(0.0959329,0.1724336){\color[rgb]{0,0,0}\makebox(0,0)[lt]{\lineheight{1.25}\smash{\begin{tabular}[t]{l}$T_1$\end{tabular}}}}%
    \put(0.40434242,0.17290668){\color[rgb]{0,0,0}\makebox(0,0)[lt]{\lineheight{1.25}\smash{\begin{tabular}[t]{l}$T_2$\end{tabular}}}}%
    \put(0.24905687,0.09909149){\color[rgb]{0,0,0}\makebox(0,0)[lt]{\lineheight{1.25}\smash{\begin{tabular}[t]{l}$+$\end{tabular}}}}%
    \put(0.5760196,0.09909149){\color[rgb]{0,0,0}\makebox(0,0)[lt]{\lineheight{1.25}\smash{\begin{tabular}[t]{l}$=$\end{tabular}}}}%
    \put(0.72304732,0.1724336){\color[rgb]{0,0,0}\makebox(0,0)[lt]{\lineheight{1.25}\smash{\begin{tabular}[t]{l}$T_1$\end{tabular}}}}%
    \put(0.87547146,0.1724336){\color[rgb]{0,0,0}\makebox(0,0)[lt]{\lineheight{1.25}\smash{\begin{tabular}[t]{l}$T_2$\end{tabular}}}}%
  \end{picture}%
\endgroup%

  \caption{A belted tangle\index{belted tangle} is a link with an unknotted component bounding an embedded 2-punctured disk; two are shown on the left. On the right, a belted sum\index{belted sum} is obtained from two belted tangles via the tangle addition shown.}
  \label{Fig:BeltedTangle}
\end{figure}

\begin{definition}\label{Def:BeltSum}
  A \emph{belted tangle}\index{belted tangle} is a link in $S^3$ with one link component unknotted in $S^3$, bounding a disk meeting other components of the link exactly two times; see \reffig{BeltedTangle}, left. 
  
  The \emph{belted sum}\index{belted sum} of two belted tangles is the belted tangle obtained by tangle addition as in \reffig{BeltedTangle}, right.
\end{definition}

\begin{corollary}\label{Cor:BeltedSum}
  If $L_1$ and $L_2$ are belted tangles\index{belted tangle} that are hyperbolic, then their belted sum\index{belted sum} $L$ is a hyperbolic link with volume satisfying $\vol(L) = \vol(L_1)+\vol(L_2)$.
\end{corollary}

\begin{proof}
Note that $S^3-L_1$ and $S^3-L_2$ each contain an embedded 3-punctured sphere, namely the 2-punctured disk whose boundary is on the unknotted component of the link. Because these two link complements are hyperbolic, the 3-punctured sphere must be incompressible\index{incompressible} in both cases (easy exercise). The result then follows from \refcor{Gluing3PunctSpheres}.
\end{proof}

\subsection{4-punctured spheres and mutation}
Note that to prove \refcor{BeltedSum}, we glued isometric 3-punctured spheres. Unfortunately, the 3-punctured sphere is the only hyperbolic surface with a unique hyperbolic structure. All others have infinitely many hyperbolic structures, and so a gluing homeomorphism will not necessarily give an isometry, and volume will not necessarily be additive. However, in certain cases we may still cut and glue along an essential surface and still ensure that geometry is well behaved. One way to do this is a process called \emph{mutation}, which applies to 4-punctured spheres.

This section gives a condition on the diagrams of two knots and links that will guarantee that the geometries of their complements are similar; in particular they will have the same hyperbolic volume. This was first discovered by Ruberman~\cite{Ruberman:mutation}, who proved the result using minimal surfaces. Because the full proof requires more background on minimal surfaces than we wish to include here, we will refer to his paper for the complete result. However, we will provide full details in the special case that an embedded essential 4-punctured sphere is isotopic to an embedded pleated\index{pleated surface} 4-punctured sphere isometric to the boundary of an ideal tetrahedron. 

\begin{definition}\label{Def:ConwaySphere}
A \emph{Conway sphere}\index{Conway sphere} is a 4-punctured sphere obtained from the diagram of a knot or link $K$ as follows. Let $\gamma$ be a simple closed curve in the plane of projection of the diagram of $K$ that meets the diagram exactly four times, transversely in edges of the diagram. Let $\overline{S}$ be the sphere embedded in $S^3-K$ obtained by attaching two disks to $\gamma$, one on either side of the plane of projection. Let $S$ denote the corresponding 4-punctured sphere in $S^3-K$. In the case that $S$ is essential,\index{essential} we say that it is a \emph{Conway sphere}\index{Conway sphere} for $K$. 
\end{definition}

We will put geometric structures on Conway spheres, and cut and reglue via isometry of the spheres. Note that a hyperbolic ideal tetrahedron has boundary a pleated\index{pleated surface} 4-punctured sphere. This gives us a special case of a hyperbolic structure on a 4-punctured sphere and an isometry preserving it. 

\begin{lemma}\label{Lem:IsomTetr}
  Let $T$ be a hyperbolic ideal tetrahedron. For each pair of opposite edges, there is an axis in $\HH^3$ meeting the two edges orthogonally. Rotation by $\pi$ about such an axis is an isometry of the ideal tetrahedron. 
\end{lemma}

\begin{proof}
Rotation by $\pi$ through an axis is an isometry of $\HH^3$ that maps the tetrahedron back to itself. 
\end{proof}

\begin{corollary}\label{Cor:Isom4PunctSphere}
  Let $S$ be a pleated\index{pleated surface} 4-punctured sphere with hyperbolic structure identical to the boundary of an ideal hyperbolic tetrahedron. Then any of the three rotations of \reflem{IsomTetr} gives an isometry of $S$.
\end{corollary}

\begin{figure}
\includegraphics{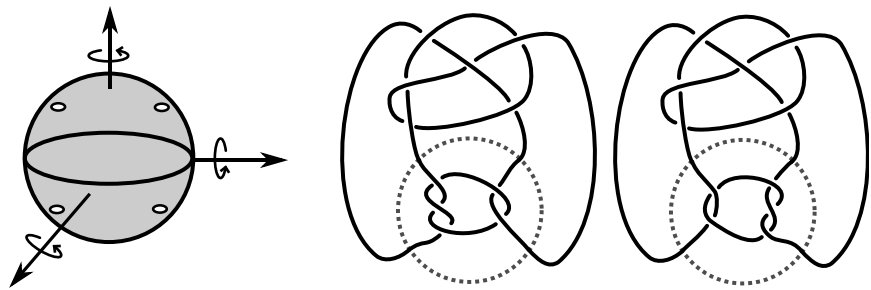}
  \caption{Mutation\index{mutation} cuts along a Conway sphere,\index{Conway sphere} performs one of the involutions shown on the left, and then reglues. Shown on the right is an example of two distinct knots related by mutation.}
  \label{Fig:Mutation}
\end{figure}

\begin{definition}\label{Def:Mutation}
A \emph{mutation}\index{mutation} of a knot or link is obtained by cutting along a Conway sphere,\index{Conway sphere} rotating by $\pi$ along one of three axes shown in \reffig{Mutation}, and then regluing.
\end{definition}

\begin{theorem}\label{Thm:Mutation}
  Let $K$ be a hyperbolic knot or link admitting an embedded essential Conway sphere.\index{Conway sphere} Let $K^\mu$ be any mutation\index{mutation} of $K$. Then $K^\mu$ is hyperbolic, and $\vol(K) = \vol(K^\mu)$. 
\end{theorem}

\begin{proof}
Let $S$ be the essential Conway sphere,\index{Conway sphere} and pleat $S$. If the pleating is embedded, isometric to the boundary of an ideal tetrahedron, then we may cut along the pleated\index{pleated surface} surface to obtain two hyperbolic manifolds whose boundaries are isometric pleated 4-punctured spheres, and isometric to the boundary of a hyperbolic ideal tetrahedron. By \refcor{Isom4PunctSphere}, rotation by $\pi$ through an axis orthogonal to opposite edges of the tetrahedron gives an isometry of the pleated surface $S$. Thus we may apply this isometry to one of the pieces and reglue, to obtain a complete hyperbolic manifold with an embedded essential Conway sphere,\index{Conway sphere} and volume equal to the volume of $S^3-K$. Note this is exactly a mutation.\index{mutation}

In the case that the pleating is not embedded, then Ruberman shows that $S$ is still isotopic to an embedded 4-punctured sphere that is a minimal surface with respect to the hyperbolic metric, and that mutation\index{mutation} is an isometry of this minimal surface \cite{Ruberman:mutation}. Thus the same argument applies to show the volumes agree.
\end{proof}

\section{Fuchsian, quasifuchsian, and accidental surfaces}

We now return to more general essential surfaces embedded in a hyperbolic 3-manifold. 

Let $M$ be hyperbolic, such that $M$ is the interior of a compact 3-manifold $\overline{M}$ with boundary.
Since $M$ is hyperbolic, we know there is a discrete, faithful representation $\rho\from \pi_1(M)\to \PSL(2,\CC)$ (\refprop{FreePropDisc}). If $S$ is a surface properly embedded in $M$, then the restriction of $\rho$ to $\pi_1(S)$ will be a discrete subgroup of $\PSL(2,\CC)$. We will consider properties of this subgroup. 

Let $\Gamma\leq\PSL(2,\CC)$ be a discrete group. Recall from \refdef{LimitSet} that the limit set\index{limit set} of $\Gamma$ is the set of accumulation points on $\bdy\HH^3$ of the orbit $\Gamma(x)$ for any point $x\in\HH^3$.

\begin{definition}\label{Def:QFuchsian}
A discrete group $\Gamma\leq\PSL(2,\CC)$ is said to be \emph{Fuchsian}\index{Fuchsian} if its limit set is a geometric circle on $\bdy\HH^3$. If its limit set is a Jordan curve and no element of $\Gamma$ interchanges the complementary components of the limit set, then $\Gamma$ is said to be \emph{quasifuchsian}.\index{quasifuchsian group}
\end{definition}

\begin{example}\label{Example:Quasifuchsian}
Let $\Gamma \leq \PSL(2,\RR)$ be the image of a discrete faithful representation of the fundamental group of a hyperbolic surface $S$ that is either closed or punctured without boundary. Then the limit set of $\Gamma$ in $\HH^2$ is all of $\bdy\HH^2$. Now view $\Gamma\leq\PSL(2,\RR)\leq\PSL(2,\CC)$ as acting on a hyperplane $H$ in $\HH^3$. When we extend the action of $\Gamma$ to all of $\HH^3$, the limit set is the geometric circle that is the boundary of the hyperplane $\bdy H$. Thus $\Gamma$ is Fuchsian.\index{Fuchsian}

Now adjust the representation very slightly, to $\Gamma_\epsilon \leq \PSL(2,\CC)$. The limit set also adjusts slightly. If $\Gamma_\epsilon$ is no longer a subgroup of $\PSL(2,\RR)$, then its limit set is no longer a geometric circle. However, it will be a topological circle. Thus $\Gamma_\epsilon$ is quasifuchsian. An example is shown in \reffig{Fractal}; this figure first appeared in \cite{thurston:bulletin}. 

\begin{figure}
\begin{center}
\includegraphics{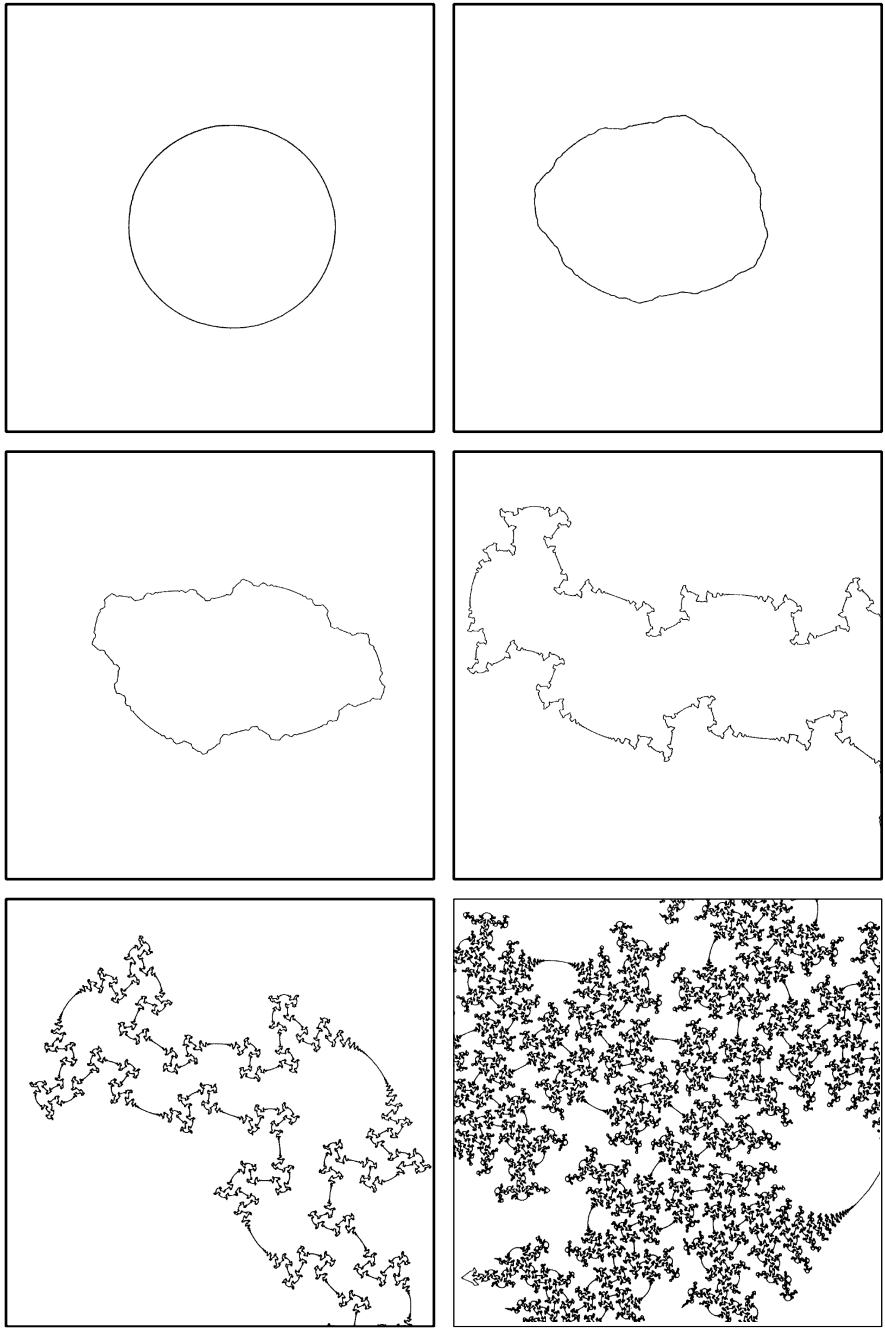}
\end{center}
\caption{The limit set of a Fuchisan group, and various limit sets of quasifuchsian groups obtained by deforming the Fuchsian\index{Fuchsian} group slightly. Figures are from \cite{thurston:bulletin}.}
\label{Fig:Fractal}
\end{figure}

The examples of \reffig{Fractal} were created by computer. Adjusting deformations of Fuchsian group by computer leads to beautiful fractal images. See, for example, \cite{IndrasPearls}. Software to visualize limit sets has also been developed by Wada \cite{Wada:Opti}. Yamashita has written a note to help users create their own software \cite{Yamashita:Software}. As a first step for the interested reader, we suggest working through Yamashita's example in \refex{Software}. For further work, the book \cite{IndrasPearls} includes direction on creating and exploring limit sets by computer.
\end{example}

\begin{definition}\label{Def:QFuchsianSfce}
Let $M$ be a hyperbolic 3-manifold and $S\subset M$ a properly embedded essential\index{essential} surface. Let $\rho\from \pi_1(M)\to\PSL(2,\CC)$ be a discrete, faithful representation. The surface $S$ is \emph{totally geodesic},\index{totally geodesic surface}\index{surface!totally geodesic} if, under the induced representation, the image $\rho(\pi_1(S))\leq\PSL(2,\CC)$ is Fuchsian. Sometimes a totally geodesic surface is also called \emph{Fuchsian}.\index{Fuchsian}\index{surface!Fuchsian} The surface $S$ is \emph{quasifuchsian}\index{quasifuchsian surface}\index{surface!quasifuchsian} if $\rho(\pi_1(S))\leq \PSL(2,\CC)$ is quasifuchsian.
\end{definition}

\begin{definition}\label{Def:Accidental}
Let $S$ be a surface properly embedded in $M$. A nontrivial loop $\gamma$ that is not freely homotopic into $\bdy S$ in $S$ is called an \emph{accidental parabolic}\index{accidental parabolic} if $\rho(\gamma)$ is parabolic\index{parabolic} in $\PSL(2,\CC)$. The surface $S$ is said to be \emph{accidental}\index{accidental surface} if it contains an accidental parabolic. 
\end{definition}

\begin{theorem}\label{Thm:NotAccidentalIfGeodesic}
If $S$ is a totally geodesic or quasifuchsian surface\index{quasifuchsian surface}\index{surface!quasifuchsian} properly embedded in the hyperbolic 3-manifold $M$, then $S$ is not accidental.\index{accidental surface}
\end{theorem}

\begin{proof}
If $S$ is a totally geodesic surface, then any closed curve $\gamma$ in $S$ that is not freely homotopic into $\bdy S$ must be freely homotopic to a closed geodesic in $S$. In turn, this closed geodesic in $S$ is a closed geodesic in $M$. Thus $\rho(\gamma)$ has a geodesic axis, and cannot be parabolic.\index{parabolic} Thus $S$ has no accidental parabolics in this case. 
 
If $S$ is quasifuchsian, the same argument does not immediately apply, However, it is known that a quasifuchsian group\index{quasifuchsian group} cannot contain an accidental parabolic\index{accidental parabolic} element. See Chapter IX, proposition D.17 of \cite{Maskit:KleinianGroups}.
Thus $S$ is not accidental. 
\end{proof}

\begin{corollary}\label{Cor:AltNoGeodesicSfce}
  Let $K$ be a knot or link with a prime, connected, alternating diagram,\index{alternating diagram}\index{alternating knot or link} and suppose $K$ is not a $(2,q)$-torus knot. Then the hyperbolic manifold $S^3-K$ contains no closed embedded totally geodesic surface.
\end{corollary}

\begin{proof}
Suppose $S$ is a closed, embedded, totally geodesic surface in the link complement $S^3-K$. By \reflem{Meridional}, any closed essential\index{essential} surface contains a closed curve that encircles a meridian of $K$. In particular, $S$ must contain a closed curve $\gamma$ encircling a meridian. But then $\gamma$ is freely isotopic to a meridian of $K$, meaning $\gamma$ is an accidental parabolic.\index{accidental parabolic} This contradicts \refthm{NotAccidentalIfGeodesic}. 
\end{proof}

In \cite{MenascoReid}, Menasco and Reid were the first to observe \refcor{AltNoGeodesicSfce}. Based on their observation, they made the following conjecture, which is still open at the time of writing this book.

\begin{conjecture}[Menasco--Reid conjecture]\label{Conj:MenascoReid}\index{Menasco--Reid conjecture}
  Let $K$ be a knot in $S^3$ such that $S^3-K$ is hyperbolic. Then $S^3-K$ admits no closed, embedded, totally geodesic surface. 
\end{conjecture}

Since the conjecture was proposed in the 1990s, evidence has developed both for and against the Menasco--Reid conjecture. As evidence for the conjecture, Menasco and Reid showed that in addition to alternating knots,\index{alternating knot or link} additional classes of hyperbolic knots cannot contain a closed embedded totally geodesic surface (closed 3-braids and tunnel number one knots) \cite{MenascoReid}. Since then, even more classes of knots and links have been shown to contain no closed, embedded, totally geodesic surfaces; a summary of such results can be found in the survey \cite{Adams:HyperbolicKnots}.

On the other hand, \refconj{MenascoReid} is known to be false for link complements, shown first in \cite{MenascoReid}. Leininger showed that there exists a sequence of hyperbolic knots whose complements contain closed embedded essential surfaces with principal curvatures converging to zero \cite{leininger}; if the principal curvatures were known to be zero the surfaces would be totally geodesic. DeBlois showed that \refconj{MenascoReid} does not hold for knots in rational homology spheres \cite{deblois-surfaces}. And Adams and Schoenfeld showed that \refconj{MenascoReid} is false if surfaces are allowed to have punctures \cite{AdamsSchoenfeld}. For example, they showed that the checkerboard surface\index{checkerboard surface} of certain pretzel knots, such as the surface shown in \reffig{TwistedBand}, is totally geodesic. 

Most of the evidence in support of \refconj{MenascoReid} is obtained by showing that any closed surface properly embedded in a particular type of knot complement must contain an accidental parabolic,\index{accidental parabolic} similar to \reflem{Meridional}. Thus there is interest in finding examples of essential surfaces without accidental parabolics. 

If we consider surfaces with (parabolic) boundary, we already have most of the tools in place to prove the following. 

\begin{theorem}\label{Thm:CheckerboardNotAccidental}
  Let $K$ be a link with a connected, prime, reduced alternating diagram,\index{alternating diagram}\index{alternating knot or link!checkerboard surface} and let $\Sigma$ be one of its checkerboard surfaces. Then $\Sigma$ is not accidental.\index{accidental surface}\index{checkerboard surface!not accidental}
\end{theorem}

Before we prove the theorem, we need a definition and a lemma.

\begin{definition}\label{Def:ParabolicLocus}
  Let $M$ be a hyperbolic 3-manifold, such that $M$ is the interior of a compact 3-manifold $\overline{M}$ with boundary. 
  The \emph{parabolic locus}\index{parabolic locus} $P$ of $M$ consists of tori and annuli in $\bdy\overline{M}$ such that each simple curve in $P$ lifts to a parabolic element\index{parabolic} of $\pi_1(M)\leq\PSL(2,\CC)$.
\end{definition}

\begin{lemma}\label{Lem:AccidentalAnnulus}
Suppose $S$ is a $\pi_1$-essential\index{$\pi_1$-essential} surface properly embedded in an irreducible,\index{irreducible} boundary irreducible\index{boundary irreducible} 3-manifold $M$, and suppose $S$ is accidental.\index{accidental surface} Then there is an essential annulus $A$ embedded in $M\cut S$ with one boundary component on the parabolic locus\index{parabolic locus} $P$ of $M$ and one boundary component an essential closed curve on $\widetilde{S}$.
\end{lemma}

The proof of the lemma uses the annulus theorem of Jaco \cite[Theorem~VIII.13]{Jaco:3Manifold} stated below. Briefly, it ensures we can replace an immersion of an annulus into a compact 3-manifold with an embedding. For a proof of the annulus theorem, see \cite{Jaco:3Manifold}. Compare to \refthm{LoopTheorem}, the loop theorem.

\begin{theorem}[Annulus theorem]\label{Thm:AnnulusTheorem}\index{annulus theorem}\index{Jaco annulus theorem}
  Let $M$ be a compact, irreducible\index{irreducible} 3-manifold with incompressible boundary. Suppose $f\from (A, \bdy A) \to (M,\bdy M)$ is a proper map, i.e.\ $f$ takes $\bdy A$ to $\bdy M$. Suppose also that $f$ is nondegenerate, i.e.\ that $f$ cannot be homotoped to the boundary of $M$. Then there exists an embedding $g\from (A, \bdy A) \to (M,\bdy M)$ that is nondegenerate. Furthermore, if the restriction of $f$ to $\bdy A$ is an embedding, then $g$ may be chosen so that its restriction to $\bdy A$ is the same embedding. 
\end{theorem}

\begin{proof}[Proof of \reflem{AccidentalAnnulus}]
If $S$ is accidental,\index{accidental surface} then there exists a nontrivial closed curve on $S$ that is freely homotopic into $\bdy M$ through $M$. Note if $S$ is nonorientable, then $\widetilde{S}$, the boundary of a regular neighborhood of $S$, is also accidental,\index{accidental surface} with accidental parabolic\index{accidental parabolic} a double cover of the curve on $S$.
So we may assume there is a nontrivial closed curve $\gamma$ on $\widetilde{S}$ that is freely homotopic into $\bdy M$ through $M$. The free homotopy defines a map of an annulus $A'$ into $M$; one boundary component of $A'$ lies on $\gamma$ and one on $\bdy M$. Adjust $A'$ so all intersections with $\widetilde{S}$ are transversal, and move the component of $A'$ on $\widetilde{S}$ in a bicollar of $S$ to be disjoint from $\widetilde{S}$.

Now consider intersections of the interior of $A'$ with $\widetilde{S}$. Consider first a closed curve of intersection that bounds a disk on $A'$. Since $\widetilde{S}$ is incompressible\index{incompressible} (because $S$ is $\pi_1$-essential),\index{$\pi_1$-essential} an innermost such curve also bounds a disk in $\widetilde{S}$. Since $M$ is irreducible,\index{irreducible} the union of the disk on $A'$ and that on $\widetilde{S}$ bounds a ball in $M$, and we may isotope $A'$ through the ball to remove the intersection. Thus we may assume there are no closed curves of intersection that bound disks on $A'$. Suppose there is an arc of intersection $A'\cap\widetilde{S}$ with both endpoints on $\bdy M$. This arc co-bounds a disk on $A'$ along with an arc on $\bdy A'\subset \bdy M$. Because $\widetilde{S}$ is boundary incompressible\index{boundary incompressible} and $M$ is boundary irreducible,\index{boundary irreducible} an innermost such arc may be isotoped away. So we assume there are no such arcs of intersection. Finally, because we have isotoped $A'$ away from $\widetilde{S}$ in a bicollar of the curve $\gamma\subset\widetilde{S}$, there are no arcs of intersection $A'\cap\widetilde{S}$ with an endpoint on $\widetilde{S}$. Thus there are no arcs of intersection of $A'\cap\widetilde{S}$. The only remaining possibility is that $A'\cap \widetilde{S}$ is a collection of essential closed curves on $A'$.

Apply a homotopy to minimize the number of closed curves of intersection. There is a sub-annulus $A''\subset A'$ that is outermost: it has one boundary component on $\bdy M$ and one on $\widetilde{S}$ and interior disjoint from $\widetilde{S}$. Thus we may consider $A''$ as an immersion of an annulus into $M\cut\widetilde{S}$. It is nondegenerate, else we could have reduced the number of closed curves of intersection of $A'$. Now we apply \refthm{AnnulusTheorem}, the annulus theorem.\index{annulus theorem}\index{Jaco annulus theorem} There exists a nondegenerate embedding of an annulus $A$ into $M\cut \widetilde{S}$ with one boundary component on $\widetilde{S}$ and one on the parabolic locus\index{parabolic locus} $\bdy M$. To finish the proof, we need to show that the embedding lies in $M\cut S$.

Note that $M\cut \widetilde{S}$ consists of two components, one homeomorphic to $M\cut S$ and one to a regular neighborhood of $S$. The regular neighborhood of $S$ only meets $\bdy M$ in a neighborhood of $\bdy S$. Since $\bdy A$ has a component on $\bdy M$, if $A$ is embedded in the neighborhood of $S$, it has a component running parallel to $\bdy S$ on $\bdy M$. But then the retraction of $A$ to $S$ defines a free homotopy of the closed curve $\bdy A\cap S$ to $\bdy S$, contradicting the definition of accidental.\index{accidental surface} Thus $A$ is embedded in the component of $M\cut \widetilde{S}$ that is homeomorphic to $M\cut S$. 
\end{proof}

\begin{proof}[Proof of \refthm{CheckerboardNotAccidental}]
Let $\Sigma$ be a checkerboard surface, and suppose by way of contradiction that it is accidental.\index{accidental surface} Then by \reflem{AccidentalAnnulus}, there exists an essential annulus $A$ embedded in $(S^3-N(K))\cut\Sigma$ with one boundary component on $\widetilde{\Sigma}$ and one on $N(K)$.

Consider the bounded polyhedral decomposition of $(S^3-N(K))\cut\Sigma$ (\reflem{BoundedPolyAlt}).
By \refthm{BddNormalForm}, we may take $A$ to be in normal form\index{normal form}\index{normal} with respect to the polyhedra. Note that because components of $\bdy A$ lie entirely in $\widetilde{S}$ and in $\bdy N(K)$, respectively, there are no arcs of $\bdy A$ intersecting a boundary edge adjacent to a surface face. Thus by \reflem{BddGaussBonnet}, the combinatorial area\index{combinatorial area} of $A$ is $0$. It follows that each normal disk\index{normal} of $A$ has combinatorial area $0$. This is possible only if the disk has one of three forms: each normal disk either meets exactly two boundary faces and no edges, or it meets exactly one boundary face and exactly two edges, or it meets exactly four edges. Because one component of $\bdy A$ lies on $\bdy N(K)$, there must be a normal disk meeting a boundary face. If the normal disk meets two boundary faces, then there is an arc of intersection of $A$ with a white face that runs from the boundary component of $\bdy A$ on $\bdy M$ back to the same boundary component, cutting off a disk in $A$. Because the white checkerboard surface is boundary incompressible,\index{boundary incompressible} such an arc bounds a disk on the white checkerboard surface. By normality, the disk cannot be contained in a single white face: the arc would run from one boundary edge back to the same edge. But an innermost intersection with a shaded face would give a crossing arc\index{crossing arc} cutting off a boundary compression disk\index{boundary compression disk} for the link. This is also impossible in a reduced diagram. 
So we may assume that there is a normal disk of $A$ meeting exactly one boundary face and exactly two edges of the polyhedron. On $A$, such an arc runs from the component $\bdy A\cap \bdy M$ of $\bdy A$ to $A\cap S$, which is the other component of $\bdy A$.

Then \reflem{21Annulus} implies that all normal disks\index{normal} of $A$ have this form, and that $K$ is a $(2,q)$-torus link. The surface $S$ is the annulus lying between the two strands of the link. Moreover, the boundary component $\gamma$ of $\bdy A$ on $\widetilde{S}$ must run along the core of the annulus $S$. It follows that $\gamma$ is boundary parallel. But then $\gamma$ is not accidental.\index{accidental surface} This is a contradiction. 
\end{proof}

\section{Fibers and semifibers}

Consider again the limit set of a group $\rho(\pi_1(S))$ where $S$ is a surface embedded in a hyperbolic 3-manifold $M$ and $\rho\from \pi_1(M)\to\PSL(2,\CC)$ is the holonomy\index{holonomy} representation. \Reffig{Fractal} shows the limit set of Fuchsian and quasifuchsian examples.\index{quasifuchsian surface}\index{surface!quasifuchsian} There is an additional option: the limit set of a discrete group isomorphic to $\pi_1(S)$ might be a space-filling curve. In this section, we will analyze surfaces with this property. First we present two topological definitions. 

\begin{definition}\label{Def:Fiber}
  Let $S$ be a surface properly embedded in a 3-manifold $M$. We say $S$ is a \emph{fiber}\index{fiber}\index{surface!fiber} if $M$ can be written as a fiber bundle over $S^1$, with fiber the surface $S$. Equivalently, there is a homeomorphism $f\from S\to S$ such that $M$ is the mapping torus
  \[ M \cong (S\times [0,1])/ (0,x)\sim(1,f(x)). \]
\end{definition}

\begin{definition}\label{Def:IBundle}
  An \emph{$I$-bundle}\index{$I$-bundle} is a 3-manifold homeomorphic to $S\times I$, where $S$ is a surface, possibly with boundary. The \emph{vertical boundary}\index{vertical boundary}\index{$I$-bundle!vertical boundary} is $\bdy S\times I$; note it is a collection of annuli. The \emph{horizontal boundary}\index{horizontal boundary}\index{$I$-bundle!horizontal boundary} consists of $S\times \bdy I$. If $S$ is orientable, this consists of the disjoint union of $S\times\{0\}$ and $S\times\{1\}$. If $S$ is nonorientable, we say that the $I$-bundle is \emph{twisted},\index{twisted $I$-bundle}\index{$I$-bundle!twisted} and the horizontal boundary is homeomorphic to the oriented double cover of $S$. We often denote a twisted $I$-bundle by $S\widetilde{\times} I$. See \refex{TwistedIBundle}.
\end{definition}

\begin{definition}\label{Def:Semifiber}
  A surface $S$ properly embedded in a 3-manifold $M$ is a \emph{semifiber}\index{semifiber}\index{surface!semifiber} if it is either a fiber, or if $S$ is the boundary of an $I$-bundle $S'\times I$ over a nonorientable surface $S'$, and $M$ is obtained by gluing two copies of this $I$-bundle by the identity on $S$. In the latter case, sometimes $S$ is called a \emph{strict semifiber}.\index{strict semifiber}\index{surface!strict semifiber}
\end{definition}

A strict semifiber is an example of a \emph{virtual fiber}, defined below. See \refex{SemifiberVirtual}.

\begin{definition}\label{Def:VirtualFiber}
A surface $S$ properly embedded in a 3-manifold $M$ is called a \emph{virtual fiber}\index{virtual fiber}\index{fiber!virtual} if there is a finite index cover of $M$ in which $S$ lifts to a fiber. 
\end{definition}

The following is due to Thurston \cite{thurston} and Bonahon \cite{Bonahon:Ends}. See also \cite{CanaryEpsteinGreen:Notes}.

\begin{theorem}\label{Thm:Trichotomy}
Let $S$ be an essential\index{essential} surface in a hyperbolic 3-manifold $M$. Then $S$ has exactly one of three forms:
  \begin{enumerate}
  \item $S$ is Fuchsian\index{Fuchsian} or quasifuchsian,\index{quasifuchsian surface}\index{surface!quasifuchsian}
  \item $S$ is accidental,\index{accidental surface} or
  \item $S$ is a virtual fiber.
  \end{enumerate}
\end{theorem}

The proof of the theorem is obtained by analyzing surfaces that are not accidental, and whose limit set is not a circle or topological circle.

\begin{example}\label{Example:Fig8Fiber}
The figure-8 knot complement contains a surface that is a fiber, namely the punctured torus shown in \reffig{Fig8Seifert}.
\begin{figure}
  \includegraphics{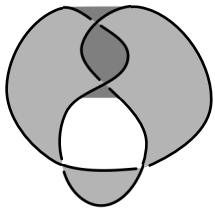}
  \caption{The Seifert surface of the figure-8 knot complement.}
  \label{Fig:Fig8Seifert}
\end{figure}

A portion of the limit set of this surface was computed by S.~Schleimer, following W.~Thurston, and is shown in \reffig{Fig8Fiber2D}.
\begin{figure}
  \includegraphics{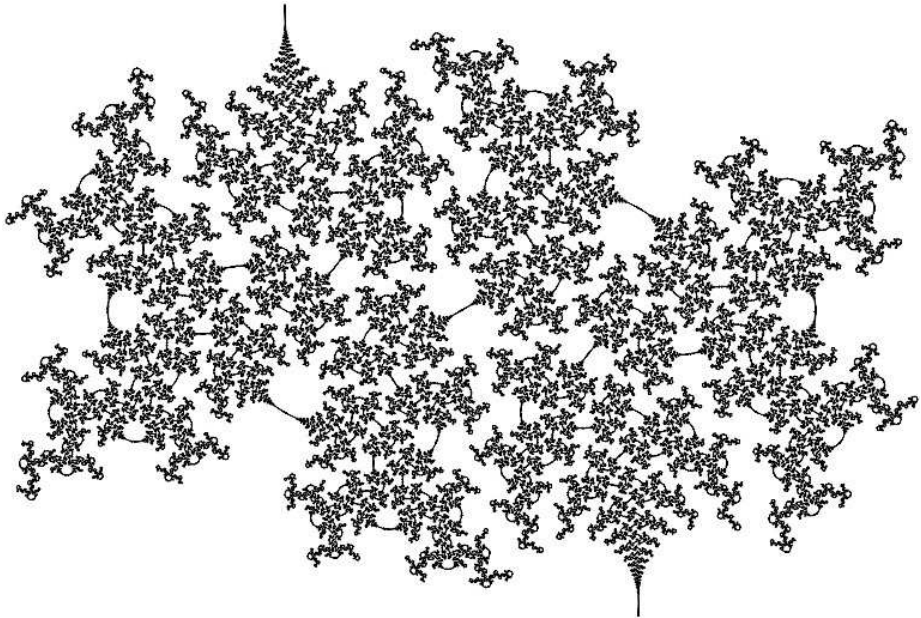}
\caption{The limit set of the Seifert surface of the figure-8 knot complement, created by S.~Schleimer.}
  \label{Fig:Fig8Fiber2D}
\end{figure}
Its lift to the universal cover, given a pleating, is shown in \reffig{Fig8Fiber3D}, due to S.~Schleimer and H.~Segerman.
\begin{figure}
  \includegraphics{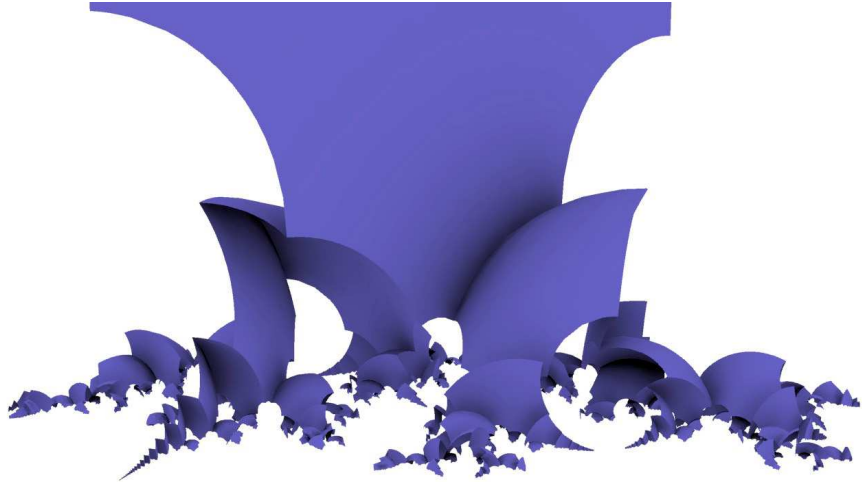}
\caption{The lift of the figure-8 knot Seifert surface to the universal cover $\HH^3$, with a pleating. Created by S.~Schleimer and H.~Segerman.}
  \label{Fig:Fig8Fiber3D}
\end{figure}

Another view of the surface is shown in \reffig{Fig8Fiber3D2}, due to D.~Bachman, S.~Schleimer, and H.~Segerman. Note that in \reffig{Fig8Fiber3D}, the cusps of the surface have been cut off to show a larger view of the pleating. On the other hand, \reffig{Fig8Fiber3D2} gives a better view of the surface near infinity, without cusps cut off. 
\begin{figure}
  \includegraphics{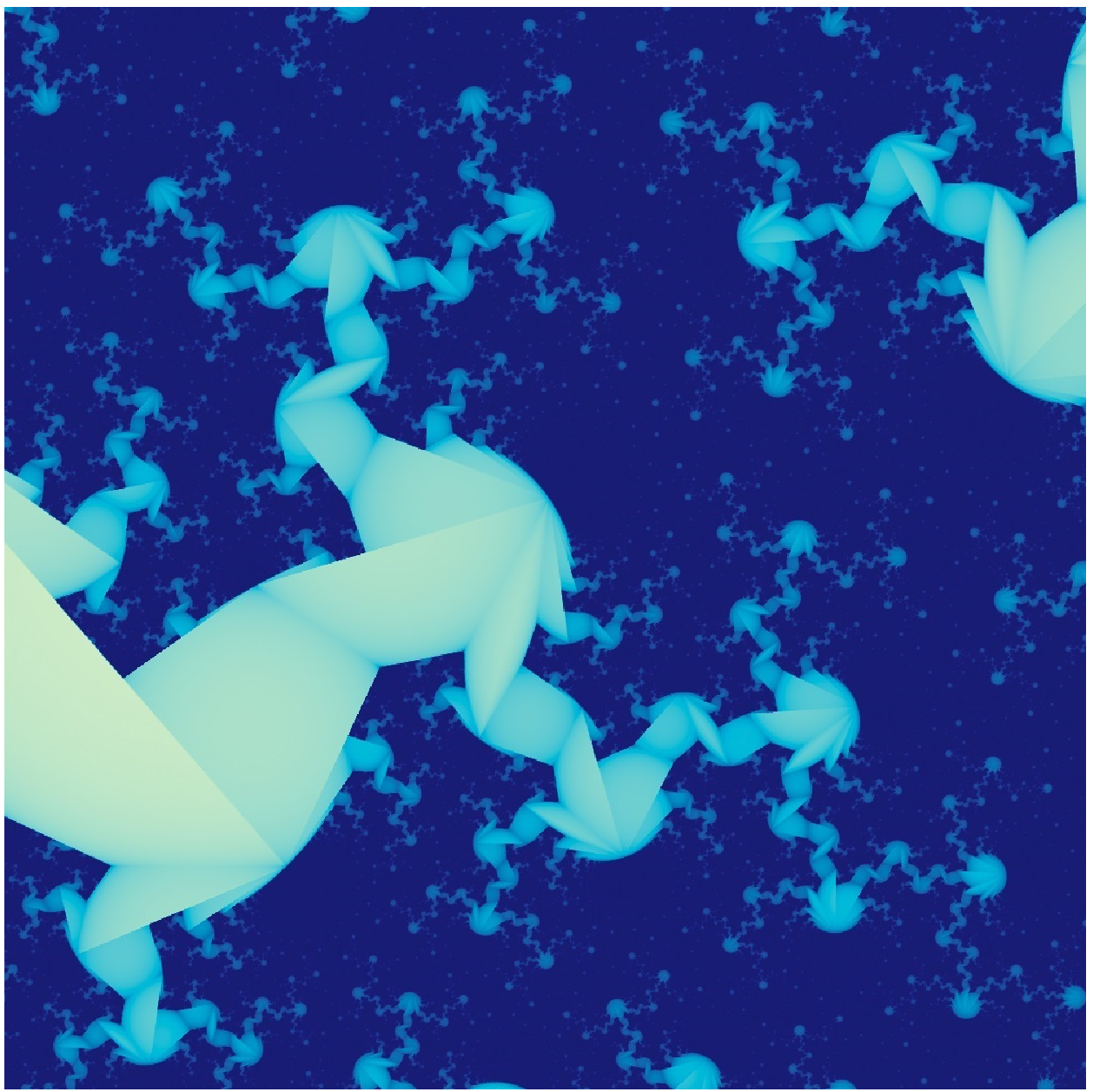}
  \caption{More of the lift of the figure-8 knot Seifert surface to $\HH^3$, without cusps cut off. Created by D.~Bachman, S.~Schleimer and H.~Segerman.}
  \label{Fig:Fig8Fiber3D2}
\end{figure}
\end{example}

We will be dealing only with embedded surfaces. In the case a surface is embedded, the virtual fiber case of the trichotomy reduces to a simpler situation.

\begin{lemma}\label{Lem:TrichotomyEmbedded}
Suppose $S$ is a properly embedded surface in a 3-manifold $M$. Then $S$ is a virtual fiber if and only if $S$ is a semifiber. 
\end{lemma}

\begin{proof}
\Refex{TrichotomyEmbedded}.
\end{proof}

As an additional example of a fibered surface in a link complement, consider the checkerboard surfaces of the $(2,q)$-torus knot or link.

\begin{lemma}\label{Lem:2qTorusFiber}
One of the checkerboard surfaces of a standard diagram of a $(2,q)$-torus knot or link is a fiber.\index{checkerboard surface}\index{checkerboard surface!fiber}
\end{lemma}

\begin{proof}
One of the checkerboard surfaces is an annulus or M\"obius band running between the two strands of the link. Let that be the white checkerboard surface. We will show the shaded surface $\Sigma$ is a fiber. Note that $\Sigma$ is orientable, as it is built of two disks with a sequence of (singly) twisted bands between them. Thus the cut manifold $(S^3-N(K))\cut\Sigma$ has boundary consisting of parabolic locus\index{parabolic locus} and $\widetilde{\Sigma}$, which is two copies of $\Sigma$ in the orientable case. 
The surface $\Sigma$ is a semifiber if and only if the cut manifold $(S^3-N(K))\cut \Sigma$ is an $I$-bundle:\index{$I$-bundle}
\[ (S^3-N(K))\cut \Sigma \cong \Sigma \times I. \]

Consider the polyhedral decomposition of the cut manifold in this case. The two polyhedra consist of a chain of adjacent white bigons\index{bigon} along with two shaded disks, one inside and one outside the chain of bigons. Note each of these polyhedra is an $I$-bundle over the shaded face, of the form $D\times I$ where $D$ is a shaded disk. White bigon faces are of the form $\alpha_i\times I$, where $\alpha_i$ is an arc with endpoints on edges of the polyhedron (shaded faces). The parabolic locus\index{parabolic locus} consists of boundary squares, which are also products $\mbox{arc}\times I$, parallel to $\alpha_i\times I$ on their sides meeting white faces, with endpoints of arcs on shaded faces. 

To obtain $(S^3-N(K))\cut \Sigma$, glue white faces. The gluing takes each bigon\index{bigon} face $\alpha_i\times I$ to another bigon face $\alpha_j\times I$, matching the $I$-bundle structure. Thus $(S^3-N(K))\cut \Sigma$ is an $I$-bundle. So $\Sigma$ is a semifiber.\index{semifiber}

To see that $\Sigma$ is actually a fiber, note that the gluing of the two polyhedra along white faces matches $D_1\times\{0\}$ in one polyhedron to $D_2\times\{0\}$ in the other, and $D_1\times\{1\}$ to $D_2\times\{1\}$ in the other. Thus the boundary of the $I$-bundle has two components, so it is not an $I$-bundle\index{$I$-bundle} over a nonorientable surface, and cannot be a strict semifiber. 
\end{proof}

\begin{theorem}\label{Thm:AltFiber}
Let $K$ be a knot or link with a connected, twist-reduced, prime, alternating diagram,\index{alternating diagram}\index{alternating knot or link} and let $\Sigma$ be an associated checkerboard surface. Then $\Sigma$ is a semifiber if and only if $K$ is a $(2,q)$-torus link and $\Sigma$ is the checkerboard surface of \reflem{2qTorusFiber} that is a fiber. \index{checkerboard surface}\index{checkerboard surface!fiber}
\end{theorem}

Before proving the theorem, we give a lemma. Its proof is very similar to Lemma~4.17 of \cite{fkp:guts}; see also \cite{HowiePurcell}.

\begin{lemma}\label{Lem:ProductRectangles}
Let $K$ be a knot or link as in the statement of \refthm{AltFiber}, and let $\Sigma$ be its shaded checkerboard surface. 
Let $B$ be an $I$-bundle\index{$I$-bundle} embedded in $M_\Sigma = (S^3-N(K))\cut \Sigma$, with horizontal boundary on $\widetilde{\Sigma}$, and suppose the vertical boundary of $B$ is essential. Let $W$ be a white face of the polyhedral decomposition of the cut manifold. Then $B\cap W$ is isotopic in $M_\Sigma$ to a collection of product rectangles $\alpha\times I$, where $\alpha\times\{0\}$ and $\alpha\times\{1\}$ are arcs of ideal edges on the boundary of $W$. 
\end{lemma}

\begin{proof}
First suppose $B=Q\times I$ is a product $I$-bundle\index{$I$-bundle} over an orientable base. 
Consider a component of $\bdy(B\cap W)$. If it lies entirely in the interior of $W$, then it lies in the vertical boundary $V=\bdy Q\times I$.
The intersection $V\cap W$ then contains a closed curve component; an innermost one bounds a disk in $W$. Since the vertical boundary is essential, we may isotope $B$ to remove such intersections. So assume each component of $\bdy(B\cap W)$ meets $\widetilde{\Sigma}$. Note that it follows that each component of $B\cap W$ is a disk. 

Note $W\cap \widetilde{\Sigma}$ consists of ideal edges on the boundary of the face $W$. It follows that the boundary of each component of $B\cap W$ consists of arcs $\alpha_1$, $\beta_1$, $\dots$, $\alpha_n$, $\beta_n$ with $\alpha_i$ an arc in an ideal edge of $W\cap \widetilde{\Sigma}$ and $\beta_i$ in the vertical boundary of $B$, in the interior of $W$.
We may assume that each arc $\beta_i$ runs between distinct ideal edges, else isotope $B$ through the disk bounded by $\beta_i$ and an ideal edge to remove $\beta_i$, and merge $\alpha_i$ and $\alpha_{i+1}$.

We may assume that $\beta_i$ runs from $Q\times\{0\}$ to $Q\times\{1\}$, for if not, then $\beta_i \subset W$ is an arc from $\bdy Q\times\{1\}$ to $\bdy Q\times\{1\}$, say, in an annulus component of $\bdy Q \times I$. Such an arc bounds a disk in $\bdy Q\times I$. This disk has boundary consisting of the arc $\beta_i$ in $W$ and an arc on $\bdy Q\times\{1\} \subset\widetilde{\Sigma}$. If the disk were essential, it would give a contradiction to \refprop{NoNormalBigons}. So it is inessential, and we may isotope $B$ to remove $\beta_i$, merging $\alpha_i$ and $\alpha_{i+1}$.

Finally we show that $n=2$, i.e.\ that each component of $B\cap W$ is a quadrilateral with arcs $\alpha_1, \beta_1, \alpha_2, \beta_2$. For if not, there is an arc $\gamma\subset W$ with endpoints on $\alpha_1$ and $\alpha_3$. By sliding along the disk $W$, we may isotope $B$ so $\gamma$ lies in $B\cap W$. Then note that $\gamma$ lies in $Q\times I$ with endpoints on $Q\times\{1\}$. It must be parallel vertically to an arc $\delta\subset Q\times\{1\} \subset \widetilde{\Sigma}$. This gives another disk with boundary consisting of an arc on $W$ and an arc on $\widetilde{\Sigma}$. Either the disk contradicts \refprop{NoNormalBigons}, or $\alpha_1$ and $\alpha_3$ lie on the same ideal edge in the boundary of $W$. But then $\beta_1$ and $\beta_2$ are arcs running from $\alpha_2$ on one ideal edge on the boundary of $W$ to the same ideal edge on the boundary of $W$ containing $\alpha_1$ and $\alpha_3$. The only way that the boundary of this component of $B\cap W$ bounds a disk in $W$ is if $n=3$, and $\beta_3$ runs from an endpoint of $\alpha_1$ to an endpoint of $\alpha_3$. But then $\beta_3$ runs from $Q\times\{1\}$ to $Q\times\{1\}$, which we ruled out in the previous paragraph. So $n=2$ and $B\cap W$ is a product rectangle $\alpha_1\times I$.

Next suppose $B$ is a twisted $I$-bundle\index{$I$ bundle!twisted} $B=Q \widetilde{\times} I$ where $Q$ is non-orientable. Let $\gamma_1, \dots, \gamma_m$ be a maximal collection of orientation reversing closed curves on $Q$. Let $A_i\subset B$ be the $I$-bundle over $\gamma_i$. Each $A_i$ is a M\"obius band. The bundle $B_0 = B\setminus (\cup_i A_i)$ is then a product bundle $B_0 = Q_0\times I$ where $Q_0=Q\setminus (\cup_i \gamma_i)$ is an orientable surface. Our work above then implies that $B_0\cap W$ is a product rectangle for each white region $W$. To obtain $B\cap W$, we attach the vertical boundary of such a product rectangle to the vertical boundary of a product rectangle of $A_i$. This procedure respects the product structure of all rectangles, hence the result is a product rectangle. 
\end{proof}

\begin{proof}[Proof of Theorem~\ref{Thm:AltFiber}]
One direction is \reflem{2qTorusFiber}: If the link diagram is the standard diagram of a $(2,q)$-torus link, then a checkerboard surface is a fiber. 
  
Conversely, if the checkerboard surface $\Sigma$ is a semifiber, then $M_\Sigma = (S^3-N(K))\cut \Sigma$ is an $I$-bundle.\index{$I$-bundle} In this case, \reflem{ProductRectangles} implies $M_\Sigma$ intersects each white face $W$ in a product rectangle of the form $\alpha\times I$, where $\alpha\times\{0\}$ and $\alpha\times\{1\}$ lie on ideal edges of $W$. Since $W \subset M_\Sigma$, the face $W$ is a product rectangle, with exactly two ideal edges $\alpha\times\{0\}$ and $\alpha\times\{1\}$. Thus $W$ is a bigon.\index{bigon} So every white face is a bigon. Thus the diagram of $K$ is a chain of bigons lined up end to end. This is a $(2,q)$-torus link. The white checkerboard surface is obtained by gluing sides of those bigons, and so forms the annulus or M\"obius band between the link components. The shaded checkerboard surface must therefore be the fiber of \reflem{2qTorusFiber}. 
\end{proof}

\begin{corollary}\label{Cor:Quasifuchsian}
  Let $K$ be a knot or link with a connected, twist-reduced, prime, alternating diagram,\index{alternating diagram}\index{alternating knot or link}\index{alternating knot or link!checkerboard surface} and let $\Sigma$ be an associated checkerboard surface. If $K$ is hyperbolic, then $\Sigma$ is quasifuchsian.\index{quasifuchsian surface}\index{surface!quasifuchsian}\index{checkerboard surface!quasifuchsian}
\end{corollary}

\begin{proof}
  If $K$ is hyperbolic, it cannot be a $(2,q)$-torus link. Then \refthm{AltFiber} implies that $\Sigma$ cannot be a semifiber. Because $\Sigma$ is an embedded surface, \reflem{TrichotomyEmbedded} implies that $\Sigma$ is not a virtual fiber. \Refthm{CheckerboardNotAccidental} implies that $\Sigma$ is not accidental.\index{accidental surface} By \refthm{Trichotomy}, it must be quasifuchsian. 
\end{proof}

\section{Exercises}

\begin{exercise}\label{Ex:BeltTangleIncompress}
(Easy) Show that the 3-punctured sphere bounded by the crossing circle in a hyperbolic belted tangle\index{belted tangle} must be incompressible.\index{incompressible}
\end{exercise}

\begin{exercise}
A chain link is a link that has the form of a circular chain, as in \reffig{ChainLink}, left.\index{chain link} Note that the link components of the chain can be twisted. We define the minimally twisted chain link with an even number of components to be the chain link with every other link component lying flat in the plane of projection, and alternate link components to be perpendicular to the plane of projection.

A minimally twisted chain link may be augmented by adding a crossing circle encircling the circular chain, as in \reffig{ChainLink}, right.

\begin{figure}
\includegraphics{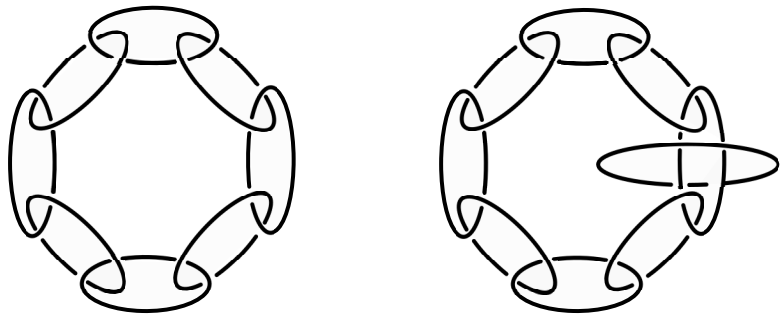}
\caption{Left: The minimally twisted chain link with eight link components. Right: the augmented minimally twisted chain link.}
\label{Fig:ChainLink}
\end{figure}

\begin{enumerate}
\item Using belted sums\index{belted sum} and the volume of the Whitehead link, find the volume of any augmented minimally twisted chain link with an even number of chain components. 
\item Find a belted tangle\index{belted tangle} $T$ such that repeatedly taking belted sums\index{belted sum} of $T$ with the Whitehead link gives the augmented minimally twisted chain link with an odd number of chain components. What is the volume of the augmented minimally twisted chain link with an odd number of chain components?
\end{enumerate}
\end{exercise}

\begin{exercise}\label{Ex:Software}
Following Yamashita's instructions, create a python program that allows us to visualize the limit set as hyperbolic structures are varied on a punctured torus. Print an example with a Fuchsian limit set, and three examples of quasifuchsian limit sets. See \cite{Yamashita:Software}.
\end{exercise}

\begin{exercise}\label{Ex:TwistedIBundle}
Suppose $S'$ is a closed nonorientable surface. Consider $S'\times I$ (also frequently denoted $S'\widetilde{\times}I$). Prove that its boundary is a closed orientable surface $S$ homeomorphic to the oriented double cover of $S'$. 
\end{exercise}

\begin{exercise}\label{Ex:SemifiberVirtual}
Prove that a strict semifiber is a virtual fiber. 
\end{exercise}

\begin{exercise}
Prove that a nonorientable surface can never be a fiber in a link complement $S^3-L$. That is, there are no strict semifibers for links in $S^3$. 
\end{exercise}

\begin{exercise}\label{Ex:TrichotomyEmbedded}
Prove \reflem{TrichotomyEmbedded}: that a properly embedded surface that is a virtual fiber in a 3-manifold must be a semifiber.
%% From Genevieve Walsh, Incompressible surfaces and spun normal form:
%% \begin{quote}
%% An embedded virtual fiber of $M$ is either a fiber of a fibration of $M$ over $S^1$, or a generic fiber of an orbifold-fibration of $M$ over an interval with mirrored endpoints. In this latter case, the surface is called a semi-fiber, and $M$ is the union of two twisted $I$-bundles over the surface.
%% \end{quote}
\end{exercise}

%    Part III
% Part III: Hyperbolic knot invariants
\part{Hyperbolic Knot Invariants}\label{Part:Invariants}
  
\chapter{Estimating Volume}\label{Chap:Volume}
\blfootnote{Jessica S. Purcell, Hyperbolic Knot Theory}

We have seen that hyperbolic 3-manifolds have finite volume if and only if they are compact or the interior of a compact manifold with finitely many torus boundary components (\refthm{FteVolIffTorusBdy}).
However, it is not completely straightforward to estimate volumes of large classes of manifolds, including knot complements. There are many open questions concerning the relationship of volume of a hyperbolic manifold to other invariants, such as knot invariants. 
In this chapter, we discuss different ways to estimate volumes of hyperbolic 3-manifolds that are defined topologically or combinatorially, such as knot complements.

%%%%%%%%%%%%%%%%%%%%%%%%%%%%%%%%%%%%%%%%%%%%%%%%%%%%%%%%%%%%%%%%%
\section{Summary of bounds encountered so far}

\subsection{Upper bounds}

It is usually an easier problem to give upper bounds on the volume of a hyperbolic 3-manifold than lower bounds, although there are exceptions, especially when sharp upper bounds are needed. Here we review two methods we have already encountered that can give upper bounds on volume.

\subsubsection{Volume bounds from polyhedra}
Recall from \refthm{MaxVolTet} that the maximal volume tetrahedron is the regular ideal tetrahedron. Its volume is the value $3\Lambda(\pi/3) := \vtet = 1.0149\dots$. In various chapters, we have found decompositions of several different knot and link complements into ideal tetrahedra. The volume of such a knot or link is therefore bounded by $\vtet$ times the number of tetrahedra in its decomposition.

For example, this can be used to show the following theorem, originally proved by Agol and D.~Thurston in the appendix to \cite{lackenby:alt-volume}.

\begin{theorem}\label{Thm:FullyAugUpperBound}
A fully augmented link\index{fully augmented link} $L$ with $t(L)$ crossing circles has volume at most $10\vtet(t(L)-1)$.\index{fully augmented link!volume bound}
\end{theorem}

\begin{proof}
In \refchap{TwistKnots}, we saw that a fully augmented link has a decomposition into two right angled ideal polyhedra $P_1$ and $P_2$, with white and shaded faces, where shaded faces are ideal triangles\index{ideal triangle} coming from 2-punctured disks bounded by crossing circles, and white faces come from the plane of projection.

For the polyhedron $P_1$, add a finite vertex $v_1$ in the interior and cone to the faces of the polyhedra. Do the same for $P_2$, adding vertex $v_2$ and coning. Each shaded triangle in $\bdy P_1$ gives rise to a tetrahedron. There are two shaded triangles per crossing circle in each of the two polyhedra, so $4t(L)$ tetrahedra arise in this way.

The white faces are coned to pyramids. Glue a pair of pyramids in $P_1$ and $P_2$ together across a matching white face, and perform stellar subdivision. That is, add an edge running from the finite vertex in one polyhedron through the center of the face to the finite vertex in the other polyhedron, then add triangles around the edge to divide the pyramids into tetrahedra. If the face has $d$ edges, it is subdivided into $d$ tetrahedra. Note each crossing circle contributes six edges to each polyhedron. Thus the total number of edges of the white faces will be $6t(L)$, and thus the white faces contribute $6t(L)$ tetrahedra to the decomposition.

This gives us $10t(L)$ tetrahedra, but these have finite vertices. We can improve the bound by choosing an ideal vertex $w_1$ in $P_1$, and collapsing the edge from $w_1$ to $v_1$. Similarly, choose the corresponding vertex $w_2$ in $P_2$, and collapse the edge from $w_2$ to $v_2$. Now simplify the triangulation by collapsing monogons to vertices, bigons\index{bigon} to a single edge, and parallel triangles to a single triangle. Note that all tetrahedra adjacent to $w_1$ and $w_2$ are collapsed to triangles under this procedure. We count the number of these. 

The ideal vertex $w_1$ is adjacent to two shaded triangles and two white faces. The white faces each have at least three edges, and so give rise to at least three tetrahedra to be collapsed, running between the two polyhedra. Thus there are at least six such tetrahedra arising from white faces. Each shaded face gives rise to one tetrahedron to be collapsed in each $P_i$, or four total. Thus there is an ideal triangulation with at most $10t(L)-10$ tetrahedra. 

Now apply \refthm{MaxVolTet}. The volume of each of the $10t(L)-10$ tetrahedra is at most $\vtet$. Thus the volume of the fully augmented link is at most $10\vtet(t(L)-1)$. 
\end{proof}

The bound of \refthm{FullyAugUpperBound} is \emph{asymptotically sharp},\index{asymptotically sharp volume bound} in the sense that there is a sequence of fully augmented links whose volumes approach the upper bound; this is proved in the first part of \refex{AsympSharpFullyAug}. 

\subsubsection{Dehn filling}
Recall Thurston's theorem on volume change under Dehn filling,\index{Dehn filling} \refthm{VolumeDF}: If $M$ is hyperbolic with cusps $C_1, \dots, C_n$, and $s_1, \dots, s_n$ are slopes, one on each $\bdy C_j$, such that $M(s_1, \dots, s_n)$ is hyperbolic, then
\[ \vol(M)>\vol(M(s_1, \dots, s_n)). \]
This result can be combined with the previous to give an upper bound on the volume of knots in terms of the twist number, first observed in \cite{lackenby:alt-volume}.

\begin{theorem}\label{Thm:VolUpperTwistRegions}
Suppose $K$ is a knot or link with a prime, twist-reduced diagram with twist number $\tw(K)\geq 2$. Then $S^3-K$ is hyperbolic, and the volume of $S^3-K$ satisfies
\[ \vol(S^3-K) < 10\vtet(\tw(K)-1). \]
\end{theorem}

Again the bound of \refthm{VolUpperTwistRegions} is asymptotically sharp;\index{asymptotically sharp volume bound} this is proved in the second part of \refex{AsympSharpFullyAug}. 

\begin{proof}[Proof of \refthm{VolUpperTwistRegions}]
Because the diagram of $K$ is prime and twist-reduced, when we fully augment $K$ by adding a crossing circle to each twist region, and then remove pairs of crossings to form a fully augmented link,\index{fully augmented link} the resulting link $L$ is hyperbolic; see \reflem{TwReducedGivesReducedAug} and \reflem{ExistsRightAngled}. It will have $t(L)=\tw(K)$ crossing circles. By \refthm{FullyAugUpperBound}, the volume of the fully augmented link is at most $10\vtet(\tw(K)-1)$.

Now, we obtain $S^3-K$ from $S^3-L$ by Dehn filling the crossing circles, filling the $i$-th one along a slope $1/n_i$ where $n_i$ is an integer such that $2n_i$ crossings were removed at that twist region to go from the diagram of $K$ to that of $L$. By Thurston's theorem on volume change under Dehn filling, \refthm{VolumeDF},
\[ \vol(S^3-K) < \vol(S^3-L) \leq 10\vtet(\tw(K)-1). \qedhere \]
\end{proof}

\subsection{Lower bounds via angle structures}

In addition to previously obtaining results that lead to upper bounds on volume, we have also built the tools to give lower bounds on hyperbolic volume in special cases, in particular when the manifold admits an angle structure.\index{angle structure} Recall \refthm{VolAngleStructs}: if the maximum of the volume functional\index{volume functional} over the set of all angle structures on a manifold $M$ occurs in the interior of the set of angle structures, then that angle structure gives the unique complete hyperbolic metric on $M$. We have the following corollary.

\begin{corollary}\label{Cor:VolAngleStructs}
Suppose $M$ is an orientable 3-manifold with boundary consisting of tori, with an ideal triangulation $\mathcal{T}$. Suppose that the set of angle structures\index{angle structure} $\calA(\calT)$ is nonempty, and that the volume functional\index{volume functional} takes its maximum on the interior of the set of all angle structures $\mathcal{A}(\mathcal{T})$. Let $A$ be any structure in the closure of $\mathcal{A}(\mathcal{T})$. Then $M$ is hyperbolic, and $\vol(M) \geq \mathcal{V}(A)$. 
\end{corollary}

\begin{proof}
The fact that $M$ is hyperbolic is a consequence of \refthm{HypAngleStruct}.
The volume functional\index{volume functional} is strictly concave down by \refthm{VolConcaveDown}, and is uniquely maximized at the complete hyperbolic structure in the interior by \refthm{VolAngleStructs} and \refthm{VolAngleStructsConverse}. Thus any structure in the closure of $\mathcal{A}(\mathcal{T})$ gives a volume at most that of $M$.
\end{proof}

\Refcor{VolAngleStructs} was used by Futer and Gu{\'e}ritaud to obtain bounds on the volumes of 2-bridge knots in terms of their continued fraction decompositions \cite{GueritaudFuter:2bridge}. They proved the following.

\begin{theorem}\label{Thm:Vol2Bridge}
Let $K$ be a reduced alternating diagram\index{alternating diagram} of a hyperbolic 2-bridge link\index{2-bridge knot or link}\index{2-bridge knot or link!volume bound} $K$ with $\tw(K)$ twist regions. Then
\[ \vol(S^3-K) > 2\vtet \tw(K) - 2.7066. \]
Moreover, the lower bound is asymptotically sharp.\index{asymptotically sharp volume bound}
\end{theorem}

\begin{proof}
We may suppose that the diagram of $K$ is determined by a continued fraction $[0,a_{n-1}, \dots, a_1]$, where $\tw(K)=n-1$, as in \refdef{2BridgeNotation}. The link complement has a geometric triangulation\index{geometric triangulation} discussed in \refchap{TwoBridge}, determined by real numbers $(z_1, z_2, \dots, z_{C-2},z_{C-1})$, where $C$ is the crossing number of the diagram of $K$; see the proof of \refprop{2BridgeAngleNonempty}. We will choose explicit values of the $z_i$ that give a structure in the boundary of the space of angle structures.\index{angle structure} By \refcor{VolAngleStructs} the volume of the result gives a lower bound on the actual hyperbolic volume. 

First assume that $\tw(K)\geq 3$, so $n\geq 2$.
We let $z_1 = z_{C-1}=0$, and we will choose $z_i=\pi/3$ for indices $i$ such that $a_1\leq i \leq C-a_{n-1}$.

These choices satisfy the hinge equation of \refeqn{2BridgeZConditions}: $|z_{i+1}-z_{i-1}| = 0<\pi-z_i=2\pi/3$, for appropriate $i$. They do not satisfy the strict inequality of the convexity equation of \refeqn{2BridgeZConditions}, but only satisfy the weak inequality: $2z_i \leq z_{i-1}+z_{i+1}$. When $a_1<i<C-a_{n-1}$, these will assign values to $x_i$ and $y_i$ using \reftable{LabelConditions}. Note the angles will take values of $\pi/3$ in the hinge case, but $2\pi/3$ and $0$ in the non-hinge case, and thus there will be flat tetrahedra. However, our choices so far give rise to a structure on the boundary of $\mathcal{A}(\mathcal{T})$, which will be sufficient for our purposes. Each hinge index contributes volume $\vtet$ to the structure, while non-hinges contribute nothing to volume. Note hinges occur between twist regions; there are $\tw(K)-3$ hinge indices between $a_1$ and $C-a_{n-1}$.

We cannot choose $z_i=\pi/3$ for all the indices in the first and last fans, else even the weak inequality of the convexity equations \eqref{Eqn:2BridgeZConditions} will not be satisfied near $i=1$ or $i=C-1$. Instead, in the first and last fans, interpolate between $0$ and $\pi/3$ in a way that satisfies the weak versions of \refeqn{2BridgeZConditions}. Then again angles $x_i$ and $y_i$ will be determined by \reftable{LabelConditions}. At the hinge indices $i=a_1$ or $i=C-a_{n-1}$, the angles will be:
\[ \frac{\pi}{2} - z_{i-1}, \quad \frac{\pi}{6}+z_{i-1}, \quad \frac{\pi}{3}. \]
The volume defined by these angles is smallest when $z_{i-1}=0$, which occurs when the three angles are $\pi/2$, $\pi/6$, $\pi/3$, and the volume is $0.84578\dots$. Thus the four tetrahedra $T_{a_1}^1$, $T_{a_1}^2$, $T_{C-a_{n-1}}^1$, and $T_{C-a_{n-1}}^2$ each have volume at least $0.84578$, and the volume of this structure satisfies
\[ \mathcal{V} > 2\vtet (\tw(K)-3) + 4(0.84578) > 2\vtet \tw(K) - 2.7066. \]

Finally, check that when $\tw(K) = 2$, $\mathcal{V} > 2(0.84578)$ still satisfies the theorem.
\end{proof}

%%%%%%%%%%%%%%%%%%%%%%%%%%%%%%%%%%%%%%%%%%%%%%%%%%%%%%%%%%%%%%%%%
\section{Negatively curved metrics and Dehn filling}

We now turn our attention to a new technique for bounding volume from below that has not arisen in previous chapters. This bound comes from differential geometric methods, in particular from finding volumes of hyperbolic manifolds under families of metrics, and showing that the hyperbolic metric maximizes volume among such a family. 

A hyperbolic manifold has constant sectional curvature equal to $-1$. The metrics we will consider will have negative sectional curvature, not necessarily constant. If a 3-manifold admits such a metric, it actually follows from the Geometrization theorem that the manifold also admits a hyperbolic metric; see \refthm{NegCurvedIsHyp}. However, the proof of that fact gives no information on how the hyperbolic metric relates to the negatively curved one. Often we can build an explicit negatively curved metric, and we will use this metric to make conclusions about the hyperbolic geometry of the manifold. This section presents a number of results along these lines, particularly relating to volume.

\begin{theorem}\label{Thm:NegCurvedIsHyp}
Suppose $M$ is a compact orientable 3-manifold whose interior admits a Riemannian metric with negative sectional curvature. Then $M$ admits a hyperbolic metric.
\end{theorem}

\begin{proof}
Because sectional curvature is negative, the Cartan--Hadamard theorem implies that the universal cover of $M$ is homeomorphic to $\RR^3$ (see for example \cite[Theorem~3.87]{GallotHulinLafontaine}). It follows that the fundamental group of $M$ is infinite. It also follows that $M$ is irreducible,\index{irreducible} for any sphere in $M$ lifts to a sphere in $\RR^3$, which bounds a ball. Then the image of the sphere in $M$ bounds the image of that ball in $M$, which is a ball.

In the case that $M$ is closed, because $M$ has strictly negative curvature, it is known that every abelian subgroup of its fundamental group is cyclic~\cite{Preissmann}. Thus there is no $\ZZ\times\ZZ$ subgroup of its fundamental group. So in this case, $M$ is irreducible,\index{irreducible} with $\pi_1(M)$ infinite, containing no $\ZZ\times\ZZ$ subgroup. By the Geometrization Theorem, \refthm{Geometrization}, $M$ admits a hyperbolic structure.

In the case that $M$ has boundary, there may be a $\ZZ\times\ZZ$ subgroup of $\pi_1(M)$, but the fact that the curvature is strictly negative in the interior implies that the subgroup is peripheral, hence $M$ is atoroidal~\cite{BallmannGromovSchroeder}.\index{atoroidal} If $M$ is a Seifert fibered space with infinite fundamental group, then $\pi_1(M)$ contains a cyclic normal subgroup (\refthm{SeifertFibered}). But again, a complete negatively curved finite volume Riemannian manifold cannot have a cyclic normal subgroup~\cite{BallmannGromovSchroeder}. It follows that $M$ is hyperbolic. 
\end{proof}

Notice that the proof of \refthm{NegCurvedIsHyp} gives no information on the relationship between the negatively curved metric and the hyperbolic metric. For example, we can make no conclusions about the difference in volumes of the manifolds. In the closed case, work of Besson, Courtois, and Gallot \cite{BessonCourtoisGallot} can be used to bound the volume of a manifold under one negatively curved metric in terms of the volume of another. This was extended to the finite volume case by Boland, Connell, and Souto \cite{BolandConnellSouto}. In 3-dimensions, a special case of their work is the following theorem. 

\begin{theorem}\label{Thm:BolandConnellSouto}
Let $\sigma$ and $\sigma'$ be two complete, finite volume Riemannian metrics on the same 3-manifold $N$. Suppose the Ricci curvature of $\sigma$ satisfies ${\mathrm{Ric}_\sigma} \geq -2\sigma$, and suppose the sectional curvatures of $\sigma'$ lie in the interval $[-a,-1]$ for some constant $a\geq 1$. Then
\[ \vol(N,\sigma) \geq \vol(N,\sigma'), \]
with equality if and only if both metrics are hyperbolic. \qed
\end{theorem}

Our goal in this section is to bound the change in volume of a hyperbolic manifold under Dehn filling by constructing a negatively curved metric on the Dehn filling of a hyperbolic manifold, and then applying \refthm{BolandConnellSouto}.
This volume estimate was first obtained in \cite{fkp:dfvjp}. Along the way we will obtain additional important consequences, for example the $2\pi$-theorem \cite{bleiler-hodgson}. 

%%%%%%%%%%%%%%%%%%%%%%%%%%%%%%%%%%%%%%%%%%%%%%%%%%%%%%%%%%%%%%%%%
\subsection{Negatively curved metrics on a solid torus}
We will construct metrics in this subsection, and use them to bound volume. The arguments here require a little more familiarity with Riemannian geometry than the rest of the book so far. However, these arguments are only needed in this section and will not be required elsewhere in the book. Thus a reader disinclined to work carefully through the calculations in Riemannian geometry at this time may accept the statements of the main results here and skip ahead to their applications in \refsubsec{ApplicationsNegCurved}. 

The metrics we construct will have constant sectional curvatures away from a collection of solid tori, namely those we glue to perform Dehn filling. Within a solid torus, we will use cylindrical coordinates.

\begin{definition}\label{Def:CylindricalCoords}
  Let $V$ be a solid torus, and let $\widetilde{V}$ be its universal cover. The \emph{cylindrical coordinates}\index{cylindrical coordinates} on $\widetilde{V}$ are given by $(r,\mu,\lambda)$, where $r\leq 0$ is the radial distance measured outward from $\bdy V$, $0\leq \mu\leq 1$ is measured around each meridional circle, and $-\infty<\lambda<\infty$ is measured in the longitudinal direction, orthogonal to $\mu$. Normalize the coordinates so that the generator of the deck transformation group on $\widetilde{V}$ changes the $\lambda$ coordinate by $1$. These coordinates descend to \emph{cylindrical coordinates} on $V$, where $r\leq 0$ is radial distance measured outward from $\bdy V$, $0\leq \mu \leq 1$ is in the meridional direction, and $0\leq \lambda \leq 1$ is measured orthogonal to $\mu$.
\end{definition}

\begin{lemma}\label{Lem:BleilerHodgson}
Let $(r,\mu,\lambda)$ be cylindrical coordinates on a solid torus or its universal cover. Then a metric of the form
\begin{equation}\label{Eqn:NegCurvedSolidTorus}
  ds^2 = dr^2 + f(r)^2d\mu^2 + g(r)^2d\lambda^2,
\end{equation}
where $f\from \RR\to \RR$ and $g\from \RR \to \RR$ are smooth functions of $r$,
satisfies the property that all sectional curvatures are convex combinations of
\[ -\frac{f''}{f}, \quad -\frac{g''}{g}, \quad -\frac{f'g'}{fg}. \]
Moreover, the metric is nonsingular if $f'(r_0)=2\pi$, where $r_0<0$ is the root of $f$ nearest $0$ (if it exists). 
\end{lemma}

\begin{proof}
The proof will be a standard calculation from Riemannian geometry, following \cite{bleiler-hodgson}.

For notational convenience, set $r=x_1$, $\mu=x_2$, $\lambda=x_3$. Our Riemannian metric can be written in coordinates as
\[ (g_{ij}) = \left(\begin{array}{ccc}
  1& 0 & 0 \\
  0 & f(r)^2 & 0 \\
  0 & 0 & g(r)^2
\end{array}\right) \quad \mbox{and} \quad
(g^{ij}) = \left(\begin{array}{ccc}
  1&0&0\\
  0&f(r)^{-2}& 0\\
  0&0& g(r)^{-2}
\end{array}\right)
\]
The Christoffel symbols $\Gamma_{ij}^k = \sum_\ell \Gamma_{ij\ell}g^{\ell k}$ can be computed using
\[ \Gamma_{ijk}= \frac{1}{2}\left(\frac{\partial g_{jk}}{\partial x_i}+ \frac{\partial g_{ik}}{\partial x_j}- \frac{\partial g_{ij}}{\partial x_k}\right).\]
Most of the 27 $\Gamma_{ijk}$ are zero; the non-zero ones are $\Gamma_{122} = f\cdot f'$, $\Gamma_{133} = g\cdot g'$, $\Gamma_{212}=f\cdot f'$, $\Gamma_{221}=-f\cdot f'$, $\Gamma_{313} = g\cdot g'$, and $\Gamma_{331}= -g\cdot g'$.
We then obtain the connection
$\nabla_{\partial/\partial x_i}(\partial/\partial x_j) = \sum_k \Gamma_{ij}^k\cdot \partial/\partial x_k$ as follows.

\begin{tabular}{cc|c cc}
  & && $j$ & \\
  & $\nabla_{\partial/\partial x_i}(\partial/\partial x_j)$ & 1& 2& 3\\
  \hline
  & 1 &0& $f'/f\cdot\partial/\partial x_2$ & $g'/g\cdot\partial/\partial x_3$\\
  $i$& 2& $f'/f\cdot\partial/\partial x_2$ & $-f\cdot f'\cdot\partial/\partial x_1$ & 0 \\
  & 3& $g'/g\cdot\partial/\partial x_3$ & 0 & $-f\cdot f'\cdot \partial/\partial x_1$
\end{tabular}

The Riemannian curvature tensor is given by
\[R(X,Y,Z) = \nabla_Y\nabla_XZ-\nabla_X\nabla_YZ+\nabla_{[X,Y]}Z,\]
and the sectional curvatures
\[ K(X,Y) = -\frac{\langle R(X,Y,X),Y\rangle}{|X|^2|Y|^2-\langle X,Y\rangle^2} \]
are all convex combinations of the three sectional curvatures
\[ K_{ij}=K(\partial/\partial x_i, \partial/\partial x_j) \]
for $\{i,j\}\subset\{1,2,3\}$.
We compute
\begin{align*}
  K_{12} &=  -\frac{\langle R(\partial/\partial x_1, \partial/\partial x_2, \partial/\partial x_1), \partial/\partial x_2\rangle}{\langle\partial/\partial x_1, \partial/\partial x_1\rangle\langle\partial/\partial x_2, \partial/\partial x_2\rangle-\langle\partial/\partial x_1, \partial/\partial x_2\rangle^2}
  \\[2mm]
  &= -\frac{\langle\nabla_{\partial/\partial x_1}\cdot\nabla_{\partial/\partial x_2}(\partial/\partial x_1) - \nabla_{\partial/\partial x_2}\cdot\nabla_{\partial/\partial x_1}(\partial/\partial x_1),\partial/\partial x_2\rangle}{1\cdot f^2 - 0^2}\\[2mm]
  &= -\frac{\langle\nabla_{\partial/\partial x_1}(f'/f\cdot\partial/\partial x_2),\partial/\partial x_2\rangle}{f^2} \\[2mm]
  &=-\frac{f''/f\cdot\langle\partial/\partial x_2, \partial/\partial x_2\rangle}{f^2} \\[2mm]
  &=-f''/f.
\end{align*}
A symmetric calculation shows $K_{13}= -g''/g$.
Finally
\begin{align*}
  K_{23} &= -\frac{\langle\nabla_{\partial/\partial x_2}\cdot\nabla_{\partial/\partial x_3}(\partial/\partial x_2) - \nabla_{\partial/\partial x_3}\cdot\nabla_{\partial/\partial x_2}(\partial/\partial x_2),\partial/\partial x_3\rangle}{f^2 \cdot g^2}\\[2mm]
  &= -\frac{\langle-\nabla_{\partial/\partial x_3}(-f'\cdot f\cdot\partial/\partial x_1),\partial/\partial x_3\rangle}{f^2\cdot g^2} \\[2mm]
  &=-\frac{f'\cdot f\cdot g'/g\langle\partial/\partial x_3, \partial/\partial x_3\rangle}{f^2} \\[2mm]
  &=-\frac{f'\cdot g'}{f\cdot g}.
\end{align*}

Finally, to ensure the metric is nonsingular, it must have a cone angle\index{cone angle} of $2\pi$ along the core, i.e.\ at the point $r=r_0<0$ nearest $0$ such that $f(r_0)=0$. If $f(r_0)=0$, then
\[ f'(r_0) = \lim_{r\to r_0} \frac{1}{r-r_0} \int_0^1 f(r)d\mu \]
gives the cone angle\index{cone angle} along the core circle of the solid torus. Thus we must ensure that $f'(r_0)=2\pi$.
\end{proof}

\begin{lemma}\label{Lem:SolidTorusMetric}
  Suppose $V$ is a solid torus with a prescribed Euclidean metric on $\bdy V$ such that a Euclidean geodesic representing a meridian has length $\ell_1>2\pi$. Then there exists a smooth Riemannian metric on $V$ that is hyperbolic on a collar neighborhood of $\bdy V$, has negative sectional curvature elsewhere, and the restriction of the metric to $\bdy V$ gives the prescribed Euclidean metric. 
\end{lemma}

\begin{proof}
Let $\widetilde{V}$ denote the universal cover of $V$. We will assign a metric to $\widetilde{V}$ that has the form of \refeqn{NegCurvedSolidTorus}. The functions $f$ and $g$ must satisfy a number of properties. 

To obtain the prescribed Euclidean metric, we must have $f(0)=\ell_1$ and $g(0)=\ell_2$ where $\ell_2 = \area(V)/\ell_1$. Then the deck transformation group on $\widetilde{V}$ is generated by
\[ (r,\mu,\lambda)\mapsto (r,\mu+\theta,\lambda+1), \]
where $\theta\in[0,1)$ is appropriately chosen so that the fundamental domain of $\bdy V$ has the correct shape. The metric on $\widetilde{V}$ descends to give a smooth metric on $V$.

In order for the metric to be hyperbolic near $\bdy V$, $f$ and $g$ must give sectional curvatures equal to $-1$ near $r=0$. This will hold if $f(r)=\ell_1e^r$ and $g(r)=\ell_2e^r$ near $r=0$. To ensure it is nonsingular, we need to ensure $f'(r_0)=2\pi$ where $r_0$ is the negative root of $f$ nearest $0$.

To finish the proof of the lemma, we need to show that there exist functions $f$ and $g$ satisfying the above properties. For purposes of this lemma, it suffices to choose $r_0$ such that $-\ell_1/2\pi < r_0 < -1$ and define $f$ and $g$ near $r=r_0$ by $f(r)=2\pi\sinh(r-r_0)$ and $g(r)=b\cosh(r-r_0)$, for $0<b<\ell_2$. Note that $f$ and $g$ give a metric of constant curvature $-1$ near $r_0$, and that $f'(r_0)=2\pi$, so the metric will be nonsingular. To see that definitions of $f$ and $g$ can be extended, note that the tangent line to $f$ at $r=r_0$ runs through the points $(r_0,0)$ and $(0,-2\pi r_0)$, and the tangent line to $f$ at $r=0$ runs through points $(0, \ell_1)$ and $(-1,0)$. Because $-2\pi r_0<\ell_1$ and $r_0<-1$, the function $f$ can be extended to be strictly convex, increasing, positive and smooth on $r_0<r<0$; see \reffig{2PiGraph}.
\begin{figure}
  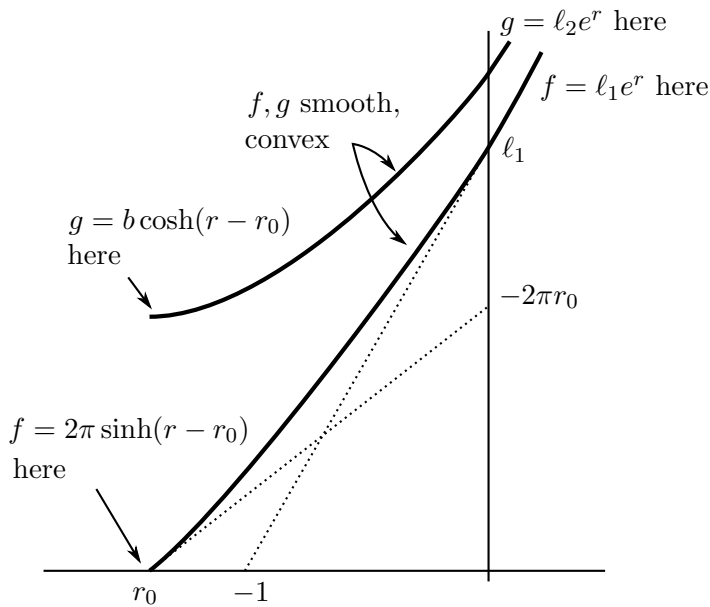
  \caption{Extending $f$ and $g$ to be strictly convex, increasing, positive, smooth functions on $r_0<r<0$}
  \label{Fig:2PiGraph}
\end{figure}
Similarly, $g'(r_0)=0<\ell_2=g'(0)$, and $0<b=g(r_0)<\ell_2=g(0)$, so $g$ can be extended to be strictly convex, increasing, positive and smooth on $r_0<r<0$. This gives the desired negatively curved metric. 
\end{proof}

\begin{theorem}[$2\pi$-Theorem]\label{Thm:2PiThm}\index{$2\pi$-theorem}
  Suppose $M$ is a hyperbolic 3-manifold with disjoint embedded cusps $C_1, \dots, C_n$ and slopes $s_j$ on $C_j$ such that a geodesic representative of each $s_j$ on $\bdy C_j$ has length strictly greater than $2\pi$ in the induced Euclidean metric. Then the Dehn filled\index{Dehn filling} manifold $M(s_1, \dots, s_n)$ admits a metric of negative curvature. Thus it is hyperbolic. 
\end{theorem}

\begin{proof}
Remove the cusps $C_1, \dots, C_n$. By \reflem{SolidTorusMetric}, there exists a negatively curved metric on a solid torus $V_j$ such that the Euclidean metric on $\bdy V_j$ agrees with that of $\bdy C_j$, and such that the metric is hyperbolic on a collar neighborhood of $\bdy V_j$. Then put a metric on $M(s_1, \dots, s_n)$ by taking the hyperbolic metric on $M-(\bigcup_{i=1}^n C_i)$, and gluing in solid tori with the metric from \reflem{SolidTorusMetric}. 
\end{proof}

Note that there is flexibility in choosing the metric of \reflem{SolidTorusMetric}. For example, the value of $b$ in the proof can be anything in a range of values. If we take a little more care to determine the metric, we can obtain bounds on additional geometric information. For example, we can determine the curvature more explicitly, using the following lemma, and bound the volume of the negatively curved solid torus, as in \reflem{VolDiffEqn}.

\begin{lemma}\label{Lem:CurvatureDiffEqn}
Let $k(r)$ be a smooth, increasing function that lies in $[0,1]$ for all $r$.
Define $f$ and $g$ to be solutions to the differential equations $f''/f=k$ and $(f'g')/(fg) = k$, subject to initial conditions $f(0)= f'(0)=\ell_1$ and $g(0)=\ell_2$. Then the function $g$ satisfies $g''/g = k + (f/f')k'$.
\end{lemma}

\begin{proof}
To check the formula for $g''/g$, note $g'/g = k f/f'$, and differentiate both sides of this equation.
\[ \frac{g''}{g} - \left(\frac{g'}{g}\right)^2 = k\left(1 + \frac{f''}{f} \left(\frac{f}{f'}\right)^2\right) + \frac{f}{f'}k'.\]
Using the fact that $g'/g = k f/f'$ and $f''/f=k$, this simplifies to the desired equation:
\[ \frac{g''}{g} = k + \frac{f}{f'}k'. \qedhere \]
\end{proof}

\begin{lemma}\label{Lem:VolDiffEqn}
Let $\ell_1>2\pi$, let $k$ be a constant function $k(r)=t\in(0,1)$, and let $f$ and $g$ be defined by the differential equations in \reflem{CurvatureDiffEqn}. Let $V$ be a solid torus with metric of \refeqn{NegCurvedSolidTorus}. Then:
\begin{enumerate}
\item Letting $r_0 = -\arctanh(\sqrt{t})/\sqrt{t}$, $f$ and $g$ have the form:
  \begin{align*}
    f(r) & = \frac{\ell_1\sqrt{1-t}}{\sqrt{t}}\sinh(\sqrt{t}(r-r_0)) \\
    g(r) & = \ell_2\sqrt{1-t}\cosh(\sqrt{t}(r-r_0))
  \end{align*}
\item At $r_0$, $f(r_0)=0$ and $f'(r_0) = \ell_1\sqrt{1-t}$. Thus the solid torus $V$ has a nonsingular metric of negative curvature $-t$ when $t=1-(2\pi/\ell_1)^2$.
\item For any $t\in(0,1)$, the volume of the (possibly singular) solid torus $V$ with metric $ds^2=dr^2+f(r)^2d\mu^2+g(r)^2d\lambda^2$ is given by
  \[ \vol(V) = \int_{r_0}^0 f(r)g(r)dr = \frac{\ell_1\ell_2}{2}. \]
\end{enumerate}
\end{lemma}

\begin{proof}
By \reflem{CurvatureDiffEqn} and \reflem{BleilerHodgson}, the metric will have negative sectional curvature. We need to show the additional properties. 
The proof is a series of calculations. First, solving the differential equation $f''/f=t$, the function $f$ has the form
\[ f(r) = c_1 e^{\sqrt{t} r} + c_2 e^{-\sqrt{t} r}. \]
Given initial conditions $f(0)=f'(0)=\ell_1$, we find
\begin{align*}
  f(r) &= \frac{\ell_1}{2}\left(1+\frac{1}{\sqrt{t}}\right) e^{\sqrt{t}r} + \frac{\ell_1}{2}\left(1-\frac{1}{\sqrt{t}}\right) e^{-\sqrt{t} r} \\
  &= \ell_1\left(\cosh(r\sqrt{t})+\frac{1}{\sqrt{t}}\sinh(r\sqrt{t})\right) \\
  &= \frac{\ell_1\sqrt{1-t}}{\sqrt{t}}\sinh(\sqrt{t}(r-r_0)),
\end{align*}
where $r_0=-\arctanh(\sqrt{t})/\sqrt{t}$, as claimed. Note that
\[ f'(r) = \ell_1\sqrt{1-t}\cosh(\sqrt{t}(r-r_0)),\]
so when $r=r_0$, we have $f(r_0)=0$ and $f'(r_0)=\ell_1\sqrt{1-t}$. Thus $f'(r_0)=2\pi$ if $t= 1-(2\pi/\ell_1)^2$, and the metric will be nonsingular in this case. 

As for $g$, we may solve $g'/g = t f'/f = \sqrt{t}\tanh(\sqrt{t}(r-r_0))$ by integration, to obtain
\[ g(r) = c_2\cosh(\sqrt{t}(r-r_0)) = \ell_2\sqrt{1-t}\cosh(\sqrt{t}(r-r_0)), \]
using the initial condition $g(0)=\ell_2$ to determine the constant $c_2$.

Finally we compute the volume of a solid torus $V$ with metric as in \refeqn{NegCurvedSolidTorus}. 
\begin{align*}
\vol(V) &= \int_{r_0}^0 f(r)g(r)\,dr \\
&= \int_{r_0}^0\frac{\ell_1\ell_2(1-t)}{\sqrt{t}}\sinh(\sqrt{t}(r-r_0))\cosh(\sqrt{t}(r-r_0)) \\
&=\left[ \frac{\ell_1\ell_2(1-t)}{2t}\sinh^2(\sqrt{t}(r-r_0)) \right]_{r=r_0}^0 \\
&=\frac{\ell_1\ell_2(1-t)}{2t}\sinh^2(\arctanh(\sqrt{t})) \\
&= \frac{\ell_1\ell_2(1-t)}{2t}\cdot\frac{t}{1-t} \\
&= \frac{\ell_1\ell_2}{2}. \qedhere
\end{align*}
\end{proof}

Now we would like to use the metric on the solid torus $V$ obtained from \reflem{VolDiffEqn}, along with \refthm{BolandConnellSouto}, to bound the volume of Dehn filled manifolds. However, at this point we have a problem. Although we have constructed a nonsingular Riemannian metric on the solid torus with nice curvature and volume, note that the metric does not give a hyperbolic metric, with sectional curvatures $-1$, on a collar neighborhood of the boundary of $V$. Thus we cannot glue the metric of \reflem{VolDiffEqn} to the metric of the cusped manifold with horoball neighborhoods removed to obtain a negatively curved metric on the Dehn filled manifold, as we did in \refthm{2PiThm}. The way to fix this problem is to do a little deeper analysis, which is done in the following lemma.

\begin{lemma}\label{Lem:FKPVolFilling}
Suppose $V$ is a solid torus with a prescribed Euclidean metric on $\bdy V$, such that a Euclidean geodesic representing a meridian has length $\ell_1>2\pi$. Let $\zeta\in (0,1)$ be a constant. Then there exists a smooth, negatively curved Riemannian metric on $V$ that satisfies the following properties.
\begin{enumerate}
\item The metric is hyperbolic on a collar neighborhood of $\bdy V$, and its restriction to $\bdy V$ gives the prescribed Euclidean metric.
\item The sectional curvatures are bounded above by $-\zeta(1-(2\pi/\ell_1)^2)$.
\item The volume of $V$ in this metric is at least $\frac{1}{2}\zeta\area(\bdy V)$.
\end{enumerate}
\end{lemma}

\begin{proof}
We use the ideas of \reflem{CurvatureDiffEqn} and \reflem{VolDiffEqn} to define $f$ and $g$ by differential equations. However, we do not choose $k(r)$ to be constant. We need $k(r)=1$ near $r=0$ to obtain the appropriate curvature estimates on the boundary of $V$. We have seen in \reflem{VolDiffEqn} that we may obtain nice volume and curvature results when $k(r)=t$ for some $t\in(0,1)$ for $r<0$. So we define $k$ to be a smooth bump function, depending on $r$, $t$, and $\epsilon>0$, as follows. If $r\leq -\epsilon$, set $k_{t,\epsilon}(r)=t$. If $r\geq -\epsilon/2$, set $k_{t,\epsilon}(r)=1$. For $r$ between $-\epsilon$ and $-\epsilon/2$, the function $k_{t,\epsilon}(r)$ is smooth and strictly increasing. See \reffig{BumpFunction} for a typical graph.

\begin{figure}
%% Accompanies image file 'F13-02-kter.eps' (pdf, eps, ps)
\begingroup%
  \makeatletter%
  \providecommand\color[2][]{%
    \errmessage{(Inkscape) Color is used for the text in Inkscape, but the package 'color.sty' is not loaded}%
    \renewcommand\color[2][]{}%
  }%
  \providecommand\transparent[1]{%
    \errmessage{(Inkscape) Transparency is used (non-zero) for the text in Inkscape, but the package 'transparent.sty' is not loaded}%
    \renewcommand\transparent[1]{}%
  }%
  \providecommand\rotatebox[2]{#2}%
  \newcommand*\fsize{\dimexpr\f@size pt\relax}%
  \newcommand*\lineheight[1]{\fontsize{\fsize}{#1\fsize}\selectfont}%
  \ifx\svgwidth\undefined%
    \setlength{\unitlength}{174.37249374bp}%
    \ifx\svgscale\undefined%
      \relax%
    \else%
      \setlength{\unitlength}{\unitlength * \real{\svgscale}}%
    \fi%
  \else%
    \setlength{\unitlength}{\svgwidth}%
  \fi%
  \global\let\svgwidth\undefined%
  \global\let\svgscale\undefined%
  \makeatother%
  \begin{picture}(1,0.4927441)%
    \lineheight{1}%
    \setlength\tabcolsep{0pt}%
    \put(0,0){\includegraphics[width=\unitlength]{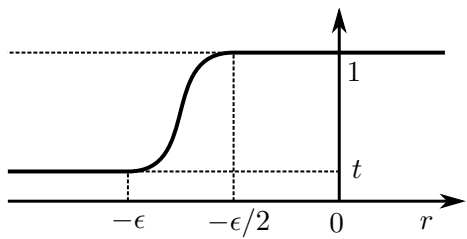}}%
    \put(0.73545451,0.32974792){\color[rgb]{0,0,0}\makebox(0,0)[lt]{\lineheight{1.25}\smash{\begin{tabular}[t]{l}$1$\end{tabular}}}}%
    \put(0.74631641,0.11903559){\color[rgb]{0,0,0}\makebox(0,0)[lt]{\lineheight{1.25}\smash{\begin{tabular}[t]{l}$t$\end{tabular}}}}%
    \put(0.22626091,0.01172578){\color[rgb]{0,0,0}\makebox(0,0)[lt]{\lineheight{1.25}\smash{\begin{tabular}[t]{l}$-\epsilon$\end{tabular}}}}%
    \put(0.43328314,0.00760466){\color[rgb]{0,0,0}\makebox(0,0)[lt]{\lineheight{1.25}\smash{\begin{tabular}[t]{l}$-\epsilon/2$\end{tabular}}}}%
    \put(0.89272972,0.00933556){\color[rgb]{0,0,0}\makebox(0,0)[lt]{\lineheight{1.25}\smash{\begin{tabular}[t]{l}$r$\end{tabular}}}}%
    \put(0.69765018,0.00369141){\color[rgb]{0,0,0}\makebox(0,0)[lt]{\lineheight{1.25}\smash{\begin{tabular}[t]{l}$0$\end{tabular}}}}%
  \end{picture}%
\endgroup%

  \caption{Graph of $k_{t,\epsilon}(r)$.}
  \label{Fig:BumpFunction}
\end{figure}
Then $k$ is continuous in the three variables $t, \epsilon,r$. 
We also define $k_{t,0}(r)$ to be the step function
\[ k_{t,0}(r) = \lim_{\epsilon\to 0^+} k_{t,\epsilon}(r) =
\begin{cases}
  t & \mbox{ if } r<0, \\ 1 & \mbox{ if } r\geq 0.
\end{cases}
\]

Now for $\epsilon\geq 0$ and $t\in (0,1)$, define $f_{t,\epsilon}$ and $g_{t,\epsilon}$ by the differential equations
\[ \frac{f''_{t,\epsilon}(r)}{f_{t,\epsilon}(r)} = k_{t,\epsilon}(r), \quad
\frac{g'_{t,\epsilon}(r)}{g_{t,\epsilon}(r)} = k_{t,\epsilon}(r)\frac{f_{t,\epsilon}(r)}{f'_{t,\epsilon}(r)}. \]

The family of functions $f_{t,\epsilon}(r)$ and $g_{t,\epsilon}(r)$ can be shown to have a number of nice properties. Away from $\epsilon=0$, these mostly follow by standard facts in differential equation. As $\epsilon\to 0^+$, a little more analysis is required, which we will omit here. For full details see \cite{fkp:dfvjp}. In particular, the following hold.

\textbf{Nonsingularity:}
For all $t\in(0,1)$ and $\epsilon\geq 0$, $f_{r,\epsilon}(r)$ has a unique root $r_0(t,\epsilon)$. The function $f'_{t,\epsilon}(r_0(t,\epsilon))$ is continuous in $t$ and $\epsilon$, and strictly decreasing in both variables. For every $t$ between $0$ and $1-(2\pi/\ell_1)^2<1$, there is a unique value $\epsilon(t)>0$ such that $f'_{t, \epsilon(t)}(r_0(t, \epsilon(t)) = 2\pi$. This gives a nonsingular metric for every $t$. Moreover, as $t\to 1-(2\pi/\ell_1)^2$, $\epsilon(t)\to 0$.

Now let $t\in (0,1)$ and define $\tau(t)$ to be the nonsingular Riemannian metric given by the functions $f_t(r) = f_{t, \epsilon(t)}(r)$ and $g_t(r)=g_{t,\epsilon(t)}(r)$. 

\textbf{Sectional curvatures:}
The metric $\tau(t)$ has all sectional curvatures bounded above by $-t$.
This follows from \reflem{CurvatureDiffEqn}, along with the fact that
the function $f_t(r)$ is positive and increasing. Thus $f_t(r)/f'_t(r)$ is positive. Moreover, $k'_{t,\epsilon(t)}(r)$ is positive, since $k_{t,\epsilon(t)}(r)$ is increasing with $r$. Thus \reflem{CurvatureDiffEqn} implies that $g''_t(r)/g_t(r)$ is least $k_{t,\epsilon(t)}(r)\geq t$. By definition, $f''_t(r)/f_t(r)$ and $(f'_t(r)g'_t(r))/(f_t(r)g_t(r))$ are equal to $k_{t,\epsilon(t)}(r)\geq t$. So all sectional curvatures are bounded above by $-t$. 

\textbf{Volumes:}
Recall we have fixed $\zeta>0$. For notational purposes, define $t_0$ to be $t_0=1-(2\pi/\ell_1)^2$. Let $t$ lie in the interval $(\zeta t_0,t_0)$. Then for the metric $\tau(t)$, we have
\[\lim_{t\to t_0} \vol(V, \tau(t)) = \frac{\ell_1\ell_2}{2} = \frac{1}{2}\area{\bdy V}. \]
This follows from the fact that $f_t$ and $g_t$ converge uniformly to $f_{t_0,0}$ and $g_{t_0,0}$ as $t\to t_0$. Moreover, $r_0(t,\epsilon(t))$ converges to $r_0(t_0,0)$. Then the limit must be the limit of the differential equation in the case $k$ is constant, which we computed in \reflem{VolDiffEqn}. In particular, we have
\[ \lim_{t\to t_0} \vol(V, \tau(t)) = \vol(V, t_0) = \frac{\ell_1\ell_2}{2}.\]

To finish the proof of the lemma, select $t \in (\zeta t_0, t_0)$ near enough to $t_0$ so that $\vol(V, \tau(t))\geq \frac{1}{2}\zeta\area(\bdy V)$. For this metric, sectional curvatures are bounded above by $-t \leq -\zeta t_0 = -\zeta(1-(2\pi/\ell_1)^2)$. Finally, the metric is nonsingular, and by choice of bump function and initial conditions, on a collar neighborhood of $\bdy V$ it is hyperbolic, with metric agreeing with the prescribed metric on $\bdy V$. 
\end{proof}

We are now ready to prove the main result of this section. 

\begin{theorem}[Volume change under Dehn filling]\label{Thm:VolumeChangeUnderFilling}
Let $M$ be a complete, finite volume hyperbolic manifold with cusps. Suppose $C_1, \dots, C_n$ are disjoint embedded cusps with slopes $s_j$ on $C_j$ such that a geodesic representative of $s_j$ on $\bdy C_j$ has length strictly greater than $2\pi$. Denote the minimal slope length by $\ell_{\min}$. Then the Dehn filled\index{Dehn filling}\index{volume bound!lower bound on Dehn filling} manifold $M(s_1, \dots, s_n)$ is a hyperbolic manifold with
\[ \vol(M(s_1, \dots, s_n) \geq \left( 1 - \left(\frac{2\pi}{\ell_{\min}}\right)^2\right)^{3/2} \vol(M). \]
\end{theorem}

\begin{proof}
Fix an arbitrary constant $\zeta \in (0,1)$. Replace each cusp $C_j$ by a solid torus $V_j$ whose meridian is $s_j$. By \reflem{FKPVolFilling}, the smooth Riemannian metric $\tau_j$ on $V_j$ agrees with the hyperbolic metric on $\bdy C_j$, so this gives a smooth Riemannian metric $\tau$ on $M(s_1, \dots, s_n)$. Additionally, for each $j$, sectional curvatures on $V_j$ are at most $-\zeta(1-(2\pi/\ell_{\min})^2)$, and the volume of $V_j$ is at least $\zeta\area(\bdy V_j)/2 = \zeta\vol(C_j)$ where $\vol(C_j)$ is the cusp volume in the hyperbolic metric. Note that sectional curvatures in $V_j$ are also bounded below by some constant, since $V_j$ is compact.

Thus the metric $\tau$ on $M(s_1,\dots, s_n)$ has sectional curvatures bounded above by $-\zeta(1-(2\pi/\ell_{\min})^2)$ and bounded below by some constant. Moreover
\begin{align*}
  \vol(M(s_1, \dots, s_n))& \geq\vol(M-\cup_{j=1}^n C_j) + \zeta\sum\vol(C_j) \\
  &\geq \zeta\vol(M).
\end{align*}

Rescale the metric to obtain a metric with sectional curvatures bounded above by $-1$. To do this, replace $\tau$ by $\sigma = \sqrt{\zeta (1-(2\pi/\ell_{\min})^2)}\tau$. Note this multiplies all sectional curvatures by $(\zeta(1-(2\pi/\ell_{\min})^2))^{-1}$. The volume is rescaled by a factor of $(\zeta(1-(2\pi/\ell_{\min})^2))^{3/2}$. Thus under the metric $\sigma$, sectional curvatures of $M(s_1, \dots, s_n)$ lie in $[-a, 1]$ for some $a\geq 1$, and $\vol(M(s_1, \dots, s_n), \sigma) \geq \zeta^{5/2}(1-(2\pi/\ell_{\min})^2)^{3/2}\vol(M)$.

Now let $S$ denote the set of all metrics on $M(s_1, \dots, s_n)$ whose sectional curvatures lie in the interval $[-a,-1]$. Since $\zeta$ is arbitrary, by the above work the supremum of volumes of $M(s_1, \dots, s_n)$ over all metrics in $S$ satisfies:
\[ \sup_{\sigma \in S} \vol(M(s_1, \dots, s_n), \sigma) \geq \left(1-\left(\frac{2\pi}{\ell_{\min}}\right)^2\right)^{3/2}\vol(M). \]
Here $\vol(M(s_1, \dots, s_n),\sigma)$ denotes volume under the metric $\sigma$, and $\vol(M)$ denotes the volume of $M$ under its given hyperbolic metric.

Now, \refthm{BolandConnellSouto} implies that the hyperbolic metric $\sigma_{\mathrm{hyp}}$ on the Dehn filled manifold $M(s_1, \dots, s_n)$ uniquely maximizes volume over the set $S$ of all metrics whose sectional curvatures lie in the interval $[-a,1]$. Thus:
\begin{align*}
\vol(M(s_1, \dots, s_n), \sigma_{\mathrm{hyp}}) & \geq
\sup_{\sigma \in S} \vol(M(s_1, \dots, s_n), \sigma) \\
& \geq \left(1-\left(\frac{2\pi}{\ell_{\min}}\right)^2\right)^{3/2}\vol(M). \qedhere
\end{align*}
\end{proof}

%%%%%%%%%%%%%%%%%%%%%%%%%%%%%%%%%%%%%%%%%%%%%%%%%%%%%%%%%%%%%%%%%
\subsection{Applications to knots}\label{Subsec:ApplicationsNegCurved}

Recall the definitions of \emph{twist-reduced}\index{twist-reduced}, \refdef{TwReduced} or \refdef{TwistReduced}, and \emph{twist-number}\index{twist-number}, \refdef{TwistNumber}. We will denote the twist-number of a twist-reduced diagram $K$ by $\tw(K)$. An application of \refthm{VolumeChangeUnderFilling} is the following result, which first appeared in \cite{fkp:dfvjp}.

\begin{theorem}[Volume bounds for highly twisted links]\label{Thm:VolBoundHighlyTwisted}\index{highly twisted!volume bounds}\index{highly twisted}
Let $K\subset S^3$ be a link with a prime, twist-reduced diagram. Assume the diagram has $\tw(K)\geq 2$ twist regions, and that each twist region contains at least seven crossings. Then $K$ is a hyperbolic link satisfying
\[ 0.70734 (\tw(K) -1) \leq \vol(S^3-K) < 10 \vtet (\tw(k)-1), \]
where $\vtet = 1.0149\dots$ is the volume of a regular ideal tetrahedron.
\end{theorem}

The proof of \refthm{VolBoundHighlyTwisted} uses a theorem due to Miyamoto, which in full generality gives a lower bound on the volume of an $n$-dimensional hyperbolic manifold with geodesic boundary \cite{Miyamoto}. Here, we state only the 3-dimensional case, which is the case we will use. 

\begin{theorem}[Miyamoto]\label{Thm:Miyamoto}\index{Miyamoto's theorem}
If $N$ is a hyperbolic 3-manifold with totally geodesic boundary, then $\vol(N)\geq -\voct\chi(N)$,
where $\voct =3.66\dots$ is the volume of a regular ideal octahedron,\index{ideal octahedron, regular}\index{regular ideal octahedron}\index{regular ideal octahedron!volume} with equality if and only if $N$ decomposes into $-\chi(N)$ ideal octahedra.
\end{theorem}

%%%%%%%%%%%%%%%%%%%%%%%%%%%%%%%%%%%%%%%%%%%%%%%%%%%%%%%%%%%%%%%%%
\begin{proof}[Proof sketch]
Suppose $N$ has totally geodesic boundary. Then the lift of $N$ to the universal cover $\HH^3$ is a subspace $\widetilde{N}$ of $\HH^3$ bounded by disjoint hyperplanes, where each hyperplane is a lift of the geodesic surface $\bdy N$.

Pick one such hyperplane $O$, and consider the set $D_O$ consisting of points in $\widetilde{N}$ closer to $O$ than to any other hyperplane in the lift of $\bdy N$. The set $D_O$ is convex, with boundary consisting of faces $F$ made up of points equidistant from the hyperplane $O$ and some other hyperplane. Define the truncated cone $C_F$ to be all points in $D_O$ that lie on a line running from $F$ to meet $O$ orthogonally. We may decompose all of $\widetilde{N}$ into truncated cones.

Now project to $N$. Because $\bdy \widetilde{N}$ is invariant under the action of the covering transformations, and distances along geodesics are preserved, the decomposition projects to a decomposition of $N$. The volume of $N$ is obtained by summing the volumes of all the truncated cones decomposing $N$.

The main step of the proof is to bound the ratio $\vol(C_F)/\area(O\cap C_F)$, with notation as above. To do so, \emph{truncated tetrahedra} are introduced.

Consider a combinatorial polyhedron $P$ obtained by removing from a tetrahedron a small open neighborhood of each vertex. The faces of the tetrahedra become hexagonal faces of $P$, and four new triangular faces are added. A truncated tetrahedron,\index{truncated tetrahedron} also called a \emph{hyperideal tetrahedron}\index{hyperideal tetrahedron} in the literature, is a compact polyhedron in hyperbolic space that realizes $P$, such that the triangle faces and hexagonal faces are totally geodesic, and such that hexagons meet triangles at right angles. See \reffig{TruncTet}, left.

\begin{figure}
\includegraphics{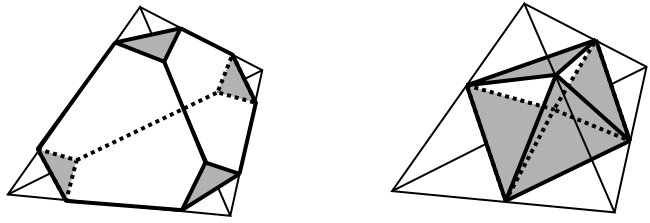}
  \caption{Left: a truncated tetrahedron, or hyperideal tetrahedron. Right: when lengths of edges between triangles go to zero, the truncated tetrahedron becomes a regular ideal octahedron.\index{ideal octahedron, regular}\index{regular ideal octahedron}}
  \label{Fig:TruncTet}
\end{figure}

It can be shown that a truncated tetrahedron is determined by the lengths of the six edges between its triangular faces. When each of these is the same length, we say the truncated tetrahedron is regular, and we denote the regular truncated tetrahedron of edge length $r$ by $T_r$. Denote its four triangular faces by $\tau_1, \dots, \tau_4$. We will be considering the ratio
\[ \rho(r) = \frac{\vol(T_{2r})}{\area(\tau_1\cup\tau_2\cup\tau_3\cup\tau_4)}. \]
Observe that $r\geq 0$, and when $r$ approaches $0$, the triangles $\tau_j$ each become ideal triangles,\index{ideal triangle} and the truncated tetrahedron $T_0$ becomes a regular ideal octahedron.\index{ideal octahedron, regular}\index{regular ideal octahedron} See \reffig{TruncTet}, right. Miyamoto shows that $\rho(r)$ increases with $r$~\cite[Lemma~2.1]{Miyamoto}. Thus $\rho(r)\geq\rho(0)$. 

Now consider again the truncated cone $C_F$ and its truncation face $O$. Consider geodesics in $N$ with both endpoints on $\bdy N$ and orthogonal to $\bdy N$. Such a geodesic is called a \emph{return path}. Let $\ell$ be the length of the shortest return path in $N$; note $\ell\geq 0$. The return path lifts to a collection of geodesics in $\widetilde{N}$, each running between hyperplane lifts of $\bdy N$, each perpendicular to the hyperplane, and each of length $\ell$. If such a geodesic meets the cone $C_F$, its intersection with $C_F$ has length $\ell/2$.

The main technical result in \cite{Miyamoto} is a proof that
\begin{equation}\label{Eqn:Miyamoto}
  \frac{\vol(C_F)}{\area(O\cap C_F)} \geq \frac{\vol(T_{\ell})}{\area(\tau_1\cup\tau_2\cup\tau_3\cup\tau_4)} = \rho(\ell/2).
\end{equation}
The proof is obtained by subdividing $C_F$ into pieces, matching pieces making up $T_{\ell}$, and observing relationships between volume and edge lengths for such pieces.

Assuming \refeqn{Miyamoto}, we complete the proof. When $N$ has shortest return path of length at least $\ell$,
\[
\vol(N) = \sum_C \vol(C) \geq \rho\left(\frac{\ell}{2}\right)\sum_C \area(C\cap O) =  \rho\left(\frac{\ell}{2}\right)\area(\bdy N),\]
where the sum is over all truncated cones $C$ in the decomposition, and $C\cap O\subset \bdy N$ denotes the portion of $C$ on $\bdy N$.

Finally, note that the shortest return path always has length at least $\ell=0$. Using the fact that $\rho$ is increasing, the above equation becomes
\[ \vol(N) \geq \rho(0)\area(\bdy N) = \frac{\voct}{4\pi}\cdot (-2\pi\chi(\bdy N)) = -\voct\chi(N).
\]
Here we are using the fact that $\vol(T_0)=\voct$, that the area of an ideal triangle\index{ideal triangle} is $\pi$, the Gauss--Bonnet formula $\area(\bdy N) = -2\pi\chi(\bdy N)$, and the fact that for a 3-manifold $N$ with boundary, $\chi(\bdy N) =2\chi(N)$.
\end{proof}
%%%%%%%%%%%%%%%%%%%%%%%%%%%%%%%%%%%%%%%%%%%%%%%%%%%%%%%%%%%%%%%%%

Given Miyamoto's theorem, we prove \refthm{VolBoundHighlyTwisted}.

\begin{proof}[Proof of \refthm{VolBoundHighlyTwisted}]
The fact that the link $K$ is hyperbolic follows from \refthm{FuterPurcellFilling}. The upper bound on volume comes from \refthm{VolUpperTwistRegions}.

To obtain the lower bound, we will consider fully augmented links.\index{fully augmented link} Let $L$ be the fully augmented link obtained by adding a crossing circle encircling each twist region of $K$. By \refthm{AugSlopeLengths}, $S^3-K$ is obtained from $S^3-L$ by performing Dehn fillings on crossing circles, along a slopes of length at least $\sqrt{7^2+1}=\sqrt{50} > 2\pi$.

We will find a lower bound on the volume of $S^3-L$. To do so, first remove all half-twists from the diagram of $L$. That is, recall $L$ may have single crossings at twist regions. Replace $L$ with a new fully augmented link $L'$ that has no crossing at twist regions. Note the complement of $L'$ is obtained from that of $L$ by cutting along 2-punctured disks bounded by crossing circles and regluing, thus it follows from \refcor{Gluing3PunctSpheres} that the volume of $S^3-L'$ is identical to the volume of $S^3-L$.

Now cut $S^3-L'$ along the plane of projection, separating it into two identical pieces, each with totally geodesic boundary coming from the white surface. Call one of these $M$. By Miyamoto's theorem, $\vol(M) \geq -\voct\chi(M)$. Note that $M$ is homeomorphic to a ball in $S^3$ with a tube drilled out for each crossing circle, and there are $\tw(K)$ crossing circles. 
Thus the Euler characteristic\index{Euler characteristic} of $M$ is $\chi(M) = (1-\tw(K))$. Because we form $S^3-L'$ by taking two copies of $M$, the volume satisfies
\[ \vol(S^3-L) = \vol(S^3-L') \geq -2\voct\chi(M) = 2\voct(\tw(K)-1). \]

Now by \refthm{VolumeChangeUnderFilling} (volume change under Dehn filling), the volume of $S^3-K$ satisfies:
\begin{align*}
  \vol(S^3-K) &\geq \left( 1-\left(\frac{2\pi}{\sqrt{50}}\right)^2\right)^{3/2} \vol(S^3-L) \\
  & \geq \left(1-\frac{2\pi^2}{25}\right)^{3/2} 2\voct(\tw(K)-1) \\
  & \geq 0.70735(\tw(K)-1). \qedhere
\end{align*}
\end{proof}

%%%%%%%%%%%%%%%%%%%%%%%%%%%%%%%%%%%%%%%%%%%%%%%%%%%%%%%%%%%%%%%%%
\section{Volume, guts, and essential surfaces}

\Refthm{VolBoundHighlyTwisted} gives volume bounds highly twisted\index{highly twisted} knots and links, but only with at least seven crossings per twist region. A similar result holds for alternating knots and links,\index{alternating knot or link} without a restriction on the number of crossings for twist regions, originally due to Lackenby \cite{lackenby:alt-volume}. The method of proof is different, but illustrates another tool for bounding hyperbolic volume from below, developed by Agol, Storm, and Thurston \cite{ast:volumes}. In this section, we explain the tool, and use it to bound volumes of alternating links.

The main theorem of the section is the following. 

\begin{theorem}[Volume bounds for alternating links]\label{Thm:VolAlt}\index{alternating knot or link!volume bounds}
  Let $K$ be a hyperbolic knot or link with a twist-reduced alternating diagram\index{alternating diagram} with twist number $\tw(K)$. Then
  \[ \vol(S^3-K) \geq \frac{\voct}{2}(\tw(K) - 2). \]
\end{theorem}

We will prove \refthm{VolAlt} by considering again the checkerboard surfaces of the alternating link, and the bounded polyhedral decomposition of the link complement cut along those surfaces, \refthm{PolyAltKnot} and \reflem{BoundedPolyAlt}. To describe our main tool, we need additional terminology.

First, recall the JSJ-decomposition of a 3-manifold, \refthm{JSJ} and \refdef{JSJ}. We will apply a special form of this decomposition to a 3-manifold $M$ cut along an essential\index{essential} surface $S$. Recall that $M\cut S$ is the closure of the manifold obtained by removing a regular neighborhood of $S$ (\refdef{Cut}). Its boundary consists of components of the parabolic locus,\index{parabolic locus} which are remnants of the torus boundary components of $M$, and $\widetilde{S}$.

\begin{definition}\label{Def:Double}
Let $M$ be a compact 3-manifold with torus boundary components, and let $S$ be an essential\index{essential} surface properly embedded in $M$. The \emph{double of $M\cut S$},\index{double of $M\cut S$} denoted $D(M\cut S)$ is the manifold obtained by taking two copies of $M\cut S$ and gluing them by the identity map on $\widetilde{S}$. 
\end{definition}

Note the double of $M\cut S$ will have torus boundary components coming from the parabolic locus of the boundary of $M\cut S$. We will consider the JSJ-decomposition of the double, as in \refdef{JSJ}.\index{JSJ-decomposition}

\begin{lemma}\label{Lem:AnnulusDecomposition}
  Let $M$ be a hyperbolic 3-manifold, homeomorphic to the interior of a compact manifold with torus boundary. Let $S$ be a properly embedded essential\index{essential} surface in $M$. Consider the double $D(M\cut S)$, and let $\calT$ denote the  JSJ-decomposition of $D(M\cut S)$. Finally, slice $\calT$ and $D(M\cut S)$ along $\widetilde{S}$, obtaining two copies of $M\cut S$. The following hold. 
  \begin{enumerate}
  \item The tori in the collection $\calT$ can be isotoped to be preserved by the reflection of $D(M\cut S)$ in the surface $\widetilde{S}$; thus cutting $D(M\cut S)$ along $\widetilde{S}$ cuts $\calT$ into two identical pieces. 
  \item Each essential torus $T\in \calT$ is sliced into essential annuli in $M\cut S$ with boundary on $\widetilde{S}$.
  \item The characteristic submanifold\index{characteristic submanifold} of $D(M\cut S)$ intersects $M\cut S$ in components that are either $I$-bundles\index{$I$-bundle} over a subsurface of $\widetilde{S}$, or Seifert fibered solid tori. 
  \end{enumerate}
\end{lemma}

\begin{proof}
  The first item follows from the equivariant torus theorem, due to Holzmann \cite{holzmann}. It is an exercise to prove the remaining two items; \refex{AnnularJSJ}. 
\end{proof}

\begin{definition}\label{Def:Guts}
  Let $M$, $S$, $\calT$ be as in \reflem{AnnulusDecomposition}. We say that the intersection of the characteristic submanifold\index{characteristic submanifold} of $D(M\cut S)$ with $M\cut S$ is the \emph{characteristic submanifold of $M\cut S$}. Its complement in $M\cut S$ is the \emph{guts}\index{guts} of $S$, denoted $\guts(M\cut S)$ or sometimes simply $\guts(S)$. 
\end{definition}

\begin{example}
  Recall from \refexamp{Fig8Fiber} that the figure-8 knot complement $M$ contains a surface $S$ that is a fiber, shown in \reffig{Fig8Seifert}. This surface $S$ is essential\index{essential} and properly embedded. The manifold $M\cut S$ is homeomorphic to $S\times I$. Thus in this case, all of $M\cut S$ is an $I$-bundle.\index{$I$-bundle} Thus $\guts(M\cut S)$ is empty. 
\end{example}

By contrast, later in this section we will find examples of alternating knot complements $M$ and surfaces $S$ such that $\guts(M\cut S) = M\cut S$, that is the characteristic submanifold is empty.

\begin{lemma}\label{Lem:GutsHyperbolic}
Let $M$, $S$, and $\calT$ be as in \reflem{AnnulusDecomposition}. Then the manifold $\guts(M\cut S)$ admits a hyperbolic metric with totally geodesic boundary.\index{guts}
\end{lemma}

\begin{proof}
  If we double $\guts(M\cut S)$ along the portion of the boundary on $\widetilde S$, i.e. take two copies of $\guts(M\cut S)$ and glue by the identity along their common boundary on $\widetilde{S}$, we obtain the complement of the characteristic submanifold of the manifold $D(M\cut S)$. This admits a finite volume hyperbolic metric. It also admits an involution fixing the surface $\guts(M\cut S)\cap \widetilde{S}$ pointwise. It follows from the proof of the Mostow--Prasad rigidity theorem that an embedded surface fixed pointwise by an involution of a finite volume hyperbolic 3-manifold must be totally geodesic. Hence cutting along it yields a hyperbolic structure on $\guts(M\cut S)$ with totally geodesic boundary. 
\end{proof}

The following theorem, from \cite{ast:volumes}, gives us a tool to bound volumes from below using the guts of surfaces.

\begin{theorem}[Agol, Storm, and Thurston]\label{Thm:astguts}
Let $S$ be a $\pi_1$-essential\index{$\pi_1$-essential} surface properly embedded in an orientable hyperbolic $3$-manifold $M$. Then
\[\vol(M)\geq-\voct\chi(\guts(M\cut S)),\]
where $\voct =3.66\dots$ is the volume of a regular ideal octahedron.\index{ideal octahedron, regular}\index{regular ideal octahedron}\index{regular ideal octahedron!volume}\index{guts}
\end{theorem}

\begin{proof}[Proof sketch.]
The essential surface $S$ can be isotoped to be minimal in $M$. Cut along the minimal surface isotopic to $S$, and denote $M\cut S$ by $N$. Note that $N$ inherits from $M$ a Riemannian metric for which the mean curvature on its boundary $\widetilde{S}$ is $0$. Denote this metric by $g$. Then $\vol(M)= \vol(N,g)$. 

Let $D(N)$ denote the double of $N$, i.e.\ the manifold obtained by taking two copies of $N$ and identifying them along their common boundary. By \reflem{GutsHyperbolic}, $D(\guts(N)) \subset D(N)$ inherits a complete hyperbolic metric with $\widetilde{S}\cap \guts(N)$ a totally geodesic surface embedded in $D(\guts(N))$; denote this metric by $h$. On the other hand, $D(N)$ inherits a singular Riemannian metric with singularities on $S$ obtained from the metric $g$; denote this metric by $g$ as well. Agol, Storm, and Thurston show that this singular metric $g$ can be approximated by smooth Riemannian metrics $\{g_i\}$ with restricted curvature properties. The volumes $\vol(D(N),g_i)$ under the smooth metrics $g_i$ converge to the volume $\vol(D(N),g)$ under its singular metric, with $\vol(D(N),g) = 2\vol(N,g) = 2\vol(M)$.

First suppose $M$ is closed. Use Ricci flow with surgery to evolve the metric, as in Perelman's proof of the geometrization theorem \cite{perelman02,perelman03}. The evolution will give $D(\guts(N))$ the hyperbolic metric $h$. Perelman's techniques imply a monotonicity result, in particular that
\[ \vol(D(N),g) \geq \vol(D(\guts(N)), h), \]
with equality if and only if $S$ is totally geodesic in $(D(N),g)$. 
Then
\[ \vol(M) = \frac{1}{2}\vol(D(N),g) \geq \vol(\guts(M\cut S)), \]
with equality if and only if $S$ is totally geodesic in $M$.

If $M$ has torus boundary, then similar techniques can be used to approximate the metric on compact sets, giving the same result.

Finally, $(\guts(M\cut S), h)$ is a hyperbolic manifold with totally geodesic boundary. Thus by Miyamoto's theorem, \refthm{Miyamoto},
\[ \vol(\guts(M\cut S)) \geq -\voct\chi(\guts(M\cut S)).\qedhere \]
\end{proof}

We will be considering the guts of the checkerboard surfaces in an alternating link. By \refthm{CheckerboardEssential}, the checkerboard surfaces are essential,\index{essential} and thus satisfy the hypotheses of \refthm{astguts}.

\begin{lemma}\label{Lem:BigonCharacteristic}
  Let $S$ be one of the checkerboard surfaces of a link with a twist-reduced alternating diagram\index{alternating diagram}\index{alternating knot or link!checkerboard surface}\index{checkerboard surface} $K$, whose hyperbolic complement we denote by $S^3-K=M$. Without loss of generality say $S$ is shaded, with $W$ the other checkerboard surface, colored white. Suppose in the polyhedral decomposition of $M$ that there are white bigon\index{bigon} faces. Then a neighborhood of each white bigon\index{bigon} face is part of an $I$-bundle\index{$I$-bundle} in $M\cut S$, and thus any bigon face of $W$ does not belong to $\guts(M\cut S)$. \index{guts}
\end{lemma}

\begin{proof}
  The $I$-bundle components\index{$I$-bundle} of $M\cut S$ have the form $Y\times I$, with $Y\times \{0\}$ and $Y\times\{1\}$ subsets of $\widetilde{S} \subset \bdy(M\cut S)$. Recall that $Y\times\{0\}$ and $Y\times\{1\}$ form the \emph{horizontal boundary}\index{horizontal boundary}\index{$I$-bundle!horizontal boundary} of the $I$-bundle. The subset $\bdy Y\times I$ forms the \emph{vertical boundary}.\index{vertical boundary}\index{$I$-bundle!vertical boundary}

  A white bigon\index{bigon} region of the polyhedral decomposition of $M$ is bounded by two edges and two vertices; recall that edges are crossing arcs,\index{crossing arc} and lie in the intersection $W\cap S$, and the vertices are ideal, forming portions of the parabolic locus.\index{parabolic locus} Thus a regular neighborhood of a white bigon face can be visualized as a thickened square, with two sides of its boundary on $S$ and two sides of its boundary on the parabolic locus. Note such a thickened square is a portion of an $I$-bundle, with horizontal boundary a neighborhood in $S$ of the two crossing arcs of $W\cap S$ that form the boundary of the bigon, and vertical boundary on the parabolic locus. We can complete this thickened square to an $I$-bundle over a subsurface of $S$ with boundary by attaching a neighborhood of the annulus that forms the parabolic locus. This annulus has boundary components on $S$, and is parallel to a link component. Its neighborhood can be given the structure of an $I$-bundle with fibers $I$ parallel to those of the bigon.\index{bigon} Thus the union of the neighborhood of the bigon and the neighborhood of this annulus (or possibly two annuli in the case of a link) forms an $I$-bundle\index{$I$-bundle} in $M\cut S$. 
\end{proof}

\begin{corollary}\label{Cor:RemoveCrossingsTwist}
  Let $K$ be a twist-reduced diagram of a hyperbolic alternating link.\index{alternating knot or link!checkerboard surface}\index{checkerboard surface} Let $K'$ be the diagram obtained from $K$ by removing all crossings but one in each twist region of $K$. Let $S$ denote a checkerboard surface of $K$, and let $S'$ denote the corresponding checkerboard surface of $K'$. Then\index{guts}
  \[ \guts((S^3-K)\cut S) = \guts((S^3-K')\cut S'). \qedhere\]
\end{corollary}

\subsection{Essential annuli}

Our method of proving the volume bound on alternating links, \refthm{VolAlt}, is to apply the volume bound via guts, \refthm{astguts}, to the modified diagram $K'$ of $K$, as in \refcor{RemoveCrossingsTwist}. We will determine the Euler characteristic\index{Euler characteristic} of the guts of checkerboard surfaces of $K'$. Recall that to identify guts,\index{guts} we must first cut along essential annuli. Thus the next step in the proof is to find essential\index{essential} annuli in the cut manifold.

Suppose there is an essential annulus. Then the proof most easily breaks into two cases, depending on whether the annulus is \emph{parabolically compressible} or not, in the sense of the following definition.

\begin{figure}
%% Accompanies image file 'F13-04-ParCom.eps' (pdf, eps, ps)
\begingroup%
  \makeatletter%
  \providecommand\color[2][]{%
    \errmessage{(Inkscape) Color is used for the text in Inkscape, but the package 'color.sty' is not loaded}%
    \renewcommand\color[2][]{}%
  }%
  \providecommand\transparent[1]{%
    \errmessage{(Inkscape) Transparency is used (non-zero) for the text in Inkscape, but the package 'transparent.sty' is not loaded}%
    \renewcommand\transparent[1]{}%
  }%
  \providecommand\rotatebox[2]{#2}%
  \newcommand*\fsize{\dimexpr\f@size pt\relax}%
  \newcommand*\lineheight[1]{\fontsize{\fsize}{#1\fsize}\selectfont}%
  \ifx\svgwidth\undefined%
    \setlength{\unitlength}{222.16092682bp}%
    \ifx\svgscale\undefined%
      \relax%
    \else%
      \setlength{\unitlength}{\unitlength * \real{\svgscale}}%
    \fi%
  \else%
    \setlength{\unitlength}{\svgwidth}%
  \fi%
  \global\let\svgwidth\undefined%
  \global\let\svgscale\undefined%
  \makeatother%
  \begin{picture}(1,0.45405616)%
    \lineheight{1}%
    \setlength\tabcolsep{0pt}%
    \put(0,0){\includegraphics[width=\unitlength]{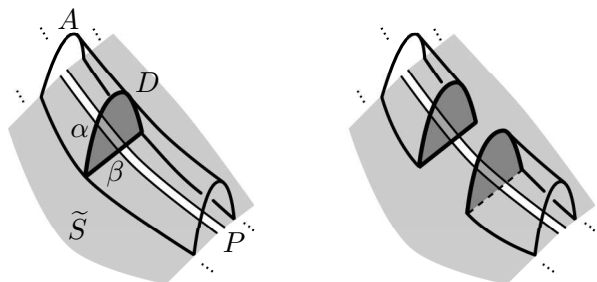}}%
    \put(0.36377093,0.05459141){\color[rgb]{0,0,0}\makebox(0,0)[lt]{\lineheight{0}\smash{\begin{tabular}[t]{l}$P$\end{tabular}}}}%
    \put(0.10089314,0.07679299){\color[rgb]{0,0,0}\makebox(0,0)[lt]{\lineheight{0}\smash{\begin{tabular}[t]{l}$\widetilde{S}$\end{tabular}}}}%
    \put(0.08361476,0.43327307){\color[rgb]{0,0,0}\makebox(0,0)[lt]{\lineheight{0}\smash{\begin{tabular}[t]{l}$A$\end{tabular}}}}%
    \put(0.10544008,0.25048597){\color[rgb]{0,0,0}\makebox(0,0)[lt]{\lineheight{0}\smash{\begin{tabular}[t]{l}$\alpha$\end{tabular}}}}%
    \put(0.16696641,0.1751627){\color[rgb]{0,0,0}\makebox(0,0)[lt]{\lineheight{0}\smash{\begin{tabular}[t]{l}$\beta$\end{tabular}}}}%
    \put(0.21511045,0.32597627){\color[rgb]{0,0,0}\makebox(0,0)[lt]{\lineheight{0}\smash{\begin{tabular}[t]{l}$D$\end{tabular}}}}%
  \end{picture}%
\endgroup%

  \caption{A portion of a parabolically compressible\index{parabolically compressible} annulus on the left, and a parabolic compression\index{parabolic compression} on the right.}
  \label{Fig:ParabolicCompression}
\end{figure}

\begin{definition}\label{Def:ParabolicCompress}
Let $M$ be a hyperbolic 3-manifold with a properly embedded essential\index{essential} surface $S$. Let $P$ denote the parabolic locus of $M\cut S$, as in \refdef{Cut}.\index{parabolic locus} An annulus $A$, properly embedded in $M\cut S$ with $\bdy A\subset \widetilde{S}$, is \emph{parabolically compressible}\index{parabolically compressible}\index{compressible!parabolically compressible} if there exists a disk $D$ with interior disjoint from $A$, with $\bdy D$ meeting $A$ in an essential arc $\alpha$ on $A$, and with $\beta = \bdy D - \alpha$ lying on $\widetilde{S} \cup P$, with $\beta$ meeting $P$ transversely exactly once.
We may surger along such a disk; this is called a \emph{parabolic compression},\index{parabolic compression} and it turns the annulus $A$ into a disk meeting $P$ transversely exactly twice, with boundary otherwise on $\widetilde{S}$. See \reffig{ParabolicCompression}.
\end{definition}

\begin{definition}\label{Def:EPD}
Let $D$ be a disk properly embedded in $M\cut S$ with boundary consisting of two arcs on $\widetilde{S}$ and two arcs on the parabolic locus\index{parabolic locus} $P$. We say $D$ is a \emph{essential product disk (EPD)}.\index{essential product disk (EPD)}\index{EPD}
\end{definition}

A proof nearly identical to that of \reflem{BigonCharacteristic} shows that EPDs belong to the $I$-bundle of $M\cut S$ (\refex{EPD}).

\begin{lemma}\label{Lem:NoParabComprAnnulus}
Let $S$ be the shaded checkerboard surface of a link with a twist-reduced alternating diagram\index{alternating diagram!checkerboard surface}\index{checkerboard surface} $K$, whose hyperbolic complement we denote by $S^3-K=M$.
Suppose there are no white bigon\index{bigon} regions, and suppose $A$ is an essential annulus properly embedded in $M\cut S$, disjoint from the parabolic locus\index{parabolic locus} and not parallel to the parabolic locus, with $\bdy A\subset\widetilde{S}$. Then $A$ is not parabolically compressible.\index{parabolically compressible} 
\end{lemma}

The proof of \reflem{NoParabComprAnnulus} is completed by considering how a parabolically compressible annulus intersects the polyhedra in the decomposition of an alternating link, similar to several proofs in \refchap{Alternating}. The following lemma will be useful.

\begin{lemma}\label{Lem:MarcLemma7}
Let $K$ be a link with a prime, twist-reduced alternating diagram,\index{alternating diagram} with corresponding ideal polyhedral decomposition. Let $D_1$ and $D_2$ be normal disks\index{normal} in the polyhedra such that $\bdy D_1$ and $\bdy D_2$ meet exactly four interior edges. Isotope $\bdy D_1$ and $\bdy D_2$ to minimize intersections $\bdy D_1\cap\bdy D_2$ in faces. If $\bdy D_1$ intersects $\bdy D_2$, then $\bdy D_1$ intersects $\bdy D_2$ exactly twice, in two faces of the same color. 
\end{lemma}

\begin{proof}
The boundaries $\bdy D_1$ and $\bdy D_2$ are quadrilaterals, with sides of $\bdy D_i$ between intersections with interior edges. Note that $\bdy D_1$ can intersect $\bdy D_2$ at most once in any of its sides by the requirement that the number of intersections be minimal (else isotope through a face). 
Thus there are at most four intersections of $\bdy D_1$ and $\bdy D_2$. If $\bdy D_1$ meets $\bdy D_2$ four times, then the two quads run through the same faces, both bounding disks, and can be isotoped off each other using the fact that the diagram is prime. Since the quads intersect an even number of times, there are either zero or two intersections. If zero intersections, we are done.

So suppose there are two intersections. Suppose $\bdy D_1$ intersects $\bdy D_2$ exactly twice in faces of the opposite color. Then an arc $\alpha_1\subset \bdy D_1$ has both endpoints on $\bdy D_1\cap \bdy D_2$ and meets only one intersection of $\bdy D_1$ with an interior edge of the polyhedron. Similarly, an arc $\alpha_2\subset\bdy D_2$ has both endpoints on $\bdy D_1\cap \bdy D_2$ and meets only one intersection of $\bdy D_2$ with an interior edge of the polyhedral decomposition. Then $\alpha_1\cup \alpha_2$ is a closed curve on the boundary of the polyhedron meeting exactly two interior edges. This gives a curve in the diagram of $K$ meeting $K$ exactly twice. Because the diagram is prime, there must be no crossings on one side of the curve. In the polyhedron, this means the arcs $\alpha_1$ and $\alpha_2$ are parallel, and we can isotope them to remove the intersections of $D_1$ and $D_2$. 
\end{proof}

\begin{proof}[Proof of \reflem{NoParabComprAnnulus}]
Suppose by way of contradiction that $A$ is an essential, parabolically compressible annulus properly embedded in $M\cut S$ with $\bdy A\subset\widetilde{S}$. Perform a parabolic compression\index{parabolic compression} to obtain an EPD\index{essential product disk (EPD)}\index{EPD} $E$. Put $E$ into normal form\index{normal form}\index{normal} with respect to the polyhedral decomposition of $M\cut S$.

Suppose $E$ intersects a white face $V$ of $W$. Consider the arcs $E\cap V$; such an arc has both endpoints on $\widetilde{S}$. If one cuts off a disk on $E$ that does not meet the parabolic locus,\index{parabolic locus} then there will be an innermost such disk. Its boundary consists of an arc in $W$ and an arc in $S$. We may sketch the boundary of the disk on the diagram of $K$, since the graph on the polyhedra is identical to the projection graph of the diagram (see \refthm{PolyAltKnot}). Thus this innermost disk has boundary intersecting the link diagram exactly twice. Because the diagram of $K$ is prime, this disk bounds a region containing no crossings. But then the original innermost disk in the polyhedron it is not normal:\index{normal} its boundary runs from a single edge back to that edge. This is a contradiction. 

So $E\cap W$ consists of arcs running from $\widetilde{S}$ to $\widetilde{S}$, cutting off disks on either side meeting the parabolic locus.\index{parabolic locus} Thus the white surface cuts $E$ into normal quadrilaterals\index{normal} $\{E_1, \dots, E_n\}$, with $n\geq 2$ by assumption that $E$ intersects a white face. On the end of $E$, the quadrilateral $E_1$ has one side on $W$, two sides on $\widetilde{S}$ and the final side on the parabolic locus (a boundary face). Isotope slightly off the boundary face into the adjacent white face so that $E_1$ remains normal. Do the same for $E_n$. Then all quadrilaterals $E_1, \dots, E_n$ have two sides on $S$ and two sides on $W$. 

Superimpose $E_1$ and $E_2$ onto the boundary of one of the (identical) polyhedra. An edge of $E_1$ in a white face $V$ is glued to an edge of $E_2$ in the same white face, but by a rotation in the face $V$. Thus when we superimpose, $\bdy E_2\cap V$ is obtained from $\bdy E_2\cap V$ by a rotation in $V$.

If $E_1\cap V$ is not parallel to a single boundary edge, then $\bdy E_1\cap V$ must intersect $\bdy E_2\cap V$; see \reffig{NoEPD}. Then \reflem{MarcLemma7} implies that $\bdy E_1$ and $\bdy E_2$ also intersect in another white face. But $\bdy E_1$ is parallel to a single boundary edge in its second white face, so $\bdy E_2$ cannot intersect it. This is a contradiction.

\begin{figure}
%% Accompanies image file 'F13-05-NoEPD.eps' (pdf, eps, ps)
\begingroup%
  \makeatletter%
  \providecommand\color[2][]{%
    \errmessage{(Inkscape) Color is used for the text in Inkscape, but the package 'color.sty' is not loaded}%
    \renewcommand\color[2][]{}%
  }%
  \providecommand\transparent[1]{%
    \errmessage{(Inkscape) Transparency is used (non-zero) for the text in Inkscape, but the package 'transparent.sty' is not loaded}%
    \renewcommand\transparent[1]{}%
  }%
  \providecommand\rotatebox[2]{#2}%
  \newcommand*\fsize{\dimexpr\f@size pt\relax}%
  \newcommand*\lineheight[1]{\fontsize{\fsize}{#1\fsize}\selectfont}%
  \ifx\svgwidth\undefined%
    \setlength{\unitlength}{231.21491432bp}%
    \ifx\svgscale\undefined%
      \relax%
    \else%
      \setlength{\unitlength}{\unitlength * \real{\svgscale}}%
    \fi%
  \else%
    \setlength{\unitlength}{\svgwidth}%
  \fi%
  \global\let\svgwidth\undefined%
  \global\let\svgscale\undefined%
  \makeatother%
  \begin{picture}(1,0.37352353)%
    \lineheight{1}%
    \setlength\tabcolsep{0pt}%
    \put(0,0){\includegraphics[width=\unitlength]{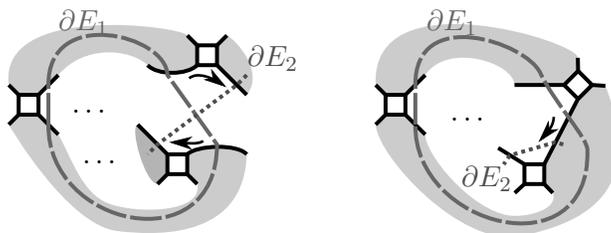}}%
    \put(0.39154324,0.28143778){\color[rgb]{0.30196078,0.30196078,0.30196078}\makebox(0,0)[lt]{\lineheight{1.25}\smash{\begin{tabular}[t]{l}$\bdy E_2$\end{tabular}}}}%
    \put(0.08338838,0.33992376){\color[rgb]{0.4,0.4,0.4}\makebox(0,0)[lt]{\lineheight{1.25}\smash{\begin{tabular}[t]{l}$\bdy E_1$\end{tabular}}}}%
    \put(0.68347943,0.33992376){\color[rgb]{0.4,0.4,0.4}\makebox(0,0)[lt]{\lineheight{1.25}\smash{\begin{tabular}[t]{l}$\bdy E_1$\end{tabular}}}}%
    \put(0.74030787,0.0852984){\color[rgb]{0.30196078,0.30196078,0.30196078}\makebox(0,0)[lt]{\lineheight{1.25}\smash{\begin{tabular}[t]{l}$\bdy E_2$\end{tabular}}}}%
    \put(0.10567882,0.20495067){\color[rgb]{0,0,0}\makebox(0,0)[lt]{\lineheight{1.25}\smash{\begin{tabular}[t]{l}$\dots$\end{tabular}}}}%
    \put(0.12285325,0.12247509){\color[rgb]{0,0,0}\makebox(0,0)[lt]{\lineheight{1.25}\smash{\begin{tabular}[t]{l}$\dots$\end{tabular}}}}%
    \put(0.72248342,0.19210424){\color[rgb]{0,0,0}\makebox(0,0)[lt]{\lineheight{1.25}\smash{\begin{tabular}[t]{l}$\dots$\end{tabular}}}}%
  \end{picture}%
\endgroup%

\caption{Left: $\bdy E_1$ is not parallel to a boundary edge in $U$, hence $\bdy E_1$ meets $\bdy E_2$ in $U$. Right: $\bdy E_1$ is parallel to a boundary edge.}
\label{Fig:NoEPD}
\end{figure}

So $E_1\cap V$ is parallel to a single boundary edge (and hence so is $E_2\cap V$). But then $E_1$ meets both white faces in arcs parallel to boundary edges. Isotoping $\bdy E_1$ slightly into these boundary faces and transfer the curve to the diagram of the link. This gives a closed curve in the link diagram meeting the projection graph of $K$ in exactly two crossings, running to opposite sides of the crossings. If the crossings are distinct, then because the diagram is twist-reduced, the two crossings must bound white bigons\index{bigon} between them, contradicting the fact that there are no white bigon regions in the diagram. So the crossings are not distinct. Returning to the polyhedron, $\bdy E_1$ encircles a single ideal vertex of the polyhedron. Repeating the argument with $E_2$ and $E_3$, and so on, we find that each $\bdy E_i$ encircles a single ideal vertex. Gluing these together, the original annulus $A$ is parallel to the parabolic locus.\index{parabolic locus} This contradicts our assumption on $A$.

So if there is an EPD\index{essential product disk (EPD)}\index{EPD} $E$, it cannot meet $W$. Then it lies completely in a single polyhedron of the decomposition. Its boundary runs through two shaded faces and two boundary faces. Transfer to the link diagram; its boundary defines a curve meeting the link diagram in exactly two crossings, running to opposite sides of the crossings. Because the diagram is twist-reduced, the curve $\bdy E$ encloses a string of white bigons.\index{bigon} But there are no white bigons, so $\bdy E$ must run in and out of the same boundary face. This contradicts the fact that it was normal.\index{normal}
\end{proof}

\begin{lemma}\label{Lem:ParabIncomp}
Let $K$ be a hyperbolic alternating link with a prime, twist-reduced diagram and corresponding polyhedral decomposition. Let $M$ denote $S^3-K$ and let $S$ denote the shaded checkerboard surface.\index{checkerboard surface} Suppose that there are no white bigons\index{bigon} in the polyhedra. Suppose $A$ is an essential\index{essential} annulus embedded in $M\cut S$, disjoint from the parabolic locus\index{parabolic locus} and not parallel to it, with $\bdy A\subset\widetilde{S}$. Then $A$ bounds a Seifert fibered solid torus. 
\end{lemma}

\begin{proof}
By \reflem{NoParabComprAnnulus} we may assume that $A$ is not parabolically compressible.\index{parabolically compressible} Put it into normal form\index{normal form} with respect to the polyhedral decomposition. Because the Euler characteristic\index{Euler characteristic} of an annulus is $0$, each normal disk\index{normal} making up $A$ must have combinatorial area\index{combinatorial area} $0$ by the Gauss--Bonnet lemma, \reflem{GaussBonnet}. Because $A$ does not meet the parabolic locus,\index{parabolic locus} each such disk must meet exactly four interior edges; see \refdef{CombinatorialArea}. Thus the white surface $W$ cuts $A$ into squares $E_1, \dots, E_n$. Note that if a component of intersection of $E_i\cap W$ is parallel to a boundary edge, then the disk of $W$ bounded by $E_i\cap W$, the boundary edge, and portions of edges of $\widetilde{S}\cap W$ defines a parabolic compression\index{parabolic compression} disk for $A$, contradicting the fact that $A$ cannot be parabolically compressible. So no component of $E_i\cap W$ is parallel to a boundary edge.

Again superimpose all squares $E_1, \dots, E_n$ on one of the polyhedra. The squares are glued in white faces, and cut off more than a single boundary edge in each white face, so $\bdy E_i$ must intersect $\bdy E_{i+1}$ in a white face; see again \reffig{NoEPD}. Then \reflem{MarcLemma7} implies $\bdy E_i$ intersects $\bdy E_{i+1}$ in both of the white faces it meets. Similarly, $\bdy E_i$ intersects $\bdy E_{i-1}$ in both its white faces. Because $E_{i-1}$ and $E_{i+1}$ lie in the same polyhedron, they are disjoint (or $E_{i-1} = E_{i+1}$, but this makes $A$ a M\"obius band rather than an annulus; see \refex{Mobius}).
This is possible only if $E_{i-1}$, $E_i$, and $E_{i+1}$ line up as in \reffig{FusedUnits} left, bounding portions of the polyhedron as shown. These transfer to the link diagram to bound tangles; Lackenby calls such tangles \emph{units} in \cite{lackenby:alt-volume}. Then all $E_j$ form a cycle of such tangles, as in \reffig{FusedUnits} right.

\begin{figure}
%% Accompanies image file 'F13-06-FusUnt.eps' (pdf, eps, ps)
\begingroup%
  \makeatletter%
  \providecommand\color[2][]{%
    \errmessage{(Inkscape) Color is used for the text in Inkscape, but the package 'color.sty' is not loaded}%
    \renewcommand\color[2][]{}%
  }%
  \providecommand\transparent[1]{%
    \errmessage{(Inkscape) Transparency is used (non-zero) for the text in Inkscape, but the package 'transparent.sty' is not loaded}%
    \renewcommand\transparent[1]{}%
  }%
  \providecommand\rotatebox[2]{#2}%
  \newcommand*\fsize{\dimexpr\f@size pt\relax}%
  \newcommand*\lineheight[1]{\fontsize{\fsize}{#1\fsize}\selectfont}%
  \ifx\svgwidth\undefined%
    \setlength{\unitlength}{329.52296448bp}%
    \ifx\svgscale\undefined%
      \relax%
    \else%
      \setlength{\unitlength}{\unitlength * \real{\svgscale}}%
    \fi%
  \else%
    \setlength{\unitlength}{\svgwidth}%
  \fi%
  \global\let\svgwidth\undefined%
  \global\let\svgscale\undefined%
  \makeatother%
  \begin{picture}(1,0.34611023)%
    \lineheight{1}%
    \setlength\tabcolsep{0pt}%
    \put(0,0){\includegraphics[width=\unitlength]{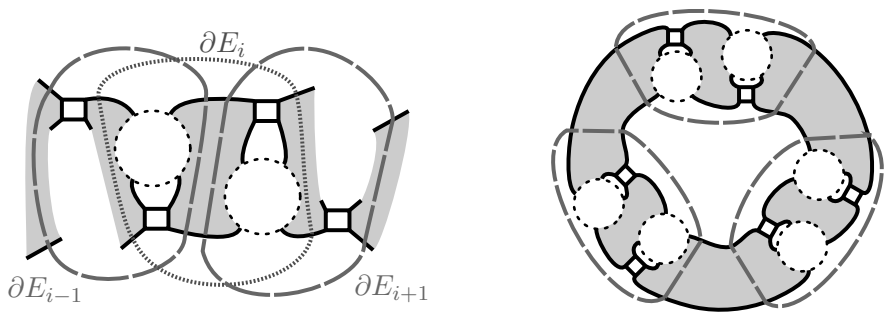}}%
    \put(-0.00236159,0.01945726){\color[rgb]{0.4,0.4,0.4}\makebox(0,0)[lt]{\lineheight{0}\smash{\begin{tabular}[t]{l}$\bdy E_{i-1}$\end{tabular}}}}%
    \put(0.39399183,0.02175615){\color[rgb]{0.4,0.4,0.4}\makebox(0,0)[lt]{\lineheight{0}\smash{\begin{tabular}[t]{l}$\bdy E_{i+1}$\end{tabular}}}}%
    \put(0.21673686,0.29993596){\color[rgb]{0.30196078,0.30196078,0.30196078}\makebox(0,0)[lt]{\lineheight{0}\smash{\begin{tabular}[t]{l}$\bdy E_i$\end{tabular}}}}%
  \end{picture}%
\endgroup%

  \caption{Left: $E_{i-1}$, $E_i$, and $E_{i+1}$ must intersect as shown. Right: cycle of three such tangles.}
  \label{Fig:FusedUnits}
\end{figure}

Observe from \reffig{FusedUnits} that each disk $E_i$ encircles two units, with a band of shaded surface between $E_i$ and $E_{i+2}$ in the same polyhedron. 
Then in each polyhedron, these disks of $A$ bound a solid cylinder (a ball) with top and bottom on white faces --- one the central region of \reffig{FusedUnits}, right, and one the unbounded region --- and sides along disks $E_j$ and shaded faces. The two solid cylinders glue across white faces with a twist, to form a solid torus. As each cylinder can be written as $D^2\times I$, with $D^2\times\{0\}$ and $D^2\times\{1\}$ on white faces, the gluing by a twist in the white face gives the solid torus a Seifert fibering. Thus $A$ bounds a Seifert fibered solid torus. 
\end{proof}

\begin{theorem}[Lackenby, \cite{lackenby:alt-volume}]\label{Thm:ChiGuts}
Let $K$ be a link with a prime, twist-reduced alternating diagram,\index{alternating diagram}\index{alternating knot or link!checkerboard surface}\index{checkerboard surface} and corresponding polyhedral decomposition. Let $M$ denote the complement of $K$, let $S$ and $W$ denote the checkerboard surfaces, and let $r_S$ and $r_W$ denote the number of non-bigon\index{bigon} regions of $S$ and $W$ respectively. Then\index{guts}
\[\chi(\guts(M\cut S))=2-r_W,\quad \chi(\guts(M\cut W))=2-r_S.\]
\end{theorem}

\begin{proof}
  Suppose first that the diagram has no white bigon\index{bigon} regions. Then \reflem{NoParabComprAnnulus} implies there is no embedded essential annulus that is parabolically compressible,\index{parabolically compressible} and \reflem{ParabIncomp} implies any parabolically incompressible annulus bounds a Seifert fibered solid torus. Thus
  \[ \chi(\guts(M\cut S)) = \chi(M\cut S).\]
  Since $M\cut S$ is obtained by gluing two balls along white faces, $\chi(M\cut S) = 2-r_W$.

If the diagram contains white bigon regions, then replace each string of white bigons\index{bigon} in the diagram by a single crossing, obtaining a new link $K'$. Let $M'$ denote $S^3-K'$ and let $S'$ be the checkerboard surface coming from the same shaded regions as $S$ in $K$. By \refcor{RemoveCrossingsTwist}, $\guts(M\cut S) = \guts(M'\cut S')$. Hence $\chi(\guts(M\cut S)) = \chi(\guts(M'\cut S')) = 2-r_W$.

An identical argument applies to $M\cut W$, replacing $S$ with $W$. 
\end{proof}

\Refthm{VolAlt} is now almost an immediate consequence of \refthm{ChiGuts} and \refthm{astguts}.

\begin{proof}[Proof of \refthm{VolAlt}]
Let $\Gamma$ be the $4$-regular diagram graph associated to $K$ by replacing each twist-region with a vertex. 
Let $|v(\Gamma)|$ denote the number of vertices of $\Gamma$, and $|f(\Gamma)|$ the number of regions. Because $\Gamma$ is 4-valent, the number of edges is $2|v(\Gamma)|$, so
\[\chi(S^2) = 2 =-|v(\Gamma)|+|f(\Gamma)|=-\tw(K)+r_S+r_W.\]
Then applying Theorems~\ref{Thm:astguts} and~\ref{Thm:ChiGuts} gives
\begin{align*}
  \vol(S^3-K) \ &\geq \
  -\frac{1}{2}\voct\chi(\guts(M\cut S))-\frac{1}{2}\voct\chi(\guts(M\cut W)) \\
\ &= \ -\frac{1}{2}\voct(2-r_S-r_W) \\
\ &= \ \frac{1}{2}\voct(\tw(K)-2). \qedhere
\end{align*}
\end{proof}

%%%%%%%%%%%%%%%%%%%%%%%%%%%%%%%%%%%%%%%%%%%%%%%%%%%%%%%%%%%%%%%%%
\section{Exercises}
\begin{exercise}\label{Ex:AsympSharpFullyAug}
  Show the upper bound of \refthm{VolUpperTwistRegions} is asymptotically sharp, in two steps.\index{asymptotically sharp volume bound} First, show there is a sequence of fully augmented links\index{fully augmented link} $L_i$ with $t(L_i)$ crossing circles such that $\vol(S^3-L_i)/t(L_i)$ approaches $10\vtet$ as $i$ goes to infinity.\index{fully augmented link!volume bound} (Hint: take white faces to be regular hexagons.) Then show that there is a sequence of links $K_i$ with twist number $t(K_i)$ such that $\vol(S^3-K_i)/t(K_i)$ approaches $10\vtet$ as $i$ goes to infinity.
\end{exercise}

\begin{exercise}
Use \refthm{MaxVolTet} to give an upper bound on the volume of a 2-bridge knot with continued fraction expansion $[0, a_{n-1}, \dots, a_1]$. 
Find an example of a 2-bridge knot such that your upper bound becomes $2\vtet(\tw(K)-1)$. Use this to show that the lower bound of \refthm{Vol2Bridge} is asymptotically sharp:\index{asymptotically sharp volume bound} The ratio of upper and lower bounds goes to $1$ as $\tw(K)\to\infty$. 
\end{exercise}

\begin{exercise}
The volume of a regular ideal octahedron\index{ideal octahedron, regular}\index{regular ideal octahedron}\index{regular ideal octahedron!volume} is denoted by $\voct$. In \refex{FullyAug2Bridge}, it was shown that a fully augmented 2-bridge link decomposes into regular ideal octahedra. Use this to prove that the volume of a 2-bridge link with twist number $\tw(K)$ is at most $2\voct(\tw(K)-1)$. 
\end{exercise}

\begin{exercise}
  Suppose $K$ is a link complement that admits a rotational symmetry about an axis, with order $p$. That is, suppose there is a curve $\gamma$ in $S^3$ such that a rotation of order $p$ about $\gamma$ preserves $K$. Show that if $p\geq 7$,
  \[ \vol(S^3-K) \geq \left( 1 - \frac{4\pi^2}{49} \right)^{3/2}\vol(S^3-(K\cup\gamma)). \]
\end{exercise}

\begin{exercise}\label{Ex:AnnularJSJ}
  Prove \reflem{AnnulusDecomposition}. 
\end{exercise}

\begin{exercise}\label{Ex:EPD}(EPDs lie in the characteristic submanifold)
  Prove that an essential product disk\index{essential product disk (EPD)}\index{EPD} in $M\cut S$ is a subset of the $I$-bundle\index{$I$-bundle} of $M\cut S$, thus cannot be part of the guts.\index{guts}\index{characteristic submanifold}
\end{exercise}

\begin{exercise}\label{Ex:Mobius}
Let $K$ be a hyperbolic alternating knot\index{alternating knot or link} with a prime twist-reduced diagram, and corresponding polyhedral decomposition. Let $S$ denote the shaded checkerboard surface, and suppose that $A$ is an essential\index{essential} surface with boundary on $\widetilde{S}$ such that $A$ is the union of exactly two normal squares\index{normal} $E_1$ and $E_2$, and such that the sides of $\bdy E_1$ and $\bdy E_2$ in white faces are not parallel to boundary edges of the polyhedra. Then prove that $A$ is a M\"obius band. 
\end{exercise}

\begin{exercise}
  We obtain lower bounds on volumes of highly twisted\index{highly twisted} 2-bridge knots with at least seven crossings per twist region from three theorems in this chapter, namely \refthm{Vol2Bridge}, \refthm{VolBoundHighlyTwisted}, and \refthm{VolAlt}. Compare the bounds coming from each theorem. Which gives the best volume estimate? 
\end{exercise}

\begin{exercise}
  How sharp are \refthm{VolBoundHighlyTwisted} and \refthm{VolAlt}? By tracing through the proofs, find conditions that must be satisfied for the lower bound on volume to be sharp.
\end{exercise}

\chapter{Ford Domains and Canonical Polyhedra}\label{Chap:Canonical}
\blfootnote{Jessica S. Purcell, Hyperbolic Knot Theory}

We have noted that there is (currently) no guarantee that every finite volume cusped hyperbolic 3-manifold admits a decomposition into positively oriented ideal tetrahedra. However, we can guarantee that every cusped hyperbolic 3-manifold admits a decomposition into convex ideal polyhedra. This is the canonical decomposition, first studied by \cite{EpsteinPenner}, which we describe in this chapter.
%% Epstein and Penner build the canonical decomposition using a different model of $\HH^3$ than we have been using, namely the \emph{hyperboloid model}\index{hyperboloid model}. While we strongly recommend learning about the hyperboloid model separately, we will not take the time to introduce it here. 
The canonical decomposition is dual to another decomposition, the Ford domain (sometimes called the Ford--Voronoi domain)\index{Ford--Voronoi domain}, which we will describe first. Our exposition is similar to that of \cite{LackenbyPurcell}, also \cite{aswy}, and \cite{Bonahon:LowDimGeom}.

Before we begin, we give a few words motivating the canonical decomposition. In the case of a hyperbolic knot, if two knot complements have the same canonical decomposition, then they must necessarily be isometric, and hence equivalent by \refthm{GordonLuecke}, the Gordon--Luecke theorem. This result follows from \refthm{CanonicalUnique}, below. Thus for hyperbolic knots, the canonical decomposition is a \emph{complete} knot invariant. Unfortunately, it is not easy to compute in general. However, in this chapter we will explain how it is defined, and give a few examples.

Both Ford domains and canonical polyhedra arise from natural geometric ideas. However, they are somewhat difficult to describe in words because of various choices that must be made. If a hyperbolic 3-manifold has more than one cusp, they depend on a choice of horoball neighborhood of the cusp. For this reason, we begin this chapter by describing choices of horoball neighborhoods and the sets equidistant from horoballs, in \refsec{Horoballs}.

Moreover, the definition of a Ford domain differs in the literature, although all definitions are closely related. Perhaps the simplest definition is the one given by Gu{\'e}ritaud and Schleimer~\cite{GueritaudSchleimer}: they define a Ford domain (or Ford--Voronoi domain)\index{Ford domain}\index{Ford--Voronoi domain} to be the set $S$ of points in $M$ that have a unique shortest path to the fixed horoball neighborhood. The drawback to this definition is that the resulting set $S$ is not a fundamental domain\index{fundamental domain} for the manifold in the sense that not every point of $M$ has a preimage in $S$. Additionally, the components of $S$ are not simply connected, because each component admits a  deformation retraction to a horoball neighborhood, which is homeomorphic to the thickened torus. We will give a closely related definition of the Ford domain in \refsec{FordDomain} that overcomes these difficulties, but at the cost of being slightly more complicated and dependent upon an additional choice (of fundamental domain for the horoball neighborhood). Still, throughout the discussion it is useful to keep Gu{\'e}ritaud and Schleimer's definition in mind. 

However, all the different definitions of Ford domain in the literature still have the same geometric dual: the canonical polyhedral decomposition.\index{canonical decomposition} This convex cell decomposition does depend on choice of horoball neighborhood, but it is independent of all other choices involved in defining the Ford domain. We describe the canonical polyhedral decomposition in \refsec{Canonical}.

\section{Horoballs and isometric spheres}\label{Sec:Horoballs}

Throughout, our setup is the following. We let $M$ be an orientable 3-manifold admitting a complete hyperbolic structure with at least one cusp. The universal cover of $M$ is then $\HH^3$. We may apply an isometry so that the point at infinity $\infty \in \bdy_\infty \HH^3$ maps to a cusp of $M$ under the covering map. Then $M \cong \HH^3/\Gamma$, where $\Gamma \leq \PSL(2,\CC)$ is a discrete group of isometries isomorphic to $\pi_1(M)$ via the holonomy\index{holonomy} representation $\rho\from \pi_1(M)\to \PSL(2,\CC)$.

Because the point at infinity projects to a cusp of $M$, there will be a parabolic\index{parabolic} subgroup $\Gamma_\infty$ of $\Gamma$ fixing the point at infinity. If the cusp of $M$ is a rank-1 cusp, then $\Gamma_\infty$ will be isomorphic to $\ZZ$. If it is a rank-2 cusp, $\Gamma_\infty$ will be isomorphic to $\ZZ\times\ZZ$; see \refdef{RankOneRankTwoCusp}. Only a rank-2 cusp can occur in a finite volume hyperbolic manifold such as a knot complement, and so this is the case we will consider in this chapter. However, much of the discussion here generalizes to the infinite volume case. 

\begin{proposition}\label{Prop:EmbeddedHoroballNbhd}
  A complete hyperbolic 3-manifold contains an embedded horoball neighborhood. That is, there is an embedded neighborhood $N$ of the cusps of $M$ such that $N$ lifts to a disjoint collection of embedded horoballs in $\HH^3$. 
\end{proposition}

\begin{proof}
This result follows immediately from the structure of the thin part, \refthm{ThinPart}. By that theorem, for any $0<\epsilon\leq\epsilon_3$, where $\epsilon_3$ is a universal constant, the $\epsilon$-thin part of $M$ consists of tubes around short geodesics and rank-1 and rank-2 cusps. Ignore the tubes; the cusps are embedded. Their lift to $\HH^3$ consists of disjoint embedded horoballs as required. 
\end{proof}

\begin{lemma}\label{Lem:CountableHoroballs}
Suppose $N$ is an embedded horoball neighborhood of a cusp of $M$ that lifts to the horoball about $\infty\in\bdy_\infty\HH^3$. Then all the lifts of $N$ to $\HH^3$ give countably many horoballs in $\HH^3$, with centers at the points
\[ \{ g(\infty) \mid g \in \Gamma\}. \]
\end{lemma}

\begin{proof}
Let $H_\infty$ denote the horoball about $\infty$ that projects to $N$. Note that for all $g\in\Gamma$, the horoball $g(H_\infty)$ must also project to $N$, and its center is $g(\infty)$. On the other hand, if $H$ is a horoball in $\HH^3$ that projects to $N$, then there must exist $h\in\Gamma$ such that $h(H) = H_\infty$, and so $H$ has center $h^{-1}(\infty)$. Thus the set $\{g(H_\infty) \mid g\in \Gamma\}$ is exactly the set of horoballs projecting to $N$. Because $\Gamma$ is a discrete group, it has countably many elements (\refex{CountableDiscrete}). Thus all lifts of $N$ to $\HH^3$ is a countable set of horoballs. 
\end{proof}

\begin{corollary}\label{Cor:CountableCuspNbhd}
Let $M$ be finite volume. Any embedded horoball neighborhood about all cusps of $M$ lifts to countably many disjoint horoballs in $\HH^3$.
\end{corollary}

\begin{proof}
\Reflem{CountableHoroballs} shows that an embedded horoball neighborhood of one cusp of $M$ lifts to countably many disjoint horoballs.  \Refprop{EmbeddedHoroballNbhd} implies that all horoball lifts from all cusps are embedded. Since $M$ has finite volume, it has only finitely many cusps. Hence the collection of all lifts of horoball neighborhoods is countable. 
\end{proof}

\begin{figure}
  \includegraphics{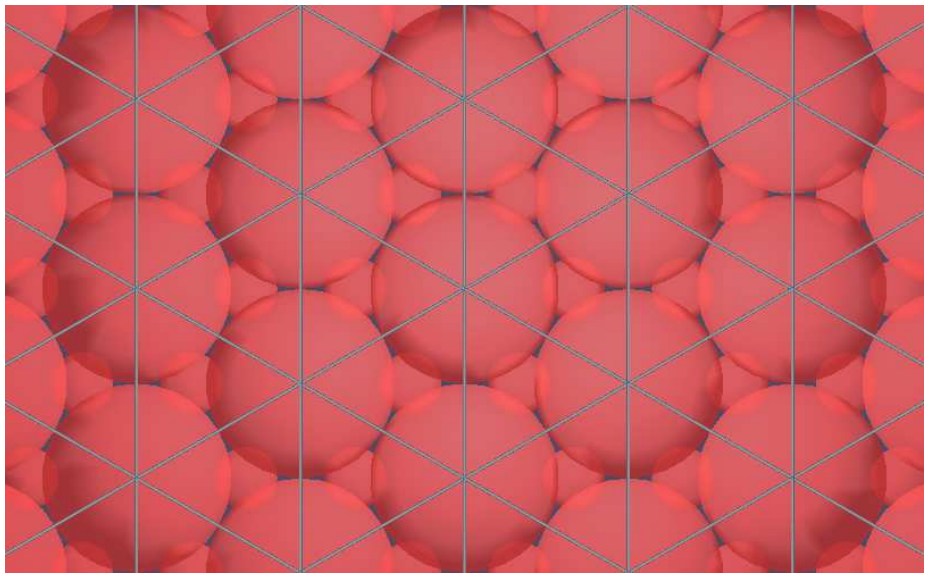}
  \caption{The horoballs in the lift of an embedded cusp of the figure-8 knot complement, from SnapPy \cite{SnapPy}.}
  \label{Fig:Fig8KnotSnappyCusp}
\end{figure}

The software SnapPy \cite{SnapPy} has a feature `Cusp Neighborhoods' that shows the horoballs making up the lift of an embedded horoball neighborhood of $M$ --- or at least those that have Euclidean diameter larger than some specified lower bound. For example, the pattern arising from the figure-8 knot is shown in \reffig{Fig8KnotSnappyCusp}.

Now consider an embedded cusp neighborhood of a manifold with a single cusp, for example the figure-8 knot complement. \Refprop{EmbeddedHoroballNbhd} and \ref{Lem:CountableHoroballs} imply that the cusp lifts to a countable collection of embedded horoballs in $\HH^3$. If we adjust the size of the initial cusp neighborhood of $M$, the sizes of the horoball lifts will also be adjusted. For example, if we shrink the cusp of $M$, the horoball about infinity $H_\infty$ will also shrink, which we see as an increase in Euclidean height of the horoball. All its translates $\gamma(H_\infty)$ will also shrink, which we see as a shrinking of the Euclidean diameters of the horoballs with centers away from $\infty$.

On the other hand, if we increase the size of the cusp of $M$, the horoballs will grow. We can increase their sizes, keeping the cusp embedded, up until the point where the cusp neighborhood becomes tangent to itself. In the lift of the horoball neighborhood, at this point two horoballs are tangent.

If $M$ has multiple cusps, then we can grow and shrink the sizes of cusps independently. However, we can still increase sizes of embedded cusps only until each is tangent either to itself or to another cusp.

\begin{definition}\label{Def:MaximalCusp}
  A \emph{maximal cusp neighborhood}\index{maximal cusp neighborhood}\index{cusp!maximal cusp neighborhood} is an (open) embedded cusp neighborhood for $M$ that is maximal in the sense that no cusp can be expanded while keeping the set of cusps embedded and disjoint. 
\end{definition}

\begin{definition}\label{Def:FullSized}
  Consider the lift of an embedded maximal cusp neighborhood to $\HH^3$, with one cusp lifting to a horoball at infinity. A \emph{full-sized horoball}\index{full-sized horoball}\index{horoball!full-sized} is a horoball in this pattern that is tangent to the horoball at infinity. Viewed from infinity, it has maximal Euclidean diameter. 
\end{definition}

\begin{example}\label{Example:Fig8FullSized}
  The complement of the figure-8 knot admits four full-sized horoballs, distinct up to translation by $\Gamma_\infty$; see again \reffig{Fig8KnotSnappyCusp}.
  For this manifold, each full-sized horoball is tangent to the other three, in a pattern that is known to be the densest possible horoball packing. (This follows from a theorem of B{\"o}r{\"o}czky \cite{boroczky}: the 3-dimensional analogue of \refthm{Boroczky}.)\index{B{\"o}r{\"o}czky cusp density theorem!3-dimensional}\index{cusp density theorem!3-dimensional}
\end{example}

%%%%%%%%%%%%%%%%%%%%%%%%%%%%%%%%%%%%%%%%%%%%%%%%%%%%%%%%%%%%%%%%%

Recall that a \emph{fundamental domain}\index{fundamental domain} for the action of a group on a space is a subset of the space that contains a point from each orbit, whose interior contains exactly one point from each orbit. In this chapter, we will restrict to complete $(G,X)$-structures\index{$(G,X)$-structure} on a manifold $M$, where $X$ is the metric space $\EE^2$ or $\HH^3$ and $G$ acts by isometries. In this case, we require a fundamental domain $R$ to be cut out by geodesic planes. Distinct points in the interior of $R$ project to distinct points in the manifold. The boundary of $R$ is made up of faces intersecting in edges and vertices, and the interior of each face is paired by an isometry of $G$ to exactly one other face; this is called a \emph{face-pairing isometry}\index{face-pairing isometry}. Finally, the quotient of the fundamental domain under the action of face-pairing isometries, which agrees with the restriction of the covering map $X\to M$, is all of $M$.

For example, a Euclidean structure\index{Euclidean structure} on a torus has fundamental domain a single (closed) parallelogram. It follows that the boundary of a cusp of a hyperbolic 3-manifold has a fundamental domain that is a parallelogram on a horosphere in $\HH^3$.

\begin{lemma}\label{Lem:FullSized}
Let $M \cong \HH^3/\Gamma$ be a hyperbolic 3-manifold with at least one cusp. In a horoball pattern in $\HH^3$ given by lifting an embedded maximal cusp neighborhood\index{maximal cusp neighborhood}\index{cusp!maximal cusp neighborhood} for $M$, apply any isometry taking a desired horoball to the one at infinity. Then in the new pattern obtained by applying this isometry, there is at least one full-sized horoball meeting a fundamental domain for the boundary of the horoball about infinity.

Moreover, if $M$ has only one cusp, then there are at least two full-sized horoballs in a fundamental domain. The second is often called the \emph{Adams horoball}.\index{Adams horoball}\index{horoball!Adams}
\end{lemma}

\begin{proof}
By definition, an embedded maximal cusp neighborhood\index{maximal cusp neighborhood}\index{cusp!maximal cusp neighborhood} cannot be expanded or the cusp will no longer be embedded. Thus its lift to $\HH^3$ will no longer consist of disjoint horoballs. That means that for each cusp, one horoball in its fundamental domain must be tangent to another. Hence when we apply an isometry taking a horoball projecting to that cusp to a horoball about infinity, another horoball becomes tangent to the one at infinity, hence full-sized.

In the case that $M$ has exactly one cusp, let $H_\infty$ denote the horoball at infinity and let $H_f$ denote the full-sized horoball. Because there is only one cusp of $M$, the two horoballs must project to the same cusp of $M$. Thus there must be a covering transformation, i.e.\ an isometry $g \in \Gamma$, taking $H_f$ to $H_\infty$. Consider the image $g(H_\infty)$. This is a horoball tangent to $g(H_f)=H_\infty$, hence it must be a full-sized horoball. Apply an isometry of $w\in \Gamma_\infty\leq \Gamma$ fixing $\infty$, if necessary, so that $wg(H_\infty)$ lies in the same fundamental domain of the cusp as $H_f$, and replace $g$ with $wg$. Now either $H_f$ and $g(H_\infty)$ are disjoint full-sized horoballs, as desired, or possibly $H_f=g(H_\infty)$. We now rule out the latter case.

Suppose $g(H_\infty)= H_f$. Consider the effect of $g$ on the geodesic from the center of $H_f$ to $\infty$. If $g$ takes $\infty$ to the center of $H_f$, then this geodesic is mapped to itself, with the point of tangency between $H_f$ and $H_\infty$ mapped to the point of tangency between the two horoballs; hence $g$ has a fixed point in the interior of $\HH^3$, so it is elliptic.\index{elliptic} But $M$ is a manifold, hence \refprop{FreePropDisc} implies the action of $\Gamma$ is fixed point free, so $g$ cannot have a fixed point. Thus $H_f$ and $g(H_\infty)$ are disjoint full-sized horoballs. 
\end{proof}

The Adams horoball is so-called because it appears prominently in work of Adams, e.g.\ \cite{adams:waist}.\index{Adams horoball}\index{horoball!Adams}

\begin{corollary}\label{Cor:CuspVolume}
  The volume of any cusp component in a maximal cusp neighborhood\index{maximal cusp neighborhood}\index{cusp!maximal cusp neighborhood} of $M$ is at least $\sqrt{3}/4$. If $M$ has only one cusp, the volume of a maximal cusp neighborhood is at least $\sqrt{3}/2$. 
\end{corollary}

\begin{proof}
  \Refex{CuspVolumeMaxCusp}.
\end{proof}

We learn a great deal of information from a 3-manifold by considering its cusp neighborhoods, and lifts of maximal cusp neighborhoods\index{maximal cusp neighborhood}\index{cusp!maximal cusp neighborhood} to $\HH^3$. However, because there are countably infinitely many horoballs in such a pattern, it is difficult to compute this pattern and it can be difficult to work with. We can reduce the difficulty of the problem by considering points closer to one lift than another, and their boundaries in $\HH^3$.

\begin{lemma}\label{Lem:DistHoroballs}
  The set of points equidistant from two horoballs in $\HH^3$ forms a geodesic plane in $\HH^3$. 
\end{lemma}

\begin{proof}
  Let $H_1$ and $H_2$ be the two horoballs, and consider the geodesic $\gamma$ between their centers. There exists a unique point $p$ on $\gamma$ equidistant from the boundaries of $H_1$ and $H_2$; see \reffig{Equidistant}.

  \begin{figure}
%% Accompanies image file 'F14-02-Equid.eps' (pdf, eps, ps)
\begingroup%
  \makeatletter%
  \providecommand\color[2][]{%
    \errmessage{(Inkscape) Color is used for the text in Inkscape, but the package 'color.sty' is not loaded}%
    \renewcommand\color[2][]{}%
  }%
  \providecommand\transparent[1]{%
    \errmessage{(Inkscape) Transparency is used (non-zero) for the text in Inkscape, but the package 'transparent.sty' is not loaded}%
    \renewcommand\transparent[1]{}%
  }%
  \providecommand\rotatebox[2]{#2}%
  \newcommand*\fsize{\dimexpr\f@size pt\relax}%
  \newcommand*\lineheight[1]{\fontsize{\fsize}{#1\fsize}\selectfont}%
  \ifx\svgwidth\undefined%
    \setlength{\unitlength}{301.68752289bp}%
    \ifx\svgscale\undefined%
      \relax%
    \else%
      \setlength{\unitlength}{\unitlength * \real{\svgscale}}%
    \fi%
  \else%
    \setlength{\unitlength}{\svgwidth}%
  \fi%
  \global\let\svgwidth\undefined%
  \global\let\svgscale\undefined%
  \makeatother%
  \begin{picture}(1,0.24176512)%
    \lineheight{1}%
    \setlength\tabcolsep{0pt}%
    \put(0,0){\includegraphics[width=\unitlength]{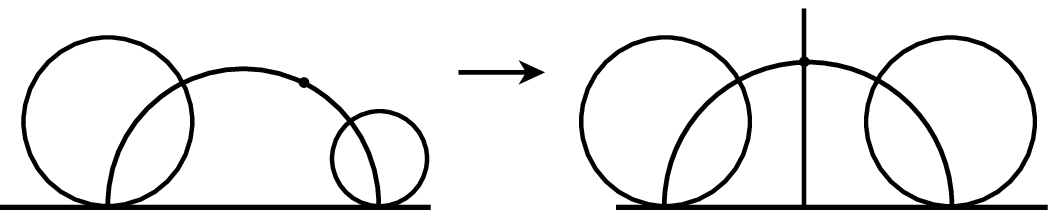}}%
    \put(0.19639522,0.18953566){\color[rgb]{0,0,0}\makebox(0,0)[lt]{\lineheight{1.25000012}\smash{\begin{tabular}[t]{l}$\gamma$\end{tabular}}}}%
    \put(0.30018644,0.17829079){\color[rgb]{0,0,0}\makebox(0,0)[lt]{\lineheight{1.25}\smash{\begin{tabular}[t]{l}$p$\end{tabular}}}}%
    \put(0.6146675,0.01918561){\color[rgb]{0,0,0}\makebox(0,0)[lt]{\lineheight{1.25}\smash{\begin{tabular}[t]{l}$-1$\end{tabular}}}}%
    \put(0.89869474,0.01918561){\color[rgb]{0,0,0}\makebox(0,0)[lt]{\lineheight{1.25}\smash{\begin{tabular}[t]{l}$1$\end{tabular}}}}%
    \put(0.05407084,0.12305794){\color[rgb]{0,0,0}\makebox(0,0)[lt]{\lineheight{1.25}\smash{\begin{tabular}[t]{l}$H_1$\end{tabular}}}}%
    \put(0.35674328,0.10627722){\color[rgb]{0,0,0}\makebox(0,0)[lt]{\lineheight{1.25}\smash{\begin{tabular}[t]{l}$H_2$\end{tabular}}}}%
    \put(0.7725295,0.21441886){\color[rgb]{0,0,0}\makebox(0,0)[lt]{\lineheight{1.25}\smash{\begin{tabular}[t]{l}$P$\end{tabular}}}}%
  \end{picture}%
\endgroup%

  \caption{Map horoballs $H_1$ and $H_2$ to lie over $-1$ and $1$, respectively, with the point equidistant from them on the geodesic between their centers mapped to the point over $0$ with height $1$.}
  \label{Fig:Equidistant}
\end{figure}

  Apply the hyperbolic isometry $\phi$ taking the center of $H_1$ to $-1 \in \CC$, taking the center of $H_2$ to $1\in \CC$, and taking $p$ to the point lying over $0\in\CC$ of height $1$. Under this isometry, $H_1$ and $H_2$ are mapped to horoballs of the same Euclidean diameter, centered at $-1$ and $1$. We claim that the set of points equidistant from two horoballs of the same Euclidean diameter and centers $-1$ and $1$ is the vertical plane $P$ that meets $\CC$ in the imaginary axis. This can be seen as follows. There is a reflection isometry fixing $P$ pointwise and exchanging the two horoballs. Thus the shortest path from a point $q\in P$ to one of the two horoballs will be mapped under the reflection to the shortest path from $q$ to the other horoball. Thus the horoballs are equidistant from $P$.

Now apply $\phi^{-1}$ to this picture. The geodesic plane $P$ is mapped to a geodesic plane in $\HH^3$ that is equidistant to the original horoballs. 
\end{proof} 

In the case that the two horoballs $H_1$ and $H_2$ project to the same cusp, the totally geodesic plane is known as an isometric sphere, as in the following definition. 

\begin{definition}\label{Def:IsoSphere}
  Let $g\in \PSL(2,\CC)$ be an element that does not fix $\infty$. Let $H$ denote a horosphere about $\infty$ in $\HH^3$. Then $g^{-1}(H)$ is a horosphere centered at a point of $\CC \subset (\CC\cup\{\infty\}) = \bdy\HH^3$. Define the set $I(g)$\index{$I(g)$} to be the set of points in $\HH^3$ equidistant from $H$ and $g^{-1}(H)$:
  \[ I(g) = \{x\in\HH^3\mid d(x,H) = d(x,g^{-1}(H))\} \]
The set $I(g)$ is the \emph{isometric sphere}\index{isometric sphere} of $g$. 
\end{definition}

Note that $I(g)$ is well-defined, independent of $H$, even if $H$ and $g^{-1}(H)$ overlap (\refex{DefIsoSphereDoesNotDependOnH}).

\begin{lemma}\label{Lem:IsoSphereImage}
For $g\in\Gamma-\Gamma_\infty$, $g$ maps $I(g)$\index{$I(g)$} isometrically to $I(g^{-1})$, taking the half ball bounded by $I(g)$ to the exterior of the half ball bounded by $I(g^{-1})$.
\end{lemma}

\begin{figure}
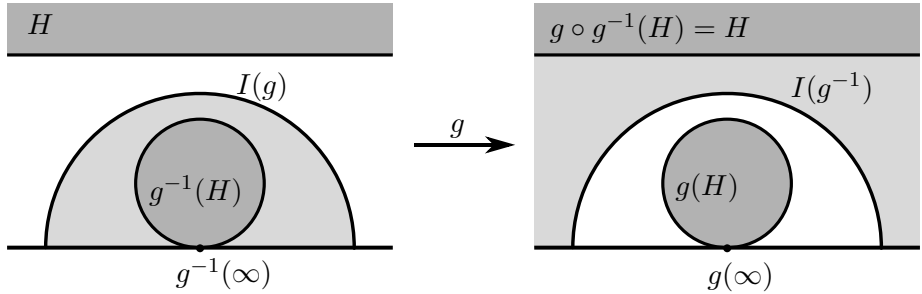
\caption{The horoballs $H$ and $g^{-1}(H)$ are shown, along with the isometric sphere\index{isometric sphere} $I(g)$.\index{$I(g)$} The effect of applying $g$ to this picture is shown on the right: $g\circ g^{-1}(H)$ maps to $H$, $H$ to $g(H)$, and $I(g)$ maps to $I(g^{-1})$. }
\label{Fig:IsoSphere}
\end{figure}

See \reffig{IsoSphere}.

\begin{proof}[Proof of \reflem{IsoSphereImage}]
Let $H$ denote a horosphere about $\infty$ in $\HH^3$. 
Note that $g$ takes the horoball $g^{-1}(H)$ to $H$, and takes $H$ to $g(H)$. Thus $g$ maps $I(g)$ isometrically to the set of points equidistant from $H$ and $g(H)$. This is $I(g^{-1})$. The half-space bounded by $I(g)$,\index{$I(g)$} which contains $g^{-1}(H)$, is mapped to the exterior of the half-space bounded by $I(g^{-1})$, which contains $H$.
\end{proof}

\begin{lemma}\label{Lem:IsoSphereRadius}
As a set, $I(g^{-1})$ (and hence $I(g)$)\index{$I(g)$} is a Euclidean hemisphere orthogonal to $\CC$. If $ g=\mat{a&b\\c&d}\in {\rm PSL}(2,\CC),$ then the center of the Euclidean hemisphere $I({g^{-1}})$ is $g(\infty) = a/c$. Its Euclidean radius is $1/|c|$.
\end{lemma}

\begin{proof}
The fact that $I(g^{-1})$ is a Euclidean hemisphere follows immediately from \reflem{DistHoroballs}: it must be a geodesic plane in $\HH^3$. Moreover, it cannot meet the point $\infty$, since $H$ is centered at that point. Thus it is a Euclidean hemisphere. 

As for the center and radius of the hemisphere, note that $g(\infty)=a/c$, and this must be the center. Consider the geodesic running from $\infty$ to $g(\infty)$.  It consists of points of the form $(a/c, t)$ in $\CC \times \RR^+ \cong \HH^3$.  It will meet the horosphere $H$ about infinity at some height $t=h_1$, and the horosphere $g(H)$ at some height $t=h_0$.  The radius of the isometric sphere\index{isometric sphere} $I({g^{-1}})$ is the height of the point equidistant from points $(a/c, h_0)$ and $(a/c, h_1)$.

Note that $g^{-1}( g(H)) = H$, and hence $h_1$ is given by the height of $g^{-1}(a/c, h_0)$, which can be computed to be $(-d/c, 1/(|c|^2h_0))$.  Thus $h_1 =1/(|c|^2h_0)$.  Then the point equidistant from $(a/c, h_0)$ and $(a/c,1/(|c|^2h_0))$ is the point of height $h = 1/|c|$.
\end{proof}

\begin{lemma}\label{Lem:IsoSphereHeight}
If $p=(x+iy, t)\in\HH^3$ lies in $I(g)$, then $g(p)$ has third coordinate $t$. That is, $g$ preserves the heights of points on $I(g)$.\index{$I(g)$} 
\end{lemma}

\begin{proof}
Let $p\in I(g)$. If $p\in I(g)$ lies on the geodesic from $\infty$ to $g^{-1}(\infty)$, then the third coordinates of $p$ and $g(p)$ can both be determined to be $1/|c|$ from \reflem{IsoSphereRadius}, and so they agree. If $p$ is another point on $I(g)$,\index{$I(g)$} construct a 2/3-ideal triangle\index{$2/3$-ideal triangle} with vertices $\infty$, $g^{-1}(\infty)$, and $p$. Then $g$ maps this 2/3-ideal triangle to a triangle with the same area, hence the same angle at its finite vertex (see \refex{2/3IdealTriangle}). Since $p$ and $g(p)$ also both lie on Euclidean hemispheres of the same radius (again by \reflem{IsoSphereRadius}), it follows that $p$ and $g(p)$ have the same third coordinate. 
\end{proof}

\begin{lemma}\label{Lem:LocallyFinite}
Let $\Gamma\leq \PSL(2,\CC)$ be a nonelementary discrete group with a parabolic\index{parabolic} subgroup $\Gamma_\infty$ fixing the point at infinity. Then the set of all isometric spheres\index{isometric sphere} $\{ I(g) \mid g\in \Gamma-\Gamma_\infty \}$ is \emph{locally finite},\index{locally finite} meaning that for any $x\in\HH^3$, there exists $\epsilon>0$ such that the ball of radius $\epsilon$ centered at $x$ meets only finitely many isometric spheres\index{isometric sphere} $I(g)$\index{$I(g)$} for $g\in\Gamma-\Gamma_\infty$.
\end{lemma}

In fact we show that for all $x\in\HH^3$, and all $\epsilon>0$, the ball $B_\epsilon(x)$ of radius $\epsilon$ centered at $x$ meets only finitely many $I(g)$.\index{$I(g)$}

\begin{proof}[Proof of \reflem{LocallyFinite}]
Suppose there exists $x\in\HH^3$ and $\epsilon>0$ so that $B_\epsilon(x)$ meets infinitely many distinct isometric spheres $I(g_n)$, for elements $g_n\in \Gamma-\Gamma_\infty$. Let $q_n\in B_\epsilon(x)\cap I(g_n)$, and let $H$ be a horoball about infinity. By definition of $I(g_n)$, the point $q_n$ is equidistant from $H$ and $g_n^{-1}(H)$. Then $g_n(q_n)$ is equidistant from $g_n(H)$ and $H$, and by \reflem{IsoSphereHeight} $g_n(q_n)$ has the same third coordinate as $q_n$, hence its third coordinate lies in an interval of length at most $2\epsilon$ centered at the third coordinate of $x$.

Consider next the first and second coordinates of $g_n(q_n)$. The group $\Gamma_\infty$ is isomorphic to $\ZZ\times\ZZ$, generated by two parabolics\index{parabolic} translating along the Euclidean plane $\bdy H$. Choose a (closed) parallelogram on $\bdy H$ that forms a fundamental domain for the action of $\Gamma_\infty$. There exists some $w_n\in\Gamma_\infty$ taking $g_n(q_n)$ to lie in this parallelogram. Then the points $w_n g_n (q_n)$ have first and second coordinates lying within this parallelogram. Since $w_n$ does not affect height, the height of $w_n g_n (q_n)$ agrees with that of $q_n$, and thus also lies in a bounded region. So all points $\{ w_n g_n (q_n) \}$ lie within a bounded parallelopiped in $\HH^3$. Thus they all lie within some bounded distance of our original point $x$, say $d(x, w_n g_n (q_n)) \leq R$ for some $R>0$. 

Now consider the points $\{ (w_n g_n)^{-1} (x)\}$.
We have
\begin{align*}
d(x, (w_n g_n)^{-1} x) & \leq d(x, q_n) + d(q_n, (w_n g_n)^{-1} x) \\
& = d(x, q_n) + d(w_n g_n(q_n), x) \\
& \leq \epsilon + R.
\end{align*}
Let $B$ denote the closed ball of radius $R+\epsilon$ centered at $x$. The above calculation shows that each point $(w_n g_n)^{-1}(x)$ lies within this ball. 

Then for each $n$, $(w_n g_n)^{-1}(B) \cap B$ contains $(w_n g_n)^{-1}(x)$, so is nonempty. Because each of the $g_n$ are distinct, each of the $(w_n g_n)^{-1}$ must be distinct, and therefore we have found an infinite set of elements of $\Gamma$ which take $B$ to a ball intersecting $B$. It follows that $\Gamma$ is not properly discontinuous, as in \refdef{ProperlyDiscont}. But $\Gamma$ is a discrete group, contradicting \reflem{DiscretePropDisc}. 
\end{proof}

%%%%%%%%%%%%%%%%%%%%%%%%%%%%%%%%%%%%%%%%%%%%%%%%%%%%%%%%%%%%%%%%%
\section{Ford domain}\label{Sec:FordDomain}

In this section, we define a special fundamental domain for a hyperbolic 3-manifold, called a Ford domain. We will build this fundamental domain for a hyperbolic 3-manifold with at least one cusp. It will not be unique or canonical, although in the case $M$ has only one cusp, a cover will be unique and canonical. Because of the non-uniqueness of domains for multiple cusps, and the consequent additional difficulties to keep track of in that case, we will first treat the case that $M$ has exactly one cusp.

\subsection{The case of one cusp}

When $M$ has a unique cusp, we define a fundamental domain in terms of isometric spheres\index{isometric sphere} of $M$. 

\begin{definition}\label{Def:EquivariantFordDomain}
Define $B(g)$ to be the open half ball bounded by $I(g)$\index{$I(g)$} in $\HH^3$, and let $\calF(\Gamma)$ be the set
\[ \calF(\Gamma) = \HH^3 - \bigcup_{g\in\Gamma-\Gamma_\infty} B(g) = \bigcap_{g\in\Gamma-\Gamma_\infty} (\HH^3 - B(g)). \]
We call $\calF(\Gamma)$ the \emph{equivariant Ford domain}.\index{equivariant Ford domain}\index{Ford domain!equivariant Ford domain}
Notice that $\calF(\Gamma)$ is invariant under the action of $\Gamma_\infty$. 
\end{definition}

\begin{figure}
\includegraphics{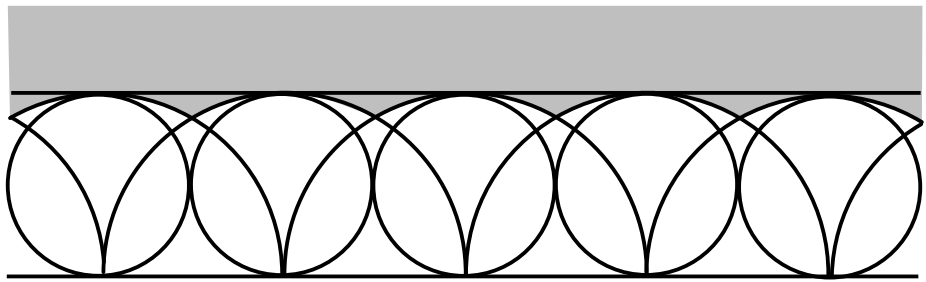}
\caption{The shaded region is 2-dimensional cross section of the equivariant Ford domain corresponding to the maximal cusp neighborhood\index{maximal cusp neighborhood}\index{cusp!maximal cusp neighborhood} of the Figure-8 knot complement.}
\label{Fig:EqFordDomainCrossSection}
\end{figure}

An example is shown in \reffig{EqFordDomainCrossSection}.

\begin{lemma}\label{Lem:EquivFordDomain}
  Fix a maximal cusp neighborhood\index{maximal cusp neighborhood}\index{cusp!maximal cusp neighborhood} of $M$, and let $H$ be the horoball that is lift of the maximal cusp neighborhood to $\HH^3$ with center at infinity. Then
  \[ B(g) = x\in \HH^3 | d(x,g^{-1}(H))< d(x,H), \quad \mbox{and} \]
  \[ \calF(\Gamma) = \{ x\in \HH^3 | d(x,H) \leq d(x,g(H)) \mbox{ for all } g\in \Gamma-\Gamma_\infty \}. \]
\end{lemma}

\begin{proof}
This follows from the definitions. By \refdef{IsoSphere}, $I(g)$\index{$I(g)$} is the set of points equidistant from $H$ and $g^{-1}(H)$. Thus $B(g)$ consists of points strictly closer to $g^{-1}(H)$ than to $H$. Then $\calF(\Gamma)$ consists of points at least as close to $H$ as to any of its translates under $\Gamma-\Gamma_\infty$. 
\end{proof}

\begin{lemma}\label{Lem:EquFordDomainProperties}
The equivariant Ford domain $\calF(\Gamma)$ satisfies the following.
\begin{enumerate}
\item $\calF(\Gamma)$ is a convex subset of $\HH^3$.
\item $\bdy \calF(\Gamma)$ consists of points on $I(g)$\index{$I(g)$} for at least one $g \in\Gamma-\Gamma_\infty$, and admits a decomposition into convex faces, edges, and vertices.
\item  $\calF(\Gamma)$ is invariant under the action of $\Gamma_\infty$. In particular, $\Gamma_\infty$ takes faces, edges, and vertices of $\calF(\Gamma)$ to faces, edges, and vertices, respectively. 
\end{enumerate}
\end{lemma}

\begin{proof}
For convexity, $\calF(\Gamma)$ is the intersection of half-spaces in $\HH^3$, which are convex, thus $\calF(\Gamma)$ is convex. 

Any point on the boundary of $\calF(\Gamma)$ must lie on the boundary of $B(g)$ for some $g\in \Gamma-\Gamma_\infty$. But $\bdy B(g) = I(g)$. Thus the decomposition into faces, edges, and vertices is via isometric spheres\index{isometric sphere} and their intersections: The faces of $\bdy\calF(\Gamma)$ are those points that lie on $I(g)$\index{$I(g)$} for a fixed $g$. The interior of the face consists of points that do not lie on any other $I(h)$ for $h\in\Gamma-\Gamma_\infty$. Edges are points in the intersection of $I(g_1)$ and $I(g_2)$, for some $g_1, g_2\in\Gamma-\Gamma_\infty$. Vertices lie in the intersection of three or more isometric spheres.\index{isometric sphere} Because faces are subsets of Euclidean hemispheres cut out by other Euclidean hemispheres, they are convex. 

Finally, let $w\in \Gamma_\infty$. A point $x\in\HH^3$ lies in $\calF(\Gamma)$ if and only if $x$ lies in the exterior of all open half balls $B(h)$ for $h\in \Gamma-\Gamma_\infty$. This holds if and only if $w(x)$ lies in the exterior of all open half balls $B(hw^{-1})$ for $hw^{-1} \in \Gamma-\Gamma_\infty$. This is the same set of open half balls. Thus $w(x)$ lies in $\calF(\Gamma)$ if and only if $x$ does, and so $\calF(\Gamma)$ is invariant under $\Gamma_\infty$.

Suppose that $x$ lies on a face of $\calF(\Gamma)$, which is a subset of $I(g)$\index{$I(g)$} for some $g\in\Gamma-\Gamma_\infty$. Then $x$ is equidistant from a horoball $H$ at infinity and its translate $g^{-1}(H)$. Thus $w(x)$ is equidistant from $w(H)=H$ and $w(g^{-1}(H)) = (gw^{-1})^{-1}(H)$. It follows that $w(x)$ lies on the isometric sphere\index{isometric sphere} $I(gw^{-1})$, and so $w$ takes isometric spheres to isometric spheres. Because $w$ preserves $\calF(\Gamma)$, $w$ takes the face to a subset of the isometric sphere $wI(g)$, which must be a face. Similarly, if $x$ lies on an edge or vertex, then it lies on the intersection of isometric spheres, and so does $w(x)$. Since $w$ preserves $\calF(\Gamma)$, $w(x)$ lies on an edge or vertex. 
\end{proof}

The faces of $\calF(\Gamma)$, which are contained in isometric spheres\index{isometric sphere} $I(g)$,\index{$I(g)$} can be glued in pairs using the group elements $g$, as in the following lemma.

\begin{lemma}\label{Lem:FaceToFace}
Suppose a subset $f_g$ of $I(g)$\index{$I(g)$} is a face of $\calF(\Gamma)$. Then $g(f_g)$ is a face of $\calF(\Gamma)$.
\end{lemma}

\begin{proof}
Any point $x$ in the interior of $f_g$ is equidistant from $H$ and $g^{-1}(H)$, and because $x$ is in the interior, $x$ lies further away from $h(H)$ for any $h\neq g^{-1}$ in $\Gamma-\Gamma_\infty$. By \reflem{IsoSphereImage}, $g$ maps $x$ to $g(x) \in I(g^{-1})$. Then $g(x)$ is equidistant from $H$ and $g(H)$, but further away from $gh(H)$ for any $h\neq g^{-1}$ in $\Gamma-\Gamma_\infty$. Equivalently, $g(x)$ is equidistant from $H$ and $g(H)$ but further from $k(H)$ for any $k\neq g$ in $\Gamma-\Gamma_\infty$. It follows that $g(x)$ lies in the interior of a face of $\calF(\Gamma)$. 
\end{proof}

\Reflem{FaceToFace} implies that if $I(g)\cap \calF(\Gamma)$ is a face for some $g$, then $g$ is a face-pairing isometry\index{face-pairing isometry} of $\calF(\Gamma)$ in the sense that it maps a face isometrically to a face. At this point, we could take the quotient of $\calF(\Gamma)$ by its face pairing isometries and obtain a manifold that is a covering space of $M$. However, we really want a fundamental domain of $M$, so we restrict $\calF(\Gamma)$ further.

%%%%%%%%%%%%%%%%%%%%%%%%%%%%%%%%%%%%%%%%%%%%%%%%%%%%%%%%%%%%%%%%%

\begin{definition}\label{Def:VerticalFundDomain}
A \emph{vertical fundamental domain}\index{vertical fundamental domain} for $\Gamma_\infty$ is a connected convex fundamental domain for the action of $\Gamma_\infty$ on $\HH^3$ that is cut out by finitely many vertical geodesic planes in $\HH^3$.
\end{definition}

For example, a vertical fundamental domain for the figure-8 knot is cut out by four vertical planes whose boundary on $\CC$ at infinity is the parallelogram bounding the triangles shown in \reffig{Fig8CuspTriang}.

\begin{definition}\label{Def:FordDomain}
A \emph{Ford domain}\index{Ford domain} for $M$ is the intersection of $\calF(\Gamma)$ and a vertical fundamental domain for $\Gamma_\infty$. 
\end{definition}

A Ford domain is not canonical; that is, it is not uniquely defined for the manifold, because the choice of vertical fundamental domain is not unique. However, the equivariant Ford domain is canonical. For this reason, sometimes in the literature the Ford domain is actually defined to be $\calF(\Gamma)$; this is the definition in \cite{Bonahon:LowDimGeom}, for example. However, $\calF(\Gamma)$ is not a finite sided region, and its interior maps to $M$ in an infinite-to-one manner rather than a one-to-one manner, meaning it is not a fundamental domain for $M$. Thus we have chosen to define the Ford domain as in \refdef{FordDomain}. 

\begin{proposition}\label{Prop:FordFiniteFaces}
Let $\overline{M}$ be a compact orientable 3-manifold with a single torus boundary component whose interior admits a complete finite-volume hyperbolic structure. Then any Ford domain $F$ for $M$ is a convex finite-sided polyhedron, cut out by finitely many geodesic planes in $\HH^3$. Moreover, it is a fundamental domain for $M$, in the sense that $M$ is obtained as a quotient of the Ford domain by face-pairing isometries.\index{face-pairing isometry}
\end{proposition}

\begin{proof}
Both the equivariant Ford domain and a vertical fundamental domain are convex, and thus their intersection is convex.

To see that $F$ is a fundamental domain for $M$, we must prove that when we restrict the covering map $\HH^3 \to \HH^3/\Gamma \cong M$ to $F\subset\HH^3$, the projection surjects onto $M$, and that no two points in the interior of $F$ project to the same point of $M$.

First, note that if $x$ is in the interior of $F$, then $x\in\calF(\Gamma)$, so $g(x) \notin \calF(\Gamma)$ for all $g\in\Gamma-\Gamma_\infty$. Since $x$ also lies in the interior of a vertical fundamental domain $V$ for the action of $\Gamma_\infty$, all $g(x) \notin V$ for all nontrivial $g\in\Gamma_\infty$. Thus $g(x)$ lies in $F$ only if $g$ is the identity, therefore no two points in the interior of $F$ project to the same point under the covering projection $\HH^3\to M$.

Now we show that the image of $F$ under the covering map surjects onto $M$. Choose a maximal cusp neighborhood\index{maximal cusp neighborhood}\index{cusp!maximal cusp neighborhood} for $M$, and let $H$ be a horoball about infinity that projects onto the cusp. Let $x\in M$. Let $\delta$ be the minimal distance from $x$ to the maximal cusp neighborhood in $M$. Then there exists a lift $\widetilde{x}$ of $x$ in $\HH^3$ that is distance $\delta$ from $H$. Since this distance is minimal, it follows that $\widetilde{x}$ lies in $\calF(\Gamma)$. There exists $w\in\Gamma_\infty$ such that $w\widetilde{x}$ lies in $V$. Thus $w\widetilde{x}$ lies in $\calF(\Gamma)\cap V= F$, and $w\widetilde{x}$ projects to $x$. So $F$ is a fundamental domain for $M$.

Next we show that the Ford domain is a finite-sided polyhedron. First, remove a small embedded horoball neighborhood from the cusp of $M$ to obtain a compact 3-manifold $\overline{M}$. The lift of the horoball neighborhood lifts to $\HH^3$; remove it from $F$, and call the result $\overline{F}$. 
By \reflem{LocallyFinite}, for any $x$ in $\HH^3$, there is a ball $B_x$ centered at $x$ meeting only finitely many isometric spheres.\index{isometric sphere} The set of all such balls for $x\in \overline{F}$ cover $\overline{F}$. Since $F$ is a fundamental domain for $M$, these balls map to a set of balls covering the compact manifold $\overline{M}$. Thus there is a finite subcollection of balls covering $\overline{M}$, which lift to give a finite cover of $\overline{F}$. Then the total collection of these balls meet only finitely many faces of $F$. Thus $F$ is finite sided.

Finally consider face-pairings.\index{face-pairing isometry} By definition, a face of $V$ is paired to another face of $V$ by an isometry $w\in\Gamma_\infty$. Any $x$ in the interior of the intersection of that face with $\calF(\Gamma)$ is mapped by $w$ to the point $w(x)$ in $V\cap \calF(\Gamma)$. So $w$ is a face-pairing isometry of $F$.

If $x\in F$ lies in the interior of a face $I(g)\cap F$, then $x$ is glued to $g(x)$ in $I(g^{-1})$ in the interior of a face in $\calF(\Gamma)$. The face may be disjoint from $V$, but because $V$ is a fundamental domain for the action of $\Gamma_\infty$, there exists some $w\in\Gamma_\infty$ such that $wg(x)\in V$. Thus $wg(x)\in I(g^{-1}w^{-1})\cap F$. By continuity, the same $w$ maps $g(y)$ to $F$ for any $y$ in a small neighborhood of $x$. Thus $wg$ is a face-pairing isometry.\index{face-pairing isometry}
\end{proof}

\begin{figure}
  %  \import{Figures/Ch14_Canonical/}{Fig8FordDomain.eps_tex}
%% Accompanies image file 'F14-05-FoMark.eps' (pdf, eps, ps)
\begingroup%
  \makeatletter%
  \providecommand\color[2][]{%
    \errmessage{(Inkscape) Color is used for the text in Inkscape, but the package 'color.sty' is not loaded}%
    \renewcommand\color[2][]{}%
  }%
  \providecommand\transparent[1]{%
    \errmessage{(Inkscape) Transparency is used (non-zero) for the text in Inkscape, but the package 'transparent.sty' is not loaded}%
    \renewcommand\transparent[1]{}%
  }%
  \providecommand\rotatebox[2]{#2}%
  \newcommand*\fsize{\dimexpr\f@size pt\relax}%
  \newcommand*\lineheight[1]{\fontsize{\fsize}{#1\fsize}\selectfont}%
  \ifx\svgwidth\undefined%
    \setlength{\unitlength}{302.11444473bp}%
    \ifx\svgscale\undefined%
      \relax%
    \else%
      \setlength{\unitlength}{\unitlength * \real{\svgscale}}%
    \fi%
  \else%
    \setlength{\unitlength}{\svgwidth}%
  \fi%
  \global\let\svgwidth\undefined%
  \global\let\svgscale\undefined%
  \makeatother%
  \begin{picture}(1,0.50574888)%
    \lineheight{1}%
    \setlength\tabcolsep{0pt}%
    \put(0,0){\includegraphics[width=\unitlength]{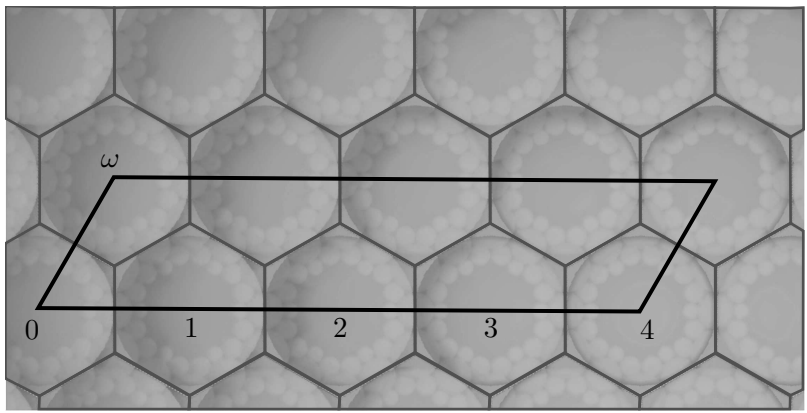}}%
    \put(0.02473593,0.08892551){\color[rgb]{0,0,0}\makebox(0,0)[lt]{\lineheight{1.25}\smash{\begin{tabular}[t]{l}$0$\end{tabular}}}}%
    \put(0.11952742,0.30483889){\color[rgb]{0,0,0}\makebox(0,0)[lt]{\lineheight{1.25}\smash{\begin{tabular}[t]{l}$\omega$\end{tabular}}}}%
    \put(0.22309556,0.09243623){\color[rgb]{0,0,0}\makebox(0,0)[lt]{\lineheight{1.25}\smash{\begin{tabular}[t]{l}$1$\end{tabular}}}}%
    \put(0.40916752,0.09243623){\color[rgb]{0,0,0}\makebox(0,0)[lt]{\lineheight{1.25}\smash{\begin{tabular}[t]{l}$2$\end{tabular}}}}%
    \put(0.59699469,0.09243623){\color[rgb]{0,0,0}\makebox(0,0)[lt]{\lineheight{1.25}\smash{\begin{tabular}[t]{l}$3$\end{tabular}}}}%
    \put(0.79184329,0.08892551){\color[rgb]{0,0,0}\makebox(0,0)[lt]{\lineheight{1.25}\smash{\begin{tabular}[t]{l}$4$\end{tabular}}}}%
  \end{picture}%
\endgroup%

  \caption{A Ford domain for the figure-8 knot complement. The vertical fundamental $V$ domain meets $\CC$ in the parallelogram shown with vertices $0$, $4$, $\omega$ and $4+\omega$. The isometric spheres\index{isometric sphere} intersect to form hexagon faces of the equivariant Ford domain $\calF(\Gamma)$. The ford domain is the intersection $\calF(\Gamma)\cap V$.}
  \label{Fig:Fig8FordDomain}
\end{figure}

\begin{example}\label{Example:Fig8FordDomain}
We can compute explicitly a Ford domain for $M$ the figure-8 knot complement. From \refexamp{Fig8GroupDiscrete} of \refchap{Margulis}, we have a description of three generators of the holonomy group\index{holonomy group} $\Gamma$ of the figure-8 knot complement, namely
\[
T_B = \frac{i}{\sqrt{\omega}}\mat{1&1\\1&-\omega^2}, \quad T_C=\mat{1&\omega\\0&1}, \quad T_D = \mat{2&-1\\1&0},
\]
where recall $\omega = \half + i\frac{\sqrt{3}}{2}$ is a cube root of unity. The transformation $T_C$ fixes the point at infinity, so it plays a part in defining a vertical fundamental domain $V$, but it does not give isometric spheres.\index{isometric sphere} Note isometric spheres corresponding to $T_B^{\pm 1}$ and $T_D^{\pm 1}$ all have radius $1$.
The center of $I(T_B^{-1})$ is $1$, that of $I(T_B)$ is $\omega^2$, that of $I(T_D)$ is $2$, and that of $I(T_D^{-1})$ is $0$.
These are equidistant to the full-sized horoballs of \refexamp{Fig8FullSized}, up to translation in $\Gamma_\infty$. Take a vertical fundamental domain $V$ for $M$ cut out by planes meeting $\CC$ as shown in \reffig{Fig8CuspTriang}. Then $V$ has face-pairing isometries\index{face-pairing isometry} $T_C$ and $T=\mat{1&4\\0&1}$. The ten isometric spheres $I(T_D^{-1})$, $I(T_B^{-1})$, $I(T_D)$, $I(T_C^{-1}TT_B)$, $I(T(T_D^{-1}))$, and their translates under $T_C$ all intersect $V$. In fact, along with $V$ they cut out a Ford domain. See \reffig{Fig8FordDomain}.
\end{example}

\begin{definition}\label{Def:CombinatoriallyEquivalent}
Let $M \cong \HH^3/\Gamma$ and $M'\cong\HH^3/\Gamma'$ be one-cusped hyperbolic 3-manifolds, and let $\calF(\Gamma)$ and $\calF(\Gamma')$ be their respective equivariant Ford domains. Suppose there is a bijection between faces, edges, and vertices of $\calF(\Gamma)$ and faces, edges, and vertices of $\calF(\Gamma')$ such that:
\begin{enumerate}
\item An edge of $\calF(\Gamma)$ is contained in a given face if and only if the corresponding edge of $\calF(\Gamma')$ is contained in the corresponding face, and a vertex of $\calF(\Gamma)$ is contained in a given edge if and only if the corresponding vertex of $\calF(\Gamma')$ is contained in the corresponding edge.
\item A face $f_1$ of $\calF(\Gamma)$ is mapped to a face $f_2$ by a parabolic\index{parabolic} translation in $\Gamma_\infty\leq \Gamma$ if and only if the face of $\calF(\Gamma')$ corresponding to $f_1$ is mapped by a parabolic translation in $\Gamma_\infty'\leq\Gamma'$ to the face corresponding to $f_2$.
\item A face pairing isometry of $\calF(\Gamma)$ matches faces, edges, and vertices if and only if a face pairing isometry of $\calF(\Gamma')$ matches corresponding faces, edges, and vertices. 
\end{enumerate}
Then the equivariant Ford domains $\calF(\Gamma)$ and $\calF(\Gamma')$ are said to be \emph{combinatorially equivalent}.\index{combinatorially equivalent}\index{Ford domain!equivariant!combinatorially equivalent}\index{equivariant Ford domain!combinatorially equivalent}
\end{definition}

\begin{theorem}\label{Thm:CombinatoriallyEquivalent}
  Suppose $M\cong \HH^3/\Gamma$ and $M'\cong \HH^3/\Gamma'$ are one-cusped hyperbolic 3-manifolds with combinatorially equivalent Ford domains $\calF(\Gamma)$ and $\calF(\Gamma')$. Then $M$ and $M'$ are isometric.

  In particular, if $M\cong S^3-K$ and $M'\cong S^3-K'$ are knot complements, then $K$ and $K'$ are isomorphic knots, up to reflection.
\end{theorem}

\begin{proof}
  Because $\calF(\Gamma)$ and $\calF(\Gamma')$ are combinatorially equivalent, the quotients $\calF(\Gamma)/\Gamma$ and $\calF(\Gamma')/\Gamma'$ are homeomorphic as 3-manifolds. Mostow--Prasad rigidity,\index{Mostow--Prasad rigidity} \refthm{MostowGeom}, then implies that the quotients are actually isometric.

  If $M$ and $M'$ are knot complements, then the fact that they come from isomorphic knots follows from Gordon and Luecke's knot complement theorem, \refthm{GordonLuecke}~\cite{gordon-luecke}.\index{Gordon--Luecke knot complement theorem}\index{knot complement theorem}
\end{proof}

\subsection{The case of multiple cusps}

When there are multiple cusps, we still build a Ford domain by considering points closer to one cusp than another. However, there will be a choice involved. We will first give the definitions, then explain how the choices affect the domain.

Let $M\cong \HH^3/\Gamma$ be a complete hyperbolic 3-manifold, and let $C_0, \dots, C_k$ denote its cusps.
Start by choosing a horoball neighborhood of all the cusps of $M$. The neighborhood does not necessarily need to be embedded for the definitions to work. In practice, however, we often consider a choice of maximal cusp neighborhood.\index{maximal cusp neighborhood}\index{cusp!maximal cusp neighborhood}

Lift the horoball neighborhood to the universal cover $\HH^3$. This gives a countable collection of horoballs in $\HH^3$, which will be disjoint if and only if the choice of horoball neighborhood is embedded in $M$. Apply an isometry of $\HH^3$ so that the horoball $H_0$ about infinity projects to the cusp $C_0$ under the covering map.

Let $\Gamma_\infty\leq \Gamma$ denote the subgroup fixing the cusp at infinity. We may choose a vertical fundamental domain $V_0$ for the action of $\Gamma_\infty$, and we still have isometric spheres\index{isometric sphere} $I(g)$\index{$I(g)$} for $g\in \Gamma-\Gamma_\infty$. These will form some of the faces of the Ford domain, but not all. In particular, isometric spheres only give points equidistant from lifts of the cusp $C_0$. We also need to consider points equidistant from the horoball $H_0$ and lifts of cusps $C_1, \dots, C_k$. These are not isometric spheres, so must be defined separately. 

For each $C_j$, $j=1, \dots, k$, choose a horoball $H_j$ such that the distance from $H_0$ to $H_j$ is the distance from the cusp $C_0$ to $C_j$ in $M$. For convenience, we may choose $H_j$ such that its center lies inside the vertical fundamental domain $V_0$. For a horoball $H$ with center on $\CC$, let $P(H)$ denote the set of points equidistant from $H_0$ and $H$. (Thus for $g\in\Gamma-\Gamma_\infty$, the isometric sphere\index{isometric sphere} $I(g)$\index{$I(g)$} is the plane $P(g^{-1}(H_0))$.) Let $B(H)$ denote the open half ball bounded by $P(H)$, so $B(g^{-1}(H_0)) = B(g)$ in the notation of the previous subsection.

Define $\calF_0$ to be the set
\[
\calF_0 = \HH^3 - \left( \bigcup_{g\in \Gamma-\Gamma_\infty} B(g) \cup \bigcup_{j=1}^k\bigcup_{g\in\Gamma} B(g(H_j)) \right).
\]

Define $F_0$ to be the set $\calF_0 \cap V_0$. 

Now repeat the entire construction above, only replacing $C_0$ with $C_j$. Thus we start by applying an isometry so that $H_j$ is a horoball about infinity projecting to $C_j$. Define a vertical fundamental domain $V_j$, and obtain sets $\calF_j$ and $F_j = \calF_j \cap V_j$. Note that under the isometry taking $C_j$ to the cusp at infinity, the sets $\calF_0$ and $F_0$ created before are mapped to some other region of $\HH^3$ whose interior will be disjoint from $\calF_j$ and $F_j$. 

\begin{definition}\label{Def:FordDomainMultiple}
Let $M$ be a complete hyperbolic 3-manifold with cusps $C_0, \dots, C_k$. For each cusp, construct subsets $\calF_j$ and $F_j$ of $\HH^3$ as above.
The (disjoint) union of the sets $\calF_j$ is the \emph{equivariant Ford domain} $\calF$.\index{equivariant Ford domain!multiple cusps}
The (disjoint) union of the sets $F_j$ is the \emph{Ford domain}\index{Ford domain!multiple cusps} of $M$.\index{Ford domain!multiple cusps}
\end{definition}

Observe that each $F_j$ is a convex polyhedron. Faces of $\calF$ are paired by isometries of $\HH^3$ that map the appropriate horoballs $H_i, H_j$ to be equidistant from the given face.

\begin{figure}
\begin{center}
  \includegraphics{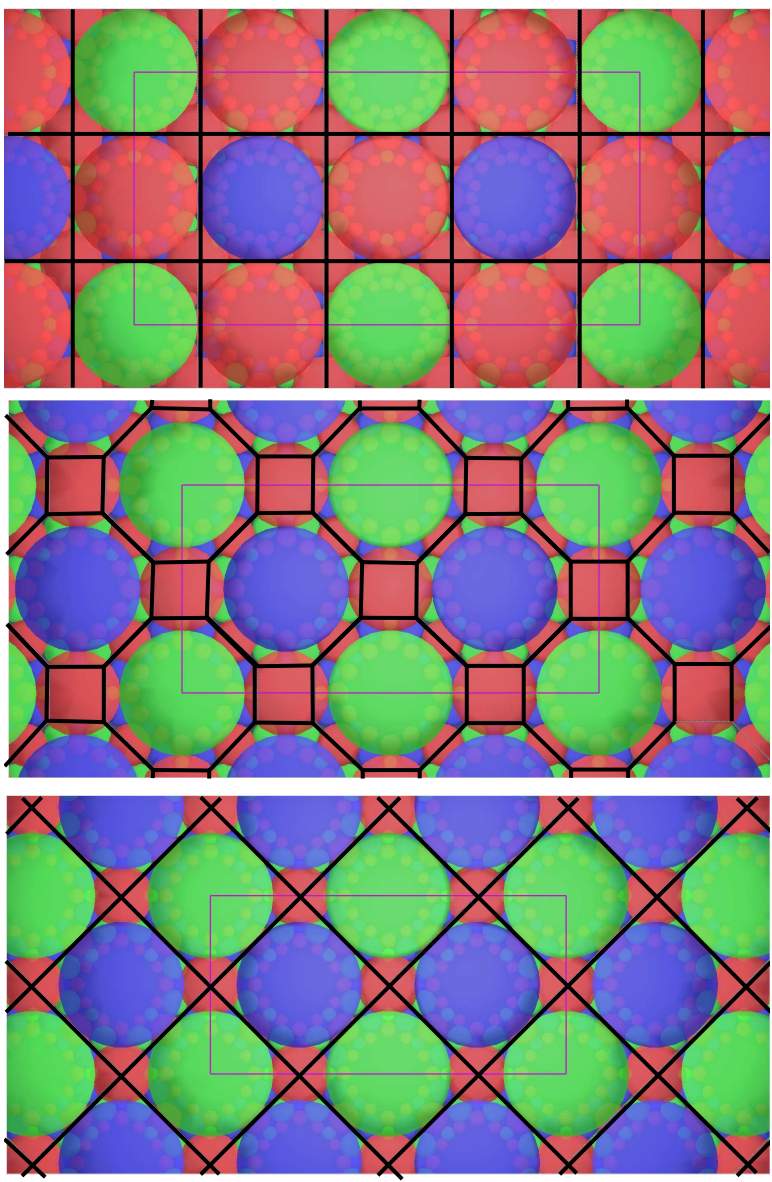}
\end{center}
\caption{Three different Ford domains for the Borromean rings.\index{Borromean rings} The boundary of a vertical fundamental domain $V_0$ is in magenta. Faces of $\calF_0$ have boundaries shown in black.}
\label{Fig:BorroRingsFord}
\end{figure}

Also note that the polyhedra $F_j$ will depend on the choice of expansion of horoballs. An example is shown in \reffig{BorroRingsFord} for the complement of the Borromean rings.\index{Borromean rings} The figures are adapted from SnapPy \cite{SnapPy}. Three different choices of maximal cusp neighborhood\index{maximal cusp neighborhood}\index{cusp!maximal cusp neighborhood} are given, and their effect on the combinatorics of the corresponding subset $\calF_0$ of the equivariant Ford domain is shown. At the top, the cusp $C_0$ has been chosen to be as large as possible. That is, first $C_0$ was expanded until it bumped itself, then $C_1$ and $C_2$ were expanded to meet $C_0$. For this choice of maximal cusp neighborhood,\index{maximal cusp neighborhood}\index{cusp!maximal cusp neighborhood} the faces of $\calF_0$ are squares oriented horizontally. In the middle, $C_0$ still has larger volume than that of $C_1$ and $C_2$, but is not as large as possible. The faces of $\calF_0$ are now octagons and squares. At the bottom, all three cusps have been chosen to have the same volume. The faces of $\calF_0$ are squares again, but oriented on a diagonal.
Akiyoshi has shown there are at most finitely many combinatorially inequivalent Ford domains for any given manifold \cite{Akiyoshi:Finiteness}. 

%%%%%%%%%%%%%%%%%%%%%%%%%%%%%%%%%%%%%%%%%%%%%%%%%%%%%%%%%%%%%%%%%
\section{Canonical polyhedra}\label{Sec:Canonical}

We now describe how to construct the canonical polyhedral decomposition of a finite volume 3-manifold.

Throughout this section, let $M \cong \HH^3/\Gamma$ be a finite volume hyperbolic 3-manifold with a choice of maximal cusp neighborhood,\index{maximal cusp neighborhood}\index{cusp!maximal cusp neighborhood} and corresponding equivariant Ford domain $\calF$. 

We begin by describing the 1-cells, or edges of the polyhedral decomposition. Take a component $\calF_j$ of the equivariant Ford domain, embedded in $\HH^3$ such that it contains a horoball $H$ about infinity. Consider a face $f$ of $\calF_j$ with nonempty interior. Points in the interior of $f$ are exactly those points in $\HH^3$ that are equidistant from $H$ and from another horoball lift $H'$ of a cusp of $M$. For each such face, take the geodesic from the center of $H'$ on $\bdy_\infty\HH^3$ to $\infty$. This geodesic is the \emph{geometric dual} to the face $f$.\index{geometric dual} For each face $f$ of each $\calF_j$, the geometric dual will be identified to a 1-cell of the canonical polyhedral decomposition. Two such geodesics are identified to the same 1-cell if they are identified by an element of $\Gamma$. 

\begin{example}\label{Example:Fig8DualEdges}
For the figure-8 knot complement, with equivariant Ford domain $\calF(\Gamma)$ shown in \reffig{Fig8FordDomain}, the geometric dual edges are those edges running from the points on $\CC$ at the centers of the hexagon faces to the point at infinity, intersecting the faces of the Ford domains in their centers.
\end{example}

\begin{remark}\label{Rmk:DualIntersection}
Note: A geometric dual edge does not necessarily intersect the face of $\calF$ that it is dual to, although the dual edges in the case of the figure-8 knot intersect their corresponding faces in \refexample{Fig8DualEdges}.
For example, a geodesic may be the geometric dual to a face $f_1$ of the Ford domain, but the highest point of the Euclidean hemisphere containing $f_1$ might be covered by another face $f_2$. Then the geometric dual to $f_1$, which runs through this highest point, would not intersect $f_1$. This phenomenon is not common, especially in the small examples that one sees using SnapPy, but it does occur in practice. 
\end{remark}

\begin{figure}
\includegraphics{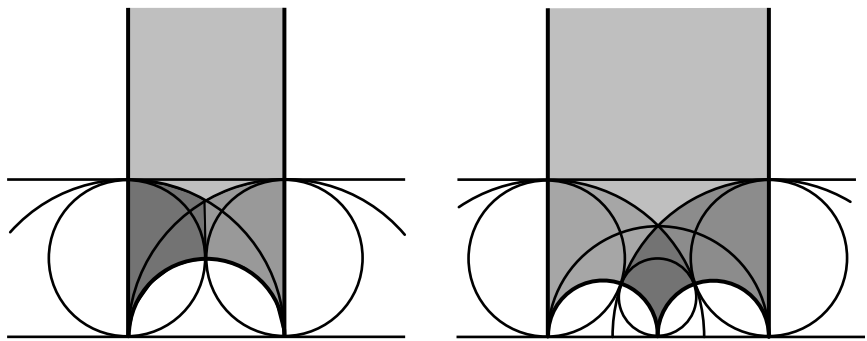}
  \caption{Left: Cross-sectional view of a 2-cell dual to an edge of the equivariant Ford domain for the figure-8 knot complement. Right: In different examples, the dual 2-cell may have more sides.}
  \label{Fig:DualFaces}
\end{figure}

We now construct the 2-cells. 
Consider an edge $e$ of $\calF_j$ with nonempty interior. The interior points on such an edge lie on faces $f'$ and $f''$ of $\calF_j$, where $f'$ consists of points equidistant from $H$ and a horoball $H'$, and where $f''$ consists of points equidistant from $H$ and a horoball $H''$. Denote the geometric dual of the face $f'$ by $\gamma'$, and the geometric dual of the face $f''$ by $\gamma''$. Thus $\gamma'$ and $\gamma''$ are infinite geodesics running from $H'$ to $H$ and from $H''$ to $H$, respectively. Consider the vertical plane containing $\gamma'$ and $\gamma''$. Form a portion $P$ of a 2-cell by taking the region of this plane lying between $\gamma'$ and $\gamma''$ and its intersection with the exterior of the faces $f'$ and $f''$. For example, in \reffig{DualFaces}, the lightly shaded region lying above two isometric spheres\index{isometric sphere} and running into infinity forms the portion $P$ of the 2-cell. Two such portions of planes are identified if they are identified by a parabolic\index{parabolic} translation of $\Gamma$ fixing the point at infinity.

The region $P$ and its translates under $\Gamma_\infty$ do not form the entire 2-cell; the 2-cell is formed by gluing portions of such regions $P$ along their intersections with isometric spheres.\index{isometric sphere} Note that each $P$ meets the faces $f'$ and $f''$ at right angles. Form a 2-cell by gluing portions of planes via the face-pairing isometries\index{face-pairing isometry} of $\calF_j$. That is, if $f'$ is glued to $\bar{f}'$ by a face-pairing isometry, then $P\cap f'$ is glued to $\bar{f}'$ by the same isometry. Similarly for $P\cap f''$. In \reffig{DualFaces}, the darker shaded regions form additional portions of the 2-cell; they are obtained by such gluings. 

\begin{lemma}\label{Lem:Dual2Cell}
  The 2-cells constructed as above are totally geodesic ideal polygons with $n\geq 3$ sides, where $n$ is the number of cusps equidistant from the image of the edge $e$ in the component of the Ford domain $\calF_j$ of $M$. 
\end{lemma}

\begin{proof}
\Refex{Dual2Cell}.
\end{proof}

The totally geodesic ideal polygon thus constructed is the \emph{geometric dual} of the edge $e$.\index{geometric dual}

\begin{example}\label{Example:Fig8DualFaces}
For the figure-8 knot complement, with equivariant Ford domain shown in \reffig{Fig8FordDomain}, each portion $P$ of a 2-cell lies over two faces of the Ford domain, running between their centers. A cross sectional view is shown on the left of \reffig{DualFaces}. In this example, three distinct portions $P$ glue to form an ideal triangle.\index{ideal triangle} Thus the 2-cells in this case are all ideal triangles.
\end{example}

An example for an edge equidistant from more than three horoballs is shown on the right in \reffig{DualFaces}. 

As in the case of the 1-cells, two 2-cells are identified if they differ by an element of $\Gamma$. 

Finally we construct the 3-cells. For each vertex of $\calF_j$, there are finitely many adjacent edges of $\calF_j$. The intersection of $\calF_j$ with the region bounded by the corresponding 2-cells forms a portion of a 3-cell $C$. The full 3-cell is formed by gluing $C\cap \bdy\calF_j$ by face-pairing isometries.\index{face-pairing isometry} The 3-cell is the \emph{geometric dual} to the vertex.\index{geometric dual} Again 3-cells are identified if they differ by an element of $\Gamma$. Note each 3-cell is bounded by a finite number of ideal polygons, thus it is an ideal polyhedron. 

\begin{definition}\label{Def:CanonicalDecomp}
The \emph{canonical polyhedral decomposition} of $M$, or simply the \emph{canonical decomposition}\index{canonical decomposition} of $M$ is the disjoint union of the 3-cells dual to the vertices of each $\calF_j$, along with their ideal faces and ideal edges. The faces are paired by isometries of $\Gamma$. 
\end{definition}

\begin{theorem}\label{Thm:CanonicalConvex}
If $M$ is a finite volume, orientable, cusped 3-manifold, with fixed embedded horoball neighborhood $H$ of all cusps, then the canonical polyhedral decomposition associated with $H$ decomposes $M$ uniquely into a finite number of convex ideal polyhedra. That is, $M$ is the quotient of the polyhedra with faces identified via face-pairing isometries,\index{face-pairing isometry} and the interiors of the polyhedra are mapped in a one-to-one manner to a subset of $M$. 
\end{theorem}

\begin{proof}[Proof sketch]
  The fact that there are finitely many such polyhedra follows from the fact that Ford domains are finite sided. Convexity follows from the fact that the polyhedra are cut out by finitely many ideal polygons dual to the convex Ford domain. All translates of the polyhedra under $\Gamma$ cover each $\calF_j$, by their definition, thus the polyhedra project surjectively to $M$. If $x$ lies in the interior of a polyhedron, projecting to some $y\in M$, then any other point $\tilde{x}$ in a 3-cell projecting to $y$ differs from $x$ by an element of $\Gamma$. By our definition of the 3-cells, it follows that the two 3-cells must agree. Thus $x$ is the only point in the polyhedra that projects to $y$. 
\end{proof}

\begin{example}
The 3-cells in the canonical polyhedral decomposition of the figure-8 knot are finite sided polyhedra whose faces are ideal triangles,\index{ideal triangle} by \refexample{Fig8DualFaces}. In fact, they are two regular ideal tetrahedra. The canonical polyhedral decomposition of the figure-8 knot is exactly the decomposition we obtained in \refchap{Fig8Decomp} and \refchap{GluingCompleteness}. 
\end{example}

\begin{theorem}\label{Thm:CanonicalUnique}
  Suppose $M$ and $M'$ are hyperbolic 3-manifolds, each with a single cusp, and suppose $M$ and $M'$ have combinatorially equivalent canonical decompositions. Then $M$ and $M'$ are isometric.

  In particular, if $M$ and $M'$ are knot complements, then $K$ and $K'$ are equivalent knots, up to reflection.
\end{theorem}

\begin{proof}
The canonical decomposition is geometrically dual to the Ford domain, and so this result follows from \refthm{CombinatoriallyEquivalent}. 
\end{proof}

In his thesis, Gu{\'e}ritaud showed that the canonical polyhedral decomposition of 2-bridge knots are the tetrahedra in the triangulation we described in \refchap{TwoBridge}~\cite{Gueritaud:thesis}; see also \cite{aswy}. There are a few other families of 3-manifolds whose canonical polyhedral decompositions are now known; for example Sakuma and Weeks find canonical polyhedra for some families of link complements using symmetry~\cite{SakumaWeeks}. SnapPy computes canonical polyhedral decompositions of most reasonably sized 3-manifolds~\cite{SnapPy}, but at the date of writing this book, such decompositions were not rigorously verified. In general it seems to be a difficult problem to determine canonical polyhedral decompositions of important families of link complements. For example, as of the writing of this book, the following is unknown, asked as a question in \cite{SakumaWeeks}.

\begin{conjecture}\label{Conj:CrossingArcGeodesic}
  For any hyperbolic alternating knot\index{alternating knot or link} and any alternating diagram\index{alternating diagram} of the knot, each crossing arc\index{crossing arc} is isotopic to an edge of the canonical polyhedral decomposition of the knot complement. 
\end{conjecture}

In fact, to date it is not even known if crossing arcs of alternating knots are isotopic to geodesics in general.

%%%%%%%%%%%%%%%%%%%%%%%%%%%%%%%%%%%%%%%%%%%%%%%%%%%%%%%%%%%%%%%%%
\section{Exercises}

\begin{exercise}\label{Ex:CountableDiscrete}
  Prove that a discrete subgroup of $\PSL(2,\CC)$ is countable.
\end{exercise}

\begin{exercise}
  Prove that a complete hyperbolic surface contains an embedded neighborhood of its cusps that lifts to countably many horodisks in $\HH^2$.
\end{exercise}

\begin{exercise}\label{Ex:DefIsoSphereDoesNotDependOnH}
  Prove that the definition of the set $I(g)$\index{$I(g)$} is independent of choice of horosphere $H$.
\end{exercise}

\begin{exercise}\label{Ex:CuspVolumeMaxCusp}
  Prove \refcor{CuspVolume}, that the volume of a cusp component in a maximal cusp neighborhood\index{maximal cusp neighborhood}\index{cusp!maximal cusp neighborhood} of $M$ is at least $\sqrt{3}/4$, and that if $M$ has only one cusp, then the volume of a maximal cusp neighborhood is at least $\sqrt{3}/2$.

Hints: If desired, you may use \refex{CuspVolume}. You may also apply the following theorem of B\"or\"oczky \cite{boroczky}, generalizing \refthm{Boroczky}.

\begin{theorem}[B\"or\"oczky density of disks in the torus]\label{BoroczkyTorus}\index{B{\"o}r{\"o}czky cusp density theorem!3-dimensional}\index{cusp density theorem!2-dimensional}
Let $D$ be a collection of disks of the same radius, embedded disjointly in the torus $T$. Then 
\[ \frac{\area(T\cap D) }{\area(T)} \leq \frac{\pi}{2\sqrt{3}}. \]
\end{theorem}
\end{exercise}

\begin{exercise}
  Extend the definition of a Ford domain to complete hyperbolic surfaces, first for those with one cusp, then generalize to finitely many cusps. Using the natural extension of the definition of a fundamental domain to hyperbolic surfaces, prove that the object you have defined is a convex fundamental domain for the surface. 
\end{exercise}

\begin{exercise}
This exercise uses SnapPy \cite{SnapPy} to investigate the combinatorics of Ford domains associated with different choices of maximal cusp neighborhoods.\index{maximal cusp neighborhood}\index{cusp!maximal cusp neighborhood} The manifold {\tt{m125}} in the SnapPy census is a hyperbolic manifold with two cusps, isometric to the complement of the (-2,3,8)-pretzel link.
\begin{enumerate}
\item Using SnapPy, find at least five different combinatorially inequivalent Ford domains for {\tt{m125}}.
\item Find a Ford domain for {\tt{m125}} where the dual canonical decomposition is a triangulation. Find a Ford domain where it is not a triangulation.
\end{enumerate}
\end{exercise}

\begin{exercise}\label{Ex:Dual2Cell}
Prove \reflem{Dual2Cell}: that the 2-cell dual to an edge of $\calF$ is a totally geodesic ideal polygon with $n\geq 3$ sides, where $n$ is the number of cusps equidistant from the projection of the edge to $M$.
\end{exercise}

%% \begin{exercise}
%% In \cite{GueritaudSchleimer}, the Ford domain is defined differently from our definition. Given a complete hyperbolic 3-manifold $M$ of finite volume, and choice of embedded horoball neighborhoods $H_1, \dots, H_n$, they define the Ford--Voronoi domain to consist of all points in $M$ having a unique shortest path to the union of the $H_i$. Explain how to obtain this domain from the equivariant Ford domain of \refdef{FordDomainMultiple} of a 3-manifold with $n$ cusps. 
%% \end{exercise}

\chapter{Algebraic Sets and the $A$-Polynomial}\label{Chap:Character}
\blfootnote{Jessica S. Purcell, Hyperbolic Knot Theory}

In this chapter, we introduce a polynomial invariant of a knot that is obtained from hyperbolic geometry. There are three closely related polynomials that appear in the literature and in calculations. We introduce all three in this chapter and discuss the relationships between them. We need to introduce a small amount of algebraic geometry, to consider how hyperbolic structures deform. M.~Culler and P.~Shalen were the first to investigate this material \cite{CullerShalen:Varieties}. The $A$-polynomial was originally introduced in \cite{APoly}. However, we start with a slightly different perspective.

\section{The gluing variety}

Suppose $M$ is a 3-manifold with a topological ideal triangulation,\index{topological ideal triangulation} as in \refdef{TopIdealTriang}. We have seen that associated with each ideal tetrahedron is an edge invariant,\index{edge invariant} namely a complex number $z(e)$ corresponding to an edge of the tetrahedron, as in \refdef{EdgeInvariant}, with all edge invariants of a tetrahedron satisfying the relations of \reflem{EdgeInvariants}. Finally, we have seen in \refthm{Gluing} that a choice of edge invariants gives a hyperbolic structure on $M$ if and only if the $z(e)$ satisfy the edge gluing equations.\index{gluing equations!edge equations}\index{edge gluing equations}

Suppose that $n$ ideal tetrahedra form the topological ideal triangulation of $M$. Then we obtain $3n$ edge invariants: $z_1, \dots, z_n$, as well as $z_i'=1/(1-z_i)$ and $z_i''=(z_i-1)/z_i$ for $i=1, \dots, n$. The edge gluing equations tell us that the product of those edge invariants that glue to the same edge of $M$ must be $1$. Clearing denominators, each edge gluing equation becomes a polynomial equation in $z_1, \dots, z_n$.

\begin{definition}\label{Def:AlgebraicSet}
  An \emph{affine algebraic set}\index{affine algebraic set} is a subset of $\CC^N$ that is defined as the set of zeros of a system of polynomial equations with coefficients in $\CC$ and $N$ variables.
\end{definition}

The union of two affine algebraic sets is an affine algebraic set, and the intersection of arbitrarily many affine algebraic sets is an affine algebraic set (\refex{Zariski}). We define a topology on $\CC^N$, called the \emph{Zariski topology},\index{Zariski topology} by taking affine algebraic sets to be closed sets. 
As a very simple preliminary example, the subset of $\CC^2$ defined by $xy=0$ is an affine algebraic set. 

\begin{definition}\label{Def:ReducibleAlgebraicSet}
An affine algebraic set is \emph{reducible}\index{reducible affine algebraic set}\index{affine algebraic set!reducible, irreducible} if it can be expressed as the union of two proper affine algebraic sets. It is \emph{irreducible}\index{irreducible affine algebraic set} if not. 

It is a corollary of the Hilbert basis theorem (says Shalen \cite{shalen:survey}) that any affine algebraic set is a finite union of irreducible affine algebraic sets; we call these the \emph{irreducible components}\index{irreducible components}. 

An irreducible affine algebraic set is called an \emph{affine variety}.\index{affine variety}
\end{definition}

The affine algebraic set defined by $xy=0$ is not an affine variety;\index{affine variety} it is reducible, with irreducible affine algebraic sets defined by polynomials $x=0$ and $y=0$. These form irreducible components. Note they are not disjoint. In general, irreducible components\index{irreducible components} of affine algebraic sets are not necessarily disjoint.
%% Note this is very different from (connected) components of topological spaces. 

\begin{corollary}\label{Cor:GluingAlgebraic}
Let $\calT$ be a topological ideal triangulation of a 3-manifold $M$, with $n$ ideal tetrahedra. Then the set of points in $\CC^n$ satisfying the edge gluing equations associated with $\calT$ forms an affine algebraic set. \qed
\end{corollary}

\begin{definition}\label{Def:GluingVariety}
Suppose $M$ has ideal triangulation $\calT$ made up of $n$ ideal tetrahedra (so $n$ gluing equations by \refex{edge=tet}). 
The \emph{gluing variety}\index{gluing variety} associated with $\calT$, denoted $\calD(\calT)$, is the affine algebraic subset of $\CC^n\times \CC$ consisting of points $(z_1, \dots, z_n, t)$ satisfying the edge gluing equations as in \refcor{GluingAlgebraic}, as well as the equation
\[ t\prod_{i=1}^n z_i(1-z_i)=1. \]
This last equation is called the \emph{degeneracy equation},\index{degeneracy equation} and ensures that the parameters $z_i$ do not degenerate to $0$ or $1$. 
\end{definition}

\begin{remark}
Degeneracy is handled differently by different authors in the literature. For example, rather than including one degeneracy equation, some authors require that for each $i$, the equations
\[ z_i(1-z_i'')=1, \quad z_i'(1-z_i)=1, \quad z_i''(1-z_i')=1 \]
hold, in addition to gluing equations. This also rules out $z_i=0,1$.

Similarly, some authors require that for each $i$, the following equations hold in addition to the gluing equations:
\[ z_iz_i'z_i''=-1 \quad \mbox{and} \quad z_i + (z_i')^{-1} -1 = 0. \]
Encoded in these equations are the relationships between $z_i$, $z_i'$, and $z_i''$ from \reflem{EdgeInvariants}. Again a solution does not allow $z_i=0,1$. 
\end{remark}

Note that the gluing variety is an affine algebraic set, but it is not necessarily irreducible.\index{irreducible affine algebraic set} In particular, the gluing variety may not satisfy the definition of an affine variety,\index{affine variety} so the terminology is somewhat misleading.

In the case that $M$ has a complete hyperbolic structure, and $\calT$ can be given a collection of edge invariants to form a geometric ideal triangulation for this structure, as in \refdef{GeomTriang}, then the gluing variety is nonempty. That is, there is an irreducible component\index{irreducible component}\index{irreducible affine algebraic set} of the gluing variety that contains this complete hyperbolic structure. The irreducible component containing the complete hyperbolic structure (when it exists) really is an affine variety.\index{affine variety}

\begin{definition}\label{Def:CanonicalComponent}
Suppose $\calT$ is an ideal triangulation of $M$ that can be assigned edge invariants to form a geometric ideal triangulation, as in \refdef{GeomTriang}. Then the irreducible component\index{irreducible component} of $\calD(\calT)$ containing the geometric triangulation is called the \emph{canonical component}\index{canonical component}\index{canonical component!gluing variety} of the gluing variety, and it is denoted $\calD_0(\calT)$.
\end{definition}

Some examples are in order.

\begin{example}[Figure-8 knot]
Take the figure-8 knot complement with the ideal triangulation of \refchap{Fig8Decomp}. This has two ideal tetrahedra, with edge invariants determined by $z$ and $w$ in $\CC$. The two edge gluing equations of the figure-8 knot actually both give a single polynomial equation, which we calculated in \refeqn{Fig8Gluing}:
\begin{equation}\label{Eqn:Fig8GluingPolynomial}
z(z-1)w(w-1)=1, \quad \mbox{ or } \quad z^2w^2 - z^2w - zw^2 + zw -1 = 0. 
\end{equation}
Note that in this case, the degeneracy equation\index{degeneracy equation} $t z (1-z)w(1-w)=1$ becomes simply $t=1$, and so the gluing variety\index{gluing variety} is therefore the set of points $(z, w, 1)$ in $\CC^2\times \CC$ satisfying \refeqn{Fig8GluingPolynomial}. 

Note that the complete hyperbolic structure occurs when
\[z=w=(1+i\sqrt{3})/2,\]
as shown in \refexamp{Fig8Completeness}. The triple
\[ \left( (1+i\sqrt{3})/2, (1+i\sqrt{3})/2, 1 \right) \in \CC^2\times\CC \] satisfies \refeqn{Fig8GluingPolynomial}.
\end{example}

\begin{example}[$6_1$ knot]
In \refchap{GluingCompleteness}, we found a triangulation $\calT_1$ of the $6_1$ knot with five tetrahedra. There is also a triangulation $\calT_2$ of the $6_1$ knot  with only four tetrahedra, for example this triangulation is stored in the SnapPy database of manifold \cite{SnapPy}. Note that the gluing variety\index{gluing variety} corresponding to $\calT_1$ is a subset of $\CC^6$, while that corresponding to $\calT_2$ is a subset of $\CC^5$. Thus even though the triangulations give the same manifold, the gluing varieties are different, living in different spaces. 
\end{example}

\begin{example}[Triangulation with inessential edge]
Dunfield notes that the figure-8 knot complement admits a topological ideal triangulation $\calT$ consisting of five tetrahedra where one of the edges $E$ is 1-valent: it meets only one of the tetrahedra \cite{Dunfield:Appendix}. Moreover, this edge $E$ is \emph{inessential}\index{inessential edge}, meaning for any horoball neighborhood $H$ of the knot, the arc $E\cap (S^3-H)$ is homotopic rel endpoints $E\cap H$ into $\bdy H$. 
For this triangulation $\calT$, note that the edge gluing equation corresponding to the inessential edge is the equation $z_1=1$. Because it is impossible for this equation to hold simultaneously with the degeneracy equation,\index{degeneracy equation} for this triangulation the gluing variety\index{gluing variety} $\calD(\calT)$ is empty. 
\end{example}

\begin{remark}
Suppose $M$ is the interior of a compact manifold with torus boundary. Segerman and Tillmann showed that for any ideal triangulation $\calT$ of $M$, the gluing variety\index{gluing variety} associated with $\calT$ is nonempty if and only if all edges in the triangulation are essential \cite{SegermanTillmann}.
\end{remark}

\subsection{The $\Ahyp$-polynomial}

We now describe a two-variable polynomial $\Ahyp(\ell,m)$ associated with any triangulated knot complement, or indeed any triangulated 3-manifold that is the interior of a compact manifold with a single torus boundary component. It was introduced by Champanerkar in \cite{Champanerkar:Thesis}, and there related to a more common polynomial, which we will describe later in this chapter. For now, we have the tools to understand the $\Ahyp$-polynomial and compute examples, so we do this first.

For a fixed triangulation $\calT$, the gluing variety $\calD(\calT)$ is defined using edge gluing equations only. Recall that we also have two completeness equations for each cusp, by \refprop{CompletenessEqns}. That is, given a knot complement $S^3-K$, with meridian $[\mu]$ and longitude $[\lambda]$, for $[\mu], [\lambda] \in \pi_1(\bdy N(K))$, there are associated equations $H(\mu)$ and $H(\lambda)$, which are rational functions in the $z_i$. Define a map $\calH\from \calD \to \CC\times\CC$ by
\[ \calH(z_1, \dots, z_n) = (H(\lambda)(z_1, \dots, z_n), H(\mu)(z_1, \dots, z_n)) = (\ell, m). \]

\begin{definition}
Let $Y_i$ be an irreducible component\index{irreducible component} of $\calD(\calT)$ for which the closure of $\calH(Y_i)$ in $\CC\times\CC$ (in the Zariski topology) is an affine algebraic set, and let $Z$ be the union of all such components.
The image of $Z$ under $\calH$ is called the \emph{holonomy variety}\index{holonomy variety} with respect to the triangulation $\calT$. The defining polynomial of the closure of the image is denoted by $\AhypT(\ell, m)$. We occasionally omit the notation $\calT$ when the triangulation is understood. The polynomial is called the \emph{$\Ahyp$-polynomial}\index{$\Ahyp$-polynomial} or \emph{hyperbolic $A$-polynomial}.\index{hyperbolic $A$-polynomial}\index{$A$-polynomial!hyperbolic, $\Ahyp$}
\end{definition}

\begin{example}[Figure-8 knot]\label{Ex:HPolynomialFig8}
For $\calT$ the usual two tetrahedra triangulation of the figure-8 knot $K$, the $\AhypT$-polynomial satisfies \refeqn{Fig8GluingPolynomial} as well as two rational equations in $m$ and $\ell$ that come from the completeness equations. We computed completeness equations in \refexample{Fig8EdgeSequence}, for curves $\alpha$ and $\beta$. The curve $\beta$ is the meridian of the figure-8 knot. The curve $\alpha$ generates $\pi_1(\bdy N(K))$ along with $\beta$, but it is not the standard longitude as defined in \refrem{MeridLongitude}. Still, we may use it to compute a polynomial. That is, we require
\begin{align*}
m &= z_2^{-1}w_1 = (1-z) w, \quad \mbox{and} \\
\ell & = \left(\frac{z_2z_3}{w_2 w_3}\right)^2 = \frac{1}{(1-z)^2}\frac{(z-1)^2}{z^2}\frac{(1-w)^2}{1}\frac{w^2}{(w-1)^2} = \frac{w^2}{z^2}
\end{align*}

Finding the polynomial $\Ahyp$, as an equation in $m$ and $\ell$ only, is now a problem in elimination theory. One way to solve it is to form an ideal generated by the above equations in the ring $\ZZ[z,w,\ell,m]$ and compute a Groebner basis that eliminates $z$ and $w$ using tools from computational algebraic geometry. For our examples, we did the latter, using Mathematica. The following command returns a polynomial in $m$ and $\ell$.
\begin{quote}
\begin{verbatim}
system := {z*(z-1)*w*(w-1)==1, M==w*(1-z), L*z^2==w^2};
elim := {z,w}; keep := {M,L}; 
GB := GroebnerBasis[system, keep, elim]; 
Print[GB];
\end{verbatim}
\end{quote}

This computation yields the following 2-variable polynomial:\index{$\Ahyp$-polynomial}\index{$A$-polynomial!hyperbolic, $\Ahyp$}
\begin{multline*}
\Ahyp(\ell, m) = \ell - 2\ell m - 3 \ell m^2 - \ell^2 m^2 + 2 \ell m^3 + 6 \ell m^4 \\
+ 2 \ell m^5 - m^6-3 \ell m^6 - 2 \ell m^7 + \ell m^8.
\end{multline*}
\end{example}

There were several choices made in the last example, particularly concerning representatives for the curves that give the completeness equations. We could change these in several different ways.

\begin{figure}
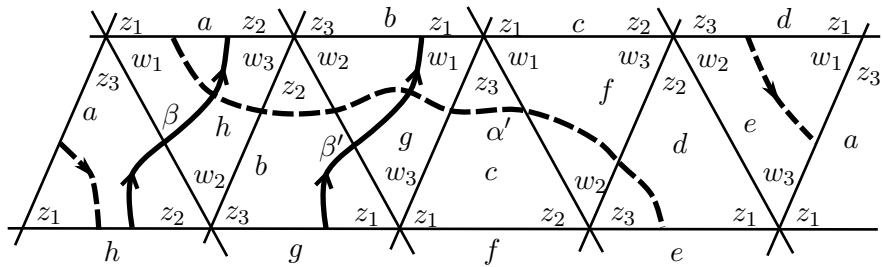
\caption{The cusp triangulation of the figure-8 knot complement, along with isotopic meridians $\beta$ and $\beta'$, and a longitude $\alpha'$.}
\label{Fig:Fig8CuspTriangAgain}
\end{figure}

In the cusp triangulation of the figure-8 knot, reproduced in \reffig{Fig8CuspTriangAgain}, 
a representative $\beta$ for the meridian was chosen to run from the bottom of the triangle labeled $a$, across the triangle labeled $h$, back to the triangle labeled $a$, giving $H(\beta) = w_1/z_2$. We could have chosen a parallel curve, for example the curve $\beta'$ running from the base of $b$, across $g$ to the base of $b$, giving the curve $H(\beta') = w_2/z_1$.

Recall that in \refchap{GluingCompleteness}, to define $H([\gamma])$ we require that a curve $\gamma$ on a cusp triangulation cuts off a single corner of each triangle it enters. Such a curve is called a \emph{normal curve}.\index{normal curve}

\begin{proposition}\label{Prop:AhypRepChoice}
Suppose $\alpha, \alpha' \in [\lambda]$ are distinct normal curves in the homotopy class of $\lambda$. Then the polynomial obtained by replacing $H(\alpha)$ by $H(\alpha')$ leaves $\Ahyp$ unchanged. Similarly for $[\mu]$. 
\index{$A$-polynomial!hyperbolic, $\Ahyp$}\index{$\Ahyp$-polynomial}
\end{proposition}

\begin{proof}
This follows from \refex{HomotopyInvarianceH}: $H([\alpha])$ is independent of choice of $\alpha$. The reason why is that if normal curves $\alpha$ and $\alpha'$ are ambient isotopic\index{ambient isotopic} on a triangulated torus, then they differ by sliding past edges of the triangulation. Because the product of all tetrahedron edge parameters $z_i$ meeting at an edge is $1$, and the edge gluing equations are satisfied for the parameters $z_i$, it follows that the two polynomials in the definition of $\Ahyp$ will differ by a factor of $1$. 
\end{proof}

We also could have chosen a different basis for $H_1(\bdy N(K))$. For example, if $\alpha'$ denotes the standard longitude of the figure-8 knot, shown in \reffig{Fig8CuspTriangAgain}, then it can be shown that
\[
H(\alpha') = z_1^{-1}\cdot w_3 \cdot z_2 \cdot w_3^{-1}\cdot z_3 \cdot w_2^{-1} \cdot z_3^{-1}\cdot w_1. 
\]

Changing the basis for $H_1(\bdy N(K))$ does affect $\Ahyp$, but in a  well-understood way.

\begin{proposition}\label{Prop:AhypChangeBasis}
Suppose $p, q, r, s$ are integers satisfying $ps-rq=1$, so that $\langle \ell^p m^q,\ell^r m^s \rangle$ is another basis for $H_1(\bdy N(K))$. Then there exist integers $a, b$ such that \index{$A$-polynomial!hyperbolic, $\Ahyp$}\index{$\Ahyp$-polynomial}
\[ \Ahyp(\ell^p m^q, \ell^r m^s) = \pm\ell^a m^b \Ahyp(\ell, m). 
\]
\end{proposition}

\begin{example}\label{Example:Fig8_AHypPoly}
The standard longitude of the figure-8 knot can be shown to be represented by the curve $\alpha'$ of \reffig{Fig8CuspTriangAgain}, which is isotopic to $\alpha\beta^2$. Replacing $\alpha$ by $\alpha'$ yields the polynomial:
\begin{multline*}
\Ahyp'(\ell,m) := \ell - 2 \ell m - 3 \ell m^2 + 2 \ell m^3 - m^4 + 6 \ell m^4  \\ - \ell^2 m^4 + 2 \ell m^5 - 3 \ell m^6 - 2 \ell m^7 + \ell m^8.
\end{multline*}
Note that this new polynomial satisfies $\Ahyp'(\ell,m) = m^{-2}\Ahyp(\ell m^2, m)$. \index{$A$-polynomial!hyperbolic, $\Ahyp$}\index{$\Ahyp$-polynomial}
\end{example}

%%%%%%%%%%%%%%%%%%%%%%%%%%%%%%%%%%%%%%%%%%%%%%%%%%%%%%%%%%%%%%%%%
\section{Representations of knots}

The original $A$-polynomial was defined in \cite{APoly} using representations of the fundamental group of a knot complement into the group $\SL(2,\CC)$. In this section, we work our way up to the original definition of the $A$-polynomial. Our goals are first, to give a taste of the interesting mathematics along the way, but more importantly, to relate the polynomial $\Ahyp$ to the original $A$-polynomial. 

\subsection{Wirtinger presentation}
We take a bit of a detour in this subsection, to recall a few important results on presentations of knot groups. The tools we need follow almost immediately from the \emph{Wirtinger presentation},\index{Wirtinger presentation} which we review now. See also \cite{rolfsen}.

\begin{figure}
\includegraphics{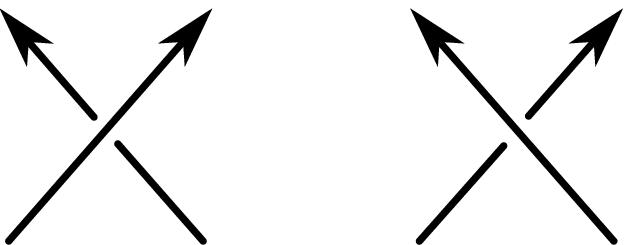}
  \caption{A positive\index{crossing!positive} crossing (left), and a negative\index{crossing!negative} crossing (right).}
  \label{Fig:PosNegWirtinger}
\end{figure}

To obtain the presentation, start with a diagram of the knot, which we view as a collection of arcs, with each arc starting and ending at undercrossings and running along only overcrossings (if any) between them. To simplify the discussion, and to fix notation, choose an orientation of the knot (either orientation is fine). We say that a crossing is \emph{positive}\index{crossing!positive} or \emph{negative}\index{crossing!negative} based on the orientation of the arcs at the crossing, as in \reffig{PosNegWirtinger}.

\begin{figure}
%% Accompanies image file 'F15-03-Wirt.eps' (pdf, eps, ps)
\begingroup%
  \makeatletter%
  \providecommand\color[2][]{%
    \errmessage{(Inkscape) Color is used for the text in Inkscape, but the package 'color.sty' is not loaded}%
    \renewcommand\color[2][]{}%
  }%
  \providecommand\transparent[1]{%
    \errmessage{(Inkscape) Transparency is used (non-zero) for the text in Inkscape, but the package 'transparent.sty' is not loaded}%
    \renewcommand\transparent[1]{}%
  }%
  \providecommand\rotatebox[2]{#2}%
  \newcommand*\fsize{\dimexpr\f@size pt\relax}%
  \newcommand*\lineheight[1]{\fontsize{\fsize}{#1\fsize}\selectfont}%
  \ifx\svgwidth\undefined%
    \setlength{\unitlength}{179.31053925bp}%
    \ifx\svgscale\undefined%
      \relax%
    \else%
      \setlength{\unitlength}{\unitlength * \real{\svgscale}}%
    \fi%
  \else%
    \setlength{\unitlength}{\svgwidth}%
  \fi%
  \global\let\svgwidth\undefined%
  \global\let\svgscale\undefined%
  \makeatother%
  \begin{picture}(1,0.37864283)%
    \lineheight{1}%
    \setlength\tabcolsep{0pt}%
    \put(0,0){\includegraphics[width=\unitlength]{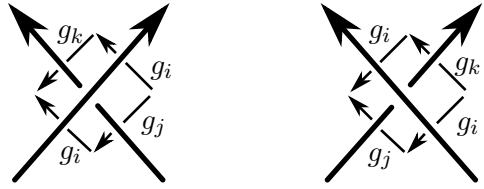}}%
    \put(0.29919061,0.2287136){\color[rgb]{0,0,0}\makebox(0,0)[lt]{\lineheight{0}\smash{\begin{tabular}[t]{l}$g_i$\end{tabular}}}}%
    \put(0.10798202,0.0502528){\color[rgb]{0,0,0}\makebox(0,0)[lt]{\lineheight{0}\smash{\begin{tabular}[t]{l}$g_i$\end{tabular}}}}%
    \put(0.28006976,0.11398841){\color[rgb]{0,0,0}\makebox(0,0)[lt]{\lineheight{0}\smash{\begin{tabular}[t]{l}$g_j$\end{tabular}}}}%
    \put(0.1032018,0.30360359){\color[rgb]{0,0,0}\makebox(0,0)[lt]{\lineheight{0}\smash{\begin{tabular}[t]{l}$g_k$\end{tabular}}}}%
    \put(0.94292624,0.11080184){\color[rgb]{0,0,0}\makebox(0,0)[lt]{\lineheight{0}\smash{\begin{tabular}[t]{l}$g_i$\end{tabular}}}}%
    \put(0.7580912,0.30519715){\color[rgb]{0,0,0}\makebox(0,0)[lt]{\lineheight{0}\smash{\begin{tabular}[t]{l}$g_i$\end{tabular}}}}%
    \put(0.75649791,0.04706568){\color[rgb]{0,0,0}\makebox(0,0)[lt]{\lineheight{0}\smash{\begin{tabular}[t]{l}$g_j$\end{tabular}}}}%
    \put(0.93814604,0.23827388){\color[rgb]{0,0,0}\makebox(0,0)[lt]{\lineheight{0}\smash{\begin{tabular}[t]{l}$g_k$\end{tabular}}}}%
  \end{picture}%
\endgroup%

  \caption{Relators for positive (left) and negative (right) crossings}
  \label{Fig:WirtingerRelators}
\end{figure}

Now each arc $a_i$ of the diagram is oriented. There will be one generator $g_i$ of the fundamental group for each arc $a_i$; the generator can be viewed as a loop beginning at a basepoint that lies high above the plane of projection, looping around $a_i$ in a positive direction, then running back to the basepoint. Relators for the group presentation come from crossings, as follows. Suppose arc $a_i$ runs along an overcrossing, meeting endpoints of arcs $a_j$ and $a_k$ at the crossing. If it is a positive crossing, then we have the relation $g_j g_i g_k^{-1} g_i^{-1} = 1$. If it is a negative crossing, we have instead $g_j g_i^{-1} g_k^{-1} g_i=1$. See \reffig{WirtingerRelators}.

\begin{theorem}[Wirtinger presentation]\label{Thm:Wirtinger}
If $K$ has a diagram with $n$ crossings, then the group $\langle g_1, \dots, g_n \mid r_1, \dots, r_n \rangle$ obtained as above forms a presentation for the fundamental group of the knot complement. Moreover, one of the relators is redundant and can be removed from the presentation. 
\end{theorem}

\begin{proof}
The proof is a standard exercise in algebraic topology, using the Seifert--Van Kampen theorem. We include it as \refex{SVK}. 
\end{proof}

We give a few immediate corollaries.

\begin{corollary}\label{Cor:Homology}
The first homology group of a knot complement is always isomorphic to $\ZZ$, generated by (the class of) the meridian.
\end{corollary}

\begin{proof}
The first homology group $H_1(S^3-K)$ is the abelianization of the fundamental group $\pi_1(S^3-K) \cong \langle g_1, \dots, g_n \mid r_1, \dots, r_n\rangle$. Under the abelianization, each relator becomes $g_j=g_k$. Thus the presentation reduces to $\langle g_1,\dots,g_n \mid g_1=\dots=g_n\rangle \cong \ZZ$, generated by the class of the meridian $g_1$. 
\end{proof}

Note that each generator of the Wirtinger presentation is a curve that bounds a disk in $S^3$, and can be isotoped to an embedded curve on the boundary of a neighborhood of the knot. 

Generators of the Wirtinger presentation are meridians (\refdef{MeridianLongitude}). 
This leads to another corollary.

\begin{corollary}\label{Cor:ParabolicGenerator}
Let $K$ be a knot whose complement admits a hyperbolic structure. Then the holonomy group\index{holonomy group} of the complement of $K$ has a presentation in which all generators are parabolic.\index{parabolic}
\end{corollary}

Recall from \refdef{holonomy} that the holonomy group is the image of the holonomy\index{holonomy} map $\rho\from \pi_1(S^3-K)\to\PSL(2,\CC)$, encoding the hyperbolic structure.

\begin{proof}[Proof of \refcor{ParabolicGenerator}]
Each generator of the Wirtinger presentation is a meridian. Each meridian is mapped by the holonomy\index{holonomy} map to an element of $\PSL(2,\CC)$ that commutes with a longitude. Then \refcor{ZxZSubgroup} ($\ZZ\times\ZZ$ subgroups) implies that meridian and longitude are mapped to parabolic\index{parabolic} elements fixing the same point on $\bdy\HH^3$. 
\end{proof}

\begin{lemma}
Given a Wirtinger presentation of the fundamental group of a knot complement, we may compute the longitude as a product of generators as follows.
\begin{enumerate}
\item Beginning on an arc $a_k$ of the knot, travel along it.
\item Write $g_i$ when traversing an undercrossing under arc $a_i$ in the positive direction, and write $g_i^{-1}$ when traversing in the negative direction.
\item Finally, multiply by $g_k^p$ where the power $p$ is chosen such that the total exponent sum is zero.
\end{enumerate}
\end{lemma}

\begin{proof}
The loop obtained as the product of generators corresponding to undercrossings as above will be homotopic to some longitude. We want to ensure it is homologically trivial. Consider its image in the abelianization of the fundamental group $H_1(S^3-K)\cong \ZZ$. Each $g_i$ maps to the generator $g$ of $\ZZ$ in homology. Therefore if the sum of the exponents in the product is zero, the image of the element in homology is trivial. 
\end{proof}

\begin{example}\label{Example:Fig8Wirtinger}
We use the Wirtinger presentation to compute the fundamental group of the figure-8 knot complement. See \reffig{Fig8Wirtinger}. \index{Wirtinger presentation}

\begin{figure}
%% Accompanies image file 'F15-04-WirtF8.eps' (pdf, eps, ps)
\begingroup%
  \makeatletter%
  \providecommand\color[2][]{%
    \errmessage{(Inkscape) Color is used for the text in Inkscape, but the package 'color.sty' is not loaded}%
    \renewcommand\color[2][]{}%
  }%
  \providecommand\transparent[1]{%
    \errmessage{(Inkscape) Transparency is used (non-zero) for the text in Inkscape, but the package 'transparent.sty' is not loaded}%
    \renewcommand\transparent[1]{}%
  }%
  \providecommand\rotatebox[2]{#2}%
  \newcommand*\fsize{\dimexpr\f@size pt\relax}%
  \newcommand*\lineheight[1]{\fontsize{\fsize}{#1\fsize}\selectfont}%
  \ifx\svgwidth\undefined%
    \setlength{\unitlength}{119.25789642bp}%
    \ifx\svgscale\undefined%
      \relax%
    \else%
      \setlength{\unitlength}{\unitlength * \real{\svgscale}}%
    \fi%
  \else%
    \setlength{\unitlength}{\svgwidth}%
  \fi%
  \global\let\svgwidth\undefined%
  \global\let\svgscale\undefined%
  \makeatother%
  \begin{picture}(1,0.94545795)%
    \lineheight{1}%
    \setlength\tabcolsep{0pt}%
    \put(0,0){\includegraphics[width=\unitlength]{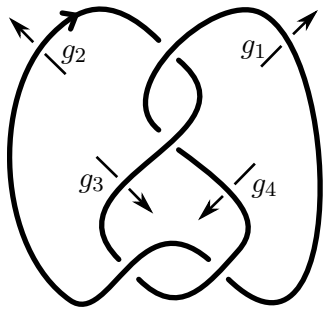}}%
    \put(0.16861528,0.78666807){\color[rgb]{0,0,0}\makebox(0,0)[lt]{\lineheight{0}\smash{\begin{tabular}[t]{l}$g_2$\end{tabular}}}}%
    \put(0.73598798,0.80076907){\color[rgb]{0,0,0}\makebox(0,0)[lt]{\lineheight{0}\smash{\begin{tabular}[t]{l}$g_1$\end{tabular}}}}%
    \put(0.2242058,0.3788402){\color[rgb]{0,0,0}\makebox(0,0)[lt]{\lineheight{0}\smash{\begin{tabular}[t]{l}$g_3$\end{tabular}}}}%
    \put(0.76991468,0.36332908){\color[rgb]{0,0,0}\makebox(0,0)[lt]{\lineheight{0}\smash{\begin{tabular}[t]{l}$g_4$\end{tabular}}}}%
  \end{picture}%
\endgroup%

\caption{Generators of the fundamental group of the figure-8 knot complement}
\label{Fig:Fig8Wirtinger}
\end{figure}

The generators are $g_1$, $g_2$, $g_3$, $g_4$. The relators are
\[ g_2 g_1^{-1}g_3^{-1}g_1=1 \quad g_4g_3^{-1}g_1^{-1}g_3=1 \quad
g_3g_2g_4^{-1}g_2^{-1}=1 \quad g_1g_4g_2^{-1}g_4^{-1}=1.\]
Use the first and third equations to eliminate $g_3$ and $g_4$, respectively. Substitute into the fourth equation to obtain the single relation
\[ g_1g_2^{-1}g_1g_2g_1^{-1}g_2 = g_2^{-1}g_1g_2g_1^{-1}g_2g_2. \]

Choose the meridian $M$ to be $g_1$.

Now starting on the arc of $g_1$ (right side of the diagram), traverse the knot by running up the arc, obtaining a longitude of the form
\begin{align*} L & = g_3g_2^{-1} g_1 g_4^{-1} \\
& = (g_1g_2g_1^{-1})g_2^{-1}g_1(g_2^{-1}g_1g_2^{-1}g_1^{-1}g_2). 
\end{align*}
Note that the sum of exponents in the example is zero, so this is the homologically trivial longitude.
\end{example}

\subsection{Representation space and character varieties}
In this subsection, we discuss representations of fundamental groups of 3-manifolds, particularly knot complements. We will only touch upon some of the details here; see \cite{shalen:survey} for a more comprehensive survey.

Suppose $M$ is a 3-manifold that is the interior of the compact manifold $\overline{M}$ with a single torus boundary component, and suppose that $M$ admits an ideal triangulation $\calT$. Then associated to $M$ is a gluing variety.
Any point in the gluing variety gives a representation of $\pi_1(M)$ into $\PSL(2,\CC)$, as follows. A point in the gluing variety gives a collection of edge parameters for the tetrahedra in $\calT$ that satisfy the edge gluing equations. These lead to a $(G,X)$-structure on $M$,\index{$(G,X)$-structure} where $G=\PSL(2,\CC)$ and $X=\HH^3$, by setting each tetrahedron in $\calT$ to be the ideal hyperbolic tetrahedron with given edge parameters. Then there is an associated developing map\index{developing map} as in \refdef{DevelopingMap}. A group element $\alpha\in\pi_1(M)$ determines a holonomy\index{holonomy} element $g_\alpha\in\PSL(2,\CC)$, as in \refdef{holonomy}. We have seen that the holonomy gives a homomorphism $\rho\from \pi_1(M)\to \PSL(2,\CC)$. 
Thus holonomy is a representation
\[ \rho\from\pi_1(M)\to\PSL(2,\CC) \cong \SL(2,\CC)/\pm\Id. \]

In fact, the original study of $A$-polynomials in \cite{APoly} considered representations to $\SL(2,\CC)$ rather than $\PSL(2,\CC)$, to avoid certain difficulties that arise. We begin with this perspective as well.

As usual, we focus on 3-manifolds that arise as knot complements. The following ensures that we have interesting representations to work with.

\begin{proposition}\label{Prop:HolonomyLiftstoSL}
Any representation $\rho\from\pi_1(S^3-K)\to\PSL(2,\CC)$ lifts to two representations $\widetilde{\rho}, \tilde{\rho}'\from \pi_1(S^3-K)\to\SL(2,\CC)$.
\end{proposition}

\begin{proof}
Let $\langle g_1, \dots, g_n \mid r_1, \dots, r_n \rangle$ be a Wirtinger presentation\index{Wirtinger presentation} for $\pi_1(S^3-K)$. Starting with $g_1$, the representation $\rho$ assigns $\rho(g_1)$ an element of $\PSL(2,\CC)$. This lifts to two matrices, say $A_1$ and $-A_1$ in $\SL(2,\CC)$. Choose $\widetilde{\rho}(g_1)=A_1$, and $\tilde{\rho}'(g_1)=-A_1$. Now consider the arc $a_1$ of the knot diagram corresponding to the generator $g_1$. Follow that arc to its endpoint. This will be an undercrossing, with an associated relator $g_1=g_j^{\pm 1}g_kg_j^{\mp 1}$. Thus there is a unique choice for $\widetilde{\rho}(g_k)$ and $\tilde{\rho}'(g_k)$. Continue along the arc associated with $g_k$, obtaining a unique choice for the representations at its endpoint. Continuing in this manner, we traverse the entire knot, and obtain a unique choice for each generator, obtaining two well-defined lifts $\widetilde{\rho}$ and $\tilde{\rho}'$. 
\end{proof}

In fact, it can be shown that for any complete hyperbolic manifold $M$ the holonomy\index{holonomy} representation always lifts to a representation into $\SL(2,\CC)$ \cite{CullerShalen:Varieties}, but we only present the result above on knot complements.

We now consider representations of a group $G$ into $\SL(2,\CC)$.

Suppose $G$ has finite presentation $\langle g_1, \dots, g_n\mid r_1, \dots , r_m\rangle$. Then a representation $\rho\from G \to \SL(2,\CC)$ is determined by $(\rho(g_1), \dots, \rho(g_n))$, each of which is a matrix in $\SL(2,\CC)$:
\[ \rho(g_i) = \mat{a_i& b_i \\ c_i & d_i}.\]
We therefore may view the space of representations $R_{\SL}(G)$ as a subset of the complex space $\CC^{4n}$.

Any representation of $G$ must satisfy the relators. Each relator $r_i$ is a word in the generators of $G$. Let $r_i(\rho(g_1), \dots, \rho(g_n))$ be the word with $\rho(g_j)^{\pm 1}$ substituted for each $g_j^{\pm 1}$ in $r_i$. Then \[(a_1, b_1, c_1, d_1, \dots, a_n, b_n, c_n, d_n)\in\CC^{4n}\]
lies in $R_{\SL}(G)$ if and only if the following hold:
\begin{equation}\label{Eqn:Determinant}
  a_id_i-b_ic_i=1 \quad \mbox{ for } i=1,\dots,n ,
\end{equation}
\begin{equation}\label{Eqn:Relator}
  r_j\left( \mat{a_1&b_1\\c_1&d_1}, \dots, \mat{a_n&b_n\\c_n&d_n} \right) =  \mat{1&0\\0&1} \quad \mbox{ for }j=1,\dots,m.
\end{equation}
The equations in \eqref{Eqn:Determinant} come from the determinant condition, to ensure each matrix $\mat{a_i&b_i\\c_i&d_i}$ lies in $\SL(2,\CC)$. 

The equations in \eqref{Eqn:Relator} come from the relators; we claim that each such equation gives four polynomial equations in the variables $\{a_i, b_i, c_i, d_i\}_{i=1}^n$. To see this, in the word on the left of \eqref{Eqn:Relator}, replace each instance of $\mat{a_i&b_i\\c_i&d_i}^{-1}$ with $\mat{d_i&-b_i\\-c_i&a_i}$. Then matrix multiplication, followed by setting each position in the matrix on the left of \eqref{Eqn:Relator} equal to the corresponding position in the identity matrix on the right of \eqref{Eqn:Relator} gives four polynomial equations for each relator.

\begin{proposition}\label{Prop:RepSpaceisAlgSet}
The space $R_{\SL}(G)$ consisting of representations of $G$ into $\SL(2,\CC)$ is an affine algebraic set.\index{affine algebraic set}
\end{proposition}

\begin{proof}
By the above discussion, a point in $\CC^{4n}$ lies in $R_{\SL}(G)$ if and only if it satisfies the polynomial equations coming from \eqref{Eqn:Determinant}, and the polynomial equations coming from each entry of each matrix equation \eqref{Eqn:Relator}. 
\end{proof}

\begin{definition}\label{Def:RepresentationVariety}
For $G$ a finitely presented group, its \emph{representation variety}\index{representation variety} is the affine algebraic set $R_{\SL}(G)$. 
\end{definition}

Notice that a different presentation of $G$ will have different relators, and hence different polynomial equations. However, we will see that the corresponding affine algebraic sets are not very different in this case, which we formalize with the following definition. 

\begin{definition}\label{Def:Morphism}
A \emph{polynomial map}\index{polynomial map} is a map completely defined by polynomials. A \emph{morphism}\index{affine algebraic set!morphism}\index{morphism of affine algebraic sets} between affine algebraic sets is a polynomial map from one to the other. It is an isomorphism if it is bijective and its inverse function is also a morphism. 
\end{definition}

\begin{proposition}\label{Prop:DifferentPresentation}
Suppose $G$ has presentations
\[ G = \langle g_1, \dots, g_n \mid r_1, \dots, r_m \rangle =
\langle h_1, \dots, h_k \mid s_1, \dots, s_\ell\rangle. \]
Then the affine algebraic sets $R_{\SL}^1(G)$ and $R_{\SL}^2(G)$ corresponding to the two presentations are isomorphic. 
\end{proposition}

\begin{proof}
We may write each $h_i$ as a word in the generating set $g_1, \dots, g_n$; denote this by $h_i=w_i(g_1, \dots, g_n)$. Define a map $\phi\from R_{\SL}^1(G)\to R_{\SL}^2(G)$ by setting $\phi(\rho)$ to be the representation defined on the generators $h_1,\dots,h_k$ by $\phi(\rho)(h_i) = w_i(\rho(g_1), \dots, \rho(g_n))$. This map only involves multiplication and addition of matrix coordinates, hence it is completely defined by polynomials, and hence the map is a morphism\index{affine algebraic set!morphism} of affine algebraic sets. By symmetry, we also have a polynomial map $\psi\from R_{\SL}^2(G)\to R_{\SL}^1(G)$ defined similarly. Then $\phi$ and $\psi$ are inverses; these are isomorphisms of affine algebraic sets. 
\end{proof}

\begin{example}\label{Example:Abelian}
The space $R_{\SL}(\pi_1(S^3-K))$ for any knot $K$ will always contain a simple irreducible affine algebraic set,\index{irreducible affine algebraic set}\index{affine algebraic set} which we now describe.

Recall that the abelianization of the fundamental group of any knot complement is $\ZZ = H_1(S^3-K)$, generated by the (homology class of the) meridian (\refcor{Homology}). Thus any representation of $\ZZ$ into $\SL(2,\CC)$ gives an induced representation of $\pi_1(S^3-K)$ into $\SL(2,\CC)$ via
\[ \pi_1(S^3-K) \to \ZZ \to \SL(2,\CC).\]
Any such representation is called an \emph{abelian representation}\index{abelian representation}.

Note that the standard longitude $\lambda$ is trivial in $H_1(S^3-K)$. Thus for any abelian representation $\rho$, we will have $\rho(\lambda)=\Id$. On the other hand, the meridian $\mu$ can be mapped to any element of $\SL(2,\CC)$. 

We claim that abelian representations form an affine algebraic set as a subset of $R_{\SL}(\pi_1(S^3-K))$, and that this irreducible component\index{irreducible component} is isomorphic to $\SL(2,\CC)$. We leave this as an exercise (\refex{AbelianReps}).
\end{example}

We only wish to consider representations into $\SL(2,\CC)$ up to conjugation. Changing the basepoint of $G$ will change the fundamental group by conjugation, and hence will change any representation by conjugation.
Thus we really want to consider two representations to be the same if they are conjugate.

A first idea would be to declare representations
$\rho_1, \rho_2\from G\to\SL(2,\CC)$ to be equivalent if they are conjugate.
That is, $\SL(2,\CC)$ acts on $R_{\SL}(G)$ via conjugation as follows. If $\rho\in R_{\SL}(G)$ and $A\in\SL(2,\CC)$, then $A\cdot\rho = i_A\circ\rho$ where $i_A$ is defined by $i_A(X)=AXA^{-1}$. We could consider the orbits of this action, which is a polynomial map from $\SL(2,\CC)\times R_{\SL}(G)$ to $R_{\SL}(G)$, and take a quotient. 
However, a point may be in the closure of several orbits, so this can lead to a non-Hausdorff space. The way to fix this is to take what is sometimes called the algebro-geometric quotient, or the algebraic quotient of invariant theory. We will not go into the details of the construction here. Instead, we give a definition of the algebro-geometric quotient that is known to be equivalent.

\begin{definition}\label{Def:Character}
  The \emph{character}\index{character!representation} of a representation $\rho\from G\to\SL(2,\CC)$ is the function $\chi_\rho\from G\to\CC$ given by $\chi_\rho(g) = \tr(\rho(g))$, where $\tr$ denotes the function that takes the trace of a matrix. 
\end{definition}

Note the character is the same for two conjugate representations. It is not quite true that representations with the same character are necessarily conjugate.

\begin{definition}\label{Def:CharacterVariety}
The \emph{character variety}\index{character variety} $X_{\SL}(G)$ is the space of characters of all representations $R_{\SL}(G)$.
\end{definition}

It is shown in \cite{CullerShalen:Varieties} that the character variety is an affine algebraic set; a more elementary proof of this fact is also given in \cite{GonzalezAcuna-Montesinos}. We refer you to those papers for the details. We note again that the standard terminology is a little misleading. In \refdef{ReducibleAlgebraicSet}, we noted that an affine variety\index{affine variety} is an irreducible affine algebraic set.\index{irreducible affine algebraic set} But the character variety\index{character variety} is typically reducible, hence not an affine variety at all. One true variety that is frequently studied is the irreducible component\index{irreducible component} of the character variety corresponding to a discrete faithful representation, i.e.\ the representation coming from a complete hyperbolic structure. This is called the \emph{canonical component}\index{canonical component}\index{character variety!canonical component}.

The canonical component of a knot complement contains a great deal of information about the manifold, including information on Dehn filling and surfaces embedded in the 3-manifold. However, it can be very difficult to compute the character variety.\index{character variety} At the time of writing this book, character varieties have been computed only for simple families of 3-manifolds and knot complements, including double twist knots and links, and some 2-bridge links \cite{MacasiebPetersenvanLuijk}, \cite{PetersenTran:TwistLinks}.

\subsection{The case of $\PSL$}

We spent the last subsection considering representations to $\SL(2,\CC)$. Much of that work can be extended to $\PSL(2,\CC)$ with a little extra effort. The extension to $\PSL(2,\CC)$ has been done carefully by \cite{BoyerZhang}. See also \cite{HeusenerPorti} for a nice exposition.

Let $G$ be a finitely presented group:
\[G = \langle g_1, \dots, g_n \mid r_1, \dots, r_m\rangle. \]
Let $R_{\PSL}(G)$ denote the space of representations from $G$ to $\PSL(2,\CC)$. 

The first thing to check is that $R_{\PSL}$ forms an algebraic set. As before, relators give polynomial equations, but now these are only defined up to multiplication by $\pm 1$. An easy way to avoid sign issues is to use the fact that the group $\PSL(2,\CC)$ is isomorphic to $\SO(3,\CC)$; this can be shown using standard tools from Lie groups, and we leave this to the reader (or see \cite[Lemma~2.1]{HeusenerPorti}).

\begin{proposition}\label{Prop:PSLRepSpaceisAlg}
The space $R_{\PSL}$ consisting of representations of $G$ into $\PSL(2,\CC)$ is an affine algebraic set.\index{affine algebraic set}
\end{proposition}

\begin{proof}
Using the fact that $\PSL(2,\CC)\cong \SO(3,\CC)$, this follows as in the proof of \refprop{RepSpaceisAlgSet}. 
\end{proof}

\begin{definition}\label{Def:PSLRepVariety}
Let $G$ be a finitely presented group. Define the \emph{$\PSL$-representation variety}\index{representation variety!$\PSL$}\index{$\PSL$-representation variety} of $G$ to be the algebraic set $R_{\PSL}(G)$. 
\end{definition}

Again we only wish to consider representations into $\PSL(2,\CC)$ up to conjugation. In the $\PSL$ case, there is a natural geometric reason: conjugate holonomy\index{holonomy} representations give isometric hyperbolic structures on a 3-manifold. We want an algebraic theory that views these structures as the same.
As in the $\SL$ case, we may take the algebro-geometric quotient, or alternatively, take the set of $\PSL$-characters, defined below. \cite{HeusenerPorti} give a proof that this is equivalent. 

\begin{definition}\label{Def:PSLCharacter}
  Let $\rho\from G\to\PSL(2\CC)$ be a representation. 
The \emph{$\PSL$-character}\index{character!$\PSL$} of  $\rho$ is the function $\xi_\rho\from G\to \CC$ given by $\xi_{\rho}(g) = \tr(\rho(g))^2$, where $\tr$ denotes the trace of the matrix. Note that the trace of an element of $\PSL(2,\CC)$ is only defined up to sign, but the square of the trace is well-defined.
\end{definition}

\begin{definition}\label{Def:PSLCharacterVariety}
The \emph{$\PSL$-character variety}\index{character variety!$\PSL$}\index{$\PSL$-character variety} $X_{\PSL}(G)$ is the space of $\PSL$-characters of all representations $R_{\PSL}(G)$. 
\end{definition}

%%%%%%%%%%%%%%%%%%%%%%%%%%%%%%%%%%%%%%%%%%%%%%%%%%%%%%%%%%%%%%%%%

\section{The $A$-polynomial}
The $A$-polynomial was first defined by Cooper, Culler, Gillet, Long, and Shalen in \cite{APoly}. The authors noticed that for manifolds with a single torus boundary component, the difficult problem of computing the $\SL(2,\CC)$ character variety could be replaced by a slightly simpler problem. Rather than considering the entire fundamental group of the knot complement, we focus on the subgroup corresponding to generators of the boundary of a regular neighborhood of the knot, $\bdy N(K)$. This reduces a complicated algebraic set to a single polynomial, called the $A$-polynomial, which still encodes a great deal of geometric information.
Our exposition in this section is based on that of Cooper and Long in \cite{CooperLong:Apoly}, with thanks to Mathews \cite{Mathews:Honours}.

\subsection{Classical definition}

Denote the meridian and longitude of a knot by $\mu$ and $\lambda$, respectively. We will be considering representations of the fundamental group $\pi_1(\bdy N(K))$ into $\SL(2,\CC)$ up to conjugation.
Since $\mu$ and $\lambda$ commute, for any representation $\rho\from \pi_1(\bdy N(K))\to\SL(2,\CC)$, the matrices $\rho(\mu)$ and $\rho(\lambda)$ must have the same fixed points on $\bdy \HH^3$ (\refex{PSL(2,C)Commute}). Hence they can be conjugated to fix either the point at infinity or the geodesic from $0$ to $\infty$ in $\bdy\HH^3$. In the former case, the two matrices have the form
$\mat{1 & * \\ 0 & 1}$, and in the latter, the matrices are diagonal. Thus there will be no loss of generality in restricting to representations for which $\rho(\mu)$ and $\rho(\lambda)$ are upper triangular.

Define $R^U_{\SL}(\pi_1(S^3-K))\subset R_{\SL}(\pi_1(S^3-K))$ to be those representations $\rho$ for which $\rho(\mu)$ and $\rho(\lambda)$ are both upper triangular. When the context is clear, we abbreviate to $R^U$. 

\begin{lemma}
The set $R^U_{\SL}(\pi_1(S^3-K))$ forms an affine algebraic set.\index{affine algebraic set}
\end{lemma}

\begin{proof}
We need to show that $R^U_{\SL}$, just like $R_{\SL}$, is defined as the set of zeros of a system of polynomial equations.

Recall that $\rho \in R_{\SL}(\pi_1(S^3-K))$ is determined by $\rho(g_1), \dots, \rho(g_n)$, where the $g_i$ are generators, and we view each $\rho(g_i)$ as a matrix with coordinates $a_i, b_i, c_i, d_i\in\CC$ satisfying \refeqn{Determinant} and \refeqn{Relator}. Then $\rho(\mu)$ and $\rho(\lambda)$ will be matrices whose coefficients are polynomials in the $\{a_i, b_i, c_i, d_i\}$. In particular, the lower left entries of $\rho(\mu)$ and $\rho(\lambda)$ are polynomials $q_\mu$ and $q_\lambda$, respectively. Thus $R^U_{\SL}(\pi_1(S^3-K))$ is obtained by adjoining $q_\mu=0$ and $q_\lambda=0$ to the polynomials defining $R_{\SL}(\pi_1(S^3-K))$. Hence $R^U_{\SL}$ is an algebraic set.
\end{proof}

Next, because the determinants of $\rho(\mu)$ and $\rho(\lambda)$ are equal to $1$, it must be the case that for some $m,\ell\in\CC^2$, the matrices $\rho(\mu)$ and $\rho(\lambda)$ have the form
\[ \rho(\mu) = \mat{m&* \\0&m^{-1}}, \quad \rho(\lambda) = \mat{\ell&*\\0&\ell^{-1}}.\]

Define functions $\xi_\mu, \xi_\lambda\from R^U_{\SL} \to \CC$ by letting $\xi_\mu(\rho)$ (respectively $\xi_\lambda(\rho)$) be the upper left entry of the matrix $\rho(\mu)$ (respectively $\rho(\lambda)$). That is, $\xi_\mu(\rho)=m$ and $\xi_\lambda(\rho)=\ell$. Each of these entries can be written as a polynomial in the $a_i, b_i, c_i, d_i$, so each map $\xi_\mu$ and $\xi_\lambda$ is a polynomial map, and thus the map $\xi\from R^U_{\SL}\to \CC^2$ given by $\xi(\rho)= (\xi_\lambda(\rho), \xi_\mu(\rho)) = (\ell, m)$ is a morphism. \index{affine algebraic set!morphism}

Consider the image $\xi(R^U_{\SL})\subset \CC^2$. This is the image of an algebraic set under a morphism,\index{affine algebraic set!morphism} hence is an algebraic set. It may have several irreducible components.\index{irreducible component} Let $C$ be an irreducible component of $R^U_{\SL}$ and let $\overline{\xi(C)}$ be its closure (in the Zariski topology).
The closure $\overline{\xi(C)}$ is the set of zeros of a family of polynomials. If there is a single such polynomial, then denote it by $F_C$.

\begin{example}\label{Example:AbelianApoly}
Let $C \subset R^U_{\SL}$ be the component of abelian representations,\index{abelian representation} as in \refexamp{Abelian}. Then $\rho(\mu)$ is an upper triangular matrix in $\SL(2,\CC)$, but $\rho(\lambda)=\Id$ always. Thus $\xi_\mu(\rho)$ can be arbitrary, but $\xi_\lambda(\rho)=1$. Thus the closure of $\xi(C)$ for $C$ the component of abelian representations gives the polynomial $F_C = \ell-1$. 
\end{example}

Now consider all polynomials $F_C$ arising from irreducible components\index{irreducible component} of $R^U_{\SL}$, aside from the abelian representation. 
Our preliminary definition of the $A$-polynomial is the polynomial obtained by multiplying all of these polynomials, or $1$ if there are no such polynomials. Observe that since $F_C$ is the set of zeros of a polynomial, it is only defined up to multiplication by a scalar and by powers of $m$ and $\ell$. It was shown in \cite{APoly} that a scalar multiple can be chosen so that all coefficients are integers. We require integer coefficients with no common factors, and multiply by $\ell^am^b$ for $a,b$ integers such that the total degree is minimal. This defines the $A$-polynomial up to sign. We summarize:

\begin{definition}\label{Def:Apoly}
Let $K$ be a knot. For every $C$ an irreducible component\index{irreducible component} of $R^U_{\SL}(\pi_1(S^3-K))$ that is not abelian, and such that $\overline{\xi(C)}$ is the set of zeros of a single polynomial, let  $F_C$ denote this polynomial. 
The \emph{$A$-polynomial}\index{$A$-polynomial} (or $\SL(2,\CC)$ $A$-polynomial)\index{$\SL(2,\CC)$ $A$-polynomial}\index{$A$-polynomial!$\SL(2,\CC)$} of $K$ is defined to be the product of $1$ and any polynomials $F_C$ as above, rescaled so that all coefficients are integers with no common factors. We normalize by multiplying by $\ell^am^b$ so that the total degree is minimized, and denote the result by $A_{\SL}(\ell,m)$. 
\end{definition}

It is easiest to make sense of the definitions if we work with examples. The simplest example is given by the following. 

\begin{proposition}\label{Prop:ApolyTrivial}
The $A$-polynomial of the unknot is $1$.\index{$A$-polynomial}
\end{proposition}

\begin{proof}
The fundamental group $G$ of the unknot is $\ZZ$, hence any representation of $G$ into $\SL(2,\CC)$ is an abelian representation.\index{abelian representation} Thus in \refdef{Apoly}, there are no polynomials $F_C$ and the $A$-polynomial is $1$. 
\end{proof}

\begin{example}[Figure-8 knot]\label{Example:Fig8_APoly}
We use the Wirtinger presentation of \refexamp{Fig8Wirtinger} to compute the polynomial $A_{SL}$ for the figure-8 knot.\index{$A$-polynomial}

There are two generators, $g_1$ and $g_2$, satisfying
\[ \rho(g_1) = \mat{a_1 & b_1 \\ c_1 & d_1}, \quad \rho(g_2) = \mat{a_2 & b_2 \\ c_2 & d_2}. \]
This initially gives the eight unknowns $a_i, b_i, c_i, d_i$ for $i=1, 2$. These must satisfy the following equations:
\begin{enumerate}
\item Two determinant equations $a_id_i-b_ic_i=1$.
\item Four equations coming from the four matrix entries of the relation
  \[ \rho(g_1g_2^{-1}g_1g_2g_1^{-1}g_2) = \rho(g_2^{-1}g_1g_2g_1^{-1}g_2g_2). \]
\item Equations coming from the matrix entries of $\rho(M)$ and $\rho(L)$:
  \[ \rho(L)= \rho(g_1^{-1}g_2 g_1g_2^{-1}g_1^{-1} g_2^{-1}g_1g_2g_1^{-1}g_2) \quad \mbox{ and } \quad \rho(M)=\rho(g_1). \]
  There are two equations to ensure $\rho(M)$ and $\rho(L)$ are upper triangular, by setting their bottom left entry to be $0$. In particular, this requires the simple equation $c_1=0$, coming from $\rho(M)$. A more complicated equation will come from the bottom left entry of $\rho(L)$, since $L$ is a more complicated product of $g_1^{\pm 1}, g_2^{\pm 1}$. 
  
  There is one equation ensuring the top left entry of $\rho(L) = \ell$, and one equation ensuring the top left entry of $\rho(M)=m$. In particular, again in the simpler $\rho(M)$ case this gives $a_1=m$. Using one of the determinant equations, we also conclude $d_1=m^{-1}$.
\end{enumerate}
In all, this gives seven equations in the unknowns $\ell$, $m$, $b_1$, $a_2$, $b_2$, $c_2$, $d_2$.

Using computer software such as Mathematica, we find a polynomial describing the curve satisfying this system:
\[ -\ell + \ell^2 +\ell m^2 - \ell^2m^2 + m^4 + \ell m^4 - \ell^2 m^4 - \ell^3m^4 + \ell m^6 - \ell^2m^6 - \ell m^8 + \ell^2m^8.  \]

Note that $(\ell-1)$ divides this polynomial, corresponding to the abelian representation.\index{abelian representation} We divide out by $(\ell-1)$, to obtain\index{$A$-polynomial}\index{$\SL(2,\CC)$ $A$-polynomial}\index{$A$-polynomial!$\SL(2,\CC)$}
\[ A_{\SL}(\ell, m) = \ell - \ell m^2 - m^4 - 2\ell m^4 - \ell^2 m^4 - \ell m^6 + \ell m^8. \]
\end{example}

\subsection{The $\APSL$-polynomial}

The work in the previous section can also be applied to representations of $\pi_1(S^3-K)$ into $\PSL(2,\CC)$ rather than $\SL(2,\CC)$. That is, for fixed generators $\mu$ and $\lambda$ of $\pi_1(\bdy N(K))$ restrict to $R^U_{\PSL} \subset R_{\PSL}$, the subset of representations $\rho$ for which $\rho(\mu)$ and $\rho(\lambda)$ are upper triangular in $\PSL(2,\CC)$. Define $\xi_{\PSL}\from R^U_{\PSL}(\pi_1(S^3-K))\to\CC^2$ by
\[ \xi_{\PSL}(\rho) = (\xi_\lambda(\rho), \xi_\mu(\rho)) = (\ell^2, m^2), \]
where $\xi_\lambda$ gives the square of the top left entry of $\rho(\lambda)$, and similarly for $\xi_\mu$.  Consider each irreducible component\index{irreducible component} $C$ of $R^U_{\PSL}$ such that the closure $\overline{\xi_{\PSL}(C)}$ is defined by a single polynomial $F_{\PSL}(C)$, and $C$ is not abelian.

\begin{definition}
The \emph{$\APSL$-polynomial}\index{$\APSL$-polynomial}\index{$A$-polynomial!$\APSL$} is defined to be the product of the polynomials $F_{\PSL}(C)$ as above, or $1$ if there are no such polynomials. 
\end{definition}

The $\APSL$-polynomial and the $A$-polynomial are related, and we describe the relationship in the hyperbolic case. 
When $K$ is hyperbolic, with $G=\pi_1(S^3-K)$, let $X_{\PSL}^0(G)$ be the irreducible component\index{irreducible component} of $X_{\PSL}(G)$ containing the complete hyperbolic structure. Similarly, let $X_{\SL}^0(G)$ be the irreducible component of $X_{\SL}(G)$ containing the lift of the complete hyperbolic structure. Then we may define polynomials $A^0_{\PSL}$ and $A^0_{\SL}$ by only considering $C$ coming from these irreducible components. Note $A^0_{\PSL}$ is a factor of the $\APSL$-polynomial, and $A^0_{\SL}$ is a factor of the $A$-polynomial.

\begin{proposition}
The polynomial $A_{\PSL}^0(\ell^2, m^2)$ divides the polynomial\index{$A$-polynomial}
\[ A^0_{\SL}(\ell, m)A^0_{\SL}(\ell, -m)A^0_{\SL}(-\ell, m)A^0_{\SL}(-\ell, -m). \]
\end{proposition}

\begin{proof}
The projection $\pi\from \SL(2,\CC) \to \PSL(2,\CC)$ induces a map on character varieties\index{character variety}
$\pi\from X_{\SL}^0(G) \to X_{\PSL}^0(G)$, which is surjective by \cite{Culler:Lifting}. Define $h\from \CC^2 \to \CC^2$ by $h(x,y) = (x^2, y^2)$. Then the following diagram commutes.
\[
\begin{tikzcd}[column sep=50pt]
  X_{\SL}^0(\pi_1(S^3-N(K)))
    \arrow[r, "r"] 
    \arrow[d, "\pi"]
  &
  X_{\SL}^0(\pi_1(\bdy N(K)))
    \arrow[d, "\pi"]
    \arrow[r, "\xi^{-1} \, \circ \, \tr", leftarrow]
  & (\CC^*)^2
    \arrow[d, "h"] \\
  X_{\PSL}^0(\pi_1(S^3-N(K)))
    \arrow[r, "r"]
  & 
  X^0_{\PSL}(\pi_1(\bdy N(K)))
    \arrow[r, "\xi_{\PSL}^{-1} \circ \tr", leftarrow]
  & (\CC^*)^2 
\end{tikzcd}
\]
Here $r$ is a map induced by the inclusion of $\pi_1(\bdy N(K))$ into $\pi_1(S^3-N(K))$.

If $D^0_{\SL}$ and $D^0_{\PSL}$ are curves defined by $A^0_{\SL}(\ell, m)$ and $A^0_{\PSL}(\ell, m)$, respectively, then $h(D^0_{\SL}) = D^0_{\PSL}$. So $h^{-1}(D^0_{\PSL})$ is the union of the curves given by $A^0_{\SL}(\pm\ell, \pm m)$. Thus $A^0_{\PSL}(\ell^2, m^2)$ divides
\[ A^0_{\SL}(\ell, m)A^0_{\SL}(\ell, -m)A^0_{\SL}(-\ell, m)A^0_{\SL}(-\ell, -m).\qedhere
\]
\end{proof}

\subsection{Relation to $\Ahyp$-polymomial}

Let $\{z_1, \dots, z_n\}$ be parameters for a triangulation coming from a point in the gluing variety. Then we obtain a well-defined developing map by attaching tetrahedra in $\HH^3$ via a face-pairing isometry. This determines a holonomy\index{holonomy} representation, hence we obtain a map $D\from \calD{\calT} \to X_{\PSL}(G)$. Champanerkar showed that this map is algebraic, i.e.\ a morphism,\index{affine algebraic set!morphism} and thus $\Ahyp(\ell, m)$ divides $\APSL(\ell,m)$ \cite{Champanerkar:Thesis}.

It follows that $A^0_{\mathrm{Hyp}}(\ell,m) = A_{\PSL}^0(\ell,m)$. 

\begin{example}
We computed above in \refexamp{Fig8_AHypPoly} that for the figure-8 knot complement,
\begin{multline*}
\Ahyp(\ell,m) = \ell - 2 \ell m - 3 \ell m^2 + 2 \ell m^3 - m^4 + 6 \ell m^4  \\ - \ell^2 m^4 + 2 \ell m^5 - 3 \ell m^6 - 2 \ell m^7 + \ell m^8.
\end{multline*}
And in \refexamp{Fig8_APoly}, its $A$-polynomial is\index{$A$-polynomial}
\[ A_{\SL}(\ell, m) = \ell - \ell m^2 - m^4 - 2\ell m^4 - \ell^2 m^4 - \ell m^6 + \ell m^8. \]

In this case, $\Ahyp(\ell^2,m^2) = -A_{\SL}(\ell,m)A_{\SL}(-\ell,m)$.\index{$\Ahyp$-polynomial}\index{$\SL(2,\CC)$ $A$-polynomial}\index{$A$-polynomial!hyperbolic, $\Ahyp$}\index{$A$-polynomial!$\SL(2,\CC)$}
\end{example}

%%%%%%%%%%%%%%%%%%%%%%%%%%%%%%%%%%%%%%%%%%%%%%%%%%%%%%%%%%%%%%%%%

\section{Exercises}

\begin{exercise}\label{Ex:Zariski}[Zariski topology]
\begin{enumerate}
\item Prove that the union of two affine algebraic sets\index{affine algebraic set} and the intersection of arbitrarily many affine algebraic sets are both affine algebraic sets.
\item Describe the following sets as affine algebraic sets: $\CC^N$, $\emptyset$, a single point $(a_1, \dots, a_N)$.
\item
The \emph{Zariski topology}\index{Zariski topology} is the topology on $\CC^N$ formed by taking affine algebraic sets to be closed sets. 
Prove that any Zariski--closed set is closed in the standard Euclidean topology. Find an example of a set that is Euclidean--closed but not Zariski--closed.
\end{enumerate}
\end{exercise}

\begin{exercise}
SnapPy finds triangulations of knot complements and allows you to print out gluing and completeness equations.
\begin{enumerate}
\item Using SnapPy, find a set of polynomial equations that lead to the $\Ahyp$ polynomials of the knots $5_2$, $6_1$, and $6_2$.
\item Using Mathematica (or other), compute the $\Ahyp$ polynomials of these knots.
\item Use SnapPy to randomize the triangulations of the knots. What happens to the $\Ahyp$ polynomial?
\item Use SnapPy to change the generators of $\pi_1(\bdy N(K))$, i.e.\ the curves giving the completeness equations. What happens to the $\Ahyp$ polynomial?
\end{enumerate}
\end{exercise}

\begin{exercise}\label{Ex:SVK}
Using the Siefert--Van Kampen theorem, or otherwise, prove that the group given by the Wirtinger presentation\index{Wirtinger presentation} is isomorphic to the fundamental group of the knot complement. 
\end{exercise}

\begin{exercise}
  Find the Wirtinger presentation\index{Wirtinger presentation} of the fundamental group of the $5_2$, $6_1$, and $6_2$ knots. Also find a presentation for their longitudes. 
\end{exercise}

\begin{exercise}
  Let $r_1, r_2, \dots, r_n$ denote the relators of the Wirtinger presentation\index{Wirtinger presentation} coming from the $n$ distinct crossings of a knot. Show that $r_n$ is always redundant.
\end{exercise}

\begin{exercise}[Abelian representations]\label{Ex:AbelianReps}
Show that abelian representations\index{abelian representation} form an affine algebraic set\index{affine algebraic set} isomorphic to $\SL(2,\CC)$. (Ensure your isomorphism is an isomorphism of affine algebraic sets,\index{affine algebraic set!morphism} i.e.\ defined by polynomial maps.)
\end{exercise}

\begin{exercise}
  Compute $A_{\SL}(\ell, m)$ for the $5_2$, $6_1$, or $6_2$ knot.
\end{exercise}

\begin{exercise}
  The \emph{support} of a polynomial $F(x,y)$ is the set $\{(a,b)\}\subset \ZZ^2$ such that the coefficient of the term $x^ay^b$ in $F(x,y)$ is nonzero. The convex hull of the support is the \emph{Newton polygon}\index{Newton polygon} of $F$. The Newton polygon of the $A$-polynomial has a remarkable relationship with essential surfaces embedded in the knot:

  \begin{theorem}[\cite{APoly}]
    Let $K$ be a knot in $S^3$ with $A$-polynomial $A_{\SL}(\ell,m)$. 
    Suppose the Newton polygon of $A_{\SL}(\ell,m)$ has a side of slope $p/q$. Then there is an essential surface\index{essential} $S$ in $S^3-N(K)$ with boundary that is a curve on the torus $\bdy N(K)$ with slope $p/q \in H_1(\bdy N(K))$.
  \end{theorem}

  Champanerkar proved that the corresponding result holds for $\Ahyp$.

  \begin{enumerate}
  \item Compute the Newton polygon for $A_{\SL}(\ell,m)$ for the figure-8 knot.
  \item Compute the Newton polygon for $\Ahyp(\ell,m)$ for the figure-8 knot.
  \end{enumerate}
\end{exercise}

\begin{exercise}
  Repeat the previous exercise for the $5_2$, $6_1$, or $6_2$ knot. 
\end{exercise}

\backmatter

\bibliographystyle{amsalpha}
\bibliography{biblio}

\printindex

\end{document}